\documentclass[11pt, a4paper]{article}
\title{Singular Contour Dynamics and Paradifferential Reduction}
\newif\ifhasstmaryrd
\IfFileExists{stmaryrd.sty}{\hasstmaryrdtrue\usepackage{stmaryrd}}{\hasstmaryrdfalse}
\author{
Xingyu Li \thanks{Università degli Studi di Trieste, Dipartimento di Matematica, Informatica e Geoscienze, Italy. \href{mailto:xingyuli92@gmail.com}{\texttt{xingyuli92@gmail.com}}}
    \and
    Emeric Roulley \thanks{Universit\`a Statale di Milano, Via Saldini 50 Milano, Italy.  \href{mailto:emeric.roulley@unimi.it}
    {\texttt{emeric.roulley@unimi.it}}}
        \and
    Stefano Scrobogna\thanks{Università degli Studi di Trieste, Dipartimento di Matematica, Informatica e Geoscienze, Italy. \href{mailto:stefano.scrobogna@units.it}
    {\texttt{stefano.scrobogna@units.it}}}
}

\usepackage{geometry}
\usepackage{empheq}
\usepackage{times}
\IfFileExists{fourier.sty}{\usepackage{fourier}}{}
\usepackage[T1]{fontenc}
\usepackage{amssymb, bm, mathrsfs, amsmath}
\usepackage{amsfonts}
\usepackage[english]{babel}
\usepackage{amssymb,amsthm}
\usepackage{mathtools}
\IfFileExists{extarrows.sty}{%
  \usepackage{extarrows}%
}{%
  \newcommand{\xlongequal}[2][]{\mathrel{\overset{##2}{=}}}%
}
\usepackage{cases}
\DeclareMathAlphabet{\mathcal}{OMS}{cmsy}{m}{n}
\DeclareMathAlphabet{\mathpzc}{OT1}{pzc}{m}{it}

\makeatletter
\newcommand{\mathpzcbold}[1]{%
  \leavevmode
  \setbox0=\hbox{$\mathpzc{#1}$}%
  \kern-.03em\copy0\kern-\wd0
  \kern.03em\copy0\kern-\wd0
  \kern0em\box0
}
\makeatother

\usepackage[hypertexnames=false]{hyperref}
\usepackage[capitalise,nameinlink,sort&compress]{cleveref}
\usepackage{cite}
\allowdisplaybreaks[1]
\usepackage{xfrac}
\usepackage[utf8]{inputenc}
\usepackage[T1]{fontenc}
\usepackage{enumerate}
\IfFileExists{accents.sty}{\usepackage{accents}}{}
\IfFileExists{bbm.sty}{\usepackage{bbm}}{}
\usepackage{bm}
\usepackage{crossreftools}
\usepackage{autonum}
\makeatletter
\autonum@generatePatchedReferenceCSL{Cref}
\makeatother

\usepackage{todonotes}

\usepackage{xcolor}

\usepackage{tikz}
\usetikzlibrary{arrows, automata, backgrounds, shadows, patterns, calc, hobby, plotmarks, shapes}
\usetikzlibrary{shapes.misc}

\tikzset{cross/.style={cross out, draw=black, minimum size=2*(#1-\pgflinewidth), inner sep=0pt, outer sep=0pt},
cross/.default={1pt}}
\allowdisplaybreaks

\DeclareSymbolFontAlphabet{\amsmathbb}{AMSb}%

\makeatletter

\newcommand{\Rmnum}[1]{\expandafter\@slowromancap\romannumeral #1@}
\makeatother

\usepackage{listings}
\usepackage{color}

\definecolor{dkgreen}{rgb}{0,0.6,0}
\definecolor{gray}{rgb}{0.5,0.5,0.5}
\definecolor{mauve}{rgb}{0.58,0,0.82}

\makeatletter
\def\maketag@@@#1{\hbox{\m@th\normalfont\normalsize#1}}
\makeatother

\newcommand{\dd}{\textnormal{d}}
\DeclareMathOperator{\dive}{div}
\DeclareMathOperator{\tr}{tr}

\DeclareMathOperator{\sgn}{sgn}
\newcommand{\QL}{\textnormal{QL}}
\newcommand{\Atw}{\mathfrak{a}}
\newcommand{\lin}{\rm lin}
\newcommand{\bil}{\rm bil}

\newcommand{\fine}{\hfill$\blacklozenge$}
\newcommand{\pare}[1]{\left( #1 \right)}
\newcommand{\bra}[1]{\left[ #1 \right]}
\newcommand{\bbra}[1]{\left\llbracket #1 \right\rrbracket}

\newcommand{\angles}[1]{\left\langle #1 \right\rangle}
\newcommand{\aangle}[1]{\left\langle\!\!\!\left\langle #1 \right\rangle\!\!\!\right\rangle}

\newcommand{\norm}[1]{\left\| #1 \right\|}
\newcommand{\av}[1]{\left| #1 \right|}

\newcommand{\pbra}[2]{\set{ #1 \big. , \  #2 }}
\newcommand{\comm}[2]{\left\llbracket #1 \,  , \,   #2 \right\rrbracket}

\newcommand{\xfrac}[2]{\left. #1 \middle/ #2 \right. }
\newcommand{\set}[1]{\left\{ #1 \right\}}
\newcommand{\system}[1]{\left\{ #1 \right.}
\newcommand{\floor}[1]{\left\lfloor #1 \right\rfloor}

\DeclareMathOperator{\OpBWPlain}{Op^{BW}}
\newcommand{\OpBW}[1]{\OpBWPlain\left(#1\right)}
\newcommand{\RBW}[2]{R^{\textnormal{BW}} \bra{#1,#2}}

\newcommand{\OpW}[1]{\textnormal{Op}^{W}\pare{#1}}

\newcommand{\supp}{\textnormal{supp}}

\newcommand{\pv}{\textnormal{p.v.}}

\newcommand{\Id}{\textnormal{Id}}
\newcommand{\defeq}{\vcentcolon=}
\newcommand{\eqdef}{=\vcentcolon}

\newcommand{\BallR}[2]{B_{#1, \bR}\pare{I;#2}}
\newcommand{\Ball}[2]{B_{#1}\pare{I;#2}}
\newcommand{\BR}[1]{\textnormal{BR}\pare{#1}}

\newcommand{\background}{\bar{\Omega}}

\newcommand{\vect}[2]{  \begin{bmatrix} #1 \\ #2 \end{bmatrix}  }
\newcommand{\vectT}[2]{  \begin{bmatrix} #1 & #2 \end{bmatrix}  }
\newcommand{\vectD}[3]{  \begin{bmatrix} #1 \\ #2 \\ #3 \end{bmatrix}  }
\newcommand{\vectDT}[3]{  \begin{bmatrix} #1 & #2 & #3 \end{bmatrix}  }

\renewcommand{\Re}{\textnormal{Re}}

\let\Im\relax
\DeclareMathOperator{\Im}{Im}
\newcommand{\RN}[1]{%
  \textup{\uppercase\expandafter{\romannumeral#1}}%
}

\newcommand{\KF}{K\cF}
\newcommand{\KR}{K\cR}

\newcommand{\cC}{\mathcal{C}}
\newcommand{\cL}{\mathcal{L}}
\newcommand{\cF}{\mathcal{F}}
\newcommand{\cR}{\mathcal{R}}
\newcommand{\bR}{\mathbb{R}}

\newcommand{\bN}{\mathbb{N}}
\newcommand{\bS}{\mathbb{S}}
\newcommand{\bC}{\mathbb{C}}

\newcommand{\ii}{{\rm i}}

\newcommand{\psc}[2]{\left\langle  #1 \ \middle\vert \ #2 \right\rangle}

\newcommand{\loc}{\textnormal{loc}}

\ifhasstmaryrd\else
  \renewcommand{\bbra}[1]{\mathopen{[\mkern-3mu[}#1\mathclose{]\mkern-3mu]}}
  \renewcommand{\comm}[2]{\mathopen{[\mkern-3mu[}#1\, , \, #2\mathclose{]\mkern-3mu]}}
  
\fi

\numberwithin{equation}{section}

\newcommand{\p}{\pare}

\newcommand{\va}{\varepsilon}

\theoremstyle{plain}
\newtheorem{theorem}{Theorem}[section]

\newtheorem*{theorem*}{Theorem}
\newtheorem{prop}[theorem]{Proposition}

\newtheorem{lemma}[theorem]{Lemma}
\newtheorem{cor}[theorem]{Corollary}

\date{}

\theoremstyle{definition}
\newtheorem{definition}[theorem]{Definition}

\newtheorem{rem}[theorem]{Remark}
\newtheorem{remark}[theorem]{Remark}
\newtheorem{notation}[theorem]{Notation}

\newtheorem{step}{Step}
\newtheorem{substep}{Substep}[step] % <-- Produce "1.1, 1.2" e si azzera in automatico

\newtheorem{proofpart}{Part}
\newtheorem{proofsubpart}{Subpart}[proofpart] % <-- Produce "1.1, 1.2" e si azzera in automatico

\makeatletter
\@addtoreset{step}{subsubsection}
\@addtoreset{proofpart}{subsubsection}
\makeatother

\crefname{prop}{Proposition}{Propositions}
\crefname{prop}{Proposition}{Propositions}
\crefname{proposition}{Proposition}{Propositions}
\crefname{proposition}{Proposition}{Propositions}
\crefname{cor}{Corollary}{Corollaries}
\crefname{cor}{Corollary}{Corollaries}
\crefname{rem}{Remark}{Remarks}
\crefname{rem}{Remark}{Remarks}
\crefname{remark}{Remark}{Remarks}
\crefname{remark}{Remark}{Remarks}
\crefname{notation}{Notation}{Notations}
\crefname{notation}{Notation}{Notations}
\crefname{example}{Example}{Examples}
\crefname{example}{Example}{Examples}
\crefname{hyp}{Hypothesis}{Hypotheses}
\crefname{hyp}{Hypothesis}{Hypotheses}
\crefname{step}{step}{steps}
\crefname{step}{Step}{Steps}
\crefname{substep}{substep}{substeps}
\crefname{substep}{Substep}{Substeps}
\crefname{proofpart}{part}{parts}
\crefname{proofpart}{Part}{Parts}
\crefname{proofsubpart}{subpart}{subparts}
\crefname{proofsubpart}{Subpart}{Subparts}

\newcounter{proofdepth}
\AtBeginEnvironment{proof}{%
  \stepcounter{proofdepth}%
  \ifnum\value{proofdepth}=1
    \setcounter{step}{0}%
    \setcounter{proofpart}{0}%
  \fi
}
\AtEndEnvironment{proof}{%
  \addtocounter{proofdepth}{-1}%
}
\newtheorem{hyp}{Hypothesis}

\numberwithin{equation}{section}

\IfFileExists{luatex85.sty}{\usepackage{pdfcomment}}{}

\begin{document}
\renewcommand{\i}{\mathrm{i}}
\maketitle

\begin{abstract}
We develop a direct paradifferential approach to nonlinear graph integral
operators arising in contour dynamics.  The main methodological result is an
abstract structure-preserving paralinearization theorem, valid in arbitrary
dimension, for a class of homogeneous kernels whose nonlinear dependence on
the graph is expressed through normalized finite differences.  The theorem
covers both critical principal-value kernels and subcritical locally
integrable kernels, gives explicit formulas for the principal symbols, and
separates the finite-order symbolic contribution generated by the
homogeneous singularity from order-zero far-field terms and arbitrarily smoothing remainders.  A
structural feature of the reduction is the cancellation of the intermediate
symbolic order; at the critical endpoint the construction also preserves the
exact skew-adjoint structure of the original singular integral.

As a main application, we consider the three-dimensional two-phase
free-boundary Euler equations with constant interfacial background vorticity.  We first
derive an autonomous contour-dynamics equation in canonical surface
variables and then apply the abstract theorem directly to the singular
integral operators generated by the Birkhoff--Rott formulation.  This yields
the complete paradifferential structure of the system without taking the
Dirichlet--Neumann reduction as a starting point.  After introducing an
Alinhac good unknown and diagonalizing the resulting system, we obtain local
well-posedness for small Sobolev perturbations in the stable
Kelvin--Rayleigh--Taylor regime, uniformly on compact subsets of the stable
parameter set.
\end{abstract}


	{\small\tableofcontents}

	\allowdisplaybreaks

\section{Presentation of the problem and main contributions}

The analysis of fluid interfaces governed by contour-dynamics formulations leads naturally to nonlinear, nonlocal evolution equations driven by singular integral operators. Such equations arise in a broad range of models in fluid dynamics, including vortex sheets, water waves, surface quasi-geostrophic sharp fronts, porous-media interfaces, and related free-boundary problems. Their analysis has traditionally proceeded along two complementary directions: through the geometric and singular-integral structure of the underlying contour equations, and through paradifferential reductions designed to capture their quasilinear character. To the best of our knowledge, however, a general framework directly connecting nonlinear singular contour-dynamics operators with their paradifferential counterparts, while retaining the structural properties of the original formulation, has not been developed.

The purpose of this work is to provide such a framework. We introduce a direct paralinearization strategy for singular contour-dynamics operators, yielding a robust procedure for reducing nonlinear, nonlocal interface equations to paradifferential systems while preserving their essential structural features. The strategy is formulated at the level of an abstract class of homogeneous kernels, in arbitrary dimension, whose nonlinear dependence on the interface enters through normalized finite differences of the graph function; it covers both critical principal-value kernels and subcritical locally integrable ones, produces explicit principal symbols, and isolates the finite-order symbolic contribution of the singularity from far-field terms and smoothing remainders. In particular, this approach makes it possible to work directly with the underlying singular integral operators, without first reformulating the problem in terms of a Dirichlet--Neumann operator, as is customary in the analysis of water-wave equations.

Because the reduction is carried out on the kernel itself rather than on a specific model, it applies uniformly to a number of contour-dynamics operators that have so far been treated by model-dependent techniques. In the corresponding small-interface regimes, the singular parts of the circular Kelvin--Helmholtz operators, the local parts of the operators governing $\alpha$-SQG sharp fronts for $\alpha\in(1,2)$, and the compensated critical kernels of the three-dimensional Muskat equation all fall within the abstract class, with periodic coefficients, non-integer orders, and compensated finite-difference structures handled by the extensions of the abstract structure-preserving paralinearization theorem. These examples, collected in \cref{rem:main_result_graph_kernel_applications_Rd}, indicate that the framework is not tied to any particular free-boundary problem but reflects a common paradifferential structure shared across contour dynamics.

As a flagship application, we consider the motion of two immiscible, incompressible, inviscid fluids separated by a smooth interface and allowing for interfacial vorticity. This setting contains both the one-phase water-wave regime and the Kelvin--Helmholtz regime as distinguished limiting configurations of the density parameters. While substantial progress has been achieved for the water-wave problem \cite{WuWW1997,WuWW1999,LannesWW2005,AlazardBurqZuily2011,CL2015}, and the two-dimensional Kelvin--Helmholtz problem has a long and well-developed theory \cite{birkhoff1962helmholtz,rott1956diffraction,saffman1992vortex}, the three-dimensional two-phase system with general density ratio and interfacial vorticity remains largely open. In this setting, the main difficulty stems from the singular and nonlocal velocity field induced on the interface by the Birkhoff--Rott integral. We address this problem within the general framework developed here: we obtain a direct paralinearization of the governing singular integral operators, identify the full quasilinear structure of the system, and use this reduction to establish local well-posedness in Sobolev spaces.

\begin{notation}
\label{notation:vector_convention}
Throughout the manuscript we adopt the following notational conventions to distinguish scalar, vector, and matrix quantities:
\begin{itemize}
    \item scalar quantities (functions, operators, symbols, etc.) are denoted by lower-case letters (e.g., $a, b, c, \ldots$);
    \item two-coordinate vector quantities (functions, operators, symbols, etc.) are denoted by upper-case letters (e.g., $A, B, C, \ldots$) and shall be written in components as $A=\vect{a_1}{a_2}$;
    \item three-coordinate vector quantities (functions, operators, symbols, etc.) are denoted by lower-case bold letters (e.g., ${\bm a}, {\bm b}, {\bm c}, \ldots$), such vectors shall be decomposed as ${\bm a}=\vect{A}{a_3}$ relatively to the first two and last component;
    \item two-by-two matrix quantities (functions, operators, symbols, etc.) are denoted by upper-case bold letters (e.g., ${\bm A}, {\bm B}, {\bm C}, \ldots$);
    \item two-by-three matrix quantities (functions, operators, symbols, etc.) are denoted by Fraktur upper-case letters (e.g., $\mathfrak{A}, \mathfrak{B}, \mathfrak{C}, \ldots$).
\end{itemize}
The principal exceptions to this rule are the spatial variable $x \in \bR^2$ and the frequency variable $\xi \in \bR^2$, which we keep in lower case to avoid excessively heavy notation. This convention leads to some nonstandard choices, such as denoting the Dirichlet--Neumann operator by $\mathpzc{g}$ instead of the usual $G$. We believe, however, that the improved clarity afforded by this unified notation system justifies these departures from standard conventions.\\
We emphasise that this convention is restricted to the unknowns and operators that compose the system under study, and shall be especially useful in \cref{sec:paralinearization}; standard mathematical objects such as Banach and Hilbert spaces, norms, inner products, and other common objects of functional analysis retain their usual notation.
\end{notation}

 We denote generic elements in the plane and space as
$x\defeq\vect{x_1}{x_2}\in \bR^2$ and ${\bm x}\defeq  \vect{x}{x_3}\in\mathbb{R}^3$. 
The corresponding gradient operators in $\bR^2$ and $\bR^3$ are defined as
$
    \nabla\defeq\vect{\partial_{x_1}}{\partial_{x_2}}$ and $ \nabla_{\bm x}\defeq\vect{\nabla}{\partial_{x_3}}$.
We define 
\begin{equation}\label{def:bJ}
    \nabla^\perp \defeq {\bm J}\nabla,\qquad{\bm J}\defeq\begin{bmatrix}
        0 &-1\\
        1 & 0
    \end{bmatrix}.
\end{equation}
Given a finite time $T>0,$ we shall work on the finite time interval $I\defeq[-T, T]$. 
The interface separating both fluids is denoted $\Gamma(t)$ and is parametrized by ${\bm \gamma}(t,x)$, namely
$$\Gamma\pare{t}\defeq\set{{\bm \gamma}\pare{t,x}\in \bR^3\middle| \ x\in \bR^2}.$$
In what follows, we shall make the assumption that the interface $\Gamma(t)$ is a graph on $ \bR^2 $ parametrized as
\begin{equation}\label{eq:surface}
\pare{t,x}\mapsto{\bm \gamma}\pare{t,x}=\pare{x,\eta\pare{t,x}}\in \cC^1\pare{I\times\bR^2 ; \bR^3}.
\end{equation}
The elevation function $\eta:I\times \bR^2\to\bR $ is one of the main quantities of interest in our study. It separates $\mathbb R^3$ into two open and simply-connected sets  
 \begin{align}
     \Omega^{+}\pare{t}&\defeq\set{\pare{t, {\bm x}} = \pare{t, x, x_3}\in I\times\bR^3 \quad\textnormal{s.t.}\quad x_3 > \eta\pare{x}},\\
     \Omega^{-}\pare{t}&\defeq\set{\pare{t, {\bm x}} = \pare{t, x, x_3}\in I\times\bR^3 \quad\textnormal{s.t.}\quad x_3< \eta\pare{x}}.
 \end{align}
We assume that the domains $\Omega^\pm\pare{t}$ are filled with inviscid, immiscible fluids, with constant densities 
$$\rho^-\ge\rho^+ \ge 0.$$
The associated velocity and pressure fields are denoted ${\bm u^\pm}$ and $p^{\pm}$, respectively. The fluid is considered to be incompressible in both phases, i.e.
\begin{align}
\label{div1}
\dive_{{\bm x}} {\bm u^\pm}\defeq \nabla_{{\bm x}}\cdot {\bm u^\pm} =0, && \text{in}\quad I\times \Omega^\pm(t).
\end{align}
Here and in the sequel, the notation $u\cdot v\defeq u^\intercal v $ stands for the usual Euclidean scalar product. The evolution of the fluids is described by the incompressible non-homogeneous Euler equations
\begin{align}
\label{Euler}
\rho^\pm \pare{\partial_t{\bm u^\pm}+ \pare{{\bm u^\pm}\cdot\nabla_{{\bm x}}} {\bm u^\pm}}=-\nabla_{{\bm x}} p^\pm- \mathfrak{g} \rho^\pm \  {\bm e}_3, && \text{in}\quad I\times \Omega^\pm(t) . 
\end{align}
where $\mathfrak{g}\in \bR$ is the gravitational constant and $ {\bm e}_3 \defeq \nabla_{{\bm x}}x_3$. {Since the fluids are immiscible, the interface $\Gamma(t)$ is transported by the normal velocity component. Specifically, we define the velocity field $\bm{u}^{\pm}$ by continuous extension to the interface, acknowledging that while the normal component is continuous, the tangential component may exhibit a jump across $\Gamma(t)$. Thus, the interface evolves according to the transport equation:}

\begin{align}
\label{boundary1}
\pare{ \partial_t{\bm \gamma} -\left. {\bm u^\pm} \right|_{\Gamma\pare{t}} } \cdot{\bm \nu}=0,&&  {\bm \nu}= {\bm \nu}\pare{{\bm \gamma}}\defeq \partial_{x_1}{\bm \gamma}\times \partial_{x_2}{\bm \gamma} =\begin{bmatrix}
    -\nabla \eta\\
    1
\end{bmatrix}. 
\end{align}
Given a vector field $ \Omega^\pm $ defined on $ \Omega^\pm\pare{t} $ with trace $ \omega^\pm $ on $ \Gamma\pare{t},$ we denote
\begin{align}
\bbra{ \omega^\pm }\defeq \omega^{-}-\omega^+,&&
\aangle{\omega^\pm}\defeq\frac{\omega^- + \omega^+}{2}\cdot 
\end{align} 
Moreover, we define the {\it Atwood number}
\begin{equation}\label{eq:Atwood}
\Atw\defeq \Atw\pare{\rho^\pm}\defeq\frac{\bbra{\rho^\pm}}{ \aangle{\rho^\pm}} = 2 \  \frac{\rho^--\rho^+}{\rho^-+\rho^+}\in[0,2].
\end{equation}
The boundary condition for the pressure is given by the equation
\begin{equation}\label{eq:boundary_condition_pressure}
\bbra{\left. p^\pm \right|_{\Gamma\pare{t}}} = \sigma\mathpzc{k}\pare{{\bm \gamma}}, 
\end{equation}
where $ \sigma\geq 0 $ is the surface tension coefficient and $\mathpzc{k}\pare{{\bm \gamma}} $ is the mean curvature of the surface $ {\bm \gamma} $ which, under the assumption \cref{eq:surface}, becomes
\begin{align}
\label{curvature2}
\left.\mathpzc{k}\pare{{\bm \gamma}}  \middle|_{{\bm \gamma}\pare{x} = \pare{x, \eta\pare{x}}}\right. = - \nabla_x \cdot \mathpzc{K}\pare{\nabla\eta},
&&
\mathpzc{K}\pare{\nabla\eta}\defeq \frac{\nabla_x \eta}{\sqrt{1+\av{\nabla_x \eta}^2}}\cdot
\end{align} 
We complement the system identified by
\cref{div1,Euler,boundary1,curvature2} by fixing a constant
background sheet strength
\(\background\in\bR^2\) and working in the centered Galilean frame in which
the common horizontal velocity at infinity vanishes, i.e., we
impose
\begin{equation}\label{eq:background_far_field}
    {\bm u}^\pm\pare{t,x,z}
    \xrightarrow{\av{z}\to\infty}
    \pm\frac12
    \begin{bmatrix}
        \background^\perp\\
        0
    \end{bmatrix}
  ,
\end{equation}
and assume that the perturbations of these constant far-field states have
the Sobolev regularity specified below.  We also assume that
\[
    \nabla_{\bm x}\times{\bm u}^\pm=0
    \qquad\textnormal{in}\qquad
    \Omega^-\cup\Omega^+.
\]
The system that we finally obtain, and that we shall study along the manuscript, is the following
\begin{equation}
\label{eq:2phaseWW3D}
\system{
\begin{aligned}
& \rho^\pm \pare{\partial_t{\bm u^\pm}+ \pare{{\bm u^\pm}\cdot\nabla_{{\bm x}}} {\bm u^\pm}}=-\nabla_{{\bm x}} p^\pm- \mathfrak{g} \rho^\pm \  {\bm e}_3, && \text{in}\quad I\times \Omega^\pm(t),\\
& \pare{ \partial_t{\bm \gamma} -\left. {\bm u^\pm} \right|_{\Gamma\pare{t}} } \cdot{\bm \nu}=0, && \text{in}\quad I \times \Gamma(t),\\
&\bbra{\left. p^\pm \right|_{\Gamma\pare{t}}} = \sigma\mathpzc{k}\pare{{\bm \gamma}}, && \text{in}\quad I \times \Gamma(t),\\
& \nabla_{{\bm x}} \cdot {\bm u^\pm} =0, && \text{in}\quad I\times \Omega^\pm(t),\\
& \nabla_{{\bm x}} \times {\bm u^\pm} =0, && \text{in}\quad I\times \Omega^\pm(t). 
\end{aligned}
}
\end{equation}
It is well known that the irrotationality assumption in \cref{eq:2phaseWW3D} implies that there exists $ \phi^\pm : I\times\Omega^\pm\pare{t}\to \bR $ such that
\begin{align}\label{eq:velocity_potential}
\nabla_{\bm x} \phi^\pm={\bm u^\pm}, 
&&
\text{in}\quad I \times\Omega^\pm\pare{t}. 
\end{align}
\begin{rem}\label{rem:endpoint_identification}
    Let us mention that the two limiting cases $\Atw=0$ and $\Atw=2$ correspond to the Kelvin-Helmholtz and deep gravity-capillary 3D water-waves systems, respectively.
    \fine
\end{rem}
\subsection{Functional setting}
Given a function $u\in L^1\pare{\bR^2;\bC}$, we denote with
\[
    \hat{u}\pare{\xi} \defeq \frac{1}{2\pi} \int _{\bR^2} u\pare{x} e^{-\ii x\cdot\xi}\dd x ,  \quad \xi\in \bR^2 , 
\]
so that, if $\hat{u}\in L^1\pare{\bR^2;\bC}$, we can recover $u$ from $\hat{u}$ via the classical inversion formula
\begin{equation}
    u\pare{x} = \frac{1}{2\pi}  \int_{\bR^2} \hat{u}\pare{\xi}e^{\ii x\cdot\xi}\dd \xi . 
\end{equation}
It is immediate that $u$ is real valued if and only if $\hat{u}\pare{\xi} = \overline{\hat{u}\pare{-\xi}}$. 
We define, for any $s\in \bR$ the non-homogeneous Sobolev space ${H}^s\defeq {H}^s\pare{\bR^2;\bC}$ as the closure of $\cC^\infty_c\pare{\bR^2;\bC}$ functions w.r.t. the norm
\begin{align}
    \norm{u}_{ H^s} \defeq \sqrt{\int_{\bR^2} \angles{\xi}^{2s}\av{\hat u \pare{\xi}}^2 \dd \xi}, && \angles{\xi}\defeq \sqrt{1+\av{\xi}^2}; 
\end{align}

We use the notation
$$
u\otimes v = 
\begin{bmatrix}
    u_1\\
    u_2
\end{bmatrix}\otimes\begin{bmatrix}
    v_1\\
    v_2
\end{bmatrix}\defeq
u v^\intercal =
\begin{bmatrix}
    u_1 v_1 & u_1v_2\\
    u_2v_1 & u_2v_2
\end{bmatrix}.$$

\begin{notation}\label{notations:intro}
Let $T, r > 0$ and $s\in \bR$.
\begin{enumerate}
    \item we denote with $I\defeq \bra{-T, T}$;
    \item for $u\in H^s$, we set
    \begin{equation}
        \norm{u}_s
        \defeq
        \norm{u}_{H^s};
    \end{equation}
    \item for $u\in\cC\pare{I;H^s}$ and
    $\vec U_p\defeq\pare{u_1,\ldots,u_p}$, we set
    \begin{equation}
        \norm{u}_{s,I}
        \defeq
        \sup_{t\in I}\norm{u\pare t}_s,
        \qquad
        \norm{\vec U_p}_{s,I}
        \defeq
        \max_{1\leq j\leq p}\norm{u_j}_{s,I};
    \end{equation}
    \item  we denote with
    \begin{align}
        \Ball{s}{r} \defeq B_X \pare{0, r} , && X = \cC\pare{I;H ^s\pare{\bR^2;\bC }} . 
    \end{align}
    \item  we denote with
    \begin{align}
        \BallR{s}{r} \defeq B_X \pare{0, r} , && X = \cC\pare{I;H ^s\pare{\bR^2;\bR}} . 
    \end{align}
\end{enumerate}
   
\end{notation}

All symbol, kernel, and smoothing-operator classes below are scalar. Their
statements and estimates apply componentwise to finite-dimensional vectors and
matrices. A separately estimated argument written in brackets does not define
a new class: it denotes the evaluation of a continuous linear map with values
in the corresponding Fr\'echet space, according to
\cref{notation:paralinearization}. Such an argument is not included among the
small backgrounds and is not required to be small. Time is suppressed in
pointwise formulas, and all estimates are uniform in \(t\in I\).

{

\subsection{Main results}
\label{sec:main_result}

We organize the main results around the three steps that connect the
free-boundary Euler system to the Cauchy theory: the autonomous canonical
contour-dynamics formulation; an abstract paralinearization theorem for
graph integral operators, together with its application to the canonical
CDE; and the local well-posedness theorem.

\subsubsection{The contour-dynamics formulation}
\label{sec:main_result_contour_dynamics_3d}

The first step is to eliminate the bulk variables and close the evolution on
the interface.  At the level of contour dynamics, one could already obtain a
closed system by working with the elevation \(\eta\) and the surface
vorticity \(\Omega\).  We do not use these variables, however, since the
natural Hamiltonian coordinates for the two-fluid problem are the elevation
and the density-weighted trace potential \(\psi\) defined by
\begin{equation}\label{def:HAMpsi_intro}
   \psi\defeq\frac{ \bbra{\rho^\pm \psi^\pm}}{\aangle{\rho^\pm}},\qquad\psi^{\pm}\defeq\left.\phi^{\pm}\right|_{\Gamma\pare{t}}. 
\end{equation}
Passing from \(\Omega\) to \(\psi\) therefore does more than
close the equations in scalar surface variables: it puts the contour dynamics
directly in the canonical coordinates which will later be compatible with
the Hamiltonian symmetrization and the energy estimates.

Since the prescribed background $\background\in\mathbb{R}^2$ is represented by an affine, and therefore
non-Sobolev, potential, we further separate it by setting
\begin{equation}\label{eq:main_result_canonical_perturbation}
    \bar\psi\pare{x}
    \defeq
    x^\perp\cdot\background,
    \qquad
    \nabla^\perp\bar\psi=\background,
    \qquad
    \tilde\psi
    \defeq
    \psi-\bar\psi.
\end{equation}
The Bernoulli gauge is fixed as in \cref{eq:c_background}, so that the flat
configuration \(\pare{\eta,\tilde\psi}=\pare{0,0}\) is stationary.  The
resulting reduction is the following.

\begin{theorem}[The vortex-sheet equations as an autonomous canonical Contour Dynamic Equation (CDE)]
\label{thm:main_result_autonomous_cde_3d}
Let \(\sigma,\rho^\pm\geq0\) with
\(\aangle{\rho^\pm}>0\), let \(\mathfrak g\in\bR\), and let
\(\background\in\bR^2\).  There exist \(s_0>0\) and \(r_0>0\)
such that, in the class of  graph configurations for which
\(\eta\) and \(\psi\) are defined by \cref{eq:surface,def:HAMpsi_intro},
\(\tilde\psi\) is given by
\cref{eq:main_result_canonical_perturbation}, and $\eta\in\BallR{s_0}{r_0}$,
the free-boundary Euler system
\cref{eq:2phaseWW3D,eq:background_far_field} is equivalent to 
\begin{subequations}\label{eq:CDE_main_intro}
\begin{empheq}[left=\empheqlbrace]{align}
    \eta_t
    &=
    \mathpzc g^\star\pare{\eta}
    \bra{\tilde\Omega\bra{\background}\pare{\eta}}
    +
    \mathpzc g\pare{\Atw\middle|\eta}
    \bra{\tilde\psi},
    \label{eq:CDEeta_main_intro}
    \\
    \tilde\psi_t
    &=
    \tilde\pi
    \pare{\rho^\pm,\background\middle|\eta,\tilde\psi}.
    \label{eq:evolution_psi_main_intro}
\end{empheq}
\end{subequations}
The objects appearing in \cref{eq:CDE_main_intro} are defined, in their
order of appearance, as follows.
\begin{itemize}
    \item The operator
    \(\mathpzc g^\star\pare{\eta}\) is the normal-velocity operator in
    vorticity variables defined by
    \begin{equation}\label{def:gstar_intro}
        \mathpzc g^\star\pare{\eta}\Omega
        \defeq
        {\bm v}^\star\bra{\Omega}\pare{\eta}
        \cdot\vect{\nabla\eta}{-1},\qquad{\bm v}^\star\bra{\Omega}\pare{\eta}
        \defeq
        \vect{
             V^{\star}\bra{\Omega}\pare{\eta}
        }{
             v_3^{\star}\bra{\Omega}\pare{\eta}
        }
        =
        \vect{
            \pare{
                \mathpzc B_0\pare{\eta}
                \bra{\nabla\eta\cdot\Omega}
                -
                \mathpzc q\pare{\eta}\Omega
            }^\perp
        }{
            -\mathpzc B_0\pare{\eta}^\perp\cdot\Omega
        },
    \end{equation}
    where the singular integral operators $\mathpzc{B}_0$ and $\mathpzc{q}$ are given by
     \begin{subequations}
    \begin{align}
\label{eq:B0_intro}
    \mathpzc{B}_0\pare{\eta} \bra{\omega} \pare{x} \defeq & \ \frac{1}{4\pi} \pv \int_{\bR^2} \frac{z}{\pare{ \av{z}^2+\pare{ \delta_{z}\eta \pare{x} }^2 }^\frac 32} \ \omega\pare{x-z} \ \dd z, \\
    \label{eq:Q_intro}
    \mathpzc{q}\pare{\eta} \bra{\omega} \pare{x} \defeq & \frac{1}{4\pi} \pv \int_{\bR^2} \frac{\delta_z\eta\pare{x}}{\pare{ \av{z}^2+\pare{ \delta_{z}\eta \pare{x} }^2 }^\frac 32} \ \omega\pare{x-z} \ \dd z,
\end{align}
\end{subequations}
with $\delta_z \eta\pare x
    \defeq
    \eta\pare x-\eta\pare{x-z}.$

    \item The canonical background vorticity is
    \begin{equation}\label{eq:main_result_background_vorticity_preview_3d}
        \tilde\Omega\bra{\background}\pare{\eta}
        \defeq
        S\pare{\Atw\middle|\eta}\bar\psi,
    \end{equation}
    where the potential-to-vorticity map
    \(S\pare{\Atw\middle|\eta}\) is defined by
    \begin{equation}\label{def:S_intro}
         \begin{aligned}
             S\pare{\Atw\middle|\eta}\psi
        &\defeq
        \nabla^\perp
        \pare{1+\Atw\tau_0\pare{\eta}}^{-1}\psi,\\
        \tau_0\pare{\eta}\psi
        &\defeq
        \Delta^{-1}\nabla^\perp\cdot
        \pare{
            \mathpzc B_0\pare{\eta}
            \bra{\nabla\eta\cdot\nabla^\perp\psi}
            -
            \mathpzc q\pare{\eta}\nabla^\perp\psi
            +
            \pare{
                \mathpzc B_0\pare{\eta}
                \cdot\nabla\psi
            }
            \nabla^\perp\eta
        },
         \end{aligned}
    \end{equation}
where the operators $\mathpzc{B}_0$ and $\mathpzc{q}$ are given in \cref{eq:B0_intro,eq:Q_intro}.
    \item The generalized Dirichlet--Neumann operator in canonical
    variables is
    \begin{equation}\label{eq:main_result_dn_factorization_preview_3d}
        \mathpzc g\pare{\Atw\middle|\eta}
        \defeq
        \mathpzc g^\star\pare{\eta}
        \circ
        S\pare{\Atw\middle|\eta},
    \end{equation}
    with $\mathpzc{g}^{\star}$ and $S$ as in \cref{def:gstar_intro,def:S_intro}.  

    \item The function
    \(\tilde\pi\pare{\rho^\pm,\background\middle|\eta,\tilde\psi}\) in \cref{eq:evolution_psi_main_intro}
    is the shifted Bernoulli function defined by
    \begin{equation}\label{def:pitilde_intro}
        \tilde\pi\pare{\rho^\pm,\background\middle|\eta,\tilde\psi}\defeq
        \pi^\star
        \pare{
            \rho^\pm
            \middle|
            \eta,
            \tilde{\Omega}\bra{\background}\pare{\eta}
            +
            S\pare{\Atw\middle|\eta}\tilde{\psi}
        }
        -
        c_{\background},
    \end{equation}
 where
 \begin{align}
     \pi^\star
        \pare{
            \rho^\pm
            \middle|
            \eta,\Omega
        }&\defeq\frac{\sigma}{\aangle{\rho^\pm}}
        \nabla\cdot\mathpzc K\pare{\nabla\eta}
        -
        \mathfrak g\Atw\eta-\Omega^{\perp}\cdot V^{\star}\bra{\Omega}\pare{\eta}\\
        &\quad-\tfrac{\Atw}{4}\left[2\av{V^{\star}\bra{\Omega}\pare{\eta}}^2-2\av{v_3^{\star}\bra{\Omega}\pare{\eta}}^2+4v_3^{\star}\bra{\Omega}\pare{\eta}V^{\star}\bra{\Omega}\pare{\eta}\cdot\nabla\eta+\tfrac{1}{2}\tfrac{\av{\Omega}^2+\pare{\Omega\cdot\nabla\eta}^2}{1+\av{\nabla\eta}^2}\right],
 \end{align}
 where $\mathpzc{K},$ $V^{\star}$ and $v_3^{\star}$ are given in \cref{curvature2,def:gstar_intro}. The Bernoulli gauge 
 $$c_{\background}\defeq-\frac{\Atw}{8}\av{\background}^2$$ gives $ \tilde\pi
        \pare{\rho^\pm,\background\middle|0,0}
        =0$. 
\end{itemize}
\end{theorem}

\subsubsection{Paradifferential reduction of the CDE}
\label{sec:main_result_abstract_paralinearization_Rd}

The second, and main technical result of the manuscript consists in the paradifferential reduction of the CDE \cref{eq:CDE_main_intro}. This is achieved via a two-step process:
\begin{enumerate}[i]
    \item we prove an abstract theorem (\cref{thm:main_result_graph_kernel_paralinearization_Rd}) that allows us to paralinearize arbitrary nonlinear convolution kernel operators in any dimension $\mathtt d \geq 1$ with suitable integrability hypothesis and;
    \item  we apply this abstract result to paralinearize the terms in the CDE for the 3D vortex-sheet equation that resist classical tools of paradifferential reduction, mainly the singular integral operators $\mathpzc{B}_0$ and $\mathpzc{q}$ in \cref{eq:B0_intro,eq:Q_intro}. Next, we perform a careful analysis of every term involved in the reduction in order to obtain the final paradifferential system. 
\end{enumerate}

At this stage we deliberately avoid introducing the full paradifferential
framework needed for a completely quantitative formulation, including the
precise symbol and smoothing-operator classes and their associated tame
estimates.  The statement below should therefore be regarded as a concise
preview of the paralinearization mechanism, whose notation and precise
meaning are developed in \cref{sec:preliminaries}.
Throughout this preview, the Sobolev estimates are understood for every
\(s\geq s_0\), where \(s_0\) is a sufficiently large regularity threshold
whose value is not made explicit.  The techniques developed in the manuscript
would in principle allow one to track an admissible value of \(s_0\).
We have chosen not to do so, since this would require a substantial amount of
additional computations and low-order estimates without introducing any new
conceptual ingredient, while considerably burdening the exposition.

\paragraph{Paralinearization of abstract nonlinear convolution operators in arbitrary dimensions.}

This is one of the main methodological contributions of the
paper.  We consider $\mathtt d \in \bN^{\star}$ and  graph integral operators on \(\bR^{\mathtt d}\) of the
schematic form
\begin{align}\label{eq:main_result_schematic_graph_operator_Rd}
    \mathpzc t^{\pare{n}}\pare\eta\bra\omega\pare x
    =
    \int_{\bR^{\mathtt d}}
   \frac{h\pare{\hat z}}{\av{z}^n}\,
    k\pare{\frac{\eta\pare x-\eta\pare{x-z}}{\av z}}\,
    \omega\pare{x-z}\dd z,
    &&
    0<n\leq \mathtt d, \quad
    \hat{z}\defeq \frac{z}{\av{z}},
\end{align}
where the functions $h$ and $k$ are smooth (cf. \cref{def:abstract-admissible-angular-pair-Rd} for the precise regularity requirements) and the integral is understood in the principal-value sense in the
critical case $n=\mathtt{d}$
whereas in the subcritical case the convolution kernel is locally
integrable around zero.  Thus the singular or weakly singular kernel depends
nonlinearly on finite differences of the graph itself.  This interaction
between the \(z\)-singularity and the nonlinear geometry prevents a direct
application of standard paralinearization techniques.

Rather than analyzing such operators on a case-by-case basis, we develop
an abstract structure-preserving paralinearization theorem which generalizes
the mechanisms introduced in \cite{BCGS2025,murgante2025long}.  The theorem
works in arbitrary dimension, covers both the critical and subcritical
regimes in \cref{eq:main_result_schematic_graph_operator_Rd}, and captures
the structural cancellation of the intermediate symbolic order.  This is
one of the main technical innovations of the paper and provides the
operator-level foundation for the paradifferential reduction of the
vortex-sheet system. The proof is correspondingly divided into two parts: the critical case is established in
\cref{sec:abstract-structure-preserving-paralinearization-Rd}; while the
subcritical case is proved in
\cref{app:abstract-locally-integrable-paralinearization-Rd}, here we present a statement that is common to both critical and subcritical cases.

\begin{theorem}[Critical and subcritical paralinearization of graph kernel operators]
\label{thm:main_result_graph_kernel_paralinearization_Rd}
Let \(\mathtt d\geq1\), \(0<n\leq\mathtt d\), and let
\(\pare{h,k}\) be an admissible homogeneous pair of order \(-n\) and
vanishing order \(\ell\in\bN^\star\) in the sense of
\cref{def:abstract-admissible-angular-pair-Rd}.  Set
\[
    \beta\defeq\mathtt d-n\geq0.
\]
When \(\beta>0\), assume in addition
\cref{eq:abstract-subcritical-orders-Rd}.
There exists \(\varrho_0\geq0\) such that, for every
\(\varrho\geq\varrho_0\), there exist
\(s_0>\mathtt d/2+1\), \(\mu\geq0\), and \(r_0>0\) such that, for every
\(s\geq s_0\) and every
\[
    \eta\in H^{s+1+\mu}\pare{\bR^{\mathtt d};\bR},
    \qquad
    \omega\in H^{s+\mu}\pare{\bR^{\mathtt d};\bR},
    \qquad
    \norm\eta_{H^{s_0+1+\mu}}<r_0,
\]
the graph operator in
\cref{eq:main_result_schematic_graph_operator_Rd} satisfies
\begin{equation}\label{eq:main_result_graph_kernel_paralinearization_Rd}
    \mathpzc t^{\pare n}\pare\eta\bra\omega
    =
    \OpBW{
        A_{1,-\beta}^{h,k}\bra\omega
        \pare{\eta;x,\xi}
        +
        \vect{
            a_{\bra0}\bra\omega\pare{\eta;x,\xi}
        }{
            a_{\bra{-\beta-2}}\pare{\eta;x,\xi}
        }
    }
    \cdot
    \vect{\eta}{\omega}
    +
    R_{\bra\varrho}^{\lin}\pare{\eta,\omega}
    \cdot
    \vect{\eta}{\omega}.
\end{equation}
The objects in
\cref{eq:main_result_graph_kernel_paralinearization_Rd} are as follows.
\begin{itemize}

    \item The principal vector symbol is
    \begin{equation}\label{eq:main_result_graph_kernel_principal_vector_Rd}
        A_{1,-\beta}^{h,k}\bra\omega
        \pare{\eta;x,\xi}
        \defeq
        \vect{
            \omega\pare x
            a_1^{h,k;\eta}
            \pare{\nabla\eta\pare x,\xi}
        }{
            a_{-\beta}^{h,k;\omega}
            \pare{\nabla\eta\pare x,\xi}
        }.
    \end{equation}
    Its components have the following properties.
    \begin{itemize}

        \item The transported symbol
        \(a_1^{h,k;\eta}\pare{\nabla\eta,\xi}
        \in\Gamma^1_{\geq\ell-1}\) is given explicitly by
        \begin{equation}\label{eq:main_result_graph_kernel_a1_Rd}
            a_1^{h,k;\eta}\pare{p,\xi}
            \defeq
            \begin{cases}
                \displaystyle
                \frac{\pi}{2}
                \int_{\bS^{\mathtt d-1}}
                \av{\vartheta\cdot\xi}\,
                h\pare\vartheta
                k'\pare{p\cdot\vartheta}
                \dd\sigma\pare\vartheta,
                & \beta=0,
                \\[1.2em]
                \displaystyle
                \int_{\bR^{\mathtt d}}
                \chi\pare z
                \frac{h\pare{\hat z}}{\av z^{n+1}}\,
                k'\pare{p\cdot\hat z}
                \pare{1-e^{-\ii z\cdot\xi}}
                \dd z,
                & \beta>0,
            \end{cases}
        \end{equation}
        where $\chi$ is the cut-off function introduced in \cref{def-chi}.

        \item The quasilinear symbol satisfies
        \[
            a_{-\beta}^{h,k;\omega}
            \pare{\nabla\eta,\xi}
            -
            a_{-\beta}^{h,k;\omega}\pare{0,\xi}
            \in
            \Gamma^{-\beta}_{\geq\ell}
        \]
        and is given explicitly by
        \begin{equation}\label{eq:main_result_graph_kernel_aomega_Rd}
            a_{-\beta}^{h,k;\omega}\pare{p,\xi}
            \defeq
            \begin{cases}
                \displaystyle
                -\frac{\ii\pi}{2}
                \int_{\bS^{\mathtt d-1}}
                \sgn\pare{\vartheta\cdot\xi}\,
                h\pare\vartheta
                k\pare{p\cdot\vartheta}
                \dd\sigma\pare\vartheta,
                & \beta=0,
                \\[1.2em]
                \displaystyle
                \int_{\bR^{\mathtt d}}
                \chi\pare z
                \frac{h\pare{\hat z}}{\av z^n}\,
                k\pare{p\cdot\hat z}
                e^{-\ii z\cdot\xi}
                \dd z,
                & \beta>0.
            \end{cases}
        \end{equation}

    \end{itemize}

    \item The lower vector symbol satisfies
    \begin{equation}\label{eq:main_result_graph_kernel_lower_orders_Rd}
        a_{\bra0}\bra\cdot
        \in
        \cL\pare{H^{s_0};\Gamma^0},
        \qquad
        a_{\bra{-\beta-2}}
        \in
        \Gamma^{-\beta-2}_{\geq\ell}.
    \end{equation}
    When \(\beta=0\), the second component can be chosen in the form
    \begin{equation}\label{eq:main_result_graph_kernel_critical_lower_symbol_Rd}
        a_{\bra{-2}}\pare{\eta;x,\xi}
        \defeq
        \ii  a_{\bra{-2}}^{\bR;\omega}
        \pare{\eta;x,\xi},
        \qquad
        a_{\bra{-2}}^{\bR;\omega}
        \in
        \Gamma^{-2}_{\geq\ell},
    \end{equation}
    where \( a_{\bra{-2}}^{\bR;\omega}\) is real-valued and odd in
    \(\xi\).  In particular, there is no independent symbol of order
    \(-\beta-1\).

    \item
    \(R_{\bra\varrho}^{\lin}\pare{\eta,\omega}\) is a generic
    real-to-real linear smoothing remainder in the sense of
    \cref{notation:paralinearization}; its first column is the evaluation
    of a map in
    \(\cL\pare{H^{s_0};\cR^\varrho}\), while its second column belongs to
    \(\cR^\varrho_{\geq\ell}\).

\end{itemize}
When \(\beta=0\), the operator is exactly skew-adjoint,
\[
    a_{-\beta}^{h,k;\omega}
    =
    \ii a_0^{h,k;\omega},
\]
with \(a_0^{h,k;\omega}\) real-valued and odd in \(\xi\), and its complete
nonsmoothing quasilinear symbol acting on \(\omega\) is
\begin{equation}\label{eq:main_result_graph_kernel_critical_complete_symbol_Rd}
    \ii\pare{
        a_0^{h,k;\omega}\pare{\nabla\eta\pare x,\xi}
        +
        a_{\bra{-2}}^{\bR;\omega}\pare{\eta;x,\xi}
    },
\end{equation}
hence is purely imaginary and odd in \(\xi\).  The two principal
components satisfy
\cref{eq:abstract-principal-symbol-differential-identity}.
The complete critical and subcritical statements are respectively
\cref{thm:abstract-structure-preserving-paralinearization-Rd,%
thm:abstract-locally-integrable-paralinearization-Rd}.
\end{theorem}

\begin{rem}[Scope of the model class]
\label{rem:main_result_graph_kernel_scope_Rd}
The class in
\cref{eq:main_result_schematic_graph_operator_Rd}
is intended to capture the basic singular structure of the contour-dynamics
operators arising in applications, rather than to provide an exhaustive
description of all possible kernels and operators.  We record four extensions
which are not needed for the applications developed in the present paper.

\begin{enumerate}[i]
    \item
    \emph{Non-homogeneous kernels.}
    Exact homogeneity of the convolution kernel is not essential.  Suppose,
    more generally, that a kernel \(g\) admits near \(z=0\) a finite
    polyhomogeneous expansion of the form
    \[
        g\pare z
        =
        \sum_{j=0}^{J}
        \frac{h_j\pare{\hat z}}{\av z^{n-j}}
        +
        g_{\mathrm{rem}}\pare z .
    \]
    Each homogeneous component may be treated separately by
    \cref{thm:main_result_graph_kernel_paralinearization_Rd}.  The expansion
    is continued only until the localized remainder
    \(\chi g_{\mathrm{rem}}\) has sufficiently many powers of \(z\) at the
    origin to belong to the locally integrable kernel classes of
    \cref{prop:int-opbw}.  At that point the remaining nonlinear
    dependence is reduced by the usual Bony decomposition and
    {
\cref{prop:int-opbw}, while the part supported away from the origin
contributes an order-zero paraproduct modulo arbitrarily smoothing
remainders by \cref{prop:tail-kernel-smoothing};
}
cf. also
\cref{lem:abstract-common-local-tail-Bony-reduction-Rd}.  Consequently, only finitely
    many homogeneous components of the kernel contribute to any prescribed
    finite portion of the paradifferential expansion.

    \item
    \emph{Periodic interfaces.}
    The same mechanism applies on \(\mathbb{T}^{\mathtt d}\), after replacing the
    Euclidean paradifferential calculus by its periodic counterpart.  The
    singular analysis is entirely local near the diagonal \(z=0\) and is
    therefore unchanged.  In fact, the global part is simpler on the torus:
{
    after localizing the singularity, the complementary geometric kernel
    is smooth. Its graph dependence contributes an order-zero paraproduct
    modulo smoothing remainders.
}  Thus no
    far-field decay assumptions, and no Euclidean tail analysis at
    \(\av z\to+\infty\), are required, as can be observed for instance in the works \cite{BCGS2025,BMM2024,murgante2025long}.

    \item
    \emph{Higher-order contour-dynamics operators.}
    The restriction \(n\leq\mathtt d\) in
    \cref{eq:main_result_schematic_graph_operator_Rd} is likewise not
    intrinsic to the paralinearization mechanism.  Higher-order
    contour-dynamics equations naturally lead to compensated singular
    integrals of the form
    \begin{align}
          \int_{\bR^{\mathtt d}}
        \frac{h\pare{\hat z}}{\av z^n}\,
        k\pare{\Delta_z\eta\pare x}\,
        \delta_z\omega\pare x
        \dd z,
        &&
        n\in \pare{\mathtt{d}, \mathtt{d}+1 }
    \end{align}
  and, for more singular kernels, to the corresponding expressions
    containing higher finite differences
    \(\delta_z^m\omega\).  The finite differences provide precisely the
    cancellations required to compensate singularities with
    \(n>\mathtt d\); for instance, an \(m\)-th order finite difference gives
    \(m\) powers of \(z\) at the origin, so that the natural locally
    integrable range becomes \(n<\mathtt d+m\), with the corresponding
    parity or moment cancellation required at the endpoints.

    We do not formulate a separate abstract theorem for these operators
    here.  Nevertheless, repeating the Bony decomposition, the homogeneous
    kernel analysis, and the compensated kernel-to-symbol reductions
    developed above yields their paralinearization by the same procedure.
    In particular, increasing the differential order of the CDE changes the
    homogeneity of the resulting symbols but not the underlying
    paralinearization mechanism.
    \item
\emph{Finite regularity of the kernel function $k$.}
Real-analyticity of $k$ in \cref{def:abstract-admissible-angular-pair-Rd} is
used only in the proof of \cref{lem:kernel}, to sum the convergent series
\cref{eq:g_power_series_kernel} by Cauchy-type estimates. Since
\cref{def:kernel-global-global} only demands, at each fixed order $N$ of
$z$-differentiation, a tame bound \cref{eq:kernelfunctionproperty-local-tame}
with a loss $\mu_N^k$ allowed to grow with $N$, this convergent expansion
may be replaced, order by order, by the Taylor formula with integral
remainder
\begin{equation}\label{eq:abstract-admissible-K-finite-Taylor}
    k\pare X
    =
    \sum_{j=0}^{N}\frac{k^{\pare j}\pare0}{j!}X^j
    +
    \frac{X^{N+1}}{N!}\int_0^1\pare{1-t}^N k^{\pare{N+1}}\pare{tX}\dd t,
\end{equation}
estimating the remainder exactly as $k_{\bra2}$ is estimated in the proof of
\cref{lem:kernel} — the same mechanism already used at order two in
\cref{eq:abstract-quasilinear-composed-Taylor-remainder}. Hence
$k\in\cC^\infty\pare{\pare{-r_k,r_k};\bR}$ suffices in place of
$k\in\cC^\omega\pare{\pare{-r_k,r_k};\bR}$; the vanishing-order hypothesis
\cref{eq:abstract-admissible-K-vanishing-order} is unaffected, since it
fixes only the first $\ell$ Taylor coefficients of $k$.
\end{enumerate}
Thus \cref{eq:main_result_schematic_graph_operator_Rd} isolates the basic
homogeneous Euclidean model in which the method can be stated most
transparently, while the procedure itself is stable under the extensions
described above.

\fine
\end{rem}

\begin{rem}[Representative contour-dynamics applications]
\label{rem:main_result_graph_kernel_applications_Rd}
In the corresponding small-interface regimes, the mechanism of
\cref{thm:main_result_graph_kernel_paralinearization_Rd}, together with the
extensions in \cref{rem:main_result_graph_kernel_scope_Rd}, applies to the
following contour-dynamics operators. In each case, we use the notations of the corresponding papers \cite{murgante2025long,BCGS2025,CG2007,AN2022}.
\begin{itemize}

    \item
    \emph{Circular Kelvin--Helmholtz.}
    With
    \(q_{\mathrm{KH}}\defeq\frac12\log\pare{1+2\eta}\), the singular
    parts of \(\mathpzc{H}_0\) and
    \(\pare{1+2\eta}\mathpzc D_0-\Pi_0\) are the periodic critical
    models associated respectively with
    \[
        \pare{h,k,\ell}
        =
        \pare{\vartheta/\pi,\pare{1+X^2}^{-1},2},
        \qquad
        \pare{h,k,\ell}
        =
        \pare{1/\pi,X\pare{1+X^2}^{-1},1}.
    \]
    The remaining periodic and geometric corrections start at quadratic
    order in the chord variable and are therefore lower order
    \cite{murgante2025long}.

    \item
    \emph{\(\alpha\)-SQG sharp fronts.}
    For \(\alpha\in\pare{1,2}\), after the analogous logarithmic change
    and the extraction of one periodic chord through
    \(\delta_zf=2\sin\pare{z/2}\Delta_zf\), the local parts of
    \(\mathcal I\) and \(\mathcal J\) have
    \[
        n=\alpha-1,
        \qquad
        \beta=2-\alpha,
        \qquad
        \ell_{\mathcal I}=1,
        \qquad
        \ell_{\mathcal J}=2.
    \]
    Their \(z\)-dependent periodic coefficients are handled by
    items~\textnormal{(i)}--\textnormal{(ii)} of
    \cref{rem:main_result_graph_kernel_scope_Rd}
    \cite{BCGS2025}.

    \item
    \emph{Three-dimensional Muskat.}
    Up to the physical prefactor, the graph equation is a sum over
    \(j\in\set{1,2}\) of compensated critical kernels with
    \[
        \mathtt d=n=2,
        \qquad
        h_j\pare\vartheta=\vartheta_j,
        \qquad
        k\pare X=\pare{1+X^2}^{-3/2},
        \qquad
        \ell=2,
        \qquad
        \omega=\partial_jf,
    \]
    in which \(\omega\pare{x-z}\) is replaced by
    \(\delta_z\omega\pare x\).  This is the lower-endpoint case of the
    compensated finite-difference mechanism described in
    item~\textnormal{(iii)} of
    \cref{rem:main_result_graph_kernel_scope_Rd}
    \cite{CG2007,AN2022}.
\end{itemize}
\fine
\end{rem}

The operators \(\mathpzc B_0\pare\eta\) and \(\mathpzc q\pare\eta\) in \cref{eq:B0_intro,eq:Q_intro} are
the \(\mathtt d=2\) critical instances, with vanishing orders \(\ell=2\)
and \(\ell=1\), respectively; see
\cref{lem:B0-q-admissible-kernels,prop:B0_paralinearization,%
prop:paralinearization_Q}.  They provide the singular-integral input for
the application to the canonical CDE.

\paragraph{Paralinearization of the 3D Vortex-Sheet CDE.}

The second part consists in the application of \cref{thm:main_result_graph_kernel_paralinearization_Rd} to the  3D Vortex-Sheet CDE, the complete statement is \cref{thm:paraCDEetapsi}.

\begin{theorem}[Conditional paradifferential reduction]
\label{thm:main_result_paradifferential_reduction}
Let \(\sigma\geq0\), \(\rho^\pm\geq0\) with
\(\aangle{\rho^\pm}>0\), \(\mathfrak g\in\bR\), and
\(\background\in\bR^2\).  There exists \(\varrho_0>0\) such that, for
every \(\varrho\geq\varrho_0\), there exist \(s_0>0\) and \(r_0>0\)
with the following property.  If
\begin{align}
   \eta\in\BallR{s_0}{r_0}, && \tilde\psi \in \cC\pare{I;H ^s\pare{\bR^2;\bR}}
\end{align}
solve the canonical CDE \cref{eq:CDE_main_intro}, then, setting
\[
    Y\defeq\vect{\eta}{\tilde\psi},
\]
 one has
\begin{equation}\label{eq:main_result_paradifferential_reduction}
\begin{aligned}
    \partial_tY
    ={}&
    \OpBW{
        {\bm Q}
        \pare{
            \Atw,\background
            \middle|
            \eta,\tilde\psi;x,\xi
        }
        +
        {\bm B}
        \pare{
            \Atw,\background
            \middle|
            \eta,\tilde\psi;x,\xi
        }
        +
        \ii
        \pare{
            \mathpzc V_{\textnormal{can}}
            \pare{
                \Atw,\background
                \middle|
                \eta,\tilde\psi;x
            }
            \cdot\xi
        }
        \Id_{\bC^2}
    }Y
    \\
    &+
    \OpBW{
        \vectT{
            A_{\bra0}^\eta
            \pare{
                \Atw,\background
                \middle|
                \eta,\tilde\psi;x,\xi
            }
        }{
            A_{\bra{-1}}^\psi
            \pare{
                \Atw,\background
                \middle|
                \eta,\tilde\psi;x,\xi
            }
        }
    }Y
    +
    {\bm R}_{\bra\varrho}
    \pare{
        \Atw,\background
        \middle|
        \eta,\tilde\psi
    }Y .
\end{aligned}
\end{equation}
The objects in \cref{eq:main_result_paradifferential_reduction} have the
following structure, in their order of appearance:
\begin{itemize}
    \item The quasilinear symbol
    \({\bm Q}\in\pare{\Gamma_{\geq0}^2}^{2\times2}\) is defined by
    \begin{equation}\label{eq:main_result_Q_preview}
        {\bm Q}
        \defeq
        \begin{bmatrix}
            0
            &
            \lambda_1+\Atw\lambda_0
            \\
            -\mathfrak g\Atw
            +
            a^{\pi_0,\eta}_2+\lambda_1q^{\tilde\pi}_0
            &
            0
        \end{bmatrix}.
    \end{equation}
    Its constituent symbols are defined explicitly in
    \cref{eq:lambda1,eq:lambda0,eq:Api_star_2,eq:qpitilde_main}.  In
    particular,
    \begin{equation}\label{eq:main_result_Q_flat_preview}
        {\bm Q}
        \pare{
            \Atw,\background
            \middle|
            0,0;\xi
        }
        =
        \begin{bmatrix}
            0 & \frac12\av{\xi}
            \\
            \mathfrak q
            \pare{\Atw,\background\middle|\xi}
            & 0
        \end{bmatrix},
    \end{equation}
    where
    \begin{equation}\label{eq:main_result_flat_q_preview}
        \mathfrak q
        \pare{\Atw,\background\middle|\xi}
        \defeq
        -\mathfrak g\Atw
        -
        \frac{\sigma}{\aangle{\rho^\pm}}\av{\xi}^2
        +
        \frac{4-\Atw^2}{8}
        \frac{
            \pare{\background^\perp\cdot\xi}^2
        }{
            \av{\xi}
        },
        \qquad
        \xi\in\bR^2\setminus\set{0}.
    \end{equation}

    \item The order-one symbol
    \({\bm B}\in\pare{\Gamma_{\geq1}^1}^{2\times2}\) has the Alinhac
    block form
    \begin{equation}\label{eq:main_result_B_preview}
        {\bm B}
        \defeq
        -
        \begin{bmatrix}
            -b_{\textnormal{can}}
            & 0
            \\
            b_{\textnormal{can}}^2
            &
            b_{\textnormal{can}}
        \end{bmatrix}
        \lambda_1.
    \end{equation}
    The coefficient \(b_{\textnormal{can}}\) is defined in
    \cref{eq:bpitilde_main}.  

    \item The transport velocity
    \(\mathpzc V_{\textnormal{can}}
    \in\pare{\cF^\bR_{\geq0}}^2\) is real-valued, independent of \(\xi\),
    and satisfies
    \begin{equation}\label{eq:main_result_transport_preview}
        \mathpzc V_{\textnormal{can}}
        \pare{\Atw,\background\middle|\eta,\tilde\psi}
        \defeq
        \frac{\Atw}{4}\background^\perp
        +
        \mathpzc V_{\textnormal{can}}^{\textnormal{nl}}
        \pare{\Atw,\background\middle|\eta,\tilde\psi},
        \qquad
        \mathpzc V_{\textnormal{can}}^{\textnormal{nl}}
        \in\pare{\cF^\bR_{\geq1}}^2.
    \end{equation}
    Its full formula is given in \cref{eq:Vpitilde_main}.

    \item The remaining paradifferential symbols satisfy
    \begin{equation}\label{eq:main_result_lower_order_preview}
        A_{\bra0}^\eta
        \in\pare{\Gamma_{\geq1}^0}^2,
        \qquad
        A_{\bra{-1}}^\psi
        \in\pare{\Gamma_{\geq1}^{-1}}^2,
    \end{equation}
    and hence form a two-by-two block of orders
    \begin{equation}
        \begin{bmatrix}
            \Gamma_{\geq1}^0 & \Gamma_{\geq1}^{-1}
            \\
            \Gamma_{\geq1}^0 & \Gamma_{\geq1}^{-1}
        \end{bmatrix}.
    \end{equation}
   
    \item The remainder
    \begin{equation}\label{eq:main_result_smoothing_preview}
        {\bm R}_{\bra\varrho}
        \in
        \pare{\cR_{\geq1}^\varrho}^{2\times2}
    \end{equation}
    is a tame smoothing operator in the sense of
    \cref{def:smoothing,notation:paralinearization}.  
\end{itemize}
\end{theorem}

\begin{rem}[Formal recovery of the one-phase water-wave structure]
\label{rem:main_result_one_phase_water_wave_structure}
It is instructive to compare
\cref{thm:main_result_paradifferential_reduction} with the classical
paradifferential reduction of the one-phase water-wave equations.  We stress
that the following discussion is only a formal identification of the
quasilinear structures: no convergence result as
\(\Atw\to2\) and \(\background\to0\) is intended, and we do not attempt to
identify all the bounded symbols appearing in
\cref{eq:main_result_lower_order_preview}.

Recall first that the endpoint \(\Atw=2\) corresponds to
\(\rho^+=0\).  By \cref{eq:psi}, at this endpoint
 $\psi
    =
    2\psi^-$. Thus, when \(\background=0\), so that \(\tilde\psi=\psi\), the standard
one-phase surface potential is
\begin{equation}\label{eq:main_result_water_wave_potential_normalization}
\psi_{\rm WW}
\defeq
\frac12\tilde\psi.
\end{equation}
This factor \(2\) is entirely due to the normalization of the canonical
two-phase potential introduced following \cite{Lannes2013}.  In particular,
quadratic expressions in the potential acquire the corresponding factor
\(4\); this is consistent with the quadratic dependence of the kinetic
Hamiltonian on the surface potential and does not represent a difference
between the two dynamics.

The comparison is especially transparent at the level of the unbounded
paradifferential terms, which are the terms relevant for the derivative
balance in the energy estimates.  By \cref{eq:lambda1},
\begin{equation}\label{eq:main_result_water_wave_DN_symbol}
2\lambda_1\pare{\eta;x,\xi}
=
\sqrt{
\pare{1+\av{\nabla\eta}^2}\av{\xi}^2
-
\pare{\nabla\eta\cdot\xi}^2
},
\end{equation}
and the right-hand side is precisely the principal symbol of the one-phase
Dirichlet--Neumann operator appearing in
\cref{eq:historical-water-wave-good-unknown}.  Accordingly, after the
normalization \cref{eq:main_result_water_wave_potential_normalization}, the
order-one upper-right entry of
\cref{eq:main_result_Q_preview} is the usual water-wave
Dirichlet--Neumann contribution.  Likewise, the order-one transport and
Alinhac blocks in
\cref{eq:main_result_paradifferential_reduction,eq:main_result_B_preview}
reduce formally to the standard transport velocity and good-unknown
structure encoded by \(V\), \(B\), and
\(U=\psi_{\rm WW}-T_B\eta\) in
\cref{eq:historical-water-wave-good-unknown}.  Notice in particular that
letting \(\background\to0\) removes the prescribed affine background flow,
but not the nonlinear water-wave transport velocity.

The lower-left entry exhibits the other essential correspondence.  Since
\(\aangle{\rho^\pm}=\rho^-/2\) at \(\Atw=2\),
\cref{eq:Api_star_2} gives
\begin{equation}\label{eq:main_result_water_wave_capillary_symbol}
\frac12
a^{\pi_0,\eta}_2
\pare{\rho^\pm\middle|\eta;x,\xi}
=
-
\frac{\sigma}{\rho^-}
\frac{
\av{\xi}^2+
\av{\nabla^\perp\eta\cdot\xi}^2
}{
\pare{1+\av{\nabla\eta}^2}^{\frac32}
},
\end{equation}
which is the usual order-two symbol of the linearized mean-curvature
operator, with the physical one-phase surface-tension coefficient
\(\sigma/\rho^-\).  After dividing the second equation by the factor \(2\)
coming from
\(\tilde\psi=2\psi_{\rm WW}\), the gravity contribution
\(-2\mathfrak g\) similarly reduces to the standard
\(-\mathfrak g\) term.

Most importantly, the additional order-one contribution specific to the
two-phase problem disappears at the one-phase endpoint.  Indeed,
\cref{eq:qpitilde_main} gives the nonlinear identity
\begin{equation}\label{eq:main_result_water_wave_KH_cancellation}
q^{\tilde\pi_0}
\pare{
2,0
\middle|
\eta,2\psi_{\rm WW};x,\xi
}
=0,
\end{equation}
because \(q^{\tilde\pi}_0\) contains the factor \(4-\Atw^2\).  Thus the
order-one lower-left term
\(\lambda_1q^{\tilde\pi}_0\) in
\cref{eq:main_result_Q_preview} vanishes identically.  At the flat level
the same mechanism is visible directly in
\cref{eq:main_result_flat_q_preview}, since

$$
    \frac{4-\Atw^2}{8}
    \frac{
        \pare{\background^\perp\cdot\xi}^2
    }{
        \av{\xi}
    }
    \longrightarrow0
    \qquad
    \textnormal{as }
    \pare{\Atw,\background}\longrightarrow\pare{2,0}.
$$

This is precisely the Kelvin--Helmholtz order-one contribution which has no
counterpart in the one-phase water-wave equations.

Consequently, at the level of all positive-order terms,
\cref{eq:main_result_paradifferential_reduction} formally reduces at
\(\pare{\Atw,\background}=\pare{2,0}\), after the normalization
\(\tilde\psi=2\psi_{\rm WW}\), to the standard nonlinear
gravity--capillary water-wave paradifferential structure recalled in
\cref{eq:historical-water-wave-paradifferential-system}; see
\cite{AlazardBurqZuily2011,ABZ2014,BD2018}.  Moreover, because the
two-phase order-one Kelvin--Helmholtz contribution has vanished before any
use of surface tension, one may formally set \(\sigma=0\) in this endpoint:
the order-two curvature term in
\cref{eq:main_result_water_wave_capillary_symbol} then disappears and one
recovers the corresponding unbounded quasilinear structure of the
pure-gravity water-wave system of \cite{ABZ2014}.  The bounded
order-zero terms, including the nonlinear Taylor coefficient in the
pure-gravity formulation, are not identified by the present comparison;
the statement concerns precisely the positive-order structure responsible
for the derivative balance. 
\fine
\end{rem}

\subsubsection{Local well-posedness}

At the flat state, the first line of
\cref{eq:main_result_paradifferential_reduction} reduces to the linear symbol
in \cref{eq:linearized_CDE_matrix_symbol}.  Fix
\(\sigma>0\), \(\mathfrak g>0\), and \(\rho_\star>0\).  For this choice of
physical parameters, define the renormalized densities
\begin{equation}\label{eq:main_density_parametrization}
    \rho^-_\Atw
    \defeq
    \rho_\star\pare{1+\frac{\Atw}{2}},
    \qquad
    \rho^+_\Atw
    \defeq
    \rho_\star\pare{1-\frac{\Atw}{2}},
    \qquad
    0<\Atw\leq2.
\end{equation}
With \(\rho^\pm=\rho^\pm_\Atw\), we isolate the parameter region in which
the two dispersive branches are separated and purely imaginary.

\begin{hyp}[Strict Atwood--Weber stability]
\label{hyp:strict_atwood_weber_stability}
\begin{equation}\label{eq:strict_atwood_weber_number_intro}
    \sigma>0,
    \qquad
    \mathfrak g\Atw>0,
    \qquad
    \mathsf W_{\Atw,\background}
    \defeq
    \frac{
        \kappa_{\rm We}^{\max}
    }{
        2\kappa_{\rm cap}
    }
    =
    \frac{4-\Atw^2}{16}
    \av{\background}^2
    \sqrt{
        \frac{
            \aangle{\rho^\pm}
        }{
            \sigma\,\mathfrak g\Atw
        }
    }
    <1 .
\end{equation}
\begin{itemize}
    \item \(\sigma\) is the surface tension coefficient.  In particular,
    capillarity is present at every nonzero frequency.

    \item \(\mathfrak g\Atw\) is the flat Rayleigh--Taylor restoring
    acceleration.  In particular, the heavier fluid lies below the lighter
    one when \(\mathfrak g\Atw>0\).

    \item \(\kappa_{\rm cap}\) is the positive capillary wavenumber defined by
    $$\kappa_{\rm cap}\defeq\sqrt{\frac{
            \pare{\rho^- - \rho^+}\mathfrak g
        }{
            \sigma
        }}\cdot$$
    In particular,
    \(\kappa_{\rm cap}>0\) under \(\mathfrak g\Atw>0\).

    \item \(\kappa_{\rm We}^{\max}\) is the maximal effective Weber
    wavenumber generated by the constant background shear and defined by
    $$\kappa_{\rm We}^{\max}\defeq\frac{
            \rho^- \rho^+
        }{
            \sigma\pare{\rho^- + \rho^+}
        }
        \av{\background}^2.$$
    In particular,
    \(\kappa_{\rm We}^{\max}=0\) if \(\background=0\) or if one of the two
    densities vanishes.
\end{itemize}
With these notations, the stability criterion \cref{eq:strict_atwood_weber_number_intro} is equivalent to
\begin{equation}\label{eq:strict_linear_stability_scales_intro}
    \pare{\kappa_{\rm We}^{\max}}^2
    <
    4\kappa_{\rm cap}^2 .
\end{equation}
By \cref{prop:purely_imaginary_linear_frequencies}, this implies that the two spectral branches of the flat linearized CDE are purely imaginary
and simple at every nonzero frequency.
\end{hyp}

For the Cauchy statement, we introduce the phase space
\begin{equation}\label{eq:main_phase_space}
    \mathcal X^s
    \defeq
    H^{s+\frac14}\pare{\bR^2;\bR}
    \times
    H^{s-\frac14}\pare{\bR^2;\bR}.
\end{equation}

\begin{definition}[Stable parameter set]
\label{def:main_stable_parameter_set}
Let \(\sigma>0\), \(\mathfrak g>0\), and \(\rho_\star>0\) be fixed. We define
the \emph{stable parameter set} by
\begin{equation}\label{eq:main_stable_parameter_set}
    \mathcal S_{\rm st}
    \defeq
    \set{
        \pare{\Atw,\background}\in\pare{0,2}\times\bR^2
        \ \middle|\
        \mathfrak g\Atw>0,
        \quad
        \mathsf W_{\Atw,\background}<1
    }
    \cup
    \pare{\set{2}\times\bR^2}.
\end{equation}
Thus \(\mathcal S_{\rm st}\) consists of the two-phase parameters satisfying
the strict Atwood--Weber stability condition, together with the one-phase
endpoint \(\Atw=2\).
\end{definition}

\begin{rem}[Linear stability and strict hyperbolicity]
\label{rem:main_stable_parameter_set_linear_stability}
The terminology in \cref{def:main_stable_parameter_set} is justified by the
linearized analysis. Indeed, by
\cref{prop:purely_imaginary_linear_frequencies}, one has
\begin{equation}\label{eq:main_stable_parameter_set_spectral_characterization}
    \pare{\Atw,\background}\in\mathcal S_{\rm st}
    \quad\Longleftrightarrow\quad
    \mathfrak d\pare{\Atw,\background\middle|\xi}<0
    \quad
    \textnormal{for every }
    \xi\in\bR^2\setminus\set{0},
\end{equation}
 where $\mathfrak d\pare{\Atw,\background\middle|\xi}$ is defined in \cref{eq:def_linearized_discriminant}.
Equivalently, the two eigenvalues
\(\lambda_\pm\pare{\Atw,\background\middle|\xi}\) of the flat linearized
CDE, computed in \cref{sec:linearization-trivial-state}, are purely imaginary
and distinct at every nonzero frequency, or, equivalently,
\(-\ii{\bm L}_0\pare{\xi}\) has two distinct real eigenvalues.
Hence \(\mathcal S_{\rm st}\) is precisely the parameter region of strict
full-frequency linear hyperbolicity. \fine
\end{rem}

The uniform gap provided by
\cref{hyp:strict_atwood_weber_stability} allows the positive-order block in
\cref{eq:main_result_paradifferential_reduction} to be symmetrized and
diagonalized; this reduction is carried out in
\cref{sec:diagonal_reduction_main}.  The resulting energy estimate yields the
following local well-posedness theorem for the canonical CDE.

\begin{theorem}[Uniform local well-posedness in the Hamiltonian scale]
\label{thm:main}
Let \(\sigma>0\), \(\mathfrak g>0\), and \(\rho_\star>0\).
There exists \(s_0>0\) such that, for every \(s\geq s_0\) and every
compact set \(K\Subset\mathcal S_{\rm st}\), there exist
\(\varepsilon_0,c_{s,K}>0\), depending only on \(s\) and \(K\), such that
the following holds.

For every \(\pare{\Atw,\background}\in K\) and every real-valued datum
\(Y_0\defeq\pare{\eta_0,\tilde\psi_0}\in\mathcal X^s\) satisfying
\[
    \norm{Y_0}_{\mathcal X^s}\leq\varepsilon<\varepsilon_0,
\]
the canonical CDE \cref{eq:CDE_main_intro}, with densities
\(\rho^\pm=\rho^\pm_\Atw\), admits a unique real-valued classical solution
\[
    Y\defeq\pare{\eta,\tilde\psi}
    \in
    \cC\pare{\bra{-T_\varepsilon,T_\varepsilon};\mathcal X^s}
    \cap
    \cC^1\pare{\bra{-T_\varepsilon,T_\varepsilon};
    \mathcal X^{s-\frac32}},
    \qquad
    Y\pare{0}=Y_0,
\]
with
\[
    T_\varepsilon\geq c_{s,K}\varepsilon^{-1}.
\]
\end{theorem}

\begin{remark}
\label{rem:main_result_comments_3d_vs}
The strict Atwood--Weber threshold
\(\mathsf W_{\Atw,\background}<1\) is the CDE spectral form of the
Kelvin--Rayleigh--Taylor stability balance isolated by Lannes
\cite[eq.~(1.1)]{Lannes2013}.  Let us spell out the comparison.  In
our sign convention \(\bbra{f^\pm}=f^- - f^+\), Lannes' nonlinear
criterion reads
\begin{equation}\label{eq:main_result_lannes_criterion_physical}
    \bbra{-\partial_z p^\pm\big|_{z=\zeta}}
    >
    \frac{1}{4\sigma}
    \frac{\pare{\rho^-\rho^+}^2}{
        \pare{\rho^-+\rho^+}^2
    }
    c\pare{\zeta}\,
    \norm{\bbra{V^\pm\big|_{z=\zeta}}}_{L^\infty}^4,
\end{equation}
where \(c\pare{\zeta}\) is the geometric prefactor entering the
two-fluid symmetrizer.  The left-hand side is the Rayleigh--Taylor
restoring force, while the right-hand side measures the destabilizing
Kelvin--Helmholtz shear, regularized by capillarity.

In the flat infinite-depth configuration one has \(c\pare{0}=1\), and
the two quantities entering \cref{eq:main_result_lannes_criterion_physical}
become exactly the constant coefficients appearing in the linearized CDE.
More precisely,
\begin{equation}\label{eq:main_result_lannes_flat_identifications}
    \bbra{-\partial_z p^\pm\big|_{z=0}}
    =
    \pare{\rho^- - \rho^+}\mathfrak g,
    \qquad
    \norm{\bbra{V^\pm\big|_{z=0}}}_{L^\infty}
    =
    \av{\background},
\end{equation}
up to the fixed \(\bm J\)-rotation built into the vorticity variable
\(\Omega=\nabla^\perp\bbra{\psi^\pm}\).  Therefore
\cref{eq:main_result_lannes_criterion_physical} reduces to
\begin{equation}\label{eq:main_result_lannes_comparison}
    \pare{\rho^- - \rho^+}\mathfrak g
    >
    \frac{1}{4\sigma}
    \frac{\pare{\rho^-\rho^+}^2}{
        \pare{\rho^-+\rho^+}^2
    }
    \av{\background}^4.
\end{equation}
On the other hand, since
\begin{equation}\label{eq:main_result_atwood_density_identity}
    \frac{4-\Atw^2}{16}
    =
    \frac{\rho^-\rho^+}{\pare{\rho^-+\rho^+}^2},
    \qquad
    \Atw=\frac{2\pare{\rho^- - \rho^+}}{\rho^-+\rho^+},
\end{equation}
the inequality \(\mathsf W_{\Atw,\background}<1\) is equivalent to
\cref{eq:main_result_lannes_comparison}.

Thus the two criteria express the same principal stability mechanism:
Rayleigh--Taylor restoring force and capillarity must dominate the
Kelvin--Helmholtz shear.  In \cite{Lannes2013} this balance is formulated
as a nonlinear stability condition in physical variables, involving the
actual interfacial velocity jump
\(\bbra{V^\pm\big|_{z=\zeta}}\).  In the present theorem the same balance
is obtained from the eigenvalue calculation of the flat CDE around a state
carrying the prescribed affine background \(\background\).  The point is
therefore not that the condition is a special case of Lannes' theorem, nor
that Lannes' criterion excludes nonzero shear; rather, the comparison shows
that our Atwood--Weber threshold is the constant-background spectral
realization, in CDE variables, of Lannes' Kelvin--Rayleigh--Taylor
criterion.  This is precisely the stability gap used below to construct the
complex symmetrizer and diagonalizer.
\fine
\end{remark}

\subsection{Ideas of the proof}
\label{sec:ideas_of_the_proof_3d}

The proof follows the three layers used to present the main results.  The
first and the third are tied to the vortex-sheet system, whereas the middle
one is an operator-level statement which is formulated and proved
independently of it.  We first reduce the Euler equations to an autonomous
contour-dynamics system in canonical variables.  We then isolate the graph
integral operators generated by this reduction and paralinearize them through
a single abstract mechanism.  Only after this singular-integral analysis has
been completed do we return to the fluid system, assemble the real
paradifferential CDE, and perform the diagonal reduction compatible with the
non-homogeneous energy at frequency zero.  Schematically,
\begin{equation}\label{eq:ideas_proof_reduction_chain_3d}
\begin{aligned}
    \pare{{\bm u}^\pm,p^\pm,{\bm\gamma}}
    &\longmapsto
    \pare{\eta,\tilde\psi},
    \\
    \pare{h,k},\ 0<n\leq \mathtt d
    &\longmapsto
    \mathpzc t^{\pare{n}}\pare\eta\bra\omega
    \longmapsto
    \textnormal{structure-preserving paradifferential reduction},
    \\
    \pare{\eta,\tilde\psi}
    &\longmapsto
    \textnormal{paradifferential CDE}
    \longmapsto
    Y_{\rm Al}
    \longmapsto
    U
    \longmapsto
    W.
\end{aligned}
\end{equation}
The second line is the model-independent layer of the argument.  It extracts
the local homogeneous structure responsible for the symbols from the
particular fluid model, and is established on \(\bR^{\mathtt d}\) before the
vortex-sheet kernels are inserted.  The critical endpoint is the case used
in the vortex-sheet application; the subcritical result in the appendix
records the robustness of the same construction beyond the principal-value regime.

\paragraph{From the Euler system to an autonomous CDE.}
The starting point is the concentration of the vorticity on the interface.
Since the velocity fields are irrotational away from the sheet, the
Biot--Savart law expresses their traces through the three-dimensional
Birkhoff--Rott integral.  Once the interface is written as a graph, the
vorticity has the form
\(\bm\omega=\vect{\Omega}{\nabla\eta\cdot\Omega}\), and the average trace of
the velocity is described entirely by the two singular operators
\(\mathpzc B_0\pare{\eta}\) and \(\mathpzc q\pare{\eta}\) of
\cref{eq:B0_intro,eq:Q_intro}.  In particular, the kinematic condition becomes
\begin{equation}\label{eq:ideas_proof_kinematic_vorticity_3d}
    \eta_t
    =
    \mathpzc g^\star\pare{\eta}\bra{\Omega}.
\end{equation}
At this stage the system is already supported on the interface, but it is not
yet closed in the canonical variables, since \(\Omega\) is still an
independent unknown.

The density-weighted trace potential \(\psi\) in \cref{def:HAMpsi_intro} is precisely the variable that
allows us to eliminate \(\Omega\) as an independent unknown.  Indeed, writing
\(\varphi\defeq\bbra{\psi^\pm}\), the Plemelj jump relations give
\[
    \psi
    =
    \pare{\Id+\Atw\tau_0\pare{\eta}}\varphi,
    \qquad
    \Omega
    =
    \nabla^\perp\varphi.
\]
This is a second-kind boundary integral equation closely related to the one
used in \cite{CCG2011} to recover the vortex-sheet amplitude in the
two-dimensional Muskat problem.  There are, however, substantial differences
between the two settings.  The analysis of \cite{CCG2011} concerns a
one-dimensional interface in a two-dimensional fluid domain and exploits
tools specific to that geometry to obtain a non-perturbative inversion
result.  Here the interface is two-dimensional and the corresponding operator
acts on functions on \(\bR^2\).  We work perturbatively around the flat
interface: the Sobolev invertibility of
\(\Id+\Atw\tau_0\pare{\eta}\) follows from the smallness of \(\eta\), while
the precise structure of its inverse required in the subsequent analysis is
recovered by a multidimensional paradifferential parametrix construction,
see \cref{lem:uniform_inverse_tau0,prop:paralinearization_S}.  Thus the
mechanism in \cite{CCG2011} provides a useful structural analogy, but does
not directly apply to the present three-dimensional problem.  We therefore set
\begin{equation}\label{eq:ideas_proof_potential_to_vorticity_3d}
    \Omega
    =
    S\pare{\Atw\middle|\eta}\psi,
    \qquad
    S\pare{\Atw\middle|\eta}\psi
    \defeq
    \nabla^\perp
    \pare{\Id+\Atw\tau_0\pare{\eta}}^{-1}\psi,
\end{equation}
and substitute \cref{eq:ideas_proof_potential_to_vorticity_3d} into the
kinematic and Bernoulli equations.  This gives a closed system for
\(\pare{\eta,\psi}\).  Equivalently, the generalized Dirichlet--Neumann
operator appearing in the first equation of \cref{eq:CDEeta_main_intro} is
not introduced as an abstract elliptic boundary operator, but is factorized
as
\begin{equation}\label{eq:ideas_proof_dn_factorization_3d}
    \mathpzc g\pare{\Atw\middle|\eta}
    =
    \mathpzc g^\star\pare{\eta}
    \circ
    S\pare{\Atw\middle|\eta}.
\end{equation}
Thus both pieces of \(\mathpzc g\) remain explicit in terms of the
Birkhoff--Rott kernels.  This is important later: rather than invoking the
symbolic theory of the Dirichlet--Neumann operator, we recover its
paradifferential structure directly from the underlying singular integrals.
The factorization in \cref{eq:ideas_proof_dn_factorization_3d} is also the
point at which the fluid problem meets the abstract theory: its nonlocal
building blocks are critical graph-kernel operators associated with the
admissible homogeneous pairs of
\cref{def:abstract-admissible-angular-pair-Rd}.

The constant background requires one further step.  Its potential
\(\bar\psi=x^\perp\cdot\background\) is affine and does not belong to the
Sobolev phase space.  We therefore split
\(\psi=\bar\psi+\tilde\psi\), absorb the background contribution into
\begin{equation}\label{eq:ideas_proof_background_vorticity_3d}
    \tilde\Omega\bra{\background}\pare{\eta}
    \defeq
    S\pare{\Atw\middle|\eta}\bar\psi,
\end{equation}
and subtract the constant Bernoulli gauge \(c_{\background}\).  The flat
state is then represented by \(\pare{\eta,\tilde\psi}=\pare{0,0}\), while
the full background vorticity remains present through
\(\tilde\Omega\bra{\background}\pare{\eta}\).  This gives exactly the
autonomous system \cref{eq:CDE_main_intro}: all bulk quantities and the
independent vorticity have disappeared, and every term in the right-hand side
is a function of \(\pare{\eta,\tilde\psi}\) and of the fixed parameters.

Let us also stress how this reduction differs from the contour formulation
used by Fefferman, Ionescu, and Lie in \cite{fefferman2016absence}.  There the
interface is a curve in a two-dimensional fluid domain, no preferred
parametrization of the curve is fixed, and the boundary vorticity is retained
in the singular-integral representation in order to study the geometry of a
possible self-intersection.  The one-dimensional character of the interface
also entails structural simplifications which are no longer available in the
present setting: there is a single tangential direction, whereas for a surface
in three-dimensional space the tangential geometry is genuinely
two-dimensional and the associated singular and paradifferential operators
carry a nontrivial angular dependence in frequency space.  Here we work in
the graph gauge and eliminate the vorticity in favour of the canonical trace
potential.  The distinction is therefore not that one formulation is a
contour equation and the other is not.  It lies both in the geometry of the
interface and in the choice of variables and purpose of the reduction: our
goal is an autonomous CDE on a fixed Sobolev phase space, suitable for a
direct quasilinear analysis.

The important output of this first reduction is not only the closed system
itself.  Its nonlocal geometry has been concentrated in a finite family of
graph operators whose dependence on the interface occurs through normalized
finite differences.  This is the model-independent structure isolated in the
next step.

\paragraph{A structure-preserving paralinearization theorem for graph kernels.}
This is the main methodological contribution of the paper.  The singular
operators produced by the contour reduction are not paralinearized one at a
time.  Instead, their local geometry is encoded by an admissible homogeneous
pair \(\pare{h,k}\) of order \(-n\) and vanishing order
\(\ell\), with \(0<n\leq\mathtt d\), in the sense of
\cref{def:abstract-admissible-angular-pair-Rd}.  Here \(h\) is the angular
profile of the homogeneous kernel, whereas \(k\) describes its nonlinear
dependence on the graph. 
At the level relevant for the symbolic analysis, every operator has the
localized homogeneous core
\begin{equation}\label{eq:ideas_proof_model_singular_kernel_3d}
\begin{aligned}
   \mathpzc t_{\loc}^{\pare{n}}\pare\eta\bra\omega\pare x
    &\defeq \text{p.v.}
    \int_{\bR^{\mathtt d}}
    \chi\pare z \frac{h\pare{\hat z}}{\av{z}^n}
    k\pare{\Delta_z\eta\pare x}
    \omega\pare{x-z}\dd z,
    \\
    \Delta_z\eta\pare x
    &\defeq
    \frac{\eta\pare x-\eta\pare{x-z}}{\av z},
    \qquad
    \hat z\defeq\frac z{\av z}\cdot
\end{aligned}
\end{equation}
The space localization is obtained through the cut-off function $\chi$ as in \cref{def-chi}. When \(\beta =  \mathtt{d}-n =0\), namely \(n=\mathtt d\), the complete flat kernel is
interpreted as a homogeneous principal-value distribution.  The opposite
parities imposed on \(h\) and \(k\) make the full graph kernel odd under
the exchange of the two spatial variables and yield the exact
skew-adjointness of
\cref{lem:abstract-skew-adjoint-graph-kernel-Rd}.  When \(\beta>0\), the
homogeneous core is locally integrable and no parity condition is required.
The global version treated in
\cref{app:abstract-locally-integrable-paralinearization-Rd} is obtained by
completing this local homogeneous kernel with a sufficiently regular
{
far-field tail. The nonlinear tail contributes an order-zero paraproduct
and smoothing remainders, without changing the local principal symbols.
}
In either regime the perturbative assumption is imposed only on the graph,
so that the difference quotients remain inside the analytic range of \(k\);
no smallness is imposed on \(\omega\).

This formulation isolates and extends the mechanisms introduced in
\cite{BCGS2025,murgante2025long}.  Its point is not merely to shorten the
proofs of the two operators appearing below.  It identifies, in arbitrary
dimension, which features of the paralinearization come from homogeneity,
which come from the vanishing of \(k-k\pare0\), and which are specific to
the critical endpoint.  The singular dependence on \(z\) and the nonlinear
dependence on \(\Delta_z\eta\) occur at the same scale, and this is precisely
why a standard paralinearization performed after the integration in \(z\)
does not reveal the correct symbolic structure.

The proof is organized around four operations which are reusable beyond the
present fluid problem.

\begin{enumerate}

 \item \emph{Isolation of the homogeneous core.}
Set \(\widetilde k\defeq k-k\pare0\).  The first step is to separate the
contribution of the kernel near \(z=0\) from its far-field tail.  This
{
distinction separates the singular principal symbols from the order-zero
tail contribution and its smoothing remainders.
}

For the critical operator, the decomposition is explicit:
{
\begin{equation}\label{eq:ideas_proof_local_tail_splitting_3d}
\begin{aligned}
    \mathpzc t^{\textnormal{flat}}\bra\omega\pare x
    &\defeq
    k\pare0\,
    \pv\int_{\bR^{\mathtt d}}
    \frac{h\pare{\hat z}}{\av z^{\mathtt d}}
    \omega\pare{x-z}\dd z,
    \\
    \mathpzc t^{\textnormal{loc}}\pare\eta\bra\omega\pare x
    &\defeq
    \int_{\bR^{\mathtt d}}
    \chi\pare z
    \frac{h\pare{\hat z}}{\av z^{\mathtt d}}
    \widetilde k\pare{\Delta_z\eta\pare x}
    \omega\pare{x-z}\dd z,
    \\
    \mathpzc t^{\textnormal{tail}}\pare\eta\bra\omega\pare x
    &\defeq
    \int_{\bR^{\mathtt d}}
    \pare{1-\chi\pare z}
    \frac{h\pare{\hat z}}{\av z^{\mathtt d}}
    \widetilde k\pare{\Delta_z\eta\pare x}
    \omega\pare{x-z}\dd z.
\end{aligned}
\end{equation}
}
Thus
{
\[
    \mathpzc t^{\pare{\mathtt d}}\pare\eta\bra\omega
    =
    \mathpzc t^{\textnormal{flat}}\bra\omega
    +
    \mathpzc t^{\textnormal{loc}}\pare\eta\bra\omega
    +
    \mathpzc t^{\textnormal{tail}}\pare\eta\bra\omega.
\]
}
The flat term retains the complete homogeneous principal-value multiplier.
{
The local term contains the singularity at \(z=0\) and produces the
principal and subprincipal symbols.
}  It is this term that is analyzed by the Bony decomposition and
the local Taylor expansion below.

By contrast, the last integral only sees the far field.  On
\(\av z\geq1\),
\(
    \av{\Delta_z\eta\pare x}
    \lesssim
    \av z^{-1}\norm\eta_{L^\infty}
\),
{
so the vanishing order $\ell$ of $\widetilde k$ supplies the decay
$\av z^{-\ell}$. By
\cref{eq:abstract-tail-integrability-calculation,prop:tail-kernel-smoothing},
for every prescribed $\varrho\ge0$,
\[
\begin{aligned}
    \mathpzc t^{\textnormal{tail}}\pare\eta\bra\omega
    ={}&
    \OpBW{f_{\geq\ell}^{\textnormal{tail}}\bra\omega\pare\eta}\eta
    \\
    &+
    r_{\bra{\varrho,\geq\ell}}^{\textnormal{tail}}\pare\eta\omega
    +
    r_{\bra{\varrho,\geq\ell}}^{\textnormal{tail}}\bra\omega\pare\eta\eta.
\end{aligned}
\]
The function is given in \cref{eq:abstract-common-tail-function}.
It enters the generic order-zero coefficient of $\eta$.
}

The subcritical case is organized in the same way: the homogeneous
near-field component produces the finite-order symbols, while the
non-homogeneous far-field contribution is treated by the corresponding
tail estimates in
\cref{thm:abstract-locally-integrable-paralinearization-Rd}.

    \item \emph{Bony decomposition at fixed convolution scale.}
    For \(z\neq0\), set
    \(\omega_z\pare x\defeq\omega\pare{x-z}\).  The localized nonlinear
    product is decomposed before integration in \(z\):
    \begin{equation}\label{eq:ideas_proof_abstract_Bony_splitting_3d}
    \begin{aligned}
        \widetilde k\pare{\Delta_z\eta}\omega_z
        ={}&
        \OpBW{\widetilde k\pare{\Delta_z\eta}}\omega_z
        +
        \OpBW{\omega_z}\bra{\widetilde k\pare{\Delta_z\eta}}
       +
        \RBW{\widetilde k\pare{\Delta_z\eta}}{\omega_z}.
    \end{aligned}
    \end{equation}
    The three terms are respectively the quasilinear branch, in which the
    high frequency remains on \(\omega\); the transported branch, in which
    the high frequency is transferred to \(\eta\); and the localized Bony
    remainder.  The latter becomes arbitrarily smoothing after integration
    in \(z\), by
    \cref{lem:abstract-localized-Bony-remainder-Rd} in the critical case and
    \cref{lem:abstract-locally-integrable-localized-Bony-remainder-Rd} in the
    subcritical case.  Performing this decomposition before the
    \(z\)-integration is essential: the powers of \(\av z\) which measure
    the Taylor order of the graph are exactly those which determine the
    order of the resulting Fourier symbols.

    \item \emph{Extraction of the symbols from the local graph geometry.}
    On the support of \(\chi\), Taylor's formula gives
    \begin{equation}\label{eq:ideas_proof_local_taylor_expansion_3d}
        \Delta_z\eta\pare x
        =
        \nabla\eta\pare x\cdot\hat z
        -
        \frac{\av z}{2}
        \hat z^\intercal{\bm H}\eta\pare x\hat z
        +
        e_{\bra2}\pare{\eta;x,z},
    \end{equation}
    where \(e_{\bra2}\) has local order two by
    \cref{eq:abstract-quasilinear-difference-Taylor-expansion}.  Therefore
    \begin{equation}\label{eq:ideas_proof_composed_taylor_expansion_3d}
    \begin{aligned}
        k\pare{\Delta_z\eta}
        ={}&
        k\pare{\nabla\eta\cdot\hat z}
        -
        \frac{\av z}{2}
        k'\pare{\nabla\eta\cdot\hat z}
        \pare{\hat z^\intercal{\bm H}\eta\hat z}
        +
        k_{\bra2}\pare{\eta;x,z}.
    \end{aligned}
    \end{equation}
    The basic order conversion is
    \begin{equation}\label{eq:ideas_proof_kernel_to_symbol_order_3d}
        \textnormal{local kernel order }n+j
        \quad\longmapsto\quad
        \textnormal{symbol order }-\beta-j.
    \end{equation}
    The first term in
    \cref{eq:ideas_proof_composed_taylor_expansion_3d} gives the
    order-\(\pare{-\beta}\) quasilinear symbol, while the transported branch gives
    an order-one symbol linear in \(\omega\).  The remainder
    \(k_{\bra2}\) begins two local orders later and is therefore treated by
    \cref{prop:int-opbw} at order \(-\beta-2\), modulo an arbitrarily
    smoothing operator.

After freezing
\(p=\nabla\eta\pare x\), the leading local kernels depend on \(z\) only
through the homogeneous radial factor and the angular variable
\(\vartheta=\hat z\).  Analyticity of \(k\) allows us to separate this
angular dependence from the \(x\)-dependent coefficient:
\begin{equation}\label{eq:ideas_proof_angular_analytic_expansion_3d}
    k\pare{p\cdot\vartheta}
    =
    \sum_{\alpha\in\bN^{\mathtt d}}
    \frac{k^{\pare{\av\alpha}}\pare0}{\alpha!}
    p^\alpha\vartheta^\alpha.
\end{equation}
Thus each term of the series is the product of a coefficient \(p^\alpha\)
and a convolution operator with a fixed homogeneous angular amplitude
\(h\pare\vartheta\vartheta^\alpha\).  For instance, the leading
quasilinear branch is reduced term by term to operators of the form
\[
    p^\alpha\,
    \pv\int_{\bR^{\mathtt d}}
    \chi\pare z
    \frac{
        h\pare{\hat z}\hat z^\alpha
    }{
        \av z^{\mathtt d}
    }
    \omega\pare{x-z}\dd z.
\]
The parity assumptions in
\cref{def:abstract-admissible-angular-pair-Rd} ensure that every nonzero
Taylor coefficient produces precisely the odd or even angular amplitude
required by the homogeneous multiplier formulas of
\cref{prop:abstract-homogeneous-angular-multipliers-Rd}.  The transported
branch is reduced in the same way, after the symmetrization by
\(\delta_z^{\pare2}\eta\), to compensated homogeneous kernels of degree
\(-\mathtt d-1\).

At this point the cutoff no longer carries symbolic information.  Indeed,
the localized kernel and its complete homogeneous model have the same
singularity at \(z=0\), and
\cref{lem:localized-homogeneous-model,%
lem:localized-compensated-homogeneous-model}
show that their high-frequency Fourier multipliers differ only by a
rapidly decreasing symbol:
%\textcolor{red}
{\begin{equation}\label{eq:ideas_proof_localized_vs_complete_homogeneous_3d}
    t_{h,\alpha}^{\textnormal{loc}}f
    =
    \dot{\mathfrak m}_{h,\alpha}\pare Df
    +
    \mathfrak s_{h,\alpha}\pare Df,
    \qquad
    \mathfrak s_{h,\alpha}
    \in
    \bigcap_{M\in\bN}\tilde\Gamma^{-M},
\end{equation}
where
$$t_{h,\alpha}^{\textnormal{loc}}\bra{\omega}\pare{x}\defeq\pv\int_{\bR^{\mathtt d}}
    \chi\pare z
    \frac{
        h\pare{\hat z}\hat z^\alpha
    }{
        \av z^{\mathtt d}
    }
    \omega\pare{x-z}\dd z,\qquad\dot{\mathfrak m}_{h,\alpha}\pare{D}\bra{\omega}\pare{x}\defeq\pv\int_{\bR^{\mathtt d}}
    \frac{
        h\pare{\hat z}\hat z^\alpha
    }{
        \av z^{\mathtt d}
    }
    \omega\pare{x-z}\dd z.$$
}
Consequently, all finite-order symbolic information can be computed from
the complete homogeneous kernels.  The degree-\(\pare{-\mathtt d}\) odd kernels
give purely imaginary multipliers of order zero, whereas the compensated
even kernels arising in the transported branch give real multipliers of
order one.  Applying
\cref{prop:abstract-homogeneous-angular-multipliers-Rd} term by term and
then summing back the analytic expansions reconstructs \(k\) and \(k'\),
and yields
\begin{equation}\label{eq:ideas_proof_angular_multipliers_3d}
\begin{aligned}
    a_0^{h,k;\omega}\pare{p,\xi}
    &\defeq
    -\frac\pi2
    \int_{\bS^{\mathtt d-1}}
    \sgn\pare{\vartheta\cdot\xi}
    h\pare\vartheta
    k\pare{p\cdot\vartheta}
    \dd\sigma\pare\vartheta,
    \\
    a_1^{h,k;\eta}\pare{p,\xi}
    &\defeq
    \frac\pi2
    \int_{\bS^{\mathtt d-1}}
    \av{\vartheta\cdot\xi}
    h\pare\vartheta
    k'\pare{p\cdot\vartheta}
    \dd\sigma\pare\vartheta,
    \\
    \partial_\xi a_1^{h,k;\eta}\pare{p,\xi}
    &=
    -\partial_p
    a_0^{h,k;\omega}\pare{p,\xi}.
\end{aligned}
\end{equation}

    The last identity is
    \cref{eq:abstract-principal-symbol-differential-identity}.  It shows
    that the two principal columns are not unrelated outputs of the
    calculation: they are different manifestations of the same graph
    kernel.  In the subcritical regime, the corresponding principal symbols
    are the localized Fourier integrals in
    \cref{eq:abstract-locally-integrable-transported-symbol-Rd,eq:abstract-locally-integrable-quasilinear-symbol-Rd}.
    Their orders follow from the same local homogeneity, without an endpoint
    parity argument.

    The angular variable survives the Fourier transform.  The resulting
    symbols have the form
    \begin{equation}\label{eq:ideas_proof_genuine_two_dimensional_symbol_3d}
        a\pare{\nabla\eta\pare x,\hat\xi}\av\xi^\nu.
    \end{equation}
    In the present application, \(\mathtt d=2\) and
    \(\hat z\in\bS^1\), so the graph geometry couples to a continuum of
    frequency directions.  This is a genuinely surface-level feature: for a
    one-dimensional contour, \(\bS^0=\set{-1,1}\) and the corresponding
    angular structure collapses to two directions and remarkably simplifies the analysis.

\item \emph{Preservation of the critical adjoint structure.}
The first lower-order contribution in the quasilinear branch comes from two
different sources.  The principal symbol produces a first Bony--Weyl
correction, while the next term in the Taylor expansion of
\(\Delta_z\eta\) produces a geometric correction of the same order.
The two contributions are equal and opposite, as shown in
\cref{eq:abstract-quasilinear-composition-subprincipal,%
eq:abstract-quasilinear-geometric-subprincipal}.  Hence no independent
symbol of order \(-1\) survives in the critical case, and the first generic
lower-order term in the \(\omega\)-column has order \(-2\).

There is a simple geometric way to understand this cancellation.  Weyl
quantization is centered at the midpoint of the interaction.  In the
corresponding coordinates the graph difference takes the form
\[
    \frac{
        \eta\pare{X+z/2}-\eta\pare{X-z/2}
    }{\av z}
    =
    \nabla\eta\pare X\cdot\hat z
    +
    O\pare{\av z^2}.
\]
There is therefore no correction of order \(\av z\).  When the calculation
is first expanded around one endpoint, such a term appears separately as
the first geometric Taylor correction and the first Weyl correction; their
cancellation simply restores this centered structure.

In the critical regime, the parity assumptions have an additional role:
they make the principal-value operator exactly skew-adjoint.  Accordingly,
the principal Bony--Weyl operator is
\(\OpBW{\ii a_0^{h,k;\omega}}\), with
\(a_0^{h,k;\omega}\) real and odd in \(\xi\).  Thus the cancellation above
is also the symbolic manifestation of the adjoint structure of the original
singular integral.  The same centered-difference mechanism explains the
corresponding gap between orders \(-\beta\) and \(-\beta-2\) in the
subcritical case; there, however, no parity assumption is needed because
the kernel is already locally integrable.

\end{enumerate}

The fluid system enters only after this operator-level construction has been
completed.  By \cref{lem:B0-q-admissible-kernels}, the two scalar components
of \(\mathpzc B_0\pare\eta\) and the operator
\(\mathpzc q\pare\eta\) are critical instances with
\(\mathtt d=2\), \(n=2\), and therefore \(\beta=0\).  The two
components of \(\mathpzc B_0\) have vanishing order \(\ell=2\), whereas
\(\mathpzc q\) has vanishing order \(\ell=1\).  The former correspond to
an odd angular profile and an even geometric nonlinearity; the latter to an
even angular profile and an odd geometric nonlinearity.  Their
paralinearizations
\cref{prop:B0_paralinearization,prop:paralinearization_Q} are therefore
applications of the same critical theorem.  The model-specific part is
reduced to verifying admissibility, evaluating the angular integrals, and
assembling the composite operators of the CDE.

We then paralinearize
\(S\pare{\Atw\middle|\eta}\), the velocity traces, the generalized
Dirichlet--Neumann operators, and the Bernoulli terms, and combine them in
\cref{thm:main_result_paradifferential_reduction}.  As a consistency check,
the graph-kernel route recovers the geometric principal symbol
\begin{equation}\label{eq:ideas_proof_geometric_principal_symbol_3d}
\begin{aligned}
    2\lambda_1\pare{\eta;x,\xi}
    &=
    \sqrt{
        \pare{1+\av{\nabla\eta\pare{x}}^2}\av{\xi}^2
        -
        \pare{\nabla\eta\pare{x}\cdot\xi}^2
    }
    =
    \av{\xi}
    \sqrt{
        1+
        \av{\nabla^\perp\eta\pare{x}\cdot\hat\xi}^2
    }.
\end{aligned}
\end{equation}
Up to the normalization of the two-fluid potential, this is exactly the principal symbol of the Dirichlet-Neumann operator. The abstract approach recovers this familiar geometry without treating the Dirichlet-Neumann calculus as the point of departure, while keeping explicit control over the lower-order vorticity and background contributions. The vortex-sheet system is therefore the main application of a graph-kernel paralinearization method formulated before, and independently of, the model.

\paragraph{From the paradifferential system to the Cauchy theory.}
From this point onward no singular integral is treated directly.  We start
from the real system in
\cref{thm:main_result_paradifferential_reduction}.  Its positive-order part
contains the Alinhac block \({\bm B}\), which would lead to a loss of
derivatives if it were treated perturbatively.  We first remove it by the
real triangular transformation
\begin{equation}\label{eq:ideas_proof_real_alinhac_unknown_3d}
    Y_{\rm Al}
    \defeq
    {\bm G}Y,
    \qquad
    {\bm G}
    \defeq
    \Id_{\bC^2}
    +
    \OpBW{
        \begin{bmatrix}
            0&0
            \\
            b_{\textnormal{can}}&0
        \end{bmatrix}
    }.
\end{equation}
The matrix symbol in \cref{eq:ideas_proof_real_alinhac_unknown_3d} is
nilpotent, so the inverse is obtained exactly by changing its sign.  The
transformation leaves \(\eta\) unchanged and replaces \(\tilde\psi\) by the
usual good unknown.  It is important that this step is performed in real
variables: it removes the geometric derivative loss before the two
dispersive branches are mixed by the complexification; see
\cref{prop:real_alinhac_good_unknown_3d}.

Under \cref{hyp:strict_atwood_weber_stability}, the two branches of the flat
linearized system are simple and purely imaginary.  Recall $\mathfrak d\pare{\Atw,\background\middle|\xi}$ defined in \cref{eq:def_linearized_discriminant}, we denote the positive
dispersive frequency by
\begin{equation}\label{eq:ideas_proof_flat_positive_frequency_3d}
    \omega_{\frac32}
    \pare{\Atw,\background\middle|\xi}
    \defeq
    \sqrt{
        -\mathfrak d_{\Atw,\background}\pare{\xi}
    }.
\end{equation}
For every \(\xi\neq0\), the corresponding flat eigenvectors suggest the
complex coordinate
\begin{equation}\label{eq:ideas_proof_homogeneous_complex_unknown_3d}
    u_{\rm hom}
    \defeq
    \frac1{\sqrt2}
    \pare{
        \mathfrak m_{-\frac14}^{-1}\pare{D}\eta_{\rm Al}
        +
        \ii
        \mathfrak m_{-\frac14}\pare{D}\tilde\psi_{\rm Al}
    },
    \qquad
    \mathfrak m_{-\frac14}
    \pare{\Atw,\background\middle|\xi}
    \defeq
    \sqrt{
        \frac{\av{\xi}}{
            2\omega_{\frac32}
            \pare{\Atw,\background\middle|\xi}
        }
    }\cdot
\end{equation}
This choice diagonalizes the flat system exactly away from \(\xi=0\) and has
the correct high-frequency orders.  The obstruction is confined to low
frequencies.  Indeed, by
\cref{eq:def_linearized_q_background,eq:def_linearized_discriminant},
\begin{equation}\label{eq:ideas_proof_infrared_asymptotics_3d}
    \mathfrak q\pare{\Atw,\background\middle|\xi}
    \xrightarrow{\av{\xi}\to 0}
    -\mathfrak g\Atw,
    \qquad
    \omega_{\frac32}
    \pare{\Atw,\background\middle|\xi}
    \sim
    \sqrt{\frac{\mathfrak g\Atw}{2}}
    \av{\xi}^{\frac12},
    \qquad
    \mathfrak m_{-\frac14}
    \pare{\Atw,\background\middle|\xi}
    \sim
    \frac{\av{\xi}^{\frac14}}{
        \pare{2\mathfrak g\Atw}^{\frac14}
    }\cdot
\end{equation}
Hence
\(\mathfrak m_{-\frac14}^{-1}
\pare{\Atw,\background\middle|\xi}\)
grows like \(\av{\xi}^{-\frac14}\) as \(\xi\to0\).  The exact eigenvector
normalization therefore makes the map
\(\eta_{\rm Al}\mapsto
\mathfrak m_{-\frac14}^{-1}\pare{D}\eta_{\rm Al}\)
singular at the origin; in particular, it does not define a bounded Fourier
multiplier from \(H^{s+\frac14}\) to \(H^s\).  Equivalently, the squared
Fourier weight carried by the \(\eta_{\rm Al}\)-component behaves like
\(\av{\xi}^{-\frac12}\).

This singularity is not present in the physical quadratic energy.  At the
flat state the latter contains
\begin{equation}\label{eq:ideas_proof_flat_hamiltonian_energy_3d}
    \int_{\bR^2}
    \left(
        -\mathfrak q
        \pare{\Atw,\background\middle|\xi}
        \av{\widehat\eta\pare{\xi}}^2
        +
        \frac{\av{\xi}}{2}
        \av{\widehat{\tilde\psi}\pare{\xi}}^2
    \right)
    \dd\xi.
\end{equation}
Since
\(-\mathfrak q\pare{\Atw,\background\middle|\xi}
\to\mathfrak g\Atw>0\),
the gravitational term controls the low-frequency part of \(\eta\) with an
ordinary \(L^2\) weight.  Thus the factor
\(\av{\xi}^{-\frac14}\) is an artifact of the exact homogeneous eigenvector
normalization, rather than a degeneracy of the flat Hamiltonian.  The
homogeneous coordinate is therefore well adapted to the high-frequency
diagonalization, but not to the natural inhomogeneous Sobolev phase space.

For this reason we switch the homogeneous diagonalizer off near the origin.
With the cut-off \(\chi\) of \cref{def-chi}, we set
\begin{equation}\label{eq:ideas_proof_nonhomogeneous_complexification_3d}
    \mathfrak m_{-\frac14}^{\rm nh}
    \defeq
    \chi
    +
    \pare{1-\chi}\mathfrak m_{-\frac14},
    \qquad
    U
    \defeq
    \vect{u}{\bar u}
    =
    {\bm M}\pare{\Atw,\background\middle|D}Y_{\rm Al},
\end{equation}
where \({\bm M}\) is defined in \cref{eq:complex_matrix_M}.  This is a
bounded and boundedly invertible map between the real inhomogeneous phase
space and the complex Sobolev scale, and it preserves the conjugation
condition \(U=\vect{u}{\bar u}\).  The price is that the flat operator is no
longer exactly diagonal in the low-frequency region.  The discrepancy is the
explicit constant-coefficient multiplier
\({\bm A}^{\bC,{\rm lin}}_{\rm lo}\) of
\cref{eq:linear_low_frequency_defect_complex_3d}.

At high frequency the remaining diagonalization has a classical algebraic
core.  An order-zero matrix first diagonalizes the quasilinear block of order
\(3/2\).  The transport symbol has order one but is scalar in the two complex
components, and hence it is already diagonal.  After this first conjugation,
the only off-diagonal positive-order term has order \(1/2\).  It is removed
by an order-minus-one generator \({\bm S}_{-1}\), chosen so that, at the
homogeneous symbolic level,
\begin{equation}\label{eq:ideas_proof_subprincipal_homological_equation_3d}
    \comm{
        \ii
        \dot\omega_{\frac32}
        \dot\theta_0
        {\bm J}_{\rm d}
    }{
        \dot{\bm S}_{-1}
    }
    +
    \pare{{\bm L}_{\frac12}}^{\textnormal{off}}
    =0.
\end{equation}
This is the standard subprincipal mechanism: commutation with the separated
principal branches divides by their spectral gap and cancels the off-diagonal
half-order term.  In the notation of the proof, the commutator generates
exactly the opposite of the term containing
\(\dot\beta^{\bC,{\rm pr}}_0\av{\xi}^{1/2}{\bm J}\); the cancellation is
recorded in
\cref{eq:subprincipal_positive_half_diagonal_result_refined_3d}.

What is less standard is making this formal cancellation into an
admissible change of coordinates for the problem at hand.  Two kinds of
difficulty arise.  The first concerns the generator \(\dot{\bm S}_{-1}\)
itself.  Being homogeneous, it has to be cut off near \(\xi=0\).  It is
also built from the leading capillary part of the frequency, which
differs from the exact frequency by symbols that are lower order but not
negligible.  The cut-off then creates further terms, localized on a fixed
transition annulus.

The second concerns the resulting conjugation.  Each conjugating operator
must have a tame inverse on the non-homogeneous Sobolev scale, and its
time derivative must remain of an acceptable order.  The full
composition, moreover, must preserve the reality condition
\(U=\vect{u}{\bar u}\) of \cref{eq:ideas_proof_nonhomogeneous_complexification_3d}.
Finally, all these bounds must be uniform as \(\pare{\Atw,\background}\)
ranges over compact subsets of the stable set.

\cref{sec:positive_order_diagonalization_3d} carries out this program in
detail and produces an equivalent unknown \(W\) for which all
positive-order symbols are real and diagonal, as stated in
\cref{prop:positive_order_diagonalization_3d}.

The final step is the energy estimate.  At high frequency, if \(d_j\) is one
of the diagonal entries of the positive-order symbol, we choose the Sobolev
weight to be a scalar function \(\lambda_s\pare{d_j}\).  The potentially
order-losing Poisson bracket then vanishes identically:
\begin{equation}\label{eq:ideas_proof_poisson_cancellation_3d}
    \pbra{\lambda_s\pare{d_j}}{d_j}
    =
    \lambda_s'\pare{d_j}
    \pbra{d_j}{d_j}
    =0.
\end{equation}
This is the reason for diagonalizing the complete positive-order part rather
than using a fixed power of \(\angles{D}\) in the energy.

At low frequency we return to the real Hamiltonian structure.  Let
\begin{equation}\label{eq:ideas_proof_low_frequency_symmetrizer_3d}
    {\bm S}^{\bR}_{\rm lo}
    \defeq
    \begin{bmatrix}
        -\mathfrak q
        \pare{\Atw,\background\middle|\xi}
        &0
        \\
        0&\dfrac{\av{\xi}}{2}
    \end{bmatrix},
    \qquad
    {\bm S}^{\bC}_{\rm lo}
    \defeq
    \pare{{\bm M}^{-1}}^*
    {\bm S}^{\bR}_{\rm lo}
    {\bm M}^{-1}.
\end{equation}
The low-frequency quadratic form in the complex variable is therefore
exactly the pullback of the real quadratic Hamiltonian, not merely an
equivalent norm.  The isolated defect
\({\bm A}^{\bC,{\rm lin}}_{\rm lo}\) is not skew-adjoint on the transition
annulus, but it is not supposed to be considered in isolation: its
symmetrizer defect is cancelled by the transition part of the flat diagonal
dispersion.  The complete flat linear symbol contributes exactly zero to the
low-frequency energy balance, as shown in
\cref{eq:low_frequency_flat_linear_cancellation_3d}.

Combining the high-frequency cancellation
\cref{eq:ideas_proof_poisson_cancellation_3d} with the exact low-frequency
Hamiltonian identity gives the modified energy estimate
\begin{equation}\label{eq:ideas_proof_modified_energy_estimate_3d}
    \av{
        \frac{\dd}{\dd t}E_s\pare{t}
    }
    \lesssim_s
    E_{s_0}^{\frac12}\pare{t}
    E_s\pare{t}.
\end{equation}
The energy \(E_s\) is equivalent to the Hamiltonian low--high norm of
\cref{eq:hamiltonian_phase_space_norm_3d}.  It controls the full
inhomogeneous elevation and the differentiated low-frequency part of
\(\tilde\psi\).  The only missing component is the undifferentiated
low-frequency potential, which is propagated directly from its equation and
grows at most like \(\varepsilon\av{t}\).  The Riccati estimate
\cref{eq:ideas_proof_modified_energy_estimate_3d} therefore closes on an
interval of size \(c\varepsilon^{-1}\).  Existence follows by Friedrichs
regularization and compactness; the same reductions applied to the equation
for the difference, at one lower Sobolev index, give uniqueness, and a
Bona--Smith approximation completes the proof of continuous dependence.

\subsection{Previous results}
\label{sec:previous_results}

The literature on water waves, vortex sheets, and moving interfaces is vast;
see, among others, the monograph of Lannes and the survey of Wu
\cite{lannes2013water,wu2006mathematical}.  We recall here the developments
that are most directly connected with the Cauchy problem and with the
reduction of a free-boundary evolution to a quasilinear equation on the
interface.  We begin with the contour-dynamics approach to incompressible
flows, which provides the closest historical precedent for the point of view
adopted in this paper: the bulk variables are eliminated and the dynamics is
encoded directly by nonlinear singular integral operators acting on the
moving interface.  We then discuss the one-phase water-wave problem, where
the modern paradifferential strategy was developed most systematically, and
finally turn to two-fluid interfaces and to the Kelvin--Helmholtz mechanism,
which is the setting of the present work.

\subsubsection{Contour dynamics in incompressible flows}

The modern contour-dynamics approach to incompressible flows emerged from
the study of two-dimensional Euler vortex patches.  If the vorticity is
piecewise constant, the Euler transport equation preserves this structure
and the Biot--Savart law allows the motion to be described entirely through
the evolution of the boundaries separating regions of constant vorticity.
Thus a bulk transport equation is replaced by a nonlocal evolution equation
on a curve.  Deem and Zabusky discovered families of steadily translating
and rotating finite-area vortices, the so-called \(V\)-states
\cite{DZ1978}, and Zabusky, Hughes, and Roberts subsequently formulated
contour dynamics as a systematic method for the two-dimensional Euler
equations \cite{ZHR1979}.  This point of view became a powerful numerical
tool for resolving increasingly fine filamentation of vorticity interfaces;
an important development was Dritschel's contour-dynamics and contour-surgery
framework \cite{Dritschel1989}.  On the analytical side, Burbea placed the
construction of rotating vortex patches in a bifurcation framework
\cite{Burbea1982}.  The global regularity theory of vortex patches was later
established by Chemin and, by a different argument, by Bertozzi and
Constantin \cite{Chemin1993,BC1993}: an initially smooth patch boundary
remains smooth for all times.  These results provided one of the first
settings in which a genuinely nonlinear singular integral evolution of an
interface could be controlled globally through geometric cancellations.

The study of coherent vortex structures generated a second major branch of
the contour-dynamics theory.  Building on Burbea's formulation, Hmidi,
Mateu, and Verdera proved smoothness of rotating Euler patches close to the
bifurcation circle \cite{HMV2013}.  Castro, C\'ordoba, and G\'omez-Serrano
extended this program to the generalized surface quasi-geostrophic family,
constructing and proving regularity of rotating patch solutions
\cite{CCGS2016_1}, and established analyticity for families of uniformly
rotating patches in active scalar equations \cite{CCGS2016_2}.  They later
constructed compactly supported smooth vorticities rotating uniformly under
the two-dimensional Euler flow \cite{CCGS2019}.  These works illustrate a
recurring feature of contour dynamics: although the geometric unknown is
lower dimensional, the resulting equation contains singular nonlocal
operators whose precise cancellation properties determine both the
regularity of the evolution and the structure of special solutions.

A particularly influential development of this perspective, due in large
part to C\'ordoba, Gancedo, and collaborators, concerns sharp fronts for the
surface quasi-geostrophic equation and its generalized variants.  Rodrigo
derived a rigorous evolution equation for an SQG sharp front and obtained
local existence and uniqueness for smooth data by a Nash--Moser argument
\cite{Rodrigo2005}.  C\'ordoba, Fontelos, Mancho, and Rodrigo subsequently
gave numerical and self-similar evidence for finite-time singularity
formation in a family of contour equations interpolating between Euler
patches and the SQG sharp front \cite{CFMR2005}.  Gancedo established
local-in-time existence in Sobolev spaces for the \(\alpha\)-patch family,
including the SQG endpoint \cite{Gancedo2008}; C\'ordoba, C\'ordoba, and
Gancedo later completed the Sobolev theory at the SQG endpoint by proving
uniqueness \cite{CCG2018}.  Further developments include the local theory
and blow-up criteria for generalized SQG patches of Gancedo and Patel
\cite{GP2021}, and the low-regularity theory of Gancedo, Nguyen, and Patel,
which allows SQG fronts in \(H^{2+s}\), for arbitrarily small \(s>0\), and
hence permits unbounded curvature \cite{GNP2022}.  In the more singular
generalized SQG regime, C\'ordoba, G\'omez-Serrano, and Ionescu proved global
stability of the half-plane patch under small localized perturbations
\cite{CGI2019}.  Complementing these existence theories, Gancedo and Strain
ruled out splash singularities for SQG sharp fronts as long as the relevant
regularity remains controlled \cite{GS2014}.  More recently, the direct
paralinearization of the singular contour operator has become a tool for
long-time analysis: Berti, Cuccagna, Gancedo, and Scrobogna derived a
paradifferential reduction of the \(\alpha\)-SQG sharp-front equation and
used it to obtain an extended lifespan for small solutions
\cite{BCGS2025}.

The same contour-dynamics philosophy has played a central role in the
analysis of porous-medium interfaces.  C\'ordoba and Gancedo derived a
contour formulation for the three-dimensional Muskat problem and proved
local well-posedness in the stable regime, together with ill-posedness in
the unstable one \cite{CG2007}.  C\'ordoba, C\'ordoba, and Gancedo developed
the corresponding two-dimensional Sobolev theory under the Rayleigh--Taylor
sign condition \cite{CCG2011} and subsequently established a
three-dimensional local theory \cite{CCG2013}.  Constantin, C\'ordoba,
Gancedo, and Strain obtained global solutions under quantitative smallness
and slope assumptions \cite{CCGS2013}.  The contour formulation also made
it possible to study the failure of the stability condition dynamically:
Castro, C\'ordoba, Fefferman, Gancedo, and L\'opez-Fern\'andez showed that
the Rayleigh--Taylor condition may break down in finite time and, by the
same mechanism, constructed turning water waves \cite{CCFGL2012}; Castro,
C\'ordoba, Fefferman, and Gancedo later constructed solutions starting from
analytic data in the stable regime which enter the unstable regime and
subsequently lose smoothness \cite{CCFG2013}.  On the other hand, splash
formation is excluded in the regular Muskat evolution by the geometric
mechanism identified in \cite{GS2014}.  In this family of problems the sign
of the principal singular operator therefore distinguishes a dissipative
stable evolution from an unstable one, making the extraction of its
principal symbolic structure particularly transparent and important.

Contour dynamics was also used directly for inviscid free-boundary Euler
flows.  C\'ordoba, C\'ordoba, and Gancedo studied the interface between two
irrotational incompressible fluids in two dimensions and proved local
Sobolev well-posedness under the Rayleigh--Taylor condition through a
contour formulation \cite{CCG2010}.  For genuine vortex sheets, Castro,
C\'ordoba, and Gancedo analyzed the Birkhoff--Rott evolution in a
parametrization adapted to finite-energy sheets, obtaining a Cauchy theory
for analytic data and exhibiting the corresponding obstruction for
nonanalytic amplitudes \cite{CCG2012}.  We return below to the
Kelvin--Helmholtz instability and to the regularizing effect of surface
tension in the discussion of the two-phase water-wave problem.

Viewed together, the Euler patch, SQG, Muskat, water-wave, and vortex-sheet
theories show that contour dynamics is not merely a convenient
reformulation of particular fluid models.  It provides a common analytic
setting in which the evolution is driven by nonlinear singular
integral operators whose kernels depend on finite differences of the
unknown geometry.  The degree of the singularity changes from model to
model, as do the parity and adjoint properties responsible for the
available cancellations, but the central problem is repeatedly the same:
to extract from the nonlinear integral operator the quasilinear part which
governs the highest-order dynamics while controlling the remaining terms.

This is precisely the aspect of contour dynamics that is relevant to the
present paper.  Our purpose is not to develop another model-specific
expansion of a Birkhoff--Rott operator, but to isolate a general
structure-preserving paralinearization mechanism for this class of
geometric singular integrals.  In this sense, the results of
\cref{thm:abstract-structure-preserving-paralinearization-Rd,thm:abstract-locally-integrable-paralinearization-Rd}
may be viewed as an abstract counterpart of the model-dependent singular
integral analyses appearing throughout the contour-dynamics literature:
the homogeneous singular part determines the paradifferential principal
and subprincipal operators, the parity of the kernel retains the relevant
{
adjoint structure, and sufficiently regular nonhomogeneous tails yield
order-zero paraproducts modulo smoothing remainders.
}

\subsubsection{The one-phase problem}

At the endpoint $\rho^+=0$, equivalently $\Atw=2$, the upper phase is replaced by vacuum and \cref{eq:2phaseWW3D} reduces to the classical water-wave problem. In the irrotational graph setting, let $\phi$ be the velocity potential in
\[
    \Omega\pare{t}
    \defeq
    \set{\pare{x,y}\in\bR^{\mathtt d}\times\bR\,\middle|\,y<\eta\pare{t,x}},
    \qquad \mathtt d\in\set{1,2},
\]
and let
\[
    \psi\pare{t,x}
    \defeq
    \phi\pare{t,x,\eta\pare{t,x}}
\]
be its trace on the free surface.  In this historical discussion only, we use the conventional notation $G\pare{\eta}$ for the Dirichlet--Neumann operator,
\[
    G\pare{\eta}\psi
    \defeq
    \left.
    \pare{
        \partial_y\phi
        -\nabla\eta\cdot\nabla\phi
    }
    \right|_{y=\eta}.
\]
After normalizing the fluid density, the Euler equations are equivalent to the Zakharov--Craig--Sulem system
\begin{equation}
\label{eq:one-phase-zakharov-craig-sulem}
\system{
\begin{aligned}
    \partial_t\eta
    & =G\pare{\eta}\psi,
    \\
    \partial_t\psi
    & +\mathfrak{g}\eta
      -\sigma\dive\pare{
          \frac{\nabla\eta}
          {\sqrt{1+\av{\nabla\eta}^2}}
      }
      +\frac12\av{\nabla\psi}^2
      -\frac12
      \frac{
          \pare{G\pare{\eta}\psi
          +\nabla\eta\cdot\nabla\psi}^2
      }{1+\av{\nabla\eta}^2}
      =0.
\end{aligned}
}
\end{equation}
This is the canonical Hamiltonian system
\begin{equation}
\label{eq:one-phase-zakharov-hamiltonian}
    \partial_t\eta
    =\frac{\delta\mathcal H}{\delta\psi},
    \qquad
    \partial_t\psi
    =-\frac{\delta\mathcal H}{\delta\eta},
    \qquad
    \mathcal H\pare{\eta,\psi}
    \defeq
    \frac12\int_{\bR^{\mathtt d}}
    \bra{
        \psi G\pare{\eta}\psi
        +\mathfrak{g}\eta^2
        +2\sigma\pare{
            \sqrt{1+\av{\nabla\eta}^2}-1
        }
    }
    \dd x,
\end{equation}
identified by Zakharov and subsequently developed in the Dirichlet--Neumann formulation by Craig and Sulem \cite{Zakharov1968,CS1993}.  The reduction of the bulk Euler equations to \cref{eq:one-phase-zakharov-craig-sulem} isolates the main analytic difficulty in the nonlinear first-order elliptic operator $G\pare{\eta}$.

The first rigorous Cauchy theories for the exact equations concerned analytic data, small data, or special geometries.  Foundational contributions were made by Nalimov and Shinbrot \cite{Nalimov1974,Shinbrot1976}, and by Reeder and Shinbrot \cite{RS1976,RS1979}.  Further perturbative theories were developed by Kano--Nishida, Yosihara, and Craig \cite{KN1979,Yosihara1982,yosihara1983capillary,Craig1985}.  These works established that the exact free-boundary system admits a local dynamics, but they did not yet provide a large-data Sobolev theory at the natural quasilinear level.

The decisive breakthrough was achieved by Wu.  For the infinite-depth pure-gravity problem, she proved local well-posedness for Sobolev data, without a smallness assumption, first for a one-dimensional interface and then for a two-dimensional interface \cite{WuWW1997,WuWW1999}.  The argument combines a boundary-integral description of the harmonic velocity field with nonlinear changes of variables adapted to the geometry of the interface.  In the one-phase irrotational setting, the Taylor sign condition is generated by the equations themselves under the standard geometric assumptions, and this positivity supplies the coercive mechanism for the energy estimates.  Wu later proved almost-global existence for small two-dimensional water waves and global existence for small three-dimensional water waves \cite{Wu2009,Wu2011}, revealing additional nonlinear cancellations that are invisible in a direct Eulerian energy estimate. We also refer to the work of Beyer and G{\"u}nther \cite{beyer1998cauchy}.

A complementary geometric and Lagrangian theory was developed for the free-boundary Euler equations, including flows with vorticity.  Christodoulou and Lindblad obtained intrinsic a priori estimates under the Taylor sign condition, Lindblad established the corresponding local existence theory, and Coutand and Shkoller treated the free-surface problem with or without surface tension in Lagrangian coordinates \cite{CL2000,Lindblad2005,CS2007}.  Shatah and Zeng subsequently formulated geometric energies for both free surfaces and fluid interfaces and established local well-posedness in a framework that does not require irrotationality or simple connectivity \cite{SZ2008_1,SZ2008_2,SZ2011}.  These approaches remain fundamental, especially when the bulk vorticity is part of the dynamics.

The surface formulation \cref{eq:one-phase-zakharov-craig-sulem} led to a different route, based on Bony's paradifferential calculus \cite{Bony1981}.  Lannes derived tame estimates and shape-derivative formulas for the Dirichlet--Neumann operator and used an Alinhac good unknown in a Nash--Moser construction of solutions under the Taylor sign condition \cite{LannesWW2005}.  Alazard and M\'etivier then obtained a precise paralinearization of $G\pare{\eta}$ and exposed the quasilinear structure of three-dimensional water waves \cite{AM2009}.  Denoting Bony's paraproduct by $T$, the notation customary in this literature is
\begin{equation}
\label{eq:historical-water-wave-good-unknown}
\begin{aligned}
    B
    &\defeq
    \frac{
        G\pare{\eta}\psi
        +\nabla\eta\cdot\nabla\psi
    }{1+\av{\nabla\eta}^2},
    &
    V
    &\defeq
    \nabla\psi-B\nabla\eta,
    \\
    U
    &\defeq
    \psi-T_B\eta,
    &
    \lambda\pare{\eta;x,\xi}
    &\defeq
    \sqrt{
        \pare{1+\av{\nabla\eta}^2}\av{\xi}^2
        -\pare{\nabla\eta\cdot\xi}^2
    }.
\end{aligned}
\end{equation}
The principal paralinearization takes the form
\begin{equation}
\label{eq:historical-dn-paralinearization}
    G\pare{\eta}\psi
    =
    T_\lambda U
    -T_V\cdot\nabla\eta
    +R\pare{\eta,\psi},
\end{equation}
where the remainder is smoother in the tame sense relevant to the energy method.  If $a$ denotes the Taylor coefficient and $\ell$ the order-two principal symbol of the linearized curvature operator, the evolution is thereby reduced to
\begin{equation}
\label{eq:historical-water-wave-paradifferential-system}
\system{
\begin{aligned}
    \pare{\partial_t+T_V\cdot\nabla}\eta
    -T_\lambda U
    &=f_1,
    \\
    \pare{\partial_t+T_V\cdot\nabla}U
    +T_a\eta
    +\sigma T_\ell\eta
    &=f_2.
\end{aligned}
}
\end{equation}
The dispersive symbol of the symmetrized system is, at principal level,
\[
    \gamma
    \sim
    \sqrt{\lambda\pare{a+\sigma\ell}};
\]
it has order $1/2$ in the pure-gravity regime and order $3/2$ in the capillary regime.  Alazard, Burq, and Zuily used this structure to establish the gravity-capillary and pure-gravity Cauchy theories by a direct paradifferential symmetrization \cite{AlazardBurqZuily2011,ABZ2014}.  M\'etivier's general theory of paradifferential nonlinear systems supplies the corresponding abstract Cauchy framework \cite{metivier2008differential}.  In this approach, the good unknown removes the apparent loss in the shape derivative of the Dirichlet--Neumann operator, while the symbolic diagonalization turns the free-boundary Euler equations into a transport-dispersive system amenable to energy and microlocal estimates.

The reduction \cref{eq:historical-water-wave-paradifferential-system} became a common starting point for several further developments.  It is important to distinguish the paradifferential reduction itself from the additional mechanisms used after it: depending on the problem, the latter may be semiclassical Strichartz estimates, Fourier resonances, normal forms, holomorphic coordinates, or KAM reducibility.

At the local level, Alazard, Burq, and Zuily developed semiclassical Strichartz estimates, a Cauchy theory for nonlocalized data, and lower-regularity gravity theories \cite{ABZ2011_2,ABZ2016,ABZ2018}.  De Poyferr\'e and Nguyen carried the gravity-capillary reduction down to the quasilinear Sobolev threshold and combined it with Strichartz estimates to allow non-Lipschitz initial velocity traces \cite{PN2017,PN2016}.  Nguyen proved Hadamard well-posedness for the gravity flow constructed by the paradifferential method \cite{Nguyen2016}.  Using quantitative paracomposition, he also established a sharp two-dimensional gravity-capillary Cauchy theory \cite{Nguyen2017}.  Ai further improved the Sobolev thresholds for gravity waves through refined parametrices and lossless Strichartz estimates \cite{Ai2019,Ai2020}.  More recently, Ai, Ifrim, and Tataru developed balanced cubic energy estimates and used them in a low-regularity local and global theory for two-dimensional gravity waves \cite{AIT2026,AIT2022}.  These works progressively move the Cauchy theory below the threshold accessible by a purely Lipschitz-coefficient energy estimate.

For small and localized data on Euclidean domains, the quasilinear analysis---written either in paradifferential form or in equivalent boundary and holomorphic coordinates---has been combined with dispersive decay and nonlinear normal forms.  Large-time control in long-wave regimes was established by Alvarez-Samaniego and Lannes \cite{ASL2008}.  In the pure-gravity problem, Germain, Masmoudi, and Shatah proved global existence in three space dimensions, while Ionescu-Pusateri, Alazard-Delort and Wang obtained global solutions and modified scattering in two space dimensions \cite{GMS2012,IP2015,AD2015_1,AD2015_2,wang2018global}. We refer to \cite{alazard2022cauchy,wang20173,wang2019global} when the system is above a flat bottom and \cite{harrop2017finite} when the system is irrotational with finite bottom. Hunter, Ifrim, and Tataru gave a holomorphic-coordinate proof of local and almost-global existence, which was followed by the global theory of Ifrim and Tataru \cite{HIT2016,IT2016}. For long-term regularity results in three space dimensions, we refer to \cite{zheng2022long,ginsberg2025long}. Moreover, Wu proved the case that the interface is not in $C^1$ \cite{wu2019wellposedness}.  For capillary waves, Germain, Masmoudi, and Shatah proved global existence in three space dimensions, Ifrim and Tataru obtained a cubic lifespan in two space dimensions, and Ionescu and Pusateri proved global regularity and modified scattering for the full two-dimensional system \cite{GMS2015,IT2017,IP2018}, and Wang proved global regularity, scattering, and the non-existence of small traveling waves for the three-dimensional  finite depth system in \cite{Wang2020}. Deng, Ionescu, Pausader, and Pusateri established the corresponding global theory for the three-dimensional gravity-capillary system \cite{DIPP2017}. In \cite{zhu2024propagation}, Zhu proved the propagation result for gravity-capillary system of general dimension.  The proofs use different nonlinear structures, but all require an accurate identification of the quasilinear principal and subprincipal dynamics before the dispersive analysis can begin.

On compact spatial domains, where linear dispersion does not yield time decay, paradifferential conjugations instead enter Nash--Moser, KAM, and Birkhoff normal-form schemes.  Alazard and Baldi constructed gravity-capillary standing waves, Baldi, Berti, Haus, and Montalto constructed time-quasiperiodic pure-gravity waves in finite depth, and Berti and Montalto obtained quasiperiodic gravity-capillary standing waves \cite{AB2015,BBHM2018,BM2020}.  Berti and Delort proved almost-global existence for periodic gravity-capillary waves, while Berti, Feola, and Pusateri derived the quartic Birkhoff normal form and the associated long-time stability for periodic pure-gravity waves \cite{BD2018,BFP2023}.  For gravity-capillary waves with constant vorticity, Berti, Maspero, and Murgante developed a Hamiltonian paradifferential Birkhoff normal form and obtained almost-global existence for generic surface tension \cite{BMM2024}.  Paradifferential reducibility has also been combined with constant-vorticity formulations in the construction of traveling quasiperiodic waves \cite{BFM2021}; for a local Hamiltonian formulation with general vorticity and its shallow-water stability, see also \cite{CL2015}.  Thus, in the one-phase theory, paradifferential calculus has become the standard device that separates the universal quasilinear geometry from the problem-specific dispersive or Hamiltonian argument. For recent results about paradifferential reducibility and KAM theory, we refer to \cite{alazard2024paracomposition,alazard2026kam}.

\subsubsection{The two-phase problem}

When both densities are positive, the interface carries two harmonic velocity potentials, or equivalently two Dirichlet--Neumann operators, and the jump of tangential velocity becomes an independent dynamical quantity.  This changes the principal symbol in an essential way: in addition to the Rayleigh--Taylor restoring mechanism, it contains the destabilizing Kelvin--Helmholtz contribution generated by the velocity jump.  The Hamiltonian structure of the two-fluid problem was clarified by Benjamin and Bridges \cite{BB1997_1}, but the analytic behavior differs sharply from the one-phase endpoint.

In the genuine vortex-sheet regime, a nonzero tangential velocity jump and the absence of surface tension produce a high-frequency Kelvin--Helmholtz instability \cite{murgante2026steady}.  Local solutions can be constructed in analytic classes, as in the work of Sulem, Sulem, Bardos, and Frisch \cite{sulem1985finite,SSBF1981}; Caflisch and Orellana subsequently obtained a long-time analytic theory for small perturbations of a flat sheet \cite{CO1986}.  The complementary construction of singular solutions by Caflisch and Orellana proves ill-posedness in Sobolev spaces \cite{CO1989}.  Lebeau later showed that, under a priori H\"older regularity assumptions, any local two-dimensional Kelvin--Helmholtz solution must in fact have an analytic interface and analytic sheet strength \cite{Lebeau2002}; see also Wu's survey \cite{wu2006mathematical}.  This rigidity explains why one cannot expect a generic Sobolev Cauchy theory for unregularized vortex sheets.

Surface tension suppresses the short-wave instability and restores a quasilinear energy theory.  Ambrose established local well-posedness for two-dimensional periodic vortex sheets with surface tension, and Ambrose and Masmoudi treated three-dimensional vortex sheets \cite{Ambrose2003,AM2007}.  Cheng, Coutand, and Shkoller proved local well-posedness for three-dimensional vortex sheets with surface tension in the presence of bulk vorticity \cite{CCS2008}.  The geometric framework of Shatah and Zeng gives a further local theory for fluid interfaces with surface tension without assuming irrotationality or simple connectivity \cite{SZ2008_2,SZ2011}.  Pusateri obtained estimates uniform in the surface-tension coefficient and in the density of the lighter phase, and used them to justify the simultaneous one-fluid and zero-surface-tension limit \cite{Pusateri2011}.

A central paradifferential contribution in the two-phase setting is Lannes' nonlinear stability criterion \cite{Lannes2013}.  Working with the two Dirichlet--Neumann operators, Lannes derived a symbolic condition that combines the Rayleigh--Taylor coefficient, the tangential velocity jump, the density contrast, and surface tension.  The criterion shows that gravity can control the shear at low and intermediate frequencies, whereas capillarity controls the high-frequency Kelvin--Helmholtz component; it also yields estimates adapted to singular asymptotic regimes and requires a symbolic decomposition of the bottom contribution into a principal part and an infinitely smoothing tail.  This work provides the closest Dirichlet--Neumann analogue of the stability mechanism used below.  Two-fluid geometry also imposes restrictions absent from the one-phase problem: for instance, locally smooth two-dimensional interfaces with two positive densities cannot develop a splash singularity \cite{fefferman2016absence}.

The energy theories of Ambrose and of Ambrose--Masmoudi are already direct interface arguments.  What remained much less developed than in the Dirichlet--Neumann setting was a systematic symbolic procedure that converts the nonlinear Birkhoff--Rott operators themselves into paradifferential operators while retaining their adjoint and parity structure.  A recent step in this direction is the work of Murgante, Roulley, and Scrobogna on vortex sheets close to a circle, where a direct paralinearization of nonlinear singular integral operators is combined with Hamiltonian Birkhoff normal forms to obtain almost-global stability below a Weber threshold \cite{murgante2025long}.

The present paper addresses the corresponding methodological problem for three-dimensional graph interfaces.  Rather than first replacing the contour dynamics by a Dirichlet--Neumann system, we paralinearize the singular integral operators themselves.  The abstract results of \cref{thm:abstract-structure-preserving-paralinearization-Rd,thm:abstract-locally-integrable-paralinearization-Rd} isolate the homogeneous singular structure, preserve the adjoint and parity information needed by the energy method, and reduce nonhomogeneous tails to order-zero paraproducts modulo smoothing remainders.  Applied to the two-phase contour-dynamics equation, this yields the complete quasilinear symbolic system and the Sobolev local well-posedness result proved in the sequel.

\subsection{Structure of the paper}

In \cref{sec:recast}, we recast \cref{eq:2phaseWW3D} as an autonomous
Contour Dynamics Equation (CDE) in canonical surface variables.
\cref{sec:linearization-trivial-state} is devoted to the linearization of
the canonical CDE at the zero equilibrium and to the analysis of the
associated dispersion and stability conditions.  In
\cref{sec:preliminaries}, we introduce the functional setting and the
paradifferential tools used throughout the paper.

The main methodological part of the analysis is contained in
\cref{sec:paralinearization}.  We first develop an abstract
structure-preserving paralinearization theorem for critical homogeneous
graph integral operators on \(\bR^{\mathtt d}\).  We then apply this result
to the singular integral operators entering the contour dynamics and derive
the complete paradifferential reduction of the canonical CDE.  In
\cref{sec:diagonal_reduction_main}, we introduce the Alinhac good unknown,
reduce the resulting system to a diagonal complex form, establish the
modified energy estimates, and conclude the proof of
\cref{thm:main}.  Finally, \cref{sec:technical lemma} is devoted to the proof of some technical lemmas and
\cref{app:abstract-locally-integrable-paralinearization-Rd} gives the
subcritical locally integrable counterpart of the abstract
paralinearization theorem.

\subsection*{Acknowledgements}

XL would like to thank Thomas Alazard, Jean-Marc Delort and Chengyang Shao, and SS would like to thank Massimiliano Berti, Alberto Maspero and Federico Murgante, for
fruitful discussions.  XL and SS are supported by PRIN 2022 ``Turbulent
effects vs Stability in Equations from Oceanography'' (TESEO), project
number 2022HSSYPN.  ER is supported by the ERC Starting Grant 2021
``Hamiltonian Dynamics, Normal Forms and Water Waves'' (HamDyWWa), project
number 101039762. The views and opinions expressed are however those of the authors only and do not necessarily reflect those of the European Union or the European Research Council. Neither the European Union nor the granting authority can be held responsible for them. ER is also supported by the INdAM--GNAMPA project
``(In)-stability in Fluid Mechanics'', CUP E53C25002010001.

\subsection*{AI-disclosure}
The main mathematical ideas contained in this article as well as the vast majority of its write-up were generated without the assistance of AI. The authors started working on this project in September 2024 and overcame all the mathematical difficulties by the beginning of 2026, which predates the most recent AI-powered breakthroughs. The authors have used Sage to check some algebraic identities in the paradifferential reduction. The authors acknowledge OpenAI through the ChatGPT for Academic Researchers program and Anthropic through the Claude for Science program. ChatGPT and Claude were used to assist with language editing and with checking the consistency of notation. All mathematical content was written and verified by the authors, who take full responsibility for it.

\subsection*{Data availability statement} All data generated or analysed during this study are included in this article.

\subsection*{Author contribution statement} All authors contributed equally to the manuscript.

\section{Recasting \cref{eq:2phaseWW3D} as a Contour Dynamic Equation (CDE)}\label{sec:recast}

It is well-known that the irrotationality assumption in \cref{eq:2phaseWW3D} allows us to recast \cref{eq:2phaseWW3D} as a system of contour dynamics equations of the interface $ \Gamma\pare{t}.$ Our purpose here is to describe this formulation. We consider
\begin{subequations}
    \begin{align}
{\bm v}^\pm\pare{t,x}&\defeq {\bm u^\pm}\pare{t,{\bm \gamma}\pare{t,x}},\label{eq:traces_boundary}\\
q^\pm\pare{t,x}&\defeq p^\pm\pare{t,{\bm \gamma}\pare{t,x}},\\
 \psi^\pm \pare{t,x}&\defeq \phi^\pm\pare{t,{\bm \gamma}\pare{t,x}},\label{eq:psipm}
\end{align}
\end{subequations}
  respectively, the trace of the velocity field, the pressure and the potential onto the interface $\Gamma\pare{t}$.\\
  
  Next, following \cite{Lannes2013}, we introduce the quantity 
\begin{equation}\label{eq:psi}
\psi\defeq \frac{ \bbra{\rho^\pm \psi^\pm}}{\aangle{\rho^\pm}}\cdot
\end{equation}

\begin{notation}\label{notation:deltagamma}
Throughout this section, we shall use the following
notations.
\begin{enumerate}[\bf i]
    \item \label{item:deltagamma1} Let us define
\begin{align}
\delta_z f\pare{x} \defeq  f\pare{x} - f\pare{x-z}, 
&&
\Delta_z f\pare{x} \defeq  \frac{f\pare{x} - f\pare{x-z}}{\av{z}}, 
&&
x, z\in \bR^2. 
\end{align}

\item  Let $\eta \in B_{\dot{W}^{1,\infty}\pare{\bR^2;\bR}}\pare{0,1} $  and a {\it scalar} function $ \omega:\bR^2\to \bR $, let us define the following operators of nonlinear-convolution type
\begin{subequations}
    \begin{align}
\label{eq:B0}
    \mathpzc{B}_0\pare{\eta} \bra{\omega} \pare{x} \defeq & \ \frac{1}{4\pi} \pv \int_{\bR^2} \frac{z}{\pare{ \av{z}^2+\pare{ \delta_{z}\eta \pare{x} }^2 }^\frac 32} \ \omega\pare{x-z} \ \dd z, \\
    \label{eq:Q}
    \mathpzc{q}\pare{\eta} \bra{\omega} \pare{x} \defeq & \frac{1}{4\pi} \pv \int_{\bR^2} \frac{\delta_z\eta\pare{x}}{\pare{ \av{z}^2+\pare{ \delta_{z}\eta \pare{x} }^2 }^\frac 32} \ \omega\pare{x-z} \ \dd z, 
\end{align}
\end{subequations}

where $\delta_z \eta$ is defined in \cref{item:deltagamma1}. We also denote by
$\mathfrak{R}$ the vector-valued Riesz
transform,
\begin{align}\label{eq:Riesz}
    \mathfrak{R}\omega
    \defeq \vect{\mathfrak{R}_1\omega}{\mathfrak{R}_2\omega},
    &&
    \mathfrak{R}_{\mathsf j}\omega\pare{x}
    \defeq
    \frac{1}{2\pi}\pv\int_{\bR^2}
    \frac{z_{\mathsf j}}{\av{z}^3}
    \omega\pare{x-z}\dd z,
\end{align}
that is a Fourier multiplier given by
$$\widehat{\mathfrak{R}_j\omega}(\xi)=-\ii\frac{\xi_j}{\av{\xi}}\hat{\omega}(\xi).$$
In particular, we have
\begin{equation}\label{eq:B0_q_flat_linearization}
    \mathpzc B_0\pare{0}
    =
    \frac{\mathfrak R}{2},
    \qquad
    \mathpzc q\pare{0} = 0.
\end{equation}
{{
As a consequence,
\begin{equation}
\label{eq:Omega_bg_can_flat}
\mathpzc B_0\pare{0}[1]\equiv 0.
\end{equation}}}
\end{enumerate}
    
\end{notation}

The main goal of this section is to prove the following theorem, which is about reformulating the system \cref{eq:2phaseWW3D} into the CDEs for $\eta$ and $\psi$. 

{
{
\begin{theorem}[The Vortex Sheet equations as a canonical CDE with background vorticity]
\label{proprecast1}
Let $\sigma,\rho^\pm\geq 0$ with
\(\aangle{\rho^\pm}>0\), $\mathfrak g\in\bR$,
$\background\in\bR^2$, and let
$\pare{{\bm u}^\pm,p^\pm,{\bm\gamma}}$ be a classical solution of
\cref{eq:2phaseWW3D}. Let $\eta$ and $\psi$ be defined as in
\cref{eq:surface} and \cref{eq:psi}. Define the affine background potential
\begin{align}\label{eq:psibar_expression}
    \bar{\psi}\pare{x}
    \defeq
    x^\perp\cdot\background ,
    &&
    \nabla^\perp\bar{\psi}=\background
\end{align}
and the canonical perturbation
\begin{equation}\label{eq:psi_can_def}
    \tilde{\psi}
    \defeq
    \psi-\bar{\psi}.
\end{equation}
There exist $s_0>0$ and $r_0>0$ such that, if
\begin{equation}
    \norm{\eta}_{s_0,I}
    \leq r_0,
\end{equation}
then the set of equations \cref{eq:2phaseWW3D} is equivalent to the canonical contour-dynamics
system
\begin{subequations}
\label{eq:CDE}
\begin{empheq}[left=\empheqlbrace]{align}
    \eta_t
    &=
    \mathpzc g^\star\pare{\eta}
    \bra{\tilde{\Omega}\bra{\background}\pare{\eta}}
    +
    \mathpzc g\pare{\Atw\middle|\eta}
    \bra{\tilde{\psi}},
    \label{eq:CDEeta}
    \\
    \tilde{\psi}_t
    &=
    \tilde{\pi}
    \pare{\rho^\pm,\background\middle|\eta,\tilde{\psi}}.
    \label{eq:evolution_psi}
\end{empheq}
\end{subequations}
The objects appearing in \cref{eq:CDE} are defined as follows.
\begin{itemize}
    \item The vorticity variable carried by the canonical slice
    $\set{\tilde{\psi}=0}$ is
    \begin{equation}\label{eq:Omega_bg_can_def}
        \tilde{\Omega}\bra{\background}\pare{\eta}
        \defeq
        S\pare{\Atw\middle|\eta}\bar{\psi},
    \end{equation}
    where
    \begin{itemize}
        \item the potential-to-vorticity  operator $S\pare{\Atw\middle|\eta}$ is defined by
    \begin{equation}\label{eq:S}
        S\pare{\Atw\middle|\eta}\psi
        \defeq
        \nabla^\perp
        \pare{1+\Atw\tau_0\pare{\eta}}^{-1}\psi.
    \end{equation}

    \item  the operator $\tau_0\pare{\eta}$ is defined as
    \begin{equation}\label{eq:Tstar}
        \tau_0\pare{\eta}\psi
        \defeq
        \Delta^{-1}\nabla^\perp\cdot
        \pare{
            \mathpzc B_0\pare{\eta}
            \bra{\nabla\eta\cdot\nabla^\perp\psi}
            -
            \mathpzc q\pare{\eta}\nabla^\perp\psi
            +
            \pare{
                \mathpzc B_0\pare{\eta}
                \cdot\nabla\psi
            }
            \nabla^\perp\eta
        },
            \end{equation}
            where $\mathpzc{B}_0\pare{\eta}$ and $\mathpzc{q}\pare{\eta}$ are defined in \cref{eq:B0} and \cref{eq:Q}, respectively. 

    \end{itemize}
   
    \item The operator $\mathpzc g^\star\pare{\eta}$ is the normal velocity
    operator acting on vorticity variables and is defined by 
    
    \begin{equation}\label{eq:gstar}
        \mathpzc g^\star\pare{\eta}\Omega
        \defeq
        {\bm v}^\star\bra{\Omega}\pare{\eta}
        \cdot\vect{\nabla\eta}{-1},
    \end{equation}
    where  the auxiliary velocity trace
    ${\bm v}^\star\bra{\Omega}\pare{\eta}$ is
    \begin{equation}\label{eq:bmvtilde}
        {\bm v}^\star\bra{\Omega}\pare{\eta}
        \defeq
        \vect{
             V^{\star}\bra{\Omega}\pare{\eta}
        }{
             v_3^{\star}\bra{\Omega}\pare{\eta}
        }
        =
        \vect{
            \pare{
                \mathpzc B_0\pare{\eta}
                \bra{\nabla\eta\cdot\Omega}
                -
                \mathpzc q\pare{\eta}\Omega
            }^\perp
        }{
            -\mathpzc B_0\pare{\eta}^\perp\cdot\Omega
        },
    \end{equation}
    
    where $\mathpzc B_0\pare{\eta}$ and $\mathpzc q\pare{\eta}$ are defined in \cref{eq:B0} and \cref{eq:Q}, respectively.

    \item The operator $\mathpzc g\pare{\Atw\middle|\eta}$ is the generalized
    Dirichlet--Neumann operator in canonical potential variables and is defined by
    \begin{equation}\label{eq:DN}
        \mathpzc g\pare{\Atw\middle|\eta}\psi
        \defeq
        \mathpzc g^\star\pare{\eta}
        \bra{S\pare{\Atw\middle|\eta}\psi}.
    \end{equation}
    Equivalently,
    \begin{equation}\label{eq:bmv_statement}
        {\bm v}\bra{\psi}\pare{\eta}
        \defeq
        \vect{
            V\bra{\psi}\pare{\Atw\middle|\eta}
        }{
            v_3\bra{\psi}\pare{\Atw\middle|\eta}
        }
        \defeq
        {\bm v}^\star
        \bra{S\pare{\Atw\middle|\eta}\psi}
        \pare{\eta},
    \end{equation}
    and
    \[
        \mathpzc g\pare{\Atw\middle|\eta}\psi
        =
        {\bm v}\bra{\psi}\pare{\eta}\cdot\vect{\nabla\eta}{-1}.
    \]

\item The shifted Bernoulli function in canonical variables is 

    \begin{equation}\label{eq:tildepsiNEW}
        \tilde{\pi}
        \pare{\rho^\pm,\background\middle|\eta,\tilde{\psi}}
        \defeq
        \pi^\star
        \pare{
            \rho^\pm
            \middle|
            \eta,
            \tilde{\Omega}\bra{\background}\pare{\eta}
            +
            S\pare{\Atw\middle|\eta}\tilde{\psi}
        }
        -
        c_{\background},
    \end{equation}
   {where $\pi^{\star}$ is given by
    $$\pi^{\star}\pare{\rho^{\pm}\middle|\eta,\Omega}\defeq\pi_0\pare{\rho^{\pm}\middle|\eta}+\pi_2^{\star}\pare{\Atw\middle|\eta}\bra{\Omega,\Omega},$$
    with}
    \begin{itemize}
        \item  the gravity-curvature term $\pi_0$ is defined as
        \begin{equation}\label{eq:pi0}
        \pi_0\pare{\rho^\pm\middle|\eta}
        \defeq
        \frac{\sigma}{\aangle{\rho^\pm}}
        \nabla\cdot\mathpzc K\pare{\nabla\eta}
        -
        \mathfrak g\Atw\eta,
    \end{equation}
    where $\mathpzc K\pare{\nabla\eta}$ is defined in \cref{curvature2}. 

    \item the
    symmetric bilinear form $\pi_2^{\star}$ is
    \begin{equation}\label{eq:pistar2}
        \pi^\star_2\pare{\Atw\middle|\eta}
        \bra{\Omega_1,\Omega_2}
        \defeq
        -
        \pi^\star_{2;\RN{1}}\pare{\eta}
        \bra{\Omega_1,\Omega_2}
        -
        \pi^\star_{2;\RN{2}}\pare{\Atw\middle|\eta}
        \bra{\Omega_1,\Omega_2},
    \end{equation}
    with
    \begin{align}
        \pi^\star_{2;\RN{1}}\pare{\eta}
        \bra{\Omega_1,\Omega_2}
        &\defeq
        \frac{
            \Omega_1^\perp\cdot V^{\star}\bra{\Omega_2}\pare{\eta}
            +
            \Omega_2^\perp\cdot V^{\star}\bra{\Omega_1}\pare{\eta}
        }{2},
        \label{eq:pi2_star_I_def_lemma}
        \\
        \pi^\star_{2;\RN{2}}\pare{\Atw\middle|\eta}
        \bra{\Omega_1,\Omega_2}
        &\defeq
        \frac{\Atw}{4}
        \Bigg[
            2 V^{\star}\bra{\Omega_1}\pare{\eta}
            \cdot
             V^{\star}\bra{\Omega_2}\pare{\eta}
            \notag
        \\
        &\quad
        -
             v_3^{\star}\bra{\Omega_1}\pare{\eta}
            \pare{
                 v_3^{\star}\bra{\Omega_2}\pare{\eta}
                -
                2 V^{\star}\bra{\Omega_2}\pare{\eta}\cdot\nabla\eta
            }
            \notag
        \\
        &\quad
        -
             v_3^{\star}\bra{\Omega_2}\pare{\eta}
            \pare{
                 v_3^{\star}\bra{\Omega_1}\pare{\eta}
                -
                2 V^{\star}\bra{\Omega_1}\pare{\eta}\cdot\nabla\eta
            }
            \notag
        \\
        &\quad
        +
            \frac{1}{2}
            \frac{
                \Omega_1^\intercal
                \pare{\Id_{\bR^2}+\nabla\eta\otimes\nabla\eta}
                \Omega_2
            }{
                1+\av{\nabla\eta}^2
            }
        \Bigg].
        \label{eq:pi2_star_II_def_lemma}
    \end{align}

    \item the constant gauge is
    \begin{equation}\label{eq:c_background}
        c_{\background}
        \defeq\pi^\star\pare{\rho^\pm\middle|0,\background}
        =
        -\frac{\Atw}{8}\av{\background}^2 .
    \end{equation}
    In particular, subtracting $c_{\background}$ in \cref{eq:tildepsiNEW} makes
    $\pare{\eta,\tilde{\psi}}=\pare{0,0}$ a stationary trivial solution.

    \end{itemize}
    Note that $\tilde{\pi}$ decomposes according to homogeneity in $\tilde{\psi}$ as
    \begin{equation}\label{eq:tildepi_decomposition}
        \tilde{\pi}
        \pare{\rho^\pm,\background\middle|\eta,\tilde{\psi}}
        =
        \tilde{\pi}_0\pare{\rho^\pm,\background\middle|\eta}
        +
        \tilde{\pi}_1\pare{\Atw,\background\middle|\eta}
        \bra{\tilde{\psi}}
        +
        \tilde{\pi}_2\pare{\Atw\middle|\eta}
        \bra{\tilde{\psi},\tilde{\psi}},
    \end{equation}
    where
    \begin{subequations}\label{eq:tildepi_homogeneous_components}
    \begin{align}
        \tilde{\pi}_0\pare{\rho^\pm,\background\middle|\eta}
        &\defeq
        \pi_0\pare{\rho^\pm\middle|\eta}
        +
        \pi^\star_2\pare{\Atw\middle|\eta}
        \bra{
            \tilde{\Omega}\bra{\background}\pare{\eta},
            \tilde{\Omega}\bra{\background}\pare{\eta}
        }
        -
        c_{\background},
        \label{eq:pi0tilde}
        \\
        \tilde{\pi}_1\pare{\Atw,\background\middle|\eta}
        \bra{\tilde{\psi}}
        &\defeq
        2
        \pi^\star_2\pare{\Atw\middle|\eta}
        \bra{
            \tilde{\Omega}\bra{\background}\pare{\eta},
            S\pare{\Atw\middle|\eta}\tilde{\psi}
        },
        \label{eq:pi1tilde}
        \\
        \tilde{\pi}_2\pare{\Atw\middle|\eta}
        \bra{\tilde{\psi}_1,\tilde{\psi}_2}
        &\defeq
        \pi^\star_2\pare{\Atw\middle|\eta}
        \bra{
            S\pare{\Atw\middle|\eta}\tilde{\psi}_1,
            S\pare{\Atw\middle|\eta}\tilde{\psi}_2
        }.
        \label{eq:pi2tilde}
    \end{align}
    \end{subequations}
    In particular, for fixed $\eta$ and $\background$, the term
    $\tilde{\pi}_j$ is $j$-homogeneous in $\tilde{\psi}$ for
    $j=0,1,2$. 
    \end{itemize}
\end{theorem}
}}

\begin{rem}
    The formulation \cref{eq:CDE} constitutes a primary novelty of this work. Departing from verifying abstract properties of the Dirichlet-to-Neumann operator, we provide in \cref{eq:DN} an explicit representation of $\mathpzc{g}\pare{\Atw\middle|\eta}$ in terms of the singular integral operators $\mathpzc{B}_0$ and $\mathpzc{q}$. While this approach introduces significant algebraic complexity, it is essential for capturing the precise paralinear structure of the system in the presence of vorticity.
    \fine
\end{rem}

Let us outline the structure of the rest of this section.
\begin{itemize}
    \item[$\bullet$]
    In \cref{secvelo}, we describe the trace of the velocity field on the interface.

    \item[$\bullet$]
    In \cref{seceta}, we derive the CDE for the elevation function $\eta$.

    \item[$\bullet$]
    In \cref{secomega}, we derive the CDE for the vorticity variable $\Omega$.

    \item[$\bullet$]
    In \cref{sec:psi}, we derive the CDE for the weighted trace of the potential $\psi$.

    \item[$\bullet$]
    
    In \cref{sec:vorticitytopotential}, we recast the potential in term of vorticity.

    \item[$\bullet$]
    In \cref{secpsi}, we prove the main result, \cref{proprecast1}.

    \item[$\bullet$]
    In \cref{sec:background}, we record the canonical background formulation
    summarized in \cref{prop:nearCvortrevisited}.
\end{itemize}

\begin{rem}\label{rem:relation_B0_Q}
    The operators $\mathpzc{B}_0\pare{\eta}$ and $\mathpzc{q}\pare{\eta}$ defined in \cref{eq:B0} and \cref{eq:Q} are crucial since they are the only elements of \cref{eq:CDEeta,eq:evolution_psi} that resist standard paralinearization techniques (Bony paraproduct decomposition, Bony paralinearization formula, etc.). Hence a major part of the present manuscript consists of the analysis of these two operators. For this reason in the present remark we lay out some informal computations that are useful in the forthcoming analysis. Notice that using the simplification $\delta_z\eta \approx \nabla\eta \cdot z$ we obtain from \cref{eq:B0} and \cref{eq:Q} that
    \begin{equation}\label{eq:Q_B0_approximation}
    \mathpzc{q}\pare{\eta} \omega \approx \nabla\eta \cdot  \mathpzc{B}_0\pare{\eta} \omega, 
    \end{equation}
    modulo lower order terms. The approximation in \cref{eq:Q_B0_approximation} is useful in the sense that it relates the operators in \cref{eq:B0,eq:Q} via a concrete relation, hence suggesting that the paralinearization of $\mathpzc{B}_0$ entails the paralinearization for $\mathpzc{q}$, thus cutting in half the computations required for the manuscript. The relation \cref{eq:Q_B0_approximation} shall be made formal, and transparent, at the paradifferential level in \cref{step:relation_B0_Q_symbols} below.
    \fine
\end{rem}

 \subsection{Trace of the velocity field}\label{secvelo}

To prove \cref{proprecast1}, we need to start by understanding the velocity field on the interface. It is a classical fact that for vortex-sheets (see for instance \cite{birkhoff1962helmholtz,rott1956diffraction,saffman1992vortex,wu2006mathematical} for the two-dimensional case) the latter can be expressed in terms of the Birkhoff-Rott operator. From \cref{boundary1}, we have
$$\bbra{ {\bm v}^\pm } \cdot {\bm \nu}=0.$$
Consider the vortex sheet strength ${\bm \omega}$ defined on $\Gamma\pare{t}$ as
\begin{align}
\label{vorticity}
{\bm \omega}\defeq - \bbra{ {\bm v}^\pm } \times {\bm \nu}:=\vect{\Omega}{\omega_3}.
\end{align}
Together with the irrotationality conditions in \cref{eq:2phaseWW3D}, the global vorticity field $\vec{\bm \omega}\pare{t, {\bm x}}$,
{{defined in the sense of distributions as the curl of the global velocity field ${\bm u}^\pm$, is supported entirely on the interface $\Gamma\pare{t}$}} and is described by the following Dirac distribution:
$$\vec{\bm \omega}\pare{t, {\bm x}}={\bm \omega}\pare{t,x}\delta_0\pare{{\bm x}-{\bm \gamma}\pare{t,x}}.$$
The definition \cref{vorticity} combined with the orthogonality properties of the cross product imply that 
\begin{equation}\label{orthoomg}
    {\bm \omega}\cdot{\bm \nu}=0.
\end{equation}
Therefore, putting together \cref{boundary1,vorticity,orthoomg} give
\begin{align}
\omega_3 = \nabla \eta \cdot\Omega\qquad\textnormal{and}\qquad
 \frac{{\bm \nu} \times {\bm \omega}}{|{\bm \nu}|^2}\times {\bm \nu}={\bm \omega}.\label{eq:omega3}
\end{align}  
Next, for  ${\bm x}\in \Omega^\pm(t)$, we define 
\begin{equation}\label{eq:Uj}
{\bm{\mathfrak{u}}}^\pm(t,{\bm x}):=\frac 1{4\pi}\pv\int_{\mathbb R^2}\frac{ {\bm x}-{\bm \gamma}\pare{t,y}}{\av{{\bm x}-{\bm \gamma}\pare{t,y}}^3} \times {\bm \omega}\pare{t,y} \ \dd y. 
\end{equation}
% {\blue which is well defined thanks to the assumption ${\bm u^\pm}\pare{\bm x}\xrightarrow{\av{\bm x}\to \infty} 0$ in \cref{eq:2phaseWW3D}.}
% \todo[inline]{\footnotesize Here we could try to expand and quite improve the result that we are studying by assuming only ${\bm u^\pm}\pare{x, z}\xrightarrow{\av{z}\to \infty} 0$, this would give $\hat{\omega}\pare{0}\neq 0$ and we would have some KH-type of background flow that produces instability}
Let us define the {\it Birkhoff-Rott integral operator}
\begin{equation}
\label{BR}
\BR{{\bm \gamma}} \bra{\omega}\pare{t,x}:=\frac{1}{4\pi}\pv\int_{\mathbb{R}^2}\frac{{\bm \gamma}\pare{t,x}-{\bm \gamma}\pare{t,y}}{\av{{\bm \gamma}\pare{t,x}-{\bm \gamma}\pare{t,y}}^3}\times{\bm \omega}\pare{t,y} \ \dd y.
\end{equation}

Our aim is to prove the following lemma about the velocity field on the interface.

\begin{lemma}[Evaluation of the trace of the velocity field]  \label{lem:trace_velocities}
Let $\pare{{\bm u^\pm}, p^\pm,{\bm \gamma}}$ be a classical solution of \cref{eq:2phaseWW3D}. 
Let $ {\bm v}^\pm,$ ${\bm \omega} $, ${\bm{\mathfrak{u}}}^\pm$, $\BR{\cdot}$ be defined as in \cref{eq:traces_boundary,vorticity,eq:Uj,BR}. Then we have, for any $\pare{t,{\bm x} }\in I\times \Omega\pare{t}^\pm$, that 
$${\bm u^\pm}\pare{t, {\bm x}} = {\bm{\mathfrak{u}}}^\pm\pare{t, {\bm x}}.$$
Moreover, we have 
\begin{align}
\label{eq:traces_velocity}
{\bm v}^\pm=\BR{{\bm \gamma}}{\bm \omega}\pm \frac 12\frac{{\bm \nu}\times{\bm \omega}}{\av{{\bm \nu}}^2}\cdot 
\end{align}
\end{lemma}

\begin{proof}
For simplicity of notation, we assume that $t\in I$ is fixed in the rest of the argument, and we drop the $t$ dependence of the functions. In the rest of the proof, the convergence is always point-wise for all $x\in\mathbb R^2$. From direct calculations, we have
$$\dive_{{\bm x}} {\bm{\mathfrak{u}}}^\pm=0, \quad\nabla_{{\bm x}}\times {\bm{\mathfrak{u}}}^\pm=0 \quad \text{on}\quad \Omega^\pm.$$
Since we assume ${\bm \gamma}$ to be at least $\cC^1,$ we have that there exists an $\varepsilon_0>0$, such that for any $\varepsilon\in\p {0,\varepsilon_0}$ and $x\in\mathbb R^2$,
\[
{\bm \gamma}\pare{x}
\pm\varepsilon {\bm \nu}\pare{x}\in\Omega^\pm. 
\]
Next, we study the limit $\varepsilon\downarrow 0$. For any $x\in\mathbb R^2$,
we calculate the limit $\displaystyle\lim_{\varepsilon\downarrow 0}{\bm{\mathfrak{u}}}^\pm\pare{{\bm \gamma}\pare{x}\pm\varepsilon {\bm \nu}\pare{x}}$. We can split 
$${\bm{\mathfrak{u}}}^-\pare{{\bm \gamma}\pare{x}-\varepsilon {\bm \nu}\pare{x}}=I_1^{\varepsilon}\pare{x}+I_2^{\varepsilon}\pare{x},$$ 
where
\begin{align}
    I_1^{\varepsilon}\pare{x}&=\frac 1{4\pi}\text{p.v.}\int_{|y-x|>\varepsilon ^\delta|{\bm \nu}\pare{x}|}\frac{{\bm \gamma}\pare{x}-\varepsilon {\bm \nu}\pare{x}-{\bm \gamma}\pare{y}}{\av{{\bm \gamma}\pare{x}-\varepsilon {\bm \nu}\pare{x}-{\bm \gamma}\pare{y}}^3}\times \omega(y)\dd y,\\
    I_2^{\varepsilon}\pare{x}&=\frac 1{4\pi}\text{p.v.}\int_{|y-x|\le\varepsilon^\delta|{\bm \nu}\pare{x}|}\frac{{\bm \gamma}\pare{x}-\varepsilon {\bm \nu}\pare{x}-{\bm \gamma}\pare{y}}{\av{{\bm \gamma}\pare{x}-\varepsilon {\bm \nu}\pare{x}-{\bm \gamma}\pare{y}}^3}\times \omega(y) \dd y.
\end{align}
Here, we have introduced a parameter $\delta\in (\frac 23,1)$. \begin{proofpart}[Estimate of $I_1^{\varepsilon}$]
 We first show that $I_1^{\varepsilon}\to \BR{{\bm \gamma}}\omega$ point-wise as $\varepsilon \downarrow 0$. Observe from \cref{eq:surface}, that
\begin{equation}
\label{eq:control1}
\av{{\bm \gamma}\pare{y}-{\bm \gamma}\pare{x}}\geq\av{y-x}.
\end{equation}
Hence, for $|y-x|>\varepsilon^\delta|{\bm \nu} \pare{x}|$, we have
\begin{equation}
\label{eq:control2}
\frac{\varepsilon |{\bm \nu} \pare{x}|}{\av{{\bm \gamma}\pare{y}-{\bm \gamma}\pare{x}}}\leq\varepsilon^{1-\delta}\to 0\quad\quad\text{as}\quad\quad \varepsilon\downarrow 0.
\end{equation}
Define the function $\mathpzc f_i\pare u:=\frac{u_i}{\av u^3}$. Then, by Taylor integral formula, we have
\[
\mathpzc f_i\pare{{\bm \gamma}\p x-\va{\bm \nu}\p x-{\bm \gamma}\p y}-\mathpzc f_i\pare{{\bm \gamma}\p x-{\bm \gamma}\p y}=-\va\int_0^1 {\bm \nu}\p x\cdot\nabla \mathpzc f_i\pare{{\bm \gamma}\p x-t\va{\bm \nu}\p x-{\bm \gamma}\p y}\dd t.
\]
Thus, denoting $\mathpzc f:=\pare{\mathpzc f_1,\mathpzc f_2},$ we obtain that
\[
I_1^\va\pare{x}\defeq I^\va_{1,1}\pare{x}+I^\va_{1,2}\pare{x},
\]
where
\[
I^\va_{1,1}\pare{x}\defeq \frac 1{4\pi}\text{p.v.}\int_{|y-x|>\varepsilon ^\delta|{\bm \nu}\pare{x}|}\frac{{\bm \gamma}\pare{x}-{\bm \gamma}\pare{y}}{\av{{\bm \gamma}\pare{x}-{\bm \gamma}\pare{y}}^3}\times \omega(y)\dd y,
\]
and
\[
I^\va_{1,2}\pare{x}\defeq -\frac \va{4\pi}\text{p.v.}\int_{|y-x|>\varepsilon ^\delta|{\bm \nu}\pare{x}|}
\int_0^1 \dd\mathpzc f\pare{{\bm \gamma}\p x-t\va{\bm \nu}\p x-{\bm \gamma}\p y}\bra{{\bm \nu}\p x}\dd t\times{\bm \omega}\p y
\dd y.
\]
By definition of the principal value, we have $I_{1,1}^{\varepsilon}\pare{x}\to \BR{{\bm \gamma}}\bra{\omega}\pare{x}$ as $\varepsilon \downarrow 0.$ For the term $I^\va_{1,2}$, first we have for $i,j=1,2,3$,  $\partial_j\mathpzc f_i\p u=O\pare{\av u^{-3}}$. So together with \cref{eq:control1}, \cref{eq:control2} and $t\in [0,1]$, we obtain that
\[
\norm{\dd\mathpzc f\pare{{\bm \gamma}\p x-t\va{\bm \nu}\p x-{\bm \gamma}\p y}}_{\mathcal{L}(\mathbb{R}^3)}\lesssim \av{{\bm \gamma}\p x-t\va{\bm \nu}\p x-{\bm \gamma}\p y}^{-3}\lesssim \av{{\bm \gamma}\p x-{\bm \gamma}\p y}^{-3} \lesssim \av{y-x}^{-3}.
\]
Notice that
\begin{equation}
\label{eq:control3}
\int_{\va^\delta}^\infty\frac{\va}{r^2}dr=\va^{1-\delta}\to 0,\quad\text{as}\quad\va\downarrow 0.
\end{equation}
Therefore, passing to the polar coordinates $ y-x=r\vect{\cos\theta}{\sin\theta}$, we obtain that
$I_{1,2}^{\varepsilon}\pare{x}\to 0$ as $\va\downarrow 0$, which means that $I_{1}^{\varepsilon}\pare{x}\to \BR{{\bm \gamma}}\bra{\omega}\pare{x}$ as $\va\downarrow 0$. \end{proofpart}\begin{proofpart}[Estimate of $I_2^{\varepsilon}$] We  write 
\[
I_2^\va\p x \defeq I_{2,1}^\va\p x+I_{2,2}^\va\p x,
\]
where
\[
I_{2,1}^\va\p x\defeq \frac 1{4\pi}\text{p.v.}\int_{\mathbb{R}^2}\mathbf{1}_{|y-x|\le \varepsilon ^\delta|{\bm \nu}\pare{x}|}\frac{{\bm \gamma}\pare{x}-{\bm \gamma}\pare{y}}{\av{{\bm \gamma}\pare{x}-\va{\bm \nu}\p x-{\bm \gamma}\pare{y}}^3}\times \omega(y)\dd y,
\]
and
\[
I_{2,2}^\va\p x\defeq -\frac \va{4\pi}{\bm \nu}\p x\times \text{p.v.}\int_{\mathbb{R}^2}\omega\p y\mathbf{1}_{|y-x|\le \varepsilon ^\delta|{\bm \nu}\pare{x}|}\frac{1}{\av{{\bm \gamma}\pare{x}-\va{\bm \nu}\p x-{\bm \gamma}\pare{y}}^3}\dd y.
\]
We first estimate $I_{2,1}^\va\p x$. From \cref{eq:surface} and \cref{boundary1}, we can write
\begin{equation}
\begin{aligned}
\av{{\bm \gamma}\p x-\va {\bm \nu}\p x-{\bm \gamma} 
\p y}^2&=|y-x-\varepsilon\nabla\eta\p x|^2+\pare{\eta\p y-\eta\p x+\va}^2\\
&=|y-x|^2+\va^2\av{{\bm \nu}\p x}^2
+\av{\eta\p y-\eta\p x}^2+2\varepsilon\pare{\eta\p y-\eta\p x-\nabla\eta \p x\cdot \p{y-x}}.
\end{aligned}
\end{equation}
By Taylor expansion, we have
\[\eta\p {y}=\eta\p {x}+\nabla\eta\p x\cdot\p{y-x}+O\pare{\av{y-x}^2}\quad \Longrightarrow\quad
|\eta\p y-\eta\p x|^2-\pare{\nabla\eta\p x\cdot\p{y-x}}^2=O\pare{|y-x|^3}.\]
Moreover, 
\begin{equation}
\begin{aligned}
    &\av{y-x}^2+\pare{\nabla\eta\p x\cdot\p{y-x}}^2=\sum_{i=1}^2\pare{1+\av{\partial_i\eta\p x}^2}\pare{y_i-x_i}^2+2\partial_1\eta\p x\partial_2\eta\p x\pare{y_1-x_1}\pare{y_2-x_2}.
\end{aligned}
\end{equation}
For the symmetric matrix 
$$
\begin{pmatrix}
1+\av{\partial_1\eta\p x}^2 & \partial_1\eta\p x\partial_2\eta\p x\\
\partial_1\eta\p x\partial_2\eta\p x & 1+\av{\partial_2\eta\p x}^2
\end{pmatrix},$$ 
its two eigenvalues are $\av{{\bm \nu}\p x}^2$ and $1$. Under the orthogonal transformation, we only need to consider $\av{{\bm \nu}\p x}^2\pare{y_1-x_1}^2+\pare{y_2-x_2}^2$. Let $\tilde{y}:=\pare{\tilde{y}_1,\tilde{y}_2},$ with \[\tilde y_1:=\av{{\bm \nu}\p x}^\frac 12\pare{y_1-x_1},\quad\tilde y_2:=\av{{\bm \nu}\p x}^{-\frac 12}\pare{y_2-x_2}.\]
Then, we have 
\[
\av{{\bm \nu}\p x}^2\pare{y_1-x_1}^2+\pare{y_2-x_2}^2=\av{{\bm \nu}\p x}\av{\tilde y}^2.
\]
Recall that $\eta\in H^{s_0}$, we have $\norm{\nabla\eta}_{L^\infty}<\infty$. Denote $\mathpzc c\defeq \pare{1+\norm{\nabla\eta}_{L^\infty}^2}^\frac 14$, then $\av{{\bm \nu}\p x}^\frac 12\le \mathpzc c$ for any $x\in\mathbb R^2$. Moreover,
recalling that $\av{y-x}<\va^\delta\av{{\bm \nu}\p x}$ and $\delta>\frac 23$, we obtain that 
\[
\av{{\bm \gamma}\p x-\va {\bm \nu}\p x-{\bm \gamma} 
\p y}^2\approx \av{{\bm \nu}\p x}\av{\tilde y}^2+\av{{\bm \nu} \p x}^2\va^2.
\]
 The integration region of $\tilde y$ is an ellipse, located between two circles with radii $\frac 1{c}\va^\delta\av{{\bm \nu}\p x}$ and $c\va^\delta \av{{\bm \nu}\p x}$.
Next, we denote the polar coordinates $\tilde y=r\vect{\cos\theta} {\sin\theta}$. For any constant $C>0$ that does not depend on $\va$,
from symmetry, we have
\[
\int_0^{C\va^\delta}\int_0^{2\pi}\frac{r^2\pare{\cos\theta,\sin\theta}}{\pare{\av{{\bm \nu}\p x}r^2+\av{{\bm \nu}\p x}^2\va^2}^\frac 32}\dd r\dd\theta=0.
\]
Moreover, as $\va\downarrow 0$, we have
\[
\int_0^{C\va^\delta}\frac{r^3}{\pare{\av{{\bm \nu}\p x}r^2+\av{{\bm \nu}\p x}^2\va^2}^\frac 32}\dd r<\frac 1{\av{{\bm \nu}\p x}^{3/2}}\int_0^{C\va^{\delta}} 1\dd r\le C\va^\delta\to 0.
\]
As a consequence, $I_{2,1}^\va\pare{x}\to 0$ as $\va\downarrow 0$. To compute the limit of  $I_{2,2}^\va\pare{x}$, we only need to study the limit of
\[
-\frac{{\bm \nu}\p x}{4\pi}\times \pv\int_0^{2\pi}\int_{0}^{C\varepsilon^\delta}\frac{ \varepsilon r\omega\pare{x+|{\bm \nu}(x)|r\pare{\cos\theta, \sin\theta}}}{\pare{\av{{\bm \nu}\p x}r^2+\av{{\bm \nu}\p x}^2\va^2}^\frac 32}\dd r\dd\theta
\]
as $\va\downarrow 0$. Notice that
\begin{align}
\int_{0}^{C\varepsilon^\delta}\frac{ \varepsilon r}{\pare{\av{{\bm \nu}\p x}r^2+\av{{\bm \nu}\p x}^2\va^2}^\frac 32}\dd r&\xlongequal{r=\varepsilon{\av{\bm\nu}}^\frac 12 r'}\frac 1{\av{{\bm \nu}\p x}^2}\int_0^{\frac{C\va^{\delta-1}}{\av{{\bm \nu}\p x}^{1/2}}}\frac{r'}{\pare{r'^2+1}^\frac 32}\dd r'\\&\xrightarrow{\varepsilon\downarrow 0} \frac 1{\av{{\bm \nu}\p x}^2}\int_0^{\infty}\frac{r'}{\pare{r'^2+1}^\frac 32}\dd r'=\frac 1{\av{{\bm \nu}\p x}^2}\cdot
\end{align}
So we obtain that
\[
\lim_{\va\to 0^+}I_{2}^\va\p x=\lim_{\va\to 0^+}I_{2,2}^\va\p x=-\frac 12\frac{{\bm \nu}\pare{x}\times \omega\pare{x}}{|{\bm \nu}\pare{x}|^2}.
\]
\end{proofpart}
\begin{proofpart}[Synthesis]
Combining the foregoing calculations, we conclude that
\begin{equation}
\label{eq:trace-1}
\lim_{\varepsilon\downarrow 0}{\bm{\mathfrak{u}}}^-\pare{{\bm \gamma}\pare{x}-\varepsilon {\bm \nu}\pare{x}}=\BR{{\bm \gamma}}\bra{\omega}\pare{x}-\frac 12\frac{{\bm \nu}\pare{x}\times \omega\pare{x}}{|{\bm \nu}\pare{x}|^2}\cdot
\end{equation}
Similarly, we obtain
\begin{equation}
\label{eq:trace-2}
\lim_{\varepsilon\downarrow 0}{\bm{\mathfrak{u}}}^+\pare{{\bm \gamma}\pare{x}+\varepsilon {\bm \nu}\pare{x}}=\BR{{\bm \gamma}}\bra{\omega}\pare{x}+\frac 12\frac{{\bm \nu}\pare{x}\times\omega\pare{x}}{|{\bm \nu}\pare{x}|^2}\cdot
\end{equation}
 Let \({\bm w}^\pm\defeq{\bm{\mathfrak u}}^\pm-{\bm u}^\pm\) in \(\Omega^\pm\). By the foregoing trace formulas \cref{eq:trace-1,eq:trace-2},
\[
    \bbra{\left.{\bm{\mathfrak u}}^\pm\right|_{\Gamma}}
    =
    -\frac{{\bm\nu}\times{\bm\omega}}{\av{{\bm\nu}}^2}.
\]
On the other hand, recalling that \(\bbra{{\bm v}^\pm}\cdot{\bm\nu}=0\), from \cref{vorticity} we obtain
\[
    \bbra{{\bm v}^\pm}
    =
    -\frac{{\bm\nu}\times{\bm\omega}}{\av{{\bm\nu}}^2}.
\]
Therefore \(\bbra{\left.{\bm w}^\pm\right|_{\Gamma}}=0\), and we may define the continuous vector field
\[
    {\bm w}
    \defeq
    {\bm w}^-\mathbf 1_{\Omega^-}
    +
    {\bm w}^+\mathbf 1_{\Omega^+}
    \qquad\textnormal{on}\quad\bR^3.
\]
Since both \({\bm{\mathfrak u}}^\pm\) and \({\bm u}^\pm\) are divergence-free and irrotational in \(\Omega^\pm\), and since \({\bm w}\) has no jump across \(\Gamma\), we have \(\dive_{\bm x}{\bm w}=0\) and \(\nabla_{\bm x}\times{\bm w}=0\) in the sense of distributions on \(\bR^3\). Consequently,
\[
    \Delta_{\bm x}{\bm w}
    =
    \nabla_{\bm x}\pare{\dive_{\bm x}{\bm w}}
    -
    \nabla_{\bm x}\times
    \pare{\nabla_{\bm x}\times{\bm w}}
    =0
    \qquad\textnormal{in}\quad\mathcal D'\pare{\bR^3}.
\]
Hence \({\bm w}\) is harmonic on \(\bR^3\).  By the standing
far-field normalization, \({\bm w}\) is bounded.  Liouville's theorem
therefore implies that \({\bm w}\) is constant.  The same normalization
fixes the common constant velocity to be zero, and hence
\({\bm w}\equiv0\).  Thus
\({\bm u}^\pm={\bm{\mathfrak u}}^\pm\) in \(\Omega^\pm\), and the trace
formulas computed above yield \cref{eq:traces_velocity}.
\end{proofpart}
\end{proof}

\subsection{\texorpdfstring{The CDE for the elevation $ \eta $}{The CDE for the elevation η}}\label{seceta}

Our next goal is to find the CDE for the elevation function $\eta.$ It is obtained by particularizing the transport equation \cref{boundary1} in our context via the explicit computation of the Birkhoff-Rott operator \cref{BR}. More precisely, the result reads as follows.

\begin{lemma}[Derivation of the CDE for the elevation function] \label{lem:evolution_eta}
    Let the interface $\gamma$ be a graph onto $\bR^2$ as in \cref{eq:surface}, there exist $s_0 > 0$ and $r_0 > 0$ such that if $\eta \in \BallR{s_0}{r_0}$ then we have that
    \begin{equation}\label{eq:evolution_eta_omegah}
    \eta_t = \mathpzc{g}^\star \pare{\eta}\Omega,
\end{equation}
where
\begin{itemize}
    \item $\Omega$ is defined in \cref{vorticity};
    \item the operator $\mathpzc{g}^\star \pare{\eta}$ is defined as
    \begin{equation}
        \mathpzc{g}^\star\pare{\eta}\Omega \defeq {\bm v}^\star\bra{\Omega}\pare{\eta}\cdot \vect{\nabla\eta}{-1},
    \end{equation}
    and ${\bm v}^\star\bra{\Omega}\pare{\eta}$ is defined in \cref{eq:bmvtilde}.
\end{itemize}
\end{lemma}

\begin{proof}
First observe that, by virtue of \cref{eq:surface} and \cref{boundary1}, one has    
\begin{equation}\label{gm.nu}
    {\bm \gamma}_t\cdot{\bm \nu} = \eta_t.
\end{equation}
 Besides, the relation \cref{eq:traces_velocity} gives
\begin{equation}\label{vpm.nu}
    {\bm v}^\pm\cdot{\bm \nu}=\BR{{\bm \gamma}}{\bm \omega}\cdot{\bm \nu}.
\end{equation}
But, the integral operator \cref{BR} can be simplified using  \cref{eq:surface,eq:B0,eq:Q,eq:omega3} 
as
\begin{equation}\label{eq:BR_on_graph}
    \begin{aligned}
\left.\BR{{\bm \gamma}}{\bm \omega}\p x\right|_{{\bm \gamma} \p x=\pare{x, \eta\pare{x}}}
= 
& \
\frac{1}{4\pi} \pv \int _{\bR^2}
\begin{bmatrix}
\delta_z\eta\pare{x} \Omega^\perp\pare{x-z} - z^\perp \omega_3\pare{x-z}
\\
z^\perp \cdot\Omega\pare{x-z}
\end{bmatrix}
\frac{\dd z}{\pare{ \av{z}^2+\pare{ \delta_{z}\eta \pare{x} }^2 }^\frac 32}
\\
= & 
\vect{-\mathpzc{B}_0\pare{\eta}^\perp \bra{\nabla\eta\cdot\Omega}+ \mathpzc{q}\pare{\eta}\Omega^\perp}{\mathpzc{B}_0\pare{\eta}^\perp \cdot\Omega}\pare{x}. 
\end{aligned}
\end{equation}
Inserting \cref{gm.nu}, \cref{vpm.nu} and \cref{eq:traces_boundary} into \cref{boundary1} gives the desired result.
\end{proof}

\subsection{\texorpdfstring{The CDE for the vorticity ${\bm \omega}$}{The CDE for the vorticity ω}}\label{secomega}

It is clear from \cref{lem:trace_velocities}, \cref{eq:Uj} and \cref{eq:surface} that we can describe completely the velocity flow in terms of $ \eta $ and $ {\bm \omega}$. Moreover, thanks to \cref{eq:omega3}, we obtain that the evolution of the system \cref{eq:2phaseWW3D} is completely determined by the evolution of $\eta$ and $\Omega$. We now provide the evolution equation of the quantity $ \Omega :=\pare{\omega_1, \omega_2} $ in terms of boundary contributions only.

\begin{notation}
    Let us now define the following auxiliary variables in order to simplify the computations in the present section:
\begin{align}\label{eq:BOmega}
{\bm v}_N\defeq \BR{{\bm \gamma}}{\bm \omega}, 
&&
{\bm v}_T \defeq \frac{{\bm \nu}\times{\bm \omega}}{\av{{\bm \nu}}^2}, 
\end{align}
so that, by virtue of \cref{eq:traces_velocity},
$${\bm v}^\pm = {\bm v}_N\pm \frac{{\bm v}_T}{2}\cdot$$ 
\end{notation}

\begin{notation}
    Throughout the manuscript, we denote the Jacobian matrix and the rotated Jacobian matrix of a vector field $V:\mathbb{R}^2\to\mathbb{R}^p$, with $p\in\mathbb{N}^*$, as 
    $$J V \defeq \pare{\nabla V^\intercal}^\intercal\qquad\textnormal{and}\qquad J^\perp V \defeq \pare{\nabla^\perp V^\intercal}^\intercal.$$ Notice that according to this notation, given a scalar function $v$ then $\nabla v$ is a column vector, while $J v = \pare{\nabla v}^\intercal$ is a row vector. We also denote the Hessian matrix $${\bm H} V \defeq \nabla\nabla^\intercal V.$$
        An analogous definition holds for the operators $J_{{\bm x}}$ and ${\bm H}_{{\bm x}}$.
\end{notation}

With these notations in hand, we can state our result.

\begin{lemma}[Derivation of the CDE for the vorticity]\label{lem:evolution_omega}
Let $\sigma\geq 0$, $\mathfrak{g}\in \bR$. Let $x\mapsto {\bm \gamma}\pare{x}$ be a $\cC^2$ surface and let ${\bm \omega}$ be defined as in \cref{vorticity}.
\begin{enumerate}
    \item \textbf{General Surface Case.} In the case in which ${\bm \gamma}$ is a generic surface, the vorticity ${\bm \omega}$ evolves accordingly to the evolutionary equation
    \begin{equation}
    \label{eq:simplification_eq_omega_3}
    \partial_t \pare{ \pare{1+\Atw\tilde{{\bm{m}}}\pare{{\bm \gamma}}}\bra{{\bm \omega}}} =
    -{\bm \pi}^\dagger \pare{\rho^\pm\middle| {\bm \gamma}, {\bm \omega}} +
    \tilde{{\bm{n}}} \pare{{\bm \gamma}}\bra{{\bm \omega}},
    \end{equation}
    where
    \begin{itemize}
        \item $\Atw $ is the Atwood number defined in \cref{eq:Atwood};
        \item the nonlinear terms in \cref{eq:simplification_eq_omega_3} are explicitly defined as
        \begin{subequations}
        \label{eq:TNP}
        \begin{align}
            \tilde{{\bm{m}}}\pare{{\bm \gamma}}\bra{{\bm \omega}} \defeq & \  -J^\perp{\bm \gamma}  \pare{J{\bm \gamma}}^\intercal {\bm v}_N, \label{eq:m_tilde_def}
            \\
             \tilde{{\bm{n}}} \pare{{\bm \gamma}}\bra{{\bm \omega}} \defeq & \ J^\perp {\bm \gamma}_t \pare{J{\bm \gamma}}^\intercal \pare{{\bm v}_T - \Atw {\bm v}_N}, \label{eq:n_tilde_def}
            \\
            {\bm \pi}^\dagger \pare{\rho^\pm\middle| {\bm \gamma}, {\bm \omega}}
            \defeq & \ {\bm \pi}^\dagger _0 \pare{\rho^\pm\middle| {\bm \gamma}} + {\bm \pi}^\dagger_{2} \pare{\Atw\middle| {\bm \gamma}} \bra{{\bm \omega},{\bm \omega}}, \label{eq:pi_dagger_def}
        \end{align}
        \end{subequations}
        with 
        \begin{align}
            {\bm \pi}^\dagger _0 \pare{\rho^\pm\middle| {\bm \gamma}} \defeq & \ J{\bm \gamma}\nabla^{\perp}\pare{\frac{\sigma}{\aangle{\rho^\pm}}\mathpzc{k}\pare{{\bm \gamma}} + \mathfrak{g}\Atw {\bm \gamma}_3}, \label{eq:pi0_dagger}
            \\
            {\bm \pi}^\dagger _2 \pare{\Atw\middle| {\bm \gamma}} \bra{{\bm \omega}, {\bm \omega}} \defeq &
            \underbrace{-\pare{J^\perp {\bm \gamma}\nabla \pare{{\bm v}_T\cdot {\bm \gamma}_t} + J{\bm \gamma} \nabla^\perp \pare{{\bm v}_N\cdot {\bm v}_T}} }_{\eqdef {\bm \pi}^\dagger _{2,\RN{1}} \pare{ {\bm \gamma}} \bra{{\bm \omega}, {\bm \omega}}} \nonumber
            \\
            & +
            \underbrace{ \Atw \pare{\frac{1}{2}J{\bm \gamma} \nabla^\perp \pare{\av{{\bm v}_N}^2 + \frac{\av{{\bm v}_T}^2}{4}} + J^\perp {\bm \gamma}\nabla \pare{{\bm v}_N\cdot {\bm \gamma}_t} }}_{\eqdef {\bm \pi}^\dagger _{2,\RN{2}} \pare{\Atw\middle| {\bm \gamma}} \bra{{\bm \omega}, {\bm \omega}}}, \label{eq:pi2_dagger}
        \end{align}
        where ${\bm v}_N$ and ${\bm v}_T$ are defined in \cref{eq:BOmega} and $\mathpzc{k}\pare{{\bm \gamma}}$ is the mean curvature of the surface $x\mapsto {\bm \gamma}\pare{x}$. 
    \end{itemize}

    \item \textbf{Graph Case.} In the case in which ${\bm \gamma}$ is a graph onto $\bR^2$ as per \cref{eq:surface}, there exist $s_0 >0$ and $r_0 > 0$ such that if $\eta \in \BallR{s_0}{r_0}$ the equation \cref{eq:simplification_eq_omega_3} becomes
    \begin{equation}
    \label{eq:evolution_omegah}
    \partial_t \pare{ \pare{ 1+\Atw T \pare{\eta} } \Omega  } = \nabla^\perp  \pi^\star\pare{\rho^\pm\middle| \eta , \Omega}, 
    \end{equation}
    where $\Omega$ corresponds to the first two components of $\bm{\omega}$ (cf. \cref{vorticity}) and:
    \begin{itemize}
        \item the operator $T\pare{\eta}$ is defined by 
        \begin{equation}\label{eq:T_eta_def}
            T\pare{\eta}\Omega \defeq \mathpzc{B}_0\pare{\eta} \bra{\nabla\eta\cdot\Omega} - \mathpzc{q}\pare{\eta}\Omega + \pare{\mathpzc{B}_0 \pare{\eta}^\perp \cdot \Omega} \nabla^\perp \eta,
        \end{equation}
         with $\mathpzc{B}_0$ and $\mathpzc{q}$ as in \cref{eq:B0} and \cref{eq:Q};
        
        \item the nonlinear term $\pi^\star\pare{\rho^\pm\middle|\eta, \Omega}$ is defined as
              \begin{subequations}
        \begin{equation}\label{eq:pistar}
            \pi^\star\pare{\rho^\pm\middle|\eta,\Omega} \defeq \pi_0\pare{\rho^\pm\middle|\eta} + \pi_2^\star\pare{\Atw\middle|\eta}\bra{\Omega,\Omega},
        \end{equation}
        where $\pi_0$ is defined in \cref{eq:pi0},
       and the bilinear form $\pi_2^\star\pare{\Atw\middle|\eta}\bra{\Omega_1,\Omega_2}$ is defined as 
        \begin{equation}\label{eq:pistar2_boundary_lemma}
    \pi^\star_{2}\pare{\Atw\middle| \eta }\bra{\Omega_1, \Omega_2} \defeq  - \pi^\star_{2;\RN{1}}\pare{\eta }\bra{\Omega_1, \Omega_2} -  \pi^\star_{2;\RN{2}}\pare{\Atw\middle| \eta }\bra{\Omega_1, \Omega_2} ,
\end{equation}
where
    \begin{align}
        \pi^\star_{2;\RN{1}}\pare{\eta }\bra{\Omega_1, \Omega_2} \defeq & \ 
        \frac{ \Omega_1^\perp \cdot  V^{\star}\bra{\Omega_2}\pare{\eta} + \Omega_2^\perp \cdot  V^{\star}\bra{\Omega_1}\pare{\eta}}{2} , \label{eq:pi2_star_I_boundary_lemma}
        \\
        %--------------------------------------------
        \pi^\star_{2;\RN{2}}\pare{\Atw\middle| \eta }\bra{\Omega_1, \Omega_2} \defeq & \ 
        \frac{\Atw}{4} \Big[
         2  V^{\star}\bra{\Omega_1}\pare{\eta} \cdot  V^{\star}\bra{\Omega_2}\pare{\eta} 
         \nonumber
        \\
        & \quad
        -  v_3^{\star}\bra{\Omega_1}\pare{\eta} \pare{ v_3^{\star}\bra{\Omega_2}\pare{\eta} - 2 V^{\star}\bra{\Omega_2}\pare{\eta}\cdot\nabla\eta}
        \nonumber
        \\
        & \quad
        -  v_3^{\star}\bra{\Omega_2}\pare{\eta} \pare{ v_3^{\star}\bra{\Omega_1}\pare{\eta} - 2 V^{\star}\bra{\Omega_1}\pare{\eta}\cdot\nabla\eta}
        \nonumber
        \\
        & \quad
        +
        \frac{1}{2}\frac{ \Omega_1^\intercal \pare{\Id_{\bR^2} + \nabla\eta\otimes\nabla\eta}\Omega_2}
        {\p{1+|\nabla\eta|^2}}
        \Big].
        \label{eq:pi2_star_II_boundary_lemma}
    \end{align}
\end{subequations}

        with ${\bm v}^\star\bra{\Omega}\pare{\eta} = \pare{ V^{\star},  v_3^{\star}}$ as in \cref{eq:bmvtilde}.
    \end{itemize}
\end{enumerate}
\end{lemma}

\begin{proof}

We divide the proof into two steps.
 
\begin{step}[The evolution equation for a generic parametrization ${\bm \gamma}$ of the interface]

First observe that for any vector field ${\bm f}\defeq\vectDT{f_1}{f_2}{f_3}^\intercal \in\mathbb R^3$, from \cref{boundary1}, we have
\begin{equation}
\label{operator2}
{\bm f}\times{\bm \nu}=\pare{{\bm f}\cdot\partial_2{\bm \gamma}}\partial_1{\bm \gamma}- \pare{{\bm f}\cdot\partial_1{\bm \gamma}}\partial_2{\bm \gamma}  =-J {\bm \gamma}\pare{\pare{J ^\perp{\bm \gamma}}^\intercal {\bm f}}.
\end{equation}
Therefore, in view of \cref{vorticity} and \cref{operator2}, we have
\begin{equation}
\label{vorticity3}
{\bm \omega}=-\pare{\bbra{{\bm v}^\pm}\cdot\partial_2{\bm \gamma}} \partial_1{\bm \gamma}+\pare{\bbra{{\bm v}^\pm}\cdot\partial_1{\bm \gamma}} \partial_2{\bm \gamma}.
\end{equation}
Consequently, taking the derivative in time, we obtain
\begin{align}
    \partial_t{\bm \omega}&=-\partial_{t}\pare{\bbra{{\bm v}^\pm}\cdot\partial_2{\bm \gamma}} \partial_1{\bm \gamma}-\pare{\bbra{{\bm v}^\pm}\cdot\partial_2{\bm \gamma}} \partial_1\partial_{t}{\bm \gamma}\\
    &\quad+\partial_t\pare{\bbra{{\bm v}^\pm}\cdot\partial_1{\bm \gamma}} \partial_2{\bm \gamma}+\pare{\bbra{{\bm v}^\pm}\cdot\partial_1{\bm \gamma}} \partial_2\partial_t{\bm \gamma}.
\end{align}
Using \cref{operator2}, we have
\begin{equation}
\label{eq:expansion0}
\partial_t\pare{{\bm v}^\pm\times{\bm \nu}}=\partial_t\pare{{\bm v}^\pm\cdot\partial_2{\bm \gamma}}\partial_1{\bm \gamma}-\partial_t\pare{{\bm v}^\pm\cdot\partial_1{\bm \gamma}}\partial_2{\bm \gamma}+\pare{{\bm v}^\pm\cdot\partial_2{\bm \gamma}}\partial_t\partial_1{\bm \gamma}-\pare{{\bm v}^\pm\cdot\partial_1{\bm \gamma}}\partial_t\partial_2{\bm \gamma}.
\end{equation}
From the definitions of ${\bm v}^\pm$ in \cref{eq:traces_boundary}, we have for any $\mathsf{j}=1,2 $,
\begin{equation}
\label{eq:expansion1}
\begin{aligned}
\partial_t\pare{{\bm v}^\pm\cdot\partial_{\mathsf{j}}{\bm \gamma}}&= {\bm v}^\pm \cdot \partial_t \partial_{\mathsf{j}} {\bm \gamma} + \pare{\partial_t {\bm u^\pm}}\pare{{\bm \gamma}} \cdot\partial_{\mathsf{j}}
{\bm \gamma}+ \pare{J_{\bm x} {\bm u^\pm} \pare{{\bm \gamma}} \partial_t{\bm \gamma}} \cdot \partial_{\mathsf{j}} {\bm \gamma}, 
\end{aligned}
\end{equation}
and
\begin{equation}
\label{eq:expansion2}
\begin{aligned}
\partial_{\mathsf{j}}\pare{{\bm v}^\pm\cdot\partial_t{\bm \gamma}} &={\bm v}^\pm \cdot \partial_t \partial_{\mathsf{j}} {\bm \gamma} +  \pare{J_{\bm x}{\bm u}^{\pm}\pare{{\bm \gamma}}\partial_{\mathsf{j}}{\bm \gamma}}\cdot \partial_t{\bm \gamma} , 
\end{aligned}
\end{equation}
and
\begin{equation}
\label{eq:expansion3}
   \partial_{\mathsf{j}}\pare{\frac{\av{{\bm v}^\pm}^2}{2}}= \pare{ {\bm v}^\pm }^\intercal J_{\bm x} {\bm u}^{\pm}\pare{{\bm \gamma}} \partial_{\mathsf{j}} {\bm \gamma}.
\end{equation}
Observe that from \cref{eq:expansion1} and \cref{eq:expansion2} we obtain
\begin{equation}\label{eq:expansion3.1}
    \partial_t\pare{{\bm v}^\pm\cdot\partial_{\mathsf{j}}{\bm \gamma}} - \partial_{\mathsf{j}}\pare{{\bm v}^\pm\cdot\partial_t{\bm \gamma}} = \pare{\partial_t {\bm u^\pm}}\pare{{\bm \gamma}} \cdot\partial_{\mathsf{j}} {\bm \gamma} + 
     \pare{\partial_{\mathsf{j}} {\bm \gamma}}^\intercal \pare{J_{\bm x} {\bm u}^\pm\pare{{\bm \gamma}}\partial_t {\bm \gamma}}  - \pare{\partial_t {\bm \gamma}}^\intercal \pare{J_{\bm x} {\bm u}^\pm\pare{{\bm \gamma}}\partial_{\mathsf{j}} {\bm \gamma}} . 
\end{equation}
Notice that thanks to \cref{eq:velocity_potential} we have that  $J_{{\bm x}} {\bm u^\pm}={\bm H}_{{\bm x}} \phi^\pm$. Hence by Schwarz Theorem $J_{{\bm x}} {\bm u^\pm}$ is a symmetric matrix so that
\begin{equation}
    \label{eq:expansion4}
 \pare{\partial_{\mathsf{j}} {\bm \gamma}}^\intercal \pare{J_{\bm x} {\bm u}^\pm\pare{{\bm \gamma}}\partial_t {\bm \gamma}}  - \pare{\partial_t {\bm \gamma}}^\intercal \pare{J_{\bm x} {\bm u}^\pm\pare{{\bm \gamma}}\partial_{\mathsf{j}} {\bm \gamma}} =0. 
\end{equation}
Plugging \cref{eq:expansion4} into \cref{eq:expansion3.1}, we infer
\begin{equation}
\label{eq:expansion5}
     \partial_t\pare{{\bm v}^\pm\cdot\partial_{\mathsf{j}}{\bm \gamma}} - \partial_{\mathsf{j}}\pare{{\bm v}^\pm\cdot\partial_t{\bm \gamma}} =\pare{\partial_t {\bm u^\pm}}\pare{{\bm \gamma}} \cdot\partial_{\mathsf{j}} {\bm \gamma}. 
\end{equation}
The trace of  the first equation of \cref{eq:2phaseWW3D} onto $ \Gamma\pare{t} $ is 
$$\pare{\partial_t {\bm u^\pm}}\pare{{\bm \gamma}} = -\pare{  J_{\bm x} {\bm u}^{\pm}\pare{{\bm \gamma}} {\bm v}^\pm + \frac{1}{\rho^\pm} \nabla_{\bm x}p^\pm\pare{{\bm \gamma}}     +\mathfrak{g} {\bm e}_3 }.$$
Inserting this into \cref{eq:expansion5} then combining with \cref{eq:expansion3} and the fact that the matrix $ J_{\bm x} {\bm u}^{\pm}\pare{{\bm \gamma}} $ is symmetric, we get
\begin{equation}\label{eq:expansion6}
     \partial_t\pare{{\bm v}^\pm\cdot\partial_{\mathsf{j}}{\bm \gamma}} - \partial_{\mathsf{j}}\pare{{\bm v}^\pm\cdot\partial_t{\bm \gamma}} = -  \partial_{\mathsf{j}}\pare{\frac{\av{{\bm v}^\pm}^2}{2}}
     -\pare{  \frac{1}{\rho^\pm} \nabla_{\bm x} p^\pm \pare{{\bm \gamma}} +\mathfrak{g}{\bm e}_3} \cdot\partial_{\mathsf{j}} {\bm \gamma}. 
\end{equation}
Remark that
\begin{equation}
\pare{ \frac{1}{\rho^\pm} \nabla_{\bm{x}} p^\pm\pare{{\bm \gamma}} + \mathfrak{g} {\bm e}_3 }\cdot \partial_{\mathsf{j}}{\bm \gamma} = \nabla_{x, z} \left( \frac{p^\pm\pare{{\bm \gamma}}}{\rho^\pm} + \mathfrak{g} z \right) \cdot \partial_{\mathsf{j}}{\bm \gamma} = \partial_{\mathsf{j}} \pare{  \frac{q^\pm}{\rho^\pm} + \mathfrak{g} {\bm \gamma}_3 }.
\end{equation}
Thus,
\begin{align}\label{eq:expansion6.1}
     \partial_t\pare{{\bm v}^\pm\cdot\partial_{\mathsf{j}}{\bm \gamma}} - \partial_{\mathsf{j}}\pare{{\bm v}^\pm\cdot\partial_t{\bm \gamma}} = -  \partial_{\mathsf{j}}\pare{\frac{\av{{\bm v}^\pm}^2}{2} +   \frac{q^\pm}{\rho^\pm} + \mathfrak{g} {\bm \gamma}_3 }, 
     && \mathsf{j}=1,2. 
\end{align}
Hence, \cref{eq:expansion6.1} and \cref{operator2} imply that
\begin{equation}
\pare{ \partial_t\pare{{\bm v}^\pm\cdot\partial_2{\bm \gamma}} \partial_1{\bm \gamma} - \partial_t\pare{{\bm v}^\pm\cdot\partial_1{\bm \gamma}}\partial_2{\bm \gamma}}
-
\pare{  \partial_2\pare{{\bm v}^\pm\cdot\partial_t{\bm \gamma}} \partial_1{\bm \gamma} - \partial_1\pare{{\bm v}^\pm\cdot\partial_t{\bm \gamma}}\partial_2{\bm \gamma}}
=
  J {\bm \gamma}  \nabla^\perp \pare{\frac{\av{{\bm v}^\pm}^2}{2} +   \frac{q^\pm}{\rho^\pm} + \mathfrak{g} {\bm \gamma}_3 }. 
\end{equation}
Combined with \cref{eq:expansion0} and obtain  that
\begin{equation}
\label{eq:expansion8}
\partial_t \pare{{\bm v}^\pm\times {\bm \nu}} + \mathcal{I}^\pm 
=
J {\bm \gamma}  \nabla^\perp \pare{\frac{\av{{\bm v}^\pm}^2}{2} +   \frac{q^\pm}{\rho^\pm} + \mathfrak{g} {\bm \gamma}_3 }, 
\end{equation}
where
\begin{align}   
\mathcal{I}^\pm \defeq & \ \mathcal{I}^\pm_1+\mathcal{I}^\pm_2,
\\
\mathcal{I}^\pm_1\defeq & \ \pare{{\bm v}^\pm\cdot\partial_1 {\bm \gamma}}\partial_t\partial_2{\bm \gamma}-\pare{{\bm v}^\pm\cdot\partial_2 {\bm \gamma}}\partial_t\partial_1{\bm \gamma} &=& \ -J^\perp{\bm \gamma} _t \pare{J{\bm \gamma}}^\intercal {\bm v}^\pm,
\\
\mathcal{I}^\pm_2\defeq & \ - \partial_2\pare{{\bm v}^\pm\cdot\partial_t{\bm \gamma}}\partial_1{\bm \gamma} + \partial_1\pare{{\bm v}^\pm\cdot\partial_t{\bm \gamma}}\partial_2{\bm \gamma}
&=& \ -J^\perp {\bm \gamma}\nabla\pare{{\bm v}^\pm \cdot {\bm \gamma}_t}.
\end{align}
Thus, we have
\begin{equation}
\label{equationI0}
\mathcal{I}^\pm=-\pare{J^\perp {\bm \gamma}\nabla\pare{{\bm v}^\pm \cdot {\bm \gamma}_t} + J^\perp{\bm \gamma} _t \pare{J{\bm \gamma}}^\intercal {\bm v}^\pm}.
\end{equation}
We insert \cref{equationI0} into \cref{eq:expansion8} and obtain that
\begin{equation}\label{eq:expansion9}
M_1\bra{{\bm \gamma}, {\bm \gamma}} \bra{{\bm v}^\pm} \defeq\partial_t \pare{{\bm v}^\pm\times {\bm \nu}} 
-
\pare{
J^\perp {\bm \gamma}\nabla\pare{{\bm v}^\pm \cdot {\bm \gamma}_t} + J^\perp{\bm \gamma} _t \pare{J{\bm \gamma}}^\intercal {\bm v}^\pm 
}
=
J {\bm \gamma}  \nabla^\perp \pare{\frac{\av{{\bm v}^\pm}^2}{2} +   \frac{q^\pm}{\rho^\pm} + \mathfrak{g} {\bm \gamma}_3 }. 
\end{equation}
From \cref{eq:expansion9} and \cref{eq:boundary_condition_pressure} it is immediate that 
\begin{equation}
\label{eq:difference_momentum_trace}
M_1\bra{{\bm \gamma},{\bm \gamma}}\bra{ \bbra{ \rho^\pm {\bm v}^\pm  } } 
=
J {\bm \gamma}  \nabla^\perp \pare{\frac{\bbra{ \rho^\pm\av{{\bm v}^\pm}^2 }}{2} +  \sigma \mathpzc{k}\pare{{\bm \gamma}} + \mathfrak{g} \bbra{\rho^\pm} {\bm \gamma}_3 }.
\end{equation}
From \cref{eq:BOmega}, \cref{eq:traces_velocity} and \cref{eq:Atwood}, we have that
\begin{equation}
\bbra{ \rho^\pm {\bm v}^\pm } = \bbra{ \rho^\pm } {\bm v}_N-  \aangle{\rho^\pm} {\bm v}_T.
\end{equation}
Therefore, using \cref{eq:BOmega} and \cref{eq:omega3}, we infer
\begin{equation}
\label{eq:simplification_eq_omegaQeta1}
\begin{aligned}
\partial_t \pare{\bbra{\rho^\pm {\bm v}^\pm }\times{\bm \nu}} = -& \  \aangle{\rho^\pm} \partial_t {\bm \omega}+ \bbra{ \rho^\pm } \partial_t\pare{{\bm v}_N\times {\bm \nu}}.
\end{aligned}
\end{equation}
We use now \cref{eq:BOmega} and \cref{eq:traces_velocity} in order to argue that
\begin{equation}
\label{eq:simplification_eq_omega_2}
\bbra{ \rho^\pm \av{{\bm v}^\pm}^2 } =\bbra{ \rho^\pm } \pare{ {\bm v}_N ^2 +\frac{{\bm v}_T^2}{4}}
-2\aangle{\rho^\pm} {\bm v}_N\cdot{\bm v}_T.
\end{equation}
Inserting \cref{eq:simplification_eq_omegaQeta1}, \cref{eq:simplification_eq_omega_2} and \cref{eq:Atwood} into \cref{eq:difference_momentum_trace}, we get
\begin{equation}\label{eq:simplification_eq_omega_2.1}
\begin{aligned}
    -{\bm \omega}_t + \Atw \partial_t \pare{{\bm v}_N\times{\bm \nu}}
&=
J^\perp {\bm \gamma} \nabla\pare{\pare{-{\bm v}_T + \Atw {\bm v}_N}\cdot{\bm \gamma}_t}
+
J^\perp {\bm \gamma}_t \pare{J{\bm \gamma}}^\intercal \pare{-{\bm v}_T + \Atw {\bm v}_N}
\\
&\quad+
J{\bm \gamma} \nabla^\perp 
\pare{
\frac{\Atw}{2} \pare{{\bm v}_N^2 + \frac{{\bm v}_T^2}{4}} -   {\bm v}_N \cdot{\bm v}_T +\frac{\sigma}{\aangle{\rho^\pm}}\mathpzc{k}\pare{{\bm \gamma}} + \mathfrak{g}\Atw{\bm \gamma}_3
}. 
\end{aligned}
\end{equation}With standard algebraic manipulations we obtain \cref{eq:simplification_eq_omega_3}.

\end{step}

\begin{step}[The evolution under the graph assumption \cref{eq:surface}]

Under the graph assumption \cref{eq:surface}, we have that
\begin{align}\label{eq:nablazetaspecifico}
J{\bm \gamma}=\vect{\Id_{\bR^2}}{\pare{\nabla \eta}^\intercal }, 
&&
J^\perp {\bm \gamma}=\vect{\vect{0 & 1}{-1 & 0}}{\pare{\nabla^\perp \eta}^\intercal },
\end{align}
so that
\begin{align}
J{\bm \gamma} \nabla^\perp = \vect{\nabla^\perp}{\nabla\eta\cdot\nabla^\perp},
&&
J^\perp {\bm \gamma} \nabla = \vect{-\nabla^\perp}{-\nabla\eta\cdot\nabla^\perp}.
\end{align}
Combined with the definition of $ {\bm \pi}^\dagger \pare{\rho^\pm} \defeq \vect{ \Pi^\dagger \pare{\rho^\pm}}{ \pi ^\dagger _3\pare{\rho^\pm}} $ given in \cref{eq:pi_dagger_def} implies that
\begin{equation}\label{eq:pi}
\begin{aligned}
\Pi^\dagger \pare{\rho^\pm\middle| \pare{\bullet, \eta }, {\bm \omega}} = & \
\nabla^\perp \phi_\Pi \pare{\rho^\pm\middle| \eta,{\bm \omega}},
\\
\phi_\Pi \pare{\rho^\pm\middle| \eta,{\bm \omega}} \defeq & \  \phi_{\Pi;0} \pare{\rho^\pm \middle| \eta } + \phi_{\Pi;2} \pare{\Atw \middle| \eta }\bra{ {\bm \omega}, {\bm \omega}},
\\
\phi_{\Pi;0} \pare{\rho^\pm \middle| \eta } \defeq & \ - \frac{\sigma}{\aangle{\rho^\pm}}\nabla\cdot\mathpzc{K}\pare{\nabla\eta} + \mathfrak{g}\Atw \eta ,
\\
\phi_{\Pi;2} \pare{\Atw \middle| \eta }\bra{ {\bm \omega}, {\bm \omega}}
\defeq & \ 
\underbrace{{\bm v}_T\cdot\pare{\vect{0_{\bR^2}}{\eta_t} -{\bm v}_N}}_{\eqdef \phi_{\Pi,\RN{1};2} \pare{ \eta }\bra{ {\bm \omega}, {\bm \omega}}} +  
\underbrace{\frac{\Atw }{2} \pare{{\bm v}_N^2 + \frac{{\bm v}_T^2}{4}-2v_{N;3} \eta_t} } _{\eqdef \phi_{\Pi,\RN{2};2} \pare{\Atw \middle| \eta }\bra{ {\bm \omega}, {\bm \omega}}}.
\end{aligned}
\end{equation}
Next, we have that
\begin{align}\label{eq:N=0}
J^\perp {\bm \gamma}_t =\vect{0_{\bR^{2\times 2}}}{\pare{\nabla^\perp \eta_t}^\intercal}
\ \text{and} \ \tilde{{\bm {n}}} = \vect{\tilde{N}}{\tilde{n}_3}
&& \Rightarrow && \tilde{N}\pare{\pare{\bullet,\eta}, {\bm \omega}}  \equiv 0. 
\end{align}
Moreover letting ${\bm v}_N\defeq\vect{V_N}{v_{N;3}}$ we have that
$$\pare{J {\bm \gamma}}^\intercal {\bm v}_N = V_N + v_{N;3} \nabla \eta,$$
	which we combine with \cref{eq:nablazetaspecifico} and $\tilde{{\bm{m}}}=\vect{\tilde{M}}{m_3}$ in order to obtain that
			\begin{equation}\label{eq:Th}
				\left. \tilde{M}\pare{{\bm \gamma}}{\bm \omega} \right|_{{\bm \gamma}\pare{x} = \pare{x, \eta\pare{x}}} =\pare{ V_N + v_{N;3} \nabla \eta}^\perp. 
			\end{equation}
			Now using \cref{eq:BOmega}, \cref{eq:bmvtilde} and \cref{eq:omega3}, we get
			\begin{equation}\label{eq:T(eta)}
				\begin{aligned}
					T \pare{\eta}\Omega \defeq & \ 
					\left.\tilde{M}\pare{{\bm \gamma}}{\bm \omega} \right|_{\substack{{\bm \gamma}\pare{x} = \pare{x, \eta\pare{x}}\\ \omega_3 = \nabla\eta\cdot\Omega}}
					\\
					= & \ \mathpzc{B}_0\pare{\eta} \bra{\nabla\eta\cdot\Omega} - \mathpzc{q}\pare{\eta}\Omega + \pare{\mathpzc{B}_0 \pare{\eta}^\perp \cdot \Omega} \nabla^\perp \eta
					\\
					= & -\pare{ V^{\star}\bra{\Omega}\pare{\eta} +  v_3^{\star} \bra{\Omega}\pare{\eta} \ \nabla\eta}^\perp.
				\end{aligned} 
			\end{equation}
			Plugging \cref{eq:Th}, \cref{eq:T(eta)}, \cref{eq:N=0} and \cref{eq:pi} into \cref{eq:simplification_eq_omega_3}, we infer
			\begin{equation}
				\label{eq:simplification_eq_omega_4}
				\partial_t \pare{ \pare{ 1+\Atw T \pare{\eta} } \Omega  } =  -\left. \nabla^\perp \phi_\Pi \pare{\rho^\pm\middle| \eta,{\bm \omega}} \middle|_{ \omega_3 =\nabla\eta\cdot\Omega } \right.. 
			\end{equation}
			We use now the special structure of the interface $ {\bm \gamma} $ and ${\bm \omega}$, given in \cref{eq:surface} and \cref{eq:omega3}, in order to simplify the term  $\left.  \phi_\Pi \pare{\rho^\pm\middle| \eta,{\bm \omega}} \middle|_{ \omega_3 =\nabla\eta\cdot\Omega } \right.  $ in \cref{eq:simplification_eq_omega_4}.

			\begin{enumerate}

				\item {\bf Analysis of $\phi_{\Pi;0} \pare{\rho^\pm \middle| \eta } $:} We define
				\begin{equation}\label{eq:pi3}
					\begin{aligned}
						\pi_0\pare{\rho^\pm\middle| \eta }&\defeq
						-\phi_{\Pi;0} \pare{\rho^\pm \middle| \eta } .
					\end{aligned}  
				\end{equation}
				
				\item {\bf Analysis of $\phi_{\Pi,\RN{1};2} \pare{ \eta }\bra{ {\bm \omega}, {\bm \omega}}$:} From \cref{eq:evolution_eta_omegah} and \cref{eq:gstar}, we have $\eta_t+ v_3^{\star}=\nabla\eta\cdot  V^{\star}$. Observe that for any $A,B,C\in\mathbb R^2$, we have
				\begin{equation}
					\label{eq:formulascalar1}
					\pare{A\cdot B}C+\pare{C^\perp \cdot B}A^\perp=\pare{A\cdot C}B.
				\end{equation}
				We combine  \cref{eq:formulascalar1,eq:evolution_eta_omegah,eq:bmvtilde,eq:CDEeta} and obtain that{{
\begin{align}
\pi_{2;\RN{1}}^\star\pare{\eta }\bra{\Omega, \Omega }\defeq
\left. \phi_{\Pi,\RN{1};2} \pare{ \eta }\bra{ {\bm \omega}, {\bm \omega}}
\middle|_{\omega_3 =\nabla\eta\cdot\Omega } \right.={V_T\cdot  V^{\star}+\frac{\p{\nu\times \bm\omega}_3}{1+\av{\nabla\eta}^2}\pare{\nabla\eta\cdot  V^{\star}}}.
\end{align}}}
Making appeal to \cref{eq:BOmega}, \cref{boundary1} and \cref{eq:omega3}, we find
\begin{align}
     \pi_{2;\RN{1}}^\star\pare{\eta }\bra{\Omega, \Omega }&=\frac{1}{1+\av{\nabla\eta}^2}\bra{\pare{-\omega_2-\partial_2\eta\pare{\nabla\eta\cdot \Omega}} V^{\star}_1+\pare{\omega_1+\partial_1\eta\pare{\nabla\eta\cdot \Omega}} V^{\star}_2+\pare{-\partial_1\eta\omega_2+\partial_2\eta\omega_1}\pare{\nabla\eta\cdot  V^{\star}}}
\\
&= \frac{1}{1+\av{\nabla\eta}^2}\bra{-\pare{1+\av{\nabla\eta}^2}\omega_2  V^{\star}_1+\pare{1+\av{\nabla\eta}^2}\omega_1  V^{\star}_2}\\
&= \Omega^\perp\cdot  V^{\star} \bra{\Omega}\pare{\eta}.
\end{align}
Therefore,
\begin{equation}\label{quadra:pi2Istar}
    \pi_{2;\RN{1}}^\star\pare{\eta }\bra{\Omega, \Omega }=\Omega^{\perp}\cdot V^{\star} \bra{\Omega}\pare{\eta}.
\end{equation}
By polarization, we obtain \cref{eq:pi2_star_I_def_lemma}.
Since \cref{eq:simplification_eq_omega_4} gives
\[\partial_t\pare{\pare{1+\Atw T\pare{\eta}}\Omega}=-\nabla^\perp\phi_\Pi,\]
the quantity in \cref{quadra:pi2Istar} enters the Hamiltonian-coordinate right-hand side through $\pi^\star=-\phi_\Pi$. Consequently, the contribution of $\pi^\star_{2;\RN{1}}$ to $\pi^\star_2$ carries the minus sign displayed in \cref{eq:pistar2}. After the substitution $\Omega=S\pare{\Atw\middle|\eta}\psi$, this is the origin of the minus sign in the definition of $\pi^\star_{2;\RN{1}}\pare{\Atw\middle|\eta}$, c.f. \cref{eq:pi2tilde}.

\item {\bf Analysis of $\phi_{\Pi,\RN{2};2} \pare{\Atw \middle| \eta }\bra{ {\bm \omega}, {\bm \omega}}$:} Let us recall from \cref{eq:pi} that
\begin{equation}
\phi_{\Pi,\RN{2};2} \pare{\Atw \middle| \eta }\bra{ {\bm \omega}, {\bm \omega}}
=
\left.
    \frac{\Atw }{2} \pare{{\bm v}_N^2 + \frac{{\bm v}_T^2}{4}-2{\bm v}_N\cdot{\bm \gamma}_t} 
    \right|_{\substack{{\bm \gamma}\pare{x} = \pare{x, \eta\pare{x}} \\ \omega_3 =\nabla\eta\cdot\Omega }} .
\end{equation}
From ${\bm \omega}\cdot{\bm \nu}=0$, we have \[|{\bm \nu}\times{\bm \omega}|^2=|{\bm \nu}|^2 |{\bm \omega}|^2 - \pare{{\bm \nu} \cdot {\bm \omega}}^2 \iff |{\bm \nu}\times{\bm \omega}| = |{\bm \nu}| |{\bm \omega}|.\] So\begin{equation}\label{eq:Omegasquared}
\left.     \av{{\bm v}_T}^2 \middle|_{\substack{{\bm \gamma}\pare{x} = \pare{x, \eta\pare{x}} \\ \omega_3 =\nabla\eta\cdot\Omega }} \right. =\left. \frac{|{\bm \omega}|^2}{|{\bm \nu}|^2}\middle|_{\substack{{\bm \gamma}\pare{x} = \pare{x, \eta\pare{x}} \\ \omega_3 =\nabla\eta\cdot\Omega }} \right.=\frac{|\Omega|^2+|\nabla\eta\cdot\Omega|^2}{1+|\nabla\eta|^2}\cdot
\end{equation}
For the term $\av{{\bm v}_N}^2$  using \cref{eq:bmvtilde} we have that
\begin{equation}\label{eq:Bsquared}
    \left. \av{{\bm v}_N}^2 \right|_{\substack{{\bm \gamma}\pare{x} = \pare{x, \eta\pare{x}} \\ \omega_3 =\nabla\eta\cdot\Omega }} = \av{{\bm v}^\star\bra{\Omega}\pare{\eta}}^2 = \av{ V^{\star}\bra{\Omega}\pare{\eta}}^2 + \av{ v_3^{\star}\bra{\Omega}\pare{\eta}}^2 . 
\end{equation}

For the term ${\bm v}_N\cdot{\bm \gamma}_t$ we have by virtue of \cref{eq:gstar}
\begin{equation}\label{eq:Betat}
\begin{aligned}
    \left. {\bm v}_N\cdot{\bm \gamma}_t
\middle|_{\substack{{\bm \gamma}\pare{x} = \pare{x, \eta\pare{x}} \\ \omega_3 =\nabla\eta\cdot\Omega }} \right. &=- v_3^{\star}\bra{\Omega}\pare{\eta} \ \eta_t= -  v_3^{\star}\bra{\Omega}\pare{\eta}  \  \mathpzc{g}^\star \pare{\eta} \Omega
\\
&=
-  v_3^{\star}\bra{\Omega}\pare{\eta}
\pare{
 V^{\star}\bra{\Omega}\pare{\eta}\cdot\nabla\eta -  v_3^{\star}\bra{\Omega}\pare{\eta}
}.
\end{aligned}
\end{equation}
We combine \cref{eq:Omegasquared,eq:Bsquared,eq:Betat} in the definition of $\tilde\pi_2 \pare{\Atw\middle|{\bm \gamma}, \omega}$ (cf. \cref{eq:pi2tilde}) and obtain that
%\paragraph{Proof of \cref{prop:paralinearization_pi}}
\begin{equation}\label{eq:pi2II_quadratic}
      \begin{aligned}
          \pi^\star_{2;\RN{2}}\pare{\Atw\middle| \eta }\bra{\Omega, \Omega}&\defeq
\left. \phi_{\Pi,\RN{2};2} \pare{\Atw \middle| \eta }\bra{ {\bm \omega}, {\bm \omega}}
\middle|_{ \omega_3 =\nabla\eta\cdot\Omega } \right.\\
    &=\frac{\Atw }{2}
    \pare{
    \av{ V^{\star}\bra{\Omega}\pare{\eta}}^2 
    -  v_3^{\star}\bra{\Omega}\pare{\eta} \pare{ v_3^{\star}\bra{\Omega}\pare{\eta} - 2 V^{\star}\bra{\Omega}\pare{\eta}\cdot\nabla\eta}
    +
    \frac{1}{4}\frac{|\Omega|^2+|\nabla\eta\cdot\Omega|^2}{\p{1+|\nabla\eta|^2}}
    }.
      \end{aligned}  
\end{equation}
Once again,  the associated bilinear form \cref{eq:pi2_star_II_def_lemma} is obtained  via polarization.

\end{enumerate}
\noindent \cref{eq:pistar2} follows by summation, and we conclude.
\end{step}
\end{proof}

\subsection{\texorpdfstring{The CDE for the variable $\psi$}{The CDE for the variable ψ}}\label{sec:psi}

In this subsection, we prove the following lemma dealing with the evolution of the weighted trace of the potential $\psi$.

\begin{lemma}[Evolution of the weighted trace of the velocity potential]\label{lempsi}
The function $\psi$  defined via \cref{eq:psi,eq:psipm}, satisfies
\begin{equation}\label{link:omg-psi}
    \Omega=S\pare{\Atw\middle|\eta}\psi,
\end{equation}
where $S\pare{\Atw\middle|\eta}$ is the operator defined in \cref{eq:S}. Moreover, there exist $s_0 > 0$ and $r_0 > 0$ such that for any $\eta\in \BallR{s_0}{r_0}$ the function $\psi$
solves the  evolutionary equation \cref{eq:evolution_psi}, i.e.
\begin{equation}\label{eqpsifinal}
    \psi_t=\pi\pare{\rho^{\pm}\middle|\eta,\psi},
\end{equation}
where the right hand-side is defined by
    \begin{equation}
    \label{eq:pi_as_pistar_composed_with_S}\pi\pare{\rho^\pm\middle| \eta, \psi}\defeq\pi^\star\pare{\rho^\pm\middle| \eta ,  S\pare{\Atw\middle|\eta}\psi},\end{equation} 
    with $\pi^\star\pare{\rho^\pm\middle| \eta, \Omega}$ as in \cref{eq:pistar}.
\end{lemma}

\begin{proof}
    Let us recall that $\psi$ is defined in \cref{eq:psi}. Let us define the \emph{modified vorticity}
\begin{equation}
\label{eq:zeta}
    {\bm \zeta}\defeq\vect{\mathrm{Z}}{\zeta_3} \defeq - \frac{ \bbra{\rho^\pm {\bm v}^\pm}\times {\bm \nu}}{\aangle{\rho^\pm}}\cdot 
\end{equation}
Standard computations combined with \cref{vorticity,eq:traces_velocity} and the first identity of \cref{eq:TNP} show that
\begin{equation}
\label{eq:zeta_omega_relation}
\begin{aligned}
        {\bm \zeta} = & {\bm \omega}- \Atw \BR{{\bm \gamma}}{\bm \omega}\times {\bm \nu} 
        =  \pare{1 + \Atw \tilde{\bm{m}}\pare{{\bm \gamma}}} \bra{{\bm \omega}}. 
        \end{aligned}
\end{equation}
We plug \cref{eq:zeta_omega_relation} into \cref{eq:simplification_eq_omega_3} and obtain that
\begin{equation}\label{eq:zeta_eq_generic}
    {\bm \zeta}_t = -{\bm \pi}^\dagger \pare{\rho^\pm\middle| {\bm \gamma}, {\bm \omega}} +
\tilde{{\bm{n}}} \pare{{\bm \gamma}}\bra{{\bm \omega}}. 
\end{equation}
Let us now reduce the first two components \cref{eq:zeta_eq_generic} to the case in which ${\bm \gamma}$ satisfies the graph assumption \cref{eq:surface}. The transformed r.h.s. of \cref{eq:zeta_eq_generic} is given by the r.h.s. of \cref{eq:evolution_omegah}, so that
\begin{equation}\label{eq:evolution_zeta_elevation}
      \partial_t \mathrm{Z} = \nabla^\perp \pi^\star \pare{\rho^\pm\middle|\eta, \Omega}. 
\end{equation}
Notice now that, defining $ {\bm v}^\pm \defeq\vect{V^\pm}{v_3^\pm} $, we have 
\begin{equation}
\label{eq:psiimportant}
\nabla \psi^{\pm}\pare{x} = \nabla \phi^\pm\pare{x, \eta\pare{x}} + \partial_{x_3}\phi^\pm\pare{x, \eta\pare{x}} \nabla\eta\pare{x}
=
V^\pm \pare{x} + v^\pm_3 \pare{x}\nabla\eta\pare{x}, 
\end{equation}
so that
\begin{equation}\label{eq:nablapsi}
\nabla\bbra{\psi^\pm} = \bbra{V^\pm} + \bbra{v^\pm_3} \nabla \eta. 
\end{equation}
Now from \cref{vorticity} we have that
\begin{equation}\label{eq:omegaexplicit}
{\bm \omega}=\vect{\pare{ \bbra{V^\pm} + \bbra{v^\pm_3}\nabla \eta}^\perp}{-\nabla^\perp \eta\cdot \bbra{V^\pm}}, 
\end{equation}
so that, from \cref{eq:nablapsi,eq:omegaexplicit}, we obtain that
\begin{equation}\label{eq:omega(psi)}
 \Omega = \nabla^\perp \bbra{\psi^\pm }. 
\end{equation}
Similar computations moreover give that
\begin{equation}\label{eq:zeta(psi)}
    \mathrm{Z} = \nabla^\perp \psi. 
\end{equation}
We plug \cref{eq:omega(psi),eq:zeta(psi)} into \cref{eq:evolution_zeta_elevation}  in order to express the evolution of $\psi$ in terms of $\eta$ and $\bbra{\psi^\pm}$, this gives
\begin{equation}\label{eq:psi([psi])_evolution}
 \psi_t = 
    \pi^\star \pare{\rho^\pm\middle|\eta, \nabla^\perp \bbra{\psi^\pm}}. 
\end{equation}
Let us now relate $\psi$ and $\bbra{\psi^\pm}$ via an invertible, bounded operator. We combine \cref{eq:zeta_omega_relation,eq:T(eta)} and obtain that
\begin{equation}\label{eq:zeta_omega_relation_elevation}
    \mathrm{Z} =  \pare{1+\Atw T\pare{\eta}}\bra{\Omega }.
\end{equation}
Putting together \cref{eq:omega(psi),eq:zeta(psi),eq:zeta_omega_relation_elevation}, we obtain that
\begin{equation}\label{eq:psi([psi])_0}
    \nabla^\perp\psi  =  \pare{1+\Atw T\pare{\eta}} \bra{\nabla^\perp\bbra{\psi^\pm}}. 
\end{equation}
\cref{eq:psi([psi])_0} has no curl-free component, hence we can apply the operator $\Delta^{-1}\nabla^\perp\cdot \bullet $ to it and obtain that
\begin{equation}\label{eq:psi([psi])}
    \psi =   \pare{1+\Atw \tau_0 \pare{\eta}} \bra{\bbra{\psi^\pm}}, 
\end{equation}
where (cf. \cref{eq:T(eta)})
\begin{equation}\label{eq:Tstar(T)}
    \tau_0\pare{\eta}\defeq  \Delta^{-1}\nabla^\perp\cdot T\pare{\eta}\nabla^\perp. 
\end{equation}
Notice that there exist $s_0 > 0$ and $r_0 > 0$ such that for all $\eta\in \BallR{s_0}{r_0}$ the operator $ \pare{1+\Atw \tau_0 \pare{\eta}}\in \cL\pare{H ^s}$ is invertible for all $s\geq s_0$ (the proof is postponed to \cref{lem:uniform_inverse_tau0}) ,  so we have that
\begin{equation}\label{eq:[psi](psi)}
    \bbra{\psi^\pm} = \pare{1+\Atw \tau_0\pare{\eta}}^{-1} \psi.
\end{equation}
It will be convenient to distinguish the unweighted trace jump from the
canonical weighted trace. We therefore set
\begin{equation}\label{eq:trace_jump_potential_def}
    \varphi
    \defeq
    \bbra{\psi^\pm}
    =
    \pare{1+\Atw\tau_0\pare{\eta}}^{-1}\psi .
\end{equation}
Thus
\begin{equation}\label{eq:canonical_weighted_trace_def}
    \psi
    =
    \pare{1+\Atw\tau_0\pare{\eta}}\varphi,
    \qquad
    \Omega
    =
    \nabla^\perp\varphi .
\end{equation}
All these identities are understood in homogeneous spaces, hence modulo
additive constants.
Plugging \cref{eq:[psi](psi)} into \cref{eq:psi([psi])_evolution} and defining $\pi\pare{\rho^\pm\middle| \eta , \psi }$ as in the statement we obtain \cref{eq:evolution_psi}.
\end{proof}
\subsection{Recasting the potential in term of vorticity}\label{sec:vorticitytopotential}
{{

Recall $\psi^\pm\p{t,x}$ defined in \cref{eq:psipm}. Let us now express $ \psi^\pm, \psi $ in terms of $ \Omega$.

\begin{lemma}\label{lem:psi_omegah}
    Let $\psi^\pm$ and $\psi$ be as in \cref{eq:psipm,eq:psi} respectively, there exists a $s_0 > 0$ and a $r_0 > 0$ such that if $\eta\in \BallR{s_0}{r_0}, \Omega\in\p{\BallR{s_0}{r_0}}^2$ then
    {{
    \begin{subequations}
         \begin{align}\label{eq:psipm(Omega)}
\psi^\pm\p x= &\Delta^{-1}\nabla^\perp\cdot\p{\mp\frac 12+T\p\eta}\bra{\Omega},
\\
\label{eq:psi(Omega)}
\psi\p x= & \ \Delta^{-1}\nabla^\perp\cdot\p{ 1+\Atw T\p\eta}\bra{\Omega},
\end{align}
    \end{subequations}
    }}
where $T\pare{\eta}$ is defined in \cref{eq:T(eta)}. Moreover, we have
   \begin{equation}
   \label{eq:Tdivergencefree}
T\pare{\eta}\Omega=\nabla^\perp\p{\frac{\psi^++\psi^-}2}\p x.
   \end{equation}
   
\end{lemma}

\begin{proof}
Recall $\omega=\vect{\Omega}{\nabla\eta\cdot\Omega}$. By \cref{eq:psiimportant} and \cref{eq:traces_velocity}, we get
\begin{equation}
\label{eq:psiimportant2}
\begin{aligned}
{V^\pm \pare{x} + v^\pm_3 \pare{x}\nabla\eta\pare{x}}={\p{\textnormal{BR}\p{\gamma}\omega\pm\frac 12\frac{\nu\times\omega}{\av\nu^2}}_h+\p{\textnormal{BR}\p{\gamma}\omega\pm\frac 12\frac{\nu\times\omega}{\av\nu^2}}_3\nabla\eta}.
\end{aligned}
\end{equation}
Recall \cref{eq:bmvtilde} and \cref{eq:T(eta)}, we have
\begin{equation}
\label{eq:psiimportant3}
\begin{aligned}
{\p{\textnormal{BR}\p{\gamma}\omega}_h+\p{\textnormal{BR}\p{\gamma}\omega}_3\nabla\eta}
&={-\mathpzc{B}_0\pare{\eta}^\perp \bra{\nabla\eta \cdot \Omega}+ \mathpzc{q}\pare{\eta}\Omega^\perp+\p{\mathpzc{B}_0\pare{\eta}^\perp \cdot \Omega}\nabla\eta}\\
&=-\p{ T\p\eta\Omega}^\perp,
\end{aligned}
\end{equation}
and
\begin{equation}
\label{eq:psiimportant4}
\begin{aligned}
\p{\frac{\nu\times\omega}{\av\nu^2}}_h+\p{\frac{\nu\times\omega}{\av\nu^2}}_3\nabla\eta=\vect{-\Omega_2}{\Omega_1}=\Omega^\perp.
\end{aligned}
\end{equation}
So
\begin{equation}
\label{eq:psiimportant5}
\begin{aligned}
\nabla\cdot\p{\p{\textnormal{BR}\p{\gamma}\omega}_h+\p{\textnormal{BR}\p{\gamma}\omega}_3\nabla\eta}=\nabla^\perp \cdot T\p\eta\Omega,\quad\nabla\cdot\p{\p{\frac{\nu\times\omega}{\av\nu^2}}_h+\p{\frac{\nu\times\omega}{\av\nu^2}}_3\nabla\eta}=\nabla\cdot\Omega^\perp=-\nabla^\perp\cdot \Omega,
\end{aligned}
\end{equation}
and we directly take the divergence of $\nabla\psi^\pm\p x$ and use \cref{eq:psiimportant5} to obtain \cref{eq:psipm(Omega)}. \cref{eq:psi(Omega)} is a consequence of \cref{eq:psi,eq:psipm(Omega)}.\\
Finally, we obtain from \cref{eq:psiimportant2} and 
\cref{eq:psiimportant3} that
\begin{equation}
\label{eq:psiimportant6}
\begin{aligned}
T\p\eta\Omega&=\p{\p{\textnormal{BR}\p{\gamma}\omega}_h+\p{\textnormal{BR}\p{\gamma}\omega}_3\nabla\eta}^\perp\\&=\p{\frac{V^+ \pare{x}+V^-\p x}2 + \frac{v^+_3 \pare{x}+v^-_3\p{x}}2\nabla\eta\pare{x}}^\perp.
\end{aligned}
\end{equation}
Thus we obtain \cref{eq:Tdivergencefree} from \cref{eq:psiimportant,eq:psiimportant6}.
\end{proof}

}}

\paragraph{The normal velocity in Hamiltonian coordinates.} We can now rewrite \cref{eq:bmvtilde} in the Hamiltonian variable $\psi$, indeed we have that
\begin{align}
    \Omega = & \ S\pare{\Atw\middle|\eta}\psi,\quad
    \omega_3 = \ \nabla\eta\cdot\Omega = \nabla\eta \cdot S\pare{\Atw\middle|\eta}\psi
\end{align}
so that
\begin{equation}\label{eq:bmv}
\begin{aligned}
    {\bm v}\bra{\psi}\pare{\Atw\middle|\eta} \defeq &\vect{V\bra{\psi}\pare{\Atw\middle|\eta}}{v_3\bra{\psi}\pare{\Atw\middle|\eta}} 
    \defeq
    {\bm v}^\star\bra{S\pare{\Atw\middle|\eta}\psi}\pare{\eta}\\
    = 
    &-
    \left. \BR{{\bm \gamma}}\vect{S\pare{\Atw\middle|\eta}\psi}{\nabla\eta \cdot S\pare{\Atw\middle|\eta}\psi} \right|_{{\bm \gamma}\pare{x}=\pare{x, \eta\pare{x}}} 
    \\
    =
    &\vect{\mathpzc{B}_0\pare{\eta}^\perp \bra{\nabla\eta \cdot S\pare{\Atw\middle|\eta}\psi} - \mathpzc{q}\pare{\eta}\pare{S\pare{\Atw\middle|\eta}\psi}^\perp}{-\mathpzc{B}_0\pare{\eta}^\perp \cdot S\pare{\Atw\middle|\eta}\psi}.
\end{aligned}
\end{equation}
The Hamiltonian structure of the canonical CDE is proved in
\cref{proprecast1}.

\subsection{Proof of \cref{proprecast1}}\label{secpsi}

The proof combines the CDE in vorticity variables with the canonical
background shift.

\begin{proof}
By \cref{lem:evolution_eta}, the elevation equation is
\[
    \eta_t
    =
    \mathpzc g^\star\pare{\eta}\Omega .
\]
By \cref{eqpsifinal}, the
potential equation is
\[
    \psi_t
    =
    \pi\pare{\rho^\pm\middle|\eta,\psi}.
\]
We fix once and for all the Bernoulli gauge. Recall that the equation
for the weighted trace is determined only up to an additive function of time,
namely
\[
\psi_t
=
\pi\pare{\rho^\pm\middle|\eta,\psi}
+ c(t).
\]
In the background formulation we choose this gauge by imposing
$c(t)\equiv -c_{\background}$. Equivalently, one may replace the
weighted trace by $\psi-c_{\background}t$. With this convention,
setting $\psi=\bar\psi+\tilde\psi$, and using that $\bar\psi$ is
time-independent, we obtain
\[
\tilde\psi_t
=
\pi\pare{\rho^\pm\middle|\eta,\bar\psi+\tilde\psi}
-
c_{\background}.
\]
By \cref{eq:Omega_can_decomposition}, \cref{eq:pi_as_pistar_composed_with_S}
and the definition of $\tilde\pi$, this gives
\[
\tilde\psi_t
=
\tilde\pi\pare{\rho^\pm,\background\middle|\eta,\tilde\psi}.
\]
Now set $\psi=\bar{\psi}+\tilde{\psi}$. Since $\bar{\psi}$ is independent of
time, $\psi_t=\tilde{\psi}_t$. Moreover, by \cref{eq:Omega_bg_can_def},
\begin{equation}
\label{eq:Omega_can_decomposition}
    \Omega
    =
    \tilde{\Omega}\bra{\background}\pare{\eta}
    +
    S\pare{\Atw\middle|\eta}\tilde{\psi}.
\end{equation}
Therefore
\[
    \eta_t
    =
    \mathpzc g^\star\pare{\eta}
    \bra{\tilde{\Omega}\bra{\background}\pare{\eta}}
    +
    \mathpzc g^\star\pare{\eta}
    \bra{S\pare{\Atw\middle|\eta}\tilde{\psi}}
    =
    \mathpzc g^\star\pare{\eta}
    \bra{\tilde{\Omega}\bra{\background}\pare{\eta}}
    +
    \mathpzc g\pare{\Atw\middle|\eta}
    \bra{\tilde{\psi}},
\]
which is \cref{eq:CDEeta}. On the other hand,
\begin{equation}
\label{eq:canonical_shift_differential}
    \tilde{\psi}_t
    =
    \pi\pare{\rho^\pm\middle|\eta,\bar{\psi}+\tilde{\psi}}
    =
    \tilde{\pi}
    \pare{\rho^\pm,\background\middle|\eta,\tilde{\psi}},
\end{equation}
after subtracting the constant gauge $c_{\background}$ as in
\cref{eq:tildepsiNEW}. This proves \cref{eq:evolution_psi}.

It remains only to justify the homogeneous decomposition of $\tilde{\pi}$.
Using \cref{eq:pi_as_pistar_composed_with_S},
\cref{eq:Omega_can_decomposition}, and the bilinearity of
$\pi^\star_2$, we get
\begin{align}
\pi^\star
\pare{
    \rho^\pm
    \middle|
    \eta,
    \tilde{\Omega}\bra{\background}\pare{\eta}
    +
    S\pare{\Atw\middle|\eta}\tilde{\psi}
}
-
c_{\background} &=
\pi_0\pare{\rho^\pm\middle|\eta}
+
\pi^\star_2\pare{\Atw\middle|\eta}
\bra{
    \tilde{\Omega}\bra{\background}\pare{\eta},
    \tilde{\Omega}\bra{\background}\pare{\eta}
}
-
c_{\background}
\\
&\qquad
+
2\pi^\star_2\pare{\Atw\middle|\eta}
\bra{
    \tilde{\Omega}\bra{\background}\pare{\eta},
    S\pare{\Atw\middle|\eta}\tilde{\psi}
}
\\
&\qquad
+
\pi^\star_2\pare{\Atw\middle|\eta}
\bra{
    S\pare{\Atw\middle|\eta}\tilde{\psi},
    S\pare{\Atw\middle|\eta}\tilde{\psi}
}.
\end{align}
The three lines are exactly
\cref{eq:pi0tilde}, \cref{eq:pi1tilde}, and \cref{eq:pi2tilde}. This
concludes the proof.
\end{proof}

\subsection{Canonical background variables}
\label{sec:background}

The canonical background variables used in \cref{proprecast1} are introduced
there because they are part of the final CDE formulation. In this subsection
we record the elementary identities which justify the canonical rewriting.

The affine background potential $\bar{\psi}$ and the canonical perturbation
$\tilde{\psi}$ are defined in
\cref{eq:psibar_expression,eq:psi_can_def}. Since $\bar{\psi}$ is independent
of $\eta$, \cref{eq:canonical_shift_differential} holds, and therefore the
change of variables preserves the canonical symplectic form.

By \cref{link:omg-psi}, applied with $\psi=\bar{\psi}$, the vorticity
carried by the slice $\set{\tilde{\psi}=0}$ is
\[
    \tilde{\Omega}\bra{\background}\pare{\eta}
    =
    S\pare{\Atw\middle|\eta}\bar{\psi}.
\]
Consequently, inserting $\psi=\bar{\psi}+\tilde{\psi}$ into
\cref{link:omg-psi} gives the decomposition
\cref{eq:Omega_can_decomposition}. Finally, the gauge
$c_{\background}$ in \cref{eq:c_background} removes the constant value of
the Bernoulli function at the flat state.

\begin{prop}[Canonical background formulation]\label{prop:nearCvortrevisited}
Assume the hypotheses of \cref{proprecast1}. Let
$\bar{\psi}$, $\tilde{\psi}$, $\tilde{\Omega}\bra{\background}\pare{\eta}$,
and $c_{\background}$ be defined by
\cref{eq:psibar_expression,eq:psi_can_def,eq:Omega_bg_can_def,eq:c_background}.
Then the change of variables
\[
    \psi=\bar{\psi}+\tilde{\psi}
\]
preserves the canonical two-form in the sense of
\cref{eq:canonical_shift_differential}. In these variables, the canonical
CDE and the homogeneous decomposition of its shifted Bernoulli function are
those in \cref{eq:CDE,eq:tildepi_decomposition}. In particular, for fixed
$\eta$ and $\background$, the term $\tilde{\pi}_j$ is $j$-homogeneous in
$\tilde{\psi}$ for $j=0,1,2$.
\end{prop}

\begin{proof}
The shift $\psi\mapsto\tilde{\psi}=\psi-\bar{\psi}$ is a translation by the
$\eta$-independent function $\bar{\psi}$, hence it preserves the canonical
two-form by \cref{eq:canonical_shift_differential}. By
\cref{eq:Omega_can_decomposition,eq:gstar,eq:DN}, substituting
$\psi=\bar{\psi}+\tilde{\psi}$ gives the first equation in \cref{eq:CDE}.
Since $\bar{\psi}$ is time-independent, $\psi_t=\tilde{\psi}_t$, and the
second equation follows from \cref{eq:tildepsiNEW}. Finally,
\cref{eq:tildepi_decomposition} gives the asserted homogeneous
decomposition.
\end{proof}

\begin{rem}[The equilibrium $\pare{0,0}$]
\label{rem:flat_equilibrium_canonical}
At $\pare{\eta,\tilde{\psi}}=\pare{0,0}$ one has 
\begin{equation}\label{tot-can-vort_at_zero}
    \tilde{\Omega}\bra{\background}\pare{0}=\background.
\end{equation} by \cref{eq:Omega_bg_can_flat}, $\mathpzc g\pare{\Atw\middle|0}\bar{\psi}=0$, and $\pi\pare{\rho^\pm\middle|0,\bar{\psi}}-c_{\background}=0$ by \cref{eq:c_background}. The flat configuration is therefore a stationary point of \cref{eq:CDE}.
\fine
\end{rem}

\section{Linearization at the trivial state}
\label{sec:linearization-trivial-state}

The purpose of this section is to linearize the canonical CDE \cref{eq:CDE}
at the equilibrium state $\pare{\eta,\tilde\psi} = \pare{0,0}$ and to analyze
the resulting dispersion relation. The two main results are
\cref{prop:linearization_trivial_state}, which gives the explicit linearized
matrix operator together with its eigenvalues, and
\cref{prop:purely_imaginary_linear_frequencies}, which characterizes the
purely imaginary part of the spectrum in terms of the stability scales
introduced in \cref{def:linear_stability_scales}.

\begin{prop}[Linearization of the canonical CDE at the trivial state]
\label{prop:linearization_trivial_state}
Let $\sigma\geq0$, $\rho^\pm\geq0$ with $\aangle{\rho^\pm}>0$,
$\mathfrak g\in\bR$, and $\background\in\bR^2$. The linearization of the
canonical CDE \cref{eq:CDE} at $\pare{\eta,\tilde\psi}=\pare{0,0}$ in the
direction $\pare{h,\varphi}$ is
\begin{equation}\label{eq:linearized_CDE_trivial_state}
    \partial_t
    \vect{h}{\varphi}
    =
    {\bm L}_0\pare{D}
    \vect{h}{\varphi},
\end{equation}
where
\begin{equation}\label{eq:linearized_CDE_matrix_operator}
    {\bm L}_0\pare{D}
    \defeq
    \frac{\Atw}{4}
    \background^{\perp}\cdot\nabla\,\Id_{\bC^2}
    +
    \begin{bmatrix}
        0
        &
        \tfrac{1}{2}\av{D}
        \\
        \mathfrak q\pare{\Atw,\background\middle|D}
        &
        0
    \end{bmatrix}.
\end{equation}
Here
\begin{equation}\label{eq:def_linearized_q_background}
    \mathfrak q\pare{\Atw,\background\middle|\xi}
    \defeq
    -\mathfrak g\Atw
    -
    \frac{\sigma}{\aangle{\rho^\pm}}\av{\xi}^2
    +
    \frac{4-\Atw^2}{8}
    \frac{
        \pare{\background^\perp\cdot\xi}^2
    }{
        \av{\xi}
    }.
\end{equation}
Equivalently, in Fourier variables,
\begin{equation}\label{eq:linearized_CDE_matrix_symbol}
    \partial_t
    \vect{\widehat h}{\widehat\varphi}
    =
    {\bm L}_0\pare{\xi}
    \vect{\widehat h}{\widehat\varphi},
    \qquad
    {\bm L}_0\pare{\xi}
    \defeq
    \ii\frac{\Atw}{4}
    \pare{\background^{\perp}\cdot\xi}
    \Id_{\bC^2}
    +
    \begin{bmatrix}
        0
        &
        \tfrac{1}{2}\av{\xi}
        \\
        \mathfrak q\pare{\Atw,\background\middle|\xi}
        &
        0
    \end{bmatrix}.
\end{equation}
The characteristic polynomial is
\begin{equation}\label{eq:linearized_characteristic_polynomial}
    \det\pare{
        \omega\,\Id_{\bC^2}
        -
        {\bm L}_0\pare{\xi}
    }
    =
    \pare{
        \omega
        -
        \ii\frac{\Atw}{4}
        \pare{\background^{\perp}\cdot\xi}
    }^2
    -
    \mathfrak d\pare{\Atw,\background\middle|\xi},
\end{equation}
where
\begin{equation}\label{eq:def_linearized_discriminant}
    \mathfrak d\pare{\Atw,\background\middle|\xi}
    \defeq
    \frac{\av{\xi}}{2}\,
    \mathfrak q\pare{\Atw,\background\middle|\xi}
    =
    \frac{4-\Atw^2}{16}
    \pare{\background^\perp\cdot\xi}^2
    -
    \frac{\mathfrak g\Atw}{2}\av{\xi}
    -
    \frac{\sigma}{2\aangle{\rho^\pm}}\av{\xi}^3.
\end{equation}
Consequently, the two eigenvalues are
\begin{equation}\label{eq:linearized_eigenvalues}
    \lambda_\pm\pare{\Atw,\background\,\middle|\,\xi}
    \defeq
    \ii\frac{\Atw}{4}
    \pare{\background^{\perp}\cdot\xi}
    \pm
    \sqrt{
        \mathfrak d\pare{\Atw,\background\middle|\xi}
    }.
\end{equation}
\end{prop}

\begin{proof}
We linearize the two equations in \cref{eq:CDE} separately. The argument is
organized as follows. We first record the flat-state values of the operators
$\mathpzc B_0$, $\mathpzc q$, $\tau_0$, $S$, and their first variations,
which feed into every subsequent computation. We then linearize the canonical
background vorticity $\tilde\Omega\bra{\background}\pare{\eta}$, which couples
the two equations. With this input we linearize successively the
$\eta$-equation and the $\tilde\psi$-equation, and finally we assemble the
matrix form and read off the spectrum.

The vector-valued Riesz transform $\mathfrak R$ used throughout has principal
symbol $-\ii\xi/\av{\xi}$ and is defined as a singular integral operator in
\cref{eq:Riesz}.

\begin{proofpart}[Flat identities]
Combining \cref{eq:B0_q_flat_linearization} and \cref{eq:bmvtilde}, we get
\begin{equation}\label{eq:vtilde_flat_linearization}
    {\bm v}^\star\bra{\Omega}\pare{0}
    =\vect{ V^{\star}\bra{\Omega}\pare{0}}{ v_3^{\star}\bra{\Omega}\pare{0}}=
    \vect{
        0_{\bR^2}
    }{
        -\tfrac{1}{2}\mathfrak R^\perp\cdot\Omega
    },
    \qquad
    \mathpzc g^\star\pare{0}\bra{\Omega}
    =
    \frac{1}{2}\mathfrak R^\perp\cdot\Omega.
\end{equation}
Moreover, \cref{eq:B0_q_flat_linearization} and \cref{eq:Tstar} gives
\begin{equation}\label{eq:tau0_flat_linearization}
    \tau_0\pare{0} = 0,
    \qquad
    S\pare{\Atw\,\middle|\,0} = \nabla^\perp.
\end{equation}
Consequently, by virtue of \cref{eq:DN}, \cref{eq:vtilde_flat_linearization} and \cref{eq:tau0_flat_linearization}, we find
\begin{equation}\label{eq:DN_flat_linearization}
    \mathpzc g\pare{\Atw\,\middle|\,0}\bra{\varphi}
    =
    \mathpzc g^\star\pare{0}\bra{\nabla^\perp\varphi}
    =
    \frac{1}{2}\mathfrak R^\perp\cdot\nabla^\perp\varphi
    =
    \frac{1}{2}\av{D}\varphi,
\end{equation}
the last equality following from the symbolic identity
$\mathfrak R^\perp\cdot\nabla^\perp = \av{D}$. We shall also need the first variations. Differentiating \cref{eq:B0} and \cref{eq:Q}, we find
\begin{equation}\label{eq:first_variations_B0_q_flat}
    d_\eta\mathpzc B_0\pare{0}\bra{h}\background = 0\qquad\textnormal{and}
    \qquad
    d_\eta\mathpzc q\pare{0}\bra{h}\background
    =\frac{\background}{4\pi}\pv\int_{\mathbb{R}^2}\frac{\delta_zh}{\av{z}^3}\dd z=
    \frac{1}{2}\av{D}h\,\background.
\end{equation}
The last identity of \cref{eq:first_variations_B0_q_flat} is from
\begin{equation}
\label{eq:fraclap}
\mathscr{F}_{x\to\xi}\left( \text{p.v.} \int_{\mathbb{R}^2} \frac{h(x) - h(x-z)}{|z|^3} \dd z \right)(\xi) = 2\pi |\xi| \hat{h}(\xi).
\end{equation}

\end{proofpart}

\begin{proofpart}[Linearization of the canonical background vorticity]
Recall the definition \cref{eq:Omega_bg_can_def},
\begin{equation}
    \tilde\Omega\bra{\background}\pare{\eta}
    =
    S\pare{\Atw\,\middle|\,\eta}\bar\psi.
\end{equation}
Since $S\pare{\Atw\,\middle|\,\eta}
    = \nabla^\perp\pare{1+\Atw\tau_0\pare{\eta}}^{-1}$ and
$\tau_0\pare{0}=0$, we obtain
\begin{equation}\label{eq:dOmega_bg_start}
    d_\eta\tilde\Omega\bra{\background}\pare{0}\bra{h}
    =
    -\Atw\nabla^\perp
    \pare{
        d_\eta\tau_0\pare{0}\bra{h}\bar\psi
    }.
\end{equation}
We compute the scalar factor in \cref{eq:dOmega_bg_start}. By
\cref{eq:Tstar}, \cref{eq:B0_q_flat_linearization}, and
\cref{eq:first_variations_B0_q_flat},
\begin{equation}\label{eq:dtau0_background_flat}
\begin{aligned}
    d_\eta\tau_0\pare{0}\bra{h}\bar\psi
    &=
    \Delta^{-1}\nabla^\perp\cdot
    \pare{
        \mathpzc B_0\pare{0}
        \bra{\nabla h\cdot\background}
        -
        d_\eta\mathpzc q\pare{0}\bra{h}\background
    }
    \\
    &=
    \Delta^{-1}\nabla^\perp\cdot
    \pare{
        \frac{1}{2}\mathfrak R\bra{\nabla h\cdot\background}
        -
        \frac{1}{2}\av{D}h\,\background
    }
    \\
    &=
    -\frac{1}{2}
    \background\cdot\mathfrak R^\perp h.
\end{aligned}
\end{equation}
Indeed, the first term in the second line is annihilated by
$\nabla^\perp\cdot\mathfrak R = 0$, while the second is rewritten using
$\Delta^{-1}\nabla^\perp\cdot\pare{\av{D}h\,\background}
    = \background\cdot\mathfrak R^\perp h$. Therefore
\begin{equation}\label{eq:dOmega_bg_flat}
    d_\eta\tilde\Omega\bra{\background}\pare{0}\bra{h}
    =
    \frac{\Atw}{2}
    \nabla^\perp
    \pare{
        \background\cdot\mathfrak R^\perp h
    }.
\end{equation}
Equivalently, in Fourier variables,
\begin{equation}\label{eq:dOmega_bg_flat_symbol}
    \mathscr F_{x\to\xi}\pare{
        d_\eta\tilde\Omega\bra{\background}\pare{0}\bra{h}
    }\pare{\xi}
    =
    \frac{\Atw}{2}
    \frac{
        \pare{\background\cdot\xi^\perp}
    }{
        \av{\xi}
    }
    \xi^\perp\,
    \widehat h\pare{\xi}.
\end{equation}
\end{proofpart}

\begin{proofpart}[Linearization of the $\eta$-equation]
The first equation in \cref{eq:CDE} reads
\begin{equation}
    \eta_t
    =
    \mathpzc g^\star\pare{\eta}
    \bra{\tilde\Omega\bra{\background}\pare{\eta}}
    +
    \mathpzc g\pare{\Atw\,\middle|\,\eta}\bra{\tilde\psi}.
\end{equation}
Its linearization in the direction $\pare{h,\varphi}$ is
\begin{equation}\label{eq:linearized_eta_intermediate}
\begin{aligned}
    h_t
    &=
    d_\eta
    \pare{
        \mathpzc g^\star\pare{\eta}
        \bra{\tilde\Omega\bra{\background}\pare{\eta}}
    }_{\eta=0}\!\!\bra{h}
    +
    \mathpzc g\pare{\Atw\,\middle|\,0}\bra{\varphi}
    \\
    &=
    d_\eta\mathpzc g^\star\pare{0}\bra{h}\bra{\background}
    +
    \mathpzc g^\star\pare{0}
    \bra{
        d_\eta\tilde\Omega\bra{\background}\pare{0}\bra{h}
    }
    +
    \mathpzc g\pare{\Atw\,\middle|\,0}\bra{\varphi}.
\end{aligned}
\end{equation}
Let us deal with the first contribution. Recall from \cref{eq:gstar} and \cref{eq:bmvtilde} that
$$\mathpzc g^\star\pare{\eta}\bra{\Omega}
    =  V^{\star}\bra{\Omega}\pare{\eta}\cdot\nabla\eta
    -  v_3^{\star}\bra{\Omega}\pare{\eta}.$$
    Now, according to \cref{eq:vtilde_flat_linearization}, \cref{eq:bmvtilde} and \cref{eq:first_variations_B0_q_flat}, we get
    \begin{equation}\label{tildeV and dv3 at zero}
         V^{\star}\bra{\background}\pare{0} = 0\qquad\textnormal{and}\qquad d_\eta v_3^{\star}\bra{\background}\pare{0}\bra{h} = 0.
    \end{equation}
Therefore, the first contribution vanishes, namely
\begin{equation}\label{eq:deta_gstar_fixed_background_zero}
    d_\eta\mathpzc g^\star\pare{0}\bra{h}\bra{\background} = 0.
\end{equation}
Using \cref{eq:vtilde_flat_linearization} and \cref{eq:dOmega_bg_flat} and the fact that
$$\mathfrak R^\perp\cdot\nabla^\perp = \av{D},\qquad\background\cdot\mathfrak{R}^{\perp}=-\background^{\perp}\cdot\mathfrak{R},$$
we
obtain
\begin{equation}\label{eq:gstar_dOmega_bg_flat}
\begin{aligned}
    \mathpzc g^\star\pare{0}
    \bra{
        d_\eta\tilde\Omega\bra{\background}\pare{0}\bra{h}
    }
    &=
    \frac{1}{2}
    \mathfrak R^\perp\cdot
    \bra{
        \frac{\Atw}{2}
        \nabla^\perp
        \pare{
            \background\cdot\mathfrak R^\perp h
        }
    }
    \\
    &=
    \frac{\Atw}{4}
    \av{D}
    \pare{
        \background\cdot\mathfrak R^\perp h
    }\\
    &=
    \frac{\Atw}{4}
    \background^{\perp}\cdot\nabla h.
\end{aligned}
\end{equation}
Inserting \cref{eq:deta_gstar_fixed_background_zero},
\cref{eq:gstar_dOmega_bg_flat}, and \cref{eq:DN_flat_linearization} into
\cref{eq:linearized_eta_intermediate} yields
\begin{equation}\label{eq:etahat}
    h_t
    =
    \frac{\Atw}{4}
    \background^{\perp}\cdot\nabla h
    +
    \frac{1}{2}\av{D}\varphi,
\end{equation}
which is the first component of \cref{eq:linearized_CDE_trivial_state}.
\end{proofpart}

\begin{proofpart}[Linearization of the $\tilde\psi$-equation]
The second equation in \cref{eq:CDE} reads
\begin{equation}
    \tilde\psi_t
    =
    \tilde\pi
    \pare{\rho^\pm,\background\,\middle|\,\eta,\tilde\psi}.
\end{equation}
By \cref{eq:tildepsiNEW} and \cref{eq:pistar}, this can be rewritten as
\begin{equation}\label{eq:tildepi_linearization_start}
    \tilde\pi
    \pare{\rho^\pm,\background\,\middle|\,\eta,\tilde\psi}
    =
    \pi_0\pare{\rho^\pm\,\middle|\,\eta}
    +
    \pi^\star_2\pare{\Atw\,\middle|\,\eta}
    \bra{
        \Omega_{\textnormal{\tiny{Tot}}},
        \Omega_{\textnormal{\tiny{Tot}}}
    }
    -
    c_{\background},
\end{equation}
where the canonical total vorticity is
\begin{equation}
    \Omega_{\textnormal{\tiny{Tot}}}\defeq\Omega_{\textnormal{\tiny{Tot}}}\pare{\eta,\tilde{\psi}}
    \defeq
    \tilde\Omega\bra{\background}\pare{\eta}
    +
    S\pare{\Atw\,\middle|\,\eta}\tilde\psi.
\end{equation}
Recall from \cref{tot-can-vort_at_zero} that
$\Omega_{\textnormal{\tiny{Tot}}}\pare{0,0}=\background$, we deduce that
\begin{align}
    d_{\pare{\eta,\tilde{\psi}}}\tilde\pi
    \pare{\rho^\pm,\background\,\middle|\,0,0}\bra{h,\varphi}
    =
    d_{\eta}\pi_0\pare{\rho^\pm\,\middle|\,\eta}
    \bra{h}+
    \pare{d_{\eta}\pi^\star_2\pare{\Atw\,\middle|0}}\bra{h}
    \bra{
        \background,
        \background
    }+2\pi_2^{\star}\pare{\Atw\middle|0}\bra{\background,d_{\pare{\eta,\tilde\psi}}
    \Omega_{\textnormal{\tiny{Tot}}}\pare{0,0}\bra{h,\varphi}}.
\end{align}
The gravity-capillary contribution is directly obtained from the definition of
$\pi_0$ in \cref{eq:pi0},
\begin{equation}\label{eq:dpi0_flat_direct}
    d_\eta\pi_0\pare{\rho^\pm\,\middle|\,0}\bra{h}
    =
    \frac{\sigma}{\aangle{\rho^\pm}}\Delta h
    -
    \mathfrak g\Atw\,h.
\end{equation}
The first variation of the canonical total vorticity at $\pare{\eta,\tilde\psi}=\pare{0,0}$ is, by virtue of \cref{eq:dOmega_bg_flat},
\begin{equation}\label{eq:dOmega_can_flat}
    \begin{aligned}
        d_{\pare{\eta,\tilde\psi}}
    \Omega_{\textnormal{\tiny{Tot}}}\pare{0,0}\bra{h,\varphi}
    &=
    d_\eta\tilde\Omega\bra{\background}\pare{0}\bra{h}
    +
    \nabla^\perp\varphi\\
    &=\frac{\Atw}{2}
    \nabla^\perp
    \pare{
        \background\cdot\mathfrak R^\perp h
    }+
    \nabla^\perp\varphi.
    \end{aligned}
\end{equation}
Since $\background\in\mathbb R^2$ , one gets from \cref{eq:vtilde_flat_linearization} that

\begin{equation}\label{v3 at zero}
     v_3^{\star}\bra{\background}\pare{0} = 0.
\end{equation} Therefore, coming back to the definition of $\pi_2^{\star}$ in \cref{eq:pistar2}-\cref{eq:pi2_star_I_def_lemma}-\cref{eq:pi2_star_II_def_lemma} and using \cref{tildeV and dv3 at zero}, we get for any $\Omega$,
\begin{equation}\label{eq:pi2star_flat_one_background}
    \pi^\star_2\pare{\Atw\,\middle|\,0}\bra{\background,\Omega}
    =
    -\frac{\Atw}{8}
    \background\cdot\Omega.
\end{equation}
We now substitute the two pieces of \cref{eq:dOmega_can_flat} into
\cref{eq:pi2star_flat_one_background}. For the contribution
$\Omega = \nabla^\perp\varphi$ we obtain, using
$\background\cdot\xi^\perp = -\background^\perp\cdot\xi$,
\begin{equation}\label{eq:pi2star_phi_transport}
    2\pi^\star_2\pare{\Atw\,\middle|\,0}
    \bra{\background,\nabla^\perp\varphi}
    =
    -\frac{\Atw}{4}
    \background\cdot\nabla^\perp\varphi
    =
    \frac{\Atw}{4}
    \background^{\perp}\cdot\nabla\varphi.
\end{equation}
For the contribution
$\Omega = d_\eta\tilde\Omega\bra{\background}\pare{0}\bra{h}$ we obtain, by
\cref{eq:dOmega_bg_flat},
\begin{equation}\label{eq:pi2star_dOmega_bg_contribution}
\begin{aligned}
    2\pi^\star_2\pare{\Atw\,\middle|\,0}
    \bra{
        \background,
        d_\eta\tilde\Omega\bra{\background}\pare{0}\bra{h}
    }
    &=
    -\frac{\Atw}{4}
    \background\cdot
    d_\eta\tilde\Omega\bra{\background}\pare{0}\bra{h}
    \\
    &=
    -\frac{\Atw^2}{8}
    \background\cdot\nabla^\perp
    \pare{
        \background\cdot\mathfrak R^\perp h
    },
\end{aligned}
\end{equation}
which in Fourier variables reads
\begin{equation}\label{eq:pi2star_dOmega_bg_symbol}
   \mathscr F_{x\to\xi}\pare{ 2\pi^\star_2\pare{\Atw\,\middle|\,0}
    \bra{
        \background,
        d_\eta\tilde\Omega\bra{\background}\pare{0}\bra{h}
    }}\p\xi
    =
    -\frac{\Atw^2}{8}
    \frac{\pare{\background^\perp\cdot\xi}^2}{\av{\xi}}
    \widehat h\p\xi.
\end{equation}
It remains to differentiate
$\pi^\star_2\pare{\Atw\,\middle|\,\eta}\bra{\background,\background}$ with
respect to $\eta$ at fixed background vorticity argument. By \cref{eq:pi2_star_II_def_lemma}, \cref{tildeV and dv3 at zero} and \cref{v3 at zero}, the
term $\pi^\star_{2;\RN{2}}$ has zero first variation at $\eta=0$. Now recall from \cref{eq:pi2_star_I_def_lemma} that
\begin{equation}
    \pi^\star_{2;\RN{1}}\pare{\eta}\bra{\background,\background}
    =
    \background^\perp\cdot
     V^{\star}\bra{\background}\pare{\eta}.
\end{equation}
By virtue of \cref{eq:bmvtilde} and \cref{eq:first_variations_B0_q_flat}, we find
\begin{equation}\label{eq:dVtilde_background_flat}
    d_\eta V^{\star}\bra{\background}\pare{0}\bra{h}
    =
    \pare{
        \frac{1}{2}
        \mathfrak R\bra{\nabla h\cdot\background}
        -
        \frac{1}{2}\av{D}h\,\background
    }^\perp.
\end{equation}
Consequently,
\begin{equation}\label{eq:dpi2Istar_background_flat}
\begin{aligned}
    d_\eta\pi^\star_{2;\RN{1}}\pare{0}\bra{h}
    \bra{\background,\background}
    &=
    \background^\perp\cdot
    d_\eta V^{\star}\bra{\background}\pare{0}\bra{h}
    \\
    &=
    \frac{1}{2}
    \background\cdot
    \mathfrak R\bra{\nabla h\cdot\background}
    -
    \frac{1}{2}\av{\background}^2\av{D}h.
\end{aligned}
\end{equation}
In Fourier variables,
\begin{equation}
   \mathscr F_{x\to\xi}\p{ \background\cdot
    \mathfrak R\bra{\nabla h\cdot\background}}\p\xi
    =\frac{
            \pare{\background\cdot\xi}^2
        }{
            \av{\xi}}\widehat h\p\xi,
\end{equation}
so that, using
$$\av{\background}^2\av{\xi}
    = \av{\xi}^{-1}\bra{\pare{\background\cdot\xi}^2 + \pare{\background\cdot\xi^\perp}^2},$$
we deduce
\begin{equation}\label{eq:dpi2star_fixed_background_flat}
    \mathscr F_{x\to\xi}\pare{d_\eta\pi^\star_2\pare{\Atw\,\middle|\,0}\bra{h}
    \bra{\background,\background}}\p\xi
    =\mathscr F_{x\to\xi}\p{
    -d_\eta\pi^\star_{2;\RN{1}}\pare{0}\bra{h}
    \bra{\background,\background}}\p\xi
    =
    \frac{1}{2}
    \frac{\pare{\background\cdot\xi^\perp}^2
        }{\av{\xi}}\widehat h\p\xi.
\end{equation}
Combining
\cref{eq:dpi0_flat_direct},
\cref{eq:pi2star_phi_transport},
\cref{eq:pi2star_dOmega_bg_symbol}, and
\cref{eq:dpi2star_fixed_background_flat}, we obtain
\begin{equation}\label{eq:tildepsihat}
    \varphi_t
    =
    \frac{\Atw}{4}
    \background^{\perp}\cdot\nabla\varphi
    +
    \frac{\sigma}{\aangle{\rho^\pm}}\Delta h
    -
    \mathfrak g\Atw\,h
    +
    \frac{4-\Atw^2}{8} \mathscr F^{-1}_{\xi \to x}
    \p{\frac{
            \pare{\background\cdot\xi^\perp}^2
        }{\av{\xi}}
    \widehat h\p\xi},
\end{equation}
which is the second component of \cref{eq:linearized_CDE_trivial_state}.
\end{proofpart}

\begin{proofpart}[Matrix form and spectrum]
Combining \cref{eq:etahat} and \cref{eq:tildepsihat}, we obtain
\begin{equation}
    \partial_t
    \vect{h}{\varphi}
    =
    \bra{
        \frac{\Atw}{4}
        \background^{\perp}\cdot\nabla\,\Id_{\bC^2}
        +
        \begin{bmatrix}
            0
            &
            \tfrac{1}{2}\av{D}
            \\
            \mathfrak q\pare{\Atw,\background\middle|D}
            &
            0
        \end{bmatrix}
    }
    \vect{h}{\varphi},
\end{equation}
with $\mathfrak q$ defined in
\cref{eq:def_linearized_q_background}; taking the Fourier transform yields
\cref{eq:linearized_CDE_matrix_symbol}. Since the transport part is scalar,
the characteristic polynomial is
\begin{equation}
    \det\pare{
        \omega\,\Id_{\bC^2}
        -
        {\bm L}_0\pare{\xi}
    }
    =
    \pare{
        \omega
        -
        \ii\frac{\Atw}{4}
        \pare{\background^{\perp}\cdot\xi}
    }^2
    -
    \frac{\av{\xi}}{2}
    \mathfrak q\pare{\Atw,\background\middle|\xi},
\end{equation}
and the identity
\begin{equation}
    \frac{\av{\xi}}{2}
    \mathfrak q\pare{\Atw,\background\middle|\xi}
    =
    \frac{4-\Atw^2}{16}
    \pare{\background^\perp\cdot\xi}^2
    -
    \frac{\mathfrak g\Atw}{2}\av{\xi}
    -
    \frac{\sigma}{2\aangle{\rho^\pm}}\av{\xi}^3
\end{equation}
is precisely \cref{eq:def_linearized_discriminant}. The two roots of the
characteristic polynomial are \cref{eq:linearized_eigenvalues}.
\end{proofpart}
This concludes the proof of \cref{prop:linearization_trivial_state}.
\end{proof}

\begin{rem}[Structural consistency]
\label{rem:linearization_direct_consistency}
The transport
$\tfrac{\Atw}{4}\background^{\perp}\cdot\nabla$ appears in both rows of
${\bm L}_0\pare{D}$ for two structurally distinct reasons: in the
$\eta$-equation it arises from the canonical background vorticity
$\tilde\Omega\bra{\background}\pare{\eta}
    = S\pare{\Atw\,\middle|\,\eta}\bar\psi$, whereas in the
$\tilde\psi$-equation it arises from the bilinear term
$2\pi^\star_2\pare{\Atw\,\middle|\,0}\bra{\background,\nabla^\perp\varphi}$.
That the same transport speed appears in both rows is the direct CDE-level
manifestation of the canonical Hamiltonian structure.
\fine
\end{rem}

\begin{definition}[Linear stability scales]
\label{def:linear_stability_scales}
We introduce the following stability scales
attached to the linear dispersion relation
\cref{eq:def_linearized_discriminant}.
\begin{enumerate}
    \item The \emph{reduced gravity}, or \emph{Rayleigh--Taylor acceleration},
    is
    \begin{equation}\label{eq:reduced_gravity_RT_acceleration}
        \mathsf g_{\Atw}'
        \defeq
        \frac{\mathfrak g\Atw}{2}\cdot
    \end{equation}

    \item The \emph{capillary coefficient} is
    \begin{equation}\label{eq:capillary_coefficient}
        \mathsf c_\sigma
        \defeq
        \frac{\sigma}{2\aangle{\rho^\pm}}
        =
        \frac{\sigma}{\rho^- + \rho^+}\cdot
    \end{equation}

    \item The \emph{directional Kelvin--Helmholtz shear coefficient} is
    \begin{equation}\label{eq:directional_KH_shear_coefficient}
        \mathsf{KH}_{\Atw,\background}\pare{\vartheta}
        \defeq
        \frac{4-\Atw^2}{16}
        \pare{\background\cdot\vartheta^\perp}^2,
        \qquad
        \vartheta\in\bR^2, \quad \av{\vartheta} = 1.
    \end{equation}

    \item If $\sigma>0$, the \emph{directional effective Weber wave number} is
    \begin{equation}\label{eq:directional_effective_Weber_wavenumber}
        \kappa_{\rm We}\pare{\vartheta}
        \defeq
        \frac{
            \mathsf{KH}_{\Atw,\background}\pare{\vartheta}
        }{
            \mathsf c_\sigma
        }
        =
        \frac{
            \rho^- \rho^+
        }{
            \sigma\pare{\rho^- + \rho^+}
        }
        \pare{\background\cdot\vartheta^\perp}^2.
    \end{equation}
    For a prescribed length scale $L>0$, the dimensionless quantity
    $L\,\kappa_{\rm We}\pare{\vartheta}$ is the corresponding effective
    directional Weber number.

    \item If $\sigma>0$, the \emph{signed inverse squared capillary length} is
    \begin{equation}\label{eq:signed_inverse_squared_capillary_length}
        \kappa_{\rm cap}^2
        \defeq
        \frac{
            \mathsf g_{\Atw}'
        }{
            \mathsf c_\sigma
        }
        =
        \frac{
            \pare{\rho^- - \rho^+}\mathfrak g
        }{
            \sigma
        }\cdot
    \end{equation}
    When $\pare{\rho^- - \rho^+}\mathfrak g>0$, this is positive and
    \begin{equation}
        \ell_{\rm cap}
        \defeq
        \frac{1}{\kappa_{\rm cap}}
        =
        \sqrt{
            \frac{
                \sigma
            }{
                \pare{\rho^- - \rho^+}\mathfrak g
            }
        }
    \end{equation}
    is the capillary length; for a prescribed length scale $L>0$, the
    dimensionless quantity $L^2\kappa_{\rm cap}^2$ is the corresponding
    Bond, or Eötvös, number.

    \item The maximal Kelvin--Helmholtz and Weber wavenumbers are
    \begin{equation}\label{eq:maximal_KH_shear_coefficient}
        \mathsf{KH}^{\max}_{\Atw,\background}
        \defeq
        \sup_{\av{\vartheta}=1}
        \mathsf{KH}_{\Atw,\background}\pare{\vartheta}
        =
        \frac{4-\Atw^2}{16}\av{\background}^2
        =
        \frac{\rho^- \rho^+}{\pare{\rho^- + \rho^+}^2}
        \av{\background}^2
    \end{equation}
    and, if $\sigma>0$,
    \begin{equation}\label{eq:maximal_effective_Weber_wavenumber}
        \kappa_{\rm We}^{\max}
        \defeq
        \sup_{\av{\vartheta}=1}
        \kappa_{\rm We}\pare{\vartheta}
        =
        \frac{
            \rho^- \rho^+
        }{
            \sigma\pare{\rho^- + \rho^+}
        }
        \av{\background}^2.
    \end{equation}
\end{enumerate}
\end{definition}

\begin{prop}[Purely imaginary linear frequencies]
\label{prop:purely_imaginary_linear_frequencies}
Assume the hypotheses of \cref{prop:linearization_trivial_state}, and let the
stability scales be those introduced in \cref{def:linear_stability_scales}.
Then, for every $\xi\in\bR^2\setminus\set{0}$, with
$\hat{\xi}\defeq\xi/\av{\xi}$, one has
\begin{equation}\label{eq:linear_discriminant_physical_factorization}
    \mathfrak d\pare{\Atw,\background\middle|\xi}
    =
    \av{\xi}
    \bra{
        \mathsf{KH}_{\Atw,\background}\pare{\hat{\xi}}\av{\xi}
        -
        \mathsf g_{\Atw}'
        -
        \mathsf c_\sigma\av{\xi}^2
    }.
\end{equation}
Consequently,
\begin{equation}\label{eq:purely_imaginary_frequency_condition}
    \lambda_\pm\pare{\Atw,\background\,\middle|\,\xi}\in\ii\bR
    \quad\Longleftrightarrow\quad
    \mathsf c_\sigma\av{\xi}^2
    -
    \mathsf{KH}_{\Atw,\background}\pare{\hat{\xi}}\av{\xi}
    +
    \mathsf g_{\Atw}'
    \geq 0.
\end{equation}
More precisely, the following alternatives hold.
\begin{enumerate}
    \item If $\sigma>0$, then
    \begin{equation}\label{eq:purely_imaginary_sigma_positive_condition}
        \lambda_\pm\pare{\Atw,\background\,\middle|\,\xi}\in\ii\bR
        \quad\Longleftrightarrow\quad
        \av{\xi}^2
        -
        \kappa_{\rm We}\pare{\hat{\xi}}\av{\xi}
        +
        \kappa_{\rm cap}^2
        \geq 0.
    \end{equation}
    In particular, if
    \begin{equation}\label{eq:global_purely_imaginary_sigma_positive}
        \pare{\kappa_{\rm We}^{\max}}^2
        \leq
        4\kappa_{\rm cap}^2,
    \end{equation}
    then both dispersion branches are purely imaginary for every
    $\xi\in\bR^2\setminus\set{0}$.

    Conversely, if there exists $\vartheta\in\bR^2$ with
    $\av{\vartheta}=1$ such that
    \begin{equation}\label{eq:directional_unstable_band_condition_sigma_positive}
        \kappa_{\rm We}\pare{\vartheta}^2
        >
        4\kappa_{\rm cap}^2,
    \end{equation}
    then the frequencies in the direction $\vartheta$ are not purely
    imaginary for
    \begin{equation}\label{eq:unstable_band_sigma_positive}
        \av{\xi}
        \in
        \pare{
            \max\set{0,\kappa_-\pare{\vartheta}},\,
            \kappa_+\pare{\vartheta}
        },
        \qquad
        \hat{\xi} = \vartheta,
    \end{equation}
    where
    \begin{equation}\label{eq:unstable_band_endpoints_sigma_positive}
        \kappa_\pm\pare{\vartheta}
        \defeq
        \frac{
            \kappa_{\rm We}\pare{\vartheta}
            \pm
            \sqrt{
                \kappa_{\rm We}\pare{\vartheta}^2
                -
                4\kappa_{\rm cap}^2
            }
        }{
            2
        }\cdot
    \end{equation}
    Moreover, independently of
    \cref{eq:global_purely_imaginary_sigma_positive}, capillarity stabilizes
    the high-frequency region: both branches are purely imaginary for all
    $\av{\xi}\geq\kappa_{\rm hi}$, where
    \begin{equation}\label{eq:high_frequency_cutoff_sigma_positive}
        \kappa_{\rm hi}
        \defeq
        \begin{cases}
        \displaystyle
        \frac{
            \kappa_{\rm We}^{\max}
            +
            \sqrt{
                \pare{\kappa_{\rm We}^{\max}}^2
                -
                4\kappa_{\rm cap}^2
            }
        }{
            2
        },
        &
        \textnormal{if }
        \pare{\kappa_{\rm We}^{\max}}^2
        >
        4\kappa_{\rm cap}^2,
        \\[1.2em]
        0,
        &
        \textnormal{otherwise.}
        \end{cases}
    \end{equation}

    \item If $\sigma=0$, then
    \begin{equation}\label{eq:purely_imaginary_sigma_zero_condition}
        \lambda_\pm\pare{\Atw,\background\,\middle|\,\xi}\in\ii\bR
        \quad\Longleftrightarrow\quad
        \mathsf{KH}_{\Atw,\background}\pare{\hat{\xi}}\av{\xi}
        \leq
        \mathsf g_{\Atw}'.
    \end{equation}
    Hence, if $\mathsf g_{\Atw}'<0$, no nonzero frequency has purely
    imaginary branches. If $\mathsf g_{\Atw}'=0$, the branches are purely
    imaginary precisely in the directions satisfying
    $\mathsf{KH}_{\Atw,\background}\pare{\hat{\xi}}=0$. If
    $\mathsf g_{\Atw}'>0$, then in every direction satisfying
    $\mathsf{KH}_{\Atw,\background}\pare{\hat{\xi}}>0$ the purely
    imaginary range is exactly
    \begin{equation}\label{eq:low_frequency_range_sigma_zero}
        0<\av{\xi}
        \leq
        \frac{
            \mathsf g_{\Atw}'
        }{
            \mathsf{KH}_{\Atw,\background}\pare{\hat{\xi}}
        }\cdot
    \end{equation}
    In particular, if $\mathsf g_{\Atw}'>0$ and
    $\mathsf{KH}^{\max}_{\Atw,\background}>0$, then both branches are
    purely imaginary for every frequency satisfying
    \begin{equation}\label{eq:uniform_low_frequency_range_sigma_zero}
        0<\av{\xi}
        \leq
        \frac{
            \mathsf g_{\Atw}'
        }{
            \mathsf{KH}^{\max}_{\Atw,\background}
        }\cdot
    \end{equation}
    Full-frequency stability in the capillarity-free case holds if and only
    if
    \begin{equation}\label{eq:global_purely_imaginary_sigma_zero}
        \mathsf g_{\Atw}'\geq 0,
        \qquad
        \mathsf{KH}^{\max}_{\Atw,\background} = 0 \ ;
    \end{equation}
    equivalently,
    \begin{equation}
        \mathfrak g\Atw \geq 0,
        \qquad
        \pare{4-\Atw^2}\av{\background}^2 = 0.
    \end{equation}
\end{enumerate}
\end{prop}

\begin{proof}
The two branches differ from the purely imaginary scalar
$\ii\tfrac{\Atw}{4}\pare{ \background^{\perp}\cdot\xi}$ only by
$\pm\sqrt{\mathfrak d_{\Atw,\background}\pare{\xi}}$, where
$\mathfrak d_{\Atw,\background}\pare{\xi}\in\bR$. Hence both branches are
purely imaginary if and only if
$\mathfrak d_{\Atw,\background}\pare{\xi}\leq 0$. We first rewrite the
discriminant in the dispersion-physics variables of
\cref{def:linear_stability_scales}, then separately analyze the cases
$\sigma>0$ and $\sigma=0$.

Let $\xi\in\bR^2\setminus\set{0}$ and $\hat{\xi}\defeq\xi/\av{\xi}$. Since
$\xi= \av{\xi}\hat{\xi}$,
\begin{equation}
    \pare{\background^\perp\cdot\xi}^2
    =
    \av{\xi}^2
    \pare{\background^\perp\cdot\hat{\xi}}^2.
\end{equation}
Using the definitions of $\mathsf{KH}_{\Atw,\background}$,
$\mathsf g_{\Atw}'$, and $\mathsf c_\sigma$ in
\cref{def:linear_stability_scales}, this gives
\begin{equation}
    \mathfrak d\pare{\Atw,\background\middle|\xi}
    =
    \av{\xi}
    \bra{
        \mathsf{KH}_{\Atw,\background}\pare{\hat{\xi}}\av{\xi}
        -
        \mathsf g_{\Atw}'
        -
        \mathsf c_\sigma\av{\xi}^2
    },
\end{equation}
which is \cref{eq:linear_discriminant_physical_factorization}. Since
$\av{\xi}>0$, the condition
$\mathfrak d_{\Atw,\background}\pare{\xi}\leq 0$ is equivalent to
\cref{eq:purely_imaginary_frequency_condition}.

\begin{proofpart}[The capillary case $\sigma>0$]
If $\sigma>0$, then $\mathsf c_\sigma>0$. Dividing
\cref{eq:purely_imaginary_frequency_condition} by $\mathsf c_\sigma$ and
using
\begin{equation}
    \kappa_{\rm We}\pare{\hat{\xi}}
    =
    \frac{
        \mathsf{KH}_{\Atw,\background}\pare{\hat{\xi}}
    }{
        \mathsf c_\sigma
    },
    \qquad
    \kappa_{\rm cap}^2
    =
    \frac{
        \mathsf g_{\Atw}'
    }{
        \mathsf c_\sigma
    },
\end{equation}
yields \cref{eq:purely_imaginary_sigma_positive_condition}.

Fix a direction $\vartheta\in\bR^2$ with $\av{\vartheta}=1$, and consider the
quadratic polynomial
\begin{equation}
    P_\vartheta\pare{X}
    \defeq
    X^2
    -
    \kappa_{\rm We}\pare{\vartheta}X
    +
    \kappa_{\rm cap}^2,
    \qquad
    X\in\bR,
\end{equation}
with discriminant
$\Delta_\vartheta
    \defeq
    \kappa_{\rm We}^2\pare{\vartheta} - 4\kappa_{\rm cap}^2$. If
$\pare{\kappa_{\rm We}^{\max}}^2 \leq 4\kappa_{\rm cap}^2$, then
$\Delta_\vartheta\leq 0$ for every $\vartheta$. Since $P_\vartheta$ has
positive leading coefficient, $P_\vartheta\pare{X}\geq 0$ for every
$X\in\bR$, hence for every $X>0$. Applying this with $X=\av{\xi}$ and
$\vartheta=\hat{\xi}$ proves that
both dispersion branches are purely imaginary for every $\xi\in\bR^2\setminus\set{0}$.

Conversely, assume that for some direction $\vartheta$,
$\kappa_{\rm We}\pare{\vartheta}^2 > 4\kappa_{\rm cap}^2$. Then $P_\vartheta$
has two real roots $\kappa_\pm\pare{\vartheta}$ given by
\cref{eq:unstable_band_endpoints_sigma_positive}; since the leading
coefficient is positive, $P_\vartheta$ is negative precisely between its two
roots. Intersecting with the positive half-line gives
\begin{equation}
    P_\vartheta\pare{X}<0
    \quad\Longleftrightarrow\quad
    X\in
    \pare{
        \max\set{0,\kappa_-\pare{\vartheta}},\,
        \kappa_+\pare{\vartheta}
    },
\end{equation}
which in view of
\cref{eq:purely_imaginary_sigma_positive_condition} is exactly the unstable
band
\cref{eq:unstable_band_sigma_positive}--\cref{eq:unstable_band_endpoints_sigma_positive}.

It remains to prove the high-frequency statement. Since
$\kappa_{\rm We}\pare{\vartheta}\leq\kappa_{\rm We}^{\max}$ for every
direction, we have, for every $X>0$,
\begin{equation}
    P_\vartheta\pare{X}
    \geq
    X^2 - \kappa_{\rm We}^{\max}X + \kappa_{\rm cap}^2.
\end{equation}
If $\pare{\kappa_{\rm We}^{\max}}^2\leq 4\kappa_{\rm cap}^2$, the right-hand
side is nonnegative for every $X$ and one may take $\kappa_{\rm hi}=0$.
Otherwise, the right-hand side is nonnegative for every $X$ larger than its
largest root, namely for
\begin{equation}
    X
    \geq
    \frac{
        \kappa_{\rm We}^{\max}
        +
        \sqrt{
            \pare{\kappa_{\rm We}^{\max}}^2
            -
            4\kappa_{\rm cap}^2
        }
    }{
        2
    },
\end{equation}
which proves \cref{eq:high_frequency_cutoff_sigma_positive}.
\end{proofpart}

\begin{proofpart}[The capillarity-free case $\sigma=0$]
If $\sigma=0$, then $\mathsf c_\sigma=0$, and
\cref{eq:purely_imaginary_frequency_condition} reduces to
\begin{equation}
    \mathsf{KH}_{\Atw,\background}\pare{\hat{\xi}}\av{\xi}
    \leq
    \mathsf g_{\Atw}',
\end{equation}
which is \cref{eq:purely_imaginary_sigma_zero_condition}.

Since
$\mathsf{KH}_{\Atw,\background}\pare{\hat{\xi}}\geq 0$ and $\av{\xi}>0$,
the preceding inequality has no solution if $\mathsf g_{\Atw}'<0$, proving
the first assertion. If $\mathsf g_{\Atw}'=0$, the inequality becomes
$\mathsf{KH}_{\Atw,\background}\pare{\hat{\xi}}\av{\xi}\leq 0$; since both
factors on the left-hand side are nonnegative and $\av{\xi}>0$, this is
equivalent to $\mathsf{KH}_{\Atw,\background}\pare{\hat{\xi}}=0$, proving
the second assertion.

Assume now $\mathsf g_{\Atw}'>0$. If
$\mathsf{KH}_{\Atw,\background}\pare{\hat{\xi}}>0$, solving
$\mathsf{KH}_{\Atw,\background}\pare{\hat{\xi}}\av{\xi}\leq\mathsf g_{\Atw}'$
for $\av{\xi}$ gives exactly \cref{eq:low_frequency_range_sigma_zero}. If
$\mathsf{KH}_{\Atw,\background}\pare{\hat{\xi}}=0$, the inequality is
satisfied for every $\av{\xi}>0$ in that direction.

Taking the maximum over directions gives the uniform low-frequency
criterion: if $\mathsf{KH}^{\max}_{\Atw,\background}>0$ and $\av{\xi}$
satisfies
$0<\av{\xi}\leq \mathsf g_{\Atw}'/\mathsf{KH}^{\max}_{\Atw,\background}$,
then
\begin{equation}
    \mathsf{KH}_{\Atw,\background}\pare{\hat{\xi}}\av{\xi}
    \leq
    \mathsf{KH}^{\max}_{\Atw,\background}\av{\xi}
    \leq
    \mathsf g_{\Atw}',
\end{equation}
proving \cref{eq:uniform_low_frequency_range_sigma_zero}.

Finally, full-frequency stability in the capillarity-free case means that
$\mathsf{KH}_{\Atw,\background}\pare{\vartheta}X\leq\mathsf g_{\Atw}'$ for
every direction $\vartheta$ and every $X>0$, which holds if and only if
$\mathsf{KH}_{\Atw,\background}\pare{\vartheta}=0$ for every $\vartheta$ and
$\mathsf g_{\Atw}'\geq 0$. By \cref{def:linear_stability_scales}, this is
equivalent to $\pare{4-\Atw^2}\av{\background}^2 = 0$ and
$\mathfrak g\Atw\geq 0$, proving
\cref{eq:global_purely_imaginary_sigma_zero}.
\end{proofpart}
This concludes the proof of \cref{prop:purely_imaginary_linear_frequencies}.
\end{proof}

\begin{rem}[Dimensionless Atwood--Weber threshold]
\label{rem:dimensionless_atwood_weber_threshold}
When $\mathfrak g\Atw>0$, one has $\kappa_{\rm cap}^2>0$, and the
full-frequency sufficient condition
\cref{eq:global_purely_imaginary_sigma_positive} can be recast, with
$\mathsf W_{\Atw,\background}$ defined in
\cref{eq:strict_atwood_weber_number_intro}, as
\begin{equation}\label{eq:dimensionless_atwood_weber_threshold}
    \mathsf W_{\Atw,\background}\leq 1.
\end{equation}
Thus the Atwood number has two distinct stabilizing effects: it generates
the Rayleigh--Taylor restoring term $\kappa_{\rm cap}^2$ when
$\mathfrak g\Atw>0$, and it weakens the Kelvin--Helmholtz contribution
through the factor $4-\Atw^2$. In the Kelvin--Helmholtz endpoint
$\Atw=0$, one has $\kappa_{\rm cap}^2 = 0$, so no nonzero background
$\background$ can satisfy a full-frequency stability condition on $\bR^2$.
In that endpoint one can only obtain a high-frequency stability statement,
which by \cref{eq:high_frequency_cutoff_sigma_positive} reduces to
\begin{equation}
    \av{\xi}
    \geq
    \kappa_{\rm We}^{\max}\big|_{\Atw=0}
    =
    \frac{
        \aangle{\rho^\pm}
    }{
        2\sigma
    }
    \av{\background}^2.
\end{equation}
This is the flat-interface analogue of the Weber threshold in the periodic
Kelvin--Helmholtz problem, but without the discrete low-frequency gap
present in the circular geometry, see \cite{murgante2025long}.
\fine
\end{rem}

\section{Paradifferential toolkit}
\label{sec:preliminaries}
Throughout this section, $\mathtt d\geq1$. In accordance with
\cref{notation:vector_convention}, all definitions and results are stated for
scalar functions, operators, and symbols; the vector- and matrix-valued
versions follow componentwise. Standard spatial, frequency, and integration
variables retain the usual lower-case notation; in particular,
\[
    x,\xi,z,\varsigma,\nu,\nu',\zeta\in\bR^{\mathtt d}.
\]
All multi-indices belong to $\bN^{\mathtt d}$. The symbols
$\nabla_x$, $\nabla_\xi$, $\nabla_x^2$,
$\Delta_x$, and $\Delta_z$ denote the Euclidean gradients, Hessian, and Laplacians in dimension $\mathtt d$.
Moreover, we denote
\begin{align}\label{eq:toolkit-finite-differences-Rd}
    \delta_z f\pare x
    &\defeq
    f\pare x-f\pare{x-z},
    &
    \Delta_z f\pare x
    &\defeq
    \frac{\delta_z f\pare x}{\av z},
    &
    \hat z
    &\defeq
    \frac z{\av z}
    \in\bS^{\mathtt d-1}.
\end{align}
\begin{notation}
\label{notation:abstract-Fourier-convention-Rd}
We use the standard surface measure $\dd\sigma$ on
$\bS^{\mathtt d-1}$; for $\mathtt d=1$, it is counting measure on
$\bS^0=\set{-1,1}$.  The Fourier convention is 
\[
    \widehat f\pare\xi
    =
    \frac1{\pare{2\pi}^{\mathtt d/2}}
    \int_{\bR^{\mathtt d}}e^{-\ii x\cdot\xi}f\pare x\dd x,
    \qquad
    \cF\pare{\mathfrak{m}\pare Df}\pare\xi
    =
    \mathfrak m\pare\xi\widehat f\pare\xi,
\]
with inverse transform
\[
    f\pare x
    =
    \frac1{\pare{2\pi}^{\mathtt d/2}}
    \int_{\bR^{\mathtt d}}
    e^{\ii x\cdot\xi}\widehat f\pare\xi\dd\xi.
\]
\end{notation}
Every functional space below is displayed with its domain and codomain. In
operator spaces, omission of the codomain means that domain and codomain
coincide. For brevity, in norm subscripts we suppress the ambient domain and
codomain. In particular, $H^s$ there stands for
$H^s\pare{\bR^{\mathtt d};\bC}$. Real scalar backgrounds belong to
$H^s\pare{\bR^{\mathtt d};\bR}$, and we use the convention
\begin{align}\label{eq:toolkit-real-Sobolev-ball-Rd}
    \norm u_s
    &\defeq
    \norm{u}_{H^s\pare{\bR^{\mathtt d};\bC}},
    &
    \norm u_{s,I}
    &\defeq
    \sup_{t\in I}
    \norm{u\pare t}_s,
    \\
    \norm{\vec U_p}_{s,I}
    &\defeq
    \max_{1\leq j\leq p}\norm{u_j}_{s,I},
    &
    \BallR{s}{r}
    &\defeq
    B_{\cC\pare{I;H^s\pare{\bR^{\mathtt d};\bR}}}\pare{0,r}.
\end{align}

This section collects the paradifferential formalism that will be used in
\cref{sec:paralinearization} to paralinearize the generalized water-waves
system \cref{eq:CDE}. We first remind the general notion of smoothing
operators and symbols. They depend on an external function and we keep track
of the homogeneity with respect to this latter. Then, we introduce the
Bony--Weyl quantization and recall its main para-composition rules. In a
second part, we develop an approach to paralinearize general singular
integral operators on the whole space. This uses localized and tail kernel
functions to deal with integrability near the origin and at infinity,
respectively.

\subsection{General paradifferential calculus}
In this subsection, we recall the basic definitions and propositions in paradifferential calculus, such as symbols, functions, smoothing operators, and Bony-Weyl quantization.

\paragraph{$p$-linear bounded maps on Fréchet spaces, continuously differentiable and analytic maps.}
\begin{definition}[$p$-linear bounded maps in locally convex vector spaces]\label{def:L_Frechet}
    Let $\Lambda \subset \bR$,  $X$ be a Banach space, $Y$ a locally convex topological vector space
    with topology generated by a family of seminorms $\set{p_\lambda}_{\lambda\in\Lambda}$,
    and $p \in \bN^\star$. We denote by $\cL^p\pare{X;Y}$ the space of continuous
    $p$-linear maps $T : X^p \to Y$, i.e. maps that are linear in each argument
    separately and satisfy
    \begin{align}
        p_\lambda\pare{T\pare{u_1,\ldots,u_p}}
        \leq C_\lambda\, \prod_{\mathsf{j} =1}^p \norm{u_{\mathsf{j}}}_X,
        && \forall\, u_1,\ldots,u_p \in X,\ \lambda \in \Lambda,
    \end{align}
    for constants $C_\lambda > 0$ depending on $\lambda$ but not on $u_1,\ldots,u_p$.
    We write 
    \begin{align}
        \cL\pare{X;Y} \defeq \cL^1\pare{X;Y}, &&  \cL^p\pare{X}\defeq \cL^p\pare{X;X}.
    \end{align}
\end{definition}

We shall also use real-analytic functions with values in a Fréchet space;
we recall the definition for the reader's convenience.

\begin{definition}[$\cC^k$ maps into locally convex spaces]
\label{def:Ck_lcs}
Let $X$ be a Banach space, let $Y$ be a locally convex space whose topology is
generated by a family of seminorms $\set{p_\lambda}_{\lambda\in\Lambda}$,
and let $U\subseteq X$ be open. A map $f:U\to Y$ is of class $\cC^1$ if there
exists a map $\dd f:U\to\cL\pare{X;Y}$ such that, for every $x_0\in U$ and
every seminorm $p_\lambda$,
\begin{equation}
    p_\lambda\pare{f\pare{x_0+h}-f\pare{x_0}-\dd f\pare{x_0}h}
    = o\pare{\norm h_X}
    \qquad \textnormal{as }\norm h_X\to0,
\end{equation}
and such that $x\mapsto\dd f\pare{x}$ is continuous for the topology of
uniform convergence on bounded subsets of $X$, seminorm by seminorm in $Y$.
Inductively, $f\in\cC^k\pare{U;Y}$ if $f\in\cC^1$ and
$\dd f\in\cC^{k-1}\pare{U;\cL\pare{X;Y}}$. We set
$\cC^\infty\pare{U;Y}\defeq\bigcap_{k\ge1}\cC^k\pare{U;Y}$.
\end{definition}

\begin{definition}[Real-analytic maps into locally convex spaces]
\label{def:Comega_lcs}
    Let $X$ be a Banach space, $Y$ a locally convex space with topology
    generated by seminorms $\set{p_\lambda}_{\lambda\in\Lambda}$, and
    $U \subseteq X$ open.
    We say that $f : U \to Y$ is \emph{real-analytic}, and write
    $f \in \cC^\omega\pare{U;Y}$, if for every $x_0 \in U$ there exists
    $r = r\pare{x_0} > 0$ such that $B_X\pare{x_0,r} \subseteq U$ and
    \begin{equation}
        f\pare{x_0 + h}
        = \sum_{k=0}^{+\infty} \frac{1}{k!}\, \dd^k f\pare{x_0}\bra{h}^k,
        \qquad \forall\, h \in B_X\pare{0,r},
    \end{equation}
    where $\dd^k f\pare{x_0} \in \cL^k\pare{X;Y}$ is the $k$-th Fréchet
    differential of $f$ at $x_0$ \textnormal{(}cf.\ \cref{def:Ck_lcs}\textnormal{)},
\[\dd^k f\pare{x_0}\bra{h}^k \defeq \dd^k f\pare{x_0}\bra{h,\ldots,h}\]
    denotes its evaluation on the $k$-fold diagonal, and the series converges
    in $Y$, i.e.\ for every $\lambda \in \Lambda$,
    \begin{equation}\label{eq:analytic_convergence}
        \lim_{N\to+\infty}
        p_\lambda\!\left(
            f\pare{x_0+h} - \sum_{k=0}^{N}
            \frac{1}{k!}\, \dd^k f\pare{x_0}\bra{h}^k
        \right) = 0,
        \qquad \forall\, h \in B_X\pare{0,r}.
    \end{equation}
    We say that the series converges \emph{uniformly} on $B_X\pare{0,r}$ if
    the limit in \cref{eq:analytic_convergence} is uniform in
    $h \in B_X\pare{0,r}$ for every $\lambda$.
\end{definition}

\paragraph{Tame smoothing operators.}
\begin{definition}[Tame smoothing operators with vanishing order]
\label{def:smoothing}
Let $\varrho\in\bR$ and $q\in\bN$. We say that
$r_\varrho\in\mathcal R^\varrho_{\ge q}$ if, for some fixed
$p\in\bN^\star$, there exists $s_0>0$ such that,
for every $s\ge s_0$, there exists a radius $r^{r_\varrho}_s>0$ with the following
property. For every $0<r\le r^{r_\varrho}_s$ and every
\[
    \vec U_p\in
   \pare{\Ball{s_0}{r}}^p,
\]
we have
\[
    r_\varrho\pare{\vec U_p}
    \in
    \cL\pare{
        H^{s-\varrho}\pare{\bR^{\mathtt d};\bC};
        H^s\pare{\bR^{\mathtt d};\bC}
    },
\]
and, for every
$h\in H^{s-\varrho}\pare{\bR^{\mathtt d};\bC}$,
\begin{equation}\label{eq:smoothing-tame-geq-q}
\begin{aligned}
    \norm{r_\varrho\pare{\vec U_p}h}_{H^s}
    \lesssim_{s,\varrho}\;&
    \norm{\vec U_p}_{s_0}^{q}
    \norm{h}_{H^{s-\varrho}}
 .
\end{aligned}
\end{equation}

We call $r_\varrho$ real-to-real if $r_\varrho(\vec U_p)h$ is real-valued whenever the
backgrounds and $h$ are real-valued.
We use the convention $\mathcal R^\varrho
    \defeq
    \mathcal R^\varrho_{\ge1}$.

\end{definition}

\paragraph{Fourier multipliers, symbols and functions.}

\begin{definition}[Fourier multiplier]\label{def:FM}
    Let $m\in \bR$, we denote with $\tilde \Gamma ^ m$ the space of maps 
    \begin{equation}
        \xi \mapsto \mathpzc{m}\pare{\xi} \in \cC^\infty\pare{\bR^{\mathtt d}\setminus \{0\}; \bC}, 
    \end{equation}
    such that for any $\beta \in \bN^{\mathtt d}$ and $\av{\xi}\geq 1/2$ we have that
    \begin{equation}
        \av{\partial_\xi^\beta \mathpzc{m}\pare{\xi}}\lesssim_\beta \av{\xi}^{m-\av{\beta}}. 
    \end{equation}
    
\end{definition}

\begin{definition}[Symbols with finite regularity and order of vanishing]
\label{def:sfr}
Let $m\in\bR$ and let $q\in\bN^\star$.

\begin{enumerate}
\item We denote by $\Gamma^m_{\ge q}$ the space of maps
\[
    a:
    \pare{\Ball{s_0}{r_0}}^p
    \to
    L^\infty_{\loc}\pare{
        \bR^{\mathtt d}\times
        \pare{\bR^{\mathtt d}\setminus\set{0}};\bC
    }
\]
\[
    \vec U_p\defeq\pare{u_1,\ldots,u_p}^\intercal
    \mapsto
    \set{\pare{x,\xi}\mapsto a\pare{\vec U_p;x,\xi}}
\]
defined for some $s_0>\mathtt d/2$, $r_0>0$, $p\in\bN^\star$, and
satisfying the following properties.

\begin{enumerate}
\item 
For every $M\in\bN,$ there exists $\mu^a_M\ge0$ such that, for every
$s\ge s_0$, there is $0<r^a_{s,M}\leq r_0$ such that, for every
\[
    \vec U_p\in
    \pare{
        \BallR{s_0+\mu^a_M}{r}
        \cap
        \cC\pare{I;H^{s+\mu^a_M}\pare{\bR^{\mathtt d};\bR}}
    }^p,
    \qquad 0<r\le r^a_{s,M},
\]
and every $x\in\bR^{\mathtt d}$, the map
\[
    \xi\mapsto a\pare{\vec U_p;x,\xi}
\]
belongs to
$\cC^\infty\pare{\bR^{\mathtt d}\setminus\set{0};\bC}$.

\item With the losses fixed above, for every $M\in\bN$, every $s\ge s_0$,
every multi-index $\beta\in\bN^{\mathtt d}$ with $\av{\beta}\le M$, every
\[
    \vec U_p\in
    \pare{
        \BallR{s_0+\mu^a_M}{r}
        \cap
        \cC\pare{I;H^{s+\mu^a_M}\pare{\bR^{\mathtt d};\bR}}
    }^p,
    \qquad 0<r\le r^a_{s,M},
\]
and every $\av{\xi}\ge \frac12$, one has
\begin{equation}\label{simb-pro-geq-q-geq-q}
    \norm{\partial_\xi^\beta a\pare{\vec U_p;\bullet,\xi}}_{s,I}
    \lesssim_{s,\beta,M}
    \norm{\vec U_p}_{s_0+\mu^a_M,I}^{q-1}
    \norm{\vec U_p}_{s+\mu^a_M,I}
    \av{\xi}^{m-\av{\beta}} .
\end{equation}
\end{enumerate}

The losses $\mu^a_M$ may be chosen nondecreasing in $M$.

\item We use the convention
\[
    \Gamma^m \defeq \Gamma^m_{\ge 1}.
\]

\item We define
\[
    \Gamma^m_{\ge 0}
    \defeq
    \tilde\Gamma^m+\Gamma^m_{\ge 1}.
\]
Equivalently, $a\in\Gamma^m_{\ge0}$ if it can be written as
\[
    a\pare{\vec U_p;x,\xi}
    =
    a_0\pare{\xi}+a_{\ge1}\pare{\vec U_p;x,\xi},
    \qquad
    a_0\in\tilde\Gamma^m,
    \qquad
    a_{\ge1}\in\Gamma^m_{\ge1}.
\]

\item For every $s\ge s_0$ and every $M\in\bN$, we set
\begin{equation}\label{seminorm-symbol}
    \av{a\pare{\vec U_p}}_{m,s,M}
    \defeq
    \max_{\av\beta\le M}
    \sup_{\av\xi\ge\frac12}
    \av\xi^{-m+\av\beta}
    \norm{
        \partial_\xi^\beta
        a\pare{\vec U_p;\bullet,\xi}
    }_{s, I}.
\end{equation}
In particular, the tame estimate \cref{simb-pro-geq-q-geq-q} defining $\Gamma^m_{\ge q}$ reads
\[
    \av{a\pare{\vec U_p}}_{m,s,M}
    \lesssim_{s,M}
    \norm{\vec U_p}_{s_0+\mu^a_M,I}^{q-1}
    \norm{\vec U_p}_{s+\mu^a_M,I}.
\]

For a background-independent symbol $a_0\in\widetilde\Gamma^m$, we set
\[
    \av{a_0}_{m,M}
    \defeq
    \max_{\av\beta\le M}
    \sup_{\av\xi\ge\frac12}
    \av\xi^{-m+\av\beta}
    \av{\partial_\xi^\beta a_0\pare\xi}.
\]
Accordingly, if
\begin{equation}\label{flat-nonlinear_decomp}
    a=a_0+a_{\ge1}\in\Gamma^m_{\ge0},
    \qquad
    a_0\in\widetilde\Gamma^m,
    \qquad
    a_{\ge1}\in\Gamma^m_{\ge1},
\end{equation}
we use the convention
\[
    \av{a\pare{\vec U_p}}_{m,s,M}
    \defeq
    \av{a_0}_{m,M}
    +
    \av{a_{\ge1}\pare{\vec U_p}}_{m,s,M}.
\]

\item For $N,N_1,N_2\in\bN^\star$, the notations
$\pare{\Gamma^m_{\ge q}}^{N_1\times N_2}$ and
$\pare{\Gamma^m_{\ge q}}^N$ are componentwise abbreviations, in accordance
with the scalar convention above.

\item A symbol $a\pare{\vec U_p;x,\xi}\in\Gamma^m_{\ge q}$ is
\emph{pointwise real} if the map
\[
    \pare{x,\xi}\mapsto a\pare{\vec U_p;x,\xi}
\]
is real-valued for every real-valued $\vec U_p$. It is
\emph{real-to-real} if
\[
    a\pare{\vec U_p;x,-\xi}
    =\overline{a\pare{\vec U_p;x,\xi}}.
\]
\end{enumerate}
\end{definition}

\begin{remark}[Elementary rules for the vanishing order]
\label{rem:Gamma-vanishing-rules}
Let $q_1,q_2\in\bN$ and $m_1,m_2\in\bR$.

\begin{enumerate}
\item If $m_1\leq m_2$ and $q_1\geq q_2$, then
$$\Gamma_{\geq q_1}^{m_1}\subset\Gamma_{\geq q_2}^{m_2}.$$
\item If $a\in\Gamma^{m_1}_{\ge q_1}$ and
$b\in\Gamma^{m_2}_{\ge q_2}$, then, using the Fourier multiplier--nonlinear
decomposition \cref{flat-nonlinear_decomp},
\[
    ab\in\Gamma^{m_1+m_2}_{\ge q_1+q_2}.
\]

\item If $a_0\in\tilde\Gamma^{m_1}$ and
$b\in\Gamma^{m_2}_{\ge q}$, then
\[
    a_0b\in\Gamma^{m_1+m_2}_{\ge q}.
\]

\item Sums satisfy
\[
    \Gamma^{m_1}_{\ge q_1}+\Gamma^{m_2}_{\ge q_2}
    \subseteq
    \Gamma^{\max\set{m_1,m_2}}_{\ge \min\set{q_1,q_2}}.
\]
\end{enumerate}
\fine
\end{remark}

\begin{definition}[Functions]\label{def:functions}
For $q\geq1$ we set
\[
    \cF_{\geq q}
    \defeq
    \set{f\in\Gamma^0_{\geq q}:\partial_\xi f=0},
    \qquad
    \cF\defeq\cF_{\geq1}.
\]
We denote by 
\[
    \cF_{\geq0}\defeq\bC +\cF_{\geq1}.
\]
The superscript $\bR$ denotes the pointwise real-valued subclass; thus
$\cF^\bR_{\geq q}$, $\cF^\bR$, and $\cF^\bR_{\geq0}$ have the analogous
meaning.
\end{definition}

\paragraph{Quantizations and para-composition rules.}
Let us recall the Weyl quantization formula of a symbol
\begin{equation}
\label{opW} 
\OpW{a} [u] \pare{x}  \defeq \frac 1{\pare{2\pi}^{\mathtt d}}\int_{\bR^{\mathtt d}}\int_{\bR^{\mathtt d}}e^{\ii x\cdot\xi} \hat a\pare{\xi-\varsigma,\frac{\xi+\varsigma}{2}}\hat u\pare\varsigma \ \dd \varsigma \dd\xi.
\end{equation}
Let $\epsilon \in \pare{0, \frac{1}{4}}$ and consider a 
smooth, radial cutoff
$\chi\in\cC^\infty\pare{\bR^{\mathtt d};\bra{0,1}}$ such that
\begin{align}\label{def-chi}
\chi(\xi) =  
\begin{cases}
1 & \textnormal{if } \av\xi \leq 1.1,\\
0 & \textnormal{if } \av\xi \geq 1.9,
\end{cases}  
&&
\chi_\epsilon(\xi) \defeq \chi \left( \frac{\xi}{\epsilon}\right) \, . 
\end{align}
The function $\pare{\varsigma, \xi} \mapsto
\chi_\epsilon\pare{\xfrac{\varsigma}{\av{\xi}}}$, extended by zero at
$\xi=0$.
\\

Given a symbol $a$ in $\Gamma^m_{\geq0}$
we define  the {\it regularized} symbol
\begin{equation}\label{achi}
a_{\chi_\epsilon} \pare{x, \xi} \defeq
\frac{1}{\pare{2\pi}^{\mathtt d/2}}
\int_{\bR^{\mathtt d}}e^{\ii  x\cdot\varsigma}\chi_\epsilon \pare{\frac{\varsigma}{\av\xi}}\hat a(\varsigma,\xi) \  \dd\varsigma
=
\chi_\epsilon \pare{\xfrac{D}{\av\xi}} a\pare{x, \xi},
\qquad \xi\ne0,
\end{equation}
and set $a_{\chi_\epsilon}(x,0)\defeq0$.

\begin{notation}[High-frequency extension and low-frequency remainder]
\label{notation:high-frequency-symbol-extension}
Let $m\in \bR$,  $a\in \Gamma^m_{\geq 0}$ smooth in $\xi\neq0$, and let $\chi$ be the even cutoff fixed in
\cref{def-chi}.  We decompose
\begin{equation}\label{eq:high-low-frequency-symbol-decomposition}
    a
    =
    a^\sharp
    +
    a^\flat,
    \qquad
    a^\sharp\pare{U;x,\xi}
    \defeq
    \pare{1-\chi\pare\xi}a\pare{U;x,\xi},
    \qquad
    a^\flat \pare{U;x,\xi}
    \defeq
    \chi\pare{\xi} a\pare{U;x,\xi} .
\end{equation}
\begin{itemize}
    \item The symbol $a^\sharp$ is the fixed high-frequency extension of $a$.
    It is independent of the chosen extension $a^{\textnormal{ext}}$.

    \item The symbol $a^\flat$ is the corresponding low-frequency remainder.
    For a complete homogeneous Fourier multiplier, the untruncated
    homogeneous expression remains the multiplier being defined on
    $\xi\ne0$; $^\sharp$ is used only for localized symbolic reductions,
    with the low-frequency part treated separately.
\end{itemize}
\end{notation}

\begin{definition}[Bony-Weyl quantization]
Let $\epsilon \in \pare{0, \frac{1}{4}}$, $m\in \bR$,
$a \in \Gamma^m_{\geq0}$, and
$u\in \mathcal{S}\pare{\bR^{\mathtt d};\bC}$. We define the
\emph{Bony-Weyl paradifferential operator}
$\OpBW{a} \defeq \OpW{a_{\chi_{\epsilon}}}$, i.e.
\begin{equation}
\label{BonyWeyl}
\OpBW{a}\bra{u} \pare{x} = \frac{1}{\pare{2\pi}^{\mathtt d}}\int_{\bR^{\mathtt d}}\pare{\int_{\bR^{\mathtt d}} \chi_{\epsilon}\pare{\frac{2(\xi-\varsigma)}{\av{\xi+\varsigma}} }\hat{a}\pare{\xi-\varsigma,\frac{\xi+\varsigma}{2}}\hat{u} \pare{\varsigma} \dd \varsigma} e^{\ii x\cdot\xi} \dd\xi.
\end{equation}

If
$\mathsf A=\bra{a_{ij}}_{1\leq i\leq N_1,\,1\leq j\leq N_2}$ is a
finite matrix of symbols in $\Gamma^m_{\geq0}$, then
$\OpBW{\mathsf A}$ is defined componentwise by
\begin{equation}
    \OpBW{\mathsf A}
    \defeq
    \bra{\OpBW{a_{ij}}}_{1\leq i\leq N_1,\,1\leq j\leq N_2}.
\end{equation}
\end{definition}

Notice that the diagonal $\set{\xi=-\varsigma}$ is a set of measure zero in $\bR^{2\mathtt d}$. Moreover, there is no issue of integrability near this diagonal, because the integrand of \cref{BonyWeyl} vanishes in a neighborhood of it. Moreover, if 
$ \chi_\epsilon \Big( \frac{2(\xi -\varsigma)}{\av{\xi + \varsigma} }\Big) \neq 0 $
then $ |\xi-\varsigma| \leq \epsilon \av {\varsigma + \xi }  $,
and therefore, for $ \epsilon \in \pare{0, \tfrac{1}{4}}$, 
\begin{align}\label{rela:para}
\frac{1-\epsilon}{1+\epsilon}\av\varsigma \leq
\av\xi \leq \frac{1+\epsilon}{1-\epsilon}\av\varsigma,
&& \forall\,\xi,\varsigma\in\bR^{\mathtt d}.
\end{align}

\begin{lemma}[Schur test]
\label{lem:Schur_test}
Let
\[
    K:\bR^{\mathtt d}\times\bR^{\mathtt d}\to\bC
\]
be measurable and assume that
\begin{equation}\label{eq:Schur_test_assumptions}
    A_K
    \defeq
    \sup_{\varsigma\in\bR^{\mathtt d}}
    \int_{\bR^{\mathtt d}}
    \av{K\pare{\xi,\varsigma}}\,\dd\xi
    <+\infty,
    \qquad
    B_K
    \defeq
    \sup_{\xi\in\bR^{\mathtt d}}
    \int_{\bR^{\mathtt d}}
    \av{K\pare{\xi,\varsigma}}\,\dd\varsigma
    <+\infty.
\end{equation}
Then the integral operator
\[
    T_K f\pare\xi
    \defeq
    \int_{\bR^{\mathtt d}}
    K\pare{\xi,\varsigma}f\pare\varsigma\,\dd\varsigma
\]
extends to an operator in
\(\cL\pare{L^2\pare{\bR^{\mathtt d}}}\), and
\begin{equation}\label{eq:Schur_test_bound}
    \norm{T_K}_{\cL\pare{L^2}}
    \leq
    \sqrt{A_KB_K}.
\end{equation}
\end{lemma}

\begin{proof}
For \(f\in\mathcal S\pare{\bR^{\mathtt d}}\), Cauchy--Schwarz with respect
to the measure
\(\av{K\pare{\xi,\varsigma}}\dd\varsigma\) gives
\[
    \av{T_Kf\pare\xi}^2
    \leq
    \left(
        \int_{\bR^{\mathtt d}}
        \av{K\pare{\xi,\varsigma}}\,\dd\varsigma
    \right)
    \left(
        \int_{\bR^{\mathtt d}}
        \av{K\pare{\xi,\varsigma}}
        \av{f\pare\varsigma}^2\,\dd\varsigma
    \right).
\]
Integrating in \(\xi\), using \cref{eq:Schur_test_assumptions} and
Fubini's theorem, we obtain
\[
\begin{aligned}
    \norm{T_Kf}_{L^2}^2
    &\leq
    B_K
    \int_{\bR^{\mathtt d}}
    \int_{\bR^{\mathtt d}}
    \av{K\pare{\xi,\varsigma}}
    \av{f\pare\varsigma}^2
    \,\dd\varsigma\dd\xi
    \leq
    A_KB_K\norm f_{L^2}^2.
\end{aligned}
\]
The conclusion follows by density.
\end{proof}

The following lemma is immediate from the definition of Bony-Weyl quantization.
\begin{lemma}[Adjoint of Bony-Weyl quantization]\label{lem:adjoint}
Let $m\in\bR$ and $a\in\Gamma^m_{\geq0}$, then the adjoint operator $\OpBW{a}^{\star}$ is given by
    $$\OpBW{a}^{\star}=\OpBW{\overline{a}}.$$
\end{lemma}
\begin{proof}
    Cf. \cite{BD2018}. 
\end{proof}
\begin{remark}\label{rem:action_paradiff_constants}
Let $a\in \Gamma^m_{\geq 1}$, the definition of the operator $\OpBW{a}$ can be naturally extended to tempered distributions by duality, using \cref{lem:adjoint}:
    \begin{align}
        \psc{\OpBW{a}\bra{\varphi'}}{\varphi}
        \defeq
        \psc{\varphi'}{\OpBW{\bar{a}}\varphi},
        &&
        \forall\,
        \varphi\in\mathcal S\pare{\bR^{\mathtt d};\bC},\quad
        \varphi'\in\mathcal S'\pare{\bR^{\mathtt d};\bC}
    \end{align}
    thus proving that constants belong to $\ker\pare{\OpBW{a}}$, in fact
    \begin{align}
    \psc{\OpBW{a}\bra{C}}{\varphi}
    &= \psc{C}{\OpBW{\overline{a}}\bra{\varphi}}\\
    &= \psc{\hat{C}}{\mathcal{F}\pare{\OpBW{\overline{a}}\bra{\varphi}}}\\
    &= C\int_{\bR^{\mathtt d}}\int_{\bR^{\mathtt d}}\delta_{0}(\xi)\overline{\overline{\hat{a}_{\chi_{\epsilon}}\pare{\xi-\varsigma,\frac{\xi+\varsigma}{2}}}\hat{\varphi}(\varsigma)}\dd\varsigma \dd\xi\\
    &= C\int_{\bR^{\mathtt d}}\underbrace{\hat{a}_{\chi_{\epsilon}}\pare{-\varsigma,\frac{\varsigma}{2}}}_{=0}\overline{\hat{\varphi}(\varsigma)}\dd\varsigma
    = 0.
\end{align}
    \fine
\end{remark}

\begin{theorem}[Action of paradifferential operators]
\label{thm:contS}
Let $m\in\bR$ and $a\in\Gamma^m_{\geq0}$. There exists
$N_{\mathtt d}\in\bN$, depending only on $\mathtt d$, such that
 $\OpBW{a \pare{U; x, \xi}}$ extends to a bounded operator in
\[
    \cL\pare{
        H^s\pare{\bR^{\mathtt d};\bC};
        H^{s-m}\pare{\bR^{\mathtt d};\bC}
    }
\]
for every $s\in\bR$, and
\begin{equation} \label{cont00}
\norm{\OpBW{a\pare{U}}}_{\cL\pare{H^s\pare{\bR^{\mathtt d};\bC};H^{s-m}\pare{\bR^{\mathtt d};\bC}}}
\lesssim
\av{a\pare{U}}_{m,s_0,N_{\mathtt d}} + \norm{ \OpBW{a^\flat\pare U} }_{ \cL\pare{ H^s\pare{\bR^{\mathtt d};\bC}; H^{s-m}\pare{\bR^{\mathtt d};\bC} }},
\end{equation}
where $s_0$ is as in \cref{def:sfr}.
\end{theorem}

This is the standard whole-space Bony--Weyl boundedness theorem; see, for
instance, \cite[Theorem~2.4]{berti2021local} for the periodic counterpart.

\begin{rem}
The seminorm in \cref{cont00} controls only the region \(\av{\xi}\geq1/2\); the low-frequency component \(a^\flat\) has therefore to be estimated separately.
This distinction is essential for symbols that are singular or merely defined away from \(\xi=0\). \fine
\end{rem}

We recall also the following classical results.

\begin{rem}\label{rem:functions_Op_R}
    Notice that for $q\in\bN^{\star}$ and $m,\varrho\in\mathbb{R},$ 
\begin{enumerate}
    \item If $a\in \Gamma_{\geq q}^m$ then $\OpBW{a\pare{U;x, \xi}} U \in \cF_{\geq q}$;
    \item  if $r_{\varrho}\in \cR_{\geq q}^\varrho$  then $r_{\varrho}\pare{U} U \in \cF_{\geq q} $;
    \item if $\mathpzc{m}\in \tilde\Gamma^m$, then $\OpBW{\mathpzc{m}} = \mathpzc{m}\pare{D}$;
    \item \cref{thm:contS} implies that if $a\in \Gamma_{\geq q}^m$, then $\OpBW{a}\in \cR_{\geq q}^{-m}$.
\end{enumerate}
\fine
\end{rem}

\begin{lemma}[Bony paraproduct decomposition]\label{lem:paraprod_Bony}
Let $r_1+r_2>0$ and
$u\in H^{r_1}\pare{\bR^{\mathtt d};\bC}$,
$v\in H^{r_2}\pare{\bR^{\mathtt d};\bC}$. Then
\begin{equation}
uv = \OpBW{u}  {v} + \OpBW{v}  {u} + \RBW{u}{v}
\end{equation}
where $\RBW{u}{v}$ is a symmetric bilinear map
\[
    H^{r_1}\pare{\bR^{\mathtt d};\bC}
    \times
    H^{r_2}\pare{\bR^{\mathtt d};\bC}
    \longrightarrow
    H^{r_1+r_2-\mathtt d/2}\pare{\bR^{\mathtt d};\bC},
\]
and
\[
\norm{\RBW{u}{v}}_{r_1+r_2-\mathtt d/2}\lesssim \norm {u}_{r_1}\norm v_{r_2}.
\]
\end{lemma}
\begin{proof}
    See Lemma 2.7 of \cite{berti2021local} for the periodic setting. 
\end{proof}

{
\begin{rem}\label{rem:smoothBONYremainder}
For every $\varrho\ge0$ and $s_0>\mathtt d/2+\varrho$, the operator
$R\pare u h\defeq\RBW{u}{h}$ satisfies
\[
    \norm{R\pare u h}_{H^s}
    \lesssim_s
    \norm u_{H^{s_0}}
    \norm h_{H^{s-\varrho}},
    \qquad s\ge s_0.
\]
Indeed, apply \cref{lem:paraprod_Bony} with
$r_1=s_0$ and $r_2=s-\varrho$ and use
$s_0+s-\varrho-\mathtt d/2>s$.
Thus $R\in\cR^\varrho_{\ge1}$.
\fine
\end{rem}
}

{
\begin{lemma}[Bony paralinearization formula]\label{lem paralin}
Let $\varrho\ge0$ and let $f\in\cC^\infty\pare{\bR;\bR}$ satisfy
$f\pare0=0$. There exist $s_0>\mathtt d/2+\varrho$, $r_0>0$, and
$r_f\in\cR^\varrho_{\ge1}$ such that, for every $s\ge s_0$ and
$u\in\BallR{s_0}{r_0}\cap\cC\pare{I;H^s}$,
\begin{equation}
    f\pare u=\OpBW{f'\pare u}u+r_f\pare u u.
\end{equation}
If $f^{\pare j}\pare0=0$ for $0\le j<q$, with $q\ge2$, then
$r_f$ may be chosen in $\cR^\varrho_{\ge q-1}$.
\end{lemma}
\begin{proof}
Use the dyadic paralinearization in \cite{metivier2008differential}
and the change of quantization in \cite{BD2018}.
Keeping one maximal-frequency factor as the argument gives
\[
    \norm{r_f\pare u h}_{H^s}
    \lesssim_{s,\varrho}
    \norm u_{H^{s_0}}
    \norm h_{H^{s-\varrho}},
    \qquad s\ge s_0,
\]
after increasing $s_0$ according to $\varrho$.
If $f$ vanishes to order $q\ge2$, Taylor's formula leaves $q-1$
background factors and gives the asserted vanishing order.
\end{proof}
}

\begin{definition}[Poisson bracket]
\begin{itemize}
    \item   Let $a,b\in \cC^1\pare{\bR^{\mathtt d}\times\bR^{\mathtt d};\bC}$, we define the \emph{Poisson bracket} as
 \begin{equation}
      \pbra{a}{b} \defeq \nabla_\xi a\cdot\nabla_x b-\nabla_x a\cdot\nabla_\xi b
      =
      J_\xi a \nabla_x b - J_x a \nabla_\xi b
      \ ;
\label{eq:Pbra}
 \end{equation}
 \item   
Let $N_1,N_2,N_3 \in \bN^\star$, and
\[\mathsf A\in \cC^1\pare{\bR^{\mathtt d}\times\bR^{\mathtt d};\bC^{N_1\times N_2}},\quad
\mathsf B\in \cC^1\pare{\bR^{\mathtt d}\times\bR^{\mathtt d};\bC^{N_2\times N_3}}.\]
We use the notation
\begin{align}\label{eq:Pbra_matrix-vector}
	  \pbra{\mathsf A}{\mathsf B}
      \defeq
      \bra{\sum_{j=1}^{N_2} \pbra{a_{ij}}{b_{jk}}}_{\substack{i=1,\ldots,N_1 \\ k=1,\ldots,N_3}}
      \in
      \cC^0\pare{\bR^{\mathtt d}\times\bR^{\mathtt d};\bC^{N_1\times N_3}}.
\end{align}
\end{itemize}

\end{definition}

\begin{definition}[Asymptotic expansion of composition symbol]\label{def:asympt_symbol}
Let $m_1, m_2 \in \bR$, $\varrho\geq 0$, and
\[
    \sigma\pare{D_x,D_\xi,D_y,D_\varsigma}
    \defeq
    D_\xi\cdot D_y-D_x\cdot D_\varsigma,
\]
where $D_x=\frac 1{\ii}\partial_x$ and
$D_\xi,D_y,D_\varsigma$ are similarly defined.
Let $a\in \Gamma^{m_1}_{\geq0}$ and $b\in \Gamma^{m_2}_{\geq0}$.
Choose a common
admissible base index $s_0$, the finite losses required by the displayed
expansion, and a radius below the two radii supplied by \cref{def:sfr}.
For $s>s_0+\varrho$ and $U$ in the corresponding ball, we define
 \begin{equation}
     \label{eq:compositionsymbol1}
     \pare{a\#_\varrho b}\pare{U;x,\xi}\defeq \sum_{k=0}^{\floor{\varrho}}\frac{1}{k!}\pare{\frac{\ii}{2}\sigma\pare{D_x,D_\xi,D_y,D_\varsigma}}^k\left [a\pare{U;x,\xi}b\pare{U;y,\varsigma}\right]\Big |_{x=y,\xi=\varsigma}. 
 \end{equation}
\end{definition}

\begin{rem}\label{rem comp}
    \begin{enumerate}
        \item For $\floor{\varrho}=1$, we have that
    \begin{equation}
        a \#_\varrho b = ab + \frac{1}{2\ii }\pbra{a}{b}.
    \end{equation}
    \item If $a\in\Gamma_{\geq q_1}^{m_1}$ and $b\in\Gamma_{\geq q_2}^{m_2}$ then 
    $$\pbra{a}{b}\in\Gamma^{m_1+m_2-1}_{\ge q_1+q_2}.$$
    \end{enumerate}
    \fine
\end{rem}

\paragraph{Composition rules.}

 The following result is classical, we refer the reader to \cite[Theorem 2.5]{berti2021local}
\begin{prop}[Composition with vanishing order]\label{prop:composition_paradiff}
Let $m_1,m_2\in\bR$, $\varrho\ge0$, $q_1,q_2\in\bN$, and let
$a\in\Gamma^{m_1}_{\ge q_1}$, $b\in\Gamma^{m_2}_{\ge q_2}$.
Assume moreover that the low-frequency contributions in \cref{cont00}
associated with $a$, $b$, and $a\#_\varrho b$ are bounded.
Define
\begin{equation}\label{eq:remainder_composition_vanishing_order}
    r_\varrho\bra{a,b}
    \defeq
    \OpBW{a}\circ\OpBW{b}
    -
    \OpBW{a\#_\varrho b}.
\end{equation}
Then $a\#_\varrho b\in\Gamma^{m_1+m_2}_{\ge q_1+q_2}$ and
\begin{equation}
    r_\varrho\bra{a,b}
    \in
    \cR^{\varrho-m_1-m_2}_{\ge q_1+q_2}.
\end{equation}
When $q_i=0$, the assertion is understood after splitting the corresponding
symbol into its flat and nonlinear parts.
\end{prop}

\begin{lemma}\label{lem:composition_rules}
Let $\varrho,\varrho',m\in\bR$, let $q_a,q_f,q_r,q_{r}'\in\bN$, and let
$a\in\Gamma^m_{\ge q_a}$, $f\in\cF_{\ge q_f}$,
$r_{\varrho}\in\cR^\varrho_{\ge q_r}$, and
$r_{\varrho'}\in\cR^{\varrho'}_{\ge q_{r}'}$. Then the following tame rules apply,
after choosing a common, sufficiently large baseline and reducing the
radius. For full paradifferential operators, the low-frequency bound in
\cref{cont00} is assumed with vanishing order $q_a$:
    \begin{enumerate}[\bf i]
        \item\label{item:composition_rules_paradiff_remainder}
        $\OpBW a\circ r_{\varrho},\ r_{\varrho}\circ\OpBW a
        \in\cR^{\varrho-m}_{\ge q_a+q_r}$;
        \item\label{item:composition_rules_remainder_remainder}
        $r_{\varrho}\circ r_{\varrho'}\in\cR^{\varrho+\varrho'}_{\ge q_r+q_{r}'}$;
\item\label{item:composition_rules_substitution_symbol}
Assume that $a\pare{\vec U_p,v;x,\xi}
    \in\Gamma^m_{\ge q_a}
    $, $
    f\pare{\vec U_p}\in\cF_{\ge1}$, then, after possibly reducing the radius of definition, the symbol $\widetilde a\pare{\vec U_p;x,\xi}
    \defeq
    a\pare{
        \vec U_p,
        f\pare{\vec U_p};
        x,\xi
    }$
satisfies $\widetilde a\in\Gamma^m_{\ge q_a}$.

\item\label{item:composition_rules_substitution_remainder}
Assume that $r_{\varrho}\pare{\vec U_p,v}
    \in\cR^\varrho_{\ge q_r},
    $ $
    f\pare{\vec U_p}\in\cF_{\ge1}$, then, after possibly reducing the radius of definition, the operator $\widetilde r_{\varrho}\pare{\vec U_p}
    \defeq
    r_{\varrho}\pare{
        \vec U_p,
        f\pare{\vec U_p}
    }$
satisfies $ \widetilde r_{\varrho}\in\cR^\varrho_{\ge q_r}$. 
    \end{enumerate}
\end{lemma}

\begin{proof}
Choose a common baseline $S$ large enough for all fixed regularity
indices below. By \cref{thm:contS,eq:smoothing-tame-geq-q},
\begin{align}
    \norm{\OpBW a\,r_\varrho\pare U h}_{H^s}
    &\lesssim_s
    \norm U_S^{q_a}\norm{r_\varrho\pare U h}_{H^{s+m}}
    \lesssim_s
    \norm U_S^{q_a+q_r}\norm h_{H^{s+m-\varrho}},
    \\
    \norm{r_\varrho\pare U\OpBW a\,h}_{H^s}
    &\lesssim_s
    \norm U_S^{q_r}\norm{\OpBW a\,h}_{H^{s-\varrho}}
    \lesssim_s
    \norm U_S^{q_a+q_r}\norm h_{H^{s+m-\varrho}},
    \\
    \norm{r_\varrho\pare U r_{\varrho'}\pare U h}_{H^s}
    &\lesssim_s
    \norm U_S^{q_r}\norm{r_{\varrho'}\pare U h}_{H^{s-\varrho}}
    \lesssim_s
    \norm U_S^{q_r+q_r'}\norm h_{H^{s-\varrho-\varrho'}}.
\end{align}
This proves \cref{item:composition_rules_paradiff_remainder,%
item:composition_rules_remainder_remainder}.
For the substitutions, \cref{def:functions,def:sfr} give
$\norm{f\pare U}_{s_0}\lesssim\norm U_{s_0+\mu}$ for a fixed finite
loss $\mu$. Increase $S$ and reduce the radius so that
$\pare{U,f\pare U}$ belongs to the required low-norm ball.
Substitution in \cref{simb-pro-geq-q-geq-q,eq:smoothing-tame-geq-q}
gives \cref{item:composition_rules_substitution_symbol,%
item:composition_rules_substitution_remainder}.
\end{proof}

\subsection{$z$-dependent paradifferential calculus}

The purpose of this subsection is to collect the $z$-dependent
paradifferential tools used in the paralinearization of the singular
integrals. We first introduce the scalar kernel classes and the corresponding
operator-valued smoothing classes, together with their elementary composition
and integration rules. We then prove the kernel-to-symbol reduction for
$z$-dependent Bony--Weyl operators and its compensated endpoint version. The
remaining paragraphs deal with the two special mechanisms used in the
paralinearization: far-field tails, which give order-zero paraproducts
modulo smoothing remainders, and localized homogeneous convolution models, whose cutoff
can be removed modulo smoothing Fourier multipliers. Finally, we record the
local Taylor and finite-difference estimates needed to verify that the
concrete kernels arising from $B_0$ and $q$ belong to the abstract classes
introduced above.

\paragraph{Kernel classes.}

\begin{definition}[Kernel classes]
\label{def:kernel-global-global}
Let $n\in\bR$, $q\in\bN^\star$, $p\in\bN^\star$, and let
$\vec U_p\defeq\pare{u_1,\ldots,u_p}^{\intercal}$. 

\begin{enumerate}
    \item The class $\KF^n_{\loc,\ge q}$ consists of scalar kernels
    $k=k\pare{\vec U_p;x,z}$ for which there exists
    $s_0^k>\mathtt d/2$ such that the following holds:
    for every $N\in\bN$, there exists a loss $\mu_N^k\ge0$ such that,
    for every $s\ge s_0^k$, there exists a radius $r_{s,N}^k>0$.
    For every $0<r\le r_{s,N}^k$ and every \[ \vec U_p
        \in
        \pare{
            \BallR{s_0^k+\mu_N^k}{r}
            \cap
            \cC\pare{
                I;
                H^{s+\mu_N^k}\pare{\bR^{\mathtt d};\bR}
            }
        }^p,\]
    the map $z\mapsto k\pare{\vec U_p;\bullet,z}$ is measurable on
    $\bR^{\mathtt d}$ and belongs to $\cC^N\pare{
            \bR^{\mathtt d}\setminus\set{0};
            H^s\pare{\bR^{\mathtt d};\bC}
        }$. 
    Moreover, for every $\alpha\in\bN^{\mathtt d}$ with
    $\av\alpha\le N$ and every $0<\av z\le2$, one has
    \begin{equation}\label{eq:kernelfunctionproperty-local-tame}
        \norm{
            \partial_z^\alpha
            k\pare{\vec U_p;\bullet,z}
        }_{s,I}
        \lesssim_{s,\alpha,N}
        \norm{\vec U_p}_{s_0^k+\mu_N^k,I}^{q-1}
        \norm{\vec U_p}_{s+\mu_N^k,I}
        \av z^{\,n-\av\alpha}.
    \end{equation}

    \item The class $\KF^n_{\ge q}$ consists of the kernels
    $k\in\KF^n_{\loc,\ge q}$ satisfying
    \begin{equation}\label{eq:compact-kernel-support}
        \supp_z k\pare{\vec U_p;x,\cdot}
        \subseteq
        \set{z\in\bR^{\mathtt d}:\av z\le2}.
    \end{equation}
    Under the same assumptions on $s$, $r$, and $\vec U_p$ as above,
    for every $\alpha\in\bN^{\mathtt d}$ with $\av\alpha\le N$ and every
    $z\neq0$, one has
    \begin{equation}\label{eq:kernelfunctionproperty-tame}
        \norm{
            \partial_z^\alpha
            k\pare{\vec U_p;\bullet,z}
        }_{s,I}
        \lesssim_{s,\alpha,N}
        \norm{\vec U_p}_{s_0^k+\mu_N^k,I}^{q-1}
        \norm{\vec U_p}_{s+\mu_N^k,I}
        \av z^{\,n-\av\alpha}
        \mathbb 1_{\set{0<\av z\le2}}.
    \end{equation}

\item The class $\KF_{\textnormal{tail},\ge q}$ consists of scalar maps
\[
    k_{\textnormal{tail}}:
    \BallR{s_0}{r_0}
    \longrightarrow
    L^\infty_{\loc}\pare{
        \bR^{\mathtt d}_x\times\bR^{\mathtt d}_z;\bC
    },
    \qquad
    u\longmapsto k_{\textnormal{tail}}\pare{u;x,z},
\]
for which there exist $s_0>\mathtt d/2$, $r_0>0$, $m>0$, $r_{\mathsf g}>0$,  $ b\in
    \cC^\infty\pare{
        \bR^{\mathtt d}\setminus\set0;\bC
    }$ and $
    \mathsf g\in
    \cC^\omega\pare{
        \pare{-r_{\mathsf g},r_{\mathsf g}};\bC
    },$
such that
\[
    k_{\textnormal{tail}}\pare{u;x,z}
    \defeq
    \pare{1-\chi\pare z}
    b\pare z\,
    \mathsf g\pare{\Delta_z u\pare x},
\]
where $\chi$ is the cutoff function introduced in \cref{def-chi}. The function $b$ satisfies, for every
$\alpha\in\bN^{\mathtt d}$,
\begin{equation}\label{eq:tail-kernel-profile-decay}
    \av{\partial_z^\alpha b\pare z}
    \lesssim_\alpha
    \av z^{-m-\av\alpha},
    \qquad
    \av z\ge1,
\end{equation}
while $\mathsf g$ vanishes to order at least $q$ at the origin,
\begin{equation}\label{eq:tail-kernel-analytic-vanishing}
    \mathsf g^{\pare j}\pare0=0,
    \qquad
    0\le j\le q-1.
\end{equation}
We further require $m+q>\mathtt d$.

    \item We use the conventions
    \[
        \KF^n_{\loc}\defeq\KF^n_{\loc,\ge1},
        \qquad
        \KF^n\defeq\KF^n_{\ge1},
        \qquad
        \KF_{\textnormal{tail}}
        \defeq
        \KF_{\textnormal{tail},\ge1}.
    \]
\end{enumerate}
\end{definition}

\begin{rem}[Cutoff of a local kernel]
\label{rem:kernel-localization}
For every $n\in\bR$ and $q\in\bN^\star$,
\begin{align}\label{eq:local-kernel-localization}
    k\in\KF^n_{\loc,\geq q}
    \quad\Longrightarrow\quad
    \chi\pare z k\in\KF^n_{\geq q},
    &&
    \KF^n_{\geq q}\subseteq\KF^n_{\loc,\geq q}.
\end{align}
Indeed, the inclusion on the right is immediate. For the implication on the
left, apply Leibniz' rule to $\partial_z^\alpha\pare{\chi\pare z k}$. The terms
in which no derivative falls on $\chi$ are estimated by
\cref{eq:kernelfunctionproperty-local-tame}, while the remaining terms are
supported in the fixed annulus where all powers of $\av z$ are harmless. The
support property follows from the compact support of $\chi$. Thus local Taylor
remainders become compactly supported kernels after multiplication by the fixed
cutoff $\chi\pare z$. \fine
\end{rem}

\begin{rem}[Integration of rapidly decreasing scalar kernels]
\label{lem:integration_KF}
Let $n>-\mathtt d$, $q\in\bN^\star$, and $k\in\KF^n_{\geq q}$. Then
\begin{equation}\label{eq:integration_KF}
    v_k\pare{\vec U_p;x}
    \defeq
    \int_{\bR^{\mathtt d}}k\pare{\vec U_p;x,z}\dd z
    \in\cF_{\geq q}.
\end{equation}
Moreover, for every $s\geq s_0^k$ one has
\begin{equation}\label{eq:integration_KF_tame}
    \norm{v_k\pare{\vec U_p;\bullet}}_{s,I}
    \lesssim_s
    \norm{\vec U_p}_{s_0^k+\mu,I}^{q-1}
    \norm{\vec U_p}_{s+\mu,I}
\end{equation}
for a finite loss $\mu\geq0$. Indeed, this follows directly from
Minkowski's integral inequality, \cref{eq:kernelfunctionproperty-tame}, and
$\int_0^2 r^{n+\mathtt d-1}\dd r<+\infty$. The localized case is reduced to this one by
\cref{rem:kernel-localization}. \fine
\end{rem}

\paragraph{Operator-valued smoothing kernel classes.}

\begin{definition}[Kernel-smoothing operators]
\label{def:kernel-global-smoothing}
Let $\varrho,n\in\bR$, $p\in\bN^\star$, $q\in\bN$, and let
$\vec U_p\defeq\pare{u_1,\ldots,u_p}^{\intercal}$. 

\begin{enumerate}
    \item The class $\KR^{\varrho,n}_{\loc,\ge q}$ consists of
    $z$-dependent families of linear operators
    $r_{\varrho}\pare{\vec U_p;z}$ for which there exists
    $s_0^R>\mathtt d/2$ such that, for every $s\ge s_0^R$, there exists
    a radius $r_s^R>0$ with the following property.

    For every $0<r\le r_s^R$, every \[\vec U_p\in
        \pare{\BallR{s_0^R}{r}}^p,\] and every $h\in H^{s-\varrho}\pare{\bR^{\mathtt d};\bC}$, the map
    $z\mapsto r_{\varrho}\pare{\vec U_p;z}h$ is measurable and, for
    $0<\av z\le2$, satisfies
    \begin{equation}\label{eq:kernel-smoothing-local-tame-geq-q}
    \begin{aligned}
        \norm{r_{\varrho}\pare{\vec U_p;z}h}_{s,I}
        \lesssim_{s,\varrho}{}&
        \norm{\vec U_p}_{s_0^R,I}^{q}
        \norm h_{{s-\varrho}}
        \av z^n .
    \end{aligned}
    \end{equation}

    \item The class $\KR^{\varrho,n}_{\ge q}$ consists of the families
    $r_{\varrho}\in\KR^{\varrho,n}_{\loc,\ge q}$ satisfying $\supp_z r_{\varrho}\pare{\vec U_p;\cdot}
        \subseteq
        \set{z\in\bR^{\mathtt d}:\av z\le2}$, and, for every $z\neq0$,
    \begin{equation}\label{eq:kernel-smoothing-tame-geq-q}
    \begin{aligned}
        \norm{r_{\varrho}\pare{\vec U_p;z}h}_{s,I}
        \lesssim_{s,\varrho}{}&
        \norm{\vec U_p}_{s_0^R,I}^{q}
        \norm h_{{s-\varrho}}
        \av z^n
        \mathbb 1_{\set{0<\av z\le2}},
    \end{aligned}
    \end{equation}
    under the same assumptions on $s$, $\vec U_p$, and $h$ as above.

    \item For $q=0$, the background factor in the preceding estimates is
    understood as $1$. We use the conventions
    \[
        \KR^{\varrho,n}_{\loc}
        \defeq
        \KR^{\varrho,n}_{\loc,\ge1},
        \qquad
        \KR^{\varrho,n}
        \defeq
        \KR^{\varrho,n}_{\ge1}.
    \]
\end{enumerate}
\end{definition}
}

\begin{rem}[Separately estimated parameters]
\label{rem:separately-estimated-parameters}
When a separately estimated parameter appears in brackets, as in
\(a\bra\omega\), \(k\bra\omega\), or \(r\bra\omega\), the notation
\[
    \cL\pare{H^{s_0};Y}
\]
is understood according to the structure of the target class \(Y\).

If
\[
    Y\in\set{
        \Gamma^m,\,
        \cF,\,
        \KF^n_{\loc},\,
        \KF^n
    },
\]
the notation records linear dependence on the parameter together with
the corresponding tame Sobolev-scale estimates. For instance, if
\(a_{\bra{m,\geq q}}\bra\omega\pare U\) is linear in \(\omega\) and has
vanishing order at least \(q\) in the sense of
\cref{rem:linear-parameter-vanishing-order}, then, schematically, for every
\(s\geq s_0\) and \(M\in\bN\),
\[
\begin{aligned}
    \av{
        a_{\bra{m,\geq q}}\bra\omega\pare U
    }_{m,s,M}
    \lesssim_{s,M}{}&
    \norm U_{H^{s_0+\mu_M}}^{q-1}
    \norm\omega_{H^{s+\mu_M}}
    \\
    &+
    \mathbb 1_{\set{q\geq2}}
    \norm U_{H^{s_0+\mu_M}}^{q-2}
    \norm U_{H^{s+\mu_M}}
    \norm\omega_{H^{s_0+\mu_M}},
\end{aligned}
\]
for some finite loss \(\mu_M\geq0\) independent of \(s\). The analogous
interpretation is used for functions and kernel functions.

If instead
\[
    Y\in\set{
        \cR^\varrho,\,
        \KR^{\varrho,n}_{\loc},\,
        \KR^{\varrho,n}
    },
\]
the notation refers to a fixed-baseline estimate. For instance, a remainder
\(r_{\bra{\varrho,\geq q}}\bra\omega\pare U\) linear in \(\omega\)
satisfies, for some \(S\) independent of \(s\),
\[
    \norm{
        r_{\bra{\varrho,\geq q}}\bra\omega\pare U h
    }_{H^s}
    \lesssim_{s,\varrho}
    \norm U_{H^S}^{q-1}
    \norm\omega_{H^S}
    \norm h_{H^{s-\varrho}}.
\]
In particular, no Sobolev norm of the parameter whose index grows with
the output index \(s\) is allowed in the smoothing estimates.

The convention above concerns only the Sobolev regularity of the separately
estimated parameter; its contribution to the vanishing order is specified in
\cref{rem:linear-parameter-vanishing-order}.
\fine
\end{rem}

\begin{rem}[Vanishing order for linear arguments]
\label{rem:linear-parameter-vanishing-order}
When an object in any of the classes introduced above depends linearly on a
separately estimated Sobolev argument, this argument contributes one to its
vanishing order. Thus, if \(T\bra\omega\pare U\) is linear in \(\omega\) and
vanishes to order at least \(q-1\) with respect to \(U\), we say that
\(T\bra\omega\pare U\) has vanishing order at least \(q\).
The same convention is used for functions, symbols, smoothing operators,
kernel functions, and kernel-smoothing operators.
\fine
\end{rem}

\begin{rem}[Cutoff of a local smoothing kernel]
\label{rem:kernel-smoothing-localization}
For every $\varrho,n\in\bR$ and $q\in\bN$,
\[
    R\in\KR^{\varrho,n}_{\loc,\geq q}
    \quad\Longrightarrow\quad
    \chi\pare zR\in\KR^{\varrho,n}_{\geq q},
    \qquad
    \KR^{\varrho,n}_{\geq q}
    \subseteq
    \KR^{\varrho,n}_{\loc,\geq q}.
\]
Indeed, multiplication by $\chi$ preserves the tame estimate in
\cref{eq:kernel-smoothing-local-tame-geq-q} and makes the family compactly
supported in $z$. The inclusion is immediate.
\fine
\end{rem}

\begin{lemma}[Composition rules for $z$-dependent operators]
\label{lem:composition_rules_z}
Let $q_1,q_2\in\bN^\star$, let $m,m',n,n'\in\bR$, and let
$\varrho,\varrho'\geq0$.

\begin{enumerate}[\bf i]
    \item\label{item:zcomp1}
    If
    \[
        k\in\KF^m_{\loc,\geq q_1},
        \qquad
        k'\in\KF^{m'}_{\loc,\geq q_2},
    \]
    then, for every $\varrho>0$,
    \begin{equation}\label{eq:zcomp1-local}
        \OpBW{
            kk'\pare{\vec U_p;x,z}
        }-\OpBW{k\pare{\vec U_p;x,z}}
        \circ
        \OpBW{k'\pare{\vec U_p;x,z}}\in
        \KR^{\varrho,m+m'}_{\loc,\geq q_1+q_2}.
    \end{equation}

    \item\label{item:zcomp2}
    If
    \[
        k\in\KF^m_{\loc,\geq q_1},
        \qquad
        r_{\varrho}\in\KR^{\varrho,n}_{\loc,\geq q_2},
    \]
    then
    \begin{equation}\label{eq:zcomp2-local}
        \OpBW{k\pare{\vec U_p;x,z}}
        \circ r_{\varrho}\pare{\vec U_p;z}\in
        \KR^{\varrho,m+n}_{\loc,\geq q_1+q_2},
        \qquad
        r_{\varrho}\pare{\vec U_p;z}
        \circ
        \OpBW{k\pare{\vec U_p;x,z}}
        \in
        \KR^{\varrho,m+n}_{\loc,\geq q_1+q_2}.
    \end{equation}

    \item\label{item:zcomp3}
    If
    \[
        r_{\varrho}\in\KR^{\varrho,n}_{\loc,\geq q_1},
        \qquad
        r_{\varrho'}\in\KR^{\varrho',n'}_{\loc,\geq q_2},
    \]
    then
    \begin{equation}\label{eq:zcomp3-local}
        r_{\varrho}\pare{\vec U_p;z}
        \circ
        r_{\varrho'}\pare{\vec U_p;z}
        \in
        \KR^{\varrho+\varrho',n+n'}_{\loc,\geq q_1+q_2}.
    \end{equation}
\end{enumerate}
The same conclusions hold with the local classes replaced by their
compactly supported counterparts.
\end{lemma}

\begin{proof}
For \cref{item:zcomp1}, apply
\cref{prop:composition_paradiff} for each $z\ne0$.
Choose a fixed $S$ large enough for the prescribed gain $\varrho$.
By \cref{eq:kernelfunctionproperty-local-tame}, increasing the baseline
to absorb the finite kernel losses gives
\[
    \norm{k\pare{U;\bullet,z}}_{H^S}
    \norm{k'\pare{U;\bullet,z}}_{H^S}
    \lesssim
    \norm U_{s_0}^{q_1+q_2}\av z^{m+m'}.
\]
This proves \cref{eq:zcomp1-local}.
The order-zero paraproduct bound and
\cref{eq:kernel-smoothing-local-tame-geq-q} give
\cref{item:zcomp2}. Applying the two kernel-smoothing estimates
successively gives \cref{item:zcomp3}; the smoothing, local, and
vanishing orders add. Compact support is preserved in each case.
\end{proof}

\begin{lemma}[Integration of kernel-smoothing maps]
\label{lem:integration_KR}
Let $\varrho\in\bR$, $q\in\bN$, and let
$r_{\varrho}\pare{\vec U_p;z}$ be a $z$-dependent family of operators. Assume that,
for some $n>-\mathtt d$, either
\begin{align}
     r_{\varrho}\in\KR^{\varrho,n}_{\geq q}
     &&
     \text{or}
     &&   r_{\varrho}=\chi\pare z r_0,
    \qquad
    r_0\in\KR^{\varrho,n}_{\loc,\geq q}.
\end{align}
Then
\begin{equation}\label{eq:integration_KR}
    \int_{\bR^{\mathtt d}}
    r_{\varrho}\pare{\vec U_p;z}\dd z
    \in
    \cR^\varrho_{\geq q},
    \qquad
    \int_{\bR^{\mathtt d}}
    r_{\varrho}\pare{\vec U_p;z}
    \bra{h\pare{\cdot-z}}
    \dd z
    =
    \widetilde r_{\varrho}\pare{\vec U_p}h,
    \qquad
    \widetilde r_{\varrho}\in\cR^\varrho_{\geq q}.
\end{equation}
\end{lemma}

\begin{proof}
By \cref{rem:kernel-smoothing-localization}, it is enough to consider
$r\in\KR^{\varrho,n}_{\geq q}$. Let $s_0$ be the baseline index in
\cref{def:kernel-global-smoothing}. Minkowski's inequality and
\cref{eq:kernel-smoothing-tame-geq-q} give
\[
\begin{aligned}
    \norm{
        \int_{\bR^{\mathtt d}}
        r_{\varrho}\pare{\vec U_p;z}h
        \dd z
    }_{s,I}
    \lesssim_{s,\varrho}{}&
    \norm{\vec U_p}_{s_0,I}^{q}
    \norm h_{H^{s-\varrho}}
    \int_{\av z\le2}\av z^n\dd z .
\end{aligned}
\]
The last integral is finite because $n>-\mathtt d$, and this is exactly
the tame estimate defining $\cR^\varrho_{\geq q}$.

For the translated operator, the same argument applies, since $\norm{h\pare{\cdot-z}}_{H^s}
    =
    \norm h_{H^s}$ for every $s\in\bR$. The vanishing order is unchanged by integration.
\end{proof}

\paragraph{Integration of $z$-dependent paradifferential operators.}
\begin{prop}[Kernel-to-symbol reduction for compactly supported kernels]
\label{prop:int-opbw}
Let $n>-\mathtt d$, $q\in\bN^\star$, %$\varrho>0$, 
and let
$k\in\KF^n_{\geq q}$ be a scalar kernel. Let $\sigma \in \bR$ and
$\omega \in H^\sigma\pare{\bR^{\mathtt d};\bC}$, and define
\begin{equation}\label{eq:cTK_def_recall}
    \pare{\tau_k^0\pare u\omega}\pare x
    \defeq
    \int_{\bR^{\mathtt d}}
    \OpBW{k\pare{u;\cdot,z}}
    \bra{\omega\pare{\cdot-z}}\pare x
    \dd z .
\end{equation}
Then there exists $ a_k\in\Gamma^{-\pare{n+\mathtt d}}_{\geq q}$ such that
\begin{equation}\label{eq:cTK_paralin}
    \tau_k^0\pare u\omega
    =
    \OpBW{a_k\pare{u;x,\xi}}\omega.
\end{equation}
If $k$ is real-valued for real-valued $u$, then
\begin{equation}\label{eq:int-opbw-reality}
    a_k\pare{u;x,-\xi}
    =
    \overline{a_k\pare{u;x,\xi}}.
\end{equation}
\end{prop}

\begin{proof}
By density it is enough to argue for smooth \(u\) and \(\omega\).

\begin{step}[Exact Bony--Weyl representation]
Define the centered symbol
\begin{equation}\label{eq:int-opbw-centered-symbol}
    a_k\pare{u;x,\xi}
    \defeq
    \int_{\bR^{\mathtt d}}
    k\pare{u;x+\frac z2,z}
    e^{-\ii z\cdot\xi}
    \dd z .
\end{equation}
Since \(n>-\mathtt d\) and \(k\) is compactly supported in \(z\), the
integral is well-defined. Taking the Fourier transform in \(x\) gives
\begin{equation}\label{eq:int-opbw-centered-symbol-hat}
    \widehat a_k\pare{u;\varsigma,\xi}
    =
    \int_{\bR^{\mathtt d}}
    \widehat k\pare{u;\varsigma,z}
    e^{-\ii\pare{\xi-\frac{\varsigma}{2}}\cdot z}
    \dd z .
\end{equation}
On the other hand, by the explicit Bony--Weyl formula
\cref{BonyWeyl},
\begin{equation}\label{eq:int-opbw-fourier-representation}
\begin{aligned}
    \mathcal F_{x\to\nu}
    \pare{\tau_k^0\pare u\omega}\pare\nu
    ={}&
    \frac{1}{\pare{2\pi}^{\mathtt d/2}}
    \int_{\bR^{\mathtt d}}
    \chi_\epsilon\!\pare{
        \frac{2\pare{\nu-\nu'}}{\av{\nu+\nu'}}
    }
    \int_{\bR^{\mathtt d}}
    \widehat k\pare{u;\nu-\nu',z}
    e^{-\ii\nu'\cdot z}
    \dd z\,
    \widehat\omega\pare{\nu'}
    \dd\nu' .
\end{aligned}
\end{equation}
Setting $  \varsigma\defeq\nu-\nu',
    $ and $
    \xi\defeq\frac{\nu+\nu'}{2}$
one has
\(
\nu'=\xi-\varsigma/2
\).
Hence \cref{eq:int-opbw-centered-symbol-hat} and
\cref{eq:int-opbw-fourier-representation} coincide exactly with the Fourier
representation of \(\OpBW{a_k\pare u}\). Therefore
\begin{equation}\label{eq:int-opbw-exact-centered-identity}
    \tau_k^0\pare u
    =
    \OpBW{a_k\pare{u;x,\xi}}.
\end{equation}
It remains only to prove that
\(
a_k\in\Gamma^{-\pare{n+\mathtt d}}_{\geq q}
\).
\end{step}

\begin{step}[Symbol estimate]
\label{step:symbol_class}
Fix \(s\geq s_0^k\) and \(M\in\bN\). Choose
\(N_\ast\in\bN\) such that
\begin{equation}\label{eq:int-opbw-centered-Nstar}
    N_\ast>n+M+\mathtt d,
\end{equation}
and set
\begin{equation}\label{eq:int-opbw-centered-loss}
     \tilde\mu_M^k
    \defeq
    \mu_{N_\ast}^k+N_\ast .
\end{equation}
Introduce the centered kernel
\begin{equation}\label{eq:int-opbw-centered-kernel}
    K_u\pare{x,z}
    \defeq
    k\pare{u;x+\frac z2,z}.
\end{equation}
For every multi-index \(\alpha\) with
\(\av\alpha\leq N_\ast\), the chain rule gives
\[
    \partial_z^\alpha K_u\pare{x,z}
    =
    \sum_{\gamma\leq\alpha}
    c_{\alpha,\gamma}
    \pare{
        \partial_x^{\alpha-\gamma}
        \partial_z^\gamma
        k
    }
    \pare{u;x+\frac z2,z},\qquad c_{\alpha,\gamma}\in\mathbb{R}.
\]
By translation invariance of Sobolev norms and
\cref{eq:kernelfunctionproperty-tame},
\begin{equation}\label{eq:int-opbw-centered-kernel-estimate}
    \norm{
        \partial_z^\alpha K_u\pare{\cdot,z}
    }_{H^s}
    \lesssim_{s,N_\ast}
    \norm{u}_{H^{s_0^k+ \tilde\mu_M^k}}^{q-1}
    \norm{u}_{H^{s+ \tilde\mu_M^k}}\,
    \av z^{n-\av\alpha}
    \mathbb 1_{\set{0<\av z\leq2}} .
\end{equation}
Let now \(\beta\in\bN^{\mathtt d}\), \(\av\beta\leq M\). From
\cref{eq:int-opbw-centered-symbol},
\begin{equation}\label{eq:int-opbw-centered-xi-derivative}
    \partial_\xi^\beta a_k\pare{u;x,\xi}
    =
    \pare{-\ii}^{\av\beta}
    \int_{\bR^{\mathtt d}}
    z^\beta K_u\pare{x,z}
    e^{-\ii z\cdot\xi}
    \dd z .
\end{equation}
For \(1/2\leq\av\xi\leq2\), the estimate
\begin{equation}\label{eq:int-opbw-centered-bounded-frequency}
    \norm{
        \partial_\xi^\beta a_k\pare{u;\cdot,\xi}
    }_{H^s}
    \lesssim_{s,M}
    \norm{u}_{H^{s_0^k+ \tilde\mu_M^k}}^{q-1}
    \norm{u}_{H^{s+ \tilde\mu_M^k}}
\end{equation}
follows directly from
\cref{eq:int-opbw-centered-kernel-estimate} with \(\alpha=0\), since
\(n>-\mathtt d\).

Assume henceforth that \(\av\xi\geq2\). Let
\(\chi_1\defeq\chi\) and \(\chi_2\defeq1-\chi\). We split
\cref{eq:int-opbw-centered-xi-derivative} as
\begin{equation}\label{eq:int-opbw-centered-close-far-splitting}
    \partial_\xi^\beta a_k
    =
    a_{k,\beta}^{0}
    +
    a_{k,\beta}^{\infty},
\end{equation}
where
\begin{equation}\label{eq:int-opbw-centered-close-far-definitions}
\begin{aligned}
    a_{k,\beta}^{0}\pare{u;x,\xi}
    &\defeq
    \pare{-\ii}^{\av\beta}
    \int_{\bR^{\mathtt d}}
    \chi_1\pare{\av\xi z}
    z^\beta K_u\pare{x,z}
    e^{-\ii z\cdot\xi}
    \dd z,
    \\
    a_{k,\beta}^{\infty}\pare{u;x,\xi}
    &\defeq
    \pare{-\ii}^{\av\beta}
    \int_{\bR^{\mathtt d}}
    \chi_2\pare{\av\xi z}
    z^\beta K_u\pare{x,z}
    e^{-\ii z\cdot\xi}
    \dd z.
\end{aligned}
\end{equation}

\begin{substep}[Close field]
\label{substep:near}
On the support of
\(\chi_1\pare{\av\xi z}\) one has
\begin{equation}\label{eq:int-opbw-centered-close-support}
    \av z\lesssim\av\xi^{-1}.
\end{equation}
Therefore, by Minkowski's inequality and
\cref{eq:int-opbw-centered-kernel-estimate},
\begin{equation}\label{eq:int-opbw-centered-close-estimate}
\begin{aligned}
    \norm{
        a_{k,\beta}^{0}
        \pare{u;\cdot,\xi}
    }_{H^s}
    &\lesssim
    \norm{u}_{H^{s_0^k+ \tilde\mu_M^k}}^{q-1}
    \norm{u}_{H^{s+ \tilde\mu_M^k}}
    \int_0^{c\av\xi^{-1}}
    r^{n+\av\beta+\mathtt d-1}
    \dd r
    \\
    &\lesssim
    \norm{u}_{H^{s_0^k+ \tilde\mu_M^k}}^{q-1}
    \norm{u}_{H^{s+ \tilde\mu_M^k}}
    \av\xi^{-\pare{n+\mathtt d}-\av\beta}.
\end{aligned}
\end{equation}
The radial integral is finite precisely because
\(n>-\mathtt d\).
\end{substep}

\begin{substep}[Far field]
\label{substep:far}
On the support of
\(\chi_2\pare{\av\xi z}\) one has
\begin{equation}\label{eq:int-opbw-centered-far-support}
    \av z\gtrsim\av\xi^{-1}.
\end{equation}
Recalling that $\hat{\xi}\defeq\frac{\xi}{\av\xi},$ observe that
\[
    e^{-\ii z\cdot\xi}
    =
    \pare{-\ii\av\xi}^{-N_*}
    \partial_{\hat{\xi}}^{N_*}
    e^{-\ii z\cdot\xi}.
\]
Therefore, integration by parts \(N_*\) times yields
\begin{equation}\label{eq:int-opbw-centered-far-ibp}
\begin{aligned}
    a_{k,\beta}^{\infty}
    \pare{u;x,\xi}
    =
    c_{\beta,N_*}
    \av\xi^{-N_*}
    \int_{\bR^{\mathtt d}}
    \partial_{\hat{\xi}}^{N_*}
    \left(
        \chi_2\pare{\av\xi z}
        z^\beta K_u\pare{x,z}
    \right)
    e^{-\ii z\cdot\xi}
    \dd z .
\end{aligned}
\end{equation}
There are no boundary contributions: the cutoff
\(\chi_2\pare{\av\xi z}\) removes a neighborhood of \(z=0\), while
\(K_u\) is compactly supported in \(\set{\av z\leq2}\).

Expanding the derivative in
\cref{eq:int-opbw-centered-far-ibp}, every term contains derivatives of
the cutoff, of \(z^\beta\), and of \(K_u\). Derivatives of the cutoff are
supported where
\(\av z\sim\av\xi^{-1}\), and hence satisfy
\[
    \av{
        \partial_z^\gamma
        \chi_2\pare{\av\xi z}
    }
    \lesssim_\gamma
    \av z^{-\av\gamma}.
\]
Combining this bound with
\cref{eq:int-opbw-centered-kernel-estimate} gives
\begin{equation}\label{eq:int-opbw-centered-far-amplitude}
    \norm{
        \partial_{\hat{\xi}}^{N_*}
        \left(
            \chi_2\pare{\av\xi z}
            z^\beta K_u\pare{\cdot,z}
        \right)
    }_{H^s}
    \lesssim_{s,M,N_\ast}
    \norm{u}_{H^{s_0^k+ \tilde\mu_M^k}}^{q-1}
    \norm{u}_{H^{s+ \tilde\mu_M^k}}
    \av z^{n+\av\beta-N_*}
    \mathbb 1_{\set{
        c\av\xi^{-1}\leq\av z\leq2
    }} .
\end{equation}
Consequently,
\begin{equation}\label{eq:int-opbw-centered-far-estimate}
\begin{aligned}
    \norm{
        a_{k,\beta}^{\infty}
        \pare{u;\cdot,\xi}
    }_{H^s}
    &\lesssim
    \norm{u}_{H^{s_0^k+ \tilde\mu_M^k}}^{q-1}
    \norm{u}_{H^{s+ \tilde\mu_M^k}}
    \av\xi^{-N_*}
    \int_{c\av\xi^{-1}}^2
    r^{n+\av\beta-N_*+\mathtt d-1}
    \dd r
    \\
    &\lesssim
    \norm{u}_{H^{s_0^k+ \tilde\mu_M^k}}^{q-1}
    \norm{u}_{H^{s+ \tilde\mu_M^k}}
    \av\xi^{-\pare{n+\mathtt d}-\av\beta},
\end{aligned}
\end{equation}
where in the last inequality we used
\[
    N_\ast>n+M+\mathtt d
    \geq
    n+\av\beta+\mathtt d.
\]
\end{substep}

Combining
\cref{eq:int-opbw-centered-bounded-frequency,%
eq:int-opbw-centered-close-estimate,%
eq:int-opbw-centered-far-estimate},
we obtain, for every
\(\av\beta\leq M\) and every \(\av\xi\geq1/2\),
\begin{equation}\label{eq:int-opbw-centered-final-symbol-estimate}
    \norm{
        \partial_\xi^\beta
        a_k\pare{u;\cdot,\xi}
    }_{H^s}
    \lesssim_{s,M}
    \norm{u}_{H^{s_0^k+ \tilde\mu_M^k}}^{q-1}
    \norm{u}_{H^{s+ \tilde\mu_M^k}}
    \av\xi^{-\pare{n+\mathtt d}-\av\beta}.
\end{equation}
Hence
\begin{equation}\label{eq:int-opbw-centered-symbol-class}
    a_k\in
    \Gamma^{-\pare{n+\mathtt d}}_{\geq q}.
\end{equation}
\end{step}

\begin{step}[Reality]
Finally, if \(k\) is real-valued for real-valued \(u\), then
\cref{eq:int-opbw-centered-symbol} gives
\[
    a_k\pare{u;x,-\xi}
    =
    \overline{
        a_k\pare{u;x,\xi}
    }.
\]
\end{step}
\end{proof}

\begin{rem}[Kernel and symbol regularity]
\label{rem:kernel-symbol-regularity}
The indices $N$ in \cref{def:kernel-global-global} and $M$ in
\cref{def:sfr} are finite-seminorm indices and should not be interpreted
as maximal differentiability orders. More precisely, there are two distinct
mechanisms in \cref{prop:int-opbw}.

First, since $n>-\mathtt d$ and the kernel is compactly supported in $z$,
\cref{eq:kernelfunctionproperty-local-tame} with $\alpha=0$ implies that
the centered kernel in \cref{eq:int-opbw-centered-kernel} is integrable in
$z$ with values in $H^s$. Its compact support also implies that all its
$z$-moments are integrable. Hence
\cref{eq:int-opbw-centered-xi-derivative} is well defined for every
multi-index $\beta$, and the symbol in
\cref{eq:int-opbw-centered-symbol} is
$\cC^\infty$ in $\xi$. No $z$-derivatives of the kernel are needed for
this conclusion.

The $z$-regularity enters instead in the quantitative symbol estimates.
To control the symbol seminorm of index $M$, the proof chooses a finite
number $N_\ast$ of $z$-derivatives according to
\cref{eq:int-opbw-centered-Nstar} and integrates by parts in the far
region. This yields the decay required for all frequency derivatives
$\av\beta\le M$. Because the Bony--Weyl symbol involves the centered
kernel $k\pare{u;x+z/2,z}$, differentiation in $z$ also produces
$x$-derivatives, which accounts for the loss in
\cref{eq:int-opbw-centered-loss}.

Thus every finite symbol seminorm requires only finitely many kernel
seminorms, with losses and radii allowed to depend on their respective
indices. This is precisely the reason for the quantifiers in
\cref{def:sfr,def:kernel-global-global}: $\cC^\infty$ regularity in $\xi$
is compatible with finite Sobolev regularity of the background, while
higher quantitative symbol seminorms may require correspondingly higher
Sobolev regularity. 
\fine
\end{rem}

\paragraph{Far-field tails.}

For $q\in\bN^\star$ and $k\in\KF_{\textnormal{tail},\ge q}$, set
\begin{equation}\label{eq:tail-kernel-operator}
    \tau_k^\infty\pare u\omega\pare x
    \defeq
    \int_{\bR^{\mathtt d}}k\pare{u;x,z}\omega\pare{x-z}\dd z.
\end{equation}

\begin{prop}[Far-field tail paralinearization]
\label{prop:tail-kernel-smoothing}
Let \(q\in\bN^\star\), \(k\in\KF_{\textnormal{tail},\ge q}\), and
\(\varrho\ge0\). There exists $s_0>0$ such that, for every $s\ge s_0$,
there exists $r_s>0$ for which, for every
$u\in\BallR{s_0}{r_s}\cap\cC\pare{I;H^s}$ and $\omega\in H^s$,
\begin{equation}\label{eq:tail-kernel-paralinearization}
    \tau_k^\infty\pare u\omega
    =
    \OpBW{
        f_{k,\ge q}^\infty\bra\omega\pare{u;x}
    }u
    +
    r_{k,\bra{\varrho,\ge q}}^\infty\pare u\omega
    +
    r_{k,\bra{\varrho,\ge q}}^\infty
    \bra\omega\pare u u .
\end{equation}
The objects in \cref{eq:tail-kernel-paralinearization} are as follows.
\begin{itemize}
    \item The function $f_{k,\ge q}^\infty$ satisfies
    \[
        f_{k,\ge q}^\infty\bra\cdot\pare u
        -f_{k,\ge q}^\infty\bra\cdot\pare0
        \in
        \cL\pare{H^{s_0};\cF_{\ge\max\set{1,q-1}}}
    \]
    and is defined by
    \begin{equation}\label{eq:tail-kernel-function}
        f_{k,\ge q}^\infty
        \bra\omega\pare{u;x}
        \defeq
        \int_{\bR^{\mathtt d}}
        \pare{1-\chi\pare z}
        \frac{b\pare z}{\av z}
        \mathsf g'\pare{\Delta_z u\pare x}
        \omega\pare{x-z}\dd z .
    \end{equation}

    \item $ r_{k,\bra{\varrho,\ge q}}^\infty
        \in
        \cR^\varrho_{\ge q}$, 
    \item $  r_{k,\bra{\varrho,\ge q}}^\infty
        \bra\cdot
        \in
        \cL\pare{
            H^{s_0};
            \cR^\varrho_{\ge q-1}
        }$. 
\end{itemize}
If \(k\) is real-valued for real-valued \(u\), all the objects above can
be chosen real-to-real.
\end{prop}

\begin{proof}
Fix $\varrho\ge0$ and choose $s_0$ sufficiently large, in particular
\[
    s_0>\frac{\mathtt d}{2}+\varrho.
\]
By the definition of $\KF_{\textnormal{tail},\ge q}$, there exist
$m>0$, $b$, and $\mathsf g$ satisfying
\cref{eq:tail-kernel-profile-decay,eq:tail-kernel-analytic-vanishing},
with $m+q>\mathtt d$, such that
\[
    k\pare{u;x,z}
    =
    \pare{1-\chi\pare z}
    b\pare z\,
    \mathsf g\pare{\Delta_z u\pare x}.
\]
For $u\in\BallR{s_0}{r_0}$ with $r_0>0$ sufficiently small, the Taylor expansion
\[
    \mathsf g\pare X
    =
    \sum_{j=q}^{+\infty}g_jX^j
\]
converges uniformly for
$X=\Delta_z u\pare x$, $z\in\supp\pare{1-\chi}$.
By \cref{eq:toolkit-finite-differences-Rd},
\[
    \pare{\Delta_z u\pare x}^{j}
    =
    \av z^{-j}
    \sum_{a=0}^{j}
    \pare{-1}^{a}
    \binom ja
    u^{j-a}\pare x
    u^{a}\pare{x-z}.
\]
Therefore
\begin{equation}\label{eq:tail-kernel-series-representation}
    \tau_k^\infty\pare u\omega
    =
    \sum_{j=q}^{+\infty}
    g_j
    \sum_{a=0}^{j}
    \pare{-1}^{a}
    \binom ja
    u^{j-a}
    \bigl[
        c_j*\pare{u^a\omega}
    \bigr],
\end{equation}
where
\[
    c_j\pare z
    \defeq
    \pare{1-\chi\pare z}
    b\pare z\,\av z^{-j}.
\]

For every $M\in\bN$, Leibniz' rule and
\cref{eq:tail-kernel-profile-decay} give
\[
    \av{
        \partial_z^\alpha c_j\pare z
    }
    \lesssim_M
    \pare{1+j}^{M}
    \av z^{-m-j-\av\alpha},
    \qquad
    \av\alpha\le M,
    \qquad
    \av z\ge1.9.
\]
The derivatives falling on $\chi$ are supported in the fixed annulus
$1.1\le\av z\le1.9$. Since
$m+j\ge m+q>\mathtt d$, it follows that
\begin{equation}\label{eq:tail-cj-WM1}
    \norm{c_j}_{W^{M,1}}
    \lesssim_M
    \pare{1+j}^{M}.
\end{equation}
{For every $\nu\ge0$, taking $M>\nu$ and using \cref{eq:tail-cj-WM1} give}
\begin{equation}\label{eq:tail-cj-Fourier-decay}
    \sup_{\xi\in\bR^{\mathtt d}}
    \angles{\xi}^{\nu}
    \av{\widehat c_j\pare\xi}
    \lesssim_{\nu,M}
    \pare{1+j}^{M},
\end{equation}
and consequently
\begin{equation}\label{eq:convolution-bound}
    \norm{c_j*f}_{H^{r+\nu}}
    \lesssim_{r,\nu,M}
    \pare{1+j}^{M}
    \norm f_{H^r},
    \qquad
    r\in\bR.
\end{equation}
In particular, convolution with $c_j$ gains an arbitrary number of
derivatives.

For $j\ge q$ and $0\le a\le j$, set
\[
    F_{j,a}\pare{u,\omega}
    \defeq
    c_j*\pare{u^a\omega},
    \qquad
    p\defeq j-a.
\]
By \cref{eq:convolution-bound}, the algebra property of
$H^{s_0-\varrho}$, and $s_0-\varrho>\mathtt d/2$, for every fixed
$N\in\bR$ one has, after increasing $M$,
\begin{equation}\label{eq:tail-kernel-Fja-smoothing}
    \norm{
        F_{j,a}\pare{u,h}
    }_{H^N}
    \lesssim_{N,\varrho,M}
    C^j\pare{1+j}^{M}
    \norm u_{H^{s_0}}^a
    \norm h_{H^{s_0-\varrho}}.
\end{equation}

If $p=0$, \cref{eq:tail-kernel-Fja-smoothing} directly gives
\[
    F_{j,j}\pare{u,\omega}
    =
    r_{j,j,\varrho}\pare u\omega,
\]
with
\begin{equation}\label{eq:tail-kernel-pzero-remainder}
    \norm{
        r_{j,j,\varrho}\pare u h
    }_{H^s}
    \lesssim_{s,\varrho}
    C^j\pare{1+j}^{M_s}
    \norm u_{H^{s_0}}^j
    \norm h_{H^{s-\varrho}},
    \qquad
    s\ge s_0.
\end{equation}

Let now $p\ge1$. By \cref{lem:paraprod_Bony},
\begin{equation}\label{eq:tail-kernel-Bony-product}
\begin{aligned}
    u^pF_{j,a}\pare{u,\omega}
    ={}&
    \OpBW{u^p}F_{j,a}\pare{u,\omega}
    +
    \OpBW{F_{j,a}\pare{u,\omega}}u^p
    \\
    &+
    \RBW{u^p}{F_{j,a}\pare{u,\omega}}.
\end{aligned}
\end{equation}
The first and third terms in
\cref{eq:tail-kernel-Bony-product} define an operator
$r_{j,a,\varrho}\pare u$ satisfying
\begin{equation}\label{eq:tail-kernel-first-remainder-estimate}
\begin{aligned}
    r_{j,a,\varrho}\pare u h
    \defeq{}&
    \OpBW{u^p}F_{j,a}\pare{u,h}
    +
    \RBW{u^p}{F_{j,a}\pare{u,h}},
    \\
    \norm{
        r_{j,a,\varrho}\pare u h
    }_{H^s}
    \lesssim_{s,\varrho}{}&
    C^j\pare{1+j}^{M_s}
    \norm u_{H^{s_0}}^j
    \norm h_{H^{s-\varrho}}.
\end{aligned}
\end{equation}
Indeed, the paraproduct estimate only uses
$\norm{u^p}_{H^{s_0}}$, whereas in the Bony remainder the additional
regularity is placed on $F_{j,a}$ and follows from
\cref{eq:tail-kernel-Fja-smoothing}.

For the second term in \cref{eq:tail-kernel-Bony-product},
\cref{lem paralin}, applied to $v\mapsto v^p$, gives
\begin{equation}\label{eq:tail-kernel-power-paralinearization}
    u^p
    =
    \OpBW{p\,u^{p-1}}u
    +
    r_{p,\varrho}\pare u u,
\end{equation}
where, by increasing $s_0$ with respect to $\varrho$ if necessary,
\begin{equation}\label{eq:tail-kernel-power-remainder-estimate}
    \norm{
        r_{p,\varrho}\pare u h
    }_{H^s}
    \lesssim_{s,\varrho}
    C^p
    \norm u_{H^{s_0}}^{p-1}
    \norm h_{H^{s-\varrho}}.
\end{equation}
Hence
\begin{align}
    \OpBW{F_{j,a}\pare{u,\omega}}u^p
    ={}&
    \OpBW{F_{j,a}\pare{u,\omega}}
    \OpBW{p\,u^{p-1}}u
    \notag\\
    &+
    \OpBW{F_{j,a}\pare{u,\omega}}
    r_{p,\varrho}\pare u u .
    \label{eq:tail-kernel-Bony-transported}
\end{align}

For $p=1$, take $r_{p,\varrho}=0$. For $p\ge2$, the estimate
\cref{eq:tail-kernel-power-remainder-estimate} follows by induction:
apply \cref{lem:paraprod_Bony} to $u^{p-1}u$, then use
\cref{rem:smoothBONYremainder,prop:composition_paradiff}.
Applying \cref{prop:composition_paradiff} to the
functions $F_{j,a}\pare{u,\omega}$ and $p\,u^{p-1}$ gives
\begin{equation}\label{eq:tail-kernel-Bony-composition}
\begin{aligned}
    \OpBW{F_{j,a}\pare{u,\omega}}
    \OpBW{p\,u^{p-1}}
    ={}&
    \OpBW{
        p\,u^{p-1}F_{j,a}\pare{u,\omega}
    }
    +
    r_{j,a,\varrho}^{\mathrm{comp}}
    \bra\omega\pare u ,
\end{aligned}
\end{equation}
where
\begin{equation}\label{eq:tail-kernel-composition-remainder-estimate}
    \norm{
        r_{j,a,\varrho}^{\mathrm{comp}}
        \bra\omega\pare u h
    }_{H^s}
    \lesssim_{s,\varrho}
    C^j\pare{1+j}^{M_s}
    \norm u_{H^{s_0}}^{j-1}
    \norm\omega_{H^{s_0}}
    \norm h_{H^{s-\varrho}}.
\end{equation}
Together with
\cref{eq:tail-kernel-Fja-smoothing,%
eq:tail-kernel-power-remainder-estimate},
the second term in
\cref{eq:tail-kernel-Bony-transported} satisfies the same estimate.
Thus, setting
\[
    \widetilde r_{j,a,\varrho}
    \bra\omega\pare u
    \defeq
    \OpBW{F_{j,a}\pare{u,\omega}}
    r_{p,\varrho}\pare u
    +
    r_{j,a,\varrho}^{\mathrm{comp}}
    \bra\omega\pare u ,
\]
we obtain
\begin{equation}\label{eq:tail-kernel-jth-paralinearization}
\begin{aligned}
    u^{j-a}F_{j,a}\pare{u,\omega}
    ={}&
    \OpBW{
        \pare{j-a}
        u^{j-a-1}
        F_{j,a}\pare{u,\omega}
    }u
    \\
    &+
    r_{j,a,\varrho}\pare u\omega
    +
    \widetilde r_{j,a,\varrho}
    \bra\omega\pare u u ,
\end{aligned}
\end{equation}
and
\begin{equation}\label{eq:tail-kernel-jth-smoothing-estimates}
\begin{aligned}
    \norm{
        r_{j,a,\varrho}\pare u h
    }_{H^s}
    &\lesssim_{s,\varrho}
    C^j\pare{1+j}^{M_s}
    \norm u_{H^{s_0}}^j
    \norm h_{H^{s-\varrho}},
    \\
    \norm{
        \widetilde r_{j,a,\varrho}
        \bra\omega\pare u h
    }_{H^s}
    &\lesssim_{s,\varrho}
    C^j\pare{1+j}^{M_s}
    \norm u_{H^{s_0}}^{j-1}
    \norm\omega_{H^{s_0}}
    \norm h_{H^{s-\varrho}}.
\end{aligned}
\end{equation}

It remains to identify the first term in
\cref{eq:tail-kernel-jth-paralinearization}. By
\cref{eq:tail-kernel-series-representation}, its coefficient is
\[
    \sum_{j=q}^{+\infty}
    g_j
    \sum_{a=0}^{j-1}
    \pare{-1}^{a}
    \binom ja
    \pare{j-a}
    u^{j-a-1}
    \bigl[
        c_j*\pare{u^a\omega}
    \bigr].
\]
Using
\[
    \pare{j-a}\binom ja
    =
    j\binom{j-1}{a}
\]
and the definition of $c_j$, this is equal to
\[
    {f_{k,\ge q}^\infty\bra\omega\pare{u;x}}
    =
    \int_{\bR^{\mathtt d}}
    \pare{1-\chi\pare z}
    \frac{b\pare z}{\av z}\,
    \mathsf g'\pare{\Delta_z u\pare x}
    \omega\pare{x-z}\dd z.
\]
The terms with $j\ge\max\set{q,2}$ give the asserted $\cF$ membership,
while the possible $j=1$ term is independent of $u$. Hence
$f_{k,\ge q}^\infty\bra\omega\pare u$ has total vanishing order $q$ by
\cref{rem:linear-parameter-vanishing-order}.
Finally, set
{
\[
\begin{aligned}
    r_{k,\bra{\varrho,\ge q}}^\infty\pare u
    &\defeq
    \sum_{j=q}^{+\infty}
    g_j
    \sum_{a=0}^{j}
    \pare{-1}^{a}
    \binom ja
    r_{j,a,\varrho}\pare u,
    \\
    r_{k,\bra{\varrho,\ge q}}^\infty
    \bra\omega\pare u
    &\defeq
    \sum_{j=q}^{+\infty}
    g_j
    \sum_{a=0}^{j-1}
    \pare{-1}^{a}
    \binom ja
    \widetilde r_{j,a,\varrho}
    \bra\omega\pare u .
\end{aligned}
\]
}
By \cref{eq:tail-kernel-pzero-remainder,eq:tail-kernel-jth-smoothing-estimates}
and the analyticity of $\mathsf g$, for each $s\ge s_0$ the two series
converge after reducing $r_s$. Their vanishing orders in $u$ are $q$ and $q-1$, respectively; the second
remainder is linear in $\omega$ and hence has total vanishing order $q$ by
\cref{rem:linear-parameter-vanishing-order}.
Substituting into
\cref{eq:tail-kernel-series-representation} yields the decomposition
in the statement.

If $k$ is real-valued for real-valued $u$, all the terms in the above
construction preserve real-valued functions, which proves the
real-to-real property.
\end{proof}

\paragraph{Homogeneous angular multipliers}

\begin{prop}[Homogeneous angular multipliers on $\bR^{\mathtt d}$]
\label{prop:abstract-homogeneous-angular-multipliers-Rd}
Let \(f\in\mathcal S\pare{\bR^{\mathtt d};\bC}\). The following statements hold.
\begin{enumerate}

    \item \label{item:abstract-homogeneous-angular-degree-zero}
    Let \(a\in\cC^\infty\pare{\bS^{\mathtt d-1};\bR}\) be odd. Then its
    spherical mean vanishes and there exists a real, odd Fourier multiplier
    \(\dot{\mathfrak m}_a^{\pare0}\), smooth for \(\xi\neq0\) and
    homogeneous of degree zero, explicitly given by
    \begin{equation}\label{eq:abstract-critical-angular-multiplier-Rd}
        \dot{\mathfrak m}_a^{\pare0}\pare\xi
        \defeq
        -\frac\pi2
        \int_{\bS^{\mathtt d-1}}
        \sgn\pare{\vartheta\cdot\xi}
        a\pare\vartheta\dd\sigma\pare\vartheta,
    \end{equation}
    such that
    \begin{equation}\label{eq:abstract-critical-angular-operator-Rd}
        \pv\int_{\bR^{\mathtt d}}
        \frac{a\pare{\hat z}}{\av z^{\mathtt d}}
        f\pare{\bullet-z}\dd z
        =
        \ii\dot{\mathfrak m}_a^{\pare0}\pare Df.
    \end{equation}
    Moreover, \(\ii\dot{\mathfrak m}_a^{\pare0}\pare D\) is real-to-real.
    \item \label{item:abstract-homogeneous-angular-degree-minus-one}
    Let \(a\in\cC^\infty\pare{\bS^{\mathtt d-1};\bR}\) be even and
    \(\mathtt d\geq2\). Then there exists a real, even Fourier multiplier
    \(\dot{\mathfrak m}_a^{\pare{-1}}\in\tilde{\Gamma}^{-1}\), smooth for \(\xi\neq0\) and
    homogeneous of degree \(-1\), explicitly given by
    \begin{equation}\label{eq:abstract-subprincipal-angular-multiplier-Rd}
    \begin{aligned}
        \dot{\mathfrak m}_a^{\pare{-1}}\pare\xi
        &\defeq
        \pi
        \int_{\bS^{\mathtt d-1}}
        \delta\pare{\vartheta\cdot\xi}
        a\pare\vartheta\dd\sigma\pare\vartheta
        \\
        &=
        \frac\pi{\av\xi}
        \int_{\bS^{\mathtt d-1}\cap\xi^\perp}
        a\pare\vartheta
        \dd\sigma_{\mathtt d-2}\pare\vartheta,
    \end{aligned}
    \end{equation}
    such that
    \begin{equation}\label{eq:abstract-subprincipal-angular-operator-Rd}
        \int_{\bR^{\mathtt d}}
        \frac{a\pare{\hat z}}{\av z^{\mathtt d-1}}
        f\pare{\bullet-z}\dd z
        =
        \dot{\mathfrak m}_a^{\pare{-1}}\pare Df.
    \end{equation}
    Moreover, \(\dot{\mathfrak m}_a^{\pare{-1}}\pare D\) is real-to-real.
    
    If \(\mathtt d=1\) and \(a\) is even on
    \(\bS^0=\set{-1,1}\), then \(a\pare{\hat z}\) is constant on
    \(\bR\setminus\set0\) and
    \begin{equation}\label{eq:abstract-one-dimensional-degree-zero-model}
        \supp
        \cF_{z\to\xi}
        \bra{a\pare{\hat z}}
        \subseteq
        \set0.
    \end{equation}
    In particular, its complete high-frequency homogeneous multiplier
    vanishes.

    \item \label{item:abstract-homogeneous-angular-degree-one}
    Let \(a\in\cC^\infty\pare{\bS^{\mathtt d-1};\bR}\) be even. Then
    there exists a real, even Fourier multiplier
    \(\dot{\mathfrak m}_a^{\pare1}\in\tilde{\Gamma}^1\), smooth for \(\xi\neq0\) and
    homogeneous of degree one, explicitly given by
    \begin{equation}\label{eq:abstract-compensated-angular-multiplier-Rd}
        \dot{\mathfrak m}_a^{\pare1}\pare\xi
        \defeq
        \frac\pi2
        \int_{\bS^{\mathtt d-1}}
        \av{\vartheta\cdot\xi}
        a\pare\vartheta\dd\sigma\pare\vartheta,
    \end{equation}
    such that
    \begin{equation}\label{eq:abstract-compensated-angular-operator-Rd}
        \frac12
        \int_{\bR^{\mathtt d}}
        \frac{a\pare{\hat z}}{\av z^{\mathtt d+1}}
        \delta_z^{\pare2}f\dd z
        =
        \dot{\mathfrak m}_a^{\pare1}\pare Df.
    \end{equation}
    Moreover, \(\dot{\mathfrak m}_a^{\pare1}\pare D\) is real-to-real.

\end{enumerate}

Moreover, for every \(M\in\bN\),
in each of the cases above, we have
\begin{equation}\label{eq:abstract-angular-multiplier-seminorm-bound}
    \max_{\av\beta\leq M}
    \sup_{\av\xi\geq1/2}
    \av\xi^{-j+\av\beta}
    \av{
        \partial_\xi^\beta
        \dot{\mathfrak m}_a^{\pare j}\pare\xi
    }
    \lesssim_M
    \norm{a}_{\cC^{M}\pare{\bS^{\mathtt d-1}}},
\end{equation}
where \(j\in\set{-1,0,1}\), with \(j=-1\) restricted to
\(\mathtt d\geq2\).
\end{prop}

\begin{proof}
\begin{proofpart}[Homogeneous multiplier identities]
All Fourier identities below are understood in the sense of tempered
distributions. We use throughout
\cref{eq:abstract-radial-distribution-identities}.

\begin{proofsubpart}[Proof of \cref{eq:abstract-critical-angular-operator-Rd}]
By oddness,
\[
    \int_{\bS^{\mathtt d-1}}a\dd\sigma=0.
\]
In polar coordinates \(z=r\vartheta\), pairing \(\vartheta\) with
\(-\vartheta\) and using the first identity in
\cref{eq:abstract-radial-distribution-identities}, we obtain
\[
\begin{aligned}
    \pv\int_{\bR^{\mathtt d}}
    \frac{a\pare{\hat z}}{\av z^{\mathtt d}}
    e^{-\ii z\cdot\xi}\dd z
    &=
    -\ii
    \int_{\bS^{\mathtt d-1}}
    a\pare\vartheta
    \int_0^\infty
    \frac{\sin\pare{r\vartheta\cdot\xi}}r\dd r
    \dd\sigma\pare\vartheta
    \\
    &=
    \ii\dot{\mathfrak m}_a^{\pare0}\pare\xi,
\end{aligned}
\]
where the last equality is
\cref{eq:abstract-critical-angular-multiplier-Rd}.
\end{proofsubpart}

\begin{proofsubpart}[Proof of the first identity in \cref{eq:abstract-subprincipal-angular-multiplier-Rd}]
Polar coordinates and the second identity in
\cref{eq:abstract-radial-distribution-identities} give
\begin{equation}\label{eq:abstract-subprincipal-kernel-fourier-transform}
\begin{aligned}
    \int_{\bR^{\mathtt d}}
    \frac{a\pare{\hat z}}{\av z^{\mathtt d-1}}
    e^{-\ii z\cdot\xi}\dd z
    &=
    \int_{\bS^{\mathtt d-1}}
    a\pare\vartheta
    \pare{
        \int_0^\infty
        e^{-\ii r\vartheta\cdot\xi}\dd r
    }
    \dd\sigma\pare\vartheta
    \\
    &=
    \int_{\bS^{\mathtt d-1}}
    a\pare\vartheta
    \pare{
        \pi\delta\pare{\vartheta\cdot\xi}
        -
        \ii\pv\frac1{\vartheta\cdot\xi}
    }
    \dd\sigma\pare\vartheta
    \\
    &=
    \pi
    \int_{\bS^{\mathtt d-1}}
    \delta\pare{\vartheta\cdot\xi}
    a\pare\vartheta
    \dd\sigma\pare\vartheta.
\end{aligned}
\end{equation}
The principal-value contribution vanishes under
\(\vartheta\mapsto-\vartheta\), since \(a\) is even.
\end{proofsubpart}

\begin{proofsubpart}[Proof of the second identity in \cref{eq:abstract-subprincipal-angular-multiplier-Rd} and of \cref{eq:abstract-subprincipal-angular-operator-Rd}]
For \(\mathtt d\geq2\), it remains to identify geometrically the
delta contribution in
\cref{eq:abstract-subprincipal-kernel-fourier-transform}.
Fix \(\xi\neq0\) and define
\[
    g_\xi\pare\vartheta
    \defeq
    \vartheta\cdot\xi,
    \qquad
    \vartheta\in\bS^{\mathtt d-1}.
\]
Then
\begin{equation}
    \label{eq:preimage-orthogonal}
      g_\xi^{-1}\pare0
    =
    \bS^{\mathtt d-1}\cap\angles{\xi}^\perp.
\end{equation}
Moreover, since the spherical gradient is the orthogonal projection
of the Euclidean gradient onto\(T_\vartheta\bS^{\mathtt d-1}\),
\begin{equation}
\label{eq:gradient-as-projection}
     \nabla_{\bS^{\mathtt d-1}}g_\xi\pare\vartheta
    =
    \xi-\pare{\vartheta\cdot\xi}\vartheta = \pare{\Id -\vartheta\otimes\vartheta}\xi.
\end{equation}
Let $\bar\vartheta \in g_\xi^{-1}\pare0 $, by \cref{eq:preimage-orthogonal,eq:gradient-as-projection} we have that $\bar\vartheta\cdot\xi =0$ and $\nabla_{\bS^{\mathtt d-1}}g_\xi\pare{\bar\vartheta}
    =
    \xi$, this implies that \(0\) is a regular value of \(g_\xi\) and, on
\(g_\xi^{-1}\pare0\) holds $ \av{
        \nabla_{\bS^{\mathtt d-1}}
        g_\xi\pare\vartheta
    }
    =
    \av\xi$.
To identify the action of the delta distribution, let
\(\varphi\in\cC_c^\infty\pare{\bR}\). For \(\mathtt d\geq2\), fix \(\xi\neq0\). 
We write $\vartheta$ via the decomposition $\angles{\hat\xi}\oplus \angles{\hat\xi}^\perp $ as
\begin{equation}
    \vartheta = \mathbb{P}_{\hat\xi} \vartheta + \mathbb{P}_{\hat\xi}^\perp \vartheta = \hat\xi\otimes\hat\xi \ \vartheta + \pare{\Id - \hat\xi\otimes\hat\xi } \vartheta = t\hat\xi +\pare{\vartheta - t\hat\xi}. 
\end{equation}
Explicit computations show that $\av{\vartheta - t\hat\xi}^2=1-t^2$. Thus defining $\omega\defeq \frac{\vartheta-t\hat\xi}{\sqrt{1-t^2}}$, we obtain the decomposition
\[
    \vartheta
    =
    t\hat\xi+\sqrt{1-t^2}\,\omega,
    \qquad
    t\in\pare{-1,1},
    \qquad
    \omega\in\bS^{\mathtt d-1}\cap\xi^\perp.
\]
Notice that, writing
\[
    \vartheta
    =
    \cos\phi\,\hat\xi+\sin\phi\,\omega,
    \qquad
    \omega\in\bS^{\mathtt d-1}\cap\xi^\perp,
\]
the standard spherical measure is
\[
    \dd\sigma\pare\vartheta
    =
    \sin^{\mathtt d-2}\phi\,
    \dd\phi\,\dd\sigma_{\mathtt d-2}\pare\omega.
\]
Since \(t=\cos\phi\), so that
\(\dd t=-\sin\phi\,\dd\phi\) and
\(\sin\phi=\sqrt{1-t^2}\), this becomes
\[
    \dd\sigma\pare\vartheta
    =
    \pare{1-t^2}^{\frac{\mathtt d-3}{2}}
    \dd t\,\dd\sigma_{\mathtt d-2}\pare\omega.
\]
which we combine with $\vartheta\cdot\xi=\av\xi\,t$ in order to obtain that
\[
\begin{aligned}
    \int_{\bS^{\mathtt d-1}}
    \delta\pare{\vartheta\cdot\xi}
    a\pare\vartheta
    \dd\sigma\pare\vartheta
    &=
    \int_{-1}^{1}
    \delta\pare{\av\xi t}
    \pare{1-t^2}^{\frac{\mathtt d-3}{2}}
    \int_{\bS^{\mathtt d-1}\cap\xi^\perp}
    a\pare{
        t\hat\xi+\sqrt{1-t^2}\,\omega
    }
    \dd\sigma_{\mathtt d-2}\pare\omega
    \dd t
    \\
    &=
    \frac1{\av\xi}
    \int_{\bS^{\mathtt d-1}\cap\xi^\perp}
    a\pare\omega
    \dd\sigma_{\mathtt d-2}\pare\omega,
\end{aligned}
\]
where in the last equality we used
\[
    \delta\pare{\av\xi t}
    =
    \frac1{\av\xi}\delta\pare t.
\]
Together with
\cref{eq:abstract-subprincipal-kernel-fourier-transform},
it also gives
\cref{eq:abstract-subprincipal-angular-operator-Rd}.
\end{proofsubpart}

\begin{proofsubpart}[Proof of \cref{eq:abstract-one-dimensional-degree-zero-model}]
Let \(\mathtt d=1\). Since
\(\bS^0=\set{-1,1}\) and \(a\) is even, there is \(c_a\in\bR\) such that
\[
    a\pare{\hat z}=c_a,
    \qquad
    z\in\bR\setminus\set0.
\]
The value at \(z=0\) is irrelevant as a tempered distribution, hence
\(a\pare{\hat z}\) is the constant distribution \(c_a\). Therefore
\[
    \supp
    \cF_{z\to\xi}
    \bra{
        a\pare{\hat z}
    }
    \subseteq
    \set0.
\]
\end{proofsubpart}

\begin{proofsubpart}[Proof of \cref{eq:abstract-compensated-angular-operator-Rd}]
By \cref{eq:localized-compensated-second-difference},
\[
    \cF
    \bra{
        \frac12\delta_z^{\pare2}f
    }\pare\xi
    =
    \pare{
        1-\cos\pare{z\cdot\xi}
    }
    \widehat f\pare\xi.
\]
Hence, by polar coordinates and the third identity in
\cref{eq:abstract-radial-distribution-identities},
\[
\begin{aligned}
    \int_{\bR^{\mathtt d}}
    \frac{a\pare{\hat z}}{\av z^{\mathtt d+1}}
    \pare{
        1-\cos\pare{z\cdot\xi}
    }
    \dd z
    &=
    \int_{\bS^{\mathtt d-1}}
    a\pare\vartheta
    \int_0^\infty
    \frac{
        1-\cos\pare{r\vartheta\cdot\xi}
    }{r^2}
    \dd r
    \dd\sigma\pare\vartheta
    \\
    &=
    \dot{\mathfrak m}_a^{\pare1}\pare\xi,
\end{aligned}
\]
where the last equality is
\cref{eq:abstract-compensated-angular-multiplier-Rd}.
\end{proofsubpart}
\end{proofpart}

\begin{proofpart}[Symbol properties and seminorms]
The formulas in
\cref{eq:abstract-critical-angular-multiplier-Rd,%
eq:abstract-subprincipal-angular-multiplier-Rd,%
eq:abstract-compensated-angular-multiplier-Rd}
give the stated homogeneity, reality and parity properties.
It remains to prove smoothness away from the origin and the corresponding
symbol estimates.

By homogeneity,
\begin{equation}\label{eq:abstract-angular-multiplier-polar-form}
    \dot{\mathfrak m}_a^{\pare j}\pare{\xi}
    =
    \av{\xi}^j
    \dot{\mathfrak m}_a^{\pare j}\pare{\hat{\xi}},
    \qquad
    j\in\set{-1,0,1},
\end{equation}
where \(j=-1\) is restricted to \(\mathtt d\geq2\).

If \(\mathtt d=1\), counting measure on
\(\bS^0=\set{-1,1}\) gives
\[
    \dot{\mathfrak m}_a^{\pare0}\pare\xi
    =
    -\pi a\pare1\sgn\pare\xi,
    \qquad
    \dot{\mathfrak m}_a^{\pare1}\pare\xi
    =
    \pi a\pare1\av\xi.
\]
Hence smoothness on
\(\bR\setminus\set0\) and
\cref{eq:abstract-angular-multiplier-seminorm-bound}
are immediate. We may therefore assume
\(\mathtt d\geq2\).

We first record how derivatives in \(\xi\) reduce to derivatives in the
angular variable. For
\(k\in\set{1,\ldots,\mathtt d}\), let \(e_k\) denote the \(k\)-th
standard basis vector of \(\bR^{\mathtt d}\), and define the smooth
tangential vector field 
\[
    X_k\pare{\hat{\xi}}
    \defeq
    e_k-\hat{\xi}_k\hat{\xi}
    \in
    T_{\hat{\xi}}\bS^{\mathtt d-1},\qquad\hat{\xi}_k\defeq\hat{\xi}\cdot e_k.
\]
Since
\[
    \partial_{\xi_k}\av{\xi}=\hat{\xi}_k,
    \qquad
    \partial_{\xi_k}\hat{\xi}
    =
    \frac{
        e_k-\hat{\xi}_k\hat{\xi}
    }{\av{\xi}}
    =
    \frac{
        X_k\pare{\hat{\xi}}
    }{\av{\xi}},
\]
the chain rule gives, for every
\(b\in\cC^\infty\pare{\bS^{\mathtt d-1}}\),
\begin{equation}\label{eq:abstract-first-xi-angular-derivative}
    \partial_{\xi_k}
    \bra{
        \av{\xi}^j b\pare{\hat{\xi}}
    }
    =
    \av{\xi}^{j-1}
    \pare{
        j\hat{\xi}_k b\pare{\hat{\xi}}
        +
        X_k b\pare{\hat{\xi}}
    }.
\end{equation}
Iterating
\cref{eq:abstract-first-xi-angular-derivative}
shows that, for every multi-index \(\beta\), there exists a tangential
differential operator $P_\beta^j$ on \(\bS^{\mathtt d-1}\), of order at most \(\av\beta\) and with smooth
coefficients, such that
\begin{equation}\label{eq:abstract-homogeneous-xi-derivative-form}
    \partial_\xi^\beta
    \dot{\mathfrak m}_a^{\pare j}\pare{\xi}
    =
    \av{\xi}^{j-\av\beta}
    P_\beta^j
    \bra{
        \dot{\mathfrak m}_a^{\pare j}
        \big|_{\bS^{\mathtt d-1}}
    }\pare{\hat{\xi}}.
\end{equation}
More explicitly, on any coordinate neighborhood
\(U\subset\bS^{\mathtt d-1}\),
\begin{equation}\label{eq:abstract-xi-angular-local-coordinate-form}
    P_\beta^j
    =
    \sum_{\av\gamma\leq\av\beta}
    c_{\beta,\gamma}^{j,U}\pare{\hat{\xi}}
    \partial_{\hat{\xi}}^\gamma,
    \qquad
    c_{\beta,\gamma}^{j,U}
    \in
    \cC^\infty\pare U.
\end{equation}
Thus
\cref{eq:abstract-homogeneous-xi-derivative-form,%
eq:abstract-xi-angular-local-coordinate-form}
reduce the required estimates in \(\xi\) to estimates for derivatives of
the restriction of
\(\dot{\mathfrak m}_a^{\pare j}\)
to the unit sphere. We now establish these angular estimates.

Fix coordinate neighborhoods
\(U\Subset V\subset\bS^{\mathtt d-1}\) and a smooth family of rotations
\[
    V\ni\hat{\xi}
    \longmapsto
    Q_{\hat{\xi}}\in\textnormal{SO}\pare{\mathtt d},
    \qquad
    Q_{\hat{\xi}} e_{\mathtt d}=\hat{\xi}.
\]
Such a family exists locally by smoothly completing \(\nu\) to a
positively oriented orthonormal basis. Since
\(\overline U\subset V\), all its coordinate derivatives are uniformly
bounded on \(U\).

For
\(\phi\in\bS^{\mathtt d-1}\), set
\[
    \phi_{\mathtt d}
    \defeq
    \phi\cdot e_{\mathtt d}.
\]
The change of variables
\(\vartheta=Q_{\hat{\xi}}\phi\) preserves surface measure and satisfies
\[
    \vartheta\cdot\hat{\xi}
    =
    \pare{Q_{\hat{\xi}}\phi}
    \cdot
    \pare{Q_{\hat{\xi}} e_{\mathtt d}}
    =
    \phi_{\mathtt d}.
\]
Consequently,
\cref{eq:abstract-critical-angular-multiplier-Rd,%
eq:abstract-subprincipal-angular-multiplier-Rd,%
eq:abstract-compensated-angular-multiplier-Rd}
give
\begin{equation}\label{eq:abstract-angular-multiplier-rotated-forms}
\begin{aligned}
    \dot{\mathfrak m}_a^{\pare0}\pare{\hat{\xi}}
    &=
    -\frac\pi2
    \int_{\bS^{\mathtt d-1}}
    \sgn\pare{\phi_{\mathtt d}}
    a\pare{Q_{\hat{\xi}}\phi}
    \dd\sigma\pare\phi,
    \\
    \dot{\mathfrak m}_a^{\pare{-1}}\pare{\hat{\xi}}
    &=
    \pi
    \int_{\bS^{\mathtt d-1}\cap e_{\mathtt d}^\perp}
    a\pare{Q_{\hat{\xi}}\phi}
    \dd\sigma_{\mathtt d-2}\pare\phi,
    \\
    \dot{\mathfrak m}_a^{\pare1}\pare{\hat{\xi}}
    &=
    \frac\pi2
    \int_{\bS^{\mathtt d-1}}
    \av{\phi_{\mathtt d}}
    a\pare{Q_{\hat{\xi}}\phi}
    \dd\sigma\pare\phi.
\end{aligned}
\end{equation}
The purpose of the rotation is that the integration domains and the
possibly nonsmooth weights are now independent of \(\hat{\xi}\). All dependence
on \(\hat{\xi}\) is contained in the smooth factor
\(a\pare{Q_{\hat{\xi}}\phi}\).

Let \(\gamma\) be a multi-index in the local coordinates of \(V\).
By the chain rule and the smoothness of
\(\hat{\xi}\mapsto Q_{\hat{\xi}}\),
\[
    \sup_{\substack{\hat{\xi}\in U\\
    \phi\in\bS^{\mathtt d-1}}}
    \av{
        \partial_{\hat{\xi}}^\gamma
        \bra{
            a\pare{Q_{\hat{\xi}}\phi}
        }
    }
    \lesssim_{\gamma,U,V}
    \norm{
        a
    }_{\cC^{\av\gamma}\pare{\bS^{\mathtt d-1}}}.
\]
Since the domains in
\cref{eq:abstract-angular-multiplier-rotated-forms}
have finite measure and their weights are bounded, differentiation under
the integral sign is therefore justified. We obtain
\begin{equation}\label{eq:abstract-angular-derivative-bound}
    \sup_{\hat{\xi}\in U}
    \av{
        \partial_{\hat{\xi}}^\gamma
        \dot{\mathfrak m}_a^{\pare j}\pare\nu
    }
    \lesssim_{\gamma,U,V}
    \norm{
        a
    }_{\cC^{\av\gamma}\pare{\bS^{\mathtt d-1}}},
    \qquad
    j\in\set{-1,0,1}.
\end{equation}
In particular, the restriction of each multiplier to \(U\) is smooth.
A finite collection of such neighborhoods covers
\(\bS^{\mathtt d-1}\), and hence
\cref{eq:abstract-angular-derivative-bound}
holds uniformly on the sphere.

We may now substitute the angular bounds into the derivative formula
obtained at the beginning of this step.
Indeed,
\cref{eq:abstract-homogeneous-xi-derivative-form,%
eq:abstract-xi-angular-local-coordinate-form,%
eq:abstract-angular-derivative-bound}
give, for every multi-index \(\beta\),
\[
    \av{
        \partial_\xi^\beta
        \dot{\mathfrak m}_a^{\pare j}\pare{\xi}
    }
    \lesssim_\beta
    \av{\xi}^{j-\av\beta}
    \norm{
        a
    }_{\cC^{\av\beta}\pare{\bS^{\mathtt d-1}}}.
\]
Therefore, for every \(M\in\bN\),
\[
    \max_{\av\beta\leq M}
    \sup_{\av\xi\geq1/2}
    \av\xi^{-j+\av\beta}
    \av{
        \partial_\xi^\beta
        \dot{\mathfrak m}_a^{\pare j}\pare\xi
    }
    \lesssim_M
    \norm{
        a
    }_{\cC^M\pare{\bS^{\mathtt d-1}}},
\]
which proves
\cref{eq:abstract-angular-multiplier-seminorm-bound}.
By \cref{def:FM},
\[
    \dot{\mathfrak m}_a^{\pare j}
    \in
    \tilde\Gamma^j.
\]

Finally, reality and parity give
\[
    \ii
    \dot{\mathfrak m}_a^{\pare0}\pare{-\xi}
    =
    \overline{
        \ii
        \dot{\mathfrak m}_a^{\pare0}\pare\xi
    },
\]
whereas, for \(j\in\set{-1,1}\),
\[
    \dot{\mathfrak m}_a^{\pare j}\pare{-\xi}
    =
    \overline{
        \dot{\mathfrak m}_a^{\pare j}\pare\xi
    },
\]
with \(j=-1\) restricted to \(\mathtt d\geq2\).
These are precisely the Fourier reality conditions, and therefore prove
the stated real-to-real properties.
\end{proofpart}
\end{proof}

\paragraph{Homogeneous convolution models and cutoff removal.} The following setting of the models corresponds to the critical regime of the graph integral operators, which will be used in the proof of \cref{thm:abstract-structure-preserving-paralinearization-Rd}.

\begin{definition}[Localized homogeneous convolution models]
\label{def:localized-homogeneous-convolution-models}
Let $\lambda\in\set{-\mathtt d,-\mathtt d+1}$ and let
$h\in\cC^\infty\pare{\bR^{\mathtt d}\setminus\set0;\bC}$ be homogeneous of degree
$\lambda$. If $\lambda=-\mathtt d$, we assume the cancellation condition
\begin{equation}\label{eq:localized-homogeneous-model-cancellation}
    \int_{\bS^{\mathtt d-1}}
    h\pare\vartheta\dd\sigma\pare\vartheta=0.
\end{equation}
We associate with $h$ the tempered distribution
\begin{equation}\label{eq:localized-homogeneous-model-distribution}
    \mathfrak h
    \defeq
    \begin{cases}
        \pv h, & \lambda=-\mathtt d,\\
        h, & \lambda=-\mathtt d+1.
    \end{cases}
\end{equation}
We define
\begin{equation}\label{eq:localized-homogeneous-model-localized-operator}
    t_h^{\loc}\omega
    \defeq
    \pare{\chi\mathfrak h}*\omega,
\end{equation}
and the corresponding complete homogeneous multiplier by
\begin{equation}\label{eq:homogeneous-model-multiplier}
    \dot{\mathfrak m}_h
    \defeq
    \pare{2\pi}^{\mathtt d/2}\mathcal F_{z\to\xi}\mathfrak h .
\end{equation}
\end{definition}
{{
\begin{remark}
   The cancellation condition \cref{eq:localized-homogeneous-model-cancellation} ensures that the convolution operator
   \[
   T_hf\p x\defeq \pv\int_{\mathbb R^{\mathtt d}} h\p yf\p{x-y}\dd y
   \]
   is a Calderon-Zygmund operator, which is $L^p$ bounded for any $1<p<\infty$.
   \fine
\end{remark}
}}

\begin{lemma}[Localized homogeneous models are Fourier multipliers modulo smoothing]
\label{lem:localized-homogeneous-model}
Let $\lambda\in\set{-\mathtt d,-\mathtt d+1}$ and let $h$ satisfy
\cref{def:localized-homogeneous-convolution-models}. Then,
$$\dot{\mathfrak m}_h\in\tilde\Gamma^{-\pare{\lambda+\mathtt{d}}}$$
and
\begin{equation}\label{eq:localized-homogeneous-model-comparison}
    t_h^{\loc}-\dot{\mathfrak m}_h \pare{D}
    \in
    \bigcap_{M\in\bN}
    \tilde\Gamma^{-M}.
\end{equation}
Moreover, if $h$ is real-valued, then both terms in
\cref{eq:localized-homogeneous-model-comparison} are real-to-real.
\end{lemma}

\begin{proof}
Let
\[
    \mathfrak m_{\chi h}
    \defeq
    \pare{2\pi}^{\mathtt d/2}
    \mathcal F_{z\to\xi}
    \pare{\chi\pare z\mathfrak h}.
\]
Since $\chi\mathfrak h$ is a compactly supported distribution, convolution
with $\chi\mathfrak h$ is the Fourier multiplier
$\mathfrak m_{\chi h}\pare D$. Thus
\begin{equation}\label{eq:localized-homogeneous-model-factorization}
    t_h^{\loc}
    =
    \mathfrak m_{\chi h}\pare D .
\end{equation}
The cancellation condition in
\cref{eq:localized-homogeneous-model-cancellation} guarantees that, when
$\lambda=-\mathtt d$, the principal value $\pv h$ is a homogeneous tempered
distribution of degree $-\mathtt d$. When $\lambda=-\mathtt d+1$, the function $h$ is
locally integrable at the origin and defines a homogeneous tempered
distribution of degree $-\mathtt d+1$. Hence, in both cases, the Fourier transform
$\mathcal F\mathfrak h$ is smooth on $\bR^{\mathtt d}\setminus\set0$, homogeneous of
degree $-\pare{\lambda+\mathtt d}$, and therefore satisfies the symbol estimates
\[
    \av{\partial_\xi^\beta\dot{\mathfrak m}_h\pare\xi}
    \lesssim_\beta
    \av\xi^{-\pare{\lambda+\mathtt d}-\av\beta},
    \qquad
    \av\xi\ge\frac12.
\]
This proves
$\dot{\mathfrak m}_h\in\tilde\Gamma^{-\pare{\lambda+\mathtt d}}$.

Now, by construction \cref{eq:homogeneous-model-multiplier,eq:localized-homogeneous-model-factorization}, one has
\[t_h^{\loc}-\dot{\mathfrak m}_h
=\mathfrak m_{\chi h}-\dot{\mathfrak m}_h
=\pare{2\pi}^{\mathtt d/2}\mathcal F_{z\to\xi}\bra{\pare{\chi\pare z-1}h\pare z}.\]
Note that the factor $\chi\pare z-1$ vanishes near $z=0$.
Thus
\[
    g\pare z
    \defeq
    \pare{\chi\pare z-1}h\pare z
\]
is smooth on $\bR^{\mathtt d}$ and is equal to $-h\pare z$ for large $\av z$.
Fix a multi-index $\beta$ and $M\in\bN$. Since $g$ vanishes in a
neighborhood of the origin, there is no integrability issue there. Moreover,
for $\av z$ sufficiently large one has $g\pare z=-h\pare z$. Hence for any $L\in \bN$, by the
homogeneity of $h$,
\[
    \pare{-\Delta_z}^L
    \bra{\pare{-\ii z}^\beta g\pare z}
    =
    \mathcal{O}\pare{\av z^{\lambda+\av\beta-2L}},
    \qquad
    \av z\to+\infty.
\]
Therefore
\[
    \pare{-\Delta_z}^L
    \bra{\pare{-\ii z}^\beta g\pare z}
    \in
    L^1\pare{\bR^{\mathtt d};\bC}
\]
provided $2L>\mathtt d+\lambda+\av\beta$.
Since $\lambda\in\set{-\mathtt d,-\mathtt d+1}$, such an $L$ can be
chosen arbitrarily large.
Indeed, by homogeneity the last expression is
$O\pare{\av z^{\lambda+\av\beta-2L}}$ as $\av z\to+\infty$, and it is
smooth near the origin. Therefore, for $\av\xi\ge2$,
\begin{align}
    \av{
        \partial_\xi^\beta
        \mathcal F_{z\to\xi}g\pare\xi
    }
    &\le
    \av\xi^{-2L}
    \norm{
        \pare{-\Delta_z}^L
        \bra{\pare{-\ii z}^\beta g}
    }_{L^1}
    \lesssim_{M,\beta,h}
    \av\xi^{-M-\av\beta},
\end{align}
after choosing $L$ so large that $2L\ge M+\av\beta$. Hence,
$$t_h^{\loc}-\dot{\mathfrak m}_h\pare{D}\in\bigcap_{M\in\bN}\tilde\Gamma^{-M}.$$

Finally, if $h$ is real-valued, then the distributions
$\mathfrak h$, $\chi\mathfrak h$, and $g$ are real. Their Fourier
multipliers satisfy the usual symmetry
\[
    \overline{\mathfrak m\pare\xi}
    =
    \mathfrak m\pare{-\xi},
\]
and therefore the corresponding Fourier multiplier operators map real-valued
functions to real-valued functions.
\end{proof}

\begin{rem}[Composition with smoothing multipliers]
\label{rem:localized-homogeneous-model-coefficients}
Let $q\in\bN^\star$, $m\in\bR$, and let
$a\in\Gamma^m_{\geq q}$ be independent of $z$. If
$\mathfrak s\in\bigcap_{M\in\bN}\tilde\Gamma^{-M}\cap L^\infty \pare{\av{\xi}\leq 1 ; \mathbb{C}} $, then, for every
prescribed $\varrho\geq0$,
\begin{equation}\label{eq:localized-homogeneous-model-coefficient-smoothing}
    \OpBW{a\pare{U;x,\xi}}\circ \mathfrak s\pare D
    \in
    \cR^\varrho_{\geq q}.
\end{equation}
The same statement holds for a separately estimated linear parameter. If
\[
    a\bra\cdot
    \in\cL\pare{H^{s_0};\Gamma_{\geq q}^m},
\]
then
\[
    \omega
    \longmapsto
    \OpBW{a\bra\omega\pare{U;x,\xi}}\circ\mathfrak s\pare D
    \in
    \cL\pare{H^{s_0};\cR_{\geq q}^\varrho}.
\]
The same conclusion holds
componentwise for vector-valued smoothing multipliers. This is the standard Bony--Weyl composition estimate with one
factor of arbitrarily negative order. \fine
\end{rem}

\begin{definition}[Localized compensated homogeneous convolution models]
\label{def:localized-compensated-homogeneous-convolution-models}
Let
$h\in\cC^\infty\pare{\bR^{\mathtt d}\setminus\set0;\bC}$ be even and
homogeneous of degree $-\mathtt d-1$:
\begin{equation}\label{eq:degree-minus-three-even-kernel}
    h\pare{rz}=r^{-\mathtt d-1}h\pare z,
    \qquad
    h\pare{-z}=h\pare z,
    \qquad
    r>0.
\end{equation}
For a scalar function $f$, set
\begin{equation}\label{eq:localized-compensated-second-difference}
    \delta_z^{\pare2}f
    \defeq
    \delta_z f+\delta_{-z}f
    =
    2f-f\pare{\bullet-z}-f\pare{\bullet+z}.
\end{equation}
We define the localized compensated model by
\begin{equation}\label{eq:localized-compensated-homogeneous-operator}
    t_h^{\loc,\pare2}f
    \defeq
    \frac12
    \int_{\bR^{\mathtt d}}
    \chi\pare z h\pare z
    \delta_z^{\pare2}f
    \dd z .
\end{equation}
The corresponding global homogeneous multiplier is
\begin{equation}\label{eq:global-compensated-homogeneous-multiplier}
    \dot{\mathfrak m}_h^{(2)}\pare\xi
    \defeq
    \frac12
    \int_{\bR^{\mathtt d}}
    h\pare z
    \pare{
        2-e^{-\ii\xi\cdot z}-e^{\ii\xi\cdot z}
    }
    \dd z,
    \qquad
    \xi\in\bR^{\mathtt d}\setminus\set0.
\end{equation}
We set
\begin{equation}\label{eq:localized-compensated-order-zero-correction}
    c_h
    \defeq
    -
    \int_{\bR^{\mathtt d}}
    \pare{1-\chi\pare z}h\pare z
    \dd z
    \in\bC .
\end{equation}
\end{definition}

\begin{lemma}[Localized compensated homogeneous models are Fourier multipliers modulo smoothing]
\label{lem:localized-compensated-homogeneous-model}
Let $h$ satisfy
\cref{def:localized-compensated-homogeneous-convolution-models}. Then,
\[
    \dot{\mathfrak m}_h^{\pare{2}}\in\tilde\Gamma^1
\]
and
\begin{equation}\label{eq:localized-compensated-homogeneous-comparison}
    t_h^{\loc,\pare2}-
    \dot{\mathfrak m}_h^{\pare{2}}\pare D -c_h\cdot\in
    \bigcap_{M\in\bN}
    \tilde\Gamma^{-M} .
\end{equation}
If $h$ is real-valued, then all terms in
\cref{eq:localized-compensated-homogeneous-comparison} are real-to-real.
\end{lemma}

\begin{proof}
We first discuss the global multiplier. Since $2-e^{-\ii\xi\cdot z}-e^{\ii\xi\cdot z}
    =
    2\pare{1-\cos\pare{\xi\cdot z}}$,
the kernel multiplied by the second difference is
$O\pare{\av z^{1-\mathtt d}}$ near $z=0$, and the polar measure therefore
produces a bounded radial integrand. At infinity the polar radial weight is
$O\pare{r^{-2}}$. Hence the integral defining $\dot{\mathfrak m}_h^{\pare{2}}$ is absolutely
convergent. Writing $z=r\vartheta$, with $\vartheta\in\bS^{\mathtt d-1}$, gives
\begin{equation}\label{eq:compensated-homogeneous-multiplier-polar}
    \dot{\mathfrak m}_h^{\pare{2}}\pare\xi
    =
    \frac\pi2
    \int_{\bS^{\mathtt d-1}}
    h\pare\vartheta
    \av{\vartheta\cdot\xi}
    \dd\sigma\pare\vartheta .
\end{equation}
Here we use the third identity in \cref{eq:abstract-radial-distribution-identities}.
Formula \cref{eq:compensated-homogeneous-multiplier-polar} shows that
$\dot{\mathfrak m}_h^{\pare{2}}$ is homogeneous of degree $1$. If $\mathtt d=1$, its
smoothness away from the origin and the estimates below follow directly from
the formula on the two-point sphere $\bS^0$. If $\mathtt d\geq2$, fix a
coordinate patch on $\bS^{\mathtt d-1}$ and a smooth family of rotations
$Q_\nu\in SO\pare{\mathtt d}$ satisfying $Q_\nu e_{\mathtt d}=\nu$. After
the change of variables $\vartheta=Q_\nu\phi$, all angular derivatives fall
on the smooth amplitude $h\pare{Q_\nu\phi}$. A finite cover of the sphere
therefore gives
\[
    \av{\partial_\xi^\beta\dot{\mathfrak m}_h\pare\xi}
    \lesssim_\beta
    \av\xi^{1-\av\beta},
    \qquad
    \av\xi\geq\frac12,
\]
after the fixed high-frequency truncation. Thus,
$\dot{\mathfrak m}_h^{\pare{2}}\in\tilde\Gamma^1$.

Let $\mathfrak m_{\chi h}$ denote the multiplier of the localized model:
\begin{equation}\label{eq:localized-compensated-localized-multiplier}
    \mathfrak m_{\chi h}^{\pare{2}}\pare\xi
    \defeq
    \frac12
    \int_{\bR^{\mathtt d}}
    \chi\pare z h\pare z
    \pare{
        2-e^{-\ii\xi\cdot z}-e^{\ii\xi\cdot z}
    }
    \dd z .
\end{equation}
For Schwartz function $f$, the Fourier transform of
$\delta_z^{\pare2}f$ is the multiplication by
$2-e^{-\ii\xi\cdot z}-e^{\ii\xi\cdot z}$. Hence,
\begin{equation}\label{eq:localized-compensated-bare-factorization}
    t_h^{\loc,\pare2}
    =
    \mathfrak m_{\chi h}^{\pare{2}}\pare D .
\end{equation}
Set $g\pare z
    \defeq
    \pare{1-\chi\pare z}h\pare z $.
Then $g$ is smooth on $\bR^{\mathtt d}$ and is a classical symbol of degree
$-\mathtt d-1$ as
$\av z\to+\infty$. By subtracting the global multiplier from the localized
one, we obtain
\begin{equation}\label{eq:localized-compensated-tail-decomposition}
\begin{aligned}
    t_h^{\loc,\pare{2}}
    -
    \dot{\mathfrak m}_h^{\pare{2}}\pare\xi
    =\mathfrak m_{\chi h}^{\pare{2}}\pare\xi
    -
    \dot{\mathfrak m}_h^{\pare{2}}\pare\xi
    ={}&
    -
    \int_{\bR^{\mathtt d}}g\pare z\dd z+
    \frac12
    \int_{\bR^{\mathtt d}}
    g\pare z
    \pare{
        e^{-\ii\xi\cdot z}
        +
        e^{\ii\xi\cdot z}
    }
    \dd z
    \\
    ={}&
    c_h
    +
    \mathfrak s_h^{\mathrm{tail}}\pare\xi .
\end{aligned}
\end{equation}
We claim that $\mathfrak s_h^{\mathrm{tail}}$ is rapidly decaying at high
frequency. Let $\beta\in\bN_0^{\mathtt d}$ and $M\in\bN$. Choose $L$ large enough, such that
\[
    \pare{-\Delta_z}^L
    \bra{z^\beta g\pare z}
    \in
    L^1\pare{\bR^{\mathtt d};\bC}.
\]
This is possible because $g$ has degree $-\mathtt d-1$ at infinity and is smooth near
the origin. Integrating by parts gives, for $\av\xi\geq2$,
\[
    \av{
        \partial_\xi^\beta
        \int_{\bR^{\mathtt d}}
        g\pare z e^{\pm\ii\xi\cdot z}
        \dd z
    }
    \lesssim_{\beta,L}
    \av\xi^{-2L}
    \norm{
        \pare{-\Delta_z}^L
        \bra{z^\beta g}
    }_{L^1}.
\]
Taking $L$ with $2L\geq M+\av\beta$ yields
\begin{equation}\label{eq:localized-compensated-tail-smoothing-symbol}
    \av{
        \partial_\xi^\beta
        \mathfrak s_h^{\mathrm{tail}}\pare\xi
    }
    \lesssim_{\beta,M}
    \av\xi^{-M-\av\beta},
    \qquad
    \av\xi\geq2 .
\end{equation}
This proves that
\[
    \mathfrak s_h^{\mathrm{tail}}
    \in
    \bigcap_{M\in\bN}
    \tilde\Gamma^{-M} .
\]

If $h$ is real-valued, then $c_h$ is real and the multipliers satisfy
\[
    \overline{\mathfrak m\pare\xi}
    =
    \mathfrak m\pare{-\xi}.
\]
Thus the associated Fourier multiplier operators map real-valued functions to
real-valued functions.
\end{proof}

\paragraph{Local Taylor remainders for composed nonlinearities.}
\begin{lemma}[Local Taylor remainder]
\label{lem:kernel}
Let $r_{\mathtt g}>0$ and
$\mathsf g\in\cC^\omega\pare{B_\bR\pare{0,r_{\mathtt g}};\bR}$. Define
\begin{equation}\label{eq:kernel_Taylor_remainder}
\begin{aligned}
    k_{\bra2}\pare{\eta;x,z}
    \defeq{}&
    \mathsf g\pare{\Delta_z\eta\pare x}
    -
    \mathsf g\pare{\nabla_x\eta\pare x\cdot\hat z}
    +
    \frac{\av z}{2}
    \mathsf g'\pare{\nabla_x\eta\pare x\cdot\hat z}
    \pare{\hat z^\intercal\nabla_x^2\eta\pare x\hat z}.
\end{aligned}
\end{equation}
Then
\begin{equation}\label{eq:kernel_Taylor_remainder_classes}
    k_{\bra2}\in\KF^2_{\loc,\geq1}.
\end{equation}
If, in addition, $\mathsf g'\pare0=0$, then
\begin{equation}\label{eq:kernel_Taylor_remainder_quadratic}
    k_{\bra2}\in\KF^2_{\loc,\geq2},
    \qquad
    \chi\pare z k_{\bra2}\in\KF^2_{\geq2}.
\end{equation}

\end{lemma}

\begin{proof}
The statement is purely local in the convolution variable. Choose
$s_+>\mathtt d/2+2$ and
then $r_0>0$ so small that every $\eta\in\BallR{s_+}{r_0}$ satisfies
\begin{equation}\label{eq:smallness_condition}
    \norm{\nabla_x\eta}_{L^\infty}
    +
    \sup_{0<\av z\leq2}
    \norm{\Delta_z\eta}_{L^\infty}
    <
    \frac{r_{\mathtt g}}4\cdot
\end{equation}
Thus all arguments of $\mathsf g$ remain in a compact subset of its analytic
interval.

Set
\begin{equation}\label{eq:kernel-building-blocks}
    v\pare{x,z}
    \defeq
    \nabla_x\eta\pare x\cdot\hat z,
    \qquad
    h\pare{x,z}
    \defeq
    \hat z^\intercal\nabla_x^2\eta\pare x\hat z.
\end{equation}
The second order expansion of the difference quotient reads
\begin{equation}\label{eq:Delta_z_eta_second_order_remainder}
    \Delta_z\eta\pare x
    =
    v\pare{x,z}
    -
    \frac{\av z}{2}h\pare{x,z}
    +
    e\pare{x,z}.
\end{equation}
For every $s\geq s_+$ and every multi-index $\alpha$,
\begin{equation}\label{eq:kernel_building_blocks_tame}
\begin{aligned}
    \norm{\partial_z^\alpha v\pare{\cdot,z}}_{H^s}
    &\lesssim_{s,\alpha}
    \av z^{-\av\alpha}
    \norm{\eta}_{H^{s+1+\av\alpha}},
    \\
    \norm{\partial_z^\alpha h\pare{\cdot,z}}_{H^s}
    &\lesssim_{s,\alpha}
    \av z^{-\av\alpha}
    \norm{\eta}_{H^{s+2+\av\alpha}},
    \\
    \norm{\partial_z^\alpha e\pare{\cdot,z}}_{H^s}
    &\lesssim_{s,\alpha}
    \av z^{2-\av\alpha}
    \norm{\eta}_{H^{s+3+\av\alpha}},
\end{aligned}
\qquad
0<\av z\leq2.
\end{equation}
Indeed, derivatives of $\hat z$ have homogeneity $-\av\alpha$, while $e$ is
$\av z^2$ times an averaged third derivative of $\eta$ along the segment
$x-\theta z$; differentiating in $z$ either lowers the homogeneous degree or
differentiates the translated factor, producing the finite loss displayed
above.

Let
\begin{equation}\label{eq:theta-local-taylor-def}
    \vartheta\pare{x,z}
    \defeq
    \Delta_z\eta\pare x-v\pare{x,z}
    =
    -\frac{\av z}{2}h\pare{x,z}+e\pare{x,z}.
\end{equation}
Taylor's formula for $\mathsf g$ at $v$ gives
\begin{equation}\label{eq:g_power_series_kernel}
    \mathsf g\pare{v+\vartheta}
    =
    \mathsf g\pare v
    +
    \mathsf g'\pare v\vartheta
    +
    \sum_{\ell\geq2}
    \frac{\mathsf g^{\pare\ell}\pare v}{\ell!}
    \vartheta^\ell.
\end{equation}
Hence
\begin{equation}\label{eq:kernel_remainder_decomposition}
    k_{\bra2}
    =
    \mathsf g'\pare v e
    +
    \sum_{\ell\geq2}
    \frac{\mathsf g^{\pare\ell}\pare v}{\ell!}
    \vartheta^\ell.
\end{equation}
The analytic-composition estimates and the tame product estimates, combined
with \cref{eq:kernel_building_blocks_tame}, imply that for every
$\av\alpha\leq N$,
\begin{equation}\label{eq:kernel_KF2_local_bound}
    \norm{
        \partial_z^\alpha
        k_{\bra2}\pare{\eta;\cdot,z}
    }_{H^s}
    \lesssim_{s,N}
    \av z^{2-\av\alpha}
    \norm{\eta}_{H^{s+\mu_N}},
    \qquad
    0<\av z\leq2,
\end{equation}
for a finite loss $\mu_N$.

The vanishing order follows by separating
\begin{equation}\label{eq:gprime-vanishing-split}
    \mathsf g'\pare v
    =
    \mathsf g'\pare0
    +
    v\int_0^1\mathsf g''\pare{\theta v}\dd\theta.
\end{equation}
If $\mathsf g'\pare0\neq0$, the term $\mathsf g'\pare0e$ is linear in
$\eta$, and this gives the order $\geq1$. If $\mathsf g'\pare0=0$, the
first term in \cref{eq:kernel_remainder_decomposition} contains the
additional factor $v$, while the series starts quadratically in
$\vartheta$; therefore
\begin{equation}\label{eq:kernel_KF2_local_quadratic_bound}
    \norm{
        \partial_z^\alpha
        k_{\bra2}\pare{\eta;\cdot,z}
    }_{H^s}
    \lesssim_{s,N}
    \norm{\eta}_{H^{s_++\mu_N}}
    \norm{\eta}_{H^{s+\mu_N}}
    \av z^{2-\av\alpha}.
\end{equation}
This proves the local claims in
\cref{eq:kernel_Taylor_remainder_classes} and
\cref{eq:kernel_Taylor_remainder_quadratic}. Their rapidly decreasing
versions follow from \cref{rem:kernel-localization}.
\end{proof}

\paragraph{Finite-difference Sobolev estimates.}
\begin{lemma}
\label{lem:hsnorm}
For every $s\in\bR$ one has
\begin{equation}\label{eq:techlemma}
    \norm{\delta_z\eta}_{H^s}
    \lesssim
    \av z\norm{\eta}_{H^{s+1}},
    \qquad
    \norm{\Delta_z\eta-\nabla_x\eta\cdot\hat z}_{H^s}
    \lesssim
    \av z\norm{\eta}_{H^{s+2}}.
\end{equation}
\end{lemma}

\begin{proof}
The estimates are non-homogeneous. By Plancherel,
\begin{equation}\label{eq:techlemma-first-proof}
    \norm{\delta_z\eta}_{H^s}^2
    =
    \int_{\bR^{\mathtt d}}
    \angles\xi^{2s}
    \av{\widehat\eta\pare\xi}^2
    \av{1-e^{-\ii\xi\cdot z}}^2
    \dd\xi.
\end{equation}
Since $\av{1-e^{-\ii X}}\lesssim\min\set{1,\av X}\lesssim\av X$, we get
\begin{equation}\label{eq:techlemma-first-bound}
    \norm{\delta_z\eta}_{H^s}^2
    \lesssim
    \av z^2
    \int_{\bR^{\mathtt d}}
    \angles\xi^{2s}\av\xi^2
    \av{\widehat\eta\pare\xi}^2
    \dd\xi
    \lesssim
    \av z^2
    \norm{\eta}_{H^{s+1}}^2.
\end{equation}
Similarly,
\begin{equation}\label{eq:techlemma-second-proof}
\begin{aligned}
    &\norm{\Delta_z\eta-\nabla_x\eta\cdot\hat z}_{H^s}^2
   =
    \int_{\bR^{\mathtt d}}
    \angles\xi^{2s}
    \av{\widehat\eta\pare\xi}^2
    \av{
        \frac{1-e^{-\ii\xi\cdot z}-\ii\xi\cdot z}{\av z}
    }^2
    \dd\xi.
\end{aligned}
\end{equation}
Using
$\av{1-e^{-\ii X}-\ii X}\lesssim\min\set{\av X,\av X^2}\lesssim\av X^2$
we infer
\begin{equation}\label{eq:techlemma-second-bound}
    \norm{\Delta_z\eta-\nabla_x\eta\cdot\hat z}_{H^s}^2
    \lesssim
    \av z^2
    \int_{\bR^{\mathtt d}}
    \angles\xi^{2s}\av\xi^4
    \av{\widehat\eta\pare\xi}^2
    \dd\xi
    \lesssim
    \av z^2
    \norm{\eta}_{H^{s+2}}^2.
\end{equation}
Taking square roots proves \cref{eq:techlemma}.
\end{proof}

\begin{rem}[Dimension convention]
\label{rem:paradifferential-toolkit-dimension-convention}
All the results of this section are stated on $\bR^{\mathtt d}$ for an
arbitrary $\mathtt d\geq1$. Unless explicitly stated otherwise, in the
remainder of the manuscript we set $\mathtt d=2$; the only exceptions are
\cref{sec:abstract-structure-preserving-paralinearization-Rd} and \cref{app:abstract-locally-integrable-paralinearization-Rd}, where the
dimension remains arbitrary. \fine
\end{rem}

\section{Paralinearization}\label{sec:paralinearization}
\newcommand{\tail}{\textnormal{tail}}
This section is devoted to the paralinearization of the system of interest in this study. It represents the core of the manuscript.

\begin{notation}\label{notation:paralinearization}
\begin{enumerate}

\item In the present section $1 \ll \varrho \ll  s_0 $ and the explicit value of $s_0$ may vary from line to line.

\item Let $m,n\in\bR$, $p\in\bN^\star$, and
$\vec{U}_p\in\pare{\BallR{s_0}{r_0}}^p$.
\begin{enumerate}
  \item For \(\ell\in\bN\), we denote with
  \[a_{\bra{m,\geq\ell}}\pare{\vec{U}_p;x,\xi},\quad
  A_{\bra{m,\geq\ell}}\pare{\vec{U}_p;x,\xi},\quad
  {\bm A}_{\bra{m,\geq\ell}}\pare{\vec{U}_p;x,\xi}, \quad
\mathfrak{A}_{\bra{m,\geq\ell}}\pare{\vec{U}_p;x,\xi}\]
  a generic element in \[\Gamma_{\geq\ell}^m,\quad
  \pare{\Gamma_{\geq\ell}^m}^2,\quad
  \pare{\Gamma_{\geq\ell}^m}^{2\times2},\quad
  \pare{\Gamma_{\geq\ell}^m}^{2\times3},\] respectively, whose explicit
  expression may vary from line to line.  When \(\ell=1\) we omit the homogeneity subscript. {When the corresponding symbols are pointwise
  real componentwise, we add the superscript \(\bR\), writing
  \[a^\bR_{\bra{m,\geq\ell}},\quad
  A^\bR_{\bra{m,\geq\ell}},\quad
  {\bm A}^\bR_{\bra{m,\geq\ell}}\quad, \mathfrak{A}^\bR_{\bra{m,\geq\ell}},\] respectively. This superscript
  records pointwise reality and not the real-to-real symmetry of a symbol.}
   Of these we define the subclasses
\begin{enumerate}
    \item \emph{Linear symbols} defined as
    \begin{align}
        A^{\rm lin}_{\bra{m,\geq\ell}}\pare{\eta, \omega} \defeq & \ \vect{a_{\bra{m,\geq\ell}}\bra{\omega}\pare{\eta}}{a_{\bra{m,\geq\ell}}\pare{\eta}},
        &
                {\bm A}^{\rm lin}_{\bra{m,\geq\ell}}\pare{\eta, \omega} \defeq & \  \vectT{A_{\bra{m,\geq\ell}}\bra{\omega}\pare{\eta}}{A_{\bra{m,\geq\ell}}\pare{\eta}},
                \\
                 {\bm a}^{\rm lin}_{\bra{m,\geq\ell}}\pare{\eta, \Omega} \defeq & \ \vect{a_{\bra{m,\geq\ell}}\bra{\Omega}\pare{\eta}}{A_{\bra{m,\geq\ell}}\pare{\eta}},
        &
                \mathfrak{ A}^{\rm lin}_{\bra{m,\geq\ell}}\pare{\eta, \Omega} \defeq & \  \vectT{A_{\bra{m,\geq\ell}}\bra{\Omega}\pare{\eta}}{{\bm A}_{\bra{m,\geq\ell}}\pare{\eta}}.
    \end{align}

    \item \emph{Bilinear symbols} defined as
    \begin{align}
        {\bm a}^{\rm bil}_{\bra{m,\geq\ell}}\pare{\eta, \omega_1, \omega_2} \defeq & \
            \vectD
                {a_{\bra{m,\geq\ell}}\bra{\omega_1,\omega_2}\pare{\eta}}
                {a_{\bra{m,\geq\ell}}\bra{\omega_2}\pare{\eta}}
                {a_{\bra{m,\geq\ell}}\bra{\omega_1}\pare{\eta}},
        \\
        \mathfrak{A}^{\rm bil}_{\bra{m,\geq\ell}}\pare{\eta, \omega_1, \omega_2} \defeq & \
            \vectDT
                {A_{\bra{m,\geq\ell}}\bra{\omega_1,\omega_2}\pare{\eta}}
                {A_{\bra{m,\geq\ell}}\bra{\omega_2}\pare{\eta}}
                {A_{\bra{m,\geq\ell}}\bra{\omega_1}\pare{\eta}}.
    \end{align}
    
\end{enumerate}
The above definitions are such that the application $$\omega \mapsto \OpBW{A^{\lin}_{\bra{m,\geq\ell}}\pare{\eta, \omega}}\cdot\vect{\eta}{\omega} \in \cL\pare{H^s ; H^{s-m} }$$ for all $s\geq s_0$, the same holds for all the other terms with due modifications. 

\item For \(\ell\in\bN^\star\), we denote with
\[k_{\bra{n,\geq\ell}}\pare{\vec U_p;x,z},\quad
K_{\bra{n,\geq\ell}}\pare{\vec U_p;x,z},\quad
{\bm K}_{\bra{n,\geq\ell}}\pare{\vec U_p;x,z},\quad
\mathfrak K_{\bra{n,\geq\ell}}\pare{\vec U_p;x,z}\]
a generic element of \[\KF^n_{\geq\ell},\quad
\pare{\KF^n_{\geq\ell}}^2,\quad
\pare{\KF^n_{\geq\ell}}^{2\times2},\quad
\pare{\KF^n_{\geq\ell}}^{2\times3},\] respectively, whose explicit
expression may vary from line to line. When \(\ell=1\) we omit the homogeneity subscript.
The forms
\[k^{\loc}_{\bra{n,\geq\ell}},\quad
K^{\loc}_{\bra{n,\geq\ell}},\quad
{\bm K}^{\loc}_{\bra{n,\geq\ell}}, \quad
\mathfrak K^{\loc}_{\bra{n,\geq\ell}}\]
denote generic elements of \(\KF^n_{\loc,\geq\ell}\) and its corresponding
vector and matrix classes. 
For \(n<0\) with \(\ell-n>\mathtt d\), we write
\[k^{\tail}_{\bra{n,\geq\ell}},\quad
K^{\tail}_{\bra{n,\geq\ell}},\quad
{\bm K}^{\tail}_{\bra{n,\geq\ell}},\quad
\mathfrak K^{\tail}_{\bra{n,\geq\ell}}\]
for generic elements of the existing class
\(\KF_{\textnormal{tail},\geq\ell}\) and its corresponding vector and
matrix classes, whose profiles in \cref{def:kernel-global-global} satisfy
\cref{eq:tail-kernel-profile-decay} with \(m=-n\), componentwise.
For tail kernels, \(n\) records the order of the profile \(b\), not the
decay of the full kernel.

When these kernels are real-valued componentwise for real backgrounds
and real separately estimated arguments, we add the superscript
\(\bR\), writing \[k^\bR_{\bra{n,\geq\ell}},\quad
K^\bR_{\bra{n,\geq\ell}},\quad
{\bm K}^\bR_{\bra{n,\geq\ell}},\quad
\mathfrak K^\bR_{\bra{n,\geq\ell}},\] respectively.
The corresponding local forms are
\[k^{\loc,\bR}_{\bra{n,\geq\ell}},\quad
K^{\loc,\bR}_{\bra{n,\geq\ell}},\quad
{\bm K}^{\loc,\bR}_{\bra{n,\geq\ell}}, \quad
\mathfrak K^{\loc,\bR}_{\bra{n,\geq\ell}},\]
and the tail forms are
\[k^{\tail,\bR}_{\bra{n,\geq\ell}},\quad
K^{\tail,\bR}_{\bra{n,\geq\ell}},\quad
{\bm K}^{\tail,\bR}_{\bra{n,\geq\ell}},\quad
\mathfrak K^{\tail,\bR}_{\bra{n,\geq\ell}}.\]

\item For \(\ell\in\bN^\star\), we denote with
\[r_{\bra{\varrho,\geq\ell}}\pare{\vec U_p},\quad
R_{\bra{\varrho,\geq\ell}}\pare{\vec U_p},\quad
{\bm R}_{\bra{\varrho,\geq\ell}}\pare{\vec U_p},\quad
\mathfrak R_{\bra{\varrho,\geq\ell}}\pare{\vec U_p}\]
a generic element of \[\cR^\varrho_{\geq\ell},\quad
\pare{\cR^\varrho_{\geq\ell}}^2,\quad
\pare{\cR^\varrho_{\geq\ell}}^{2\times2},\quad
\pare{\cR^\varrho_{\geq\ell}}^{2\times3},\] respectively, whose explicit
expression may vary from line to line. When the corresponding smoothing
operators are real-to-real componentwise, we add the superscript \(\bR\), writing
\[r^\bR_{\bra{\varrho,\geq\ell}},\quad
R^\bR_{\bra{\varrho,\geq\ell}},\quad
{\bm R}^\bR_{\bra{\varrho,\geq\ell}},\quad
\mathfrak R^\bR_{\bra{\varrho,\geq\ell}},\] respectively.

We denote with
\[r_{\bra{\varrho,n,\geq\ell}}\pare{\vec U_p;z},\quad
R_{\bra{\varrho,n,\geq\ell}}\pare{\vec U_p;z},\quad
{\bm R}_{\bra{\varrho,n,\geq\ell}}\pare{\vec U_p;z},\quad
\mathfrak R_{\bra{\varrho,n,\geq\ell}}\pare{\vec U_p;z}\]
a generic element of \[\KR^{\varrho,n}_{\geq\ell},\quad
\pare{\KR^{\varrho,n}_{\geq\ell}}^2,\quad
\pare{\KR^{\varrho,n}_{\geq\ell}}^{2\times2},\quad
\pare{\KR^{\varrho,n}_{\geq\ell}}^{2\times3},\] respectively, whose
explicit expression may vary from line to line. Aggregate notation with one
threshold is used only when every component has the displayed class
membership. Of these we define the subclasses
\begin{enumerate}
    \item \emph{Linear remainders} defined as
    \begin{align}
        R^{\rm lin}_{\bra{\varrho,\geq\ell}}\pare{\eta, \omega} \defeq & \ \vect{r_{\bra{\varrho,\geq\ell}}\bra{\omega}\pare{\eta}}{r_{\bra{\varrho,\geq\ell}}\pare{\eta}},
        &
                {\bm R}^{\rm lin}_{\bra{\varrho,\geq\ell}}\pare{\eta, \omega} \defeq & \  \vectT{R_{\bra{\varrho,\geq\ell}}\bra{\omega}\pare{\eta}}{R_{\bra{\varrho,\geq\ell}}\pare{\eta}},
                \\
                 {\bm r}^{\rm lin}_{\bra{\varrho,\geq\ell}}\pare{\eta, \Omega} \defeq & \ \vect{r_{\bra{\varrho,\geq\ell}}\bra{\Omega}\pare{\eta}}{R_{\bra{\varrho,\geq\ell}}\pare{\eta}},
        &
                \mathfrak{ R}^{\rm lin}_{\bra{\varrho,\geq\ell}}\pare{\eta, \Omega} \defeq & \  \vectT{R_{\bra{\varrho,\geq\ell}}\bra{\Omega}\pare{\eta}}{{\bm R}_{\bra{\varrho,\geq\ell}}\pare{\eta}}.
    \end{align}

    \item \emph{Bilinear remainders} defined as
    \begin{align}
        {\bm r}^{\rm bil}_{\bra{\varrho,\geq\ell}}\pare{\eta, \omega_1, \omega_2} \defeq & \
            \vectD
                {r_{\bra{\varrho,\geq\ell}}\bra{\omega_1,\omega_2}\pare{\eta}}
                {r_{\bra{\varrho,\geq\ell}}\bra{\omega_2}\pare{\eta}}
                {r_{\bra{\varrho,\geq\ell}}\bra{\omega_1}\pare{\eta}},
        \\
        \mathfrak{R}^{\rm bil}_{\bra{\varrho,\geq\ell}}\pare{\eta, \omega_1, \omega_2} \defeq & \
            \vectDT
                {R_{\bra{\varrho,\geq\ell}}\bra{\omega_1,\omega_2}\pare{\eta}}
                {R_{\bra{\varrho,\geq\ell}}\bra{\omega_2}\pare{\eta}}
                {R_{\bra{\varrho,\geq\ell}}\bra{\omega_1}\pare{\eta}}.
    \end{align}
    
\end{enumerate}
The real-to-real versions
\[R^{\lin,\bR}_{\bra{\varrho,\geq\ell}},\quad
{\bm R}^{\lin,\bR}_{\bra{\varrho,\geq\ell}},\quad
{\bm r}^{\lin,\bR}_{\bra{\varrho,\geq\ell}},\quad
\mathfrak R^{\lin,\bR}_{\bra{\varrho,\geq\ell}},\quad
{\bm r}^{\bil,\bR}_{\bra{\varrho,\geq\ell}},\quad
\mathfrak R^{\bil,\bR}_{\bra{\varrho,\geq\ell}}\]
are defined by the same formulas with each constituent
\[r_{\bra{\varrho,\geq\ell}},\quad
R_{\bra{\varrho,\geq\ell}}, \quad
{\bm R}_{\bra{\varrho,\geq\ell}}\]
replaced by its \(\bR\)-superscripted version. For constituents with
bracketed arguments, the real-to-real property holds for every real
choice of those arguments.

The above definitions are such that the application $$\omega \mapsto R^{\lin}_{\bra{\varrho,\geq\ell}}\pare{\eta, \omega}\cdot\vect{\eta}{\omega} \in \cL\pare{H^s ; H^{s+\varrho} }$$ for all $s\geq s_0$, the same holds for all the other terms with due modifications.

\item For \(\ell\in\bN^\star\), we denote with $a_{\bra{m,\geq\ell}}\pare{\Atw\middle|\vec{U}_p;x,\xi}, \ r_{\bra{\varrho,\geq\ell}}\pare{\Atw\middle|\vec{U}_p}$ and so on families of, respectively symbols and smoothing operators, that are real analytic w.r.t. the Atwood number $\Atw\in \bra{-2,2}$. An analogous notation is used for vector- and matrix-valued symbols and smoothing operators.

\item Given $A\in \pare{\Gamma_{\geq1}^m}^2$ and a vector-valued function $U$, we use the compact notation
\begin{equation}
    \OpBW{A}\cdot\bra{U}\defeq \OpBW{A}\cdot U \defeq\sum_{\mathsf{j}=1}^2 \OpBW{a_\mathsf{j}}\bra{u_\mathsf{j}}.
\end{equation}
    \end{enumerate}
\end{enumerate}
\end{notation}
We shall consider $\Omega\approx \nabla^{\perp}\psi.$ Therefore, the quasilinear terms are of order $0$ in $\Omega$ and $1$ in $\eta.$

The main goal of this section is to prove the following theorem.
\begin{theorem}[Paralinearization of the canonical CDE]\label{thm:paraCDEetapsi}
Let $\sigma\geq0$, $\rho^\pm\geq0$ with $\aangle{\rho^\pm}>0$,
$\mathfrak g\in\bR$, and
$\background\in\bR^2$. There exists $\varrho_0>0$ such that, for any
$\varrho\geq\varrho_0$, there exist $s_0>0$ and $r_0>0$ such that, if
\[
    \eta,\tilde{\psi}\in \BallR{s_0}{r_0}
\]
solve the canonical CDE \cref{eq:CDE}, then
$\pare{\eta,\tilde{\psi}}$ solves the paradifferential system
\begin{equation}\label{eq:paralinearazationetapsi}
\begin{aligned}
\vect{\eta_t}{\tilde{\psi}_t}=&
\OpBW{
    {\bm Q}
    \pare{\Atw,\background\middle|\eta,\tilde{\psi};x,\xi}
    +
    {\bm B}
    \pare{\Atw,\background\middle|\eta,\tilde{\psi};x,\xi}
    +
    \ii
    \pare{
        \mathpzc V_{\textnormal{can}}\pare{\Atw,\background\middle|\eta,\tilde{\psi};x}
        \cdot\xi
    }
    \Id_{\bC^2}
}
\vect{\eta}{\tilde{\psi}}
\\
+&
\OpBW{ \vectT{ A_{\bra{0,\geq1}}^\eta
    \pare{\Atw,\background\middle|\eta,\tilde{\psi};x,\xi}}{ A_{\bra{-1,\geq1}}^\psi\pare{\Atw,\background\middle|\eta,\tilde{\psi};x,\xi}}
    }
   \vect{\eta}{\tilde{\psi}}
+
{\bm R}_{\bra{\varrho,\geq1}}^{{\lin},\bR}
\pare{\Atw,\background\middle|\eta,\tilde{\psi}}
\vect{\eta}{\tilde{\psi}} .
\end{aligned}
\end{equation}
When the lower-order term in the second line is viewed as a two-by-two block,
it has the order structure
\begin{equation}\label{eq:lower-order-block-structure-main}
    \begin{bmatrix}
        \Gamma_{\geq1}^0 & \Gamma_{\geq1}^{-1}\vspace{0.1cm}\\
        \Gamma_{\geq1}^0 & \Gamma_{\geq1}^{-1}
    \end{bmatrix}.
\end{equation}
Here the objects appearing in \cref{eq:paralinearazationetapsi} are defined
as follows.
\begin{itemize}
    \item $\Atw$ is the Atwood number, defined in \cref{eq:Atwood}.

    \item ${\bm Q}\in\pare{\Gamma_{\geq0}^2}^{2\times 2}$ is the quasilinear symbol
    explicitly defined as
    \begin{multline}\label{eq:mathpzcQ}
        {\bm Q}
        \pare{\rho^\pm,\Atw,\background\middle|\eta,\tilde{\psi};x,\xi}
        \\
        \defeq
        \begin{bmatrix}
            0
            &
            \lambda_1\pare{\eta;x,\xi} + \Atw \lambda_0\pare{\eta;x,\xi}
            \\
            -\mathfrak g\Atw
            +
            a^{\pi_0,\eta}_2
            \pare{\rho^\pm\middle|\eta;x,\xi}
            +
            \lambda_1\pare{\eta;x,\xi}
            q^{\tilde{\pi}}_0
            \pare{\Atw,\background\middle|\eta,\tilde{\psi};x,\xi}
            &
            0
        \end{bmatrix}.
    \end{multline}
    More explicitly:
    \begin{itemize}
        \item $\lambda_1\in\Gamma_{\geq0}^1$ is real-valued and is explicitly defined as
        \begin{equation}\label{eq:lambda1}
            \lambda_1\pare{\eta;x,\xi}
            \defeq
            \frac{\av{\xi}}{2}
            \sqrt{
                1+
                \av{
                    \nabla^\perp\eta\cdot
                    \frac{\xi}{\av{\xi}}
                }^2
            } .
        \end{equation}

        \item $\lambda_0 \in \Gamma_{\geq1}^0$ is real valued and  defined as
                 \begin{equation}\label{eq:lambda0}
                      \lambda_0\pare{\eta;x,\xi}\defeq \pare{\frac{\av{\xi}}{4\sqrt{2} \ \lambda_1\pare{\eta;x,\xi}}}^2 \ \pare{\hat{\xi}^\perp}^\intercal {\bm H}_x \eta \hat{\xi}^\perp, 
                 \end{equation}
                 where $\lambda_1$ is defined in \cref{eq:lambda1};

        \item $a^{\pi_0,\eta}_2\in\Gamma_{\geq0}^2$ is real-valued and is explicitly
        defined as
        \begin{equation}\label{eq:Api_star_2}
            a^{\pi_0,\eta}_2
            \pare{\rho^\pm\middle|\eta;x,\xi}
            \defeq
            -
            \frac{\sigma}{\aangle{\rho^\pm}}
            \frac{
                \av{\xi}^2+\av{\nabla^\perp\eta\cdot\xi}^2
            }{
                \pare{1+\av{\nabla\eta}^2}^{\frac32}
            }{=-\frac{4\sigma}{\aangle{\rho^\pm}}
            \frac{
                \lambda_1^2\pare{\eta;x,\xi}
            }{
                \pare{1+\av{\nabla\eta\pare{x}}^2}^{\frac32}
            }}\cdot
        \end{equation}

        \item $q^{\tilde{\pi}}_0\in\Gamma_{\geq0}^0$ is real-valued, even in $\xi$ and is explicitly
        defined as
        \begin{equation}\label{eq:qpitilde_main}
            q^{\tilde{\pi}}_0
            \pare{\Atw,\background\middle|\eta,\tilde{\psi};x,\xi}
            \defeq
            \pare{4-\Atw^2}
            \left(
                A_{0}^{\mathpzc{B}_0,\omega;\perp}\pare{\eta;x,\xi}
                \cdot
                \Omega
                \pare{\background\middle|\eta,\tilde{\psi};x}
            \right)^2 ,
        \end{equation}
        where $A_{0}^{\mathpzc{B}_0,\omega}$ is the vector symbol defined as 
        \begin{equation}\label{eq:H0_definition}
        A_{0}^{\mathpzc{B}_0,\omega}\pare{\eta;x,\xi}
        \defeq
        \nabla_\xi\phi^{\mathpzc B_0}_1\pare{\eta;x,\xi}
        ,\qquad\phi^{\mathpzc B_0}_1\pare{\eta;x,\xi}\defeq
        -
        \frac{\lambda_1\pare{\eta;x,\xi}}
        {1+\av{\nabla\eta\pare x}^2};
    \end{equation}
    \end{itemize}

    \item ${\bm B}\in\pare{\Gamma_{\geq1}^1}^{2\times 2}$ has the Alinhac block form
    \begin{equation}\label{eq:mathpzcB}
        {\bm B}
        \pare{\Atw,\background\middle|\eta,\tilde{\psi};x,\xi}
        \defeq
        -
        \begin{bmatrix}
            -b_{\textnormal{can}}
            \pare{\Atw,\background\middle|\eta,\tilde{\psi};x}
            &
            0
            \\
            \pare{
                b_{\textnormal{can}}
                \pare{\Atw,\background\middle|\eta,\tilde{\psi};x}
            }^2
            &
            b_{\textnormal{can}}
            \pare{\Atw,\background\middle|\eta,\tilde{\psi};x}
        \end{bmatrix}
        \lambda_1\pare{\eta;x,\xi}.
    \end{equation}
    Here
    \begin{equation}\label{eq:bpitilde_main}
        b_{\textnormal{can}}
        \pare{\Atw,\background\middle|\eta,\tilde{\psi};x}
        \defeq
        \Atw\,
         v_3^{\star}
        \bra{
            \Omega
            \pare{\background\middle|\eta,\tilde{\psi}}
        }
        \pare{\eta;x}
       -
        \frac{
            \nabla^\perp\eta\cdot
            \Omega
            \pare{\background\middle|\eta,\tilde{\psi};x}
        }{
            1+\av{\nabla\eta}^2
        } ,
    \end{equation}
  and the canonical vorticity variable is
    \begin{equation}\label{eq:Omegacan_paraCDE}
        \Omega
        \pare{\background\middle|\eta,\tilde{\psi};x}
        \defeq
        \tilde{\Omega}\bra{\background}\pare{\eta;x}
        +
        S\pare{\Atw\middle|\eta}\tilde{\psi}\pare{x},
    \end{equation}
    where $S\pare{\Atw\middle|\eta}$ is defined in \cref{eq:S} and $\tilde{\Omega}\bra{\background}\pare{\eta}$ is defined in
    \cref{eq:Omega_bg_can_def};

    \item The transport velocity is
    \begin{equation}\label{eq:Vpitilde_main}
\begin{aligned}
        \mathpzc V_{\textnormal{can}}
        \pare{\Atw,\background\middle|\eta,\tilde{\psi};x}
        &\defeq
         V^{\star}
        \bra{
            \Omega
            \pare{\background\middle|\eta,\tilde{\psi}}
        }
        \pare{\eta;x}
        \\
        &\quad
        +
        \frac{\Atw}{4\pare{1+\av{\nabla\eta}^2}}
        \bm J
        \pare{
            \Omega
            \pare{\background\middle|\eta,\tilde{\psi};x}
            +
            \pare{
                \nabla\eta\cdot
                \Omega
                \pare{\background\middle|\eta,\tilde{\psi};x}
            }
            \nabla\eta
        } .
\end{aligned}
    \end{equation}
    where
     ${\bm v}^\star=\vect{ V^{\star}}{ v_3^{\star}}$ is defined in
    \cref{eq:bmvtilde}.
    In particular, the transport coefficient belongs to
    $\pare{\cF^\bR_{\geq0}}^2$ and splits as
    \begin{equation}\label{eq:Vpitilde_main_split}
        \mathpzc V_{\textnormal{can}}
        \pare{\Atw,\background\middle|\eta,\tilde{\psi}}
        =
        \frac{\Atw}{4}\background^\perp
        +
        \mathpzc V_{\textnormal{can}}^{\mathrm{nl}}
        \pare{\Atw,\background\middle|\eta,\tilde{\psi}},
        \qquad
        \mathpzc V_{\textnormal{can}}^{\mathrm{nl}}
        \in\pare{\cF_{\geq1}^\bR}^2.
    \end{equation}
    Moreover,
    $\mathpzc V_{\textnormal{can}}
    \pare{\Atw,\background\middle|\eta,\tilde{\psi};x}$
    is real-valued and independent of $\xi$.

    \item The lower-order symbols satisfy
    \begin{equation}\label{eq:lower-order-symbol-classes-main}
        A_{\bra{0,\geq1}}^\eta
        \in
        \pare{\Gamma_{\geq1}^0}^2,
        \qquad
        A_{\bra{-1,\geq1}}^\psi
        \in
        \pare{\Gamma_{\geq1}^{-1}}^2,
    \end{equation}
    or equivalently the associated two-by-two block has the order structure
    stated in \cref{eq:lower-order-block-structure-main}. They denote generic
    symbols as per \cref{notation:paralinearization}.

    \item ${\bm R}_{\bra{\varrho,\geq1}}^{{\lin},\bR}\in
    \pare{\cR_{\geq1}^\varrho}^{2\times 2}$ denotes a generic smoothing
    remainder as per \cref{notation:paralinearization}.
\end{itemize}
\end{theorem}

\begin{remark}[On the equivalent formulas for $\lambda_1$]\label{rem:lambdas}
Let 
\begin{equation}
\tilde{\lambda}_1(x,\xi) \defeq \sqrt{\pare{1+\av{\nabla\eta(x)}^2}\,\av{\xi}^2-\pare{\nabla\eta(x)\cdot\xi}^2},
\end{equation}
which is the principal symbol in the expansion of the DN operator in \cite{ABZ2014}. Factoring $\av{\xi}$ gives $$\tilde{\lambda}_1(x,\xi)=\av{\xi}\sqrt{\,1+\av{\nabla \eta}^2-\pare{\nabla \eta \cdot \hat{\xi}}^2\,}.$$
Notice that $\pare{\hat{\xi},\hat{\xi}^{\perp}}$ is an orthonormal basis. Hence,
\begin{equation}
\av{\nabla \eta}^2=\pare{\nabla \eta\cdot\hat{\xi}}^2+\pare{\nabla \eta\cdot\hat{\xi}^{\perp}}^2
      =\pare{\nabla \eta\cdot\hat{\xi}}^2+\pare{\nabla^\perp\eta\cdot\hat{\xi}}^2,
\end{equation}
so that $$\tilde{\lambda}_1(x,\xi)=\av{\xi}\sqrt{\,1+\av{\nabla^\perp\eta(x)\cdot\hat{\xi}}^2\,}.$$
Therefore, $\lambda_1 =\xfrac{\tilde{\lambda}_1}{2}$, and the factor $1/2$ depends on the fact that here $\psi$ is $2\psi$ of \cite{ABZ2014} (c.f. \cref{rem:main_result_one_phase_water_wave_structure}).
\fine
\end{remark}

The proof of \cref{thm:paraCDEetapsi} is organized into four logical steps:
\begin{enumerate}
    \item \textbf{Abstract structure-preserving paralinearization.}
    In \cref{sec:abstract-structure-preserving-paralinearization-Rd}, we establish the abstract paralinearization theorem for critical homogeneous graph-kernel operators on $\bR^{\mathtt d}$. The proof separates the local--tail and Bony reductions from the quasilinear and transported contributions, and identifies the corresponding homogeneous Fourier multipliers and lower-order symbols.

    \item \textbf{Paralinearization of the singular integral operators.}
    In \cref{sec:paralinearization_SIO}, we apply \cref{thm:abstract-structure-preserving-paralinearization-Rd} to the primitive singular integral operators $\mathpzc{B}_0\pare{\eta}$ and $\mathpzc{q}\pare{\eta}$. We then identify their principal symbols and combine the resulting paralinearizations in the form needed for the subsequent analysis.

    \item \textbf{Paralinearization of the evolution operators.}
    In \cref{sec:paralinearization_EO}, we use the preceding singular-integral expansions to paralinearize the operators entering the canonical evolution equations. Starting from ${\bm v}^\star$, we treat successively $T\pare{\eta}$ and $\tau_0\pare{\eta}$, the change of variable $S\pare{\Atw\middle|\eta}$, the Dirichlet--Neumann-type operators $\mathpzc{g}^\star\pare{\eta}$ and $\mathpzc{g}\pare{\Atw\middle|\eta}$, and the Bernoulli contributions $\pi^\star$, $\pi$, and $\tilde{\pi}$.

    \item \textbf{Synthesis.}
    Finally, we assemble these expansions in the canonical variables $\pare{\eta,\tilde{\psi}}$ and obtain the complete paralinear system \cref{eq:paralinearazationetapsi}, thereby concluding the proof of \cref{thm:paraCDEetapsi}.
\end{enumerate}

\subsection{Abstract structure-preserving paralinearization on \texorpdfstring{$\bR^{\mathtt d}$}{R-d}}
\label{sec:abstract-structure-preserving-paralinearization-Rd}

\paragraph{Dimension and scalar convention.}
Throughout this subsection, $\mathtt d\geq1$.  In accordance with
\cref{notation:vector_convention}, scalar functions, operators, and symbols
are denoted by lower-case letters; only $x,\xi\in\bR^{\mathtt d}$ retain the
lower-case notation reserved there for spatial and frequency variables.  We
state the result in the scalar case; vector- and matrix-valued versions follow
componentwise.  The Fourier normalization and the measure on
$\bS^{\mathtt d-1}$ are fixed in
\cref{notation:abstract-Fourier-convention-Rd} above.

\begin{definition}[Admissible homogeneous pair]
\label{def:abstract-admissible-angular-pair-Rd}
Let
\(
    0<n\leq\mathtt d
\)
and
\(
    \ell\in\bN^\star
\).
Thus \(-n\) is the degree of homogeneity of the kernel
\(\av z^{-n}\), and the critical endpoint corresponds to
\(n=\mathtt d\).
An \emph{admissible pair of order \(-n\) and vanishing order
\(\ell\)} is a pair
\begin{equation}\label{eq:abstract-admissible-angular-pair-Rd}
    \pare{h,k}
    \in
    \cC^\infty\pare{\bS^{\mathtt d-1};\bR}
    \times
    \cC^\omega\pare{\pare{-r_k,r_k};\bR}
\end{equation}
with the following properties:
\begin{itemize}

    \item the function
    \begin{equation}\label{eq:abstract-admissible-K-vanishing}
        \widetilde k\pare X
        \defeq
        k\pare X-k\pare0
    \end{equation}
    vanishes to order at least \(\ell\) at the origin, namely
    \begin{equation}\label{eq:abstract-admissible-K-vanishing-order}
        \widetilde k^{\pare j}\pare0
        =
        0,
        \qquad
        0\leq j\leq\ell-1.
    \end{equation}

    \item In the critical case \(n=\mathtt d\), there exists
    \(\varepsilon\in\set{-1,1}\) such that
    \begin{subequations}
    \begin{align}
        h\pare{-\vartheta}
        &=
        \varepsilon h\pare\vartheta,
        \label{eq:abstract-admissible-H-parity}
        \\
        k\pare{-X}
        &=
        -\varepsilon k\pare X.
        \label{eq:abstract-admissible-K-parity}
    \end{align}
    \end{subequations}
\end{itemize}
\end{definition}

\begin{rem}
      In the case $n=\mathtt{d}$, for every \(p\in\bR^{\mathtt d}\) we have that
    \begin{equation}\label{eq:abstract-admissible-derived-parities}
    \begin{aligned}
        h\pare\vartheta k\pare{p\cdot\vartheta}
        &\textnormal{ is odd in }\vartheta,
        \\
        h\pare\vartheta k'\pare{p\cdot\vartheta}
        &\textnormal{ is even in }\vartheta,
        \\
        h\pare\vartheta k''\pare{p\cdot\vartheta}
        &\textnormal{ is odd in }\vartheta.
    \end{aligned}
    \end{equation}
    \fine
\end{rem}

Fix once and for all a number \(r_{\mathrm{an}}\in\pare{0,r_k}\) strictly
inside the radius of convergence at the origin of the power series of \(k\),
and fix \(r_*\in\pare{0,r_{\mathrm{an}}}\).
All analytic expansions below are performed on
\(\av X\leq r_*<r_{\mathrm{an}}\); control of the full operator uses the
same bound for every difference quotient.

For an admissible pair as per \cref{def:abstract-admissible-angular-pair-Rd}, set for all $n\leq \mathtt{d}$
\begin{equation}\label{eq:abstract-singular-graph-operator-Rd}
\begin{aligned}
    \mathpzc t_{h,k}^{\pare{n}}\pare\eta\bra\omega\pare x
    &\defeq
    \pv\int_{\bR^{\mathtt d}}
    \frac{h\pare{\hat z}}{\av z^{n}}
    k\pare{\Delta_z\eta\pare x}
    \omega\pare{x-z}\dd z,
    \\
    \Delta_z\eta\pare x
    &\defeq
    \frac{\eta\pare x-\eta\pare{x-z}}{\av z},
    \qquad
    \hat z\defeq\frac z{\av z}\cdot
\end{aligned}
\end{equation}

\begin{lemma}[Skew-symmetric structure of admissible graph kernels in the critical case]
\label{lem:abstract-skew-adjoint-graph-kernel-Rd}
 For every real \(\eta\) for which
\cref{eq:abstract-singular-graph-operator-Rd} is defined the operator $\mathpzc t_{h,k}^{\pare{\mathtt d}}\pare\eta$ is $L^2\pare{\bR^{\mathtt{d}};\bR}$--skew-symmetric.
\end{lemma}

\begin{proof}
 We prove that for every 
\(f,g\in\mathcal S\pare{\bR^{\mathtt d};\bR}\), we have that
\begin{equation}\label{eq:abstract-skew-adjoint-graph-kernel-Rd}
    \psc{
        \mathpzc t_{h,k}^{\pare{\mathtt d}}\pare\eta\bra f}{
        g
    }_{L^2}
    =
    -
    \psc{
        f}{
        \mathpzc t_{h,k}^{\pare{\mathtt d}}\pare\eta\bra g
    }_{L^2},
\end{equation}
the generic assertion thus follows by density. 
For \(0<\epsilon<R<+\infty\), let
\[
    \mathpzc t_{h,k;\epsilon,R}^{\pare{\mathtt d}}
    \pare\eta\bra f\pare x
    \defeq
    \int_{\epsilon<\av{x-y}<R}
    \mathcal K_\eta\pare{x,y}f\pare y\dd y,
\]
where
\[
    \mathcal K_\eta\pare{x,y}
    \defeq
    \frac{h\pare{\widehat{x-y}}}{\av{x-y}^{\mathtt d}}
    k\pare{
        \frac{\eta\pare x-\eta\pare y}{\av{x-y}}
    }.
\]
The truncation removes the singularity at \(x=y\), so that the following
computations are justified by Fubini's theorem.
We first verify that \(\mathcal K_\eta\) is skew-symmetric.  Exchanging
\(x\) and \(y\) gives
\[
    \widehat{y-x}=-\widehat{x-y},
    \qquad
    \frac{\eta\pare y-\eta\pare x}{\av{y-x}}
    =
    -
    \frac{\eta\pare x-\eta\pare y}{\av{x-y}}\cdot
\]
Hence, by
\cref{eq:abstract-admissible-H-parity,eq:abstract-admissible-K-parity}, we have
\begin{equation}\label{sk_sym:calK}
    \begin{aligned}
    \mathcal K_\eta\pare{y,x}
    &=
    \frac{
        h\pare{-\widehat{x-y}}
    }{\av{x-y}^{\mathtt d}}
    k\pare{
        -
        \frac{\eta\pare x-\eta\pare y}{\av{x-y}}
    }
    \\
    &=
    \frac{
        \varepsilon h\pare{\widehat{x-y}}
    }{\av{x-y}^{\mathtt d}}
    \pare{-\varepsilon}
    k\pare{
        \frac{\eta\pare x-\eta\pare y}{\av{x-y}}
    }
    \\
    &=
    -\mathcal K_\eta\pare{x,y}.
\end{aligned}
\end{equation}
Since the truncated domain
\[
    \set{
        \pare{x,y}\in\bR^{\mathtt d}\times\bR^{\mathtt d}
        :
        \epsilon<\av{x-y}<R
    }
\]
is invariant under the exchange \(x\leftrightarrow y\), using \cref{sk_sym:calK}, we obtain
\begin{align}
    \psc{
        \mathpzc t_{h,k;\epsilon,R}^{\pare{\mathtt d}}
        \pare\eta\bra f}{
        g
    }_{L^2}
    &=
    \iint_{\epsilon<\av{x-y}<R}
    \mathcal K_\eta\pare{x,y}
    f\pare y g\pare x
    \dd y\dd x
    \\
    &=
    \iint_{\epsilon<\av{x-y}<R}
    \mathcal K_\eta\pare{y,x}
    f\pare x g\pare y
    \dd x\dd y
    \\
    &=
    -
    \iint_{\epsilon<\av{x-y}<R}
    \mathcal K_\eta\pare{x,y}
    f\pare x g\pare y
    \dd y\dd x
    \\
    &=
    -
    \psc{
        f}{
        \mathpzc t_{h,k;\epsilon,R}^{\pare{\mathtt d}}
        \pare\eta\bra g
    }_{L^2}.
\end{align}
Finally, by the definition of
\(\mathpzc t_{h,k}^{\pare{\mathtt d}}\pare\eta\) in
\cref{eq:abstract-singular-graph-operator-Rd}, we let
\(\epsilon\downarrow0\) and \(R\uparrow+\infty\).  Passing to the
principal-value limit in the preceding identity yields
\cref{eq:abstract-skew-adjoint-graph-kernel-Rd}.
\end{proof}

\begin{theorem}[Abstract structure-preserving paralinearization]
\label{thm:abstract-structure-preserving-paralinearization-Rd}
Let \(\pare{h,k}\) be an admissible homogeneous pair of critical order
\(-\mathtt d\) and vanishing order \(\ell\) in the sense of
\cref{def:abstract-admissible-angular-pair-Rd}. There exists
\(\varrho_0>0\) such that, for every \(\varrho\geq\varrho_0\), there exist
\(s_0>\mathtt d/2+1\), \(\mu\geq0\), and \(r_0>0\) such that, for every
\(s\geq s_0\), every
\(\eta\in H^{s+1+\mu}\pare{\bR^{\mathtt d};\bR}\), and 
\(\omega\in H^{s+\mu}\pare{\bR^{\mathtt d};\bR}\) satisfying
\(\norm\eta_{H^{s_0+1+\mu}}<r_0\), the operator defined in
\cref{eq:abstract-singular-graph-operator-Rd} admits the paralinearization
\begin{equation}\label{eq:abstract-structure-preserving-paralinearization-Rd}
    \mathpzc t_{h,k}^{\pare{\mathtt d}}\pare\eta\bra\omega
    =
    \OpBW{
        A^{h,k}_1\bra\omega
        \pare{\eta;x,\xi}
        +
        \vect{
            a_{\bra{0,\geq\ell}}\bra\omega\pare{\eta;x,\xi}
        }{
            \ii a_{\bra{-2,\geq\ell}}^{\bR;\omega}
            \pare{\eta;x,\xi}
        }
    }
    \cdot
    \vect{\eta}{\omega}
    +
    R^{\lin,\bR}_{\bra{\varrho,\geq\ell}}\pare{\eta,\omega}
    \cdot
    \vect{\eta}{\omega}.
\end{equation}
The objects in
\cref{eq:abstract-structure-preserving-paralinearization-Rd} are as follows.
\begin{itemize}

    \item The vector symbol
    \begin{equation}\label{eq:abstract-vector-symbol-Rd}
        A^{h,k}_1\bra\omega
        \pare{\eta;x,\xi}
        \defeq
        \vect{
            \omega\pare x
            a_1^{h,k;\eta}
            \pare{\nabla\eta\pare x,\xi}
        }{
            \ii
            a_0^{h,k;\omega}
            \pare{\nabla\eta\pare x,\xi}
        }
    \end{equation}
    has the following components.
    \begin{itemize}

        \item
        \(a_1^{h,k;\eta}\pare{\nabla\eta,\xi}
        \in\Gamma^1_{\geq\ell-1}\) is the real-valued symbol
        \begin{equation}\label{eq:abstract-explicit-principal-symbol-A1-Rd}
            a_1^{h,k;\eta}\pare{p,\xi}
            \defeq
            \frac\pi2
            \int_{\bS^{\mathtt d-1}}
            \av{\vartheta\cdot\xi}
            h\pare\vartheta
            k'\pare{p\cdot\vartheta}
            \dd\sigma\pare\vartheta.
        \end{equation}
        It is even in \(\xi\). In addition,
        \(a_1^{h,k;\eta}\pare{0,\xi}=0\) whenever \(\ell\geq2\);

        \item
        \(a_0^{h,k;\omega}\pare{\nabla\eta,\xi}
        -a_0^{h,k;\omega}\pare{0,\xi}
        \in\Gamma^0_{\geq\ell}\) is the real-valued symbol
        \begin{equation}\label{eq:abstract-explicit-principal-symbol-A0-Rd}
            a_0^{h,k;\omega}\pare{p,\xi}
            \defeq
            -\frac\pi2
            \int_{\bS^{\mathtt d-1}}
            \sgn\pare{\vartheta\cdot\xi}
            h\pare\vartheta
            k\pare{p\cdot\vartheta}
            \dd\sigma\pare\vartheta.
        \end{equation}
        It is odd in \(\xi\). In addition,
        \(a_0^{h,k;\omega}\pare{0,\xi}\) is the symbol of the flat
        homogeneous contribution, whereas its nonlinear part vanishes to
        order at least \(\ell\) at \(\eta=0\);

        \item 
        the two components satisfy
    \begin{equation}\label{eq:abstract-principal-symbol-differential-identity}
       \partial_{\xi} a_1^{h,k;\eta}\pare{p,\xi}
        =
        -\partial_p
        a_0^{h,k;\omega}\pare{p,\xi}.
    \end{equation}

    \end{itemize}

    \item The lower-order vector symbol
    \[
        \vect{
            a_{\bra{0,\geq\ell}}\bra\omega\pare{\eta;x,\xi}
        }{
            \ii a_{\bra{-2,\geq\ell}}^{\bR;\omega}
            \pare{\eta;x,\xi}
        }
    \]
    has the following components.
    \begin{itemize}

        \item $ a_{\bra{0,\geq\ell}}\bra\cdot
            \in
            \cL\pare{
                H^{s_0};
                \Gamma^0
            }$, and for every real \(\omega\),
        \begin{equation}\label{eq:abstract-transported-lower-symbol-reality-Rd}
            a_{\bra{0,\geq\ell}}\bra\omega\pare{\eta;x,-\xi}
            =
            \overline{
                a_{\bra{0,\geq\ell}}\bra\omega\pare{\eta;x,\xi}
            };
        \end{equation}
        If \(\ell\geq2\), the range improves to
        \(\Gamma^0_{\geq\ell-1}\);

        \item $  a_{\bra{-2,\geq\ell}}^{\bR;\omega}
            \in
            \Gamma^{-2}_{\geq\ell}$
        is real-valued and satisfies
        \begin{equation}\label{eq:abstract-imaginary-lower-symbol-parity-Rd}
            a_{\bra{-2,\geq\ell}}^{\bR;\omega}
            \pare{\eta;x,-\xi}
            =
            -
            a_{\bra{-2,\geq\ell}}^{\bR;\omega}
            \pare{\eta;x,\xi}.
        \end{equation}

    \end{itemize}

    \item The smoothing vector
    $R^{\lin,\bR}_{\bra{\varrho,\geq\ell}}\pare{\eta,\omega}$ is as per
    \cref{notation:paralinearization}; both components are real-to-real in the
    sense of \cref{def:smoothing}.

\end{itemize}
\end{theorem}

\begin{rem}[Skew-adjointness and the vanishing subprincipal symbol]
\label{rem:abstract-skew-adjoint-subprincipal-cancellation-Rd}
Two distinct structural facts enter the critical paralinearization.

First, the exact skew-adjointness of the critical quasilinear branch implies
that its complete Bony--Weyl symbol can be chosen purely imaginary. This is
proved in
\cref{lem:abstract-quasilinear-contribution-Rd} and yields the symbol in
\cref{eq:abstract-structure-preserving-paralinearization-Rd}.

Second, the absence of an order-\((-1)\) term is a stronger statement and does
not follow from skew-adjointness alone: a symbol of the form
\(\ii b_{-1}\), with \(b_{-1}\) real and odd in \(\xi\), would still define
a skew-adjoint Bony--Weyl operator. For \(\mathtt d\geq2\), the proof below
establishes the stronger cancellation explicitly. The first Bony--Weyl
correction generated by the principal homogeneous kernel is
\(-a_{-1}^{h,k}\), whereas the next Taylor coefficient of the graph
difference produces \(+a_{-1}^{h,k}\). Hence
\[
    -a_{-1}^{h,k}+a_{-1}^{h,k}=0,
\]
as shown in
\cref{eq:abstract-quasilinear-composition-subprincipal,%
eq:abstract-quasilinear-geometric-subprincipal}.
For \(\mathtt d=1\), the degree-zero principal multiplier is constant on
each half-line, so its derivatives in \(\xi\) vanish away from the origin;
the corresponding localized order-minus-one model also has vanishing
high-frequency multiplier. Thus no independent order-\((-1)\) symbol occurs
in any dimension.

The cancellation removes the subprincipal order, while skew-adjoint
completion controls the reality type of the entire remaining symbol. These
are complementary consequences of, respectively, the centered Weyl
structure and the opposite parities in
\cref{eq:abstract-admissible-H-parity,eq:abstract-admissible-K-parity}.
\fine
\end{rem}

\subsubsection{Local--tail and Bony reduction}

Let $\chi$ be the even radial cutoff fixed in \cref{def-chi}.  For fixed
$z\neq0$, write $\omega_z\pare x\defeq\omega\pare{x-z}$. Recall from \cref{eq:abstract-admissible-K-vanishing} that
\(\widetilde k\defeq k-k\pare0\), so that
\(\widetilde k\) vanishes to order at least \(\ell\) at the origin.
Define the flat contribution without multiplying a possibly undefined
principal value by zero:
\begin{equation}\label{eq:abstract-flat-contribution-Rd}
    \mathpzc t_{h,k}^{\mathrm{flat}}\bra\omega
    \defeq
    \begin{cases}
        k\pare0\displaystyle\pv\int_{\bR^{\mathtt d}}
        \dfrac{h\pare{\hat z}}{\av z^{\mathtt d}}
        \omega\pare{\bullet-z}\dd z,
        & k\pare0\neq0,
        \\
        0,
        & k\pare0=0.
    \end{cases}
\end{equation}
If \(k\pare0\neq0\),
\cref{eq:abstract-admissible-K-parity} forces \(h\) to be odd, so the
principal-value distribution in the first line is well-defined.

\begin{lemma}[Local-tail and paraproduct reduction]
\label{lem:abstract-common-local-tail-Bony-reduction-Rd}
With the same hypothesis of
\Cref{thm:abstract-structure-preserving-paralinearization-Rd}, one has
\begin{equation}\label{eq:abstract-common-Bony-decomposition}
\begin{aligned}
    \mathpzc t_{h,k}^{\pare{\mathtt d}}\pare\eta\bra\omega
    ={}&
    \mathpzc t_{h,k}^{\QL}\pare\eta\bra\omega
    +
    \mathpzc t_{h,k}^{\mathrm{tr}}\pare\eta\bra\omega
    +
    \mathpzc t_{h,k}^{R}\pare\eta\bra\omega
    \\
    &+
    \OpBW{
        f_{\geq\ell}^{\tail,\bR}
        \bra\omega\pare{\eta;x}
    }\eta
    +
    R^{\lin,\bR}_{\bra{\varrho,\geq\ell}}\pare{\eta,\omega}
    \cdot\vect{\eta}{\omega},
\end{aligned}
\end{equation}
where:
\begin{itemize}
    \item the quasilinear branch, including the flat contribution, is
    \begin{equation}\label{eq:abstract-common-quasilinear-branch}
    \begin{aligned}
        \mathpzc t_{h,k}^{\QL}\pare\eta\bra\omega
        \defeq{}&
        \mathpzc t_{h,k}^{\mathrm{flat}}\bra\omega
        \\
        &+
        \pv\int_{\bR^{\mathtt d}}
        \chi\pare z
        \frac{h\pare{\hat z}}{\av z^{\mathtt d}}
        \OpBW{
            \widetilde k\pare{\Delta_z\eta}
        }\omega_z\dd z;
    \end{aligned}
    \end{equation}

    \item the transported branch is
    \begin{equation}\label{eq:abstract-common-transported-branch}
        \mathpzc t_{h,k}^{\mathrm{tr}}\pare\eta\bra\omega
        \defeq
        \pv\int_{\bR^{\mathtt d}}
        \chi\pare z
        \frac{h\pare{\hat z}}{\av z^{\mathtt d}}
        \OpBW{\omega_z}
        \bra{\widetilde k\pare{\Delta_z\eta}}
        \dd z;
    \end{equation}

    \item the localized Bony remainder is
    \begin{equation}\label{eq:abstract-common-localized-Bony-remainder}
        \mathpzc t_{h,k}^{R}\pare\eta\bra\omega
        \defeq
        \pv\int_{\bR^{\mathtt d}}
        \chi\pare z
        \frac{h\pare{\hat z}}{\av z^{\mathtt d}}
        \RBW{
            \widetilde k\pare{\Delta_z\eta}
        }{\omega_z}
        \dd z;
    \end{equation}

    \item the tail function is
    \begin{equation}\label{eq:abstract-common-tail-function}
        f_{\geq\ell}^{\tail,\bR}
        \bra\omega\pare{\eta;x}
        \defeq
        \int_{\bR^{\mathtt d}}
        \pare{1-\chi\pare z}
        \frac{h\pare{\hat z}}{\av z^{\mathtt d+1}}
        k'\pare{\Delta_z\eta\pare x}
        \omega_z\pare x\dd z;
    \end{equation}

    \item the smoothing vector
    \(R^{\lin,\bR}_{\bra{\varrho,\geq\ell}}\pare{\eta,\omega}\)
    is as per \cref{notation:paralinearization}.
\end{itemize}
\end{lemma}

\begin{proof}
Let $r_0$ in
\cref{thm:abstract-structure-preserving-paralinearization-Rd}
be such that $0<r_0<r_\ast$, so that
$k\pare{\Delta_z\eta}$ is well defined.
Using $\widetilde k\defeq k-k\pare0$, we first split
\begin{equation}\label{eq:abstract-common-local-tail-splitting}
\begin{aligned}
    \mathpzc t_{h,k}^{\pare{\mathtt d}}
    \pare\eta\bra\omega
    ={}&
    \mathpzc t_{h,k}^{\mathrm{flat}}\bra\omega
    \\
    &+
    \pv\int_{\bR^{\mathtt d}}
    \chi\pare z
    \frac{h\pare{\hat z}}{\av z^{\mathtt d}}
    \widetilde k\pare{\Delta_z\eta}
    \omega_z\dd z
    \\
    &+
    \int_{\bR^{\mathtt d}}
    \pare{1-\chi\pare z}
    \frac{h\pare{\hat z}}{\av z^{\mathtt d}}
    \widetilde k\pare{\Delta_z\eta}
    \omega_z\dd z.
\end{aligned}
\end{equation}
We first treat the last term. Set
\[
    b_h\pare z
    \defeq
    \frac{h\pare{\hat z}}{\av z^{\mathtt d}}.
\]
Since $h\in\cC^\infty\pare{\bS^{\mathtt d-1}}$, for every multi-index
$\alpha$,
\begin{equation}\label{eq:abstract-common-tail-homogeneous-kernel-bound}
    \av{\partial_z^\alpha b_h\pare z}
    \lesssim_{\alpha,h}
    \av z^{-\mathtt d-\av\alpha}.
\end{equation}
Thus $b_h$ has far-field order $m=\mathtt d$. Moreover,
$\widetilde k$ vanishes to order at least $\ell\geq1$ by
\cref{eq:abstract-admissible-K-vanishing-order}; hence
$m+\ell=\mathtt d+\ell>\mathtt d$, which means
\begin{equation}
\label{eq:abstract-tail-integrability-calculation}
\int_1^\infty r^{\mathtt d-1-m-\ell}\dd r=\int_1^\infty r^{\mathtt -1-\ell}\dd r<\infty.
\end{equation}Therefore
\begin{equation}\label{eq:abstract-common-tail-kernel-class}
    k^{\tail,\bR}_{\bra{-\mathtt d,\geq\ell}}
    \pare{\eta;x,z}
    \defeq
    \pare{1-\chi\pare z}
    \frac{h\pare{\hat z}}{\av z^{\mathtt d}}
    \widetilde k\pare{\Delta_z\eta}
    \in
    \KF_{\textnormal{tail},\geq\ell}.
\end{equation}
Applying \cref{prop:tail-kernel-smoothing}, for every prescribed
$\varrho\geq0$,
\begin{equation}\label{eq:abstract-common-tail-operator}
\begin{aligned}
    \int_{\bR^{\mathtt d}}
    k^{\tail,\bR}_{\bra{-\mathtt d,\geq\ell}}
    \pare{\eta;x,z}\omega_z\dd z
    ={}&
    \OpBW{
        f_{\geq\ell}^{\tail,\bR}
        \bra\omega\pare{\eta;x}
    }\eta
    +
    R^{\lin,\bR}_{\bra\varrho}
    \pare{\eta,\omega}
    \cdot\vect{\eta}{\omega}.
\end{aligned}
\end{equation}

It remains to decompose the localized nonlinear term in
\cref{eq:abstract-common-local-tail-splitting}. For each fixed
$z\neq0$, \cref{lem:paraprod_Bony} gives
\begin{equation}\label{eq:abstract-common-local-Bony-identity}
\begin{aligned}
    \widetilde k\pare{\Delta_z\eta}\omega_z
    ={}&
    \OpBW{
        \widetilde k\pare{\Delta_z\eta}
    }\omega_z
    +
    \OpBW{\omega_z}
    \bra{
        \widetilde k\pare{\Delta_z\eta}
    }
    +
    \RBW{
        \widetilde k\pare{\Delta_z\eta}
    }{\omega_z}.
\end{aligned}
\end{equation}
Multiplying by
$\chi\pare z h\pare{\hat z}/\av z^{\mathtt d}$ and integrating in $z$
gives, by
\cref{eq:abstract-common-quasilinear-branch,%
eq:abstract-common-transported-branch,%
eq:abstract-common-localized-Bony-remainder},
\[
    \pv\int_{\bR^{\mathtt d}}
    \chi\pare z
    \frac{h\pare{\hat z}}{\av z^{\mathtt d}}
    \widetilde k\pare{\Delta_z\eta}
    \omega_z\dd z
    =
    \mathpzc t_{h,k}^{\QL}\pare\eta\bra\omega
    -
    \mathpzc t_{h,k}^{\mathrm{flat}}\bra\omega
    +
    \mathpzc t_{h,k}^{\mathrm{tr}}\pare\eta\bra\omega
    +
    \mathpzc t_{h,k}^{R}\pare\eta\bra\omega.
\]
Combining this identity with
\cref{eq:abstract-common-local-tail-splitting,%
eq:abstract-common-tail-operator}
proves \cref{eq:abstract-common-Bony-decomposition}.
\end{proof}

\begin{lemma}[Localized Bony remainder]
\label{lem:abstract-localized-Bony-remainder-Rd}
For every prescribed $\varrho\geq0$, the localized Bony remainder satisfies
\begin{equation}\label{eq:abstract-localized-Bony-remainder-smoothing}
    \mathpzc t_{h,k}^{R}\pare\eta\bra\omega
    =
    r^\bR_{\bra{\varrho,\geq\ell}}\pare\eta\omega
    .
\end{equation}
\end{lemma}

\begin{proof}
Set
\begin{align}\label{eq:abstract-quasilinear-local-building-blocks}
    \mathsf{v}\bra{\eta}\pare{x,z}
    \defeq
    \nabla\eta\pare x\cdot\hat z,
    &&
    \mathsf{h}\bra{\eta}\pare{x,z}
    \defeq
    \hat z^\intercal{\bm H}\eta\pare x\hat z.
\end{align}
By bilinearity of the Bony remainder,
\(\mathpzc t_{h,k}^{R}\pare\eta\bra\omega\) is the sum of
\begin{equation}\label{eq:abstract-product-Bony-remainder-decomposition}
\begin{aligned}
    &\pv\int_{\bR^{\mathtt d}}
    \chi\pare z
    \frac{h\pare{\hat z}}{\av z^{\mathtt d}}
    \RBW{\widetilde k\pare{\mathsf{v}\bra{\eta}}}{\omega\pare x}
    \dd z,
    \\
    &\pv\int_{\bR^{\mathtt d}}
    \chi\pare z
    \frac{h\pare{\hat z}}{\av z^{\mathtt d}}
    \RBW{\widetilde k\pare{\Delta_z\eta}}
    {\omega_z-\omega\pare x}
    \dd z,
    \\
    &\pv\int_{\bR^{\mathtt d}}
    \chi\pare z
    \frac{h\pare{\hat z}}{\av z^{\mathtt d}}
    \RBW{
        \widetilde k\pare{\Delta_z\eta}
        -
        \widetilde k\pare{\mathsf{v}\bra{\eta}}
    }{\omega\pare x}
    \dd z.
\end{aligned}
\end{equation}
By
\cref{eq:abstract-admissible-H-parity,%
eq:abstract-admissible-K-parity},
\[
    h\pare{-\vartheta}
    =
    \varepsilon h\pare\vartheta,
    \qquad
    \widetilde k\pare{-\mathsf{v}\bra{\eta}}
    =
    -\varepsilon\widetilde k\pare{\mathsf{v}\bra{\eta}}.
\]
Hence the first integrand in
\cref{eq:abstract-product-Bony-remainder-decomposition}
is odd in \(z\), and its integral over every symmetric annulus vanishes.

For the remaining two terms, the identities
\[
    \omega_z-\omega\pare x
    =
    -\av z
    \int_0^1
    \nabla\omega\pare{x-tz}\cdot\hat z
    \dd t
\]
and
\[
    \Delta_z\eta-\mathsf{v}\bra{\eta}
    =
    -\av z
    \int_0^1
    \pare{1-t}
    \hat z^\intercal
    {\bm H}\eta\pare{x-tz}
    \hat z
    \dd t
\]
give one power of \(\av z\). Moreover, by the mean value formula,
\[
\begin{aligned}
    \widetilde k\pare{\Delta_z\eta}
    -
    \widetilde k\pare{\mathsf{v}\bra{\eta}}
    =
    \pare{\Delta_z\eta-\mathsf{v}\bra{\eta}}
    \int_0^1
    \widetilde k'
    \pare{
        \mathsf{v}\bra{\eta}+t\pare{\Delta_z\eta-\mathsf{v}\bra{\eta}}
    }
    \dd t,
\end{aligned}
\]
and therefore, by
\cref{eq:abstract-admissible-K-vanishing,%
lem:composition_rules_z},
\begin{equation}\label{eq:abstract-product-Bony-difference-classes}
    \omega\longmapsto
    \omega_z-\omega\pare x
    \in
    \cL\pare{H^{s+1};\KF^1_{\loc,\geq 1}},
    \qquad
    \widetilde k\pare{\Delta_z\eta}
    -
    \widetilde k\pare{\mathsf{v}\bra{\eta}}
    \in
    \KF^1_{\loc,\geq\ell}.
\end{equation}
Thus
\cref{rem:smoothBONYremainder,lem:composition_rules_z},
using one additional derivative in the term containing
\(\omega_z-\omega\), place both nonzero operator-valued kernels in
\[
    \KR^{\varrho,1-\mathtt d}_{\loc,\geq\ell}.
\]
Since $1-\mathtt d>-\mathtt d$,
\cref{lem:integration_KR} yields
\[
    \mathpzc t_{h,k}^{R}\pare\eta\bra\omega
    =
    r^\bR_{\bra{\varrho,\geq\ell}}\pare\eta\omega,
    \qquad
    r^\bR_{\bra{\varrho,\geq\ell}}
    \in
    \cR^\varrho_{\geq\ell}.
\]
The vanishing order follows from
\cref{eq:abstract-admissible-K-vanishing}, and the operator is
real-to-real because the kernels and the Bony remainder are real-to-real.
\end{proof}

\subsubsection{Quasilinear contribution}

\begin{lemma}[Quasilinear contribution]
\label{lem:abstract-quasilinear-contribution-Rd}
For every prescribed \(\varrho\geq0\) and every dimension
\(\mathtt d\geq1\), the quasilinear branch in
\cref{eq:abstract-common-quasilinear-branch} satisfies
\begin{equation}\label{eq:abstract-quasilinear-skew-completion-Rd}
    \mathpzc t_{h,k}^{\QL}\pare\eta\bra\omega
    =
    \ii \ \OpBW{
            a_0^{h,k;\omega}\pare{\nabla\eta,\xi}
            +
            a_{\bra{-2,\geq\ell}}^{\bR;\omega}
            \pare{\eta;x,\xi}
    }\omega
    +
    r^{\textnormal{sk},\bR}_{\bra{\varrho,\geq\ell}}\pare\eta\omega.
\end{equation}
The objects in this equation have the following properties.
\begin{itemize}

    \item The symbol \(a_0^{h,k;\omega}\) is defined in
    \cref{eq:abstract-explicit-principal-symbol-A0-Rd}; it is real-valued
    and odd in \(\xi\). In particular, its nonlinear part belongs to
    \(\Gamma^0_{\geq\ell}\);

    \item $  a_{\bra{-2,\geq\ell}}^{\bR;\omega}
        \in
        \Gamma^{-2}_{\geq\ell}$
    is real-valued and odd in \(\xi\);

    \item $ r^{\textnormal{sk},\bR}_{\bra{\varrho,\geq\ell}}
        \in
        \cR^\varrho_{\geq\ell}$
    is real-to-real and skew-adjoint.

\end{itemize}
\end{lemma}

\begin{proof}
\begin{proofpart}[Local Taylor decomposition]
Recall $\mathsf{v}\bra{\eta}\pare{x,z}$ and $ \mathsf{h}\bra{\eta}\pare{x,z}$ defined in \cref{eq:abstract-quasilinear-local-building-blocks}.
Taylor's formula along the segment \(t\mapsto x-tz\) gives
\begin{equation}\label{eq:abstract-quasilinear-difference-Taylor-expansion}
\begin{aligned}
    \Delta_z\eta
    &=
    \mathsf{v}\bra{\eta}-\frac{\av z}{2}\mathsf{h}\bra{\eta}+\widetilde k^{\loc,\bR}_{\bra{2,\geq1}},
    \\
    \widetilde k^{\loc,\bR}_{\bra{2,\geq1}}\pare{\eta;x,z}
    &\defeq
    -\av z
    \int_0^1\pare{1-t}
    \hat z^\intercal
    \left[
        {\bm H}\eta\pare{x-tz}
        -
        {\bm H}\eta\pare x
    \right]
    \hat z\dd t.
\end{aligned}
\end{equation}
\Cref{eq:abstract-quasilinear-difference-Taylor-expansion,lem:kernel,lem:hsnorm} imply
\(\widetilde k^{\loc,\bR}_{\bra{2,\geq1}}\in\KF^2_{\loc,\geq1}\).
Write
\[
    \mathsf{w}\bra{\eta}\defeq\Delta_z\eta-\mathsf{v}\bra{\eta}
    =-\frac{\av z}{2}\mathsf{h}\bra{\eta}+\widetilde k^{\loc,\bR}_{\bra{2,\geq1}}
\]
and define
\begin{equation}\label{eq:abstract-quasilinear-composed-Taylor-remainder}
\begin{aligned}
    k^{\loc,\bR}_{\bra{2,\geq\ell}}\pare{\eta;x,z}
    \defeq{}&
    k\pare{\Delta_z\eta\pare x}
    -k\pare{\mathsf{v}\bra{\eta}\pare{x,z}}
    +\frac{\av z}{2}
    k'\pare{\mathsf{v}\bra{\eta}\pare{x,z}}\mathsf{h}\bra{\eta}\pare{x,z}
    \\
    ={}&
    k'\pare{\mathsf{v}\bra{\eta}} \widetilde k^{\loc,\bR}_{\bra{2,\geq1}}
    +\mathsf{w}\bra{\eta}^2\int_0^1
    \pare{1-t}k''\pare{\mathsf{v}\bra{\eta}+t\mathsf{w}\bra{\eta}}\dd t.
\end{aligned}
\end{equation}
By \cref{eq:abstract-admissible-K-vanishing}, the composition
\(k'\pare{\mathsf{v}\bra{\eta}}\) has background vanishing order at least \(\ell-1\), while
\(k''\pare{\mathsf{v}\bra{\eta}+t\mathsf{w}\bra{\eta}}\) has background vanishing order at least
\(\max\set{\ell-2,0}\). Hence every term in the last line of
\cref{eq:abstract-quasilinear-composed-Taylor-remainder} has local order two
and background vanishing order at least \(\ell\). The analytic composition
and product rules in
\cref{lem:kernel,lem:composition_rules_z} give
\begin{equation}\label{eq:abstract-quasilinear-composed-Taylor-class}
    k^{\loc,\bR}_{\bra{2,\geq\ell}}
    \in
    \KF^2_{\loc,\geq\ell},
    \qquad
    \chi\pare z
    \frac{h\pare{\hat z}}{\av z^{\mathtt d}}
    k^{\loc,\bR}_{\bra{2,\geq\ell}}
    \in
    \KF^{-\mathtt d+2}_{\geq\ell}.
\end{equation}
Accordingly,
\begin{equation}\label{eq:abstract-quasilinear-local-Taylor-splitting}
    \mathpzc t_{h,k}^{\QL}
    =
    \mathpzc t_{h,k}^{\QL,0}
    +\mathpzc t_{h,k}^{\QL,-1}
    +\mathpzc t_{h,k}^{\QL,\leq-2},
\end{equation}
where
\begin{subequations}\label{eq:abstract-quasilinear-local-Taylor-models}
\begin{align}
    \mathpzc t_{h,k}^{\QL,0}\pare\eta\bra\omega
    \defeq{}&
    \mathpzc t_{h,k}^{\mathrm{flat}}\bra\omega
    \\
    &+
    \pv\int_{\bR^{\mathtt d}}
    \chi\pare z
    \frac{h\pare{\hat z}}{\av z^{\mathtt d}}
    \OpBW{\widetilde k\pare{\mathsf{v}\bra{\eta}}}\omega_z\dd z,
    \label{eq:abstract-quasilinear-local-leading-model}
    \\
    \mathpzc t_{h,k}^{\QL,-1}\pare\eta\bra\omega
    \defeq{}&
    -\frac12
    \int_{\bR^{\mathtt d}}
    \chi\pare z
    \frac{h\pare{\hat z}}{\av z^{\mathtt d-1}}
    \OpBW{k'\pare{\mathsf{v}\bra{\eta}} \mathsf{h}\bra{\eta}}\omega_z\dd z,
    \label{eq:abstract-quasilinear-local-subprincipal-model}
    \\
    \mathpzc t_{h,k}^{\QL,\leq-2}\pare\eta\bra\omega
    \defeq{}&
    \int_{\bR^{\mathtt d}}
    \chi\pare z
    \frac{h\pare{\hat z}}{\av z^{\mathtt d}}
    \OpBW{k^{\loc,\bR}_{\bra{2,\geq\ell}}}\omega_z\dd z.
    \label{eq:abstract-quasilinear-local-lower-model}
\end{align}
\end{subequations}
We first analyze the leading and subprincipal branches, and then the lower
Taylor remainder.
\end{proofpart}

\begin{proofpart}[Analysis of \cref{eq:abstract-quasilinear-local-leading-model}]
For \(\mathtt d\geq2\), we prove
\begin{equation}\label{eq:abstract-quasilinear-order-zero-branch}
    \mathpzc t_{h,k}^{\QL,0}\pare\eta\bra\omega
    =
    \OpBW{
        \ii a_0^{h,k;\omega}\pare{\nabla\eta,\xi}
        -a_{-1}^{h,k}\pare{\eta;x,\xi}
        +a_{\bra{-2,\geq\ell}}\pare{\eta;x,\xi}
    }\omega
    +r^\bR_{\bra{\varrho,\geq\ell}}\pare\eta\omega,
\end{equation}
where
\[
    a_{\bra{-2,\geq\ell}}\in\Gamma^{-2}_{\geq\ell},
    \qquad
    r^\bR_{\bra{\varrho,\geq\ell}}\in\cR^\varrho_{\geq\ell}.
\]
We first remove the localization in the homogeneous model and then identify
the order \(0\) and order \(-1\) symbols in
\cref{eq:abstract-quasilinear-order-zero-branch}.

For every \(\alpha\in\bN^{\mathtt d}\) such that
\(\av\alpha\geq\ell\) and \(\kappa_{\av\alpha}\neq0\), set
\begin{equation}\label{eq:abstract-quasilinear-order-zero-monomial-data}
\begin{aligned}
    c_\alpha\pare{\eta;x}
    &\defeq
    \kappa_{\av\alpha}
    \binom{\av\alpha}{\alpha}
    \pare{\nabla\eta\pare x}^\alpha,
    \qquad
    h_\alpha^{\pare0}\pare z
    \defeq
    \frac{h\pare{\hat z}\hat z^\alpha}{\av z^{\mathtt d}}\cdot
\end{aligned}
\end{equation}
By \cref{eq:abstract-admissible-K-parity},
\(\kappa_{\av\alpha}\neq0\) implies
\(\pare{-1}^{\av\alpha}=-\varepsilon\). Hence
\[
    h\pare{-\vartheta}\pare{-\vartheta}^\alpha
    =
    \varepsilon\pare{-1}^{\av\alpha}
    h\pare\vartheta\vartheta^\alpha
    =
    -h\pare\vartheta\vartheta^\alpha.
\]
Thus the angular amplitude of \(h_\alpha^{\pare0}\) is odd and has zero
spherical mean. In particular, the critical principal-value distribution
in \cref{eq:abstract-quasilinear-order-zero-monomial-data} satisfies the
hypothesis of
\cref{item:abstract-homogeneous-angular-degree-zero}.

According to
\cref{def:localized-homogeneous-convolution-models}, define
\begin{equation}\label{eq:abstract-quasilinear-order-zero-localized-model}
    t_{h_\alpha^{\pare0}}^{\loc}\omega\pare x
    \defeq
    \pv\int_{\bR^{\mathtt d}}
    \chi\pare z
    \frac{h\pare{\hat z}\hat z^\alpha}{\av z^{\mathtt d}}
    \omega\pare{x-z}\dd z.
\end{equation}
Expanding the operator
\(\mathpzc t_{h,k}^{\QL,0}\pare\eta\) defined in
\cref{eq:abstract-quasilinear-local-leading-model} by means of the monomial
decomposition in
\cref{eq:abstract-quasilinear-order-zero-monomial-data} gives
\begin{equation}\label{eq:abstract-quasilinear-localized-leading-expansion}
\begin{aligned}
    \mathpzc t_{h,k}^{\QL,0}\pare\eta\bra\omega
    ={}&
    \mathpzc t_{h,k}^{\mathrm{flat}}\bra\omega
   +
    \sum_{\substack{\alpha\in\bN^{\mathtt d}\\
                     \av\alpha\geq\ell\\
                     \kappa_{\av\alpha}\neq0}}
    \OpBW{c_\alpha\pare{\eta;x}}
    \bra{t_{h_\alpha^{\pare0}}^{\loc}\omega}.
\end{aligned}
\end{equation}
Using the Fourier multiplier
\(\dot{\mathfrak m}_{h_\alpha^{\pare0}}\) from
\cref{eq:abstract-critical-angular-multiplier-Rd}, define
\begin{equation}\label{eq:abstract-quasilinear-global-monomial-expansion}
\begin{aligned}
    \mathpzc g_{h,k}^{\QL,0}\pare\eta\bra\omega
    \defeq{}&
    \mathpzc t_{h,k}^{\mathrm{flat}}\bra\omega
    +
    \sum_{\substack{\alpha\in\bN^{\mathtt d}\\
                     \av\alpha\geq\ell\\
                     \kappa_{\av\alpha}\neq0}}
    \OpBW{c_\alpha\pare{\eta;x}}
    \bra{\dot{\mathfrak m}_{h_\alpha^{\pare0}}\pare D\omega}.
\end{aligned}
\end{equation}
Comparing
\cref{eq:abstract-quasilinear-localized-leading-expansion}
and
\cref{eq:abstract-quasilinear-global-monomial-expansion},
and applying
\cref{lem:localized-homogeneous-model,rem:localized-homogeneous-model-coefficients}
term by term, and after summation for $r_0$ sufficiently small, we obtain
\begin{equation}\label{eq:abstract-quasilinear-cutoff-removal}
    \mathpzc t_{h,k}^{\QL,0}\pare\eta\bra\omega
    =
    \mathpzc g_{h,k}^{\QL,0}\pare\eta\bra\omega
    +r^\bR_{\bra{\varrho,\geq\ell}}\pare\eta\omega,
    \qquad
    r^\bR_{\bra{\varrho,\geq\ell}}\in\cR^\varrho_{\geq\ell}.
\end{equation}
Thus \cref{eq:abstract-quasilinear-cutoff-removal} identifies the smoothing
remainder in \cref{eq:abstract-quasilinear-order-zero-branch}; it remains
to identify the symbol of
\(\mathpzc g_{h,k}^{\QL,0}\pare\eta\) through order \(-1\).

{
Define, for every \(\mathtt d\geq1\), the first quasilinear
subprincipal symbol by
\begin{equation}\label{eq:abstract-quasilinear-subprincipal-symbol-Rd}
    a_{-1}^{h,k}\pare{\eta;x,\xi}
    \defeq
    \frac12
    \pare{\nabla_x\cdot\nabla_\xi}
    a_0^{h,k;\omega}\pare{\nabla\eta\pare x,\xi}.
\end{equation}
For \(\mathtt d\geq2\), it admits the explicit representation
\begin{equation}\label{eq:abstract-quasilinear-subprincipal-equatorial-Rd}
    a_{-1}^{h,k}\pare{\eta;x,\xi}
    =
    -\frac\pi2
    \int_{\bS^{\mathtt d-1}}
    \delta\pare{\vartheta\cdot\xi}
    h\pare\vartheta
    k'\pare{\nabla\eta\pare x\cdot\vartheta}
    \pare{
        \vartheta^\intercal{\bm H}\eta\pare x\vartheta
    }
    \dd\sigma\pare\vartheta.
\end{equation}
Indeed, in the sense of tempered distributions in \(\xi\),
\[
    \partial_{\xi_i}
    \sgn\pare{\vartheta\cdot\xi}
    =
    2\vartheta_i
    \delta\pare{\vartheta\cdot\xi}.
\]
Therefore, differentiating
\cref{eq:abstract-explicit-principal-symbol-A0-Rd} and evaluating at
\(p=\nabla\eta\pare x\), we obtain
\[
    \pare{\nabla_x\cdot\nabla_\xi}
    a_0^{h,k;\omega}\pare{\nabla\eta\pare x,\xi}
    =
    -\pi
    \int_{\bS^{\mathtt d-1}}
    \delta\pare{\vartheta\cdot\xi}
    h\pare\vartheta
    k'\pare{\nabla\eta\pare x\cdot\vartheta}
    \pare{
        \vartheta^\intercal{\bm H}\eta\pare x\vartheta
    }
    \dd\sigma\pare\vartheta,
\]
which proves
\cref{eq:abstract-quasilinear-subprincipal-equatorial-Rd}
by \cref{eq:abstract-quasilinear-subprincipal-symbol-Rd}.
For \(\mathtt d=1\),
\cref{eq:abstract-explicit-principal-symbol-A0-Rd} shows that
\(a_0^{h,k;\omega}\pare{p,\cdot}\) is constant on each half-line, and hence
\(a_{-1}^{h,k}\equiv0\).
}
Applying \cref{prop:composition_paradiff} to
\cref{eq:abstract-quasilinear-global-monomial-expansion} and summing the
monomial expansion gives, at orders $0$ and $-1,$
\begin{equation}\label{eq:abstract-quasilinear-composition-subprincipal}
\begin{aligned}
    a_0^{h,k;\omega}\pare{0,\xi}
    +
    \sum_{\substack{\alpha\in\bN^{\mathtt d}\\
                     \av\alpha\geq\ell\\
                     \kappa_{\av\alpha}\neq0}}
    c_\alpha
    \dot{\mathfrak m}_\alpha^{\pare0}
    &=
    a_0^{h,k;\omega}\pare{\nabla\eta,\xi},
    \\
    -\frac12
    \sum_{\substack{\alpha\in\bN^{\mathtt d}\\
                     \av\alpha\geq\ell\\
                     \kappa_{\av\alpha}\neq0}}
    \nabla_xc_\alpha\cdot
    \nabla_\xi\dot{\mathfrak m}_\alpha^{\pare0}
    &=
    -\frac12
    \pare{\nabla_x\cdot\nabla_\xi}
    a_0^{h,k;\omega}
    =
    -a_{-1}^{h,k}.
\end{aligned}
\end{equation}
The first identity in
\cref{eq:abstract-quasilinear-composition-subprincipal}, together with
\(\dot{\mathfrak m}_{h_\alpha^{\pare0}}
=\ii\dot{\mathfrak m}_\alpha^{\pare0}\), gives the principal term
\[
    \ii a_0^{h,k;\omega}\pare{\nabla\eta,\xi}
\]
in \cref{eq:abstract-quasilinear-order-zero-branch}, while the second
identity gives its subprincipal term
\[
    -a_{-1}^{h,k}\pare{\eta;x,\xi}.
\]
The remaining symbolic contribution is
\(a_{\bra{-2,\geq\ell}}\in\Gamma^{-2}_{\geq\ell}\).
Together with
\cref{eq:abstract-quasilinear-cutoff-removal}, this proves
\cref{eq:abstract-quasilinear-order-zero-branch}.
\end{proofpart}

\begin{proofpart}[Analysis of \cref{eq:abstract-quasilinear-local-subprincipal-model}]
For \(\mathtt d\geq2\), we prove
\begin{equation}\label{eq:abstract-quasilinear-order-minus-one-branch}
    \mathpzc t_{h,k}^{\QL,-1}\pare\eta\bra\omega
    =
    \OpBW{
        a_{-1}^{h,k}\pare{\eta;x,\xi}
        +a_{\bra{-2,\geq\ell}}\pare{\eta;x,\xi}
    }\omega
    +r^\bR_{\bra{\varrho,\geq\ell}}\pare\eta\omega,
\end{equation}
where
\[
    a_{\bra{-2,\geq\ell}}\in\Gamma^{-2}_{\geq\ell},
    \qquad
    r^\bR_{\bra{\varrho,\geq\ell}}\in\cR^\varrho_{\geq\ell}.
\]
We first remove the localization and then identify the order \(-1\)
symbol in \cref{eq:abstract-quasilinear-order-minus-one-branch}.

For every \(\alpha\in\bN^{\mathtt d}\) such that
\(\av\alpha+1\geq\ell\) and \(\kappa_{\av\alpha+1}\neq0\), and for
\(1\leq j,k\leq\mathtt d\), set
\begin{equation}\label{eq:abstract-quasilinear-order-minus-one-monomial-data}
\begin{aligned}
    d_{\alpha,j,k}\pare{\eta;x}
    &\defeq
    -\frac12
    \pare{\av\alpha+1}\kappa_{\av\alpha+1}
    \binom{\av\alpha}{\alpha}
    \pare{\nabla\eta\pare x}^\alpha
    \partial_{jk}\eta\pare x,
    \\
    h_{\alpha,j,k}^{\pare{-1}}\pare z
    &\defeq
    \frac{
        h\pare{\hat z}\hat z^\alpha\hat z_j\hat z_k
    }{
        \av z^{\mathtt d-1}
    }\cdot
\end{aligned}
\end{equation}
By \cref{eq:abstract-admissible-K-parity},
\(\kappa_{\av\alpha+1}\neq0\) implies
\(\pare{-1}^{\av\alpha}=\varepsilon\). Hence
\[
    h\pare{-\vartheta}
    \pare{-\vartheta}^\alpha
    \pare{-\vartheta_j}\pare{-\vartheta_k}
    =
    h\pare\vartheta
    \vartheta^\alpha\vartheta_j\vartheta_k.
\]
Thus the angular amplitude of \(h_{\alpha,j,k}^{\pare{-1}}\) is even and
satisfies the hypothesis of
\cref{item:abstract-homogeneous-angular-degree-minus-one}.

According to
\cref{def:localized-homogeneous-convolution-models}, define
\begin{equation}\label{eq:abstract-quasilinear-order-minus-one-localized-model}
    t_{h_{\alpha,j,k}^{\pare{-1}}}^{\loc}\omega\pare x
    \defeq
    \int_{\bR^{\mathtt d}}
    \chi\pare z
    \frac{
        h\pare{\hat z}\hat z^\alpha\hat z_j\hat z_k
    }{
        \av z^{\mathtt d-1}
    }
    \omega\pare{x-z}\dd z.
\end{equation}
Expanding the operator
\(\mathpzc t_{h,k}^{\QL,-1}\pare\eta\) defined in
\cref{eq:abstract-quasilinear-local-subprincipal-model} by means of the
monomial decomposition in
\cref{eq:abstract-quasilinear-order-minus-one-monomial-data} gives
\begin{equation}\label{eq:abstract-quasilinear-localized-subprincipal-expansion}
    \mathpzc t_{h,k}^{\QL,-1}\pare\eta\bra\omega
    =
    \sum_{\substack{\alpha\in\bN^{\mathtt d}\\
                     \av\alpha+1\geq\ell\\
                     \kappa_{\av\alpha+1}\neq0}}
    \sum_{j,k=1}^{\mathtt d}
    \OpBW{d_{\alpha,j,k}\pare{\eta;x}}
    \bra{
        t_{h_{\alpha,j,k}^{\pare{-1}}}^{\loc}\omega
    }.
\end{equation}
Using the Fourier multiplier
\(\dot{\mathfrak m}_{h_{\alpha,j,k}^{\pare{-1}}}\) defined in
\cref{eq:homogeneous-model-multiplier}, define
\begin{equation}\label{eq:abstract-quasilinear-global-subprincipal-model}
    \mathpzc g_{h,k}^{\QL,-1}\pare\eta\bra\omega
    \defeq
    \sum_{\substack{\alpha\in\bN^{\mathtt d}\\
                     \av\alpha+1\geq\ell\\
                     \kappa_{\av\alpha+1}\neq0}}
    \sum_{j,k=1}^{\mathtt d}
    \OpBW{d_{\alpha,j,k}\pare{\eta;x}}
    \bra{
        \dot{\mathfrak m}_{h_{\alpha,j,k}^{\pare{-1}}}
        \pare D\omega
    }.
\end{equation}
Assume \(\mathtt d\geq2\). Comparing
\cref{eq:abstract-quasilinear-localized-subprincipal-expansion}
and
\cref{eq:abstract-quasilinear-global-subprincipal-model},
set, for every admissible \(\alpha,j,k\),
\[
    \mathfrak s_{\alpha,j,k}\pare\xi
    \defeq
    \mathfrak m_{\chi h_{\alpha,j,k}^{\pare{-1}}}\pare\xi
    -
    \dot{\mathfrak m}_{h_{\alpha,j,k}^{\pare{-1}}}\pare\xi .
\]
By \cref{lem:localized-homogeneous-model},
\[
    \mathfrak s_{\alpha,j,k}
    \in
    \bigcap_{M\in\bN}\widetilde\Gamma^{-M}.
\]
Moreover,
\(\mathfrak m_{\chi h_{\alpha,j,k}^{\pare{-1}}}\) is bounded, whereas
\(\dot{\mathfrak m}_{h_{\alpha,j,k}^{\pare{-1}}}\) is homogeneous of
degree \(-1\). Hence
\[
    \av{
        \mathfrak s_{\alpha,j,k}\pare\xi
    }
    \lesssim
    1+\av\xi^{-1},
    \qquad
    0<\av\xi\leq1.
\]
The high-frequency contribution $\mathfrak s_{\alpha,j,k}^\sharp$ is therefore arbitrarily smoothing.
It remains to consider \(\mathfrak s_{\alpha,j,k}^\flat\pare D\). Since \[ \av{ \mathfrak s_{\alpha,j,k}^\flat\pare\xi } \lesssim \av\xi^{-1}, \qquad 0<\av\xi\lesssim1, \] the support condition in \cref{BonyWeyl,rela:para} and Cauchy--Schwarz give \[ \av{\mathfrak s_{\alpha,j,k}^\flat\pare\varsigma} \int_{\av\theta\lesssim\av\varsigma} \av{ \widehat{ d_{\alpha,j,k}\pare\eta }\pare\theta } \dd\theta \lesssim \av\varsigma^{\frac{\mathtt d}{2}-1} \norm{ d_{\alpha,j,k}\pare\eta }_{L^2} \lesssim \norm{ d_{\alpha,j,k}\pare\eta }_{L^2}, \] for \(\mathtt d\geq2\). Schur's lemma therefore gives the \(L^2\)-boundedness of the low-frequency contribution. Since, by \cref{rela:para}, it is supported in a fixed low-frequency region, this contribution belongs to \(\cR^\varrho_{\geq\ell}\) for every \(\varrho\geq0\).
Summation therefore gives
\begin{equation}\label{eq:abstract-quasilinear-subprincipal-cutoff-removal}
    \mathpzc t_{h,k}^{\QL,-1}\pare\eta\bra\omega
    =
    \mathpzc g_{h,k}^{\QL,-1}\pare\eta\bra\omega
    +r^\bR_{\bra{\varrho,\geq\ell}}\pare\eta\omega,
    \qquad
    r^\bR_{\bra{\varrho,\geq\ell}}\in\cR^\varrho_{\geq\ell}.
\end{equation}
Thus \cref{eq:abstract-quasilinear-subprincipal-cutoff-removal}
identifies the smoothing remainder in
\cref{eq:abstract-quasilinear-order-minus-one-branch}; it remains to
identify the symbol of
\(\mathpzc g_{h,k}^{\QL,-1}\pare\eta\) through order \(-1\).

By
\cref{item:abstract-homogeneous-angular-degree-minus-one}, for every
admissible \(\alpha,j,k\),
\[
    \dot{\mathfrak m}_{h_{\alpha,j,k}^{\pare{-1}}}
    =
    \dot{\mathfrak m}_{\alpha,j,k}^{\pare{-1}},
    \qquad
    \dot{\mathfrak m}_{\alpha,j,k}^{\pare{-1}}\pare\xi
    \defeq
    \pi
    \int_{\bS^{\mathtt d-1}}
    \delta\pare{\vartheta\cdot\xi}
    h\pare\vartheta
    \vartheta^\alpha\vartheta_j\vartheta_k
    \dd\sigma.
\]
Applying \cref{prop:composition_paradiff} to
\cref{eq:abstract-quasilinear-global-subprincipal-model} and summing the
monomial expansion gives
\begin{equation}\label{eq:abstract-quasilinear-geometric-subprincipal}
\begin{aligned}
    &\sum_{\substack{\alpha\in\bN^{\mathtt d}\\
                     \av\alpha+1\geq\ell\\
                     \kappa_{\av\alpha+1}\neq0}}
    \sum_{j,k=1}^{\mathtt d}
    d_{\alpha,j,k}
    \dot{\mathfrak m}_{\alpha,j,k}^{\pare{-1}}
    =
    -\frac\pi2
    \int_{\bS^{\mathtt d-1}}
    \delta\pare{\vartheta\cdot\xi}
    h\pare\vartheta
    k'\pare{\nabla\eta\cdot\vartheta}
    \pare{\vartheta^\intercal{\bm H}\eta\vartheta}
    \dd\sigma
    =
    a_{-1}^{h,k},
\end{aligned}
\end{equation}
where the last equality follows from
\cref{eq:abstract-quasilinear-subprincipal-equatorial-Rd}.
Thus \cref{eq:abstract-quasilinear-geometric-subprincipal} gives the
order \(-1\) term
\[
    a_{-1}^{h,k}\pare{\eta;x,\xi}
\]
in \cref{eq:abstract-quasilinear-order-minus-one-branch}.

Since
\(\dot{\mathfrak m}_{\alpha,j,k}^{\pare{-1}}\) has order \(-1\), all
non-principal terms in the composition have order at most \(-2\).
Hence the remaining symbolic contribution in
\cref{eq:abstract-quasilinear-order-minus-one-branch} is
\[
    a_{\bra{-2,\geq\ell}}
    \in
    \Gamma^{-2}_{\geq\ell}.
\]
Together with
\cref{eq:abstract-quasilinear-subprincipal-cutoff-removal}, this proves
\cref{eq:abstract-quasilinear-order-minus-one-branch}.
\end{proofpart}

\begin{proofpart}[Exceptional one-dimensional case]
Let \(\mathtt d=1\). By
\cref{eq:abstract-explicit-principal-symbol-A0-Rd},
\(a_0^{h,k;\omega}\pare{p,\xi}\) is constant on each of the two half-lines
\(\xi>0\) and \(\xi<0\). Hence
\[
    \partial_\xi^q
    a_0^{h,k;\omega}\pare{p,\xi}
    =
    0,
    \qquad
    \xi\neq0,
    \quad q\geq1,
\]
and therefore, by
\cref{eq:abstract-quasilinear-subprincipal-symbol-Rd} $ a_{-1}^{h,k}\equiv0.$
Applying \cref{prop:composition_paradiff} to
\cref{eq:abstract-quasilinear-global-monomial-expansion}, all non-principal
composition terms are supported in the fixed transition annulus of the
high-frequency extension and are therefore smoothing. Together with
\cref{eq:abstract-quasilinear-cutoff-removal}, this gives
\begin{equation}\label{eq:abstract-quasilinear-order-zero-branch-d1}
    \mathpzc t_{h,k}^{\QL,0}\pare\eta\bra\omega
    =
    \OpBW{
        \ii a_0^{h,k;\omega}\pare{\partial_x\eta,\xi}
    }\omega
    +
    r^\bR_{\bra{\varrho,\geq\ell}}\pare\eta\omega,
    \qquad
    r^\bR_{\bra{\varrho,\geq\ell}}\in\cR^\varrho_{\geq\ell}.
\end{equation}

For the order-minus-one branch,
\cref{eq:abstract-quasilinear-order-minus-one-monomial-data} reduces to
\[
    h_{\alpha,1,1}^{\pare{-1}}\pare z
    =
    h\pare{\hat z}\hat z^\alpha .
\]
Since \(\kappa_{\alpha+1}\neq0\),
\cref{eq:abstract-admissible-K-parity} gives
\(\pare{-1}^{\alpha}=\varepsilon\), and hence
\[
    h\pare{-\vartheta}\pare{-\vartheta}^{\alpha}
    =
    h\pare\vartheta\vartheta^\alpha,
    \qquad
    \vartheta\in\bS^0.
\]
Thus \(h_{\alpha,1,1}^{\pare{-1}}\) is constant on
\(\bR\setminus\set0\). By
\cref{eq:abstract-one-dimensional-degree-zero-model}, its complete
high-frequency multiplier vanishes. Therefore
\cref{eq:abstract-quasilinear-subprincipal-cutoff-removal} gives
\begin{equation}\label{eq:abstract-quasilinear-order-minus-one-branch-d1}
    \mathpzc t_{h,k}^{\QL,-1}\pare\eta\bra\omega
    =
    r^\bR_{\bra{\varrho,\geq\ell}}\pare\eta\omega,
    \qquad
    r^\bR_{\bra{\varrho,\geq\ell}}\in\cR^\varrho_{\geq\ell}.
\end{equation}
\end{proofpart}
\begin{proofpart}[Lower Taylor remainder]
Since \(-\mathtt d+2>-\mathtt d\), the generalized
\cref{prop:int-opbw} applied to
\cref{eq:abstract-quasilinear-composed-Taylor-class} gives
\begin{equation}\label{eq:abstract-quasilinear-Taylor-remainder-symbol}
    \mathpzc t_{h,k}^{\QL,\leq-2}\pare\eta\bra\omega
    =
    \OpBW{
        a_{\bra{-2,\geq\ell}}\pare{\eta;x,\xi}
    }\omega
    +r^\bR_{\bra{\varrho,\geq\ell}}\pare\eta\omega,
\end{equation}
where
\[
    a_{\bra{-2,\geq\ell}}\in\Gamma^{-2}_{\geq\ell},
    \qquad
    r^\bR_{\bra{\varrho,\geq\ell}}\in\cR^\varrho_{\geq\ell}.
\]
Because the kernel in
\cref{eq:abstract-quasilinear-composed-Taylor-class} is real,
\cref{prop:int-opbw} also gives
\[
    a_{\bra{-2,\geq\ell}}\pare{\eta;x,-\xi}
    =
    \overline{
        a_{\bra{-2,\geq\ell}}\pare{\eta;x,\xi}
    },
\]
and a real-to-real remainder.

\end{proofpart}

\begin{proofpart}[Preliminary synthesis]
For \(\mathtt d\geq2\), combining
\cref{eq:abstract-quasilinear-local-Taylor-splitting,%
eq:abstract-quasilinear-order-zero-branch,%
eq:abstract-quasilinear-order-minus-one-branch,%
eq:abstract-quasilinear-Taylor-remainder-symbol}
and cancelling the two opposite order-\((-1)\) contributions gives a principal
symbol \(\ii a_0^{h,k;\omega}\) and generic terms of order at most \(-2\).
For \(\mathtt d=1\), the same representation follows by combining
\cref{eq:abstract-quasilinear-local-Taylor-splitting,%
eq:abstract-quasilinear-order-zero-branch-d1,%
eq:abstract-quasilinear-order-minus-one-branch-d1,%
eq:abstract-quasilinear-Taylor-remainder-symbol}.
Thus, after relabeling the generic symbols and remainders,
\begin{equation}\label{eq:abstract-quasilinear-preliminary-symbol-Rd}
    \mathpzc t_{h,k}^{\QL}\pare\eta\bra\omega
    =
    \OpBW{
        \ii a_0^{h,k;\omega}\pare{\nabla\eta,\xi}
        +
        \widetilde a_{\bra{-2,\geq\ell}}\pare{\eta;x,\xi}
    }\omega
    +
    \widetilde r^\bR_{\bra{\varrho,\geq\ell}}\pare\eta\omega,
\end{equation}
where
\begin{equation}\label{eq:abstract-quasilinear-preliminary-lower-classes-Rd}
    \widetilde a_{\bra{-2,\geq\ell}}
    \in
    \Gamma^{-2}_{\geq\ell},
    \qquad
    \widetilde r^\bR_{\bra{\varrho,\geq\ell}}
    \in
    \cR^\varrho_{\geq\ell}.
\end{equation}
Since the kernels and coefficients above are real-valued, the reality
statements in \cref{lem:localized-homogeneous-model,prop:int-opbw} and the
termwise symmetry of the Bony--Weyl expansion
\cref{eq:compositionsymbol1} allow the generic terms in
\cref{eq:abstract-quasilinear-preliminary-symbol-Rd} to be chosen so that
\begin{equation}\label{eq:abstract-quasilinear-preliminary-reality-Rd}
    \widetilde a_{\bra{-2,\geq\ell}}\pare{\eta;x,-\xi}
    =
    \overline{
        \widetilde a_{\bra{-2,\geq\ell}}\pare{\eta;x,\xi}
    },
\end{equation}
and \(\widetilde r^\bR_{\bra{\varrho,\geq\ell}}\) is real-to-real.
\end{proofpart}

\begin{proofpart}[Skew-adjoint completion]
    \begin{proofsubpart}[Exact skew-adjointness]
For \(z\neq0\), set
\begin{equation}\label{eq:abstract-quasilinear-skew-building-blocks-Rd}
\begin{aligned}
    \tau_z\omega\pare x
    &\defeq
    \omega\pare{x-z},
    \\
    \mathfrak h\pare z
    &\defeq
    \chi\pare z
    \frac{h\pare{\hat z}}{\av z^{\mathtt d}},
    \\
    a_z\pare{\eta;x}
    &\defeq
    \widetilde k\pare{\Delta_z\eta\pare x}.
\end{aligned}
\end{equation}
Evaluating
\cref{eq:abstract-admissible-K-parity} at the origin gives
\(\pare{1+\varepsilon}k\pare0=0\). Consequently,
\begin{equation}\label{eq:abstract-tilde-k-derived-parity-Rd}
    \widetilde k\pare{-X}
    =
    -\varepsilon\widetilde k\pare X.
\end{equation}
Since \(\chi\) is even,
\begin{equation}\label{eq:abstract-local-critical-kernel-parity-Rd}
    \mathfrak h\pare{-z}
    =
    \varepsilon\mathfrak h\pare z.
\end{equation}
Moreover,
\begin{equation}\label{eq:abstract-quasilinear-shifted-parity-Rd}
\begin{aligned}
    \Delta_{-z}\eta\pare{x-z}
    &=
    \frac{\eta\pare{x-z}-\eta\pare x}{\av z}
    =
    -\Delta_z\eta\pare x,
    \\
    a_{-z}\pare{\eta;x-z}
    &=
    -\varepsilon a_z\pare{\eta;x}.
\end{aligned}
\end{equation}

The definition of Bony--Weyl quantization in \cref{BonyWeyl} is covariant
under translations. Since \(a_z\) is real-valued,
\cref{lem:adjoint} gives
\begin{equation}\label{eq:abstract-quasilinear-translation-adjoint-Rd}
\begin{aligned}
    \pare{\OpBW{a_z}\tau_z}^{\star}
    &=
    \tau_{-z}\OpBW{a_z}
    =
    \OpBW{a_z\pare{\eta;\cdot+z}}\tau_{-z}.
\end{aligned}
\end{equation}
For \(0<\delta<R<+\infty\), define the truncated localized branch by
\[
    \mathpzc t_{h,k;\delta,R}^{\QL,\loc}\pare\eta
    \defeq
    \int_{\delta<\av z<R}
    \mathfrak h\pare z
    \OpBW{a_z}\tau_z
    \dd z.
\]
Using
\cref{eq:abstract-quasilinear-translation-adjoint-Rd}, then changing
\(z\) into \(-z\), and finally applying
\cref{eq:abstract-local-critical-kernel-parity-Rd,%
eq:abstract-quasilinear-shifted-parity-Rd}, we obtain
\begin{align}
    \pare{
        \mathpzc t_{h,k;\delta,R}^{\QL,\loc}\pare\eta
    }^{\star}
    &=
    \int_{\delta<\av z<R}
    \mathfrak h\pare z
    \OpBW{a_z\pare{\eta;\cdot+z}}
    \tau_{-z}
    \dd z
    \\
    &=
    \int_{\delta<\av z<R}
    \mathfrak h\pare{-z}
    \OpBW{a_{-z}\pare{\eta;\cdot-z}}
    \tau_z
    \dd z
    \\
    &=
    -
    \int_{\delta<\av z<R}
    \mathfrak h\pare z
    \OpBW{a_z}
    \tau_z
    \dd z
    \\
    &=
    -
    \mathpzc t_{h,k;\delta,R}^{\QL,\loc}\pare\eta.
\end{align}
Passing to the principal-value limit gives the same identity for the
localized contribution in
\cref{eq:abstract-common-quasilinear-branch}.

It remains to consider the flat contribution in
\cref{eq:abstract-flat-contribution-Rd}. If \(k\pare0=0\), it vanishes.
If \(k\pare0\neq0\), then
\cref{eq:abstract-admissible-K-parity} forces \(\varepsilon=-1\), and
hence \(h\) is odd. Therefore
\cref{item:abstract-homogeneous-angular-degree-zero} identifies the flat
operator with a purely imaginary Fourier multiplier, so it is
skew-adjoint. We have proved
\begin{equation}\label{eq:abstract-quasilinear-exact-skew-adjointness-Rd}
    \pare{
        \mathpzc t_{h,k}^{\QL}\pare\eta
    }^{\star}
    =
    -
    \mathpzc t_{h,k}^{\QL}\pare\eta.
\end{equation}
\end{proofsubpart}

\begin{proofsubpart}[Projection of the lower-order terms]
We now combine the preliminary representation
\cref{eq:abstract-quasilinear-preliminary-symbol-Rd} with the exact
skew-adjointness identity
\cref{eq:abstract-quasilinear-exact-skew-adjointness-Rd}.
Every smoothing remainder entering
\cref{eq:abstract-quasilinear-preliminary-symbol-Rd} is a finite sum of
terms of the following two types: smoothing Fourier multipliers furnished
by \cref{lem:localized-homogeneous-model}, Bony--Weyl composition
remainders furnished by \cref{prop:composition_paradiff}. In each case, adjunction only
conjugates the symbol or exchanges the two Fourier variables of the kernel;
the symbol bounds, support properties, tame estimates, and background
vanishing orders are unchanged. Therefore
\begin{equation}\label{eq:abstract-quasilinear-adjoint-remainder-class-Rd}
    \widetilde r^{\star,\bR}_{\bra{\varrho,\geq\ell}}
    \in
    \cR^\varrho_{\geq\ell}.
\end{equation}
Since \(a_0^{h,k;\omega}\) is real-valued,
\cref{lem:adjoint,eq:abstract-quasilinear-exact-skew-adjointness-Rd} and
\cref{eq:abstract-quasilinear-preliminary-symbol-Rd} yield
\begin{align}
    \mathpzc t_{h,k}^{\QL}\pare\eta
    &=
    \frac12
    \left(
        \mathpzc t_{h,k}^{\QL}\pare\eta
        -
        \pare{
            \mathpzc t_{h,k}^{\QL}\pare\eta
        }^{\star}
    \right)
    \\
    &=
    \OpBW{
        \ii a_0^{h,k;\omega}\pare{\nabla\eta,\xi}
        +
        \frac{
            \widetilde a_{\bra{-2,\geq\ell}}
            -
            \overline{\widetilde a_{\bra{-2,\geq\ell}}}
        }2
    }
    +
    \frac{
        \widetilde r^\bR_{\bra{\varrho,\geq\ell}}
        -
        \widetilde r^{\star,\bR}_{\bra{\varrho,\geq\ell}}
    }2\cdot
\end{align}
Define
\begin{equation}\label{eq:abstract-quasilinear-skew-projected-objects-Rd}
\begin{aligned}
    a_{\bra{-2,\geq\ell}}^{\bR;\omega}
    &\defeq
    \frac{
        \widetilde a_{\bra{-2,\geq\ell}}
        -
        \overline{\widetilde a_{\bra{-2,\geq\ell}}}
    }{2\ii},\quad \quad r^{\textnormal{sk},\bR}_{\bra{\varrho,\geq\ell}}
    \defeq
    \frac{
        \widetilde r^\bR_{\bra{\varrho,\geq\ell}}
        -
        \widetilde r^{\star,\bR}_{\bra{\varrho,\geq\ell}}
    }2\cdot
\end{aligned}
\end{equation}
Then
\(a_{\bra{-2,\geq\ell}}^{\bR;\omega}\) is real-valued, belongs to
\(\Gamma^{-2}_{\geq\ell}\), and, by
\cref{eq:abstract-quasilinear-preliminary-reality-Rd}, satisfies
\[
    a_{\bra{-2,\geq\ell}}^{\bR;\omega}\pare{\eta;x,-\xi}
    =
    -a_{\bra{-2,\geq\ell}}^{\bR;\omega}\pare{\eta;x,\xi}.
\]
Moreover, since \(\widetilde r^\bR_{\bra{\varrho,\geq\ell}}\) is real-to-real, so is
\(\widetilde r^{\star,\bR}_{\bra{\varrho,\geq\ell}}\). Therefore
\(r^{\textnormal{sk},\bR}_{\bra{\varrho,\geq\ell}}\in
\cR^\varrho_{\geq\ell}\) is real-to-real and skew-adjoint. Substituting
\cref{eq:abstract-quasilinear-skew-projected-objects-Rd} into the preceding
identity gives
\cref{eq:abstract-quasilinear-skew-completion-Rd}.
\end{proofsubpart}
\end{proofpart}

\end{proof}
\subsubsection{Transported contribution}

\begin{lemma}[Transported contribution]
\label{lem:abstract-transported-contribution-Rd}
For every prescribed \(\varrho\geq0\), the transported branch in
\cref{eq:abstract-common-transported-branch} satisfies
\begin{equation}\label{eq:abstract-transported-final}
    \mathpzc t_{h,k}^{\mathrm{tr}}\pare\eta\bra\omega
    =
    \OpBW{
        \omega\pare x
        a_1^{h,k;\eta}\pare{\nabla\eta\pare x,\xi}
        +
        a_{\bra{0,\geq\ell}}\bra\omega\pare{\eta;x,\xi}
    }\eta
    +
    r^\bR_{\bra{\varrho,\geq\ell}}\bra\omega\pare\eta\eta,
\end{equation}
where
\begin{itemize}
    \item
    \(a_1^{h,k;\eta}\in\Gamma^1_{\geq\ell-1}\) is the real, even symbol
    defined in \cref{eq:abstract-explicit-principal-symbol-A1-Rd};

    \item
    \(a_{\bra{0,\geq\ell}}\bra\cdot
    \in\cL\pare{H^{s_0};\Gamma^0}\) and
    \(r^\bR_{\bra{\varrho,\geq\ell}}\bra\cdot
    \in\cL\pare{H^{s_0};\cR^\varrho}\) are generic linear
    terms as per \cref{notation:paralinearization}. If \(\ell\geq2\),
    their ranges improve to \(\Gamma^0_{\geq\ell-1}\) and
    \(\cR^\varrho_{\geq\ell-1}\), respectively.
\end{itemize}
Moreover,
\begin{equation}\label{eq:abstract-transported-lower-symbol-reality}
    a_{\bra{0,\geq\ell}}\bra\omega\pare{\eta;x,-\xi}
    =
    \overline{
        a_{\bra{0,\geq\ell}}\bra\omega\pare{\eta;x,\xi}
    },
\end{equation}
and the remainder is real-to-real.
\end{lemma}

\begin{proof}
\begin{proofpart}[Bony splitting]
Applying \cref{lem paralin} to
\(\widetilde k\pare{\Delta_z\eta}\) gives
\begin{equation}\label{eq:abstract-transported-Bony-paralinearization}
    \widetilde k\pare{\Delta_z\eta}
    =
    \OpBW{k'\pare{\Delta_z\eta}}\Delta_z\eta
    +
    \mathpzc q_k\pare{\Delta_z\eta}
    \bra{\Delta_z\eta},
\end{equation}
where \(\mathpzc q_k\pare u\defeq R_{\widetilde k}\pare u\).
Substituting this identity into
\cref{eq:abstract-common-transported-branch} yields
\begin{equation}\label{eq:abstract-transported-Bony-splitting}
    \mathpzc t_{h,k}^{\mathrm{tr}}\pare\eta\bra\omega
    =
    \mathpzc t_{h,k}^{\mathrm{tr},\mathrm I}\pare\eta\bra\omega
    +
    \mathpzc t_{h,k}^{\mathrm{tr},\mathrm {II}}\pare\eta\bra\omega,
\end{equation}
where
\begin{subequations}\label{eq:abstract-transported-Bony-branches}
\begin{align}
    \mathpzc t_{h,k}^{\mathrm{tr},\mathrm I}
    \pare\eta\bra\omega
    \defeq{}&
    \pv\int_{\bR^{\mathtt d}}
    \chi\pare z
    \frac{h\pare{\hat z}}{\av z^{\mathtt d}}
    \OpBW{\omega_z}
    \OpBW{k'\pare{\Delta_z\eta}}
    \Delta_z\eta\dd z,
    \label{eq:abstract-transported-principal-branch}
    \\
    \mathpzc t_{h,k}^{\mathrm{tr},\mathrm {II}}
    \pare\eta\bra\omega
    \defeq{}&
    \pv\int_{\bR^{\mathtt d}}
    \chi\pare z
    \frac{h\pare{\hat z}}{\av z^{\mathtt d}}
    \OpBW{\omega_z}
    \mathpzc q_k\pare{\Delta_z\eta}
    \bra{\Delta_z\eta}\dd z.
    \label{eq:abstract-transported-Bony-remainder-operator}
\end{align}
\end{subequations}
\end{proofpart}

\begin{proofpart}[Analysis of \cref{eq:abstract-transported-principal-branch}]
Let $\mathsf{v}\bra{\eta}$ and $\mathsf{h}\bra{\eta}$ be defined by
\cref{eq:abstract-quasilinear-local-building-blocks}. Since
$\av{\mathsf{v}\bra{\eta}}\leq\norm{\nabla\eta}_{L^\infty}
\lesssim\norm{\eta }_{s_0}<r_0$, the derivatives of $k$ admit the
absolutely convergent multinomial expansions
\begin{subequations}\label{eq:abstract-transported-multinomial-expansions}
\begin{align}
    k'\pare{\mathsf{v}\bra{\eta}}
    ={}&
    \sum_{\substack{\alpha\in\bN^{\mathtt d}\\
                     \av\alpha+1\geq\ell\\
                     k^{\pare{\av\alpha+1}}\pare0\neq0}}
    \frac{k^{\pare{\av\alpha+1}}\pare0}{\alpha!}
    \pare{\nabla\eta\pare x}^\alpha\hat z^\alpha,
    \label{eq:abstract-transported-kprime-multinomial-expansion}
    \\
    k''\pare{\mathsf{v}\bra{\eta}}
    ={}&
    \sum_{\substack{\alpha\in\bN^{\mathtt d}\\
                     \av\alpha+2\geq\ell\\
                     k^{\pare{\av\alpha+2}}\pare0\neq0}}
    \frac{k^{\pare{\av\alpha+2}}\pare0}{\alpha!}
    \pare{\nabla\eta\pare x}^\alpha\hat z^\alpha \ ,
    \label{eq:abstract-transported-ksecond-multinomial-expansion}
\end{align}
\end{subequations}
where $\alpha!\defeq\prod_{j=1}^{\mathtt d}\alpha_j!$ and the constraints on
the summation indices follow from
\cref{eq:abstract-admissible-K-vanishing-order}. Set
\begin{equation}\label{eq:abstract-transported-composition-defect}
    \mathpzc c_z\bra\omega\pare\eta
    \defeq
    \OpBW{\omega_z}\circ\OpBW{k'\pare{\Delta_z\eta}}
    -
    \OpBW{\omega_zk'\pare{\Delta_z\eta}} \ .
\end{equation}
Taylor's formula,
\cref{eq:abstract-quasilinear-difference-Taylor-expansion}, and
\cref{lem:kernel,lem:hsnorm,lem:composition_rules_z} give
\begin{equation}\label{eq:abstract-transported-local-expansions}
\begin{aligned}
    k'\pare{\Delta_z\eta}
    &=
    k'\pare{\mathsf{v}\bra{\eta}}
    -\frac{\av z}2k''\pare{\mathsf{v}\bra{\eta}}\mathsf{h}\bra{\eta}
    +k^{\loc,\mathrm{tr},\bR}_{\bra{2,\geq\max\set{1,\ell-1}}},
    &
    k^{\loc,\mathrm{tr},\bR}_{\bra{2,\geq\max\set{1,\ell-1}}}
    &\in
    \KF^2_{\loc,\geq\max\set{1,\ell-1}},
    \\
    \omega_z
    &=
    \omega\pare x
    -\av z\,\nabla\omega\pare x\cdot\hat z
    +e_{\bra2}\bra\omega,
    &
    e_{\bra2}
    &\in
    \cL\pare{H^{s_0};\KF^2_{\loc}} \ ,
\end{aligned}
\end{equation}
whence
\begin{equation}\label{eq:abstract-transported-product-expansion}
    \omega_zk'\pare{\Delta_z\eta}
    =
    \omega\pare xk'\pare{\mathsf{v}\bra{\eta}}
    -
    \av z\bra{
        \frac12\omega\pare xk''\pare{\mathsf{v}\bra{\eta}}
        \mathsf{h}\bra{\eta}
        +
        \pare{\nabla\omega\pare x\cdot\hat z}
        k'\pare{\mathsf{v}\bra{\eta}}
    }
    +
    e_{\bra2}^{\mathrm{tr}}\bra\omega \ ,
\end{equation}
where $e_{\bra2}^{\mathrm{tr}}\bra\cdot\in
\cL\pare{H^{s_0};\KF^2_{\loc,\geq\ell-1}}$ is defined by
\cref{eq:abstract-transported-product-expansion} itself. Accordingly, since
$\Delta_z\eta=\delta_z\eta/\av z$,
\begin{equation}\label{eq:abstract-transported-principal-splitting}
    \mathpzc t_{h,k}^{\mathrm{tr},\mathrm I}\pare\eta\bra\omega
    =
    \mathpzc t_{h,k}^{\mathrm{tr},1,\loc}\pare\eta\bra\omega
    +
    \mathpzc t_{h,k}^{\mathrm{tr},0,\loc}\pare\eta\bra\omega
    +
    \mathpzc t_{h,k}^{\mathrm{tr},\leq-1,\loc}\pare\eta\bra\omega
    +
    \mathpzc t_{h,k}^{\mathrm{tr},\mathrm{cp}}\pare\eta\bra\omega \ .
\end{equation}
Here
\begin{subequations}\label{eq:abstract-transported-principal-models}
\begin{align}
    \mathpzc t_{h,k}^{\mathrm{tr},1,\loc}\pare\eta\bra\omega
    \defeq{}&
    \frac12
    \int_{\bR^{\mathtt d}}
    \chi\pare z
    \frac{h\pare{\hat z}}{\av z^{\mathtt d+1}}
    \OpBW{\omega\pare xk'\pare{\mathsf{v}\bra{\eta}}}
    \delta_z^{\pare2}\eta\dd z,
    \label{eq:abstract-transported-local-principal-model}
    \\
    \mathpzc t_{h,k}^{\mathrm{tr},0,\loc}\pare\eta\bra\omega
    \defeq{}&
    \int_{\bR^{\mathtt d}}
    \OpBW{k^{\pare0}\bra\omega}
    \delta_z\eta\dd z,
    \label{eq:abstract-transported-order-zero-model}
    \\
    \mathpzc t_{h,k}^{\mathrm{tr},\leq-1,\loc}\pare\eta\bra\omega
    \defeq{}&
    \int_{\bR^{\mathtt d}}
    \OpBW{k^{\pare{\leq-1}}\bra\omega}
    \delta_z\eta\dd z,
    \label{eq:abstract-transported-lower-order-model}
    \\
    \mathpzc t_{h,k}^{\mathrm{tr},\mathrm{cp}}\pare\eta\bra\omega
    \defeq{}&
    \pv\int_{\bR^{\mathtt d}}
    \chi\pare z
    \frac{h\pare{\hat z}}{\av z^{\mathtt d}}
    \mathpzc c_z\bra\omega\pare\eta\bra{\Delta_z\eta}\dd z,
    \label{eq:abstract-transported-composition-model}
\end{align}
\end{subequations}
and
\begin{subequations}\label{eq:abstract-transported-expansion-kernels}
\begin{align}
    k^{\pare0}\bra\omega\pare{\eta;x,z}
    \defeq{}&
    -\chi\pare z
    \frac{h\pare{\hat z}}{\av z^{\mathtt d}}
    \bra{
        \frac12\omega\pare xk''\pare{\mathsf{v}\bra{\eta}}
        \mathsf{h}\bra{\eta}
        +
        \pare{\nabla\omega\pare x\cdot\hat z}
        k'\pare{\mathsf{v}\bra{\eta}}
    },
    \label{eq:abstract-transported-endpoint-kernel}
    \\
    k^{\pare{\leq-1}}\bra\omega\pare{\eta;x,z}
    \defeq{}&
    \chi\pare z
    \frac{h\pare{\hat z}}{\av z^{\mathtt d+1}}
    e_{\bra2}^{\mathrm{tr}}\bra\omega
    \in
    \KF^{1-\mathtt d}_{\geq\ell-1} \ .
    \label{eq:abstract-transported-lower-order-kernel}
\end{align}
\end{subequations}
In \cref{eq:abstract-transported-local-principal-model} we used that
$h\pare\vartheta k'\pare{p\cdot\vartheta}$ is even by
\cref{eq:abstract-admissible-derived-parities}: the change of variables
$z\mapsto-z$ replaces $\delta_z\eta$ by
$\frac12\delta_z^{\pare2}\eta$ and turns the principal value into an
absolutely convergent integral.

We analyze the terms in
\cref{eq:abstract-transported-principal-splitting} in decreasing symbolic
order.

\begin{proofsubpart}[Analysis of \cref{eq:abstract-transported-local-principal-model}]
For \(\av{\mathsf{v}\bra{\eta}}\) sufficiently small we use
\cref{eq:abstract-transported-kprime-multinomial-expansion}.
For every index in this sum, set
\begin{equation}\label{eq:abstract-transported-monomial-data}
    q_\alpha\bra\omega\pare{\eta;x}
    \defeq
    \frac{k^{\pare{\av\alpha+1}}\pare0}{\alpha!}
    \omega\pare x
    \pare{\nabla\eta\pare x}^{\alpha},
    \qquad
    H_\alpha^{\pare1}\pare z
    \defeq
    \frac{h\pare{\hat z}\hat z^\alpha}
    {\av z^{\mathtt d+1}}\cdot
\end{equation}
By \cref{eq:abstract-admissible-K-parity},
\(k^{\pare{\av\alpha+1}}\pare0\neq0\) implies
\(\pare{-1}^{\av\alpha}=\varepsilon\); hence
\(H_\alpha^{\pare1}\) is even. Define
\begin{equation}\label{eq:abstract-transported-local-compensated-monomial}
    \mathpzc p_\alpha^{1,\loc}\bra{f}\pare x
    \defeq
    \frac12
    \int_{\bR^{\mathtt d}}
    \chi\pare z
    H_\alpha^{\pare1}\pare z
    \delta_z^{\pare2}f\pare x\dd z.
\end{equation}
Then
\begin{equation}\label{eq:abstract-transported-local-monomial-expansion}
    \mathpzc t_{h,k}^{\mathrm{tr},1,\loc}
    \pare\eta\bra\omega
    =
    \sum_{\substack{\alpha\in\bN^{\mathtt d}\\
                     \av\alpha+1\geq\ell\\
                     k^{\pare{\av\alpha+1}}\pare0\neq0}}
    \OpBW{q_\alpha\bra\omega}
    \bra{\mathpzc p_\alpha^{1,\loc}{\eta}}.
\end{equation}

Using the complete homogeneous multiplier
\(\dot{\mathfrak m}_{H_\alpha^{\pare1}}^{\pare2}\) from
\cref{def:localized-compensated-homogeneous-convolution-models}, define
\begin{equation}\label{eq:abstract-transported-global-model}
    \mathpzc g_{h,k}^{\mathrm{tr},1}
    \pare\eta\bra\omega
    \defeq
    \sum_{\substack{\alpha\in\bN^{\mathtt d}\\
                     \av\alpha+1\geq\ell\\
                     k^{\pare{\av\alpha+1}}\pare0\neq0}}
    \OpBW{q_\alpha\bra\omega}
    \bra{
        \dot{\mathfrak m}_{H_\alpha^{\pare1}}^{\pare2}
        \pare D\eta
    }.
\end{equation}
By
\cref{eq:localized-compensated-homogeneous-comparison},
for every admissible \(\alpha\),
\[
    \mathpzc p_\alpha^{1,\loc}
    =
    \dot{\mathfrak m}_{H_\alpha^{\pare1}}^{\pare2}\pare D
    +
    c_{H_\alpha^{\pare1}}\Id
    +
    \mathfrak s_\alpha\pare D,
    \qquad
    \mathfrak s_\alpha
    \in
    \bigcap_{M\in\bN}\tilde\Gamma^{-M}.
\]
Therefore, applying
\cref{rem:localized-homogeneous-model-coefficients} term by term to
\cref{eq:abstract-transported-local-monomial-expansion} gives
\begin{equation}\label{eq:abstract-transported-cutoff-reduction}
    \mathpzc t_{h,k}^{\mathrm{tr},1,\loc}
    \pare\eta\bra\omega
    =
    \mathpzc g_{h,k}^{\mathrm{tr},1}
    \pare\eta\bra\omega
    +
    \OpBW{a_{\bra{0,\geq\ell}}\bra\omega}\eta
    +
    r^\bR_{\bra{\varrho,\geq\ell}}\bra\omega\pare\eta\eta.
\end{equation}
Here the constant corrections
\(c_{H_\alpha^{\pare1}}\) are included in
\(a_{\bra{0,\geq\ell}}\bra\omega\), whereas the multipliers
\(\mathfrak s_\alpha\) are included in
\(r^\bR_{\bra{\varrho,\geq\ell}}\bra\omega\pare\eta\).
Summability follows from
\cref{eq:abstract-angular-multiplier-seminorm-bound} and the absolute
convergence of
\cref{eq:abstract-transported-kprime-multinomial-expansion} for
\(\av{\mathsf{v}\bra{\eta}}<r_0\), after shrinking \(r_0\) if necessary.
By
\cref{eq:compensated-homogeneous-multiplier-polar},
\[
    \dot{\mathfrak m}_{H_\alpha^{\pare1}}^{\pare2}\pare\xi
    =
    \dot{\mathfrak m}_\alpha^{\pare1}\pare\xi
    \defeq
    \frac\pi2
    \int_{\bS^{\mathtt d-1}}
    \av{\vartheta\cdot\xi}
    h\pare\vartheta\vartheta^\alpha
    \dd\sigma\pare\vartheta.
\]
Applying \cref{prop:composition_paradiff} to
\cref{eq:abstract-transported-global-model} gives
\begin{equation}\label{eq:abstract-transported-principal-sum}
\begin{aligned}
    \sum_{\substack{\alpha\in\bN^{\mathtt d}\\
                     \av\alpha+1\geq\ell\\
                     k^{\pare{\av\alpha+1}}\pare0\neq0}}
    q_\alpha\bra\omega
    \dot{\mathfrak m}_\alpha^{\pare1}\pare\xi
    &={}
    \omega\pare x
    \frac\pi2
    \int_{\bS^{\mathtt d-1}}
    \av{\vartheta\cdot\xi}
    h\pare\vartheta
    k'\pare{\nabla\eta\pare x\cdot\vartheta}
    \dd\sigma\pare\vartheta
   =
    \omega\pare x
    a_1^{h,k;\eta}
    \pare{\nabla\eta\pare x,\xi}.
\end{aligned}
\end{equation}
Since the multiplier has order one, all non-principal terms in the
composition have order at most zero. Therefore
\begin{equation}\label{eq:abstract-transported-composition-final}
    \mathpzc g_{h,k}^{\mathrm{tr},1}
    \pare\eta\bra\omega
    =
    \OpBW{
        \omega\pare x
        a_1^{h,k;\eta}\pare{\nabla\eta\pare x,\xi}
        +
        a_{\bra{0,\geq\ell}}\bra\omega
    }\eta
    +
    r^\bR_{\bra{\varrho,\geq\ell}}\bra\omega\pare\eta\eta.
\end{equation}
\end{proofsubpart}

\begin{proofsubpart}[Analysis of \cref{eq:abstract-transported-order-zero-model}]
By \cref{eq:abstract-admissible-derived-parities},
\[
    h\pare\vartheta
    k''\pare{\nabla\eta\cdot\vartheta}
    \pare{\vartheta^\intercal{\bm H}\eta\vartheta}
    \quad\textnormal{is odd in }\vartheta,
\]
while
\[
    h\pare\vartheta
    k'\pare{\nabla\eta\cdot\vartheta}
    \pare{\nabla\omega\cdot\vartheta}
    \quad\textnormal{is odd in }\vartheta.
\]
Hence the angular kernel in
\cref{eq:abstract-transported-endpoint-kernel} is odd, and therefore
\begin{equation}\label{eq:abstract-transported-endpoint-cancellation}
    \int_{\bS^{\mathtt d-1}}
    k^{\pare0}\bra\omega\pare{\eta;x,r\vartheta}
    \dd\sigma\pare\vartheta
    =0,
    \qquad
    r>0 \ .
\end{equation}
For every $\alpha\in\bN^{\mathtt d}$ and every
$1\leq j,l\leq\mathtt d$ set
\begin{equation}
\label{eq:abstract-transported-order-zero-monomial-data}
\begin{aligned}
    c_{\alpha,j,l}\bra\omega\pare{\eta;x}
    \defeq{}&
    -\frac{k^{\pare{\av\alpha+2}}\pare0}{2\,\alpha!}
    \omega\pare x
    \pare{\nabla\eta\pare x}^\alpha
    \partial_{jl}\eta\pare x,
    &
    h_{\alpha,j,l}^{\pare0}\pare z
    \defeq{}&
    \frac{h\pare{\hat z}\hat z^\alpha\hat z_j\hat z_l}
    {\av z^{\mathtt d}},
    \\
    c_{\alpha,j}\bra\omega\pare{\eta;x}
    \defeq{}&
    -\frac{k^{\pare{\av\alpha+1}}\pare0}{\alpha!}
    \partial_j\omega\pare x
    \pare{\nabla\eta\pare x}^\alpha,
    &
    h_{\alpha,j}^{\pare0}\pare z
    \defeq{}&
    \frac{h\pare{\hat z}\hat z^\alpha\hat z_j}
    {\av z^{\mathtt d}}\cdot
\end{aligned}
\end{equation}

By
\cref{eq:abstract-transported-ksecond-multinomial-expansion,%
eq:abstract-transported-kprime-multinomial-expansion},
the first and second coefficient families in
\cref{eq:abstract-transported-order-zero-monomial-data}
occur only for indices satisfying
\[
    \av\alpha+2\geq\ell,
    \qquad
    \av\alpha+1\geq\ell,
\]
respectively. Since
\(c_{\alpha,j,l}\bra\omega\pare{\eta;x}\) has degree
\(\av\alpha+1\) in \(\eta\), whereas
\(c_{\alpha,j}\bra\omega\pare{\eta;x}\) has degree
\(\av\alpha\) in \(\eta\), both coefficient families vanish to order at
least \(\ell-1\) at \(\eta=0\).

Differentiating
\cref{eq:abstract-admissible-K-parity} at the origin gives
\(k^{\pare m}\pare0\neq0\Longrightarrow\pare{-1}^{m}=-\varepsilon\).
Hence
\(k^{\pare{\av\alpha+2}}\pare0\neq0\Longrightarrow
\pare{-1}^{\av\alpha}=-\varepsilon\).
Using also \cref{eq:abstract-admissible-H-parity}, we obtain
\[
    h\pare{-\vartheta}\pare{-\vartheta}^{\alpha}
    \pare{-\vartheta}_j\pare{-\vartheta}_l
    =
    \varepsilon\pare{-1}^{\av\alpha}
    h\pare\vartheta\vartheta^\alpha\vartheta_j\vartheta_l
    =
    -h\pare\vartheta\vartheta^\alpha\vartheta_j\vartheta_l.
\]
Similarly,
\(k^{\pare{\av\alpha+1}}\pare0\neq0\Longrightarrow
\pare{-1}^{\av\alpha}=\varepsilon\), and therefore
\[
    h\pare{-\vartheta}\pare{-\vartheta}^{\alpha}\pare{-\vartheta}_j
    =
    -\varepsilon\pare{-1}^{\av\alpha}
    h\pare\vartheta\vartheta^\alpha\vartheta_j
    =
    -h\pare\vartheta\vartheta^\alpha\vartheta_j.
\]
Thus the angular amplitudes of
\(h_{\alpha,j,l}^{\pare0}\) and \(h_{\alpha,j}^{\pare0}\) are odd and
therefore satisfy the hypothesis of
\cref{item:abstract-homogeneous-angular-degree-zero}.

Let $t_{h_{\alpha,j,l}^{\pare0}}^{\loc}$ and
$t_{h_{\alpha,j}^{\pare0}}^{\loc}$ be the localized homogeneous models
associated with these amplitudes as per
\cref{def:localized-homogeneous-convolution-models}. Splitting
$\delta_z\eta=\eta-\eta_z$ in
\cref{eq:abstract-transported-order-zero-model}, the unshifted contribution
vanishes by \cref{eq:abstract-transported-endpoint-cancellation}, and
\cref{eq:abstract-transported-endpoint-kernel,%
eq:abstract-transported-order-zero-monomial-data}
give
\begin{equation}\label{eq:abstract-transported-order-zero-expansion}
\begin{aligned}
    \mathpzc t_{h,k}^{\mathrm{tr},0,\loc}\pare\eta\bra\omega
    ={}&
    -\sum_{\substack{\alpha\in\bN^{\mathtt d}\\
                     \av\alpha+2\geq\ell\\
                     k^{\pare{\av\alpha+2}}\pare0\neq0}}
    \ \sum_{j,l=1}^{\mathtt d}
    \OpBW{c_{\alpha,j,l}\bra\omega}
    \bra{t_{h_{\alpha,j,l}^{\pare0}}^{\loc}\eta}
    \\
    &-\sum_{\substack{\alpha\in\bN^{\mathtt d}\\
                     \av\alpha+1\geq\ell\\
                     k^{\pare{\av\alpha+1}}\pare0\neq0}}
    \ \sum_{j=1}^{\mathtt d}
    \OpBW{c_{\alpha,j}\bra\omega}
    \bra{t_{h_{\alpha,j}^{\pare0}}^{\loc}\eta}.
\end{aligned}
\end{equation}
Applying
\cref{lem:localized-homogeneous-model,%
rem:localized-homogeneous-model-coefficients}
term by term removes the cutoff,
\cref{item:abstract-homogeneous-angular-degree-zero} identifies each
homogeneous multiplier as $\ii$ times a real multiplier homogeneous of
degree zero, and \cref{prop:composition_paradiff} recasts the resulting
compositions as paradifferential operators of order zero. Summability
follows from \cref{eq:abstract-angular-multiplier-seminorm-bound} and the
absolute convergence of
\cref{eq:abstract-transported-multinomial-expansions}, after shrinking
$r_0$ if necessary. Therefore
\begin{equation}\label{eq:abstract-transported-order-zero-reduction}
    \mathpzc t_{h,k}^{\mathrm{tr},0,\loc}\pare\eta\bra\omega
    =
    \OpBW{a_{\bra{0,\geq\ell}}\bra\omega}\eta
    +
    r^\bR_{\bra{\varrho,\geq\ell}}\bra\omega\pare\eta\eta \ .
\end{equation}
\end{proofsubpart}

\begin{proofsubpart}[Analysis of \cref{eq:abstract-transported-lower-order-model}]
By \cref{eq:abstract-transported-lower-order-kernel}, the kernel
$k^{\pare{\leq-1}}\bra\cdot$ is linear in $\omega$, has local order
$1-\mathtt d>-\mathtt d$ and background vanishing order at least $\ell-1$.
Splitting $\delta_z\eta=\eta-\eta_z$ in
\cref{eq:abstract-transported-lower-order-model},
\cref{lem:integration_KF} applies to the unshifted contribution and
\cref{prop:int-opbw} to the translated one, the latter producing a symbol of
order $-1$, which we absorb in the generic order-zero class. Hence
\begin{equation}\label{eq:abstract-transported-lower-order-reduction}
    \mathpzc t_{h,k}^{\mathrm{tr},\leq-1,\loc}\pare\eta\bra\omega
    =
    \OpBW{a_{\bra{0,\geq\ell}}\bra\omega}\eta
    +
    r^\bR_{\bra{\varrho,\geq\ell}}\bra\omega\pare\eta\eta .
\end{equation}
\end{proofsubpart}

\begin{proofsubpart}[Analysis of \cref{eq:abstract-transported-composition-model}]
Set
\begin{equation}\label{eq:abstract-transported-frozen-composition-defect}
    \mathpzc c_0\bra\omega\pare\eta
    \defeq
    \OpBW{\omega}\circ\OpBW{k'\pare{\mathsf{v}\bra{\eta}}}
    -
    \OpBW{\omega k'\pare{\mathsf{v}\bra{\eta}}}.
\end{equation}
By \cref{eq:abstract-admissible-derived-parities} the map
$z\mapsto h\pare{\hat z}\mathpzc c_0\bra\omega\pare\eta
\bra{\mathsf{v}\bra{\eta}}$ is odd, hence
\begin{equation}\label{eq:abstract-transported-frozen-composition-cancellation}
    \pv\int_{\bR^{\mathtt d}}
    \chi\pare z
    \frac{h\pare{\hat z}}{\av z^{\mathtt d}}
    \mathpzc c_0\bra\omega\pare\eta
    \bra{\mathsf{v}\bra{\eta}}\dd z
    =0,
\end{equation}
and \cref{eq:abstract-transported-composition-model} may be rewritten as
\begin{equation}\label{eq:abstract-transported-composition-frozen-subtraction}
    \mathpzc t_{h,k}^{\mathrm{tr},\mathrm{cp}}\pare\eta\bra\omega
    =
    \pv\int_{\bR^{\mathtt d}}
    \chi\pare z
    \frac{h\pare{\hat z}}{\av z^{\mathtt d}}
    \set{
        \mathpzc c_z\bra\omega\pare\eta\bra{\Delta_z\eta}
        -
        \mathpzc c_0\bra\omega\pare\eta\bra{\mathsf{v}\bra{\eta}}
    }\dd z.
\end{equation}
Since $\pare{a,b}\mapsto\OpBW a\circ\OpBW b-\OpBW{ab}$ is bilinear,
\cref{eq:abstract-transported-composition-defect,%
eq:abstract-transported-frozen-composition-defect}
give that the integrand in braces equals
\begin{equation}\label{eq:abstract-transported-composition-difference}
\begin{aligned}
    &\pare{
        \OpBW{\omega_z-\omega}\circ\OpBW{k'\pare{\Delta_z\eta}}
        -
        \OpBW{\pare{\omega_z-\omega}k'\pare{\Delta_z\eta}}
    }\Delta_z\eta
    \\
    &\quad+
    \pare{
        \OpBW{\omega}\circ
        \OpBW{k'\pare{\Delta_z\eta}-k'\pare{\mathsf{v}\bra{\eta}}}
        -
        \OpBW{\omega\pare{
            k'\pare{\Delta_z\eta}-k'\pare{\mathsf{v}\bra{\eta}}
        }}
    }\Delta_z\eta
    \\
    &\quad+
    \mathpzc c_0\bra\omega\pare\eta
    \bra{\Delta_z\eta-\mathsf{v}\bra{\eta}}.
\end{aligned}
\end{equation}
Moreover,
\begin{equation}\label{eq:abstract-transported-frozen-difference-formula}
\begin{aligned}
    \omega_z-\omega
    &=
    -\av z
    \int_0^1
    \nabla\omega\pare{x-tz}\cdot\hat z\dd t,
    \\
    \Delta_z\eta-\mathsf{v}\bra{\eta}
    &=
    -\av z
    \int_0^1\pare{1-t}
    \hat z^\intercal{\bm H}\eta\pare{x-tz}\hat z\dd t,
    \\
    k'\pare{\Delta_z\eta}-k'\pare{\mathsf{v}\bra{\eta}}
    &=
    \pare{\Delta_z\eta-\mathsf{v}\bra{\eta}}
    \int_0^1
    k''\pare{
        \mathsf{v}\bra{\eta}
        +t\pare{\Delta_z\eta-\mathsf{v}\bra{\eta}}
    }\dd t,
\end{aligned}
\end{equation}
so that every summand of
\cref{eq:abstract-transported-composition-difference} carries one factor
$\av z$. Writing
\[
    \Delta_z\eta
    =
    \int_0^1\nabla\eta\pare{x-tz}\cdot\hat z\dd t
\]
to reserve $\eta$ as the argument, \cref{lem:composition_rules_z} turns
\cref{eq:abstract-transported-composition-difference} into a
kernel-smoothing family of local order $1-\mathtt d>-\mathtt d$, linear in
$\omega$ and of background vanishing order at least $\ell-1$. Applying
\cref{lem:integration_KR} to
\cref{eq:abstract-transported-composition-frozen-subtraction} we conclude
that
\begin{equation}\label{eq:abstract-transported-composition-reduction}
    \mathpzc t_{h,k}^{\mathrm{tr},\mathrm{cp}}\pare\eta\bra\omega
    =
    r^\bR_{\bra{\varrho,\geq\ell}}\bra\omega\pare\eta\eta.
\end{equation}
\end{proofsubpart}

Combining
\cref{eq:abstract-transported-principal-splitting,%
eq:abstract-transported-composition-reduction,%
eq:abstract-transported-order-zero-reduction,%
eq:abstract-transported-lower-order-reduction}
gives
\begin{equation}\label{eq:abstract-transported-endpoint-reduction}
    \mathpzc t_{h,k}^{\mathrm{tr},\mathrm I}\pare\eta\bra\omega
    =
    \mathpzc t_{h,k}^{\mathrm{tr},1,\loc}\pare\eta\bra\omega
    +
    \OpBW{a_{\bra{0,\geq\ell}}\bra\omega}\eta
    +
    r^\bR_{\bra{\varrho,\geq\ell}}\bra\omega\pare\eta\eta.
\end{equation}
\end{proofpart}

\begin{proofpart}[Analysis of \cref{eq:abstract-transported-Bony-remainder-operator}]
Applying \cref{lem paralin} to \(F=\widetilde k\), with smoothing gain
\(\varrho+2\), gives
\begin{equation}\label{eq:abstract-transported-Bony-operator-class}
    \mathpzc q_k
    =
    R_{\widetilde k}
    \in
    \cR^{\varrho+2}_{\geq1}.
\end{equation}
Moreover,
\begin{equation}\label{eq:abstract-transported-Bony-contracted-remainder}
    \mathpzc q_k\pare u\bra u
    =
    \widetilde k\pare u
    -
    \OpBW{k'\pare u}u.
\end{equation}
By
\cref{eq:abstract-admissible-K-vanishing,rem:Gamma-vanishing-rules},
both terms on the right-hand side of
\cref{eq:abstract-transported-Bony-contracted-remainder} vanish to order
at least \(\ell\) in \(u\). Hence
\(u\mapsto\mathpzc q_k\pare u\bra u\) vanishes to order at least \(\ell\),
and its differential vanishes to order at least \(\ell-1\).

By
\cref{eq:abstract-admissible-K-parity},
\[
    \widetilde k\pare{-u}
    =
    -\varepsilon\widetilde k\pare u,
    \qquad
    k'\pare{-u}
    =
    \varepsilon k'\pare u.
\]
Using
\cref{eq:abstract-transported-Bony-contracted-remainder}, we therefore get
\[
    \mathpzc q_k\pare{-u}\bra{-u}
    =
    -\varepsilon
    \mathpzc q_k\pare u\bra u.
\]

Hence the integrand obtained from
\cref{eq:abstract-transported-Bony-remainder-operator} by replacing
\((\omega_z,\Delta_z\eta)\) with
\((\omega\pare x,\mathsf{v}\bra{\eta})\) is odd.
Subtracting it gives
\begin{equation}\label{eq:abstract-transported-Bony-remainder-difference}
\begin{aligned}
    &\OpBW{\omega_z}
    \mathpzc q_k\pare{\Delta_z\eta}\bra{\Delta_z\eta}
    -
    \OpBW{\omega\pare x}
    \mathpzc q_k\pare{\mathsf{v}\bra{\eta}}
    \bra{\mathsf{v}\bra{\eta}}
    \\
    &\quad={}
    \OpBW{\omega_z-\omega\pare x}
    \mathpzc q_k\pare{\Delta_z\eta}\bra{\Delta_z\eta}
    \\
    &\qquad+
    \OpBW{\omega\pare x}
    \mathpzc q_k\pare{\Delta_z\eta}
    \bra{\Delta_z\eta-\mathsf{v}\bra{\eta}}
    \\
    &\qquad+
    \OpBW{\omega\pare x}
    \left(
        \mathpzc q_k\pare{\Delta_z\eta}
        -
        \mathpzc q_k\pare{\mathsf{v}\bra{\eta}}
    \right)\bra{\mathsf{v}\bra{\eta}}.
\end{aligned}
\end{equation}
By \cref{eq:abstract-transported-frozen-difference-formula},
\begin{equation}\label{eq:abstract-transported-Bony-z-gains}
\begin{aligned}
    \omega_z-\omega\pare x
    &=
    -\av z
    \int_0^1
    \nabla\omega\pare{x-tz}\cdot\hat z\,\dd t,
    \\
    \Delta_z\eta-\mathsf{v}\bra{\eta}
    &=
    -\av z
    \int_0^1
    \pare{1-t}
    \hat z^\intercal{\bm H}\eta\pare{x-tz}\hat z\,\dd t,
    \\
    \mathpzc q_k\pare{\Delta_z\eta}
    -
    \mathpzc q_k\pare{\mathsf{v}\bra{\eta}}
    &=
    -\av z
    \int_0^1\int_0^1
    \pare{1-s}
    \dd \mathpzc q_k
    \pare{
        \mathsf{v}\bra{\eta}
        +t\pare{\Delta_z\eta-\mathsf{v}\bra{\eta}}
    }
    \bra{
        \hat z^\intercal{\bm H}\eta\pare{x-sz}\hat z
    }
    \,\dd s\,\dd t .
\end{aligned}
\end{equation}
Here \(\dd \mathpzc q_k\) is well defined since the remainder map
\(\mathpzc q_k=R_{\widetilde k}\) inherits the analyticity of \(k\).
Thus each term in
\cref{eq:abstract-transported-Bony-remainder-difference}
gains one power of \(\av z\).
We may therefore apply
\cref{rem:smoothBONYremainder,lem:composition_rules_z,%
lem:integration_KR}
to \cref{eq:abstract-transported-Bony-remainder-operator} and obtain
\begin{equation}\label{eq:abstract-transported-Bony-remainder-final}
    \mathpzc t_{h,k}^{\mathrm{tr},\mathrm {II}}
    \pare\eta\bra\omega
    =
    r^\bR_{\bra{\varrho,\geq\ell}}\bra\omega\pare\eta\eta.
\end{equation}
Combining
\cref{eq:abstract-transported-Bony-splitting,%
eq:abstract-transported-endpoint-reduction,%
eq:abstract-transported-Bony-remainder-final}
gives
\begin{equation}\label{eq:abstract-transported-almost-paralinearized}
    \mathpzc t_{h,k}^{\mathrm{tr}}
    \pare\eta\bra\omega
    =
    \mathpzc t_{h,k}^{\mathrm{tr},1,\loc}
    \pare\eta\bra\omega
    +
    \OpBW{a_{\bra{0,\geq\ell}}\bra\omega}\eta
    +
    r^\bR_{\bra{\varrho,\geq\ell}}\bra\omega\pare\eta\eta.
\end{equation}
\end{proofpart}

\begin{proofpart}[Proof of \cref{eq:abstract-transported-final}]
Combining
\cref{eq:abstract-transported-almost-paralinearized,%
eq:abstract-transported-cutoff-reduction,%
eq:abstract-transported-composition-final}
proves \cref{eq:abstract-transported-final}, with
\[
    a_{\bra{0,\geq\ell}}\bra\cdot
    \in
    \cL\pare{H^{s_0};\Gamma^0},
    \qquad
    r^\bR_{\bra{\varrho,\geq\ell}}\bra\cdot
    \in
    \cL\pare{H^{s_0};\cR^\varrho}.
\]
If \(\ell\geq2\), the corresponding ranges improve to
\(\Gamma^0_{\geq\ell-1}\) and
\(\cR^\varrho_{\geq\ell-1}\), respectively.
The contributions to
\(a_{\bra{0,\geq\ell}}\bra\omega\)
satisfy
\cref{eq:abstract-transported-lower-symbol-reality}
term by term: the lower-order kernel contribution by
\cref{eq:int-opbw-reality}, the homogeneous endpoint contribution by
\cref{item:abstract-homogeneous-angular-degree-zero,%
lem:localized-homogeneous-model}, and the compensated principal
contribution by
\cref{lem:localized-compensated-homogeneous-model}.
The Bony--Weyl composition formula preserves the corresponding
real-to-real symmetry. The smoothing remainder is real-to-real by the
same results and \cref{lem:integration_KR}.
\end{proofpart}
\end{proof}

\subsubsection{Proof of \cref{thm:abstract-structure-preserving-paralinearization-Rd}}

For smooth inputs, combining
\cref{eq:abstract-common-Bony-decomposition,%
eq:abstract-quasilinear-skew-completion-Rd,%
eq:abstract-transported-final,%
eq:abstract-localized-Bony-remainder-smoothing}
gives
\begin{equation}\label{eq:abstract-critical-synthesis-before-relabeling-Rd}
\begin{aligned}
    \mathpzc t_{h,k}^{\pare{\mathtt d}}\pare\eta\bra\omega
    ={}&
    \OpBW{
        \omega\pare x
        a_1^{h,k;\eta}\pare{\nabla\eta,\xi}
        +
        a_{\bra{0,\geq\ell}}\bra\omega
    }\eta
    \\
    &+
    \OpBW{
        \ii\pare{
            a_0^{h,k;\omega}\pare{\nabla\eta,\xi}
            +
            a_{\bra{-2,\geq\ell}}^{\bR;\omega}
        }
    }\omega
    \\
    &+
    r^\bR_{\bra{\varrho,\geq\ell}}\bra\omega\pare\eta\,\eta
    +
    r^\bR_{\bra{\varrho,\geq\ell}}\pare\eta\,\omega.
\end{aligned}
\end{equation}
The principal symbols satisfy
\begin{align}
    -\partial_p
    a_0^{h,k;\omega}\pare{p,\xi}\bra\xi
    &=
    \frac\pi2
    \int_{\bS^{\mathtt d-1}}
    \sgn\pare{\vartheta\cdot\xi}
    {\vartheta}
    h\pare\vartheta
    k'\pare{p\cdot\vartheta}
    \dd\sigma\pare\vartheta
    \\
    &=
    \frac\pi2 \partial_\xi
    \int_{\bS^{\mathtt d-1}}
    \av{\vartheta\cdot\xi}
    h\pare\vartheta
    k'\pare{p\cdot\vartheta}
    \dd\sigma\pare\vartheta
    \\
    &=\partial_\xi
    a_1^{h,k;\eta}\pare{p,\xi},
\end{align}
which proves
\cref{eq:abstract-principal-symbol-differential-identity}.
The symbol and remainder classes, their vanishing orders, and their
real-to-real properties follow directly from
\cref{lem:abstract-quasilinear-contribution-Rd,%
lem:abstract-transported-contribution-Rd,%
lem:abstract-localized-Bony-remainder-Rd,%
lem:abstract-common-local-tail-Bony-reduction-Rd}.
In particular,
\cref{lem:abstract-quasilinear-contribution-Rd} gives the purely
imaginary and odd complete quasilinear symbol in
\cref{eq:abstract-quasilinear-skew-completion-Rd}.
Absorbing the tail contribution into the corresponding generic
order-zero and smoothing terms, and relabeling the lower-order terms
according to \cref{notation:paralinearization},
\cref{eq:abstract-critical-synthesis-before-relabeling-Rd} becomes
\cref{eq:abstract-structure-preserving-paralinearization-Rd}.
\hfill\qed

\subsection{Paralinearization of singular integral operators}
\label{sec:paralinearization_SIO}

In this subsection, we apply \cref{thm:abstract-structure-preserving-paralinearization-Rd}
to the two singular integral operators $\mathpzc{B}_0$ and $\mathpzc{q}$ introduced in
\cref{eq:B0,eq:Q}.  For the convenience of the reader, we recall here their definitions:
\begin{subequations}
\begin{align}
    \mathpzc B_0\pare\eta\bra\omega\pare x
    &\defeq
    \frac1{4\pi}\pv\int_{\bR^2}
    \frac z{
        \pare{
            \av z^2+\pare{\delta_z\eta\pare x}^2
        }^{\frac32}
    }
    \omega\pare{x-z}\dd z,
    \label{eq:B0_recall}
    \\
    \mathpzc q\pare\eta\bra\omega\pare x
    &\defeq
    \frac1{4\pi}\pv\int_{\bR^2}
    \frac{\delta_z\eta\pare x}{
        \pare{
            \av z^2+\pare{\delta_z\eta\pare x}^2
        }^{\frac32}
    }
    \omega\pare{x-z}\dd z.
    \label{eq:Q_recall}
\end{align}
\end{subequations}
First of all, let us explain how we can rewrite these operators in a form adapted to the implementation of \cref{thm:abstract-structure-preserving-paralinearization-Rd}. This is given by the following lemma.
\begin{lemma}[Admissible kernels for $\mathpzc B_0$ and $\mathpzc q$]
\label{lem:B0-q-admissible-kernels}
For \(\vartheta\in\bS^1\) and \(X\in\bR\), set
\begin{equation}\label{eq:B0-q-admissible-kernels}
\begin{aligned}
    H^{\mathpzc B_0}\pare\vartheta
    &\defeq
    \frac{\vartheta}{4\pi},
    &
    k^{\mathpzc B_0}\pare X
    &\defeq
    \frac1{\pare{1+X^2}^{\frac32}}
    ,%=1+\mathtt K\pare X
    \\
    h^{\mathpzc q}\pare\vartheta
    &\defeq
    \frac1{4\pi},
    &
    k^{\mathpzc q}\pare X
    &\defeq
    \frac X{\pare{1+X^2}^{\frac32}}
    =Xk^{\mathpzc{B}_0}\pare{X}.%\mathtt K_1\pare X.
\end{aligned}
\end{equation}
Writing \(H^{\mathpzc B_0}=\vect{h^{\mathpzc B_0}_1}{h^{\mathpzc B_0}_2}\),
the pairs
\(
    \pare{h^{\mathpzc B_0}_{\mathsf j},k^{\mathpzc B_0}}
\), \(\mathsf j\in\set{1,2}\), are admissible of vanishing order
\(\ell=2\), whereas
\(
    \pare{h^{\mathpzc q},k^{\mathpzc q}}
\)
is admissible of vanishing order \(\ell=1\).  Moreover, recalling
\cref{eq:abstract-singular-graph-operator-Rd}, one has
\begin{equation}\label{eq:B0-q-abstract-operator-identification}
    \pare{\mathpzc B_0\pare\eta\bra\omega}_{\mathsf j}
    =
    \mathpzc t^{\pare2}_{h^{\mathpzc B_0}_{\mathsf j},
    k^{\mathpzc B_0}}\pare\eta\bra\omega,
    \qquad
    \mathpzc q\pare\eta\bra\omega
    =
    \mathpzc t^{\pare2}_{h^{\mathpzc q},
    k^{\mathpzc q}}\pare\eta\bra\omega.
\end{equation}
\end{lemma}

\begin{proof}
\begin{proofpart}[Parity and vanishing order]
By construction \cref{eq:B0-q-admissible-kernels}, the vector amplitude \(H^{\mathpzc B_0}\) is odd,
\(k^{\mathpzc B_0}\) is even, and
\[
    k^{\mathpzc B_0}\pare X-k^{\mathpzc B_0}\pare0
    =O\pare{X^2}.
\]
Thus each component pair is admissible, in the sense of \cref{def:abstract-admissible-angular-pair-Rd}, with \(\varepsilon=-1\) and
\(\ell=2\).  The amplitude \(h^{\mathpzc q}\) is even,
\(k^{\mathpzc q}\) is odd, and
$$k^{\mathpzc q}\pare 0=0,\qquad k^{\mathpzc q}\pare X=O\pare X.$$ Hence the second pair is admissible, in the sense of \cref{def:abstract-admissible-angular-pair-Rd},
with \(\varepsilon=1\) and \(\ell=1\).
\end{proofpart}

\begin{proofpart}[Exact operator identification]
Since \(\delta_z\eta=\av z\Delta_z\eta\), the definitions in
\cref{eq:B0_recall,eq:Q_recall} become
\[
\begin{aligned}
    \mathpzc B_0\pare\eta\bra\omega\pare x
    &=
    \pv\int_{\bR^2}
    \frac{H^{\mathpzc B_0}\pare{\hat z}}{\av z^2}
    k^{\mathpzc B_0}\pare{\Delta_z\eta\pare x}
    \omega\pare{x-z}\dd z,
    \\
    \mathpzc q\pare\eta\bra\omega\pare x
    &=
    \pv\int_{\bR^2}
    \frac{h^{\mathpzc q}\pare{\hat z}}{\av z^2}
    k^{\mathpzc q}\pare{\Delta_z\eta\pare x}
    \omega\pare{x-z}\dd z.
\end{aligned}
\]
This is precisely \cref{eq:B0-q-abstract-operator-identification}, with
the first identity understood componentwise.
\end{proofpart}
\end{proof}
With this in hand, we can now apply \cref{thm:abstract-structure-preserving-paralinearization-Rd} to obtain the paralinearization of the operators $\mathpzc{B}_0$ and $\mathpzc{q}$. The final results are stated in the following two propositions. Their complete proofs are given in the following subsections.

\phantomsection
\label{sec:parab0}
\begin{prop}[Paralinearization of $\mathpzc B_0\pare\eta$]
\label{prop:B0_paralinearization}
There exists \(\varrho_0>0\) such that, for every
\(\varrho\geq\varrho_0\), there exist \(s_0>2\), \(\mu\geq0\), and
\(r_0>0\) such that, for every \(s\geq s_0\) and every real
\[
    \eta\in H^{s+1+\mu}\pare{\bR^2},
    \qquad
    \omega\in H^{s+\mu}\pare{\bR^2},
    \qquad
    \norm\eta_{H^{s_0+1+\mu}}<r_0,
\]
the following identity holds in \(H^s\pare{\bR^2;\bR^2}\):
\begin{equation}\label{eq:paralinearization_B0}
    \mathpzc B_0\pare\eta\bra\omega
    =
    \OpBW{
        {\bm A}^{\mathpzc B_0}_1
        \pare{\eta,\omega;x,\xi}
        +
        \vectT{
            A_{\bra{0,\geq2}}\bra\omega
            \pare{\eta;x,\xi}
        }{
            \ii A^{\mathpzc B_0,\bR}_{\bra{-2,\geq2}}
                \pare{\eta;x,\xi}
        }
    }
    \vect{\eta}{\omega}
    +
    {\bm R}^{\mathpzc B_0,\lin,\bR}_{\bra{\varrho,\geq2}}
    \pare{\eta,\omega}
    \vect{\eta}{\omega}.
\end{equation}
The objects in \cref{eq:paralinearization_B0} are as follows.
\begin{itemize}
    \item The principal matrix symbol is
    \begin{equation}\label{eq:AB0_matrix}
        {\bm A}^{\mathpzc B_0}_1
        \pare{\eta,\omega;x,\xi}
        \defeq
        \begin{bmatrix}
            \omega\pare x
            A^{\mathpzc B_0,\eta}_1
            \pare{\eta;x,\xi}
            &
            \ii A_{0}^{\mathpzc{B}_0,\omega}\pare{\eta;x,\xi}
        \end{bmatrix}.
    \end{equation}

    \item The common potential is the real symbol
    \(\phi^{\mathpzc B_0}_1\in\Gamma^1_{\geq0}\) defined by
    \begin{equation}\label{eq:phiB0}
        \phi^{\mathpzc B_0}_1\pare{\eta;x,\xi}
        \defeq
        -
        \frac{\lambda_1\pare{\eta;x,\xi}}
        {1+\av{\nabla\eta\pare x}^2},
    \end{equation}
    where $\lambda_1$ is defined in \cref{eq:lambda1},

    \item The real vector symbol
    \(A_{0}^{\mathpzc{B}_0,\omega}\in\pare{\Gamma^0_{\geq0}}^2\), defined in
    \cref{eq:H0_definition}, satisfies
    \begin{equation}\label{eq:B0-H0-potential-identity}
        A_{0}^{\mathpzc{B}_0,\omega}\pare{\eta;x,\xi}
        =
        \nabla_\xi\phi^{\mathpzc B_0}_1
        \pare{\eta;x,\xi}.
    \end{equation}

    \item The real vector symbol
    \(A^{\mathpzc B_0,\eta}_1
    \in\pare{\Gamma^1_{\geq1}}^2\) is
    \begin{equation}\label{eq:AB0eta1}
    \begin{aligned}
        A^{\mathpzc B_0,\eta}_1
        \pare{\eta;x,\xi}
        &\defeq
        -
        \frac{\partial\phi^{\mathpzc B_0}_1}
        {\partial\pare{\nabla\eta}}
        \pare{\eta;x,\xi}
        \\
        &=
        \frac{
            A_{0}^{\mathpzc{B}_0,\omega}\pare{\eta;x,\xi}\otimes\nabla\eta\pare x
            +
            \nabla\eta\pare x\otimes A_{0}^{\mathpzc{B}_0,\omega}\pare{\eta;x,\xi}
        }{
            1+\av{\nabla\eta\pare x}^2
        }
        \xi.
    \end{aligned}
    \end{equation}
    
    \item The lower-order symbols satisfy
    \begin{equation}\label{eq:B0-lower-symbol-classes}
        A_{\bra{0,\geq2}}\bra\cdot
        \in
        \cL\pare{H^{s_0};\pare{\Gamma^0_{\geq1}}^2},
        \qquad
        A^{\mathpzc B_0,\bR}_{\bra{-2,\geq2}}
        \in
        \pare{\Gamma^{-2}_{\geq2}}^2.
    \end{equation}
    The vector symbol
    \(A^{\mathpzc B_0,\bR}_{\bra{-2,\geq2}}\) is real-valued componentwise and
    odd in \(\xi\), namely
    \begin{equation}\label{eq:B0-lower-symbol-parity}
        A^{\mathpzc B_0,\bR}_{\bra{-2,\geq2}}
            \pare{\eta;x,-\xi}
            =
            -A^{\mathpzc B_0,\bR}_{\bra{-2,\geq2}}
            \pare{\eta;x,\xi}.
    \end{equation}

    \item The smoothing matrix is
    \begin{equation}\label{eq:B0-smoothing-matrix}
        {\bm R}^{\mathpzc B_0,\lin,\bR}_{\bra{\varrho,\geq2}}
        \pare{\eta,\omega}
        \defeq
        \vectT{
            R^{\mathpzc B_0,\bR}_{\bra{\varrho,\geq2}}
            \bra\omega\pare\eta
        }{
            R^{\mathpzc B_0,\bR}_{\bra{\varrho,\geq2}}\pare\eta
        },
    \end{equation}
    where
    \begin{equation}\label{eq:B0-smoothing-classes}
        R^{\mathpzc B_0,\bR}_{\bra{\varrho,\geq2}}\bra\cdot
        \in
        \cL\pare{H^{s_0};\pare{\cR^\varrho_{\geq1}}^2},
        \qquad
        R^{\mathpzc B_0,\bR}_{\bra{\varrho,\geq2}}
        \in
        \pare{\cR^\varrho_{\geq2}}^2.
    \end{equation}
\end{itemize}
All the terms in \cref{eq:paralinearization_B0} are real-to-real.
\end{prop}

\phantomsection
\label{sec:paraq0}
\begin{prop}[Paralinearization of $\mathpzc q\pare\eta$]
\label{prop:paralinearization_Q}
There exists \(\varrho_0>0\) such that, for every
\(\varrho\geq\varrho_0\), there exist \(s_0>2\), \(\mu\geq0\), and
\(r_0>0\) such that, for every \(s\geq s_0\) and every real
\[
    \eta\in H^{s+1+\mu}\pare{\bR^2},
    \qquad
    \omega\in H^{s+\mu}\pare{\bR^2},
    \qquad
    \norm\eta_{H^{s_0+1+\mu}}<r_0,
\]
the following identity holds in \(H^s\pare{\bR^2}\):
\begin{equation}\label{eq:paralinearization_Q}
    \mathpzc q\pare\eta\bra\omega
    =
    \OpBW{
        A^{\mathpzc q}_1
        \pare{\eta,\omega;x,\xi}
        +
        \vect{
            a_{\bra{0,\geq1}}\bra\omega
            \pare{\eta;x,\xi}
        }{
            \ii a^{\mathpzc q,\bR}_{\bra{-2,\geq1}}
                \pare{\eta;x,\xi}
        }
    }
    \cdot\vect{\eta}{\omega}
    +
    R^{\mathpzc q,\lin,\bR}_{\bra{\varrho,\geq1}}
    \pare{\eta,\omega}
    \cdot\vect{\eta}{\omega}.
\end{equation}
The objects in \cref{eq:paralinearization_Q} are as follows.
\begin{itemize}
    \item The principal vector symbol is
    \begin{equation}\label{eq:Aq1_definition}
        A^{\mathpzc q}_1\pare{\eta,\omega;x,\xi}
        \defeq
        \vect{
            \omega\pare x
            a^{\mathpzc q,\eta}_1\pare{\eta;x,\xi}
        }{
            \ii a^{\mathpzc q,\omega}_0\pare{\eta;x,\xi}
        }.
    \end{equation}
    \begin{itemize}
           \item The real symbol
    \(a^{\mathpzc q,\omega}_0\in\Gamma^0_{\geq1}\) is
    \begin{equation}\label{eq:AB0omega_AQL0}
        a^{\mathpzc q,\omega}_0\pare{\eta;x,\xi}
        \defeq
        \nabla\eta\pare x\cdot
        A_{0}^{\mathpzc{B}_0,\omega}\pare{\eta;x,\xi}.
    \end{equation}

    \item The real symbol
    \(a^{\mathpzc q,\eta}_1\in\Gamma^1_{\geq0}\) is
    \begin{equation}\label{eq:relationa1A1}
        a^{\mathpzc q,\eta}_1\pare{\eta;x,\xi}
        \defeq
        \nabla\eta\pare x\cdot
        A^{\mathpzc B_0,\eta}_1\pare{\eta;x,\xi}
        -
        A_{0}^{\mathpzc{B}_0,\omega}\pare{\eta;x,\xi}\cdot\xi.
    \end{equation}
    \end{itemize}

    \item The lower-order symbols satisfy
    \begin{equation}\label{eq:q-lower-symbol-classes}
        a_{\bra{0,\geq1}}\bra\cdot
        \in
        \cL\pare{H^{s_0};\Gamma^0},
        \qquad
        a^{\mathpzc q,\bR}_{\bra{-2,\geq1}}
        \in
        \Gamma^{-2}_{\geq1}.
    \end{equation}
    The symbol
    \(a^{\mathpzc q,\bR}_{\bra{-2,\geq1}}\) is real-valued and odd in \(\xi\),
    namely
    \begin{equation}\label{eq:q-lower-symbol-parity}
        a^{\mathpzc q,\bR}_{\bra{-2,\geq1}}
            \pare{\eta;x,-\xi}
            =
            -a^{\mathpzc q,\bR}_{\bra{-2,\geq1}}
            \pare{\eta;x,\xi}.
    \end{equation}
  
    \item The smoothing vector is
    \begin{equation}\label{eq:q-smoothing-vector}
        R^{\mathpzc q,\lin,\bR}_{\bra{\varrho,\geq1}}
        \pare{\eta,\omega}
        \defeq
        \vect{
            r^{\mathpzc q,\bR}_{\bra{\varrho,\geq1}}
            \bra\omega\pare\eta
        }{
            r^{\mathpzc q,\bR}_{\bra{\varrho,\geq1}}\pare\eta
        },
    \end{equation}
    where
    \begin{equation}\label{eq:q-smoothing-classes}
        r^{\mathpzc q,\bR}_{\bra{\varrho,\geq1}}\bra\cdot
        \in
        \cL\pare{H^{s_0};\cR^\varrho},
        \qquad
        r^{\mathpzc q,\bR}_{\bra{\varrho,\geq1}}
        \in
        \cR^\varrho_{\geq1}.
    \end{equation}
\end{itemize}
All the terms in \cref{eq:paralinearization_Q} are real-to-real.
\end{prop}

\subsubsection{Identification of the principal symbols}
According to \cref{lem:B0-q-admissible-kernels}, we can apply \cref{eq:abstract-explicit-principal-symbol-A0-Rd,%
eq:abstract-explicit-principal-symbol-A1-Rd} with
\(\mathtt d=2\). This provides \cref{eq:paralinearization_B0}, with 
$${\bm A}^{\mathpzc B_0}_1
        \pare{\eta,\omega;x,\xi}
        \defeq
        \begin{bmatrix}
            \omega\pare x
            A^{\mathpzc B_0,\eta}_1
            \pare{\nabla\eta\pare{x},\xi}
            &
            \ii A^{\mathpzc B_0,\omega}_0\pare{\nabla\eta\pare{x},\xi}
        \end{bmatrix}$$
        and \cref{eq:paralinearization_Q}, with
        $$A^{\mathpzc q}_1\pare{\eta,\omega;x,\xi}
        \defeq
        \vect{
            \omega\pare x
            a^{\mathpzc q,\eta}_1\pare{\nabla\eta\pare{x},\xi}
        }{
            \ii a^{\mathpzc q,\omega}_0\pare{\nabla\eta\pare{x},\xi},
        }$$
        where
\begin{subequations}\label{eq:B0-q-angular-principal-symbols}
\begin{align}
    A^{\mathpzc B_0,\omega}_0\pare{p,\xi}
    &\defeq
    -\frac18
    \int_{\bS^1}
    \sgn\pare{\vartheta\cdot\xi}
    \vartheta
    k^{\mathpzc B_0}\pare{p\cdot\vartheta}
    \dd\sigma\pare\vartheta,
    \label{eq:AB0omega}
    \\
    A^{\mathpzc B_0,\eta}_1\pare{p,\xi}
    &\defeq
    \frac18
    \int_{\bS^1}
    \av{\vartheta\cdot\xi}
    \vartheta
    \pare{k^{\mathpzc B_0}}'\pare{p\cdot\vartheta}
    \dd\sigma\pare\vartheta,
    \label{eq:AB0eta1_int}
    \\
    a^{\mathpzc q,\omega}_0\pare{p,\xi}
    &\defeq
    -\frac18
    \int_{\bS^1}
    \sgn\pare{\vartheta\cdot\xi}
    k^{\mathpzc q}\pare{p\cdot\vartheta}
    \dd\sigma\pare\vartheta,
    \label{eq:aq0omega}
    \\
    a^{\mathpzc q,\eta}_1\pare{p,\xi}
    &\defeq
    \frac18
    \int_{\bS^1}
    \av{\vartheta\cdot\xi}
    \pare{k^{\mathpzc q}}'\pare{p\cdot\vartheta}
    \dd\sigma\pare\vartheta.
    \label{eq:AB0eta1Q}
\end{align}
\end{subequations}
Note that the factor \(1/(4\pi)\) in
\cref{eq:B0-q-admissible-kernels} turns \(\pi/2\) into \(1/8\), which
gives exactly \cref{eq:B0-q-angular-principal-symbols}. What remains is to give the explicit expressions of these integral formulae. We have the following lemma.

\begin{lemma}[Potential representation and Cartesian form of the principal symbols]
\label{lem:eformB0}
Let \(k^{\mathpzc B_0}\) be defined in
\cref{eq:B0-q-admissible-kernels}, and let
\(A^{\mathpzc B_0,\omega}_0\) and
\(A^{\mathpzc B_0,\eta}_1\) be the angular symbols defined in
\cref{eq:AB0omega,eq:AB0eta1_int}. For \(p\in\bR^2\) and
\(\xi\neq0\), define
\begin{equation}\label{eq:phiB0_int}
    \phi^{\mathpzc B_0}_{\mathrm{int}}\pare{p,\xi}
    \defeq
    -\frac18
    \int_{\bS^1}
    \av{\vartheta\cdot\xi}
    k^{\mathpzc B_0}\pare{p\cdot\vartheta}
    \dd\sigma\pare\vartheta.
\end{equation}
Then
\begin{equation}\label{eq:phiB0_int_full_circle}
    \phi^{\mathpzc B_0}_{\mathrm{int}}\pare{p,\xi}
    =
    -\frac{\av\xi}{2}
    \frac{
        \sqrt{
            1+\pare{p^\perp\cdot\xi/\av\xi}^2
        }
    }{
        1+\av p^2
    }\cdot
\end{equation}
Moreover,
\begin{equation}\label{eq:phiB0_int_differential_identifications}
\begin{aligned}
    \nabla_\xi
    \phi^{\mathpzc B_0}_{\mathrm{int}}\pare{p,\xi}
    &=
    A^{\mathpzc B_0,\omega}_0\pare{p,\xi},
    %\label{eq:phiB0_int_xi_derivative}
    \\
    -
    \partial_p
    \phi^{\mathpzc B_0}_{\mathrm{int}}\pare{p,\xi}
    &=
    A^{\mathpzc B_0,\eta}_1\pare{p,\xi}.
    %\label{eq:phiB0_int_p_derivative}
\end{aligned}
\end{equation}
For \(p=\nabla\eta\pare x\), one has
\begin{equation}\label{eq:phiB0_int_identification}
    \phi^{\mathpzc B_0}_{\mathrm{int}}
    \pare{\nabla\eta\pare x,\xi}
    =
    \phi^{\mathpzc B_0}_1\pare{\eta;x,\xi},
\end{equation}
where \(\phi^{\mathpzc B_0}_1\) is defined in \cref{eq:phiB0}.
Consequently,
\begin{equation}\label{eq:H0_Aomega_identification}
\begin{aligned}
    A_{0}^{\mathpzc{B}_0,\omega}\pare{\eta;x,\xi}\defeq A^{\mathpzc B_0,\omega}_0
    \pare{\nabla\eta\pare x,\xi}
    &=
    \nabla_\xi
    \phi^{\mathpzc B_0}_{\mathrm{int}}
    \pare{\nabla\eta\pare x,\xi}
    \\
    &=
    \nabla_\xi
    \phi^{\mathpzc B_0}_1\pare{\eta;x,\xi}
    .
\end{aligned}
\end{equation}
Furthermore,
\begin{equation}\label{eq:B0-A1-potential-identification}
    A^{\mathpzc B_0,\eta}_1
    \pare{\eta;x,\xi}\defeq A^{\mathpzc B_0,\eta}_1
    \pare{\nabla\eta\pare x,\xi}
    =
    -
    \frac{\partial\phi^{\mathpzc B_0}_1}
    {\partial\pare{\nabla\eta}}
    \pare{\eta;x,\xi}.
\end{equation}
The latter identity is equivalently the Cartesian formula
\cref{eq:AB0eta1}.
\end{lemma}

\begin{proof}
\begin{proofpart}[Evaluation of the angular potential]
By rotation covariance, it is enough to take
\(\xi=\av\xi\vect{0}{1}\) and \(p=\vect{p_1}{p_2}\).  Writing
\(\vartheta\defeq\vect{\cos\theta}{\sin\theta}\) and using the
symmetry of the integrand, we get
\begin{equation}\label{eq:B0-angular-potential-one-dimensional-integral}
\begin{aligned}
    &\int_{\bS^1}
    \frac{\av{\vartheta\cdot\xi}}
    {\pare{1+\pare{p\cdot\vartheta}^2}^{\frac32}}
    \dd\sigma\pare\vartheta
    =
    2\av\xi
    \int_0^\pi
    \frac{\sin\theta}
    {\pare{
        1+\pare{p_1\cos\theta+p_2\sin\theta}^2
    }^{\frac32}}
    \dd\theta.
\end{aligned}
\end{equation}
The antiderivative
\[
    -
    \frac{
        \pare{1+p_1^2}\cos\theta
        +p_1p_2\sin\theta
    }{
        \pare{1+p_1^2+p_2^2}
        \sqrt{
            1+\pare{p_1\cos\theta+p_2\sin\theta}^2
        }
    }
\]
gives
\[
    \int_0^\pi
    \frac{\sin\theta}
    {\pare{
        1+\pare{p_1\cos\theta+p_2\sin\theta}^2
    }^{\frac32}}
    \dd\theta
    =
    \frac{2\sqrt{1+p_1^2}}{1+\av p^2}\cdot
\]
Since \(p_1^2=\pare{p^\perp\cdot\xi/\av\xi}^2\) in the chosen frame,
substitution in \cref{eq:phiB0_int} proves
\cref{eq:phiB0_int_full_circle}. 
\end{proofpart}

\begin{proofpart}[Differential identification]
With \(p=\nabla\eta\pare x\),
\cref{eq:phiB0_int_full_circle,eq:lambda1} give
\cref{eq:phiB0_int_identification}.  Differentiating
\cref{eq:phiB0_int} under the integral sign yields
\[
    \nabla_\xi\phi^{\mathpzc B_0}_{\mathrm{int}}
    =A^{\mathpzc B_0,\omega}_0,
    \qquad
    -\partial_p\phi^{\mathpzc B_0}_{\mathrm{int}}
    =A^{\mathpzc B_0,\eta}_1.
\]
This proves \cref{eq:phiB0_int_differential_identifications}. Together with
\cref{eq:H0_definition,eq:phiB0_int_identification}, this proves
\cref{eq:H0_Aomega_identification,eq:B0-A1-potential-identification}.
Finally, differentiating the Cartesian expression
\cref{eq:phiB0_int_full_circle} gives \cref{eq:AB0eta1}.
\end{proofpart}
\end{proof}

\begin{lemma}[Relation between the principal symbols of
$\mathpzc B_0$ and $\mathpzc q$]
\label{lem:B0-q-principal-symbol-relation}
\label{step:relation_B0_Q_symbols}
For \(p=\nabla\eta\pare x\), the angular symbols in
\cref{eq:aq0omega,eq:AB0eta1Q} satisfy
\begin{equation}\label{eq:qomega-symbol-identification}
    a^{\mathpzc q,\omega}_0\pare{\eta;x,\xi}\defeq a^{\mathpzc q,\omega}_0\pare{\nabla\eta\pare{x},\xi}
    =
    \nabla\eta\pare x\cdot
    A_{0}^{\mathpzc{B}_0,\omega}\pare{\eta;x,\xi},
\end{equation}
and
\begin{equation}\label{eq:qeta-symbol-identification}
    a^{\mathpzc q,\eta}_1\pare{\eta;x,\xi}\defeq a^{\mathpzc q,\eta}_1\pare{\nabla\eta\pare{x},\xi}
    =
    \nabla\eta\pare x\cdot
    A^{\mathpzc B_0,\eta}_1\pare{\eta;x,\xi}
    -
    A_{0}^{\mathpzc{B}_0,\omega}\pare{\eta;x,\xi}\cdot\xi.
\end{equation}
Thus these are precisely the symbols in
\cref{eq:AB0omega_AQL0,eq:relationa1A1}.
\end{lemma}

\begin{proof}
\begin{proofpart}[Order-zero symbol]
The identities in \cref{eq:B0-q-admissible-kernels} give
\begin{equation}\label{eq:formmattK1}
    k^{\mathpzc q}\pare X
    =Xk^{\mathpzc B_0}\pare X.
\end{equation}
Hence \cref{eq:aq0omega,eq:AB0omega} imply
\[
\begin{aligned}
    a^{\mathpzc q,\omega}_0\pare{p,\xi}
    &=-\frac18
    \int_{\bS^1}
    \sgn\pare{\vartheta\cdot\xi}
    \pare{p\cdot\vartheta}
    k^{\mathpzc B_0}\pare{p\cdot\vartheta}
    \dd\sigma\pare\vartheta
    \\
    &=p\cdot A^{\mathpzc B_0,\omega}_0\pare{p,\xi}.
\end{aligned}
\]
Taking \(p=\nabla\eta\pare x\) and using
\cref{eq:H0_Aomega_identification} proves
\cref{eq:qomega-symbol-identification}.
\end{proofpart}

\begin{proofpart}[Order-one symbol]
Differentiating \cref{eq:formmattK1},
\[
    \pare{k^{\mathpzc q}}'\pare X
    =
    k^{\mathpzc B_0}\pare X
    +X\pare{k^{\mathpzc B_0}}'\pare X.
\]
Therefore \cref{eq:AB0eta1Q} becomes
\[
\begin{aligned}
    a^{\mathpzc q,\eta}_1\pare{p,\xi}
    ={}&
    \frac18\int_{\bS^1}
    \av{\vartheta\cdot\xi}
    k^{\mathpzc B_0}\pare{p\cdot\vartheta}
    \dd\sigma\pare\vartheta
    +
    p\cdot A^{\mathpzc B_0,\eta}_1\pare{p,\xi}
    \\
    ={}&
    -\phi^{\mathpzc B_0}_{\mathrm{int}}\pare{p,\xi}
    +
    p\cdot A^{\mathpzc B_0,\eta}_1\pare{p,\xi}.
\end{aligned}
\]
Since \(\phi^{\mathpzc B_0}_{\mathrm{int}}\) is one-homogeneous in
\(\xi\), \cref{eq:H0_Aomega_identification} and Euler's identity give
\(
    \phi^{\mathpzc B_0}_{\mathrm{int}}
    =A_{0}^{\mathpzc{B}_0,\omega}\cdot\xi
\).
This proves \cref{eq:qeta-symbol-identification}.
\end{proofpart}
\end{proof}

\begin{cor}
\label{lem:relationimportant}
One has
\begin{equation}\label{eq:relationimportant3}
    \phi^{\mathpzc B_0}_1\pare{\eta;x,\xi}
    =
    \xi\cdot A_{0}^{\mathpzc{B}_0,\omega}\pare{\eta;x,\xi}
    =
    -
    \frac{\lambda_1\pare{\eta;x,\xi}}
    {1+\av{\nabla\eta\pare x}^2},
\end{equation}
and
\begin{equation}\label{eq:relationimportant4}
    -
    \pare{\nabla\eta\pare x\cdot\xi}
    A_{0}^{\mathpzc{B}_0,\omega}\pare{\eta;x,\xi}
    +
    \pare{1+\av{\nabla\eta\pare x}^2}
    A^{\mathpzc B_0,\eta}_1\pare{\eta;x,\xi}
    =
    -
    \lambda_1\pare{\eta;x,\xi}
    \frac{\nabla\eta\pare x}
    {1+\av{\nabla\eta\pare x}^2}\cdot
\end{equation}
\end{cor}

\begin{proof}
Euler's identity applied to the one-homogeneous symbol
\(\phi^{\mathpzc B_0}_1\), together with
\cref{eq:B0-H0-potential-identity,eq:phiB0}, gives
\cref{eq:relationimportant3}.  Multiplying \cref{eq:AB0eta1} by
\(1+\av{\nabla\eta}^2\) gives
\[
    \pare{1+\av{\nabla\eta}^2}
    A^{\mathpzc B_0,\eta}_1
    =
    \pare{\nabla\eta\cdot\xi}A_{0}^{\mathpzc{B}_0,\omega}
    +
    \pare{A_{0}^{\mathpzc{B}_0,\omega}\cdot\xi}\nabla\eta.
\]
Now \cref{eq:relationimportant3} yields
\cref{eq:relationimportant4}.
\end{proof}

\subsubsection{Synthesis of the two paralinearizations}

\begin{proof}[Proof of
\cref{prop:B0_paralinearization,prop:paralinearization_Q}]
\begin{proofpart}[Application to $\mathpzc B_0$]
Apply \cref{thm:abstract-structure-preserving-paralinearization-Rd}
componentwise to the two admissible pairs of vanishing order
\(\ell=2\) in \cref{lem:B0-q-admissible-kernels}.  The angular symbols
produced by the theorem are \cref{eq:AB0omega,eq:AB0eta1_int}; by
\cref{lem:eformB0}, they are respectively \(A_{0}^{\mathpzc{B}_0,\omega}\) and
\(A^{\mathpzc B_0,\eta}_1\).  Thus the two principal columns are exactly
those of \cref{eq:AB0_matrix}.  Moreover,
\cref{eq:abstract-imaginary-lower-symbol-parity-Rd} gives the second
lower symbol componentwise in the form
\(\ii A^{\mathpzc B_0,\bR}_{\bra{-2,\geq2}}\), where
\(A^{\mathpzc B_0,\bR}_{\bra{-2,\geq2}}\) is real-valued componentwise and
odd in \(\xi\).  Since the bracketed terms are linear in \(\omega\), their
vanishing order is \(\ell=2\) by
\cref{rem:linear-parameter-vanishing-order}. This gives
\cref{eq:B0-lower-symbol-classes,eq:B0-smoothing-classes}, proving
\cref{eq:paralinearization_B0}.
\end{proofpart}

\begin{proofpart}[Application to $\mathpzc q$]
Apply the same theorem to the admissible pair of vanishing order
\(\ell=1\) in \cref{lem:B0-q-admissible-kernels}.  The leading angular
symbols are \cref{eq:aq0omega,eq:AB0eta1Q};
\cref{lem:B0-q-principal-symbol-relation} identifies them with
\cref{eq:AB0omega_AQL0,eq:relationa1A1}.  Moreover,
\cref{eq:abstract-imaginary-lower-symbol-parity-Rd} gives the second
lower symbol in the form \(\ii a^{\mathpzc q,\bR}_{\bra{-2,\geq1}}\), where
\(a^{\mathpzc q,\bR}_{\bra{-2,\geq1}}\) is real-valued and odd in \(\xi\).
The transposed terms vanish to order \(\ell-1=0\) in \(\eta\) and are
linear in \(\omega\); hence their total vanishing order is \(\ell=1\) by
\cref{rem:linear-parameter-vanishing-order}. The quasilinear terms have the
same vanishing order. This gives
\cref{eq:q-lower-symbol-classes,eq:q-smoothing-classes} and
proves \cref{eq:paralinearization_Q}.
\end{proofpart}
\end{proof}

\subsection{Paralinearization of evolution operators}\label{sec:paralinearization_EO}
Based on the paralinearization of the elementary singular integral operators $\mathpzc{B}_0$ and $\mathpzc{q}$ carried out in the previous subsection, we provide here the paralinearization of each component of the evolution equations \cref{eq:CDE}. We analyze each brick separately to finally gather, in the next subsection, the obtained results, thus providing the complete paralinearization of the generalized water-waves system.

\subsubsection{Paralinearization of \texorpdfstring{${\bm v}^\star$}{tilde v}}
Let ${\bm v}^\star$ be defined as in \cref{eq:bmvtilde}, which is 
    \begin{align}
     {\bm v}^\star\bra{\Omega}\pare{\eta}\defeq & \ \vect{ V^{\star} \bra{\Omega}\pare{\eta}}{ v_3^{\star} \bra{\Omega}\pare{\eta}} \defeq
      \vect{\pare{\mathpzc{B}_0\pare{\eta} \bra{\nabla\eta \cdot \Omega} - \mathpzc{q}\pare{\eta}\Omega}^\perp}{-\mathpzc{B}_0\pare{\eta}^\perp \cdot \Omega}.
      \label{eq:bmvtilde_recall}
\end{align}
Then, we have the following

\begin{prop}[Paralinearization of $\Omega\mapsto {\bm v}^\star\bra{\Omega}\pare{\eta}$]\label{prop:paralinearization_Vtilde}
    There exists $\varrho_0 > 0$ such that for any $\varrho \geq \varrho_0$, there exists $s_0 > 0$ and $r_0 > 0$ such that for any $\eta\in \BallR{s_0}{r_0}, \Omega \in\pare{\BallR{s_0}{r_0}}^2$  we have the following:
\begin{enumerate}
    \item For $ V^{\star}$ defined as in \cref{eq:bmvtilde_recall} the following paralinear decomposition holds true
    {}
\begin{equation}\label{eq:paralinearization_Vtilde}
           V^{\star}\bra{\Omega}\pare{\eta;x}
           = \OpBW{\mathfrak{A}^{ V^{\star}}_1 \pare{\eta,\Omega;x,\xi} + \vectT{ A_{\bra{0,\geq1}}\bra{\Omega}\pare{\eta;x,\xi}}{{\bm A}_{\bra{-2,\geq1}}\pare{\eta;x,\xi}}}\vect{\eta}{\Omega}
          + \mathfrak{R}^{\lin,\bR}_{\bra{\varrho,\geq1}}\pare{\eta, \Omega} \vect{\eta}{\Omega},
    \end{equation}
    where
    \begin{equation}\label{eq:frakAVtilde}
        \mathfrak{A}^{ V^{\star}}_1 \pare{\eta,\Omega;x,\xi} \defeq \vectT{ {\bm A}^{ V^{\star}, \eta}_1\pare{\eta;x,\xi}\Omega}{ \pare{\ii {\bm A}^{ V^{\star},\Omega}_0 + {\bm A}^{ V^{\star},\Omega}_{-1}}\pare{\eta;x,\xi}},
    \end{equation}
    with
    \begin{itemize}

\item  $ {\bm A}^{ V^{\star}, \eta}_1\in \pare{\Gamma_{\geq0}^1}^{2\times 2}$ and is explicitly defined as
      \begin{equation}\label{eq:AVtildeeta1_lem}
        {\bm A}^{ V^{\star}, \eta}_1\pare{\eta;x,\xi}\defeq 
        -\xi\otimes A_{0}^{\mathpzc{B}_0,\omega;\perp}\pare{\eta;x,\xi}
        +
        \nabla\eta  \otimes A^{\mathpzc{B}_0,\eta;\perp}_1\pare{\eta;x,\xi} ,
    \end{equation}
    where $A_{0}^{\mathpzc{B}_0,\omega}$ is defined in \cref{eq:H0_definition} and $A^{\mathpzc{B}_0,\eta}_1$ is defined in \cref{eq:AB0eta1};

        \item $  {\bm A}^{ V^{\star},\Omega}_0 \in \pare{\Gamma_{\geq1}^0}^{2\times 2}$ is the real-valued symbol explicitly defined as 
          \begin{equation}\label{eq:AVtildeg0_lem}
        {\bm A}^{ V^{\star},\Omega}_0\pare{\eta;x,\xi}\defeq\nabla\eta \otimes A_{0}^{\mathpzc{B}_0,\omega;\perp}\pare{\eta;x,\xi},
    \end{equation}
    where $A_{0}^{\mathpzc{B}_0,\omega}$ is defined in \cref{eq:H0_definition};

\item $  {\bm A}^{ V^{\star},\Omega}_{-1} \in \pare{\Gamma_{\geq1}^{-1}}^{2\times 2}$ and is explicitly defined as 
\begin{equation}\label{eq:AVtildeg-1_lem}
{\bm A}^{ V^{\star},\Omega}_{-1}\pare{\eta;x,\xi} 
 = \frac{\bm J}{2} J_\xi A_{0}^{\mathpzc{B}_0,\omega} {\bm H}_x \eta ,
      \end{equation}
     where $A_{0}^{\mathpzc{B}_0,\omega}$ is defined in \cref{eq:H0_definition};

    \end{itemize}

         \item For $ v_3^{\star}$ defined as in \cref{eq:bmvtilde_recall} the following paralinear decomposition holds true
         \begin{equation}\label{eq:paralinearization_tildev3}
           \begin{aligned}
               v_3^{\star}\bra{\Omega}\pare{\eta} 
              & = \OpBW{{\bm a}^{ v_3^{\star}}_1\pare{\eta,\Omega;x,\xi} + \vect{a_{\bra{0,\geq1}}\bra{\Omega}\pare{\eta;x, \xi} }{A_{\bra{-2,\geq1}}\pare{\eta;x, \xi}}} \cdot\vect{\eta}{\Omega}
              +
               {\bm r}^{\lin,\bR}_{\bra{\varrho,\geq1}}\pare{\eta, \Omega}\cdot\vect{\eta}{\Omega},
       \end{aligned}
    \end{equation}
    where
    \begin{equation}\label{eq:atildev3_symbol}
    {\bm a}^{ v_3^{\star}}_1\pare{\eta,\Omega;x,\xi}
    \defeq
    \vect{
        a_1^{ v_3^{\star},\eta}\bra{\Omega}\pare{\eta;x,\xi}
    }{
        \ii A^{ v_3^{\star},\Omega}_0\pare{\eta;x,\xi}
    }
    \defeq
    \vect{
        -A^{\mathpzc{B}_0,\eta;\perp}_1\pare{\eta;x,\xi}\cdot\Omega
    }{
        -\ii A_{0}^{\mathpzc{B}_0,\omega;\perp}\pare{\eta;x,\xi}
    } .
\end{equation}
Here
\(
    a_1^{ v_3^{\star},\eta}\bra{\Omega}\in\Gamma_{\geq1}^1
\)
is linear in \(\Omega\), while
\(
    A^{ v_3^{\star},\Omega}_0=-A_{0}^{\mathpzc{B}_0,\omega;\perp}\in\pare{\Gamma_{\geq0}^0}^2
\).
Equivalently, the \(\Omega\)-slot in
\cref{eq:paralinearization_tildev3} is the row symbol
\(-\ii\pare{A_{0}^{\mathpzc{B}_0,\omega;\perp}}^{\intercal}\).
\end{enumerate}
\end{prop}

\begin{proof}
We prove separately the horizontal and vertical components of
\({\bm v}^\star=\vect{ V^{\star}}{ v_3^{\star}}\).  Throughout the proof,
all generic symbols and remainders are understood in the sense of
\cref{notation:paralinearization}.  We also freely enlarge \(s_0\)
and decrease \(r_0\), so that all products below are admissible in the
classes stated in \cref{prop:B0_paralinearization,prop:paralinearization_Q}.

\begin{proofpart}[Paralinearization of
\texorpdfstring{$\Omega\mapsto V^{\star}\bra{\Omega}\pare{\eta}$}{Omega -> Vtilde}]
Set
\[
    \omega_\eta\defeq \nabla\eta\cdot\Omega .
\]
Since \(\omega_\eta\) is obtained from \((\eta,\Omega)\) by a Bony
product, \cref{lem:paraprod_Bony,rem:smoothBONYremainder} and
\cref{lem:composition_rules}\,\textup{(iii)--(iv)} allow us to substitute
\(\omega=\omega_\eta\) in the paralinearization of
\(\mathpzc{B}_0\).  Thus, by
\cref{prop:B0_paralinearization},
\begin{equation}\label{eq:toughcal1}
    \begin{aligned}
        \mathpzc{B}_0\pare{\eta}\bra{\nabla\eta\cdot\Omega}
    &=
    \OpBW{
        \ii A_{0}^{\mathpzc{B}_0,\omega}\pare{\eta;x,\xi}
        +A_{\bra{-2,\geq1}}\pare{\eta;x,\xi}
    }
    \bra{\nabla\eta\cdot\Omega}
    \\
    &+
    \OpBW{
        A^{\mathpzc{B}_0,\eta}_1\pare{\eta;x,\xi}
        \pare{\nabla\eta\cdot\Omega}
        +A_{\bra{0,\geq1}}\bra{\Omega}\pare{\eta;x,\xi}
    }\eta
    +
    \mathfrak R^{\lin,\bR}_{\bra{\varrho,\geq1}}\pare{\eta,\Omega}
    \vect{\eta}{\Omega} .
    \end{aligned}
\end{equation}
Here the symbol \(A_{\bra{0,\geq1}}\bra{\Omega}\) is linear in \(\Omega\), with
coefficients depending on \(\eta\), because the factor \(\nabla\eta\)
in \(\omega_\eta\) is absorbed in the coefficient.

We next expand the first term in \cref{eq:toughcal1}.  By
\cref{lem:paraprod_Bony},
\[
    \nabla\eta\cdot\Omega
    =
    \OpBW{\pare{\nabla\eta}^{\intercal}}\Omega
    +
    \OpBW{\Omega^{\intercal}}\nabla\eta
    +
    \sum_{j=1}^2\RBW{\partial_j\eta}{\Omega_j},
\]
and \cref{rem:smoothBONYremainder} together with
\cref{lem:composition_rules}\,\textup{(i)} shows that the contribution of
\(\sum_{j=1}^2\RBW{\partial_j\eta}{\Omega_j}\) is absorbed in
\(\mathfrak R^{\lin,\bR}_{\bra{\varrho,\geq1}}\pare{\eta,\Omega}\vect{\eta}{\Omega}\).
For the first paraproduct we use \cref{prop:composition_paradiff} with
asymptotic expansion up to order one.  Since
\(A_{0}^{\mathpzc{B}_0,\omega}\in\pare{\Gamma^0_{\geq0}}^2\) and
\(\pare{\nabla\eta}^{\intercal}\in\pare{\Gamma_{\geq1}^0}^{1\times2}\),
\[
    \pare{\ii A_{0}^{\mathpzc{B}_0,\omega}}\#_1\pare{\nabla\eta}^{\intercal}
    =
    \ii A_{0}^{\mathpzc{B}_0,\omega}\pare{\nabla\eta}^{\intercal}
    +
    \frac{1}{2}\pbra{A_{0}^{\mathpzc{B}_0,\omega}}{\pare{\nabla\eta}^{\intercal}},
\]
where the second term belongs to \(\pare{\Gamma_{\geq1}^{-1}}^{2\times2}\).
The composition defect belongs to \(\cR_{\geq1}^{\varrho}\) after increasing
\(\varrho_0\), while the contribution of
\(A_{\bra{-2,\geq1}}\) gives a symbol in \(\pare{\Gamma_{\geq1}^{-2}}^{2\times2}\).

For the second paraproduct, we first write
\(
    \nabla\eta=\OpBW{\ii\xi}\eta
\).
Another application of \cref{prop:composition_paradiff} gives
\[
    \OpBW{\Omega^{\intercal}}\nabla\eta
    =
    \OpBW{\ii\,\Omega\cdot\xi
        +A_{\bra{0,\geq1}}\bra{\Omega}\pare{\eta;x,\xi}}
    \eta
    +
    R^\bR_{\bra{\varrho,\geq1}}\pare{\eta,\Omega}\eta .
\]
Composing once more with \(\OpBW{\ii A_{0}^{\mathpzc{B}_0,\omega}+A_{\bra{-2,\geq1}}}\), and using
\cref{prop:composition_paradiff,lem:composition_rules}, the principal
part is
\(
    \ii A_{0}^{\mathpzc{B}_0,\omega}\,\ii\pare{\Omega\cdot\xi}
    =-A_{0}^{\mathpzc{B}_0,\omega}\pare{\Omega\cdot\xi}
\), while all Poisson brackets and the terms generated by
\(A_{\bra{-2,\geq1}}\) are absorbed in
\(A_{\bra{0,\geq1}}\bra{\Omega}\pare{\eta;x,\xi}\).  Therefore,
\begin{equation}\label{eq:toughcal5.1}
    \begin{aligned}
    &\OpBW{\ii A_{0}^{\mathpzc{B}_0,\omega}\pare{\eta;x,\xi}+A_{\bra{-2,\geq1}}\pare{\eta;x,\xi}}
    \bra{\nabla\eta\cdot\Omega}\\
    =&
    \OpBW{
        \ii A_{0}^{\mathpzc{B}_0,\omega}\pare{\eta;x,\xi}\pare{\nabla\eta}^{\intercal}
        +
        \frac12\pbra{A_{0}^{\mathpzc{B}_0,\omega}\pare{\eta;x,\xi}}{\pare{\nabla\eta}^{\intercal}}
        +
        {\bm A}_{\bra{-2,\geq1}}\pare{\eta;x,\xi}
    }\Omega\\
    +&
    \OpBW{
        -A_{0}^{\mathpzc{B}_0,\omega}\pare{\eta;x,\xi}\pare{\Omega\cdot\xi}
        +
        A_{\bra{0,\geq1}}\bra{\Omega}\pare{\eta;x,\xi}
    }\eta
    +
    \mathfrak R^{\lin,\bR}_{\bra{\varrho,\geq1}}\pare{\eta,\Omega}
    \vect{\eta}{\Omega} .
    \end{aligned}
\end{equation}
Inserting \cref{eq:toughcal5.1} into \cref{eq:toughcal1} yields
\begin{equation}\label{eq:toughcal6}
\begin{aligned}
    \mathpzc{B}_0\pare{\eta}\bra{\nabla\eta\cdot\Omega}
    =&
    \OpBW{
        \ii A_{0}^{\mathpzc{B}_0,\omega}\pare{\eta;x,\xi}\otimes\nabla\eta
        +
        \frac12\pbra{A_{0}^{\mathpzc{B}_0,\omega}\pare{\eta;x,\xi}}{\pare{\nabla\eta}^{\intercal}}
        +
        {\bm A}_{\bra{-2,\geq1}}\pare{\eta;x,\xi}
    }\Omega
    \\
    +&
    \OpBW{
        \pare{
            -A_{0}^{\mathpzc{B}_0,\omega}\pare{\eta;x,\xi}\otimes\xi
            +
            A^{\mathpzc{B}_0,\eta}_1\pare{\eta;x,\xi}
            \otimes\nabla\eta
        }\Omega
        +
        A_{\bra{0,\geq1}}\bra{\Omega}\pare{\eta;x,\xi}
    }\eta
    +
    \mathfrak R^{\lin,\bR}_{\bra{\varrho,\geq1}}\pare{\eta,\Omega}
    \vect{\eta}{\Omega} .
\end{aligned}
\end{equation}
We now treat \(\mathpzc{q}\pare{\eta}\Omega\).  Applying
\cref{prop:paralinearization_Q} to the two scalar components of
\(\Omega\), and collecting them into a vector identity, gives
\begin{equation}\label{eq:toughcal7}
\begin{aligned}
    \mathpzc{q}\pare{\eta}\Omega
    =&
    \OpBW{
        \ii a^{\mathpzc{q},\omega}_0\pare{\eta;x,\xi}\Id_{\bR^2}
        +
        {\bm A}_{\bra{-2,\geq1}}\pare{\eta;x,\xi}
    }\Omega
    \\
    +&
    \OpBW{
        a^{\mathpzc{q},\eta}_1\pare{\eta;x,\xi}\Omega
        +
        A_{\bra{0,\geq1}}\bra{\Omega}\pare{\eta;x,\xi}
    }\eta
    +
    \mathfrak R^{\lin,\bR}_{\bra{\varrho,\geq1}}\pare{\eta,\Omega}
    \vect{\eta}{\Omega} .
\end{aligned}
\end{equation}
Here the lower symbol in the \(\eta\)-column is, more precisely,
\[
    A_{\bra{0,\geq1}}\bra\cdot
    \in
    \cL\pare{\pare{H^{s_0}}^2;\pare{\Gamma^0}^2},
\]
and the corresponding smoothing column belongs to
\(\cL\pare{\pare{H^{s_0}}^2;\pare{\cR^\varrho}^2}\), by
\cref{eq:q-lower-symbol-classes,eq:q-smoothing-classes}. Although these
coefficients need not vanish at \(\eta=0\), they multiply \(\eta\); hence the
complete terms have total degree at least one in \((\eta,\Omega)\).
By the definition \cref{eq:bmvtilde_recall},
\(
     V^{\star}\bra{\Omega}\pare{\eta}
    =
    {\bm J}
    \pare{
        \mathpzc{B}_0\pare{\eta}\bra{\nabla\eta\cdot\Omega}
        -
        \mathpzc{q}\pare{\eta}\Omega
    }
\).
Combining \cref{eq:toughcal6} and \cref{eq:toughcal7}, we obtain
\cref{eq:paralinearization_Vtilde}, with the intermediate symbols
\begin{subequations}
\begin{align}
\label{eq:AVtildeg0}
    {\bm A}^{ V^{\star},\Omega}_0\pare{\eta;x,\xi}
    &\defeq
    {\bm J}
    \pare{
        A_{0}^{\mathpzc{B}_0,\omega}\pare{\eta;x,\xi}\otimes\nabla\eta
        -
        a^{\mathpzc{q},\omega}_0\pare{\eta;x,\xi}\Id_{\bR^2}
    },
    \\
\label{eq:AVtildeg-1}
    {\bm A}^{ V^{\star},\Omega}_{-1}\pare{\eta;x,\xi}
    &\defeq
    {\bm J}
    \pare{
        \frac12
        \pbra{A_{0}^{\mathpzc{B}_0,\omega}\pare{\eta;x,\xi}}{\pare{\nabla\eta}^{\intercal}}
    },
    \\
\label{eq:AVtildeeta1}
    {\bm A}^{ V^{\star},\eta}_1\pare{\eta;x,\xi}
    &\defeq
    {\bm J}
    \pare{
        -A_{0}^{\mathpzc{B}_0,\omega}\pare{\eta;x,\xi}\otimes\xi
        +
        A^{\mathpzc{B}_0,\eta}_1\pare{\eta;x,\xi}\otimes\nabla\eta
        -
        a^{\mathpzc{q},\eta}_1\pare{\eta;x,\xi}\Id_{\bR^2}
    } .
\end{align}
\end{subequations}
The membership of these symbols in the classes stated in the proposition follows
from
\[
A_{0}^{\mathpzc{B}_0,\omega}\in\pare{\Gamma^0_{\geq0}}^2,\quad
A^{\mathpzc{B}_0,\eta}_1\in\pare{\Gamma_{\geq1}^1}^2,\quad
    a^{\mathpzc{q},\omega}_0\in\Gamma_{\geq1}^0,\quad
    a^{\mathpzc{q},\eta}_1\in\Gamma^1_{\geq0}
,\]
from \cref{prop:B0_paralinearization,prop:paralinearization_Q}, and from
\cref{prop:composition_paradiff}.  It remains only to identify the three
intermediate expressions with the explicit symbols in the statement.

\begin{proofsubpart}[Identification of \cref{eq:AVtildeg0_lem}]
For any \(W_1,W_2\in\bR^2\),
\begin{equation}\label{eq:tensor_relation}
    W_1^\perp\otimes W_2^\perp+W_2\otimes W_1
    =
    \pare{W_1\cdot W_2}\Id_{\bR^2}.
\end{equation}
Using \cref{eq:AB0omega_AQL0}, namely
\(
    a^{\mathpzc{q},\omega}_0=\nabla\eta\cdot A_{0}^{\mathpzc{B}_0,\omega}
\), and \cref{eq:tensor_relation} with
\((W_1,W_2)=(\nabla\eta,A_{0}^{\mathpzc{B}_0,\omega})\), we get
\begin{equation}\label{eq:pi1_identity_1}
    A_{0}^{\mathpzc{B}_0,\omega}\otimes\nabla\eta
    -
    a^{\mathpzc{q},\omega}_0\Id_{\bR^2}
    =
    -\nabla^\perp\eta\otimes A_{0}^{\mathpzc{B}_0,\omega;\perp} .
\end{equation}
Moreover,
\begin{equation}\label{eq:pi1_identity_2}
    -{\bm J}
    \pare{W_1^\perp\otimes W_2^\perp}
    =
    W_1\otimes W_2^\perp .
\end{equation}
Inserting \cref{eq:pi1_identity_1} into \cref{eq:AVtildeg0} and using
\cref{eq:pi1_identity_2}, we obtain
\[
    {\bm A}^{ V^{\star},\Omega}_0\pare{\eta;x,\xi}
    =
    \nabla\eta\otimes A_{0}^{\mathpzc{B}_0,\omega;\perp}\pare{\eta;x,\xi},
\]
which is \cref{eq:AVtildeg0_lem}.
\end{proofsubpart}

\begin{proofsubpart}[Identification of \cref{eq:AVtildeg-1_lem}]
Since \(\pare{\nabla\eta}^{\intercal}\) is independent of \(\xi\), the
matrix Poisson bracket of \cref{eq:Pbra_matrix-vector} gives
\[
    \pbra{A_{0}^{\mathpzc{B}_0,\omega}}{\pare{\nabla\eta}^{\intercal}}
    =
    J_\xi A_{0}^{\mathpzc{B}_0,\omega}\,{\bm H}_x\eta .
\]
Thus \cref{eq:AVtildeg-1} becomes
\[
    {\bm A}^{ V^{\star},\Omega}_{-1}\pare{\eta;x,\xi}
    =
    \frac{{\bm J}}{2}
    J_\xi A_{0}^{\mathpzc{B}_0,\omega}\,{\bm H}_x\eta,
\]
which is \cref{eq:AVtildeg-1_lem}.  The symbol belongs to
\(\pare{\Gamma_{\geq1}^{-1}}^{2\times2}\), because one \(\xi\)-derivative falls on
\(A_{0}^{\mathpzc{B}_0,\omega}\in\pare{\Gamma^0_{\geq0}}^2\).
\end{proofsubpart}

\begin{proofsubpart}[Identification of \cref{eq:AVtildeeta1_lem}]
By \cref{eq:relationa1A1},
\[
    a^{\mathpzc{q},\eta}_1
    =
    \nabla\eta\cdot A^{\mathpzc{B}_0,\eta}_1
    -
    A_{0}^{\mathpzc{B}_0,\omega}\cdot\xi .
\]
Substituting this identity in \cref{eq:AVtildeeta1} and applying
\cref{eq:tensor_relation} twice, first with \((W_1,W_2)=(\xi,A_{0}^{\mathpzc{B}_0,\omega})\)
and then with
\((W_1,W_2)=(\nabla\eta,A^{\mathpzc{B}_0,\eta}_1)\), gives
\[
    -A_{0}^{\mathpzc{B}_0,\omega}\otimes\xi
    +
    \pare{A_{0}^{\mathpzc{B}_0,\omega}\cdot\xi}\Id_{\bR^2}
    =
    \xi^\perp\otimes A_{0}^{\mathpzc{B}_0,\omega;\perp},
\]
and
\[
    A^{\mathpzc{B}_0,\eta}_1\otimes\nabla\eta
    -
    \pare{\nabla\eta\cdot A^{\mathpzc{B}_0,\eta}_1}\Id_{\bR^2}
    =
    -\nabla^\perp\eta
    \otimes
    A^{\mathpzc{B}_0,\eta;\perp}_1.
\]
Therefore
\[
    {\bm A}^{ V^{\star},\eta}_1
    =
    {\bm J}
    \pare{
        \xi^\perp\otimes A_{0}^{\mathpzc{B}_0,\omega;\perp}
        -
        \nabla^\perp\eta\otimes
        A^{\mathpzc{B}_0,\eta;\perp}_1
    } .
\]
Using \cref{eq:pi1_identity_2} on the two tensor products, we conclude
\[
    {\bm A}^{ V^{\star},\eta}_1\pare{\eta;x,\xi}
    =
    -\xi\otimes A_{0}^{\mathpzc{B}_0,\omega;\perp}\pare{\eta;x,\xi}
    +
    \nabla\eta\otimes
    A^{\mathpzc{B}_0,\eta;\perp}_1\pare{\eta;x,\xi},
\]
which is \cref{eq:AVtildeeta1_lem}.
\end{proofsubpart}
\end{proofpart}

\begin{proofpart}[Paralinearization of
\texorpdfstring{$\Omega\mapsto v_3^{\star}\bra{\Omega}\pare{\eta}$}{Omega -> v3tilde}]
Write
\(
    \mathpzc{B}_0\pare{\eta}=(\mathpzc{B}_{0,1}\pare{\eta},
    \mathpzc{B}_{0,2}\pare{\eta})
\).  From \cref{eq:bmvtilde_recall},
\begin{equation}\label{eq:tildev3_component_identity}
     v_3^{\star}\bra{\Omega}\pare{\eta}
    =
    -\mathpzc{B}_0\pare{\eta}^\perp\cdot\Omega
    =
    \mathpzc{B}_{0,2}\pare{\eta}\bra{\Omega_1}
    -
    \mathpzc{B}_{0,1}\pare{\eta}\bra{\Omega_2} .
\end{equation}
We apply \cref{prop:B0_paralinearization} to the two scalar inputs
\(\Omega_1\) and \(\Omega_2\), and then take the combination in
\cref{eq:tildev3_component_identity}.  The \(\Omega\)-principal part is
\[
    \OpBW{\ii A_{0,2}^{\mathpzc{B}_0,\omega}}\Omega_1
    -
    \OpBW{\ii A_{0,1}^{\mathpzc{B}_0,\omega}}\Omega_2
    =
    \OpBW{-\ii\pare{A_{0}^{\mathpzc{B}_0,\omega;\perp}}^{\intercal}}\Omega,
\]
which corresponds to
\(
    A^{ v_3^{\star},\Omega}_0=-A_{0}^{\mathpzc{B}_0,\omega;\perp}\in\pare{\Gamma_{\geq0}^0}^2
\).
Likewise, the \(\eta\)-principal part is
\[
    \OpBW{
        A^{\mathpzc{B}_0,\eta}_{1,2}\pare{\eta;x,\xi}\Omega_1
        -
        A^{\mathpzc{B}_0,\eta}_{1,1}\pare{\eta;x,\xi}\Omega_2
    }\eta
    =
    \OpBW{
        -A^{\mathpzc{B}_0,\eta;\perp}_1\pare{\eta;x,\xi}
        \cdot\Omega
    }\eta .
\]
The remaining terms in \cref{prop:B0_paralinearization} are either
\(a_{\bra{0,\geq1}}\bra{\Omega}\pare{\eta;x,\xi}\eta\), smoothing remainders in
\(\mathfrak r^{\lin,\bR}_{\bra{\varrho,\geq1}}\pare{\eta,\Omega}\cdot\vect{\eta}{\Omega}\),
or symbols of order \(-2\) acting on \(\Omega\); the latter are absorbed in
the generic \(A_{\bra{-2,\geq1}}\pare{\eta;x,\xi}\)-slot of
\cref{eq:paralinearization_tildev3}.  Hence
\[
     v_3^{\star}\bra{\Omega}\pare{\eta}
    =
    \OpBW{
        \vect{
            -A^{\mathpzc{B}_0,\eta;\perp}_1\cdot\Omega
        }{
            -\ii A_{0}^{\mathpzc{B}_0,\omega;\perp}
        }
        +
        \vect{
            a_{\bra{0,\geq1}}\bra{\Omega}
        }{
            A_{\bra{-2,\geq1}}
        }
    }\cdot\vect{\eta}{\Omega}
    +
    {\bm r}^{\lin,\bR}_{\bra{\varrho,\geq1}}\pare{\eta,\Omega}
    \cdot\vect{\eta}{\Omega},
\]
which is precisely \cref{eq:paralinearization_tildev3}.  The preceding
identities also prove
\(
    a_1^{ v_3^{\star},\eta}\bra{\Omega}\in\Gamma_{\geq1}^1
\)
and
\(
    A^{ v_3^{\star},\Omega}_0\in\pare{\Gamma_{\geq0}^0}^2
\), as declared.
\end{proofpart}
\end{proof}
We record the only algebraic identity from the preceding computation that is used
below in the paralinearization of the pressure contribution.

\begin{cor}[Quadratic identity for the \(\eta\)-symbol of \( V^{\star}\)]
\label{cor:simplifyimportant}
For any \(x\in\bR^2, \xi\in\bR^2\setminus\set{0}, W\in\bR^2\), one has
\begin{equation}\label{eq:relationimportant5}
\pare{
    {\bm A}^{ V^{\star},\eta}_1\pare{\eta;x,\xi}W
}
\cdot W^\perp
=
-\lambda_1\pare{\eta;x,\xi}
\left[
    \pare{2A_{0}^{\mathpzc{B}_0,\omega;\perp}\pare{\eta;x,\xi}\cdot W}^2
    -
    \pare{
        \frac{
            \nabla^\perp\eta\pare{x}\cdot W
        }{
            1+\av{\nabla\eta\pare{x}}^2
        }
    }^2
\right].
\end{equation}
\end{cor}

\begin{proof}
Both sides of \cref{eq:relationimportant5} are one-homogeneous in \(\xi\)
and invariant under simultaneous rotations of \(\xi\), \(\nabla\eta\pare{x}\),
and \(W\). It is therefore enough to prove the identity for
\(\av{\xi}=1\) and \(\xi=e_1\). Write
\[
    \nabla\eta\pare{x}=\vect{p}{q},
    \qquad
    W=\vect{\alpha}{\beta}.
\]
Then
\[
    \angles{\nabla\eta}^2=1+p^2+q^2,
    \qquad
    1+q^2=4\lambda_1^2\pare{\eta;x,e_1},
\]
and
\[
    A_{0}^{\mathpzc{B}_0,\omega}\pare{\eta;x,e_1}
    =
    -\frac{1}{4\lambda_1\pare{\eta;x,e_1}
    \angles{\nabla\eta}^2}
    \vect{1+q^2}{-pq}.
\]
Hence
\[
    A_{0}^{\mathpzc{B}_0,\omega;\perp}\pare{\eta;x,e_1}\cdot W
    =
    -\frac{
        pq\alpha+\pare{1+q^2}\beta
    }{
        4\lambda_1\pare{\eta;x,e_1}
        \angles{\nabla\eta}^2
    }\cdot
\]
Moreover,
\[
    \nabla^\perp\eta\pare{x}\cdot W
    =
    -q\alpha+p\beta,
    \qquad
    e_1\cdot W^\perp=-\beta,
    \qquad
    \nabla\eta\pare{x}\cdot W^\perp
    =
    q\alpha-p\beta .
\]
Using \cref{eq:AVtildeeta1_lem,eq:relationimportant4}, we get
\[
\begin{aligned}
\pare{
    {\bm A}^{ V^{\star},\eta}_1\pare{\eta;x,e_1}W
}
\cdot W^\perp
&=
\beta\,A_{0}^{\mathpzc{B}_0,\omega;\perp}\pare{\eta;x,e_1}\cdot W
\\
&\quad
+
\pare{q\alpha-p\beta}
\left[
    \frac{p}{\angles{\nabla\eta}^2}
    A_{0}^{\mathpzc{B}_0,\omega;\perp}\pare{\eta;x,e_1}\cdot W
    -
    \frac{\lambda_1\pare{\eta;x,e_1}}
         {\angles{\nabla\eta\pare{x}}^4}
    \pare{-q\alpha+p\beta}
\right].
\end{aligned}
\]
Since
\[
    \beta
    +
    \frac{p\pare{q\alpha-p\beta}}
         {\angles{\nabla\eta}^2}
    =
    \frac{
        pq\alpha+\pare{1+q^2}\beta
    }{
        \angles{\nabla\eta}^2
    },
\]
we infer
\[
\begin{aligned}
\pare{
    {\bm A}^{ V^{\star},\eta}_1\pare{\eta;x,e_1}W
}
\cdot W^\perp
&=
-\frac{
    \pare{pq\alpha+\pare{1+q^2}\beta}^2
}{
    4\lambda_1\pare{\eta;x,e_1}
    \angles{\nabla\eta\pare{x}}^4
}
+
\frac{\lambda_1\pare{\eta;x,e_1}}
     {\angles{\nabla\eta\pare{x}}^4}
\pare{-q\alpha+p\beta}^2
\\
&=
-\lambda_1\pare{\eta;x,e_1}
\left[
    \pare{
        2A_{0}^{\mathpzc{B}_0,\omega;\perp}\pare{\eta;x,e_1}\cdot W
    }^2
    -
    \pare{
        \frac{
            \nabla^\perp\eta\pare{x}\cdot W
        }{
            \angles{\nabla\eta}^2
        }
    }^2
\right].
\end{aligned}
\]
This proves \cref{eq:relationimportant5}.
\end{proof}

\subsubsection{\texorpdfstring{Paralinearization of operators $T\pare{\eta}$ and $\tau_0\pare{\eta}$}{Paralinearization of operators T(η) and τ₀(η)}}\label{Subsec:paraTtau0}

 Consider the operator $T\p\eta$ of nonlinear convolution type, that has been introduced in \cref{eq:T(eta)}, and that can be recast as follows using \cref{eq:bmvtilde}:
\begin{equation}
\label{eq:T(eta)_recall}
T\p\eta\bra{\Omega} = -\pare{ V^{\star}\bra{\Omega}\pare{\eta} +  v_3^{\star} \bra{\Omega}\pare{\eta}\nabla\eta}^\perp .
\end{equation}
The paralinearization of the operator $T\pare{\eta}$ is given in the following proposition.

\begin{prop}[Paralinearization of $T\pare{\eta}$]\label{prop:T_paralinearization}
There exists $\varrho_0 > 0$ such that for any $\varrho \geq \varrho_0$, there exists $s_0 > 0$ and $r_0 > 0$ such that:
  For any $\eta\in \BallR{s_0}{r_0}, \Omega \in\pare{\BallR{s_0}{r_0}}^2$, the operator $T\pare{\eta}\bra{ \Omega}$ as in \cref{eq:T(eta)_recall} admits the following paralinearization:
        \begin{equation}
        \label{eq:paraoperatorT}
                T\pare{\eta} \Omega=
                \OpBW{
                \mathfrak{A}^T_1\pare{\eta,\Omega;x,\xi}
                +
                \vectT{A_{\bra{0,\geq1}}\bra{\Omega}\pare{\eta;x, \xi}}{{\bm A}_{\bra{-2,\geq1}}\pare{\eta;x,\xi}}
                }\vect{\eta}{\Omega}
                 +\mathfrak{R}^{\lin,\bR}_{\bra{\varrho,\geq1}}\pare{\eta, \Omega} \vect{\eta}{\Omega},
        \end{equation}
        where
        \begin{itemize}
       \item $ \mathfrak{A}^T_1\in \pare{\Gamma_{\geq1}^1}^{2\times 3}$ and is defined as
       \begin{equation}
           \mathfrak{A}^T_1\pare{\eta,\Omega;x,\xi} \defeq \vectT{A^{T,\eta}_1\bra{\Omega}\pare{\eta;x,\xi}}{{\bm A}^{T,\Omega}_{-1}\pare{\eta;x, \xi}},
       \end{equation}
       where
       \begin{itemize}

        \item $A^{T,\eta}_1 \in \cL\pare{ \pare{H^{s_0}}^2 ; \pare{\Gamma_{\geq0}^1}^2}$ is defined as
                \begin{equation}
                \label{eq:AT1extra}
                 A^{T,\eta}_1\bra{\Omega}\pare{\eta;x, \xi}\defeq    \pare{A_{0}^{\mathpzc{B}_0,\omega;\perp}\pare{\eta;x,\xi}\cdot \Omega - \ii \  v_3^{\star}\bra{\Omega}\pare{\eta;x}}  \xi^\perp,
                \end{equation}
                where $ v_3^{\star}$ is defined in \cref{eq:bmvtilde_recall} and $A_{0}^{\mathpzc{B}_0,\omega}$ is defined in \cref{eq:H0_definition};
           
        \item $ {\bm A}^{T,\Omega}_{-1}\in \pare{\Gamma_{\geq1}^{-1}}^{2\times 2}$ and has the explicit form
           \begin{equation}\label{eq:ATOmega-1}
           \begin{aligned}
           {\bm A}^{T,\Omega}_{-1}\pare{\eta;x,\xi}
           &\defeq{}
       \frac{1}{2} \pare{\tr \pare{{\bm H}_\xi \phi^{\mathpzc{B}_0}_1 {\bm J} {\bm H}_x\eta}  \ {\bm J}
        + \Delta_\xi \phi^{\mathpzc{B}_0}_1 \,{\bm H}_x\eta
   + \Delta_x\eta \,{\bm H}_\xi \phi^{\mathpzc{B}_0}_1
   - \Delta_\xi \phi^{\mathpzc{B}_0}_1 \Delta_x\eta \,\Id 
   };
           \end{aligned}
        \end{equation}

       \end{itemize}

             \item ${\bm A}_{\bra{m,\geq1}}\in\pare{\Gamma_{\geq1}^m}^{2\times2}, A_{\bra{m,\geq1}}\in\pare{\Gamma_{\geq1}^m}^2,$ $m\in\bR$, denote generic symbols as per \cref{notation:paralinearization};
             \item $\mathfrak{R}^{\lin,\bR}_{\bra{\varrho,\geq1}}\in\pare{\cR_{\geq1}^\varrho}^{2\times3}$ denotes generic smoothing remainders as per \cref{notation:paralinearization}.
        \end{itemize}
\end{prop}

\begin{proof}

    Since the paralinearization for $ V^{\star}\bra{\Omega}\pare{\eta}$ is given in \cref{prop:paralinearization_Vtilde} we  turn our attention to the term
\begin{equation}
    - v_3^{\star}\bra{\Omega}\pare{\eta}\nabla^\perp\eta = \pare{\mathpzc{B}_0\p\eta^\perp\cdot \Omega}\nabla^\perp\eta
\end{equation}
in \cref{eq:T(eta)_recall}. \cref{lem:paraprod_Bony} gives
\begin{equation}
\label{eq:paralinearization_tildev3.1}
     \pare{- v_3^{\star}\bra{\Omega}\pare{\eta}}\nabla^\perp\eta = \OpBW{- v_3^{\star}\bra{\Omega}\pare{\eta}}\nabla^\perp\eta + \OpBW{\nabla^\perp\eta}\bra{- v_3^{\star}\bra{\Omega}\pare{\eta}} + \RBW{\nabla^\perp\eta}{- v_3^{\star}\bra{\Omega}\pare{\eta}}. 
\end{equation}
Using the paralinearization in \cref{eq:paralinearization_tildev3}, combined with \cref{rem:functions_Op_R} we obtain that 
\begin{equation}
    \label{eq:B0cdotg_function}
    - v_3^{\star}\bra{\Omega}\pare{\eta}\in \cF^\bR\subset \Gamma_{\geq1}^0.
\end{equation}
So we can apply \cref{prop:composition_paradiff} and obtain that
\begin{equation}
\label{eq:paralinearization_tildev3.2}
     \OpBW{- v_3^{\star}\bra{\Omega}\pare{\eta}}\nabla^\perp\eta =  \OpBW{-\ii \  v_3^{\star}\bra{\Omega}\pare{\eta} \xi^\perp + A_{\bra{0,\geq1}}\bra{\Omega}\pare{\eta;x, \xi}  }\eta + R^\bR_{\bra{\varrho,\geq1}}\bra{\Omega}\pare{\eta}  \eta . 
\end{equation}
Next, we invoke \cref{eq:paralinearization_tildev3,lem:composition_rules,prop:composition_paradiff,eq:Pbra_matrix-vector} to get 
    \begin{equation}
        \label{eq:paralinearization_tildev3.3}
            \begin{aligned}
                    \OpBW{\nabla^\perp\eta}\bra{- v_3^{\star}\bra{\Omega}\pare{\eta}}
                =
                & \ 
                 \OpBW{  \ii \ \nabla^\perp\eta \otimes A_{0}^{\mathpzc{B}_0,\omega;\perp}\pare{\eta;x,\xi} 
                 + \frac{1}{2}\pbra{\nabla^\perp\eta}{\pare{A_{0}^{\mathpzc{B}_0,\omega;\perp}\pare{\eta;x,\xi}}^\intercal}+  {\bm A}_{\bra{-2,\geq1}}\pare{\eta;x, \xi} }  \Omega
                \\
                 +& \ 
                 \OpBW{\nabla^\perp\eta \pare{A_{1}^{\mathpzc{B}_0,\eta;\perp}\pare{\eta;x,\xi}\cdot \Omega} + { A}_{\bra{0}}\bra{\Omega}\pare{\eta;x, \xi} }\eta
                +\mathfrak{R}^{\lin,\bR}_{\bra{\varrho,\geq1}}\pare{\eta, \Omega}\vect{\eta}{\Omega}. 
            \end{aligned}
    \end{equation}
An application of \cref{rem:smoothBONYremainder,lem:composition_rules,eq:B0cdotg_function} shows
\begin{equation}
\label{eq:paralinearization_tildev3.4}
    \RBW{\nabla^\perp\eta}{- v_3^{\star}\bra{\Omega}\pare{\eta}} = R^\bR_{\bra{\varrho,\geq1}}\pare{\eta, \Omega} \eta.
\end{equation}
Plugging \cref{eq:paralinearization_tildev3.4,eq:paralinearization_tildev3.3,eq:paralinearization_tildev3.2} into \cref{eq:paralinearization_tildev3.1}, we obtain that
\begin{equation}\label{eq:toughcal8}
         \begin{aligned}
             - v_3^{\star}\bra{\Omega}\pare{\eta} \ \nabla^\perp\eta
        =&
         \OpBW{ -\ii \  v_3^{\star}\bra{\Omega}\pare{\eta} \xi^\perp + \nabla^\perp\eta \pare{A_{1}^{\mathpzc{B}_0,\eta;\perp}\pare{\eta;x,\xi}\cdot \Omega} + A_{\bra{0,\geq1}}\bra{\Omega}\pare{\eta;x, \xi}  }\eta
         \\
          +&  \OpBW{ \ii \   \nabla^\perp\eta \otimes A_{0}^{\mathpzc{B}_0,\omega;\perp}\pare{\eta;x,\xi} 
          +  \frac{1}{2}\pbra{\nabla^\perp\eta}{\pare{A_{0}^{\mathpzc{B}_0,\omega;\perp}\pare{\eta;x,\xi}}^\intercal}+  {\bm A}_{\bra{-2,\geq1}}\pare{\eta;x, \xi} }  \Omega\\+&\mathfrak{R}^{\lin,\bR}_{\bra{\varrho,\geq1}}\pare{\eta, \Omega}\vect{\eta}{\Omega} .
         \end{aligned}
\end{equation}
Now, we put \cref{eq:paralinearization_Vtilde,eq:toughcal8} into \cref{eq:T(eta)_recall} and get
\begin{equation}\label{eq:paralinearization_Tetag_1}
    \begin{aligned}
        T\pare{\eta}\bra{\Omega}  &= \OpBW{ {\bm A}^{T,\Omega}_{-1}\pare{\eta;x, \xi} + {\bm A}_{\bra{-2,\geq1}}\pare{\eta;x, \xi} }\Omega\\
    &+
    \OpBW{ A^{T,\eta}_1\bra{\Omega}\pare{\eta;x, \xi} + A_{\bra{0,\geq1}} \pare{\eta, \Omega;x, \xi} }\eta
   +{\bm R}^\bR_{\bra{\rho,\geq1}}\p\eta \Omega+R^{\lin,\bR}_{\bra{\rho,\geq1}}\pare{\eta,\Omega}\eta,
    \end{aligned}
\end{equation}
where
\begin{align}
\label{eq:ATQL-1_1}
{\bm A}^{T,\Omega}_{-1}\pare{\eta;x, \xi}
\defeq \ &\frac{1}{2}
\pbra{\nabla^\perp\eta}
{\pare{A_{0}^{\mathpzc{B}_0,\omega;\perp}\pare{\eta;x,\xi}}^\intercal}
-
{\bm J}\,{\bm A}^{ V^{\star},\Omega}_{-1}\pare{\eta;x, \xi}
=
\frac{1}{2}
\pare{
{\bm H}_\xi \phi^{\mathpzc{B}_0}_1 \  {\bm H}_x\eta
- {\bm J} {\bm H}_x\eta \  \pare{{\bm J} {\bm H}_\xi \phi^{\mathpzc{B}_0}_1}^\intercal
}
\end{align}
and 
\begin{equation}\label{eq:AT1}
    \begin{aligned}
         A^{T,\eta}_1\bra{\Omega}\pare{\eta ;x, \xi}
         \defeq & \ 
        - \ii \  v_3^{\star}\bra{\Omega}\pare{\eta} \  \xi^\perp 
        +  \nabla^\perp\eta \pare{A^{\mathpzc{B}_0,\eta;\perp}_1\pare{\eta;x,\xi}\cdot \Omega}\\&-{\bm J} \pare{-\xi\otimes A_{0}^{\mathpzc{B}_0,\omega;\perp} \pare{\eta,x,\xi} + \nabla\eta \otimes A^{\mathpzc{B}_0, \eta;\perp}_1\pare{\eta;x,\xi}}\Omega . 
    \end{aligned}
\end{equation}
We now focus in simplifying the expressions in \cref{eq:ATQL-1_1,eq:AT1}.

\begin{step}[Simplification of \cref{eq:ATQL-1_1}]
   We use the identity
   \[
   {\bm J} {\bm A} {\bm J} = {\bm A}^\intercal - \tr\pare{{\bm A}}\,\Id,
   \]
   valid for any $2\times2$ matrix ${\bm A}$, so that if ${\bm S}_1, {\bm S}_2$ are symmetric
   \begin{equation}
       \label{eq:symm_matrices}
       {\bm S}_1 {\bm S}_2 - {\bm J}{\bm S}_2\pare{{\bm J}{\bm S}_1}^\intercal
       =
       {\bm S}_1 {\bm S}_2 + {\bm J} {\bm S}_2 {\bm S}_1 {\bm J}
       =
       2\,{\bm S}_1 {\bm S}_2 - \tr\pare{{\bm S}_1 {\bm S}_2}\,\Id .
   \end{equation}
   Applying \cref{eq:symm_matrices} we obtain
   \begin{equation}       
   {\bm H}_\xi \phi^{\mathpzc{B}_0}_1 \, {\bm H}_x\eta
   - {\bm J} {\bm H}_x\eta \, \pare{{\bm J} {\bm H}_\xi \phi^{\mathpzc{B}_0}_1}^\intercal
   =
   2\,{\bm H}_\xi \phi^{\mathpzc{B}_0}_1 \, {\bm H}_x\eta
   - \tr\pare{{\bm H}_x\eta \ {\bm H}_\xi \phi^{\mathpzc{B}_0}_1}\,\Id .
   \end{equation}
   Next, decomposing $2{\bm S}_1{\bm S}_2 = \comm{{\bm S}_1}{{\bm S}_2} + ({\bm S}_1{\bm S}_2+{\bm S}_2{\bm S}_1)$ and using that for symmetric $2\times2$ matrices
   \[
   {\bm S}_1{\bm S}_2+{\bm S}_2{\bm S}_1
   =
   \tr\pare{{\bm S}_1}\,{\bm S}_2
   +
    \tr\pare{{\bm S}_2}\,{\bm S}_1
    +
    \tr\pare{{\bm S}_1{\bm S}_2}\,\Id
    -
    \tr\pare{{\bm S}_1}\tr\pare{{\bm S}_2}\,\Id,
   \]
   we obtain the skew-symmetric/symmetric decomposition
   \begin{equation}
   \begin{aligned}
       {\bm H}_\xi \phi^{\mathpzc{B}_0}_1 \, {\bm H}_x\eta
   - {\bm J} {\bm H}_x\eta \, \pare{{\bm J} {\bm H}_\xi \phi^{\mathpzc{B}_0}_1}^\intercal
   &=
   \comm{{\bm H}_\xi \phi^{\mathpzc{B}_0}_1}{{\bm H}_x\eta}\\
   &\quad+ \tr\pare{{\bm H}_\xi \phi^{\mathpzc{B}_0}_1}\,{\bm H}_x\eta
   + \tr\pare{{\bm H}_x\eta}\,{\bm H}_\xi \phi^{\mathpzc{B}_0}_1
   - \tr\pare{{\bm H}_\xi \phi^{\mathpzc{B}_0}_1}\tr\pare{{\bm H}_x\eta}\,\Id .
   \end{aligned}
   \end{equation}
   Next, since both matrices are $2\times 2$ and symmetric
   \begin{equation}
       \comm{{\bm H}_\xi \phi^{\mathpzc{B}_0}_1}{{\bm H}_x\eta}
       =
       \tr\pare{{\bm H}_\xi \phi^{\mathpzc{B}_0}_1 {\bm J} {\bm H}_x\eta}\,{\bm J}.
   \end{equation}
   This yields the desired identity \cref{eq:ATOmega-1}.
\end{step}

\begin{step}[Equivalence of \cref{eq:AT1extra,eq:AT1}]
   Explicit computations using in particular \cref{eq:pi1_identity_2} show that
   \begin{equation}\label{eq:AT1extraextra}
        -{\bm J} \pare{-\xi\otimes A_{0}^{\mathpzc{B}_0,\omega;\perp}  + \nabla\eta \otimes A^{\mathpzc{B}_0, \eta;\perp}_1}\Omega
        =
        \xi^\perp \ A_{0}^{\mathpzc{B}_0,\omega;\perp}\cdot\Omega - \nabla^\perp\eta \  A^{\mathpzc{B}_0, \eta;\perp}_1\cdot\Omega.
   \end{equation}
   Thus, plugging \cref{eq:AT1extraextra} in \cref{eq:AT1} we obtain 
\cref{eq:AT1extra}. 
\end{step}
\end{proof}

Next, we give the paralinearization of the operator $\tau_0\pare{\eta}$ introduced in \cref{eq:Tstar}, which takes the form
\begin{equation}\label{eq:Tstar_recall}
\tau_0 \pare{\eta}\psi
\defeq
\Delta^{-1}\nabla^\perp\cdot
\pare{ T\pare{\eta}\bra{\nabla^\perp \psi}},
\end{equation}
where $T\pare{\eta}$ is the trace operator defined in \cref{eq:T(eta)_recall}.
The paralinearization result reads as follows.

\begin{prop}[Paralinearization of $\tau_0\pare{\eta}$]\label{prop:Tstar_paralinearization}
    There exists $\varrho_0 > 0$ such that for any $\varrho \geq \varrho_0$ there exist $s_0>0$ and $r_0>0$ such that for any $\eta,\psi\in \BallR{s_0}{r_0},$ the operator $\tau_0\pare{\eta}$ defined in \cref{eq:Tstar_recall} admits the paradifferential decomposition
   \begin{equation}
    \tau_0\pare{\eta}\psi
    = \OpBW{ \vect{a^{\tau_0,\eta}_0\bra{\psi}\pare{\eta;x,\xi}}{a^{\tau_0,\psi}_{-1}\pare{\eta;x,\xi}}
     + 
    \vect{a_{\bra{-1,\geq1}}\bra{\psi}\pare{\eta;x,\xi}}{a_{\bra{-2,\geq1}}\pare{\eta;x,\xi}}    }\cdot\vect{\eta}{\psi} 
    +
   R^{\lin,\bR}_{\bra{\varrho,\geq1}}\pare{\eta,\psi}\cdot\vect{\eta}{\psi},
\end{equation}
where 
\begin{itemize}
\item 
$a^{\tau_0,\eta}_0\in\cL\pare{H^{s_0};\Gamma_{\geq0}^0}$ is explicitly defined as
\begin{equation}\label{def a0star}
    a^{\tau_0,\eta}_0\bra{\psi}\pare{\eta;x,\xi}\defeq - v_3^{\star}\bra{\nabla^\perp\psi}\pare{\eta;x}
    -\ii \ \nabla \psi \pare{x} \cdot A_{0}^{\mathpzc{B}_0,\omega}\pare{\eta;x,\xi},
\end{equation}
where $ v_3^{\star}$ is the vertical velocity component defined in \cref{eq:bmvtilde_recall} and $A_{0}^{\mathpzc{B}_0,\omega}\pare{\eta;x,\xi}$ is the paradifferential symbol defined in \cref{eq:H0_definition};
\item 
$a^{\tau_0,\psi}_{-1}\in\Gamma_{\geq1}^{-1}$ is explicitly defined as
\begin{equation}\label{eq:atau0}
a^{\tau_0,\psi}_{-1}\pare{\eta;x,\xi}\defeq  -\frac{\av{\xi}^2}{32\,\lambda_1\pare{\eta;x,\xi}^3}
   \hat\xi^\perp \cdot {\bm H}_x \eta \,\hat\xi^\perp
   =-\frac{\lambda_0\pare{\eta;x,\xi}}{\lambda_1\pare{\eta;x,\xi}},
\end{equation}

As a joint symbol of \((\eta,\psi)\),
\(a^{\tau_0,\eta}_0\bra{\psi}\in\Gamma^0_{\geq1}\) and it is linear in
\(\psi\). For fixed \(\psi\), it belongs to
\(\Gamma^0_{\geq0}\) and need not vanish at \(\eta=0\).

    \item $a_{\bra{m,\geq1}}\in\Gamma_{\geq1}^m$ for $m\in\bR$ are lower-order symbol corrections (cf. \cref{notation:paralinearization}); 
    \item $R^{\lin,\bR}_{\bra{\varrho,\geq1}}\in\pare{\cR_{\geq1}^\varrho}^2$ are remainder terms satisfying the estimates in \cref{notation:paralinearization}. 
\end{itemize} 
\end{prop}

\begin{proof}
Recall $ \Delta^{-1}\nabla^{\perp}\cdot \bullet=\OpBW{-\frac{\ii\xi^\perp}{\av\xi^2}\cdot}\bullet$. Applying \cref{prop:T_paralinearization} we obtain
    \begin{equation}\label{eq:simpolificationsymbolsTstar-1}
        \begin{aligned}
            \Delta^{-1}\nabla^{\perp}\cdot T\pare{\eta}\bra{\nabla^{\perp}\psi} 
        &= \OpBW{\frac{-\ii\xi^{\perp} }{\av{\xi}^2}} \cdot \OpBW{{\bm A}_{-1}^{T,\Omega}\pare{\eta;x,\xi} + {\bm A}_{\bra{-2,\geq1}}\pare{\eta;x,\xi}}\OpBW{\ii \xi^\perp} \psi\\
        &+ \OpBW{\frac{-\ii\xi^{\perp} }{\av{\xi}^2}} \cdot \OpBW{ A^{T,\eta}_1\bra{\nabla^\perp \psi} \pare{\eta ; x, \xi} +A_{\bra{0,\geq1}}\pare{\eta,\nabla^{\perp}\psi;x, \xi}}\eta\\
        &+ \OpBW{\frac{-\ii\xi^{\perp} }{\av{\xi}^2}} \cdot \pare{ {\bm R}^{\lin,\bR}_{\bra{\varrho,\geq1}}\pare{\eta}\OpBW{\ii \xi^\perp}\psi+R^{\lin,\bR}_{\bra{\varrho,\geq1}}\pare{\eta,\nabla^{\perp}\psi}\eta}.
        \end{aligned}
    \end{equation}
  Thus, an application of \cref{prop:composition_paradiff,lem:composition_rules} yields
  \begin{align}
      \Delta^{-1}\nabla^{\perp}\cdot T\pare{\eta}\bra{\nabla^{\perp}\psi} 
       &=\OpBW{\hat{\xi}^\perp\cdot {\bm A}_{-1}^{T,\Omega}\  \hat{\xi}^\perp+a_{\bra{-2,\geq1}}\pare{\eta;x,\xi}}\psi\\
        &+\OpBW{-\frac{\ii\hat{\xi}^\perp}{\av{\xi}}\cdot A_{1}^{T,\eta}\bra{\nabla^{\perp}\psi}\pare{\eta;x,\xi}+a_{\bra{-1,\geq1}}\pare{\eta;x,\xi}}\eta+r^{\lin,\bR}_{\bra{\varrho,\geq1}}\pare{\eta,\psi}\cdot\begin{bmatrix}
            \eta\\
            \psi
        \end{bmatrix}.
  \end{align}
  Making appeal to \cref{eq:ATOmega-1}, we can write
\begin{equation}
   \hat{\xi}^\perp\cdot {\bm A}_{-1}^{T,\Omega}\  \hat{\xi}^\perp
   =
   \frac{1}{2}\,
   \hat\xi^\perp \cdot
   \pare{
   \Delta_\xi \phi^{\mathpzc{B}_0}_1 \,{\bm H}_x\eta
   - \Delta_\xi \phi^{\mathpzc{B}_0}_1 \Delta_x\eta \,\Id
   + \Delta_x\eta \,{\bm H}_\xi \phi^{\mathpzc{B}_0}_1
   }
   \hat\xi^\perp
   = \mathsf{b}_1 + \mathsf{b}_2 + \mathsf{b}_3,
\end{equation}
where
\begin{align}
    \mathsf{b}_1
    &:=
    \frac{\Delta_\xi \phi^{\mathpzc{B}_0}_1}{2}\,
    \hat\xi^\perp \cdot {\bm H}_x\eta \,\hat\xi^\perp,
    &
    \mathsf{b}_2
    &:=
    -\frac{\Delta_\xi \phi^{\mathpzc{B}_0}_1 \Delta_x\eta}{2},
    &
    \mathsf{b}_3
    &:=
    \frac{\Delta_x\eta}{2}\,
    \hat\xi^\perp \cdot {\bm H}_\xi \phi^{\mathpzc{B}_0}_1 \hat\xi^\perp.
\end{align}
Since $\phi^{\mathpzc{B}_0}_1$ is $1$-homogeneous in $\xi$, Euler's formula implies that
\begin{equation}
    {\bm H}_\xi \phi^{\mathpzc{B}_0}_1 \,\xi = 0.
\end{equation}
Therefore, in the orthonormal basis $(\hat\xi,\hat\xi^\perp)$,
\begin{equation}
    \Delta_\xi \phi^{\mathpzc{B}_0}_1
    =
    \tr \pare{{\bm H}_\xi \phi^{\mathpzc{B}_0}_1}
    =
    \hat\xi^{\perp \intercal} {\bm H}_\xi \phi^{\mathpzc{B}_0}_1 \hat\xi^\perp.
\end{equation}
In particular, this implies that
\begin{equation}
    \mathsf{b}_3
    =
    \frac{\Delta_x\eta}{2}\,\Delta_\xi \phi^{\mathpzc{B}_0}_1
    =
    -\mathsf{b}_2.
\end{equation}
Hence,
\begin{equation}
\hat{\xi}^\perp {\bm A}_{-1}^{T,\Omega}\ \hat\xi^\perp
   =
   \mathsf{b}_1
   =
   \frac{\Delta_\xi \phi^{\mathpzc{B}_0}_1}{2}\,
   \hat\xi^\perp \cdot {\bm H}_x\eta \,\hat\xi^\perp.
\end{equation}
Finally, recalling that
\begin{equation}
   \phi^{\mathpzc{B}_0}_1\pare{\eta;x,\xi}
   =
   -\frac{\lambda_1\pare{\eta;x,\xi}}{1+\av{\nabla\eta(x)}^2},
   \qquad
   \lambda_1\pare{\eta;x,\xi}
   =
   \frac{\sqrt{\av{\xi}^2+\pare{\nabla^\perp\eta(x)\cdot\xi}^2}}{2},
\end{equation}
we can compute explicitly
\begin{equation}
   \Delta_\xi \phi^{\mathpzc{B}_0}_1
   =
   -\frac{\av{\xi}^2}{16\,\lambda_1^3\pare{\eta;x,\xi}}\cdot
\end{equation}
Therefore,
\begin{equation}
   \hat{\xi}^\perp {\bm A}_{-1}^{T,\Omega}\ \hat\xi^\perp
   =
   -\frac{\av{\xi}^2}{32\,\lambda_1^3\pare{\eta;x,\xi}}
   \hat\xi^\perp \cdot {\bm H}_x\eta \,\hat\xi^\perp.
\end{equation}
Besides, according to \cref{eq:AT1extra}, one readily has
\begin{align}
    -\frac{\ii\hat{\xi}^\perp}{\av{\xi}}\cdot A_{1}^{T,\eta}\bra{\nabla^{\perp}\psi}\pare{\eta;x,\xi}&=- v_3^{\star}\bra{\nabla^\perp\psi}\pare{\eta;x}
    -\ii \ \nabla \psi \pare{x} \cdot A_{0}^{\mathpzc{B}_0,\omega}\pare{\eta;x,\xi}.
\end{align}
This gives the desired result.
\end{proof}

\subsubsection{\texorpdfstring{Paralinearization of $S\pare{\Atw\middle|\eta}$}{Paralinearization of S(a|η)}}
Let us recall that the operator $S\pare{\Atw\middle|\eta}$ was first defined in \cref{eq:S}, which we repeat here for the sake of readability:
\begin{equation}\label{eq:S_recall}
    S\pare{\Atw\middle|\eta}\psi\defeq \nabla^{\perp}\pare{1+\Atw \tau_0\pare{\eta}}^{-1}\psi ,
\end{equation}
where $\tau_0\pare\eta$ is defined in \cref{eq:Tstar}.
We shall also use the projector onto divergence-free vector fields
\begin{equation}
\label{eq:Psigma_S_def}
    \mathcal{P}
    \defeq
    \nabla^\perp\Delta^{-1}\nabla^\perp\cdot,
    \qquad
    \mathcal{P}\pare{\xi}
    =
    \frac{\xi^\perp\otimes\xi^\perp}{\av{\xi}^2}\cdot
\end{equation}

\begin{lemma}[Uniform operator inverse for the canonical trace map]
\label{lem:uniform_inverse_tau0}
There exist $s_0>0$ and $r_0>0$ such that, for every $s\geq s_0$,
every $\Atw\in[-2,2]$, and every
$\eta\in\BallR{s_0}{r_0}\cap\cC\pare{I;H^s}$, the operators
\[
    \Id+\Atw\tau_0\pare{\eta}:H^s\longrightarrow H^s,
    \qquad
    \Id+\Atw\mathcal P T\pare{\eta}:
    \pare{H^{s-1}}^2\longrightarrow\pare{H^{s-1}}^2
\]
are invertible, uniformly with respect to $\Atw$.  More precisely, for
$f\in H^s$,
\begin{equation}
\label{eq:uniform_inverse_tau0_tame}
\begin{aligned}
\sup_{\Atw\in[-2,2]}
\norm{\pare{\Id+\Atw\tau_0\pare{\eta}}^{-1}f}_{H^s}
\lesssim_s{}&
\norm{f}_{H^s}
+\norm{\eta}_{H^s}\norm{f}_{H^{s_0}},
\\
\sup_{\Atw\in[-2,2]}
\norm{\pare{\pare{\Id+\Atw\tau_0\pare{\eta}}^{-1}-\Id}f}_{H^s}
\lesssim_s{}&
\norm{\eta}_{H^{s_0}}\norm{f}_{H^s}
+\norm{\eta}_{H^s}\norm{f}_{H^{s_0}}.
\end{aligned}
\end{equation}
For $F\in\pare{H^{s-1}}^2$,
\begin{equation}\label{eq:S_vector_inverse_tame}
\begin{aligned}
\sup_{\Atw\in[-2,2]}
\norm{\pare{\Id+\Atw\mathcal P T\pare\eta}^{-1}F}_{H^{s-1}}
\lesssim_s{}&
\norm F_{H^{s-1}}+\norm\eta_{H^s}\norm F_{H^{s_0-1}},
\\
\sup_{\Atw\in[-2,2]}
\norm{\pare{\pare{\Id+\Atw\mathcal P T\pare\eta}^{-1}-\Id}F}_{H^{s-1}}
\lesssim_s{}&
\norm\eta_{H^{s_0}}\norm F_{H^{s-1}}
+\norm\eta_{H^s}\norm F_{H^{s_0-1}}.
\end{aligned}
\end{equation}
Moreover, for every $\background\in\bR^2$,
\begin{equation}\label{eq:S_constant_trace_bound}
    \norm{T\pare\eta\background}_{H^{s-1}}
    \lesssim_s\av\background\norm\eta_{H^s}.
\end{equation}
The same uniform tame bounds hold for differences and derivatives in
$\eta$, and the inverses depend real-ana\-lytically on $\Atw\in[-2,2]$ and
on $\eta$ in these seminorms.  Consequently, the substitution
$f\mapsto\allowbreak\pare{\Id+\Atw\tau_0\pare{\eta}}^{-1}f$ in the coefficient arguments of the symbol and
remainder maps of \cref{prop:Tstar_paralinearization} preserves their
$\Gamma^m_{\geq q}$ and $\mathcal R^\varrho$ classes, with the vanishing
index $q$ counted jointly in all displayed variables.
\end{lemma}

\begin{proof}
At low frequency, integration by parts in \cref{eq:Tstar_recall} gives
\[
\tau_0\pare\eta f\pare x
=\frac1{4\pi}\int_{\bR^2}
\frac{\pare{x-y}\cdot\nabla\eta\pare y-\eta\pare x+\eta\pare y}
{\pare{\av{x-y}^2+\pare{\eta\pare x-\eta\pare y}^2}^{3/2}}
f\pare y\dd y.
\]
The local kernel is bounded by
$C\norm{{\bm H}_x\eta}_{L^\infty}\av{x-y}^{-1}$; at infinity its
gradient term is Hilbert--Schmidt and its remaining terms satisfy
Schur's test. Thus $\norm{\tau_0\pare\eta}_{\cL\pare{L^2}}
\lesssim\norm\eta_{H^{s_0}}$.
On high frequencies, the tame bounds for $T\pare{\eta}$ from
\cref{eq:T(eta),prop:B0_paralinearization,prop:paralinearization_Q,thm:contS}
and the definition
$\tau_0=\Delta^{-1}\nabla^\perp\!\cdot T\pare{\eta}\nabla^\perp$ give
\[
\norm{\tau_0\pare{\eta}f}_{H^s}
\lesssim_s
\norm{\eta}_{H^{s_0}}\norm{f}_{H^s}
+\norm{\eta}_{H^s}\norm{f}_{H^{s_0}}.
\]
In particular, after decreasing $r_0$, one has
$\sup_{\Atw\in[-2,2]}\norm{\Atw\tau_0\pare{\eta}}_{\cL(H^{s_0})}<1/2$.
The Neumann series gives the inverse uniformly on $H^{s_0}$, and the tame
estimate above, with the order-$(-1)$ input term in
\cref{prop:Tstar_paralinearization}, lifts it to every $H^s$, yielding
\cref{eq:uniform_inverse_tau0_tame}.  Differentiating the resolvent identity
\[
d_\eta\pare{\Id+\Atw\tau_0}^{-1}
=-\pare{\Id+\Atw\tau_0}^{-1}
\Atw\pare{d_\eta\tau_0}
\pare{\Id+\Atw\tau_0}^{-1}
\]
and iterating proves the difference and higher tame estimates.  These
estimates, together with the locally convergent Neumann expansion, give the
stated analytic dependence and stability of the symbol and remainder
classes and preserve their joint vanishing order.

The kernel of $T\pare\eta$ in \cref{eq:T(eta)} has the same local and
tail bounds. Thus \cref{prop:T_paralinearization,thm:contS} and the
boundedness of $\mathcal P$ give
\[
\norm{\mathcal P T\pare\eta F}_{H^{s-1}}
\lesssim_s
\norm\eta_{H^{s_0}}\norm F_{H^{s-1}}
+\norm\eta_{H^s}\norm F_{H^{s_0-1}}.
\]
The same Neumann-series and resolvent argument, using the order-$(-1)$
input term in \cref{prop:T_paralinearization}, proves
\cref{eq:S_vector_inverse_tame}.

For constant inputs, integration by parts in \cref{eq:B0,eq:Q} gives
\begin{equation}\label{eq:S_constant_kernel_identities}
\begin{aligned}
\mathpzc B_0\pare\eta\bra1
&=-\mathpzc q\pare\eta\nabla\eta,\\
\mathpzc q\pare\eta\bra1
&=\mathpzc B_0\pare\eta\cdot\nabla\eta,
\\
T\pare\eta\background
&=\mathpzc B_0\pare\eta\bra{\background\cdot\nabla\eta}
-\background\pare{\mathpzc B_0\pare\eta\cdot\nabla\eta}
-\pare{\pare{\mathpzc q\pare\eta\nabla\eta}^{\perp}
\cdot\background}\nabla^\perp\eta.
\end{aligned}
\end{equation}
The boundary terms vanish by oddness and decay, subtracting the flat kernel.
Now \cref{prop:B0_paralinearization,prop:paralinearization_Q,thm:contS}
give \cref{eq:S_constant_trace_bound}; the same calculation after
differentiation gives the asserted analytic dependence.
\end{proof}

{

\begin{prop}[Paralinearization of $S\pare{\Atw\middle|\eta}$]
\label{prop:paralinearization_S}
There exists $\varrho_0 > 0$ such that for any $\varrho \geq \varrho_0$, there exists $s_0 > 0$ and $r_0 > 0$ such that for all $\eta,\psi,\tilde{\psi}\in \BallR{s_0}{r_0}$ and all $\background\in\bR^2$ the following statements hold.

\begin{enumerate}
\item
The operator $\pare{1+\Atw \tau_0\pare{\eta}}^{-1}\psi$ admits the paralinearization
\begin{equation}
\label{eq:inv-paralinearaTstar}
\begin{aligned}
    \pare{1+\Atw \tau_0\pare{\eta}}^{-1}\psi
    =
    &\OpBW{
    \vect{0}{1}
    +
    \Atw
    \pare{
    -
    \vect{
    a^{\textnormal{inv},\eta}_0\bra{\psi}
    \pare{\Atw \middle| \eta;x,\xi}
    }{
    a^{\tau_0,\psi}_{-1}\pare{\eta;x,\xi}
    }
    +
    \vect{
    a_{\bra{-1,\geq1}}\pare{\Atw\middle|\eta,\psi;x,\xi}
    }{
    a_{\bra{-2,\geq1}}\pare{\Atw\middle|\eta;x,\xi}
    }
    }}
    \cdot
    \vect{\eta}{\psi}\\
    +
    &\Atw
    R^{\lin,\bR}_{\bra{\varrho,\geq1}}\pare{\Atw \middle| \eta, \psi}
    \cdot
    \vect{\eta}{\psi},
\end{aligned}
\end{equation}
where
\begin{itemize}
    \item $\Atw$ is the Atwood number defined in \cref{eq:Atwood};

    \item $a^{\textnormal{inv},\eta}_0\in \cL\pare{H^{s_0};\Gamma_{\geq0}^0}$ is defined by
    \begin{equation}
    \label{eq:aSeta0}
        a^{\textnormal{inv},\eta}_0\bra{\psi}
        \pare{\Atw\middle|\eta;x,\xi}
        \defeq
        -v_3\bra{\psi}\pare{\Atw\middle|\eta;x}
        -
        \ii A_{0}^{\mathpzc{B}_0,\omega;\perp}\pare{\eta;x,\xi}
        \cdot
        S\pare{\Atw\middle|\eta}\psi,
    \end{equation}
    where $v_3\bra{\psi}\pare{\Atw\middle|\eta;x}$ is defined in \cref{eq:bmv}, $A_{0}^{\mathpzc{B}_0,\omega}$ is defined in \cref{eq:H0_definition}, and $S\pare{\Atw\middle|\eta}$ is defined in \cref{eq:S_recall};

    \item $a^{\tau_0,\psi}_{-1}\in\Gamma_{\geq1}^{-1}$ is defined in \cref{eq:atau0};

    \item $a_{\bra{m,\geq1}}\in\Gamma_{\geq1}^m$, $m\in\bR$, denotes a generic symbol as per \cref{notation:paralinearization};

    \item $R^{\lin,\bR}_{\bra{\varrho,\geq1}}\in\pare{\cR_{\geq1}^\varrho}^2$ denotes a generic smoothing remainder as per \cref{notation:paralinearization}.
\end{itemize}

\item
The operator $S\pare{\Atw\middle|\eta}\psi$ admits the paralinearization
\begin{equation}
\label{eq:paralinearization_S}
\begin{aligned}
    S\pare{\Atw\middle|\eta}\psi
    =
    &\OpBW{
    {\bm A}^S_1
    \pare{\Atw\middle|\eta,\psi;x,\xi}
    +
    \Atw
    \vectT{
    A_{\bra{0,\geq1}}\pare{\Atw\middle|\eta,\psi;x,\xi}
    }{
    A_{\bra{-1,\geq1}}\pare{\Atw\middle|\eta;x,\xi}
    }}
    \vect{\eta}{\psi}
    +
    \Atw
    {\bm R}^{\lin,\bR}_{\bra{\varrho,\geq1}}
    \pare{\Atw \middle| \eta, \psi}
    \vect{\eta}{\psi},
\end{aligned}
\end{equation}
where
\begin{itemize}
    \item ${\bm A}^S_1\in \pare{\Gamma_{\geq0}^1}^{2\times 2}$ is defined by
    \begin{equation}
        {\bm A}^S_1
        \pare{\Atw\middle|\eta,\psi;x,\xi}
        \defeq
        \vectT{
        A^{S,\eta}_1\bra{\psi}
        \pare{\Atw\middle|\eta;x,\xi}
        }{
        A^{S,\psi}_1
        \pare{\Atw\middle|\eta;x,\xi}
        },
    \end{equation}
    with
    \begin{equation}
    \label{eq:ASeta1_sobolev}
        A^{S,\eta}_1\bra{\psi}
        \pare{\Atw\middle|\eta;x,\xi}
        \defeq
        -\ii\Atw\xi^\perp
        a^{\textnormal{inv},\eta}_0\bra{\psi}
        \pare{\Atw\middle|\eta;x,\xi},
    \end{equation}
    and
    \begin{equation}
    \label{eq:ASpsi}
        A^{S,\psi}_1
        \pare{\Atw\middle|\eta;x,\xi}
        \defeq
        \ii\xi^\perp
        \pare{
        1-\Atw a^{\tau_0,\psi}_{-1}\pare{\eta;x,\xi}
        };
    \end{equation}

    \item $A_{\bra{m,\geq1}}\in\pare{\Gamma_{\geq1}^m}^2$, $m\in\bR$, denotes a generic vector symbol as per \cref{notation:paralinearization};

    \item ${\bm R}^{\lin,\bR}_{\bra{\varrho,\geq1}}\in\pare{\cR_{\geq1}^\varrho}^{2\times2}$ denotes a generic smoothing remainder as per \cref{notation:paralinearization}.
\end{itemize}

\item
Let $\bar{\psi}$ be the affine background potential defined in
\cref{eq:psibar_expression}. The canonical background vorticity
$\tilde{\Omega}\bra{\background}\pare{\eta}$, defined in
\cref{eq:Omega_bg_can_def}, satisfies the exact identity
\begin{equation}
\label{eq:Omega_bg_can_S_inverse}
   \tilde{\Omega}\bra{\background}\pare{\eta}
    =
    \pare{\Id+\Atw\mathcal{P}T\pare{\eta}}^{-1}\background.
\end{equation}
Moreover,
\begin{equation}
\label{eq:paralinearization_S_background}
   \tilde{\Omega}\bra{\background}\pare{\eta}
    =
    \background
    +
    \OpBW{
    A^{S,\eta}_{1,\textnormal{bg}}\bra{\background}
    \pare{\Atw\middle|\eta;x,\xi}
    +
    A_{\bra{0,\geq1}}\bra{\background}
    \pare{\Atw\middle|\eta;x,\xi}
    }\eta
    +
    R^{\lin,\bR}_{\bra{\varrho,\geq1}}\bra{\background}
    \pare{\Atw\middle|\eta}\eta,
\end{equation}
where
\begin{itemize}
        \item $A^{S,\eta}_{1,\textnormal{bg}}\in\cL\pare{\bR^2;\pare{\Gamma_{\geq0}^1}^2}$ is linear in the separately estimated variable and is defined by
    \begin{equation}\label{eq:ASeta1_background}
        \begin{aligned}
            A^{S,\eta}_{1,\textnormal{bg}}\bra{\background}
        \pare{\Atw\middle|\eta;x,\xi}
        &\defeq
        -\Atw
        \mathcal{P}\pare{\xi}
        A^{T,\eta}_1\bra{
        \tilde{\Omega}\bra{\background}\pare{\eta}
        }
        \pare{\eta;x,\xi}
        \\
        &=
        -\Atw
        \pare{
        A_{0}^{\mathpzc{B}_0,\omega;\perp}\pare{\eta;x,\xi}\cdot
        \tilde{\Omega}\bra{\background}\pare{\eta;x}
        -
        \ii v_3^{\star}\bra{
        \tilde{\Omega}\bra{\background}\pare{\eta}
        }\pare{\eta;x}
        }
        \xi^\perp;
        \end{aligned}
\end{equation}
    \item $A_{\bra{0,\geq1}}\bra{\background}\in\cL\pare{\bR^2;\pare{\Gamma_{\geq1}^0}^2}$ is linear in $\background$;

    \item $R^{\lin,\bR}_{\bra{\varrho,\geq1}}\bra{\background}\in\cL\pare{\bR^2;\pare{\cR_{\geq1}^\varrho}^2}$ is linear in $\background$.
\end{itemize}

\end{enumerate}
\end{prop}

}

\begin{lemma}[Low-frequency trace symbols]
\label{lem:S_low_frequency_smoothing}
For every $\varrho\geq0$, one has
\begin{equation}\label{eq:S_order_minus_one_low_frequency_smoothing}
\begin{aligned}
    \OpBW{\pare{a^{\tau_0,\psi}_{-1}}^\flat}
    &\in\cR^\varrho_{\geq1},
    &
    \OpBW{\pare{{\bm A}^{T,\Omega}_{-1}}^\flat},
    \OpBW{\pare{\mathcal P\pare\xi{\bm A}^{T,\Omega}_{-1}}^\flat}
    &\in\pare{\cR^\varrho_{\geq1}}^{2\times2},
\end{aligned}
\end{equation}
with $^\flat$ as in \cref{notation:high-frequency-symbol-extension}
and the symbols defined in \cref{eq:atau0,eq:ATOmega-1,eq:Psigma_S_def}.
These operators are real-to-real.
\end{lemma}

\begin{proof}
By \cref{eq:atau0}, for \(\hat{\xi}\in\bS^1\),
\[
    a^{\tau_0,\psi}_{-1}\pare{\eta;x,\hat{\xi}}
    =
    -\frac{
        \hat{\xi}^\perp\cdot{\bm H}_x\eta\,\hat{\xi}^\perp
    }{
        4\pare{
            1+\pare{\nabla^\perp\eta\cdot\hat{\xi}}^2
        }^{3/2}
    }.
\]
Differentiating with respect to the angular variable on \(\bS^1\) gives,
pointwise in \(x\),
\[
    \norm{
        a^{\tau_0,\psi}_{-1}
        \pare{\eta;x,\bullet}
    }_{H^1\pare{\bS^1}}
    \lesssim
    \av{{\bm H}_x\eta}
    \pare{1+\av{\nabla\eta}}.
\]
Hence, fixing \(s_0>2\) and
\(\eta\in\BallR{s_0}{r_0}\), Sobolev embedding in dimension two yields
\begin{equation}\label{eq:S_low_frequency_angular_bound}
    \norm{
        a^{\tau_0,\psi}_{-1}
        \pare{\eta;\bullet,\bullet}
    }_{L^2_xH^1_{\hat{\xi}}\pare{\bS^1}}
    \lesssim
    \norm\eta_{H^{s_0}}.
\end{equation}

Since \(a^{\tau_0,\psi}_{-1}\) is homogeneous of degree \(-1\) in its
frequency variable, denoting by
\(\widehat{\phantom{a}}\) the Fourier transform in \(x\) only, one has
\begin{equation}\label{eq:S_low_frequency_homogeneity}
    \widehat{a^{\tau_0,\psi}_{-1}}
    \pare{\eta;\theta,\zeta}
    =
    \av\zeta^{-1}
    \widehat{a^{\tau_0,\psi}_{-1}}
    \pare{\eta;\theta,\hat\zeta},
    \qquad
    \zeta\neq0.
\end{equation}
By \cref{BonyWeyl,notation:high-frequency-symbol-extension},
\[
    \widehat{
        \OpBW{\pare{a^{\tau_0,\psi}_{-1}}^\flat\pare\eta}h
    }\pare\xi
    =
    \frac1{2\pi}
    \int_{\bR^2}
    K_\eta\pare{\xi,\varsigma}
    \widehat h\pare\varsigma\,\dd\varsigma,
\]
where
\begin{equation}\label{eq:S_low_frequency_Fourier_kernel}
\begin{aligned}
    K_\eta\pare{\xi,\varsigma}
    \defeq{}&
    \chi_\epsilon\pare{
        \frac{2\pare{\xi-\varsigma}}
        {\av{\xi+\varsigma}}
    }
    \chi\pare{\frac{\xi+\varsigma}{2}}
    \widehat{a^{\tau_0,\psi}_{-1}}
    \pare{
        \eta;
        \xi-\varsigma,
        \frac{\xi+\varsigma}{2}
    }.
\end{aligned}
\end{equation}
On the support of this kernel,
\cref{BonyWeyl,rela:para} gives
\begin{equation}\label{eq:S_low_frequency_kernel_support_relations}
    \av{\xi-\varsigma}
    \lesssim
    \av\varsigma,
    \qquad
    \av{\xi+\varsigma}
    \sim
    \av\varsigma,
\end{equation}
and the analogous estimates with \(\xi\) and \(\varsigma\) interchanged.

We first verify the first condition in
\cref{eq:Schur_test_assumptions}. For \(\varsigma\neq0\), set
\(\theta=\xi-\varsigma\). Using
\cref{eq:S_low_frequency_homogeneity,eq:S_low_frequency_kernel_support_relations},
Cauchy--Schwarz, and the fact that
\(\av{B\pare{0,C\av\varsigma}}^{1/2}\lesssim\av\varsigma\) in dimension
two, we obtain
\[
\begin{aligned}
    \int_{\bR^2}
    \av{K_\eta\pare{\xi,\varsigma}}
    \,\dd\xi
    &\lesssim
    \av\varsigma^{-1}
    \int_{\av\theta\lesssim\av\varsigma}
    \sup_{\hat{\xi}\in\bS^1}
    \av{
        \widehat{a^{\tau_0,\psi}_{-1}}
        \pare{\eta;\theta,\hat{\xi}}
    }
    \,\dd\theta
    \lesssim
    \left(
        \int_{\bR^2}
        \sup_{\hat{\xi}\in\bS^1}
        \av{
            \widehat{a^{\tau_0,\psi}_{-1}}
            \pare{\eta;\theta,\hat{\xi}}
        }^2
        \,\dd\theta
    \right)^{1/2}.
\end{aligned}
\]
Using
\(H^1\pare{\bS^1}\hookrightarrow L^\infty\pare{\bS^1}\),
Plancherel's identity, and
\cref{eq:S_low_frequency_angular_bound}, this gives
\begin{equation}\label{eq:S_low_frequency_first_Schur_bound}
    \sup_{\varsigma\in\bR^2}
    \int_{\bR^2}
    \av{K_\eta\pare{\xi,\varsigma}}
    \,\dd\xi
    \lesssim
    \norm\eta_{H^{s_0}}.
\end{equation}
For \(\varsigma=0\) the kernel vanishes by the support of
\(\chi_\epsilon\), so the same estimate holds without restriction.
Fixing instead \(\xi\neq0\) and making the same change of variables
\(\theta=\xi-\varsigma\), the support condition gives
\[
    \av\theta\lesssim\av\xi,
    \qquad
    \av{\xi+\varsigma}\sim\av\xi.
\]
The same argument therefore gives the second Schur bound
\begin{equation}\label{eq:S_low_frequency_second_Schur_bound}
    \sup_{\xi\in\bR^2}
    \int_{\bR^2}
    \av{K_\eta\pare{\xi,\varsigma}}
    \,\dd\varsigma
    \lesssim
    \norm\eta_{H^{s_0}}.
\end{equation}
Again, the estimate is trivial for \(\xi=0\).
Applying \cref{lem:Schur_test} to
\cref{eq:S_low_frequency_Fourier_kernel}, with
\cref{eq:S_low_frequency_first_Schur_bound,eq:S_low_frequency_second_Schur_bound},
and then using Plancherel's identity, we obtain
\begin{equation}\label{eq:S_low_frequency_L2_bound}
    \norm{
        \OpBW{\pare{a^{\tau_0,\psi}_{-1}}^\flat\pare\eta}h
    }_{L^2}
    \lesssim
    \norm\eta_{H^{s_0}}
    \norm h_{L^2}.
\end{equation}
The cutoff
\(\chi\pare{\frac{\xi+\varsigma}{2}}\) implies $\av{\xi+\varsigma}\lesssim1$
on the support of \(K_\eta\). Together with
\cref{eq:S_low_frequency_kernel_support_relations}, this yields $  \av\xi+\av\varsigma\lesssim1
   $on $\supp K_\eta$. 
Consequently, for every \(s\geq s_0\) and every \(\varrho\geq0\) $ \angles\xi^s
    \angles\varsigma^{-(s-\varrho)}
    \lesssim_{s,\varrho}1
     $ on $ \supp K_\eta .$
Hence the weighted kernel $ K_{\eta,s,\varrho}\pare{\xi,\varsigma}
    \defeq
    \angles\xi^s
    K_\eta\pare{\xi,\varsigma}
    \angles\varsigma^{-(s-\varrho)}$
satisfies, by
\cref{eq:S_low_frequency_first_Schur_bound,eq:S_low_frequency_second_Schur_bound},
the same Schur estimates:
\[
    \sup_{\varsigma\in\bR^2}
    \int_{\bR^2}
    \av{K_{\eta,s,\varrho}\pare{\xi,\varsigma}}
    \,\dd\xi
    +
    \sup_{\xi\in\bR^2}
    \int_{\bR^2}
    \av{K_{\eta,s,\varrho}\pare{\xi,\varsigma}}
    \,\dd\varsigma
    \lesssim_{s,\varrho}
    \norm\eta_{H^{s_0}}.
\]
Applying \cref{lem:Schur_test} to
\(K_{\eta,s,\varrho}\) and using Plancherel's identity gives
\begin{equation}\label{eq:low_frequency_trace_smoothing_bound}
    \norm{
        \OpBW{\pare{a^{\tau_0,\psi}_{-1}}^\flat\pare\eta}h
    }_{H^s}
    \lesssim_{s,\varrho}
    \norm\eta_{H^{s_0}}
    \norm h_{H^{s-\varrho}}.
\end{equation}
We now consider \({\bm A}^{T,\Omega}_{-1}\). By
\cref{eq:ATOmega-1,eq:phiB0}, its entries are homogeneous of degree
\(-1\) in \(\xi\) and, after restriction to \(\bS^1\), are linear in
\({\bm H}_x\eta\), with coefficients depending smoothly on
\(\pare{\nabla\eta,\hat{\xi}}\). Consequently,
\begin{equation}\label{eq:S_low_frequency_ATOmega_angular_bound}
    \norm{
        {\bm A}^{T,\Omega}_{-1}
        \pare{\eta;\bullet,\bullet}
    }_{L^2_xH^1_{\hat{\xi}}\pare{\bS^1}}
    \lesssim
    \norm\eta_{H^{s_0}}.
\end{equation}
The argument leading from
\cref{eq:S_low_frequency_angular_bound} to
\cref{eq:low_frequency_trace_smoothing_bound} therefore applies
componentwise and gives
\[
    \OpBW{\pare{{\bm A}^{T,\Omega}_{-1}}^\flat}
    \in
    \pare{\cR^\varrho_{\geq1}}^{2\times2}.
\]
Finally, \(\mathcal P\pare\xi\) is smooth and homogeneous of degree zero
on \(\bR^2\setminus\set{0}\). Hence
\(\mathcal P\pare\xi{\bm A}^{T,\Omega}_{-1}\) satisfies the same
homogeneity and angular estimate as
\({\bm A}^{T,\Omega}_{-1}\), and the same argument gives
\[
    \OpBW{
        \pare{
            \mathcal P\pare\xi
            {\bm A}^{T,\Omega}_{-1}
        }^\flat
    }
    \in
    \pare{\cR^\varrho_{\geq1}}^{2\times2}.
\]

All the estimates above have the fixed baseline required in
\cref{def:smoothing} and vanish to order at least one in \(\eta\).
This proves \cref{eq:S_order_minus_one_low_frequency_smoothing}.
The three symbols are real and even in \(\xi\); since the cutoffs in
\cref{BonyWeyl,def-chi} are real and even, the corresponding
Bony--Weyl operators are real-to-real.
\end{proof}

\begin{proof}[Proof of \cref{prop:paralinearization_S}]

\begin{proofpart}[Sobolev input]
     Applying \cref{prop:Tstar_paralinearization}, we have
    \begin{equation}
    \begin{aligned}
        &\pare{1+\Atw \tau_0\pare{\eta}}\psi\\
        =& \OpBW{1+\Atw
         \pare{a^{\tau_0, \psi}_{-1} + a_{\bra{-2,\geq1}}}\pare{\eta;x,\xi}}\psi+\Atw\OpBW{a^{\tau_0,\eta}_0\bra{\psi}\pare{\eta;x,\xi}+a_{\bra{-1,\geq1}}\pare{\eta,\psi;x,\xi}}\eta\\
        &\quad+ \Atw R^{\lin,\bR}_{\bra{\varrho,\geq1}}\pare{\eta, \psi}\cdot\vect{\eta}{\psi}.
    \end{aligned}
    \end{equation}
By \cref{lem:uniform_inverse_tau0}, at a fixed Sobolev index and after
increasing $s_0$, for every $s\geq s_0$,
\begin{equation}\label{eq:S_low_output_trace_smoothing}
    \norm{\chi\pare{D/4}\tau_0\pare\eta h}_{H^s}
    +\norm{\chi\pare{D/4}T\pare\eta H}_{H^s}
    \lesssim_{s,\varrho}\norm\eta_{H^{s_0}}
    \pare{\norm h_{H^{s-\varrho}}+\norm H_{H^{s-\varrho}}}.
\end{equation}
By \cref{lem:S_low_frequency_smoothing}, the contributions obtained by applying
\(\chi\pare{D/4}\) to the concrete symbols below are smoothing. We may therefore
consider their contribution after applying \(1-\chi\pare{D/4}\).

Indeed, by \cref{rela:para}, on the Fourier support of
\[
    \pare{1-\chi\pare{D/4}}\OpBW a
\]
one has
\[
    \av{\frac{\xi+\varsigma}{2}}
    \gtrsim
    \av\xi
    \gtrsim1.
\]
Hence the frequency argument of \(a\) stays away from the origin, and
\cref{notation:high-frequency-symbol-extension} gives
\[
    \pare{1-\chi\pare{D/4}}\OpBW a
    =
    \pare{1-\chi\pare{D/4}}\OpBW{a^\sharp}.
\]
Thus the symbols appearing in the high-frequency contribution may be replaced
by their fixed high-frequency extensions, while retaining the notation of
\cref{notation:paralinearization}.

By \cref{lem:S_low_frequency_smoothing}, replace $a^{\tau_0,\psi}_{-1}$
by $\pare{a^{\tau_0,\psi}_{-1}}^\sharp$ in the input symbol.
By \cref{eq:atau0}, the latter vanishes
for $\av\xi\leq1.1$ and, after decreasing $r_0$,
\[
\sup_{\av\Atw\leq2,x,\xi}
\av{\Atw\pare{a^{\tau_0,\psi}_{-1}}^\sharp}\leq\frac12,
\qquad
\frac{1}{1+\Atw\pare{a^{\tau_0,\psi}_{-1}}^\sharp}-1
\in\Gamma^{-1}_{\geq1}.
\]
    Observe that \[\pbra{\frac{1}{1+\Atw \pare{a^{\tau_0,\psi}_{-1}}^\sharp}}{1+\Atw \pare{a^{\tau_0,\psi}_{-1}}^\sharp} =\pbra{\frac{1}{1+\Atw \pare{a^{\tau_0,\psi}_{-1}}^\sharp}-1}{1+\Atw \pare{a^{\tau_0,\psi}_{-1}}^\sharp} =0.\]
   Therefore, composing the previous identity on the left with $\OpBW{\frac{1}{1+\Atw \pare{a^{\tau_0,\psi}_{-1}}^\sharp\pare{\eta;x,\xi}}}$ gives, by virtue of \cref{prop:composition_paradiff} and \cref{lem:composition_rules},
    \begin{equation}
        \begin{aligned}
        &\OpBW{\frac{1}{1+\Atw \pare{a^{\tau_0,\psi}_{-1}}^\sharp}\pare{\eta;x,\xi}}\pare{1+\Atw \tau_0\pare{\eta}}\psi\\=&\psi+\Atw \OpBW{a_{\bra{-2,\geq1}}\pare{\Atw\middle|\eta;x,\xi}}\psi+\Atw\OpBW{a^{\tau_0,\eta}_0\bra{\psi}\pare{\eta;x,\xi}+a_{\bra{-1,\geq1}}\bra{\psi}\pare{\Atw\middle|\eta;x,\xi}}\eta\\ +&\Atw R^{\lin,\bR}_{\bra{\varrho,\geq1}}\pare{\Atw \middle| \eta, \psi}\cdot\vect{\eta}{\psi}.
        \end{aligned}
    \end{equation}
We deduce that
     \begin{equation}
         \label{eq:Tinverse0}
          \begin{aligned}
    \psi= & \ \pare{1+\Atw \tau_0\pare{\eta}}^{-1}\pare{1+\Atw \tau_0\pare{\eta}}\psi\\
        = & \ \OpBW{\frac{1}{1+\Atw \pare{a^{\tau_0,\psi}_{-1}}^\sharp\pare{\eta;x,\xi}}}\pare{1+\Atw \tau_0\pare{\eta}}\psi+\Atw \OpBW{a_{\bra{-2,\geq1}}\pare{\Atw\middle|\eta;x,\xi}}\psi\\
        &+\Atw \OpBW{- a^{\tau_0,\eta}_0\bra{\psi}\pare{\eta;x,\xi}+a_{\bra{-1,\geq1}}\bra{\psi}\pare{\Atw\middle|\eta;x,\xi}}\eta +\Atw R^{\lin,\bR}_{\bra{\varrho,\geq1}}\pare{\Atw \middle| \eta, \psi}\cdot\vect{\eta}{\psi}. 
    \end{aligned}
     \end{equation}
    Denoting
    \begin{equation}
        \label{eq:psi1}
        \psi_1\defeq \pare{1+\Atw \tau_0\pare{\eta}}\psi,
    \end{equation}
    then using \cref{lem:composition_rules,eq:psi1} in \cref{eq:Tinverse0} we infer that
    \begin{equation}\label{init-inv0}
            \begin{aligned}
       &\pare{1+\Atw \tau_0\pare{\eta}}^{-1}\psi_1\\
       =& \ \OpBW{\frac{1}{1+\Atw \pare{a^{\tau_0,\psi}_{-1}}^\sharp \pare{\eta;x,\xi}}}\psi_1 + \Atw \OpBW{a_{\bra{-2,\geq1}}\pare{\Atw\middle|\eta;x,\xi}}\pare{1+\Atw \tau_0\pare{\eta}}^{-1}\psi_1\\
        &+\Atw \OpBW{- a^{\tau_0,\eta}_0\pare{\eta,\pare{1+\Atw \tau_0\pare{\eta}}^{-1}\psi_1;x,\xi}+a_{\bra{-1,\geq1}}\pare{\Atw\middle|\eta,\pare{1+\Atw \tau_0\pare{\eta}}^{-1}\psi_1;x,\xi}}\eta
        \\
       & 
        +\Atw R^{\lin,\bR}_{\bra{\varrho,\geq1}} \pare{\Atw\middle|\eta,\pare{1+\Atw \tau_0\pare{\eta}}^{-1}\psi_1}\cdot\vect{\eta}{\pare{1+\Atw \tau_0\pare{\eta}}^{-1}\psi_1}.
    \end{aligned}
    \end{equation}
     Notice that \cref{lem:uniform_inverse_tau0}, the definition of symbols and smoothing remainders (cf. \cref{def:sfr,def:smoothing}) imply that, substituting only in the coefficient arguments, using \cref{notation:deltagamma}
    % \begin{small}
            \begin{equation}\label{init-inv1}
        \begin{aligned}
           a^{\textnormal{inv},\eta}_0\bra{\psi_1} \pare{\Atw\middle| \eta;x,\xi}
            \defeq & \ 
            a^{\tau_0,\eta}_0\pare{\eta, \pare{1+\Atw \tau_0\pare{\eta}}^{-1}\psi_1;x,\xi} 
            & \in \Gamma^0_{\geq0},
            \\
            a_{\bra{m,\geq1}}\pare{\Atw\middle|\eta,\psi_1 ;x,\xi}
            \defeq & \ 
            a_{\bra{m,\geq1}}\pare{\Atw\middle|\eta, \pare{1+\Atw \tau_0\pare{\eta}}^{-1}\psi_1;x,\xi}
            & \in \Gamma_{\geq1}^{m},
            \\
            R^{\lin,\bR}_{\bra{\varrho,\geq1}}\pare{\Atw\middle|\eta, \psi_1}
            \defeq & \ 
               R^{\lin,\bR}_{\bra{\varrho,\geq1}}\pare{\Atw\middle|\eta,\pare{1+\Atw\tau_0\pare\eta}^{-1}\psi_1}
             & \in \pare{\cR_{\geq1}^\varrho}^2.
        \end{aligned}
    \end{equation}
    Recalling \cref{def a0star,eq:S}  we obtain \cref{eq:aSeta0} from the first equation of \cref{init-inv1}.
        % \end{small}
        Next, plugging \cref{init-inv1} into \cref{init-inv0} and moving the smoothing term acting on $\pare{1+\Atw\tau_0\pare\eta}^{-1}\psi_1$ to the left, we get
         \begin{equation}\label{init-inv}
    \begin{aligned}
        &\pare{\Id+\Atw r^\bR_{\bra{\varrho,\geq1}}\pare{\Atw\middle|\eta}}\pare{1+\Atw \tau_0\pare{\eta}}^{-1}\psi_1\\
        = & \ \OpBW{\frac{1}{1+\Atw \pare{a^{\tau_0,\psi}_{-1}}^\sharp\pare{\eta;x,\xi}}}\psi_1 + \Atw \OpBW{a_{\bra{-2,\geq1}}\pare{\Atw\middle|\eta;x,\xi}}\pare{1+\Atw \tau_0\pare{\eta}}^{-1}\psi_1\\
        +&\Atw \OpBW{- a^{\textnormal{inv},\eta}_0\bra{\psi_1}\pare{\Atw\middle|\eta;x,\xi}+a_{\bra{-1,\geq1}}\pare{\Atw\middle|\eta,\psi_1;x,\xi}}\eta
        \\
       +&\Atw r^\bR_{\bra{\varrho,\geq1}}\bra{\psi_1}\pare{\Atw\middle|\eta}\eta .
    \end{aligned}
    \end{equation}
    We prove by induction that for any $n\in\mathbb{N}^*,$ there exists $b_{\bra{-2^{n}}}\in\Gamma_{\geq1}^{-2^{n}}$, $c_{\bra{-1}}^{[n]}, a_{\bra{-1,\geq1}}^{[n]}\in\Gamma_{\geq1}^{-1}$, $r^{[n],\bR}_{\bra{\varrho,\geq1}}\in\cR^\varrho_{\geq1}$ and $R^{\lin,[n],\bR}_{\bra{\varrho,\geq1}}\in \pare{\cR_{\geq1}^\varrho}^2$ such that
    \begin{equation}\label{relat induc}
    \begin{aligned}
        &\pare{\OpBW{1+ b_{\bra{-2^{n}}}\pare{\Atw|\eta;x,\xi}}+\Atw r^{[n],\bR}_{\bra{\varrho,\geq1}}\pare{\Atw|\eta}}\pare{1+\Atw \tau_0\pare{\eta}}^{-1}\psi_1\\
=&\OpBW{1+c_{\bra{-1}}^{[n]}\pare{\Atw|\eta;x,\xi}}\psi_1+\Atw\OpBW{- a^{\textnormal{inv},\eta}_0\bra{\psi_1}\pare{\Atw|\eta;x,\xi}+a_{\bra{-1,\geq1}}^{[n]}\pare{\Atw|\eta,\psi_1;x,\xi}}\eta\\
        +&\Atw R^{\lin,[n],\bR}_{\bra{\varrho,\geq1}}\pare{\Atw|\eta, \psi_1}\cdot\vect{\eta}{\psi_1}.
    \end{aligned}
    \end{equation}
    Here $b_{\bra{-2^n}}=\Atw a_{\bra{-2^n,\geq1}}$ vanishes near $\xi=0$. In addition, the symbol $c_{\bra{-1}}^{[n]}$ admits the expansion
    $$c_{\bra{-1}}^{[n]}\pare{\Atw|\eta;x,\xi}=-\Atw \pare{a^{\tau_0,\psi}_{-1}}^\sharp\pare{\eta;x,\xi}+\Atw a_{\bra{-2,\geq1}}\pare{\Atw|\eta;x,\xi}.$$
    The initialization $n=1$ is a consequence of \cref{init-inv}, with
    \[
        b_{\bra{-2}}\defeq-\Atw a_{\bra{-2,\geq1}},\qquad
        a_{\bra{-1,\geq1}}^{[1]}\defeq a_{\bra{-1,\geq1}},\qquad
        r^{[1],\bR}_{\bra{\varrho,\geq1}}\defeq r^\bR_{\bra{\varrho,\geq1}},\qquad
        R^{\lin,[1],\bR}_{\bra{\varrho,\geq1}}
        \defeq\vect{r^\bR_{\bra{\varrho,\geq1}}\bra{\psi_1}}{0},
    \]
    and
    \[
    c_{\bra{-1}}^{[1]}\pare{\Atw|\eta;x,\xi}
    \defeq-\frac{\Atw\pare{a^{\tau_0,\psi}_{-1}}^\sharp\pare{\eta;x,\xi}}
    {1+\Atw\pare{a^{\tau_0,\psi}_{-1}}^\sharp\pare{\eta;x,\xi}}
    =-\Atw\pare{a^{\tau_0,\psi}_{-1}}^\sharp\pare{\eta;x,\xi}
    +\Atw a_{\bra{-2,\geq1}}\pare{\Atw|\eta;x,\xi}.
    \]
    At each of the finitely many steps, decrease $r_0$ so that $\av{b_{\bra{-2^n}}}\leq1/2$. Then \cref{prop:composition_paradiff} implies
    $$\OpBW{\frac{1}{1+b_{\bra{-2^{n}}}\pare{\Atw|\eta;x,\xi}}}\OpBW{1+b_{\bra{-2^{n}}}\pare{\Atw|\eta;x,\xi}}=\OpBW{1+b_{\bra{-2^{n+1}}}\pare{\Atw|\eta;x,\xi}}+\Atw r^\bR_{\bra{\varrho,\geq1}},$$
    where $r^\bR_{\bra{\varrho,\geq1}}\in\mathcal{R}_{\geq1}^{\varrho}$ and
    $$b_{\bra{-2^{n+1}}}\defeq \pare{\frac{1}{1+b_{\bra{-2^{n}}}}}\#_{\varrho}\pare{1+b_{\bra{-2^{n}}}}-1\in\Gamma^{-2^{n+1}}_{\geq1}.$$
    This defect still vanishes near $\xi=0$.
    Composing \cref{relat induc} on the left with this reciprocal proves the induction argument by \cref{prop:composition_paradiff,lem:composition_rules}, retaining all smoothing operators acting on $\pare{1+\Atw\tau_0\pare\eta}^{-1}\psi_1$ on the left.
    By taking $n=\underline n>\log_2\pare\varrho$, one has
    $\OpBW{b_{\bra{-2^{\underline n}}}}\in\cR^\varrho_{\geq1}$.
    Thus the left-hand side of \cref{relat induc} becomes
    \[\pare{\Id+\Atw r^\bR_{\bra{\varrho,\geq1}}\pare{\Atw|\eta}}
    \pare{1+\Atw\tau_0\pare\eta}^{-1}\psi_1.\]
    For each $s$, decrease the radius in \cref{def:smoothing} so that
    $\norm{\Atw r^\bR_{\bra{\varrho,\geq1}}}_{\cL\pare{H^s}}<1/2$.
    The Neumann series and
    \begin{equation}\label{eq:S_smoothing_inverse}
    \begin{aligned}
    \pare{\Id+\Atw r^\bR_{\bra{\varrho,\geq1}}}^{-1}-\Id
    &=-\Atw\pare{\Id+\Atw r^\bR_{\bra{\varrho,\geq1}}}^{-1}
    r^\bR_{\bra{\varrho,\geq1}},
    \\
    \norm{\pare{\pare{\Id+\Atw r^\bR_{\bra{\varrho,\geq1}}}^{-1}-\Id}h}_{H^s}
    &\lesssim_{s,\varrho}\av\Atw\norm\eta_{H^{s_0}}
    \norm h_{H^{s-\varrho}}
    \end{aligned}
    \end{equation}
    give a correction in $\cR^\varrho_{\geq1}$, also for finite matrix
    operators. Inverting this factor in \cref{relat induc} and using
    \cref{lem:composition_rules,lem:S_low_frequency_smoothing} with
    \[
        \OpBW{\pare{a^{\tau_0,\psi}_{-1}}^\sharp}
        =\OpBW{a^{\tau_0,\psi}_{-1}}
        -\OpBW{\pare{a^{\tau_0,\psi}_{-1}}^\flat}
    \]
    proves \cref{eq:inv-paralinearaTstar}.
 To prove \cref{eq:paralinearization_S} with target index $\varrho$, we apply \cref{eq:inv-paralinearaTstar} first with index $\varrho+1$ and then use \cref{eq:S} and \cref{prop:composition_paradiff}. Indeed, $\nabla^{\perp}=\OpBW{\ii\xi^{\perp}}$ maps $\mathcal R_{\geq1}^{\varrho+1}$ into $\pare{\mathcal R_{\geq1}^\varrho}^2$; relabeling the index gives the stated remainder.
\end{proofpart}

\begin{proofpart}[Affine background]
Set
\begin{equation}\label{eq:S_affine_perturbation}
\begin{aligned}
    \delta\Omega\bra{\background}\pare\eta
    &\defeq\tilde\Omega\bra{\background}\pare\eta-\background,
    \\
    \pare{\Id+\Atw\mathcal P T\pare\eta}
    \delta\Omega\bra{\background}\pare\eta
    &=-\Atw\mathcal P T\pare\eta\background,
    \\
    \norm{\delta\Omega\bra{\background}\pare\eta}_{H^{s-1}}
    &\lesssim_s\av\background\norm\eta_{H^s},
\end{aligned}
\end{equation}
by \cref{eq:S_vector_inverse_tame,eq:S_constant_trace_bound}. This gives
the unique solution in $\background+\pare{H^{s-1}}^2$ of
\begin{equation}
\label{eq:Omega_background_exact_equation}
   \tilde{\Omega}\bra{\background}\pare{\eta}
    +
    \Atw
    \mathcal{P}T\pare{\eta}
   \tilde{\Omega}\bra{\background}\pare{\eta}
    =
    \background.
\end{equation}
Since $\mathcal P\delta\Omega=\delta\Omega$, this is the
$\nabla^\perp$ derivative of the canonical trace equation modulo
constants, by \cref{eq:Tstar}; it proves
\cref{eq:Omega_bg_can_S_inverse}. Equivalently,
\begin{equation}
\label{eq:Omega_background_Neumann_series}
   \tilde{\Omega}\bra{\background}\pare{\eta}
    =
    \background
    -
    \Atw
    \mathcal{P}T\pare{\eta}
   \tilde{\Omega}\bra{\background}\pare{\eta}.
\end{equation}
Since $T\pare\eta$ is linear in its vorticity argument,
$\tilde\Omega\bra{\background}\pare\eta$ is linear in $\background$.
It suffices to work with $\av\background\leq1$ and extend by linearity.
After decreasing $r_0$, \cref{eq:S_affine_perturbation} places
$\delta\Omega\bra{\background}\pare\eta$ in the ball required by
\cref{prop:T_paralinearization}. We use smoothing index $\varrho+1$ below.

For the constant part, apply
\cref{prop:B0_paralinearization,prop:paralinearization_Q} to the Sobolev
inputs in \cref{eq:S_constant_kernel_identities}, with index $\varrho+2$,
and use \cref{eq:AT1,eq:AT1extra}. By
\cref{eq:first_variations_B0_q_flat,eq:T(eta)},
\[
    d_\eta T\pare0\bra h\background
    =-\frac12\OpBW{\frac{\xi^\perp\pare{\background\cdot\xi^\perp}}{\av\xi}}h
    =\OpBW{A^{T,\eta}_1\bra{\background}\pare{0;x,\xi}}h.
\]
Subtracting the flat terms therefore gives lower-order symbols and
remainders vanishing at $\eta=0$. In particular,
\[
    \chi\pare{D/4}\pare{
    T\pare\eta\background-
    \OpBW{A^{T,\eta}_1\bra{\background}\pare\eta}\eta}
    =R^\bR_{\bra{\varrho+1,\geq1}}\bra{\background}\pare\eta\eta:
\]
integrate the derivative along $t\eta$, $0\leq t\leq1$, and use the
fixed-index analytic bounds in \cref{lem:uniform_inverse_tau0} and
compact output support, increasing the baseline.

Apply \cref{eq:paraoperatorT} to $\delta\Omega\bra{\background}\pare\eta$
and add the constant part. By
\cref{lem:S_low_frequency_smoothing,eq:S_low_output_trace_smoothing},
\begin{equation}
\label{eq:T_constant_background}
\begin{aligned}
    T\pare{\eta}
    \tilde{\Omega}\bra{\background}\pare{\eta}
    &= 
    \OpBW{
    \pare{{\bm A}^{T,\Omega}_{-1}}^\sharp\pare{\eta;x,\xi}
    +
    {\bm A}_{\bra{-2,\geq1}}\pare{\eta;x,\xi}
    }
    \delta\Omega\bra{\background}\pare{\eta}
    \\
    &+
    \OpBW{
    A^{T,\eta}_1\bra{
    \tilde{\Omega}\bra{\background}\pare{\eta}
    }\pare{\eta;x,\xi}
    +
    A_{\bra{0,\geq1}}\bra{\background}\pare{\eta;x,\xi}
    }
    \eta
    \\
    &+
    {\bm R}^\bR_{\bra{\varrho+1,\geq1}}\pare{\Atw\middle|\eta}
    \delta\Omega\bra{\background}\pare\eta
    +
    R^{\lin,\bR}_{\bra{\varrho+1,\geq1}}\bra{\background}\pare{\eta}\eta,
\end{aligned}
\end{equation}
The matrix symbols and remainders depend only on $\eta$; the vector
symbols and remainders are linear in $\background$. Applying $\mathcal P$
to \cref{eq:T_constant_background}, with both factors localized at high
frequency in each symbolic composition, yields the following identity.
The low-frequency terms are controlled by
\cref{lem:S_low_frequency_smoothing} and the boundedness of $\mathcal P$;
for the order-one column use the same Schur estimate, retaining its
flat multiplier exactly.
\begin{equation}
\label{eq:PT_constant_background}
\begin{aligned}
    \mathcal{P}T\pare{\eta}
    \tilde{\Omega}\bra{\background}\pare{\eta}
    &=
    \OpBW{
    \pare{\mathcal{P}\pare{\xi}{\bm A}^{T,\Omega}_{-1}}^\sharp\pare{\eta;x,\xi}
    +
    {\bm A}_{\bra{-2,\geq1}}\bra{\Atw}\pare{\eta;x,\xi}
    }
    \delta\Omega\bra{\background}\pare{\eta}
    \\
    &+
    \OpBW{
    \mathcal{P}\pare{\xi}
    A^{T,\eta}_1\bra{
    \tilde{\Omega}\bra{\background}\pare{\eta}
    }\pare{\eta;x,\xi}
    +
    A_{\bra{0,\geq1}}\bra{\background}\pare{\Atw\middle|\eta;x,\xi}
    }
    \eta
    \\
    &+
    {\bm R}^\bR_{\bra{\varrho+1,\geq1}}\pare{\Atw\middle|\eta}
    \delta\Omega\bra{\background}\pare\eta
    +
    R^{\lin,\bR}_{\bra{\varrho+1,\geq1}}\bra{\background}\pare{\Atw\middle|\eta}\eta.
\end{aligned}
\end{equation}
Plugging \cref{eq:PT_constant_background} into \cref{eq:Omega_background_exact_equation} gives the corrected first reduction
\begin{equation}
\label{eq:second_background_iterate}
\begin{aligned}
&\pare{\OpBW{
\Id
+
\Atw\pare{\mathcal{P}\pare{\xi}{\bm A}^{T,\Omega}_{-1}}^\sharp\pare{\eta;x,\xi}
}
+{\bm R}^\bR_{\bra{\varrho+1,\geq1}}\pare{\Atw\middle|\eta}}
\delta\Omega\bra{\background}\pare{\eta}
\\
&=-
\Atw\OpBW{
\mathcal{P}\pare{\xi}
A^{T,\eta}_1\bra{
\tilde{\Omega}\bra{\background}\pare{\eta}
}\pare{\eta;x,\xi}
+
A_{\bra{0,\geq1}}\bra{\background}\pare{\Atw\middle|\eta;x,\xi}
}
\eta
\\
&
+
\Atw\OpBW{
{\bm A}_{\bra{-2,\geq1}}\bra{\Atw}\pare{\eta;x,\xi}
}
\delta\Omega\bra{\background}\pare{\eta}
+
R^{\lin,\bR}_{\bra{\varrho+1,\geq1}}\bra{\background}\pare{\Atw\middle|\eta}\eta .
\end{aligned}
\end{equation}
Define
\begin{equation}
    {\bm P}_{\bra{0}}
    \pare{\Atw\middle|\eta;x,\xi}
    \defeq
    \pare{
    \Id+
    \Atw\pare{\mathcal{P}\pare{\xi}{\bm A}^{T,\Omega}_{-1}}^\sharp\pare{\eta;x,\xi}
    }^{-1}
    \in\pare{\Gamma_{\geq0}^0}^{2\times2}.
\end{equation}
After decreasing $r_0$,
$\sup_{\av\Atw\leq2,x,\xi}
\norm{\Atw\pare{\mathcal P\pare\xi{\bm A}^{T,\Omega}_{-1}}^\sharp}
\leq1/2$. Thus ${\bm P}_{\bra0}-\Id\in\pare{\Gamma^{-1}_{\geq1}}^{2\times2}$
vanishes near $\xi=0$. By \cref{prop:composition_paradiff,lem:composition_rules},
composing \cref{eq:second_background_iterate} on the left with
$\OpBW{{\bm P}_{\bra0}}$ gives
\begin{equation}
\label{eq:first_background_order_one_term}
\begin{aligned}
&\pare{\Id+{\bm R}^{[1],\bR}_{\bra{\varrho+1,\geq1}}\pare{\Atw\middle|\eta}}
\delta\Omega\bra{\background}\pare{\eta}
\\
=&\ 
\OpBW{
A^{S,\eta}_{1,\textnormal{bg}}\bra{\background}
\pare{\Atw\middle|\eta;x,\xi}
+
A_{\bra{0,\geq1}}\bra{\background}\pare{\Atw\middle|\eta;x,\xi}
}
\eta
\\
&\quad\,+
\Atw\OpBW{
{\bm A}_{\bra{-2,\geq1}}^{[1]}
\bra{\Atw}\pare{\eta;x,\xi}
}
\delta\Omega\bra{\background}\pare{\eta}
+
R^{\lin,[1],\bR}_{\bra{\varrho+1,\geq1}}\bra{\background}\pare{\Atw\middle|\eta}\eta,
\end{aligned}
\end{equation}
where ${\bm R}^{[1],\bR}_{\bra{\varrho+1,\geq1}}\in\pare{\cR^{\varrho+1}_{\geq1}}^{2\times2}$ and
\begin{equation}
    A^{S,\eta}_{1,\textnormal{bg}}\bra{\background}
    \pare{\Atw\middle|\eta;x,\xi}
    \defeq
    -\Atw
    \mathcal{P}\pare{\xi}
    A^{T,\eta}_1\bra{
    \tilde{\Omega}\bra{\background}\pare{\eta}
    }
    \pare{\eta;x,\xi}.
\end{equation}
Indeed, the factor ${\bm P}_{\bra{0}}-\Id\in\pare{\Gamma_{\geq1}^{-1}}^{2\times2}$ only produces order-zero corrections when it is composed with the order-one symbol
$\mathcal{P}\pare{\xi}A^{T,\eta}_1\bra{\tilde{\Omega}\bra{\background}\pare{\eta}}$.

We now remove the residual term
$\OpBW{{\bm A}_{\bra{-2,\geq1}}^{[1]}}\delta\Omega\bra{\background}\pare{\eta}$ by the same parametrix iteration used in the Sobolev-input part. We prove by induction that, for every $n\in\bN^*$, there exist
\[{\bm P}_{\bra{-2^n}}\in\pare{\Gamma_{\geq1}^{-2^n}}^{2\times2},\quad
{\bm R}^{[n],\bR}_{\bra{\varrho+1,\geq1}}\in\pare{\cR^{\varrho+1}_{\geq1}}^{2\times2},\]
\[
A_{\bra{0,\geq1}}^{[n]}\bra{\background}\in\cL\pare{\bR^2;\pare{\Gamma_{\geq1}^0}^2},\quad
R^{[n],\bR}_{\bra{\varrho+1,\geq1}}\bra{\background}\in\cL\pare{\bR^2;\pare{\cR_{\geq1}^{\varrho+1}}^2},\] with ${\bm P}_{\bra{-2^n}}=0$ near $\xi=0$, such that
\begin{equation}
\label{eq:higher_background_iterates_induction}
\begin{aligned}
&\pare{\OpBW{
\Id+{\bm P}_{\bra{-2^n}}
\pare{\Atw\middle|\eta;x,\xi}
}
+{\bm R}^{[n],\bR}_{\bra{\varrho+1,\geq1}}\pare{\Atw\middle|\eta}}
\delta\Omega\bra{\background}\pare{\eta}
\\
&\quad=
\OpBW{
A^{S,\eta}_{1,\textnormal{bg}}\bra{\background}
\pare{\Atw\middle|\eta;x,\xi}
+
A_{\bra{0,\geq1}}^{[n]}\bra{\background}
\pare{\Atw\middle|\eta;x,\xi}
}
\eta
+
R^{\lin,[n],\bR}_{\bra{\varrho+1,\geq1}}\bra{\background}\pare{\Atw\middle|\eta}\eta .
\end{aligned}
\end{equation}
The initialization $n=1$ is obtained from \cref{eq:first_background_order_one_term} by moving the residual
$\Atw\OpBW{{\bm A}_{\bra{-2,\geq1}}^{[1]}}\delta\Omega\bra{\background}\pare{\eta}$ to the left. Assume \cref{eq:higher_background_iterates_induction} holds at rank $n$. Decreasing $r_0$ at each of the finitely many steps, \cref{prop:composition_paradiff} gives
\begin{equation}
\OpBW{
\pare{\Id+{\bm P}_{\bra{-2^n}}}^{-1}
}
\OpBW{
\Id+{\bm P}_{\bra{-2^n}}
}
=
\OpBW{
\Id+{\bm P}_{\bra{-2^{n+1}}}
}
+{\bm R}^\bR_{\bra{\varrho+1,\geq1}},
\end{equation}
where ${\bm P}_{\bra{-2^{n+1}}}\in\pare{\Gamma_{\geq1}^{-2^{n+1}}}^{2\times2}$.
Composing \cref{eq:higher_background_iterates_induction} on the left with
$\OpBW{\pare{\Id+{\bm P}_{\bra{-2^n}}}^{-1}}$ proves the induction step: the new defect still vanishes near $\xi=0$, the smoothing terms acting on $\delta\Omega$ remain on the left, and the extra composition of\\
$\pare{\Id+{\bm P}_{\bra{-2^n}}}^{-1}-\Id\in\pare{\Gamma_{\geq1}^{-2^n}}^{2\times2}$ with
$A^{S,\eta}_{1,\textnormal{bg}}\in\pare{\Gamma_{\geq0}^1}^2$ has order at most $0$ for $n\geq1$, hence is absorbed into $A_{\bra{0,\geq1}}^{[n+1]}\bra{\background}$.

Taking $\underline n>\log_2\pare{\varrho+1}$, one has
$\OpBW{{\bm P}_{\bra{-2^{\underline n}}}}\in
\pare{\cR^{\varrho+1}_{\geq1}}^{2\times2}$. The left-hand side of
\cref{eq:higher_background_iterates_induction} becomes
$\pare{\Id+{\bm R}^\bR_{\bra{\varrho+1,\geq1}}\pare{\Atw\middle|\eta}}
\delta\Omega\bra{\background}\pare\eta$.
Invert this factor by \cref{eq:S_smoothing_inverse}. Its smoothing
correction composed with the order-one column belongs to
$\pare{\cR^\varrho_{\geq1}}^2$ by \cref{lem:composition_rules}.
Adding $\background$ gives
\begin{equation}
\label{eq:higher_background_iterates}
   \tilde{\Omega}\bra{\background}\pare{\eta}
    =
    \background
    +
    \OpBW{
    A^{S,\eta}_{1,\textnormal{bg}}\bra{\background}
    \pare{\Atw\middle|\eta;x,\xi}
    +
    A_{\bra{0,\geq1}}\bra{\background}
    \pare{\Atw\middle|\eta;x,\xi}
    }\eta
+R^{\lin,\bR}_{\bra{\varrho,\geq1}}\bra{\background}
    \pare{\Atw\middle|\eta}\eta.
\end{equation}
This proves \cref{eq:paralinearization_S_background}. Finally, using \cref{eq:AT1extra} and
$\mathcal{P}\pare{\xi}\xi^\perp=\xi^\perp$, the principal symbol in \cref{eq:higher_background_iterates} is exactly \cref{eq:ASeta1_background}.
\end{proofpart}

\end{proof}

\subsubsection{Paralinearization of Dirichlet-Neumann-type operators $\mathpzc{g}^\star\p\eta$ and $\mathpzc{g}\pare{\Atw\middle|\eta}$}
Let $\Omega\defeq\vect{\omega_1}{\omega_2}$, where $\omega_1, \omega_2$ are scalar functions that map from $\mathbb R^2$ to $\mathbb R$. Consider the operator $\mathpzc{g}^{\star}\p\eta$ of nonlinear convolution type, that has been introduced in \cref{eq:gstar} and that we recall here for simplicity
\begin{equation}
\label{eq:Gstar_recall}
    \mathpzc{g}^\star\p\eta\bra \Omega =  - {\bm v}^{\star} \bra{\Omega}\pare{\eta} \cdot {\vect{-\nabla\eta}{1}} , 
\end{equation}
where ${\bm v}^\star$ is defined in \cref{eq:bmvtilde}.

The paralinearization of the operator $\mathpzc g^{\star}\pare{\eta}$ is given in the following proposition.

\begin{prop}[Paralinearization of $\mathpzc{g}^\star\pare{\eta}$ ]\label{prop:Gstar_paralinearization}
There exists $\varrho_0 > 0$ such that for any $\varrho \geq \varrho_0$, there exists $s_0 > 0$ and $r_0 > 0$ such that for any $\eta\in \BallR{s_0}{r_0}, \Omega \in\pare{\BallR{s_0}{r_0}}^2$, the operator $\mathpzc{g}^\star\pare{\eta}\bra \Omega$ as in \cref{eq:Gstar_recall} admits the following paralinearization:
\begin{equation}\label{eq:paralinearization_gstar}
    \begin{aligned}
        \mathpzc{g}^\star\pare{\eta}\bra{\Omega}&= \OpBW{ \ii \ A^{\mathpzc{g}^\star,\Omega}_{0}\pare{\eta;x,\xi}{{+A^{\mathpzc{g}^\star,\Omega}_{-1}\pare{\eta;x,\xi}}}+A_{\bra{-2,\geq1}}\pare{\eta;x,\xi}}\cdot \Omega\\
&\quad+\OpBW{ a^{\mathpzc{g}^\star,\eta}_1\bra{\Omega}\pare{\eta;x,\xi}+a_{\bra{0,\geq1}}\pare{\eta,\Omega;x,\xi}}\eta+R^{\lin,\bR}_{\bra{\rho,\geq1}}\p\eta \cdot \Omega+r^{\lin,\bR}_{\bra{\rho,\geq1}}\pare{\eta,\Omega}\eta ,
    \end{aligned}
\end{equation}
where
\begin{itemize}
    \item $A_0^{\mathpzc{g}^\star,\Omega}\in\pare{\Gamma_{\geq0}^0}^2$ is explicitly defined as
        \begin{equation}
        \label{eq:AGstarW0}
         A^{\mathpzc{g}^\star,\Omega}_0\pare{\eta;x,\xi}\defeq \pare{1+\av{\nabla\eta}^2}A_{0}^{\mathpzc{B}_0,\omega;\perp}\pare{\eta;x,\xi} 
         {= - \nabla_\xi^\perp \lambda_1\pare{\eta;x,\xi}}
         ,
        \end{equation}
        where $A_{0}^{\mathpzc{B}_0,\omega}$ is defined in \cref{eq:H0_definition};
\item 
$A_{-1}^{\mathpzc{g}^\star,\Omega}\in\pare{\Gamma_{\geq1}^{-1}}^2$ is explicitly defined as
        \begin{equation}
         \label{eq:AGstarW-1}
         \pare{A^{\mathpzc{g}^\star,\Omega}_{-1}\pare{\eta;x,\xi}}^\intercal\defeq  \frac{1}{2}\pare{\nabla\eta}^\intercal \nabla_\xi^\perp J_\xi \phi^{\mathpzc{B}_0}_1 {\bm H}_x \eta - \frac{1}{4}J_x\av{\nabla\eta}^2 \nabla_\xi J_\xi^\perp \phi^{\mathpzc{B}_0}_1,
        \end{equation}
        where $\phi^{\mathpzc{B}_0}_1$ is defined in \cref{eq:phiB0};
    \item $a^{\mathpzc{g}^\star,\eta}_1\in \cL\pare{\pare{H^{s_0}}^2 ; \Gamma_{\geq1}^1}$ is a symbol $1$-positive homogeneous in $\xi$, admitting the following expression:
        \begin{equation}\label{eq:aGstar_1_corrected} 
           a^{\mathpzc{g}^\star,\eta}_1\bra{\Omega}\pare{\eta;x,\xi}
           \defeq
           \ii \xi\cdot  V^{\star}\bra{\Omega}\pare{\eta;x}
            -\frac{\lambda_1\p{\eta;x,\xi}}{1+\av{\nabla\eta}^2} \ \nabla^\perp\eta\cdot \Omega ,
        \end{equation}
        where $ V^{\star}\bra{\Omega}$ is defined in \cref{eq:bmvtilde_recall};

    \item $a_{\bra{m,\geq1}}\in\Gamma_{\geq1}^m, A_{\bra{m,\geq1}}\in\pare{\Gamma_{\geq1}^m}^2,$ $m\in\bR$, denote generic symbols as per \cref{notation:paralinearization};
    \item $r^{\lin,\bR}_{\bra{\varrho,\geq1}} \in \cR_{\geq1}^\varrho $ and $R^{\lin,\bR}_{\bra{\varrho,\geq1}} \in\pare{\cR_{\geq1}^\varrho}^2$ denote generic smoothing remainders as per \cref{notation:paralinearization}.
\end{itemize}
\end{prop}

\begin{proof}
% \begin{enumerate}
The part $-v_3\bra{\Omega}\pare{\eta}$ has already been analyzed in \cref{eq:paralinearization_tildev3}, so we only need to paralinearize the term 
\begin{equation}
 V^{\star}\bra{\Omega}\pare{\eta}\cdot \nabla \eta
\end{equation}
in \cref{eq:Gstar_recall}. Indeed \cref{lem:paraprod_Bony} give us that
\begin{equation}\label{eq:Gstar1}
\displaystyle V^{\star}\bra{\Omega}\pare{\eta}\cdot \nabla\eta=\OpBW{ V^{\star}\bra{\Omega}\pare{\eta}}\cdot\nabla\eta+
     \OpBW{\nabla\eta}\cdot \bra{ V^{\star}\bra{\Omega}\pare{\eta}}
     +
     \sum_{i=1}^2\RBW{ V^{\star}\bra{\Omega}\pare{\eta}_i}{\partial_i\eta}.
\end{equation}
First, applying \cref{prop:composition_paradiff}, we have
\begin{equation}
    \label{eq:GBony1}
\OpBW{ V^{\star}\bra{\Omega}\pare{\eta}}\cdot\nabla\eta=\OpBW{\ii\xi\cdot   V^{\star}\bra{\Omega}\pare{\eta;x}  +a_{\bra{0,\geq1}}\bra{\Omega}\pare{\eta;x,\xi}}\eta
+r^{\lin,\bR}_{\bra{\varrho,\geq1}}\pare{\Omega}\pare{\eta}\eta.
\end{equation}
We use now \cref{eq:paralinearization_Vtilde,prop:composition_paradiff} and obtain that
\begin{equation}\label{eq:toughcal7changed}
\begin{aligned}
    &\OpBW{\nabla\eta}\cdot \bra{ V^{\star}\bra{\Omega}\pare{\eta}}
    \\
    &=
    \OpBW{\ii \ \nabla\eta ^\intercal {\bm A}^{ V^{\star},\Omega}_0\pare{\eta;x,\xi} +\frac{1}{2}\pbra{\nabla\eta ^\intercal}{ {\bm A}^{ V^{\star},\Omega}_0\pare{\eta;x,\xi} } +\ \nabla\eta ^\intercal {\bm A}^{ V^{\star},\Omega}_{-1}\pare{\eta;x,\xi} + A_{\bra{-2,\geq1}} \pare{\eta,\Omega;x,\xi}}\Omega
    \\
    &\quad+
    \OpBW{\nabla\eta^\intercal {\bm A}^{ V^{\star},\eta}_1\pare{\eta;x,\xi}\Omega + a_{\bra{0,\geq1}}\pare{\eta,\Omega;x,\xi} }\eta 
    + R^{\lin,\bR}_{\bra{\varrho,\geq1}} \pare{\eta}\cdot \Omega + r^{\lin,\bR}_{\bra{\varrho,\geq1}}\pare{\eta, \Omega}\eta .
    \end{aligned}
\end{equation}
Next, using \cref{rem:smoothBONYremainder,eq:paralinearization_tildev3,lem:composition_rules} we obtain that
\begin{equation}\label{eq:Gstar2}
    \sum_{i=1}^2\RBW{ V^{\star}\bra{\Omega}\pare{\eta}_i}{\partial_i\eta} = r^{\lin,\bR}_{\bra{\varrho,\geq1}}\pare{\eta, \Omega} \eta.
\end{equation}
We now plug the results of \cref{eq:Gstar2,eq:toughcal7changed,eq:paralinearization_tildev3,eq:GBony1} in \cref{eq:Gstar1,eq:Gstar_recall} 
and obtain the paradifferential decomposition stated in \cref{eq:paralinearization_gstar} where the symbols $A^{\mathpzc{g}^\star,\Omega}_0$ and $A^{\mathpzc{g}^\star,\eta}_1$ assume the explicit form

\begin{subequations}
\begin{align}\label{eq:Gstar4}
       & \pare{A^{\mathpzc{g}^\star,\Omega}_{0}\pare{\eta;x,\xi}}^\intercal \defeq  \  A_{0}^{\mathpzc{B}_0,\omega;\perp}\pare{\eta;x,\xi} + \nabla \eta ^\intercal {\bm A}^{ V^{\star},\Omega}_0\pare{\eta;x,\xi},
      \\
      %---------------------------------------
       \label{eq:Gstar4.1}
         & \pare{A^{\mathpzc{g}^\star,\Omega}_{-1}\pare{\eta;x,\xi}}^\intercal \defeq %\pare{ \p{A^{\mathpzc B_0,\omega}_{-1}\pare{\eta;x,\xi}}^\perp}^\intercal +
         \nabla \eta ^\intercal {\bm A}^{ V^{\star},\Omega}_{-1}\pare{\eta;x,\xi}+\frac{1}{2}\pbra{\nabla\eta ^\intercal}{ {\bm A}^{ V^{\star},\Omega}_0\pare{\eta;x,\xi} } ,
      \\
      %---------------------------------------
       \label{eq:Gstar4.1-1}&a^{\mathpzc{g}^\star,\eta}_{1}\bra{\Omega}\pare{\eta;x,\xi} \defeq \  \ii \xi\cdot  V^{\star}\bra{\Omega}\pare{\eta;x}
         + 
          \underbrace{\nabla \eta ^\intercal {\bm A}^{ V^{\star},\eta}_1\pare{\eta;x,\xi} \Omega + \pare{A^{\mathpzc{B}_0,\eta}_1\pare{\eta;x,\xi}}^\perp \cdot \Omega}_{\defeq B^{\mathpzc{g}^\star,\eta}_1\pare{\eta;x,\xi}\cdot \Omega}.
\end{align}
\end{subequations}
Let us notice that using \cref{eq:Pbra_matrix-vector,eq:AVtildeg0_lem,eq:AVtildeg-1_lem} we obtain that
\begin{small}
    \begin{equation}\label{eq:simplification_-1_gstar}
        \begin{aligned}
        & 
        \frac{1}{2}\pbra{\nabla\eta ^\intercal}{ {\bm A}^{ V^{\star},\Omega}_0 }  = \bra{-\frac{1}{2} \partial_{x_ix_j} \eta \ \partial_{\xi_j} a^{ V^{\star},\Omega}_{0;ik}}_{k=1,2}
        =
        \bra{-\frac{1}{2} \partial_{x_ix_j} \eta \partial_{x_i}\eta \ \partial_{\xi_j} A_{0;k}^{\mathpzc{B}_0,\omega;\perp} }_{k=1,2} ^\intercal
        =
        -\frac{1}{4} J_x \av{\nabla_x\eta}^2 \nabla_\xi \pare{A_{0}^{\mathpzc{B}_0,\omega;\perp}}^\intercal,
        \\
        %-----------------------------------------------------------
        &
        \nabla \eta ^\intercal {\bm A}^{ V^{\star},\Omega}_{-1}
        = \frac{1}{2} \pare{\nabla\eta}^\intercal J_\xi A_{0}^{\mathpzc{B}_0,\omega;\perp} {\bm H}_x \eta .
    \end{aligned}
    \end{equation}
\end{small}
Thus, inserting \cref{eq:simplification_-1_gstar} in \cref{eq:Gstar4.1} proves \cref{eq:AGstarW-1}, and we obtain \cref{eq:aGstar_1_corrected} from \cref{eq:Gstar4.1-1}.
It remains hence only to prove that \cref{eq:Gstar4,eq:Gstar4.1} can be expressed as in \cref{eq:AGstarW0,eq:aGstar_1_corrected} respectively.
\\
The equivalence of \cref{eq:AGstarW0} and \cref{eq:Gstar4} follows from \cref{eq:AVtildeg0_lem,eq:Gstar4}, while using \cref{eq:AVtildeeta1_lem,eq:AB0eta1,eq:relationimportant3} we obtain that
\begin{equation}
    B^{\mathpzc{g}^\star,\eta}_1\pare{\eta;x,\xi}
    =
    -\frac{\lambda_1\p{\eta;x,\xi}}{1+\av{\nabla\eta}^2} \ \nabla^\perp \eta,
\end{equation}
and so we conclude.

\end{proof}

Let us recall that the operator $\mathpzc{g}\pare{\Atw\middle|\eta}\psi$ is defined through \cref{eq:DN,eq:bmv,eq:gstar} as:
\begin{equation}\label{eq:DN_recall}
         \mathpzc{g}\pare{\Atw\middle|\eta}\psi \defeq \mathpzc{g}^\star \pare{\eta}\bra{S\pare{\Atw\middle|\eta}\bra\psi } =
        - {\bm v}\bra{\psi}\pare{\Atw\middle|\eta}\cdot \vect{-\nabla\eta}{1}
           ,
\end{equation}
where $\mathpzc{g}^\star\pare{\eta}$ is defined in \cref{eq:Gstar_recall} and $S\pare{\Atw\middle|\eta}$ is defined in \cref{eq:S_recall}.

\begin{prop}[Paralinearization of the Dirichlet--Neumann operator $\mathpzc{g}\pare{\Atw\middle|\eta}\psi$]\label{prop:G_paralinearization}
 There exists $\varrho_0 > 0$ such that for any $\varrho \geq \varrho_0$, there exists $s_0 > 0$ and $r_0 > 0$ such that for any $\eta,\psi\in \BallR{s_0}{r_0}$, the operator $\mathpzc{g} \pare{\Atw\middle| \eta}\psi$ defined as in \cref{eq:DN_recall} admits the following paralinearization:
\begin{equation}\label{eq:paraoperatorG}
                    \mathpzc{g}\pare{\Atw\middle|\eta}\psi =\OpBW{ A^{\mathpzc{g}}_1\pare{\Atw \middle| \eta, \psi; x, \xi} + 
                    \vect{a_{\bra{0,\geq1}}\bra{\psi}\pare{\eta;x, \xi}}{a_{\bra{-1,\geq1}}\pare{\eta;x, \xi}}}\cdot\vect{\eta}{\psi} + R^{\lin,\bR}_{\bra{\varrho,\geq1}}\pare{\eta, \psi}\cdot\vect{\eta}{\psi},
\end{equation}
    where
\begin{itemize}
\item $\Atw$ is the Atwood number defined in \cref{eq:Atwood};
        \item $A^{\mathpzc{g}}_1 \in \pare{\Gamma_{\geq0}^1}^2$ is explicitly defined as
                 \begin{equation}\label{eq:AG1}
                    A^{\mathpzc{g}}_1\pare{\Atw \middle| \eta, \psi; x, \xi} \defeq \vect{a^{\mathpzc{g}, \eta}_1\bra{\psi}\pare{\Atw\middle|\eta;x,\xi}}{\lambda_1\pare{\eta;x,\xi} + \Atw \lambda_0\pare{\eta;x,\xi}},
                 \end{equation}
                 where
                \begin{itemize}
                     \item $\lambda_1\in\Gamma_{\geq0}^1$ and   $\lambda_0 \in \Gamma_{\geq1}^0$ 
                     are defined in \cref{eq:lambda1} and \cref{eq:lambda0} respectively;
                      %\item $\lambda_0 \in \Gamma_{\geq1}^0$ is defined in \cref{eq:lambda0}
               
                 \item $a^{\mathpzc{g}, \eta}_1\in\cL\pare{H^{s_0};\Gamma_{\geq0}^1}$ is a symbol $1$-positive homogeneous in $\xi$, admitting the following expression:
\begin{equation}\label{eq:Getafinal2}
a^{\mathpzc{g}, \eta}_1\bra{\psi}\pare{\Atw \middle| \eta;x,\xi} \defeq
 \lambda_1\pare{\eta;x,\xi}b\bra{\psi}\pare{\Atw\middle|\eta;x} + \ii \mathpzc{V}\bra{\psi}\pare{\Atw\middle|\eta;x}\cdot\xi ,
\end{equation}
where $\lambda_1$ is defined in \cref{eq:lambda1}.  The coefficient
$\mathpzc V\bra{\psi}\pare{\Atw\middle|\eta;x}$ is real-valued,
$\xi$-independent, and of order zero; it is the product
$\ii\mathpzc V\bra{\psi}\pare{\Atw\middle|\eta;x}\cdot\xi$ that has order one.
                \begin{itemize}
                    \item $b\bra\psi\in
                    \cL\pare{H^{s_0};\cF^\bR}$ is the real-valued
                    coefficient defined in
                    \cref{eq:definitionBb_pi_section};
                    \item $\mathpzc V\bra\psi\pare{\Atw\middle|\eta}\in
                    \cL\pare{H^{s_0};\pare{\cF^\bR}^2}$ is the
                    real-valued coefficient defined in
                    \cref{eq:definitionVb_pi_section}.
                \end{itemize}
                 \end{itemize}

\item $a_{\bra{0,\geq1}}\in {\Gamma_{\geq1}^0}, a_{\bra{-1,\geq1}}\in {\Gamma_{\geq1}^{-1}}$ denote  generic symbols as per \cref{notation:paralinearization};
 \item $R^{\lin,\bR}_{\bra{\varrho,\geq1}} \in\pare{\cR_{\geq1}^\varrho}^2$ denotes a generic smoothing remainder as per \cref{notation:paralinearization}.
\end{itemize}

\end{prop}

\begin{proof}
Choose the baseline regularity one derivative above the threshold in
\cref{prop:Gstar_paralinearization}, so that $S\psi$ has the required input
regularity. Applying that proposition with smoothing index $\varrho+1$ and
the composition rules stated in \cref{lem:composition_rules}, we infer that
\begin{equation}\label{Gstar of S}
    \begin{aligned}
        \mathpzc{g}^\star\pare\eta S\pare{\Atw\middle|\eta}\psi &= \OpBW{ \ii A^{\mathpzc{g}^\star,\Omega}_0\pare{\eta;x, \xi}{{+A^{\mathpzc{g}^\star,\Omega}_{-1}\pare{\eta;x, \xi}} }+ A_{\bra{-2,\geq1}}\pare{\eta;x, \xi}} \cdot\bra{S\pare{\Atw\middle|\eta}\psi} 
     \\
     &\quad+
     \OpBW{ a^{\mathpzc{g}^\star, \eta}_1 \bra{S\pare{\Atw\middle|\eta}\psi} \pare{\eta;x, \xi} + a_{\bra{0,\geq1}}\pare{\eta, \psi;x, \xi}}\eta
     + R^{\lin,\bR}_{\bra{\varrho+1,\geq1}}\pare{\eta, \psi}\cdot\vect{\eta}{\psi}.
    \end{aligned}
\end{equation}
Next, we combine \cref{eq:paralinearization_S} and \cref{Gstar of S} together with the composition rule \cref{prop:composition_paradiff}. Since $S$ has order one, this loses one derivative on the preceding remainder; after relabeling $\varrho+1$ as $\varrho$, we obtain
\begin{align}
 \mathpzc{g}^\star\pare\eta S\pare{\Atw\middle| \eta}\psi &=\OpBW{ \lambda_1\pare{\eta;x,\xi} +  \Lambda_0^{\mathrm{DN}}\pare{\eta;x,\xi} +a_{\bra{-1,\geq1}}\pare{\eta;x,\xi}}\psi\\
&\quad+\OpBW{a^{\mathpzc{g}, \eta}_1\bra{\psi}\pare{\Atw| \eta;x,\xi}+a_{\bra{0,\geq1}}\pare{\eta,\psi;x,\xi}}\eta+ R^{\lin,\bR}_{\bra{\varrho,\geq1}}\pare{\eta, \psi}\cdot\vect{\eta}{\psi},
\end{align}
where
\begin{subequations}
    \begin{align}
    \label{eq:lambda1_wip}
 \lambda_1\pare{\eta;x,\xi}
 \defeq & \ 
 - \xi^\perp\cdot A^{\mathpzc{g}^\star,\Omega}_0\pare{\eta;x,\xi},
 \\
  \label{eq:Lambda0DN_wip}
\Lambda_0^{\mathrm{DN}}\pare{\eta;x,\xi} \defeq & \ 
\Atw A^{\mathpzc{g}^\star, \Omega}_0\pare{\eta;x,\xi} \cdot \xi^\perp \ a^{\tau_0, \psi}_{-1}\pare{\eta;x,\xi} +
\ii \pare{\frac 1{2}\pbra{ \pare{A^{\mathpzc{g}^\star,\Omega}_0}^\intercal }{\xi^\perp }+ A^{\mathpzc{g}^\star,\Omega}_{-1}\pare{\eta;x, \xi}\cdot\xi^\perp },\qquad\\
a^{\mathpzc{g}, \eta}_1\bra{\psi}\pare{\Atw| \eta;x,\xi}\defeq & \ a^{\mathpzc{g}^\star, \eta}_1 \bra{S\pare{\Atw\middle|\eta}\psi} \pare{\eta;x, \xi}
+\ii A^{\mathpzc{g}^\star,\Omega}_0\pare{\eta;x,\xi}\cdot A^{S,\eta}_1\bra{\psi}\pare{\Atw\middle|\eta;x,\xi}.\label{eq:a1g_wip}
\end{align}
\end{subequations}
We now simplify the expressions in \cref{eq:Lambda0DN_wip,eq:lambda1_wip}. 

\begin{proofpart}[Equivalence of \cref{eq:lambda1_wip,eq:lambda1}]
    From \cref{eq:AGstarW0} and \cref{eq:relationimportant3}, we have   
\begin{equation}
\label{eq:Getafinal21}
\xi^\perp\cdot A^{\mathpzc{g}^\star,\Omega}_0\pare{\eta;x,\xi}=\p{1+\av{\nabla\eta}^2}\p{A_{0}^{\mathpzc{B}_0,\omega}\p{\eta;x,\xi}\cdot\xi}=-\lambda_1\p{\eta;x,\xi}. 
\end{equation}
This justifies a posteriori the choice of notation.
\end{proofpart}

\begin{proofpart}[Identification of $\Lambda_0^{\mathrm{DN}}=\Atw\lambda_0$]
We decompose this step in two subparts.
\begin{proofsubpart}[Proof that $\Im \Lambda_0^{\mathrm{DN}}\equiv 0$]\label{part:ImLambda0DN=0}
      We decompose $A^{\mathpzc{g}^\star,\Omega}_{-1}\pare{\eta;x, \xi}\cdot\xi^\perp$
using \cref{eq:AGstarW-1} as
\begin{subequations}
    \begin{align}
&A^{\mathpzc{g}^\star,\Omega}_{-1}\cdot\xi^\perp
= \  b_0\pare{\eta;x,\xi} + c_0\pare{\eta;x,\xi},
\label{eq:b0+c0}
\\
    &b_0\pare{\eta;x,\xi}
    \defeq
    \frac{1}{2}\pare{\nabla\eta}^\intercal
    \pare{\nabla_\xi^\perp J_\xi \phi^{\mathpzc{B}_0}_1}
    {\bm H}_x\eta\;\xi^\perp,
    \label{eq:b0_def}
    \\
    &c_0\pare{\eta;x,\xi}
    \defeq
    -\frac{1}{4}J_x\av{\nabla\eta}^2\;
    \pare{\nabla_\xi J_\xi^\perp \phi^{\mathpzc{B}_0}_1}\;\xi^\perp.
    \label{eq:c0_def}
\end{align}
\end{subequations}

We show that $c_0 = 0$ and that $b_0$ cancels
$\frac{1}{2}\pbra{\pare{A^{\mathpzc{g}^\star,\Omega}_0}^\intercal}{\xi^\perp}$.

\begin{step}[$c_0 = 0$]
A componentwise expansion gives for any $i=1,2$
\begin{equation}
    \bra{\pare{\nabla_\xi J_\xi^\perp \phi^{\mathpzc{B}_0}_1}\;\xi^\perp}_i
    =
    -\partial_{\xi_i}\partial_{\xi_2}\phi^{\mathpzc{B}_0}_1\;\pare{-\xi_2}
    +\partial_{\xi_i}\partial_{\xi_1}\phi^{\mathpzc{B}_0}_1\;\xi_1
    =
    \bra{J_\xi A_{0}^{\mathpzc{B}_0,\omega}\cdot\xi}_i.
\end{equation}
Since $A_{0}^{\mathpzc{B}_0,\omega} = \nabla_\xi\phi^{\mathpzc{B}_0}_1$ is $0$-homogeneous in $\xi$,
Euler's theorem yields
$J_\xi A_{0}^{\mathpzc{B}_0,\omega}\cdot\xi = 0$,
whence $c_0 = 0$.
\end{step}

\begin{step}[Cancellation of $b_0$ with the Poisson bracket]
From \cref{eq:AGstarW0} and \cref{eq:Pbra}, the Poisson bracket contribution reads
\begin{equation}\label{eq:PBterm_div}
 \frac 1{2}\pbra{ \pare{A^{\mathpzc{g}^\star,\Omega}_0}^\intercal }{\xi^\perp } = -\frac{1}{2} \nabla^\perp \cdot A^{\mathpzc{g}^\star,\Omega}_0 = -\frac{1}{2} \nabla_x \cdot \pare{\pare{1+\av{\nabla\eta}^2} A_{0}^{\mathpzc{B}_0,\omega}}
    \xlongequal{\cref{eq:phiB0}} \frac{\nabla_x\cdot\nabla_\xi}{2} \lambda_1.
\end{equation}
Moreover, since $\av{\nabla\eta}^2$ is independent of $\xi$,
\begin{align}\label{eq:scaling_H0_lambda}
    % \pare{1+\av{\nabla\eta}^2} A_{0}^{\mathpzc{B}_0,\omega}
    % =
    % -\nabla_\xi\lambda_1, &&
    \pare{1+\av{\nabla\eta}^2} J_\xi A_{0}^{\mathpzc{B}_0,\omega}
    =
    -{\bm H}_\xi\lambda_1.
\end{align}

\begin{substep}[Reduction of $b_0$]\label{substep:b0_reduction}
A direct matrix computation shows
$\nabla_\xi^\perp J_\xi\phi^{\mathpzc{B}_0}_1
= {\bm J}\,J_\xi A_{0}^{\mathpzc{B}_0,\omega}$,
where ${\bm J}\defeq\begin{psmallmatrix}0&-1\\1&0\end{psmallmatrix}$
is the standard symplectic matrix.
Since $\pare{\nabla\eta}^\intercal{\bm J} = -\pare{\nabla^\perp\!\eta}^\intercal$,
one obtains
\begin{equation}\label{eq:algebraic_reduction_b0}
    \pare{\nabla\eta}^\intercal
    \pare{\nabla_\xi^\perp J_\xi \phi^{\mathpzc{B}_0}_1}
    =
    -\pare{\nabla^\perp\!\eta}^\intercal J_\xi A_{0}^{\mathpzc{B}_0,\omega}
    \xlongequal{\cref{eq:scaling_H0_lambda}}
    \frac{\pare{\nabla^\perp\!\eta}^\intercal{\bm H}_\xi\lambda_1}
         {1+\av{\nabla\eta}^2},
\end{equation}
whence
\begin{equation}\label{eq:b0_via_lambda}
    b_0
    =
    \frac{\pare{\nabla^\perp\!\eta}^\intercal
          {\bm H}_\xi\lambda_1}
         {2\pare{1+\av{\nabla\eta}^2}}
    \cdot {\bm H}_x\eta\;\xi^\perp.
\end{equation}
We compute $\pare{\nabla^\perp\!\eta}^\intercal{\bm H}_\xi\lambda_1$
explicitly.
From \cref{eq:lambda1},
$4\lambda_1^2 = \av{\xi}^2
+ \pare{\nabla^\perp\!\eta\cdot\xi}^2$,
which gives
\[4\lambda_1\,\partial_{\xi_k}\lambda_1
= \xi_k + \pare{\nabla^\perp\!\eta\cdot\xi}\pare{\nabla^\perp\!\eta}_k.\]
Differentiating once more in $\xi_i$ and contracting with
$\pare{\nabla^\perp\!\eta}_i$, we find, using
$\nabla^\perp\!\eta\cdot\xi^\perp = \nabla\eta\cdot\xi$
and
$\av{\nabla^\perp\!\eta}^2 = \av{\nabla\eta}^2$,
\begin{equation}\label{eq:perp_contraction_Hess}
    \pare{\nabla^\perp\!\eta}^\intercal{\bm H}_\xi\lambda_1
    =
    \frac{1+\av{\nabla\eta}^2}{16\lambda_1^3}
    \pare{\av{\xi}^2\pare{\nabla^\perp\!\eta}
    - \pare{\nabla^\perp\!\eta\cdot\xi}\,\xi}^\intercal.
\end{equation}
Substituting \cref{eq:perp_contraction_Hess} into
\cref{eq:b0_via_lambda}, we obtain
\begin{equation}\label{eq:b0_explicit}
    b_0
    =
    \frac{1}{32\lambda_1^3}
    \pare{\av{\xi}^2\pare{\nabla^\perp\!\eta}
    - \pare{\nabla^\perp\!\eta\cdot\xi}\,\xi}^\intercal
    {\bm H}_x\eta\;\xi^\perp
    \xlongequal{{\bm H}_x\eta = \pare{{\bm H}_x\eta}^\intercal}
    \frac{\pare{\xi^\perp}^\intercal{\bm H}_x\eta
    \pare{\av{\xi}^2\nabla^\perp\!\eta
    - \pare{\nabla^\perp\!\eta\cdot\xi}\,\xi}}{32\lambda_1^3}\cdot
\end{equation}
\end{substep}

\begin{substep}[Reduction of the Poisson bracket]\label{substep:PB_reduction}
From $4\lambda_1^2 = \av{\xi}^2
+ \pare{\nabla^\perp\!\eta\cdot\xi}^2$
and the identities
\[\nabla_{\nabla\eta}\pare{\nabla^\perp\!\eta\cdot\xi} = -\xi^\perp,
\quad\nabla_\xi\pare{\nabla^\perp\!\eta\cdot\xi} = \nabla^\perp\!\eta,\]
we compute
\begin{equation}
    \partial_{\pare{\nabla_j\eta}}\partial_{\xi_k}\lambda_1
    =
    \frac{-\pare{\xi^\perp}_j}{16\lambda_1^3}
    \pare{\av{\xi}^2\pare{\nabla^\perp\!\eta}_k
    - \pare{\nabla^\perp\!\eta\cdot\xi}\,\xi_k}
    -
    \frac{\nabla^\perp\!\eta\cdot\xi}{4\lambda_1}\,J_{jk},
\end{equation}
where $J_{jk}$ are the entries of ${\bm J}$.
Applying the chain rule
\[\frac{\nabla_x\cdot\nabla_\xi}{2}\lambda_1
= \frac{1}{2}\sum_{j,k}\pare{{\bm H}_x\eta}_{jk}\,
\partial_{\pare{\nabla_j\eta}}\partial_{\xi_k}\lambda_1,\]
the antisymmetric term vanishes by the symmetry of ${\bm H}_x\eta$:
\begin{equation}
    \sum_{j,k}\pare{{\bm H}_x\eta}_{jk}\,J_{jk}
    =
    \operatorname{tr}\pare{{\bm H}_x\eta\,{\bm J}^\intercal}
    =
    \pare{{\bm H}_x\eta}_{12}
    -
    \pare{{\bm H}_x\eta}_{21}
    =
    0,
\end{equation}
so that
\begin{equation}\label{eq:PB_explicit}
    \frac{\nabla_x\cdot\nabla_\xi}{2}\lambda_1
    =
    -\frac{\pare{\xi^\perp}^\intercal{\bm H}_x\eta
    \pare{\av{\xi}^2\nabla^\perp\!\eta
    - \pare{\nabla^\perp\!\eta\cdot\xi}\,\xi}}{32\lambda_1^3}\cdot
\end{equation}
\end{substep}

\begin{step}[Synthesis]
    Combining \cref{eq:Lambda0DN_wip,eq:b0_explicit,eq:PB_explicit,eq:b0+c0,eq:PBterm_div}, we conclude
\begin{equation}\label{eq:agpsi0_zero}
   \Im \Lambda_0^{\mathrm{DN}}
    =
    \frac{\nabla_x\cdot\nabla_\xi}{2}\lambda_1
    + b_0
    =
    0.
\end{equation}
\end{step}

\end{step}
\end{proofsubpart}

\begin{proofpart}[Simplification of $\Re \Lambda_0^{\mathrm{DN}}$]
Considering the result proved in \cref{part:ImLambda0DN=0} we have that
\begin{equation}
    \Lambda_0^{\mathrm{DN}} =\Re \Lambda_0^{\mathrm{DN}} = \Atw A^{\mathpzc{g}^\star, \Omega}_0 \cdot \xi^\perp \ a^{\tau_0, \psi}_{-1}.
\end{equation}
We can now apply \cref{eq:AGstarW0,eq:relationimportant3,eq:atau0} and obtain
$\Lambda_0^{\mathrm{DN}}=\Atw\lambda_0$, with $\lambda_0$ defined in
\cref{eq:lambda0}.
    
\end{proofpart}

\begin{proofpart}[Proof of \cref{eq:Getafinal2}]
Using \cref{eq:aGstar_1_corrected,eq:bmv}, we obtain\begin{equation}
\begin{aligned}
\label{eq:Getafinal22}
a^{\mathpzc{g}^\star, \eta}_1&=\ii V\bra\psi\p\eta\cdot\xi-\lambda_1\p{\eta;x,\xi}\frac{\nabla^\perp\eta\cdot S\pare{\Atw\middle|\eta}\psi}{1+\av{\nabla\eta}^2}\cdot
\end{aligned}
\end{equation}
Moreover, \cref{eq:ASeta1_sobolev,eq:aSeta0,eq:Getafinal21} give the
additional contribution from the composition of the $\Omega$-column of
$\mathpzc g^\star$ with the $\eta$-column of $S$:
\begin{equation}
\label{eq:g_eta_missing_composition}
\begin{aligned}
\ii A^{\mathpzc{g}^\star,\Omega}_0\cdot A^{S,\eta}_1\bra{\psi}
&=
\Atw\lambda_1 v_3\bra{\psi}
+\ii\Atw\lambda_1
A_{0}^{\mathpzc{B}_0,\omega;\perp}\cdot S\pare{\Atw\middle|\eta}\psi.
\end{aligned}
\end{equation}
We plug \cref{eq:Getafinal21,eq:Getafinal22,eq:g_eta_missing_composition}
into \cref{eq:a1g_wip}; using \cref{eq:bmv} yields exactly
\cref{eq:Getafinal2}.
\end{proofpart}

\end{proofpart}

\end{proof}

\subsubsection{Paralinearization of $\pi^\star$}

Let $\pi^\star$ be defined as in \cref{eq:pistar}, i.e.
\begin{equation}\label{eq:pistar_recall}
            \pi^\star\pare{\rho^\pm\middle|\eta,\Omega} \defeq \pi_0\pare{\rho^\pm\middle|\eta} + \pi_2^\star\pare{\Atw\middle|\eta}\bra{\Omega,\Omega},
        \end{equation}
        where $\pi_0$ and $\pi^\star_2$ are respectively defined in \cref{eq:pi0,eq:pistar2}. The result we prove is the following 

     \begin{prop}[Paralinearization of $\pi^\star$]\label{prop:paralinearization_pistar}
Let $\sigma\geq0$, $\rho^\pm\geq0$ with $\aangle{\rho^\pm}>0$, and
$\mathfrak g\in\bR$. There exists $\varrho_0 > 0$ such that for any
$\varrho \geq \varrho_0$, there exist
$s_0>0$ and $r_0>0$ such that for any
$\eta\in \BallR{s_0}{r_0}$ and
$\Omega\in \pare{\BallR{s_0}{r_0}}^2$, the function
$\pi^\star\pare{\rho^\pm\middle|\eta,\Omega}$ defined in
\cref{eq:pistar_recall} admits the following paralinearization:
\begin{equation}\label{eq:paralinearization_pistar}
    \begin{aligned}\pi^\star\pare{\rho^\pm\middle|\eta,\Omega}
    =&
    -\mathfrak{g}\Atw\,\eta
    +
    \OpBW{
        {\bm a}^{\pi^\star}_2
        \pare{\rho^\pm\middle|\eta;x,\xi}
        +
        {\bm a}^{\pi^\star}_1
        \pare{\Atw\middle|\eta,\Omega;x,\xi}
        +
        \vect{
            a_{\bra{0,\geq1}}\pare{\sigma,\rho^\pm,\Atw\middle|\eta,\Omega;x,\xi}
        }{
            A_{\bra{-1,\geq1}}\bra{\Omega}\pare{\Atw\middle|\eta;x,\xi}
        }
    }\cdot\vect{\eta}{\Omega}
    \\
    &+
    r^{\pi_0,\bR}_{\bra{\varrho,\geq1}}
    \pare{\sigma,\rho^\pm\middle|\eta}\eta
    +
    {\bm r}^{\bil,\bR}_{\bra{\varrho,\geq1}}
    \pare{\Atw\middle|\eta,\Omega}
    \cdot\vect{\eta}{\Omega},
\end{aligned}
\end{equation}
where
\begin{itemize}
    \item ${\bm a}^{\pi^\star}_2\in\pare{\Gamma^2_{\geq0}}^3$ and is explicitly defined as
    \begin{equation}\label{eq:bmapistar2_order_two}
        {\bm a}^{\pi^\star}_2
        \pare{\rho^\pm\middle|\eta;x,\xi}
        \defeq
        \vect{
            a^{\pi^\star,\eta}_2\pare{\rho^\pm\middle|\eta;x,\xi}
        }{
            0_{\bR^2}
        },
    \end{equation}
    where
     $a^{\pi^\star,\eta}_2\in \Gamma^2_{\geq0}$ is real-valued and is explicitly defined as
        \begin{equation}\label{eq:apistar_eta2}
            a^{\pi^\star,\eta}_2\pare{\rho^\pm\middle|\eta;x,\xi}
            \defeq
            -\frac{\sigma}{\aangle{\rho^\pm}}
            \frac{
                \av{\xi}^2+\av{\nabla^\perp\eta\cdot\xi}^2
            }{
                \pare{1+\av{\nabla\eta}^2}^{\frac32}
            } .
        \end{equation}

    \item ${\bm a}^{\pi^\star}_1\in\pare{\Gamma_{\geq1}^1}^3$ and is explicitly defined as
    \begin{equation}\label{eq:bmapistar1}
       {\bm a}^{\pi^\star}_1
        \pare{\Atw\middle|\eta,\Omega;x,\xi}
        \defeq
        \vect{
            a^{\pi^\star,\eta}_1\bra{\Omega,\Omega}
            \pare{\Atw\middle|\eta;x,\xi}
        }{
            A^{\pi^\star,\Omega}_0\bra{\Omega}
            \pare{\Atw\middle|\eta;x,\xi}
        },
    \end{equation}
    where
    \begin{itemize}
        \item $a^{\pi^\star,\eta}_1
        \in \cL^2\pare{\pare{H^{s_0}}^2;\Gamma^1}$ and is explicitly defined as
        \begin{equation}\label{eq:apistar_eta1}
    \begin{aligned}&a^{\pi^\star,\eta}_1\bra{\Omega,\Omega}\pare{\Atw\middle|\eta;x,\xi}\\
      \defeq \ 
&
\lambda_1\pare{\eta;x,\xi}
\left[
    4\pare{A_{0}^{\mathpzc{B}_0,\omega;\perp}\pare{\eta;x,\xi}\cdot\Omega}^2
    -
    \left(
        \frac{\nabla^\perp\eta\cdot\Omega}
        {1+\av{\nabla\eta}^2}
    \right)^2
\right]
\\
+&\Atw\,
\pare{A_{0}^{\mathpzc{B}_0,\omega;\perp}\pare{\eta;x,\xi}\cdot\Omega}
\,
 V^{\star}\bra{\Omega}\pare{\eta;x}\cdot\xi
\\
+&\Atw\,\lambda_1\pare{\eta;x,\xi}
\, v_3^{\star}\bra{\Omega}\pare{\eta;x}
\frac{\nabla^\perp\eta\cdot\Omega}
{1+\av{\nabla\eta}^2}
\\
-&\ii\Atw
\left[
     v_3^{\star}\bra{\Omega}\pare{\eta;x}
     V^{\star}\bra{\Omega}\pare{\eta;x}
     +
\frac{1}{4\pare{1+\av{\nabla\eta}^2}}
    \left(
        \pare{\Omega\cdot\nabla\eta}\Omega
        -
        \frac{\av{\Omega}^2+\av{\Omega\cdot\nabla\eta}^2}
        {1+\av{\nabla\eta}^2}
        \nabla\eta
    \right)
\right]\cdot\xi ;
\end{aligned}
        \end{equation}
        where $\lambda_1 \in \Gamma^1_{\geq0}$ is defined in \cref{eq:lambda1},
        $A_{0}^{\mathpzc{B}_0,\omega} \in \pare{\Gamma^0_{\geq0}}^2$ and $A_{0}^{\mathpzc{B}_0,\omega;\perp} \defeq {\bm J}A_{0}^{\mathpzc{B}_0,\omega}$ are defined
        in \cref{eq:H0_definition}, and $ V^{\star}\bra{\Omega}\pare{\eta}$,
        $ v_3^{\star}\bra{\Omega}\pare{\eta}$ are the two components of
        ${\bm v}^\star\bra{\Omega}\pare{\eta}$ defined in \cref{eq:bmvtilde};

        \item $A^{\pi^\star,\Omega}_0
        \in \cL\pare{\pare{H^{s_0}}^2;\pare{\Gamma^0}^2}$ and is explicitly defined as
        \begin{equation}\label{eq:Apistar_Omega0}
\begin{aligned}
A^{\pi^\star,\Omega}_0\bra{\Omega}
\pare{\Atw\middle|\eta;x,\xi}
&\defeq
\left(
     V^{\star}\bra{\Omega}\pare{\eta;x}
    +
    \frac{\Atw}{4\pare{1+\av{\nabla\eta}^2}}
    {\bm J}
    \pare{
        \Omega+\pare{\nabla\eta\cdot\Omega}\nabla\eta
    }
\right)^\perp
\\
&\quad
-\ii\,\pare{1+\av{\nabla\eta}^2}
\left(
    \Atw\, v_3^{\star}\bra{\Omega}\pare{\eta;x}
    -
    \frac{\nabla^\perp\eta\cdot\Omega}
    {1+\av{\nabla\eta}^2}
\right)
A_{0}^{\mathpzc{B}_0,\omega;\perp}\pare{\eta;x,\xi},
\end{aligned}
        \end{equation}
        where $ V^{\star}\bra{\Omega}\pare{\eta}$ and $ v_3^{\star}\bra{\Omega}\pare{\eta}$
        are the two components of ${\bm v}^\star\bra{\Omega}\pare{\eta}$ defined in
        \cref{eq:bmvtilde}, $A_{0}^{\mathpzc{B}_0,\omega;\perp} \defeq {\bm J}A_{0}^{\mathpzc{B}_0,\omega}$ with $A_{0}^{\mathpzc{B}_0,\omega}$ defined in
        \cref{eq:H0_definition}, and ${\bm J}$ the
        $\pi/2$-rotation matrix.

    \end{itemize}

    \item $a_{\bra{0,\geq1}}\pare{\sigma,\rho^\pm,\Atw\middle|\eta,\Omega;x,\xi}
    \in\Gamma_{\geq1}^0$ denotes the aggregate generic order-zero $\eta$-symbol: it is the sum of
    the term coming from $\pi_0$ and a term bilinear in $\Omega$. Moreover,
    $A_{\bra{-1,\geq1}}\bra{\Omega} \in \cL\pare{\pare{H^{s_0}}^2;\pare{\Gamma^{-1}}^2}$
    denotes a generic symbol linear in $\Omega$, as per
    \cref{notation:paralinearization};

    \item $r^{\pi_0,\bR}_{\bra{\varrho,\geq1}}
    \pare{\sigma,\rho^\pm\middle|\eta}\in\cR_{\geq1}^\varrho$ is the smoothing
    remainder linear in the final $\eta$ input produced by $\pi_0$;

    \item ${\bm r}^{\rm bil,\bR}_{\bra{\varrho,\geq1}} \in \pare{\cR_{\geq1}^\varrho}^3$ denotes a
    bilinear smoothing remainder, whose $\eta$-component is bilinear in
    $\Omega$ and whose $\Omega$-component is linear in $\Omega$, in the sense of
    \cref{rem:separately-estimated-parameters,rem:linear-parameter-vanishing-order}.
\end{itemize}
\end{prop}

\paragraph{Paralinearization of $\pi_0$}

In this section, we paralinearize the term $\pi_0$, which represents the contributions from surface tension and gravity. For the reader's convenience, we recall from \cref{eq:pi0} its definition:
\begin{equation}\label{eq:pi0star_recall}
    \pi_0\pare{\rho^\pm\middle| \eta} \defeq \frac{\sigma}{\aangle{\rho^\pm}}\nabla \cdot\mathpzc{K}\pare{\nabla\eta}-\mathfrak{g}\Atw\eta,
\end{equation}
where $\mathpzc{K}\pare{\nabla\eta}$ is defined in \cref{curvature2} as $\mathpzc{K}\pare{\nabla\eta} \defeq \frac{\nabla \eta}{\sqrt{1+\av{\nabla \eta}^2}}$.

\begin{lemma}[Paralinearization of $\pi_0$]
\label{lem:paralinearization_pi0star}
Let $\sigma\geq0$, $\rho^\pm\geq0$ with $\aangle{\rho^\pm}>0$, and
$\mathfrak g\in\bR$. There exists $\varrho_0 > 0$ such that for any
$\varrho \geq \varrho_0$, there exist $s_0>0$ and $r_0>0$ such that for any
$\eta\in \BallR{s_0}{r_0}$, the term $\pi_0\pare{\eta}$ defined in
\cref{eq:pi0star_recall} admits the following paralinearization:
\begin{equation}\label{eq:paralinearization_pi0}
    \pi_0\pare{\eta} =
    -\mathfrak{g}\Atw \eta 
    +
    \OpBW{a_{2}^{\pi_0, \eta}\pare{\eta;x,\xi} + a_{\bra{0,\geq1}}\pare{\eta;x, \xi}}\eta
    + r^{\pi_0,\bR}_{\bra{\varrho,\geq1}}\pare{\sigma,\rho^\pm\middle|\eta}\eta,
\end{equation}
where the symbols are defined as follows:
\begin{itemize}
    \item $a_{2}^{\pi_0, \eta}\in\Gamma^2_{\geq0}$ is the real-valued principal symbol, given by
    \begin{equation}\label{eq:symbol_a2_pi0star}
        a_{2}^{\pi_0, \eta}\pare{\eta;x,\xi} \defeq - \frac{\sigma}{\aangle{\rho^{\pm}}} \frac{\av{\xi}^2 + \av{\nabla^\perp\eta\cdot \xi}^2}{\pare{1+\av{\nabla\eta}^2}^{\frac{3}{2}}};
    \end{equation}
    
       \item $a_{\bra{0,\geq1}}\in\Gamma_{\geq1}^0$, cf. \cref{notation:paralinearization};

    \item $r^{\pi_0,\bR}_{\bra{\varrho,\geq1}}
    \pare{\sigma,\rho^\pm\middle|\eta}\in\cR_{\geq1}^\varrho$, cf.
    \cref{notation:paralinearization}.
\end{itemize}
\end{lemma}

\begin{proof}
    Recall from \cref{eq:pi0} and \cref{curvature2} that
\begin{align}\label{eq:pi0star_something}
    \pi_0\pare{\eta}=\frac{\sigma}{\aangle{\rho^{\pm}}}\nabla\cdot\mathpzc{K}\pare{\nabla\eta}-\mathfrak{g}\Atw\eta,
    &&
    \mathpzc{K}\pare{\Omega}=\frac{\Omega}{\sqrt{1+\av{\Omega}^2}}\cdot
\end{align}
By vectorial paralinearization formula
\begin{align}\label{eq:Bony_vectorial_paralinearization_formula}
    F\pare{\Omega} - F\pare{0} =\OpBW{J F\pare{\Omega}} \Omega + R^\bR_{\bra{\varrho,\geq1}}\pare{\Omega}\Omega, 
    &&
    F \in \cC^\infty\pare{\bR^2}, 
\end{align}
and composition rule \cref{prop:composition_paradiff}, we get
\begin{equation}
\begin{aligned}
        \mathpzc{K}\pare{\nabla\eta} = & \  \OpBW{J \mathpzc{K}\pare{\nabla\eta}}\bra{\nabla\eta} + R^\bR_{\bra{\varrho,\geq1}}\pare{\eta}\eta
        \\
        = & \ \OpBW{\ii J \mathpzc{K}\pare{\nabla\eta} \xi + \frac{1}{2\ii} \pbra{J \mathpzc{K}\pare{\nabla\eta}}{\ii \xi} + A_{\bra{-1,\geq1}}\pare{\eta;x, \xi} }\eta + R^\bR_{\bra{\varrho,\geq1}}\pare{\eta}\eta,
\end{aligned}
\end{equation}
where we used the notation \cref{eq:Pbra_matrix-vector}.
Notice in particular that explicit computations show that
\begin{equation}\label{eq:Jk1}
    J\mathpzc{K}\pare{\Omega} = \frac{1}{\pare{1+\av{\Omega}^2}^{3/2}} \ \pare{\Id + \Omega^\perp\otimes \Omega^\perp}.
\end{equation}
Let us now compute
\begin{equation}\label{eq:Jk1.1}
\begin{aligned}
     \nabla\cdot  \mathpzc{K}\pare{\nabla\eta} = & \ \OpBW{\ii \xi^\intercal} \OpBW{\ii J \mathpzc{K}\pare{\nabla\eta} \xi + \frac{1}{2\ii} \pbra{J \mathpzc{K}\pare{\nabla\eta}}{\ii \xi} + A_{\bra{-1,\geq1}}\pare{\eta;x, \xi} }\eta + r^\bR_{\bra{\varrho,\geq1}}\pare{\eta}\eta
     \\
     = & \ 
     \OpBW{
     - \xi^\intercal J \mathpzc{K}\pare{\nabla\eta} \xi
     +
     \frac{\ii}{2}
     \pare{
     \pbra{\xi^\intercal}{J \mathpzc{K}\pare{\nabla\eta} \xi}  
     + \xi^\intercal \pbra{J \mathpzc{K}\pare{\nabla\eta}}{\xi}
         }
 +
     a_{\bra{0,\geq1}}\pare{\eta;x, \xi}
     }\eta
     + r^\bR_{\bra{\varrho,\geq1}}\pare{\eta}\eta . 
\end{aligned}
\end{equation}
We use now \cref{eq:Pbra_matrix-vector} and Einstein's convention for repeated indexes
\begin{align}\label{eq:Jk2}
     \pbra{\xi^\intercal}{J \mathpzc{K}\pare{\nabla\eta} \xi}  
     + \xi^\intercal \pbra{J \mathpzc{K}\pare{\nabla\eta}}{\xi}
     =
      \pbra{\xi_i}{J \mathpzc{K}\pare{\nabla\eta}_{ij} \xi_j}  
     + \xi _i \pbra{J \mathpzc{K}\pare{\nabla\eta}_{ij}}{\xi_j}
\end{align}
and Leibniz's formula $\pbra{fg}{h} = \pbra{f}{h}g + f\pbra{g}{h}$ to obtain that
\begin{equation}\label{eq:Jk3}
    \pbra{\xi_i}{J \mathpzc{K}\pare{\nabla\eta}_{ij} \xi_j}  = - \xi_j \pbra{J \mathpzc{K}\pare{\nabla\eta}_{ij}}{\xi_i} - J \mathpzc{K}\pare{\nabla\eta}_{ij} \underbrace{\pbra{\xi_j}{\xi_i}}_{\equiv 0}. 
\end{equation}
We plug \cref{eq:Jk3} into \cref{eq:Jk2}, and use the symmetry of $J \mathpzc{K}\pare{\nabla\eta}$, which is evident from \cref{eq:Jk1}, to obtain that
\begin{equation}\label{eq:Jk4}
    \pbra{\xi^\intercal}{J \mathpzc{K}\pare{\nabla\eta} \xi}  
     + \xi^\intercal \pbra{J \mathpzc{K}\pare{\nabla\eta}}{\xi}
     \equiv 0. 
\end{equation}
We insert now \cref{eq:Jk4} into \cref{eq:Jk1.1} and use \cref{eq:Jk1} to obtain that
\begin{equation}
 \nabla\cdot  \mathpzc{K}\pare{\nabla\eta} = 
 \OpBW{
    - \frac{\av{\xi}^2 + \av{\nabla^\perp\eta\cdot \xi}^2}{\pare{1+\av{\nabla\eta}^2}^{\frac{3}{2}}} +
     a_{\bra{0,\geq1}}\pare{\eta;x, \xi}
     }\eta
     + r^\bR_{\bra{\varrho,\geq1}}\pare{\eta}\eta.
\end{equation}
After multiplication by $\sigma/\aangle{\rho^\pm}$, we denote the resulting
smoothing operator by
$r^{\pi_0,\bR}_{\bra{\varrho,\geq1}}
\pare{\sigma,\rho^\pm\middle|\eta}$. The claim thus follows using
\cref{eq:pi0star_something}.

\end{proof}

\paragraph{Paralinearization of $\pi^\star_2$}
Here we paralinearize the term $\pi^\star_2$ defined in \cref{eq:pistar2} and whose definition we recall here for the reader's convenience

\begin{subequations}
    \begin{align}
        \label{eq:pistar2_recall}
    \pi^\star_{2}\pare{\Atw\middle| \eta }\bra{\Omega_1, \Omega_2} \defeq & \   - \pi^\star_{2;\RN{1}}\pare{\eta }\bra{\Omega_1, \Omega_2} -  \pi^\star_{2;\RN{2}}\pare{\Atw\middle| \eta }\bra{\Omega_1, \Omega_2},
    \\
    %-----------------------------------------------------------
     \pi^\star_{2;\RN{1}}\pare{\eta }\bra{\Omega_1, \Omega_2} \defeq & \ 
        \frac{ \Omega_1^\perp \cdot  V^{\star}\bra{\Omega_2}\pare{\eta} + \Omega_2^\perp \cdot  V^{\star}\bra{\Omega_1}\pare{\eta}}{2} , \label{eq:pi2_star_I_recall}
        \\
        %--------------------------------------------
        \pi^\star_{2;\RN{2}}\pare{\Atw\middle| \eta }\bra{\Omega_1, \Omega_2} \defeq & \ 
        \frac{\Atw}{4} \Big[
         2  V^{\star}\bra{\Omega_1}\pare{\eta} \cdot  V^{\star}\bra{\Omega_2}\pare{\eta}
         \nonumber
        \\
        & \quad
        -  v_3^{\star}\bra{\Omega_1}\pare{\eta} \pare{ v_3^{\star}\bra{\Omega_2}\pare{\eta} - 2 V^{\star}\bra{\Omega_2}\pare{\eta}\cdot\nabla\eta}
        \nonumber
        \\
        & \quad
        -  v_3^{\star}\bra{\Omega_2}\pare{\eta} \pare{ v_3^{\star}\bra{\Omega_1}\pare{\eta} - 2 V^{\star}\bra{\Omega_1}\pare{\eta}\cdot\nabla\eta}
        \nonumber
        \\
        & \quad
        +
        \frac{1}{2}\frac{\Omega_1^\intercal \pare{\Id_{\bR^2} + \nabla\eta\otimes\nabla\eta}\Omega_2}
        {{1+|\nabla\eta|^2}}
        \Big].
        \label{eq:pi2_star_II_recall}
    \end{align}
\end{subequations}

The result we prove is the following

\begin{lemma}[Paralinearization of $\pi^\star_2$]\label{lem:paralinearization_pistar2}
There exists $\varrho_0 > 0$ such that for any $\varrho \geq \varrho_0$, there exist
$s_0>0$ and $r_0>0$ such that for any
$\eta\in \BallR{s_0}{r_0}$ and
$\Omega\in \pare{\BallR{s_0}{r_0}}^2$, the term
$\pi^\star_2\pare{\Atw\middle|\eta}\bra{\Omega,\Omega}$ defined in
\cref{eq:pistar2_recall} admits the following paralinearization:
\begin{equation}\label{eq:paralinearization_pistar2}
\begin{aligned}
    \pi^\star_2\pare{\Atw\middle|\eta}\bra{\Omega,\Omega}
    &=
    \OpBW{
        {\bm a}^{\pi^\star_2}_1
        \pare{\Atw\middle|\eta,\Omega;x,\xi}
        +
        \vect{
            a_{\bra{0,\geq1}}\bra{\Omega,\Omega}\pare{\Atw\middle|\eta;x,\xi}
        }{
            A_{\bra{-1,\geq1}}\bra{\Omega}\pare{\Atw\middle|\eta;x,\xi}
        }
    }\cdot\vect{\eta}{\Omega}
  +
   {\bm r}^{\bil,\bR}_{\bra{\varrho,\geq1}} \pare{\Atw\middle| \eta, \Omega}\cdot\vect{\eta}{\Omega},
\end{aligned}
\end{equation}
where
\begin{itemize}
    \item $ {\bm a}^{\pi^\star_2}_1\in \pare{\Gamma_{\geq1}^1}^3$ and is explicitly defined as
    \begin{equation}
\label{eq:bmapistar2}
       {\bm a}^{\pi^\star_2}_1
    \pare{\Atw\middle|\eta,\Omega;x,\xi}
    \defeq
    \vect{
        a^{\pi^\star_2,\eta}_1\bra{\Omega,\Omega}
        \pare{\Atw\middle|\eta;x,\xi}
    }{
        A^{\pi^\star_2,\Omega}_0\bra{\Omega}
        \pare{\Atw\middle|\eta;x,\xi}
    } ,
\end{equation}
where
\begin{itemize}

    \item $a^{\pi^\star_2,\eta}_1\in \cL^2\pare{\pare{H^{s_0}}^2; \Gamma^1 }$ and is explicitly defined as
    \begin{equation}\label{eq:pistar2_eta_symbol}
\begin{aligned}
a^{\pi^\star_2,\eta}_1\bra{\Omega,\Omega}
\pare{\Atw\middle|\eta;x,\xi}
&\defeq
\lambda_1\pare{\eta;x,\xi}
\left[
    4\pare{A_{0}^{\mathpzc{B}_0,\omega;\perp}\pare{\eta;x,\xi}\cdot\Omega}^2
    -
    \left(
        \frac{\nabla^\perp\eta\cdot\Omega}
        {1+\av{\nabla\eta}^2}
    \right)^2
\right]
\\
&\quad
+\Atw\,
\pare{A_{0}^{\mathpzc{B}_0,\omega;\perp}\pare{\eta;x,\xi}\cdot\Omega}
\,
 V^{\star}\bra{\Omega}\pare{\eta;x}\cdot\xi
\\
&\quad
+\Atw\,\lambda_1\pare{\eta;x,\xi}
\, v_3^{\star}\bra{\Omega}\pare{\eta;x}
\frac{\nabla^\perp\eta\cdot\Omega}
{1+\av{\nabla\eta}^2}
\\
&\quad
-\ii\Atw
\left[
     v_3^{\star}\bra{\Omega}\pare{\eta;x}
     V^{\star}\bra{\Omega}\pare{\eta;x}
\right.
\\
&\qquad\qquad\left.
    +
    \frac{1}{4\pare{1+\av{\nabla\eta}^2}}
    \left(
        \pare{\Omega\cdot\nabla\eta}\Omega
        -
        \frac{\av{\Omega}^2+\av{\Omega\cdot\nabla\eta}^2}
        {1+\av{\nabla\eta}^2}
        \nabla\eta
    \right)
\right]\cdot\xi ,
\end{aligned}
\end{equation}
  where $\lambda_1 \in \Gamma^1_{\geq0}$ is defined in \cref{eq:lambda1},
        $A_{0}^{\mathpzc{B}_0,\omega} \in \pare{\Gamma^0_{\geq0}}^2$ is defined in
        \cref{eq:H0_definition}, and $ V^{\star}\bra{\Omega}\pare{\eta}$,
        $ v_3^{\star}\bra{\Omega}\pare{\eta}$ are the two components of
        ${\bm v}^\star\bra{\Omega}\pare{\eta}$ defined in \cref{eq:bmvtilde};

    \item $A^{\pi^\star_2,\Omega}_0\in \cL\pare{\pare{H^{s_0}}^2 ; \pare{\Gamma^0}^2 }$ and is explicitly defined as
    \begin{equation}\label{eq:pistar2_Omega_symbol}
\begin{aligned}
A^{\pi^\star_2,\Omega}_0\bra{\Omega}
\pare{\Atw\middle|\eta;x,\xi}
&\defeq
\left(
     V^{\star}\bra{\Omega}\pare{\eta;x}
    +
    \frac{\Atw}{4\pare{1+\av{\nabla\eta}^2}}
    {\bm J}
    \pare{
        \Omega+\pare{\nabla\eta\cdot\Omega}\nabla\eta
    }
\right)^\perp
\\
&\quad
-\ii\,\pare{1+\av{\nabla\eta}^2}
\left(
    \Atw\, v_3^{\star}\bra{\Omega}\pare{\eta;x}
    -
    \frac{\nabla^\perp\eta\cdot\Omega}
    {1+\av{\nabla\eta}^2}
\right)
A_{0}^{\mathpzc{B}_0,\omega;\perp}\pare{\eta;x,\xi} \ ;
\end{aligned}
\end{equation}
        where $ V^{\star}\bra{\Omega}\pare{\eta}$ and $ v_3^{\star}\bra{\Omega}\pare{\eta}$
        are the two components of ${\bm v}^\star\bra{\Omega}\pare{\eta}$ defined in
        \cref{eq:bmvtilde}, $A_{0}^{\mathpzc{B}_0,\omega;\perp} \defeq {\bm J}A_{0}^{\mathpzc{B}_0,\omega}$ with $A_{0}^{\mathpzc{B}_0,\omega}$ defined in
        \cref{eq:H0_definition}, and ${\bm J}$ the
        $\pi/2$-rotation matrix.

\end{itemize}

    \item $a_{\bra{0,\geq1}}\bra{\Omega,\Omega} \in \cL^2\pare{\pare{H^{s_0}}^2;\Gamma^0}$
    denotes a generic symbol bilinear in $\Omega$, and
    $A_{\bra{-1,\geq1}}\bra{\Omega} \in \cL\pare{\pare{H^{s_0}}^2;\pare{\Gamma^{-1}}^2}$
    denotes a generic symbol linear in $\Omega$, as per
    \cref{notation:paralinearization};

    \item ${\bm r}^{\rm bil,\bR}_{\bra{\varrho,\geq1}} \in \pare{\cR_{\geq1}^\varrho}^3$ denotes a
    bilinear smoothing remainder, whose $\eta$-component is bilinear in
    $\Omega$ and whose $\Omega$-component is linear in $\Omega$, in the sense of
    \cref{rem:separately-estimated-parameters,rem:linear-parameter-vanishing-order}.

\end{itemize}

\end{lemma}

\begin{cor}[Bilinear paralinearization of $\pi^\star_2$]\label{cor:paralinearization_pistar2_bilinear}
There exists $\varrho_0 > 0$ such that for any $\varrho \geq \varrho_0$, there exist
$s_0>0$ and $r_0>0$ such that for any
$\eta\in \BallR{s_0}{r_0}$ and
$\Omega_1,\Omega_2\in \pare{\BallR{s_0}{r_0}}^2$, the bilinear form
$\pi^\star_2\pare{\Atw\middle|\eta}\bra{\Omega_1,\Omega_2}$ defined in
\cref{eq:pistar2_recall} admits the following paralinearization:
\begin{equation}\label{eq:paralinearization_pistar2_bilinear}
\begin{aligned}
    \pi^\star_2\pare{\Atw\middle|\eta}\bra{\Omega_1,\Omega_2}
    &=
    \OpBW{
        a^{\pi^\star_2,\eta}_1
        \bra{\Omega_1,\Omega_2}
        \pare{\Atw\middle|\eta;x,\xi}
        +
        a_{\bra{0,\geq1}}\bra{\Omega_1,\Omega_2}
        \pare{\Atw\middle|\eta;x,\xi}
    }\eta
    \\
    &\quad
    +
    \frac{1}{2}
    \OpBW{
        A^{\pi^\star_2,\Omega}_0
        \bra{\Omega_2}
        \pare{\Atw\middle|\eta;x,\xi}
        +
        A_{\bra{-1,\geq1}}\bra{\Omega_2}
        \pare{\Atw\middle|\eta;x,\xi}
    }\cdot\Omega_1
    \\
    &\quad
    +
    \frac{1}{2}
    \OpBW{
        A^{\pi^\star_2,\Omega}_0
        \bra{\Omega_1}
        \pare{\Atw\middle|\eta;x,\xi}
        +
        A_{\bra{-1,\geq1}}\bra{\Omega_1}
        \pare{\Atw\middle|\eta;x,\xi}
    }\cdot\Omega_2
    \\
    &\quad
    +
    r^{\bil,\bR}_{\bra{\varrho,\geq1}}
    \bra{\Omega_1,\Omega_2}
    \pare{\Atw\middle|\eta}\eta
    \\
    &\quad
    +
    \frac{1}{2}
    R^\bR_{\bra{\varrho,\geq1}}\bra{\Omega_2}
    \pare{\Atw\middle|\eta}\cdot\Omega_1
    +
    \frac{1}{2}
    R^\bR_{\bra{\varrho,\geq1}}\bra{\Omega_1}
    \pare{\Atw\middle|\eta}\cdot\Omega_2,
\end{aligned}
\end{equation}
where
\begin{itemize}
    \item $a^{\pi^\star_2,\eta}_1
    \in \cL^2\pare{\pare{H^{s_0}}^2;\Gamma^1}$ and is explicitly defined as
    \begin{equation}\label{eq:pistar2_eta_symbol_bilinear}
\begin{aligned}
&a^{\pi^\star_2,\eta}_1
\bra{\Omega_1,\Omega_2}
\pare{\Atw\middle|\eta;x,\xi}\\
\defeq &
\lambda_1\pare{\eta;x,\xi}
\left[
    4
    \pare{A_{0}^{\mathpzc{B}_0,\omega;\perp}\pare{\eta;x,\xi}\cdot\Omega_1}
    \pare{A_{0}^{\mathpzc{B}_0,\omega;\perp}\pare{\eta;x,\xi}\cdot\Omega_2}
-
    \frac{
        \pare{\nabla^\perp\eta\cdot\Omega_1}
        \pare{\nabla^\perp\eta\cdot\Omega_2}
    }{
        \pare{1+\av{\nabla\eta}^2}^2
    }
\right]
\\
+&
\frac{\Atw}{2}
\pare{A_{0}^{\mathpzc{B}_0,\omega;\perp}\pare{\eta;x,\xi}\cdot\Omega_1}
\,
 V^{\star}\bra{\Omega_2}\pare{\eta;x}\cdot\xi
+
\frac{\Atw}{2}
\pare{A_{0}^{\mathpzc{B}_0,\omega;\perp}\pare{\eta;x,\xi}\cdot\Omega_2}
\,
 V^{\star}\bra{\Omega_1}\pare{\eta;x}\cdot\xi
\\
+&
\frac{
    \Atw\,\lambda_1\pare{\eta;x,\xi}
}{
    2\pare{1+\av{\nabla\eta}^2}
}
\left[
     v_3^{\star}\bra{\Omega_1}\pare{\eta;x}
    \pare{\nabla^\perp\eta\cdot\Omega_2}
%\right.
%\\
%&\qquad\qquad\left.
    +
     v_3^{\star}\bra{\Omega_2}\pare{\eta;x}
    \pare{\nabla^\perp\eta\cdot\Omega_1}
\right]
\\
-&\ii\Atw
\left[
    \frac{1}{2}
     v_3^{\star}\bra{\Omega_1}\pare{\eta;x}
     V^{\star}\bra{\Omega_2}\pare{\eta;x}
+
    \frac{1}{2}
     v_3^{\star}\bra{\Omega_2}\pare{\eta;x}
     V^{\star}\bra{\Omega_1}\pare{\eta;x}
\right.
\\
   &\qquad \left. +
    \frac{1}{4\pare{1+\av{\nabla\eta}^2}}
    \left(
        \frac{1}{2}
        \pare{\Omega_1\cdot\nabla\eta}\Omega_2
        +
        \frac{1}{2}
        \pare{\Omega_2\cdot\nabla\eta}\Omega_1
    \right.
\left.
        -
        \frac{
            \Omega_1\cdot\Omega_2
            +
            \pare{\Omega_1\cdot\nabla\eta}
            \pare{\Omega_2\cdot\nabla\eta}
        }{
            1+\av{\nabla\eta}^2
        }
        \nabla\eta
    \right)
\right]\cdot\xi ;
\end{aligned}
    \end{equation}
    In particular, setting $\Omega_1=\Omega_2=\Omega$ gives
    \cref{eq:pistar2_eta_symbol}.

    \item $A^{\pi^\star_2,\Omega}_0
    \in \cL\pare{\pare{H^{s_0}}^2;\pare{\Gamma^0}^2}$ and is explicitly defined as
    \begin{equation}\label{eq:pistar2_Omega_symbol_bilinear}
\begin{aligned}
A^{\pi^\star_2,\Omega}_0\bra{\Omega_j}
\pare{\Atw\middle|\eta;x,\xi}
&\defeq
\left(
     V^{\star}\bra{\Omega_j}\pare{\eta;x}
    +
    \frac{\Atw}{4\pare{1+\av{\nabla\eta}^2}}
    {\bm J}
    \pare{
        \Omega_j+\pare{\nabla\eta\cdot\Omega_j}\nabla\eta
    }
\right)^\perp
\\
&\quad
-\ii\,\pare{1+\av{\nabla\eta}^2}
\left(
    \Atw\, v_3^{\star}\bra{\Omega_j}\pare{\eta;x}
    -
    \frac{\nabla^\perp\eta\cdot\Omega_j}
    {1+\av{\nabla\eta}^2}
\right)
A_{0}^{\mathpzc{B}_0,\omega;\perp}\pare{\eta;x,\xi},
\end{aligned}
    \end{equation}
    for $j=1,2$. In particular, setting $\Omega_1=\Omega_2=\Omega$ in
    \cref{eq:paralinearization_pistar2_bilinear} gives the $\Omega$-component
    in \cref{eq:pistar2_Omega_symbol};

    \item $a_{\bra{m,\geq1}}$ and $A_{\bra{m,\geq1}}$ are generic symbols as per
    \cref{notation:paralinearization};

    \item $r^{\bil,\bR}_{\bra{\varrho,\geq1}}$ and $R^\bR_{\bra{\varrho,\geq1}}$ are smoothing remainders as per
    \cref{notation:paralinearization}.
\end{itemize}
\end{cor}

\begin{proof}
The claim follows from the polarization identity
\begin{equation}\label{eq:polarization_pistar2}
\begin{aligned}
\pi^\star_2\pare{\Atw\middle|\eta}\bra{\Omega_1,\Omega_2}
&=
\frac{1}{2}
\Big(
    \pi^\star_2\pare{\Atw\middle|\eta}
    \bra{\Omega_1+\Omega_2,\Omega_1+\Omega_2}
    -
    \pi^\star_2\pare{\Atw\middle|\eta}
    \bra{\Omega_1,\Omega_1}
    -
    \pi^\star_2\pare{\Atw\middle|\eta}
    \bra{\Omega_2,\Omega_2}
\Big).
\end{aligned}
\end{equation}
Applying \cref{lem:paralinearization_pistar2} to each term in
\cref{eq:polarization_pistar2}, the $\eta$-symbol is the polarization of
\cref{eq:pistar2_eta_symbol}, which gives
\cref{eq:pistar2_eta_symbol_bilinear}. The $\Omega$-symbol
$A^{\pi^\star_2,\Omega}_0$ is linear in its argument, hence
\[
\frac{1}{2}
\left[
    A^{\pi^\star_2,\Omega}_0\bra{\Omega_1+\Omega_2}
    \cdot\pare{\Omega_1+\Omega_2}
    -
    A^{\pi^\star_2,\Omega}_0\bra{\Omega_1}\cdot\Omega_1
    -
    A^{\pi^\star_2,\Omega}_0\bra{\Omega_2}\cdot\Omega_2
\right]
\]
equals
\[
    \frac{1}{2}
    A^{\pi^\star_2,\Omega}_0\bra{\Omega_2}\cdot\Omega_1
    +
    \frac{1}{2}
    A^{\pi^\star_2,\Omega}_0\bra{\Omega_1}\cdot\Omega_2 .
\]
The generic symbols and the smoothing remainders polarize in the same way,
yielding \cref{eq:paralinearization_pistar2_bilinear}.
\end{proof}

\begin{proof}[Proof of \cref{lem:paralinearization_pistar2}]
All along the proof we write
\begin{align}
\alpha \simeq \beta
&&
\iff
&&
    \alpha \equiv \beta  \mod \OpBW{  \vect{a_{\bra{0,\geq1}}\bra{\Omega,\Omega}\pare{\eta;x,\xi}}{A_{\bra{-1,\geq1}}\bra{\Omega}\pare{\eta;x,\xi}}}\cdot \vect{\eta}{\Omega}
        +
        \vect{r^\bR_{\bra{\varrho,\geq1}}\bra{\Omega,\Omega}\pare{\eta}}{R^\bR_{\bra{\varrho,\geq1}}\bra{\Omega}\pare{\eta}}\cdot\vect{\eta}{\Omega}. 
\end{align}

    \begin{proofpart}[Paralinearization of $\pi^\star_{2;\RN{1}}$]
This is a straightforward application of \cref{lem:paraprod_Bony,lem:composition_rules,prop:paralinearization_Vtilde,prop:composition_paradiff}, which gives that
    \begin{equation}\label{eq:paralin_pistar2I}
        \pi^\star_{2;\RN{1}}\pare{\eta}\bra{\Omega, \Omega} \simeq \OpBW{\vect{a^{\pi^\star_{2;\RN{1}},\eta}_1\bra{\Omega,\Omega}\pare{\eta;x,\xi}}{A^{\pi^\star_{2;\RN{1}},\Omega}_0\bra{\Omega}\pare{\eta;x,\xi}}} \cdot\vect{\eta}{\Omega}
      ,
    \end{equation}
    where
    \begin{align}
        a^{\pi^\star_{2;\RN{1}},\eta}_1\bra{\Omega,\Omega}\pare{\eta;x,\xi} \defeq & \ - \Omega^\intercal {\bm J}{\bm A}^{ V^{\star}, \eta}_1\pare{\eta;x,\xi} \Omega,
        \\
        A^{\pi^\star_{2;\RN{1}},\Omega}_0\bra{\Omega}\pare{\eta;x,\xi} \defeq & \ - V^{\star}\bra{\Omega}\pare{\eta}^\perp + \ii  \pare{\Omega^\perp}^\intercal   {\bm A}^{ V^{\star},\Omega}_0\pare{\eta;x,\xi}.
    \end{align}
        
    \end{proofpart}

\begin{proofpart}[Paralinearization of $\pi^\star_{2;\RN{2}}$]
    Let us write
    \begin{equation}
        \pi^\star_{2;\RN{2}}\pare{\eta }\bra{\Omega , \Omega} = \frac{\Atw}{2} \pare{\av{ V^{\star}\bra{\Omega}\pare{\eta}}^2 -  v_3^{\star}\bra{\Omega}\pare{\eta} \pare{ v_3^{\star}\bra{\Omega}\pare{\eta} - 2  V^{\star}\bra{\Omega}\pare{\eta}\cdot\nabla\eta} + \frac{1}{4}\frac{\av{\Omega}^2 + \av{\nabla\eta\cdot\Omega}^2}{1+\av{\nabla\eta}^2} }.
    \end{equation}
    Using \cref{lem:paraprod_Bony,prop:paralinearization_Vtilde} we obtain that
    \begin{equation}\label{eq:pistar2II_1}
        \av{ V^{\star}\bra{\Omega}\pare{\eta}}^2 \simeq \OpBW{ 2  V^{\star}\bra{\Omega}\pare{\eta}^\intercal \mathfrak{A}^{ V^{\star}}_1\pare{\eta, \Omega;x,\xi} } \vect{\eta}{\Omega},
    \end{equation}
    where $\mathfrak{A}^{ V^{\star}_1}$ is defined in \cref{eq:frakAVtilde}. Again \cref{lem:paraprod_Bony,prop:paralinearization_Vtilde} give us that
    \begin{footnotesize}
        \begin{equation}
        \begin{aligned}\label{eq:pistar2II_2}
        &-  v_3^{\star}\bra{\Omega}\pare{\eta} \pare{ v_3^{\star}\bra{\Omega}\pare{\eta} - 2  V^{\star}\bra{\Omega}\pare{\eta}\cdot\nabla\eta}
        \\
        &\simeq \OpBW{
        -  v_3^{\star}\bra{\Omega}\pare{\eta}  \pare{{\bm a}_1^{ v_3^{\star}}\pare{\eta,\Omega} - 2 \pare{ \pare{\mathfrak{A}^{ V^{\star}}_1\pare{\eta, \Omega} }^\intercal \nabla\eta+ \vect{\ii \xi \cdot V^{\star}\bra{\Omega}\pare{\eta}}{0_{\bR^2}}}} 
        -
        \pare{ v_3^{\star}\bra{\Omega}\pare{\eta} - 2  V^{\star}\bra{\Omega}\pare{\eta}\cdot\nabla\eta} {\bm a}_1^{ v_3^{\star}}\pare{\eta,\Omega} 
        }
        \cdot
        \vect{\eta}{\Omega}.
        \end{aligned}
    \end{equation}
    \end{footnotesize}
Let $h: \bR^4 \to \bR_+$ defined as $h\pare{X,Y}\defeq \frac{1}{4}\frac{{X}^2 + \av{X\cdot Y}^2}{1+\av{Y}^2}$, with such notation we have that
\begin{equation}
    \frac{1}{4}\frac{\av{\Omega}^2 + \av{\nabla\eta\cdot\Omega}^2}{1+\av{\nabla\eta}^2} =  h\pare{\Omega , \ \nabla\eta}.
\end{equation}
We apply \cref{eq:Bony_vectorial_paralinearization_formula,lem:composition_rules} we obtain that
\begin{equation}\label{eq:pistar23_1}
      \frac{1}{4}\frac{\av{\Omega}^2 + \av{\nabla\eta\cdot\Omega}^2}{1+\av{\nabla\eta}^2}  \simeq \OpBW{J _{X,Y} h\pare{\Omega , \ \nabla\eta} }  \vect{\Omega}{\nabla\eta} ,
\end{equation}
where
\begin{equation}
    J _{X,Y} h \pare{X,Y}
    = 
    \frac{1}{2}
    \vect{\frac{X + (X\cdot Y)Y}{1+\av{Y}^2}}{\frac{(1+\av{Y}^2)(X\cdot Y)X - (\av{X}^2+(X\cdot Y)^2)Y}{(1+\av{Y}^2)^2}}^\intercal
\end{equation}
so that, applying \cref{prop:composition_paradiff} we obtain that
\begin{equation}
     \frac{1}{4}\frac{\av{\Omega}^2 + \av{\nabla\eta\cdot\Omega}^2}{1+\av{\nabla\eta}^2}  \simeq \OpBW{{\bm a}^h_1\pare{\eta,\Omega;x,\xi}}\cdot\vect{\eta}{\Omega},
\end{equation}
where
\begin{equation}\label{eq:pistar2II_3}
    {\bm a}^h_1\pare{\eta,\Omega} \defeq
    \vect{\frac{\ii}{2} \frac{1}{1+\av{\nabla\eta}^2}\pare{\Omega\cdot\nabla\eta \ \Omega - \frac{\av{\Omega}^2 + \av{\Omega\cdot\nabla\eta}^2}{1+\av{\nabla\eta}^2} \ \nabla\eta}\cdot\xi}{\frac{1}{2}\frac{\Omega + \Omega\cdot\nabla\eta \ \nabla\eta}{1+\av{\nabla\eta}^2}}.
\end{equation}

\end{proofpart}

\begin{proofpart}[Synthesis]
We keep all principal quantities expanded in the notation of the previous
steps. Let $ {\bm a}^{\pi^\star_2}_1$ be as in \cref{eq:bmapistar2} so that
\[
    \pi^\star_2\pare{\Atw\middle|\eta}\bra{\Omega,\Omega}
    \simeq
    \OpBW{
        {\bm a}^{\pi^\star_2}_1\bra{\Omega}
        \pare{\Atw\middle|\eta;x,\xi}
    }\cdot\vect{\eta}{\Omega}.
\]

We first compute the $\Omega$-component. Since
$\pi^\star_2=-\pi^\star_{2;\RN{1}}-\pi^\star_{2;\RN{2}}$, from
\cref{eq:paralin_pistar2I,eq:pistar2II_1,eq:pistar2II_2,eq:pistar2II_3}
we obtain
\begin{align}
A^{\pi^\star_2,\Omega}_0\bra{\Omega}
\pare{\Atw\middle|\eta;x,\xi}
&=
 V^{\star}\bra{\Omega}\pare{\eta;x}^{\perp}
-\ii\pare{\Omega^\perp}^{\intercal}
{\bm A}^{ V^{\star},\Omega}_0\pare{\eta;x,\xi}
\notag
\\
&\quad
-\frac{\Atw}{2}
\Bigg[
    2\ii\,
     V^{\star}\bra{\Omega}\pare{\eta;x}^{\intercal}
    {\bm A}^{ V^{\star},\Omega}_0\pare{\eta;x,\xi}
\notag
\\
&\qquad
- v_3^{\star}\bra{\Omega}\pare{\eta;x}
\left(
    \ii A^{ v_3^{\star},\Omega}_0\pare{\eta;x,\xi}
    -2\ii
    \pare{{\bm A}^{ V^{\star},\Omega}_0\pare{\eta;x,\xi}}^{\intercal}
    \nabla\eta
\right)
\notag
\\
&\qquad
-\left(
     v_3^{\star}\bra{\Omega}\pare{\eta;x}
    -2 V^{\star}\bra{\Omega}\pare{\eta;x}\cdot\nabla\eta
\right)
\ii A^{ v_3^{\star},\Omega}_0\pare{\eta;x,\xi}
\notag
\\
&\qquad
+
\frac{1}{2}
\frac{
    \Omega+\pare{\Omega\cdot\nabla\eta}\nabla\eta
}{
    1+\av{\nabla\eta}^2
}
\Bigg].
\label{eq:pistar2_Omega_synthesis_raw}
\end{align}
Using
\[
    {\bm A}^{ V^{\star},\Omega}_0\pare{\eta;x,\xi}
    =
    \nabla\eta\otimes A_{0}^{\mathpzc{B}_0,\omega;\perp}\pare{\eta;x,\xi},
    \qquad
    A^{ v_3^{\star},\Omega}_0\pare{\eta;x,\xi}
    =
    -A_{0}^{\mathpzc{B}_0,\omega;\perp}\pare{\eta;x,\xi},
\]
we have
\begin{equation}\label{eq:pistar2_Omega_aux1}
    \pare{\Omega^\perp}^{\intercal}
    {\bm A}^{ V^{\star},\Omega}_0\pare{\eta;x,\xi}
    =
    -\pare{\nabla^\perp\eta\cdot\Omega}
    A_{0}^{\mathpzc{B}_0,\omega;\perp}\pare{\eta;x,\xi},
\end{equation}
and
\begin{equation}
  \label{eq:pistar2_Omega_aux2}
    \begin{aligned}
    &2\ii\,
     V^{\star}\bra{\Omega}\pare{\eta;x}^{\intercal}
    {\bm A}^{ V^{\star},\Omega}_0\pare{\eta;x,\xi}\\
-& v_3^{\star}\bra{\Omega}\pare{\eta;x}
\left(
    \ii A^{ v_3^{\star},\Omega}_0\pare{\eta;x,\xi}
    -2\ii
    \pare{{\bm A}^{ V^{\star},\Omega}_0\pare{\eta;x,\xi}}^{\intercal}
    \nabla\eta
\right)
\\
-&\left(
     v_3^{\star}\bra{\Omega}\pare{\eta;x}
    -2 V^{\star}\bra{\Omega}\pare{\eta;x}\cdot\nabla\eta
\right)
\ii A^{ v_3^{\star},\Omega}_0\pare{\eta;x,\xi}
\\
=
&2\ii
\pare{1+\av{\nabla\eta}^2}
 v_3^{\star}\bra{\Omega}\pare{\eta;x}
A_{0}^{\mathpzc{B}_0,\omega;\perp}\pare{\eta;x,\xi}.
\end{aligned}
\end{equation}
Inserting \cref{eq:pistar2_Omega_aux1} and
\cref{eq:pistar2_Omega_aux2} into
\cref{eq:pistar2_Omega_synthesis_raw} gives
\begin{equation}\label{eq:pistar2_Omega_synthesis_final}
\begin{aligned}
A^{\pi^\star_2,\Omega}_0\bra{\Omega}
\pare{\Atw\middle|\eta;x,\xi}
&=
 V^{\star}\bra{\Omega}\pare{\eta;x}^{\perp}
-
\frac{\Atw}{4\pare{1+\av{\nabla\eta}^2}}
\pare{
    \Omega+\pare{\Omega\cdot\nabla\eta}\nabla\eta
}\\
&\quad
-\ii
\pare{1+\av{\nabla\eta}^2}
\left(
    \Atw\, v_3^{\star}\bra{\Omega}\pare{\eta;x}
    -
    \frac{\nabla^\perp\eta\cdot\Omega}
    {1+\av{\nabla\eta}^2}
\right)
A_{0}^{\mathpzc{B}_0,\omega;\perp}\pare{\eta;x,\xi}.
\end{aligned}
\end{equation}
Since
\begingroup
\[
    \left(
        {\bm J}
        \pare{
            \Omega+\pare{\nabla\eta\cdot\Omega}\nabla\eta
        }
    \right)^\perp
    =
    -\pare{
        \Omega+\pare{\nabla\eta\cdot\Omega}\nabla\eta
    },
\]
\endgroup
the identity \cref{eq:pistar2_Omega_synthesis_final} is exactly
\cref{eq:pistar2_Omega_symbol}.

We now compute the $\eta$-component. From
\cref{eq:paralin_pistar2I} and \cref{cor:simplifyimportant}, the
contribution of $-\pi^\star_{2;\RN{1}}$ is
\begin{equation}\label{eq:pistar2_eta_I_synthesis}
\begin{aligned}
    -a^{\pi^\star_{2;\RN{1}},\eta}_1
    \bra{\Omega,\Omega}\pare{\eta;x,\xi}
    &=
    -\pare{{\bm A}^{ V^{\star},\eta}_1\pare{\eta;x,\xi}\Omega}
    \cdot\Omega^\perp
    \\
    &=
    \lambda_1\pare{\eta;x,\xi}
    \left[
        \pare{2A_{0}^{\mathpzc{B}_0,\omega;\perp}\pare{\eta;x,\xi}\cdot\Omega}^2
        -
        \left(
            \frac{\nabla^\perp\eta\cdot\Omega}
            {1+\av{\nabla\eta}^2}
        \right)^2
    \right].
\end{aligned}
\end{equation}

For the contribution of $-\pi^\star_{2;\RN{2}}$, we first record the two
pointwise identities
\begin{equation}\label{eq:pistar2_aux_identity_1}
    {\bm A}^{ V^{\star},\eta}_1\pare{\eta;x,\xi}\Omega
    +
    \nabla\eta\,
    a^{ v_3^{\star},\eta}_1\bra{\Omega}
    \pare{\eta;x,\xi}
    =
    -\xi
    \pare{A_{0}^{\mathpzc{B}_0,\omega;\perp}\pare{\eta;x,\xi}\cdot\Omega},
\end{equation}
and
\begin{equation}\label{eq:pistar2_aux_identity_2}
\begin{aligned}
    \pare{1+\av{\nabla\eta}^2}
    a^{ v_3^{\star},\eta}_1\bra{\Omega}
    \pare{\eta;x,\xi}
    &=
    -\pare{\nabla\eta\cdot\xi}
    \pare{A_{0}^{\mathpzc{B}_0,\omega;\perp}\pare{\eta;x,\xi}\cdot\Omega}
    +
    \lambda_1\pare{\eta;x,\xi}
    \frac{\nabla^\perp\eta\cdot\Omega}
    {1+\av{\nabla\eta}^2}\cdot
\end{aligned}
\end{equation}
Indeed, \cref{eq:pistar2_aux_identity_1} follows from
\cref{eq:AVtildeeta1_lem} and the definition of
$a^{ v_3^{\star},\eta}_1$, while \cref{eq:pistar2_aux_identity_2}
follows from \cref{eq:relationimportant4}.

Using \cref{eq:pistar2_aux_identity_1} and
\cref{eq:pistar2_aux_identity_2}, the $\eta$-symbol coming from
\cref{eq:pistar2II_1} and \cref{eq:pistar2II_2} satisfies
\begin{equation}\label{eq:pistar2_eta_II12_synthesis}
\begin{aligned}
&-\frac{\Atw}{2}
\Big[
    2 V^{\star}\bra{\Omega}\pare{\eta;x}\cdot
    {\bm A}^{ V^{\star},\eta}_1\pare{\eta;x,\xi}\Omega
    -2 v_3^{\star}\bra{\Omega}\pare{\eta;x}
    a^{ v_3^{\star},\eta}_1\bra{\Omega}\pare{\eta;x,\xi}
    +2 v_3^{\star}\bra{\Omega}\pare{\eta;x}
    \nabla\eta\cdot
    {\bm A}^{ V^{\star},\eta}_1\pare{\eta;x,\xi}\Omega
\\
    &+2\pare{
         V^{\star}\bra{\Omega}\pare{\eta;x}\cdot\nabla\eta
    }
    a^{ v_3^{\star},\eta}_1\bra{\Omega}\pare{\eta;x,\xi}
    +2\ii\,
     v_3^{\star}\bra{\Omega}\pare{\eta;x}
     V^{\star}\bra{\Omega}\pare{\eta;x}\cdot\xi
\Big]
\\
&=
\Atw
\pare{A_{0}^{\mathpzc{B}_0,\omega;\perp}\pare{\eta;x,\xi}\cdot\Omega}
 V^{\star}\bra{\Omega}\pare{\eta;x}\cdot\xi
+
\Atw\,\lambda_1\pare{\eta;x,\xi}
 v_3^{\star}\bra{\Omega}\pare{\eta;x}
\frac{\nabla^\perp\eta\cdot\Omega}
{1+\av{\nabla\eta}^2}
\\
&\quad-\ii\Atw\,
 v_3^{\star}\bra{\Omega}\pare{\eta;x}
 V^{\star}\bra{\Omega}\pare{\eta;x}\cdot\xi .
\end{aligned}
\end{equation}
Finally, the contribution of \cref{eq:pistar2II_3} to
$-\pi^\star_{2;\RN{2}}$ is
\begin{equation}\label{eq:pistar2_eta_II3_synthesis}
\begin{aligned}
&-\frac{\Atw}{2}
\left[
    \frac{\ii}{2}
    \frac{1}{1+\av{\nabla\eta}^2}
    \left(
        \pare{\Omega\cdot\nabla\eta}\Omega
        -
        \frac{\av{\Omega}^2+\av{\Omega\cdot\nabla\eta}^2}
        {1+\av{\nabla\eta}^2}
        \nabla\eta
    \right)\cdot\xi
\right]
=
-\frac{\ii\Atw}{4\pare{1+\av{\nabla\eta}^2}}
\left(
    \pare{\Omega\cdot\nabla\eta}\Omega
    -
    \frac{\av{\Omega}^2+\av{\Omega\cdot\nabla\eta}^2}
    {1+\av{\nabla\eta}^2}
    \nabla\eta
\right)\cdot\xi .
\end{aligned}
\end{equation}
Combining
\cref{eq:pistar2_eta_I_synthesis},
\cref{eq:pistar2_eta_II12_synthesis}, and
\cref{eq:pistar2_eta_II3_synthesis} gives precisely
\cref{eq:pistar2_eta_symbol}. Therefore the paralinearization
\cref{eq:paralinearization_pistar2} follows.
\end{proofpart}

\end{proof}

\begin{proof}[Proof of \cref{prop:paralinearization_pistar}]
By \cref{eq:pistar_recall},
\[
    \pi^\star\pare{\rho^\pm\middle|\eta,\Omega}
    =
    \pi_0\pare{\rho^\pm\middle|\eta}
    +
    \pi^\star_2\pare{\Atw\middle|\eta}\bra{\Omega,\Omega}.
\]
The term \(\pi_0\) is paralinearized in
\cref{lem:paralinearization_pi0star}.  Its principal symbol is
\cref{eq:symbol_a2_pi0star}, which is exactly
\cref{eq:apistar_eta2}, and it has no \(\Omega\)-component.  The quadratic term
\(\pi^\star_2\) is paralinearized in
\cref{lem:paralinearization_pistar2}; its order-one symbol is the one displayed
in \cref{eq:bmapistar1}, namely \cref{eq:apistar_eta1} and
\cref{eq:Apistar_Omega0}.

Combining these two decompositions gives \cref{eq:paralinearization_pistar}:
the order-zero \(\eta\)-symbols and the order \(-1\) \(\Omega\)-symbols are
absorbed into the generic symbols, the smoothing remainder of \(\pi_0\) gives
the separate linear term
\(r^{\pi_0,\bR}_{\bra{\varrho,\geq1}}
\pare{\sigma,\rho^\pm\middle|\eta}\eta\), and the
smoothing remainders of \(\pi_2^\star\) give the bilinear/linear terms displayed
in the last line of \cref{eq:paralinearization_pistar}.  The stated symbol
classes follow from
\cref{lem:paralinearization_pi0star,lem:paralinearization_pistar2} and from
\cref{lem:composition_rules}.
\end{proof}

\subsubsection{Paralinearization of $\pi$}
\label{sec:paralinearization_pi}

Let the function $\pi\pare{\rho^\pm\middle|\eta,\psi}$ be defined as in
\cref{eq:pidef_prop_recall}, namely
\begin{equation}\label{eq:pidef_prop_recall}
    \pi\pare{\rho^\pm\middle|\eta,\psi}
    \defeq
    \pi_0\pare{\rho^\pm\middle|\eta}
    +
    \pi_2\pare{\Atw\middle|\eta}\bra{\psi,\psi},
\end{equation}
where $\pi_0$ is defined in \cref{eq:pi0} and $\pi_2$ is defined as
\begin{equation}\label{eq:pi2_as_pistar2_composition}
    \pi_2\pare{\Atw\middle|\eta}\bra{\psi_1,\psi_2}
    =
    \pi^\star_2\pare{\Atw\middle|\eta}
    \bra{
        S\pare{\Atw\middle|\eta}\psi_1,
        S\pare{\Atw\middle|\eta}\psi_2
    } .
\end{equation}

\begin{prop}[Paralinearization of $\pi$]\label{prop:paralinearization_pi}
Let $\sigma\geq0$, $\rho^\pm\geq0$ with $\aangle{\rho^\pm}>0$, and
$\mathfrak g\in\bR$. There exists $\varrho_0 > 0$ such that for any
$\varrho \geq \varrho_0$, there exist
$s_0>0$ and $r_0>0$ such that for any
$\eta,\psi\in \BallR{s_0}{r_0}$, the function
$\pi\pare{\rho^\pm\middle|\eta,\psi}$ defined in
\cref{eq:pidef_prop_recall} admits the following paralinearization:
\begin{equation}\label{eq:paralinearization_pi}
\begin{aligned}
    \pi\pare{\rho^\pm\middle|\eta,\psi}
    =&
    -\mathfrak{g}\Atw\,\eta
    +
    \OpBW{
        A^\pi_2
        \pare{\rho^\pm\middle|\eta;x,\xi}
        +
        A^\pi_1
        \pare{\Atw\middle|\eta,\psi;x,\xi}
        +
        A_{\bra{0,\geq1}}
        \pare{\sigma,\rho^\pm,\Atw\middle|\eta,\psi;x,\xi}
    }\cdot\vect{\eta}{\psi}
    \\
    &+
    R^{\lin,\bR}_{\bra{\varrho,\geq1}}
    \pare{\sigma,\rho^\pm,\Atw\middle|\eta,\psi}
    \cdot\vect{\eta}{\psi}.
\end{aligned}
\end{equation}
where
\begin{itemize}
    \item $A^\pi_2\in\pare{\Gamma^2_{\geq0}}^2$ and is explicitly defined as
    \begin{equation}\label{eq:Api_2}
        A^\pi_2
        \pare{\rho^\pm\middle|\eta;x,\xi}
        \defeq
        \vect{
            a^{\pi,\eta}_2
            \pare{\rho^\pm\middle|\eta;x,\xi}
        }{
            0
        },
    \end{equation}
    where
    \begin{itemize}
        \item $a^{\pi,\eta}_2\in\Gamma^2_{\geq0}$ is real-valued and is explicitly defined as
        \begin{equation}\label{eq:Api_eta2}
            a^{\pi,\eta}_2
            \pare{\rho^\pm\middle|\eta;x,\xi}
            \defeq
            -\frac{\sigma}{\aangle{\rho^\pm}}
            \frac{
                \av{\xi}^2+\av{\nabla^\perp\eta\cdot\xi}^2
            }{
                \pare{1+\av{\nabla\eta}^2}^{\frac32}
            }\cdot
        \end{equation}
        In particular,
        \[
            a^{\pi,\eta}_2
            \pare{\rho^\pm\middle|0;x,\xi}
            =
            -\frac{\sigma}{\aangle{\rho^\pm}}\av{\xi}^2 .
        \]
    \end{itemize}

    \item $A^\pi_1\in\pare{\Gamma_{\geq1}^1}^2$ and is explicitly defined as
    \begin{equation}\label{eq:Api_star_1}
        A^\pi_1
        \pare{\Atw\middle|\eta,\psi;x,\xi}
        \defeq
        \vect{
            a^{\pi,\eta}_1
            \bra{\psi,\psi}
            \pare{\Atw\middle|\eta;x,\xi}
        }{
            a^{\pi,\psi}_1
            \bra{\psi}
            \pare{\Atw\middle|\eta;x,\xi}
        },
    \end{equation}
    where
    \begin{itemize}
        \item $a^{\pi,\eta}_1
        \in\cL^2\pare{H^{s_0};\Gamma^1}$ is real-valued and is explicitly defined as
        \begin{equation}\label{eq:Api_eta1}
            a^{\pi,\eta}_1
            \bra{\psi,\psi}
            \pare{\Atw\middle|\eta;x,\xi}
            \defeq
            \lambda_1\pare{\eta;x,\xi}
            \left[
                q^\pi_0\bra{\psi,\psi}
                \pare{\Atw\middle|\eta;x,\xi}
                -
                \pare{
                    b\bra{\psi}
                    \pare{\Atw\middle|\eta;x}
                }^2
            \right].
        \end{equation}
        In particular,
        $a^{\pi,\eta}_1\bra{\psi,\psi}\pare{\Atw\middle|\eta;x,\xi}$
        is even in $\xi$;

        \item $a^{\pi,\psi}_1
        \in\cL\pare{H^{s_0};\Gamma^1}$ is explicitly defined as
        \begin{equation}\label{eq:Api_psi1}
            a^{\pi,\psi}_1
            \bra{\psi}
            \pare{\Atw\middle|\eta;x,\xi}
            \defeq
            -\lambda_1\pare{\eta;x,\xi}
            b\bra{\psi}
            \pare{\Atw\middle|\eta;x}
            +
            \ii\,
            \mathpzc V
            \bra{\psi}
            \pare{\Atw\middle|\eta;x}
            \cdot\xi .
        \end{equation}
        In particular, the first term in \cref{eq:Api_psi1} is real-valued and even in $\xi$, while the second one is the transport symbol;

        \item $q^\pi_0
        \in\cL^2\pare{H^{s_0};\Gamma^0}$ is real-valued and is explicitly defined as
        \begin{equation}\label{eq:qpi0_definition}
            q^\pi_0\bra{\psi_1,\psi_2}
            \pare{\Atw\middle|\eta;x,\xi}
            \defeq
            \pare{4-\Atw^2}
            \pare{
                A_{0}^{\mathpzc{B}_0,\omega;\perp}\pare{\eta;x,\xi}
                \cdot
                S\pare{\Atw\middle|\eta}\psi_1
            }
            \pare{
                A_{0}^{\mathpzc{B}_0,\omega;\perp}\pare{\eta;x,\xi}
                \cdot
                S\pare{\Atw\middle|\eta}\psi_2
            } .
        \end{equation}
        In particular, $q^\pi_0$ is even in $\xi$;

        \item $b\bra{\psi}
        \in\cL\pare{H^{s_0};\cF^\bR}$ is real-valued and is explicitly defined as
        \begin{equation}\label{eq:definitionBb_pi_section}
            b\bra{\psi}
            \pare{\Atw\middle|\eta;x}
            \defeq
            \Atw\,
            v_3\bra{\psi}
            \pare{\Atw\middle|\eta;x}
            -
            \frac{
                \nabla^\perp\eta\cdot
                S\pare{\Atw\middle|\eta}\psi
            }{
                1+\av{\nabla\eta}^2
            }\cdot
        \end{equation}
        In particular, $b\bra{\psi}\pare{\Atw\middle|\eta;x}$ is independent of $\xi$;

        \item $\mathpzc V\bra{\psi}\pare{\Atw\middle|\eta}
        \in\cL\pare{H^{s_0};\pare{\cF^\bR}^2}$ is real-valued and is explicitly defined as
        \begin{equation}\label{eq:definitionVb_pi_section}
            \mathpzc V
            \bra{\psi}
            \pare{\Atw\middle|\eta;x}
            \defeq
            V\bra{\psi}\pare{\Atw\middle|\eta;x}
            +
            \frac{\Atw}{4\pare{1+\av{\nabla\eta}^2}}
            \bm J
            \pare{
                S\pare{\Atw\middle|\eta}\psi
                +
                \pare{
                    \nabla\eta\cdot
                    S\pare{\Atw\middle|\eta}\psi
                }\nabla\eta
            } .
        \end{equation}
        In particular, $\mathpzc V\bra{\psi}\pare{\Atw\middle|\eta;x}$ is independent of $\xi$.
    \end{itemize}

    \item $A_{\bra{0,\geq1}}$ is a generic vector symbol as per
    \cref{notation:paralinearization};

    \item $R^{\lin,\bR}_{\bra{\varrho,\geq1}}$ is a smoothing remainder as per
    \cref{notation:paralinearization}; its aggregate signature is the one
    displayed in \cref{eq:paralinearization_pi}.
\end{itemize}
\end{prop}

\begin{proof}
Fix the target gain $\varrho$. Apply
\cref{lem:paralinearization_pistar2,prop:paralinearization_S} with smoothing
index $\varrho+1$. If $s_\star$ is the Sobolev threshold in
\cref{lem:paralinearization_pistar2}, choose the present threshold so that
$s_0\geq s_\star+1$. Since $S\pare{\Atw\middle|\eta}$ has order one,
$\psi\in H^{s_0}$ then gives
$S\pare{\Atw\middle|\eta}\psi\in H^{s_\star}$; composing a remainder of gain
$\varrho+1$ with $S$ leaves gain $\varrho$. We relabel all such final
remainders as $R^{\lin,\bR}_{\bra{\varrho,\geq1}}$ below.

Set
\begin{equation}\label{eq:Omega_psi_definition_pi}
    \Omega_\psi
    \defeq
    S\pare{\Atw\middle|\eta}\psi .
\end{equation}
Then, by \cref{eq:pi2_as_pistar2_composition},
\[
    \pi_2\pare{\Atw\middle|\eta}\bra{\psi,\psi}
    =
    \pi^\star_2\pare{\Atw\middle|\eta}
    \bra{\Omega_\psi,\Omega_\psi}.
\]

\begin{proofpart}[Composition with $S\pare{\Atw\middle|\eta}$]
By \cref{lem:paralinearization_pistar2}, applied with
$\Omega=\Omega_\psi$, and by \cref{prop:paralinearization_S}, we have
\begin{equation}\label{eq:pi2_composition_raw}
\begin{aligned}
    \pi_2\pare{\Atw\middle|\eta}\bra{\psi,\psi}
    &=
    \OpBW{
        \vect{
            a^{\pi^\star_2,\eta}_1
            \bra{\Omega_\psi,\Omega_\psi}
            \pare{\Atw\middle|\eta;x,\xi}
        }{
            0
        }
    }\cdot\vect{\eta}{\psi}
    \\
    &
    +
    \OpBW{
        A^{\pi^\star_2,\Omega}_0
        \bra{\Omega_\psi}
        \pare{\Atw\middle|\eta;x,\xi}
    }\cdot
    S\pare{\Atw\middle|\eta}\psi
    \\
    &
    +
    \OpBW{
        A_{\bra{0,\geq1}}
        \pare{\Atw\middle|\eta,\psi;x,\xi}
    }\cdot\vect{\eta}{\psi}
    \\
    &
    +
    R^{\lin,\bR}_{\bra{\varrho,\geq1}}
    \pare{\Atw\middle|\eta,\psi}
    \cdot\vect{\eta}{\psi}.
\end{aligned}
\end{equation}
Moreover, \cref{eq:paralinearization_S} gives, at the principal orders needed here,
\begin{equation}\label{eq:S_columns_for_pi}
    \begin{aligned}
    S\pare{\Atw\middle|\eta}\psi
    &=
    \OpBW{
        \vectT{
            A^{S,\eta}_1\bra{\psi}
            \pare{\Atw\middle|\eta;x,\xi}
        }{
            \ii\xi^\perp
        }
        +
        \vectT{
            A_{\bra{0,\geq1}}\bra{\psi}
            \pare{\Atw\middle|\eta;x,\xi}
        }{
            A_{\bra{0,\geq1}}
            \pare{\Atw\middle|\eta;x,\xi}
        }
    }\vect{\eta}{\psi}\\&+
    R^{\lin,\bR}_{\bra{\varrho,\geq1}}
    \pare{\Atw\middle|\eta,\psi}
    \cdot\vect{\eta}{\psi},
    \end{aligned}
\end{equation}
because the correction
$-\ii\Atw\xi^\perp a^{\tau_0,\psi}_{-1}\pare{\eta;x,\xi}$
in \cref{eq:ASpsi} has order $0$.

From \cref{eq:ASeta1_sobolev} and \cref{eq:aSeta0}, the $\eta$-column in
\cref{eq:S_columns_for_pi} is
\begin{equation}\label{eq:ASeta_pi_rewrite}
    A^{S,\eta}_1\bra{\psi}
    \pare{\Atw\middle|\eta;x,\xi}
    =
    -\Atw
    \left(
        A_{0}^{\mathpzc{B}_0,\omega;\perp}\pare{\eta;x,\xi}\cdot\Omega_\psi
        -
        \ii\,
        v_3\bra{\psi}\pare{\Atw\middle|\eta;x}
    \right)\xi^\perp .
\end{equation}
Composing \cref{eq:pi2_composition_raw} with \cref{eq:S_columns_for_pi} and using the symbolic composition formula yields
\begin{equation}\label{eq:pi2_raw_vector_symbol}
\begin{aligned}
\pi_2\pare{\Atw\middle|\eta}\bra{\psi,\psi}
&=
    \OpBW{
        \vect{
            a^{\pi^\star_2,\eta}_1
            \bra{\Omega_\psi,\Omega_\psi}
            +
            A^{\pi^\star_2,\Omega}_0
            \bra{\Omega_\psi}
            \cdot
            A^{S,\eta}_1\bra{\psi}
        }{
            A^{\pi^\star_2,\Omega}_0
            \bra{\Omega_\psi}
            \cdot\ii\xi^\perp
        }
        +
        A_{\bra{0,\geq1}}
    }\cdot\vect{\eta}{\psi}\\&
    +
    R^{\lin,\bR}_{\bra{\varrho,\geq1}}
    \pare{\Atw\middle|\eta,\psi}
    \cdot\vect{\eta}{\psi}.
    \end{aligned}
\end{equation}
All the symbols in \cref{eq:pi2_raw_vector_symbol} are evaluated at
$\pare{\Atw\middle|\eta;x,\xi}$.
\end{proofpart}

\begin{proofpart}[$\psi$-component]
Using \cref{eq:pistar2_Omega_symbol}, \cref{eq:definitionBb_pi_section},
and \cref{eq:definitionVb_pi_section}, we compute
\begin{equation}\label{eq:pi2_psi_component_computation}
\begin{aligned}
    A^{\pi^\star_2,\Omega}_0
    \bra{\Omega_\psi}
    \pare{\Atw\middle|\eta;x,\xi}
    \cdot\ii\xi^\perp
    &=
    \ii\,
    \mathpzc V
    \bra{\psi}\pare{\Atw\middle|\eta;x}
    \cdot\xi+
    \pare{1+\av{\nabla\eta}^2}
    b\bra{\psi}
    \pare{\Atw\middle|\eta;x}
    A_{0}^{\mathpzc{B}_0,\omega;\perp}\pare{\eta;x,\xi}\cdot\xi^\perp .
\end{aligned}
\end{equation}
Since
\begin{equation}\label{eq:H0_xiperp_identity_pi}
    \pare{1+\av{\nabla\eta}^2}
    A_{0}^{\mathpzc{B}_0,\omega;\perp}\pare{\eta;x,\xi}\cdot\xi^\perp
    =
    \pare{1+\av{\nabla\eta}^2}
    A_{0}^{\mathpzc{B}_0,\omega}\pare{\eta;x,\xi}\cdot\xi
    =
    -\lambda_1\pare{\eta;x,\xi},
\end{equation}
we obtain
\begin{equation}\label{eq:pi2_psi_component_final}
    A^{\pi^\star_2,\Omega}_0
    \bra{\Omega_\psi}
    \pare{\Atw\middle|\eta;x,\xi}
    \cdot\ii\xi^\perp
    =
    -\lambda_1\pare{\eta;x,\xi}
    b\bra{\psi}
    \pare{\Atw\middle|\eta;x}
    +
    \ii\,
    \mathpzc V
    \bra{\psi}\pare{\Atw\middle|\eta;x}
    \cdot\xi .
\end{equation}
This is exactly \cref{eq:Api_psi1}.
\end{proofpart}

\begin{proofpart}[$\eta$-component]
We now simplify the first component of \cref{eq:pi2_raw_vector_symbol}. First,
by \cref{eq:pistar2_Omega_symbol} and \cref{eq:ASeta_pi_rewrite},
\begin{equation}\label{eq:eta_component_composition_pi}
\begin{aligned}
& A^{\pi^\star_2,\Omega}_0
    \bra{\Omega_\psi}
    \cdot
    A^{S,\eta}_1\bra{\psi}
 =
-\Atw
\left(
    A_{0}^{\mathpzc{B}_0,\omega;\perp}\pare{\eta;x,\xi}\cdot\Omega_\psi
    -
    \ii\,
    v_3\bra{\psi}\pare{\Atw\middle|\eta;x}
\right)
\left(
    \mathpzc V\bra{\psi}\pare{\Atw\middle|\eta;x}\cdot\xi
    +
    \ii\lambda_1\pare{\eta;x,\xi}
    b\bra{\psi}\pare{\Atw\middle|\eta;x}
\right),
\end{aligned}
\end{equation}
where \cref{eq:H0_xiperp_identity_pi} has been used in the second factor.

We shall also use the two pointwise identities
\begin{equation}\label{eq:pi_composition_identity_1}
    \frac{
        \bm J
        \pare{
            \Omega_\psi
            +
            \pare{\nabla\eta\cdot\Omega_\psi}\nabla\eta
        }\cdot\xi
    }{
        4\pare{1+\av{\nabla\eta}^2}
    }
    =
    \lambda_1\pare{\eta;x,\xi}
    \pare{
        A_{0}^{\mathpzc{B}_0,\omega;\perp}\pare{\eta;x,\xi}\cdot\Omega_\psi
    },
\end{equation}
and
\begin{equation}\label{eq:pi_composition_identity_2}
\begin{aligned}
\frac{1}{4\pare{1+\av{\nabla\eta}^2}}
\left(
    \pare{\Omega_\psi\cdot\nabla\eta}\Omega_\psi
    -
    \frac{
        \av{\Omega_\psi}^2+
        \av{\Omega_\psi\cdot\nabla\eta}^2
    }{
        1+\av{\nabla\eta}^2
    }
    \nabla\eta
\right)\cdot\xi
=
\lambda_1\pare{\eta;x,\xi}
\pare{
    A_{0}^{\mathpzc{B}_0,\omega;\perp}\pare{\eta;x,\xi}\cdot\Omega_\psi
}
\frac{
    \nabla^\perp\eta\cdot\Omega_\psi
}{
    1+\av{\nabla\eta}^2
}\cdot
\end{aligned}
\end{equation}
Both follow directly from the explicit formula for $A_{0}^{\mathpzc{B}_0,\omega}$ and from
\cref{eq:lambda1}.

Combining \cref{eq:pistar2_eta_symbol} with
\cref{eq:eta_component_composition_pi}, and using
\cref{eq:pi_composition_identity_1}--\cref{eq:pi_composition_identity_2},
the imaginary part cancels exactly. More precisely, the imaginary part is
\begin{equation}\label{eq:pi_eta_imaginary_cancellation}
\ii\Atw
\Bigg[
-\frac{1}{4\pare{1+\av{\nabla\eta}^2}}
\left(
    \pare{\Omega_\psi\cdot\nabla\eta}\Omega_\psi
    -
    \frac{
        \av{\Omega_\psi}^2+
        \av{\Omega_\psi\cdot\nabla\eta}^2
    }{
        1+\av{\nabla\eta}^2
    }
    \nabla\eta
\right)\cdot\xi
+
\lambda_1\pare{\eta;x,\xi}
\pare{
    A_{0}^{\mathpzc{B}_0,\omega;\perp}\pare{\eta;x,\xi}\cdot\Omega_\psi
}
\frac{
    \nabla^\perp\eta\cdot\Omega_\psi
}{
    1+\av{\nabla\eta}^2
}
\Bigg]
=0.
\end{equation}
The remaining real part is
\begin{equation}\label{eq:pi_eta_real_component_final}
\begin{aligned}
& a^{\pi^\star_2,\eta}_1
    \bra{\Omega_\psi,\Omega_\psi}
    \pare{\Atw\middle|\eta;x,\xi}
    +
    A^{\pi^\star_2,\Omega}_0
    \bra{\Omega_\psi}
    \pare{\Atw\middle|\eta;x,\xi}
    \cdot
    A^{S,\eta}_1\bra{\psi}
    \pare{\Atw\middle|\eta;x,\xi}
\\
=&
\lambda_1\pare{\eta;x,\xi}
\left[
    \pare{4-\Atw^2}
    \pare{
        A_{0}^{\mathpzc{B}_0,\omega;\perp}\pare{\eta;x,\xi}\cdot
        S\pare{\Atw\middle|\eta}\psi
    }^2
    -
    \left(
        \Atw\,
        v_3\bra{\psi}
        \pare{\Atw\middle|\eta;x}
        -
        \frac{
            \nabla^\perp\eta\cdot
            S\pare{\Atw\middle|\eta}\psi
        }{
            1+\av{\nabla\eta}^2
        }
    \right)^2
\right].
\end{aligned}
\end{equation}
By the definitions \cref{eq:qpi0_definition} and
\cref{eq:definitionBb_pi_section}, \cref{eq:pi_eta_real_component_final}
is precisely \cref{eq:Api_eta1}.
\end{proofpart}

\begin{proofpart}[Conclusion]
By \cref{lem:paralinearization_pi0star},
\begin{equation}\label{eq:pi0_inserted_in_pi_proof}
\begin{aligned}
    \pi_0\pare{\rho^\pm\middle|\eta}
    =&
    -\mathfrak{g}\Atw\,\eta
    +
    \OpBW{
        a^{\pi,\eta}_2
        \pare{\rho^\pm\middle|\eta;x,\xi}
        +
        a_{\bra{0,\geq1}}
        \pare{\rho^\pm\middle|\eta;x,\xi}
    }\eta\\&
    +
    r^{\pi_0,\bR}_{\bra{\varrho,\geq1}}
    \pare{\sigma,\rho^\pm\middle|\eta}\eta.
    \end{aligned}
\end{equation}
Combining \cref{eq:pi0_inserted_in_pi_proof},
\cref{eq:pi2_raw_vector_symbol},
\cref{eq:pi2_psi_component_final}, and
\cref{eq:pi_eta_real_component_final}, we obtain exactly
\cref{eq:paralinearization_pi}. The parity and reality properties stated in
the proposition follow from the fact that $\lambda_1$ is real and even,
$A_{0}^{\mathpzc{B}_0,\omega;\perp}$ is real and odd, while $b\bra{\psi}$ and
$\mathpzc V\bra{\psi}\pare{\Atw\middle|\eta}$ are real-valued and independent of $\xi$.
\end{proofpart}
\end{proof}

\subsubsection[Paralinearization of \texorpdfstring{$\tilde{\pi}$}{~pi}]{Paralinearization of $\tilde{\pi}$}

We now paralinearize the term $\tilde{\pi}$ which was defined in \cref{eq:tildepsiNEW}, this is a generalization of the paralinearization of $\pi$ provided in \cref{prop:paralinearization_pi} since $\pi$ can be recovered from $\tilde{\pi}$ in the limit $\bar{\Omega}\to 0$. The result we prove is the following

\begin{prop}[Paralinearization of $\tilde{\pi}$]\label{prop:paralinearization_pitilde}
Let $\sigma\geq0$, $\rho^\pm\geq0$ with $\aangle{\rho^\pm}>0$, and
$\mathfrak g\in\bR$. There exists $\varrho_0 > 0$ such that for any
$\varrho \geq \varrho_0$, there exist
$s_0>0$ and $r_0>0$ such that for any
$\eta,\tilde{\psi}\in \BallR{s_0}{r_0}$ and any $\background\in\bR^2$, the shifted
Bernoulli function
$\tilde{\pi}\pare{\rho^\pm,\background\middle|\eta,\tilde{\psi}}$ defined in
\cref{eq:tildepsiNEW} admits the following paralinearization:
\begin{equation}\label{eq:paralinearization_pitilde}
    \begin{aligned}\tilde{\pi}\pare{\rho^\pm,\background\middle|\eta,\tilde{\psi}}
    =
    &-\mathfrak{g}\Atw\,\eta
    +
    \OpBW{
        A^{\tilde{\pi}}_2
        \pare{\rho^\pm\middle|\eta;x,\xi}
        +
        A^{\tilde{\pi}}_1
        \pare{\Atw,\background\middle|\eta,\tilde{\psi};x,\xi}
        +
        A_{\bra{0,\geq1}}
        \pare{\sigma,\rho^\pm,\Atw,\background\middle|\eta,\tilde{\psi};x,\xi}
    }\cdot\vect{\eta}{\tilde{\psi}}
    \\
    &+
    R^\bR_{\bra{\varrho,\geq1}}
    \pare{\sigma,\rho^\pm,\Atw,\background\middle|\eta,\tilde{\psi}}
    \cdot\vect{\eta}{\tilde{\psi}},
\end{aligned}
\end{equation}
where
\begin{itemize}
    \item $A^{\tilde{\pi}}_2\in\pare{\Gamma^2_{\geq0}}^2$ and is explicitly defined as
    \begin{equation}\label{eq:Apitilde_2}
        A^{\tilde{\pi}}_2
        \pare{\rho^\pm\middle|\eta;x,\xi}
        \defeq
        \vect{
            a^{\pi_0,\eta}_2
            \pare{\rho^\pm\middle|\eta;x,\xi}
        }{
            0
        },
    \end{equation}
    where
    \begin{itemize}
        \item $a^{\pi_0,\eta}_2\in\Gamma^2_{\geq0}$ is real-valued and is explicitly defined in
        \cref{eq:symbol_a2_pi0star}. In particular,
        \[
            a^{\pi_0,\eta}_2
            \pare{\rho^\pm\middle|0;x,\xi}
            =
            -\frac{\sigma}{\aangle{\rho^\pm}}\av{\xi}^2 .
        \]
    \end{itemize}

    \item $A^{\tilde{\pi}}_1\in\pare{\Gamma^1_{\geq0}}^2$ and is explicitly defined as
    \begin{equation}\label{eq:Apitilde_1}
        A^{\tilde{\pi}}_1
        \pare{\Atw,\background\middle|\eta,\tilde{\psi};x,\xi}
        \defeq
        \vect{
            a^{\tilde{\pi},\eta}_1
            \pare{\Atw,\background\middle|\eta,\tilde{\psi};x,\xi}
        }{
            a^{\tilde{\pi},\tilde{\psi}}_1
            \pare{\Atw,\background\middle|\eta,\tilde{\psi};x,\xi}
        },
    \end{equation}
    where
    \begin{itemize}
        \item $a^{\tilde{\pi},\eta}_1\in\Gamma^1_{\geq0}$ is real-valued and is explicitly defined as
        \begin{equation}\label{eq:apitilde_eta1}
\begin{aligned}
a^{\tilde{\pi},\eta}_1
\pare{\Atw,\background\middle|\eta,\tilde{\psi};x,\xi}
&\defeq
\lambda_1\pare{\eta;x,\xi}
\left[
    q^{\tilde{\pi}}_0
    \pare{\Atw,\background\middle|\eta,\tilde{\psi};x,\xi}
    -
    \pare{
        b_{\textnormal{can}}
        \pare{\Atw,\background\middle|\eta,\tilde{\psi};x}
    }^2
\right].
\end{aligned}
        \end{equation}
        In particular,
        $a^{\tilde{\pi},\eta}_1
        \pare{\Atw,\background\middle|\eta,\tilde{\psi};x,\xi}$
        is even in $\xi$;

        \item $a^{\tilde{\pi},\tilde{\psi}}_1\in\Gamma^1_{\geq0}$ is explicitly defined as
        \begin{equation}\label{eq:apitilde_psitilde1}
            a^{\tilde{\pi},\tilde{\psi}}_1
            \pare{\Atw,\background\middle|\eta,\tilde{\psi};x,\xi}
            \defeq
            -\lambda_1\pare{\eta;x,\xi}
            b_{\textnormal{can}}
            \pare{\Atw,\background\middle|\eta,\tilde{\psi};x}
            +
            \ii\,
            \mathpzc V_{\textnormal{can}}
            \pare{\Atw,\background\middle|\eta,\tilde{\psi};x}
            \cdot\xi .
        \end{equation}
        In particular, the first term in \cref{eq:apitilde_psitilde1} is real-valued and even in $\xi$, while the second one is the transport symbol;

        \item $q^{\tilde{\pi}}_0\in\Gamma^0_{\geq0}$ is real-valued and is explicitly defined as
        \begin{equation}\label{eq:qpitilde0_definition}
            q^{\tilde{\pi}}_0
            \pare{\Atw,\background\middle|\eta,\tilde{\psi};x,\xi}
            \defeq
            \pare{4-\Atw^2}
            \left(
                A_{0}^{\mathpzc{B}_0,\omega;\perp}\pare{\eta;x,\xi}
                \cdot
                \Omega
                \pare{\background\middle|\eta,\tilde{\psi};x}
            \right)^2 .
        \end{equation}
        In particular, $q^{\tilde{\pi}}_0$ is even in $\xi$;

        \item $\Omega\pare{\background\middle|\eta,\tilde{\psi}}$ is explicitly defined as
        \begin{equation}\label{eq:Omegacan_pitilde}
            \Omega
            \pare{\background\middle|\eta,\tilde{\psi};x}
            \defeq
            \tilde{\Omega}\bra{\background}\pare{\eta;x}
            +
            S\pare{\Atw\middle|\eta}\tilde{\psi}\pare{x}.
        \end{equation}

        \item $b_{\textnormal{can}}\in\cF^\bR$ is real-valued and is explicitly defined as
        \begin{equation}\label{eq:bpitilde_can}
\begin{aligned}
            b_{\textnormal{can}}
            \pare{\Atw,\background\middle|\eta,\tilde{\psi};x}
            &\defeq
            \Atw\,
             v_3^{\star}
            \bra{
                \Omega
                \pare{\background\middle|\eta,\tilde{\psi}}
            }
            \pare{\eta;x}-
            \frac{
                \nabla^\perp\eta\cdot
                \Omega
                \pare{\background\middle|\eta,\tilde{\psi};x}
            }{
                1+\av{\nabla\eta}^2
            }\cdot
\end{aligned}
        \end{equation}
        In particular,
        \[
            b_{\textnormal{can}}
            \pare{\Atw,\background\middle|0,0;x}
            =
            0 ;
        \]

        \item $\mathpzc V_{\textnormal{can}}\pare{\Atw,\background\middle|\eta,\tilde{\psi}}\in
        \pare{\cF^\bR_{\geq0}}^2$ is real-valued and is explicitly defined as
        \begin{equation}\label{eq:Vpitilde_can}
\begin{aligned}
            \mathpzc V_{\textnormal{can}}
            \pare{\Atw,\background\middle|\eta,\tilde{\psi};x}
            &\defeq
             V^{\star}
            \bra{
                \Omega
                \pare{\background\middle|\eta,\tilde{\psi}}
            }
            \pare{\eta;x}
            \\
            &\quad
            +
            \frac{\Atw}{4\pare{1+\av{\nabla\eta}^2}}
            \bm J
            \pare{
                \Omega
                \pare{\background\middle|\eta,\tilde{\psi};x}
                +
                \pare{
                    \nabla\eta\cdot
                    \Omega
                    \pare{\background\middle|\eta,\tilde{\psi};x}
                }\nabla\eta
            } .
\end{aligned}
        \end{equation}
        In particular, its nonzero flat part is separated as
        \begin{equation}\label{eq:Vpitilde_can_split}
            \mathpzc V_{\textnormal{can}}
            \pare{\Atw,\background\middle|\eta,\tilde{\psi}}
            =
            \frac{\Atw}{4}\background^\perp
            +
            \mathpzc V_{\textnormal{can}}^{\mathrm{nl}}
            \pare{\Atw,\background\middle|\eta,\tilde{\psi}},
            \qquad
            \mathpzc V_{\textnormal{can}}^{\mathrm{nl}}
            \in\pare{\cF^\bR}^2.
        \end{equation}
        Thus
        $\mathpzc V_{\textnormal{can}}
        \pare{\Atw,\background\middle|\eta,\tilde{\psi};x}$
        is independent of $\xi$ and 
        $$\mathpzc V_{\textnormal{can}}
        \pare{\Atw,\background\middle|0,0;x}=\frac{\Atw}{4}\background^{\perp};$$
    \end{itemize}

    \item $A_{\bra{0,\geq1}}$ is a generic vector symbol as per
    \cref{notation:paralinearization};

    \item $R^\bR_{\bra{\varrho,\geq1}}$ is a smoothing remainder as per
    \cref{notation:paralinearization}.
\end{itemize}
\end{prop}

\paragraph[Paralinearization of pi-tilde 0]{Paralinearization of $\tilde{\pi}_0$}

We isolate the background contribution to the zero-th order canonical pressure. Namely, with
$\tilde{\Omega}\bra{\background}\pare{\eta}$ defined in \cref{eq:Omega_bg_can_def}, we set
\begin{equation}\label{eq:pitilde0_recall}
    \tilde{\pi}_0\bra{\background,\background}\pare{\Atw\middle|\eta}
    \defeq
    \pi^\star_2\pare{\Atw\middle|\eta}
    \bra{
        \tilde{\Omega}\bra{\background}\pare{\eta},
        \tilde{\Omega}\bra{\background}\pare{\eta}
    } .
\end{equation}
Equivalently, the full zero-th order canonical term in \cref{eq:pi0tilde} is
\[
    \pi_{\mathrm{can},0}\pare{\rho^\pm,\background\middle|\eta}
    =
    \pi_0\pare{\rho^\pm\middle|\eta}
    +
    \tilde{\pi}_0\bra{\background,\background}\pare{\Atw\middle|\eta}
    -
    c_{\background}.
\]
The local background term has value $c_{\background}$ at $\eta=0$; this
constant must be retained in its paralinearization and cancels the explicit
$-c_{\background}$ in the full canonical term.

\begin{lemma}[Paralinearization of $\tilde{\pi}_0$]\label{lem:paralinearization_pitilde0}
Let $\sigma\geq0$, $\rho^\pm\geq0$ with $\aangle{\rho^\pm}>0$, and
$\mathfrak g\in\bR$. There exists $\varrho_0 > 0$ such that for any
$\varrho \geq \varrho_0$, there exist
$s_0>0$ and $r_0>0$ such that for any
$\eta\in \BallR{s_0}{r_0}$ and any $\background\in\bR^2$, the term
$\tilde{\pi}_0\bra{\background,\background}\pare{\Atw\middle|\eta}$ defined in
\cref{eq:pitilde0_recall} admits the following paralinearization:
\begin{equation}\label{eq:paralinearization_pitilde0}
    \tilde{\pi}_0\bra{\background,\background}\pare{\Atw\middle|\eta}
    =
    c_{\background}
    +
    \OpBW{
        a^{\tilde{\pi}_0,\eta}_1
        \bra{\background,\background}
        \pare{\Atw\middle|\eta;x,\xi}
        +
        a_{\bra{0,\geq1}}
        \bra{\background,\background}
        \pare{\Atw\middle|\eta;x,\xi}
    }\eta
    +
    r^\bR_{\bra{\varrho,\geq1}}
    \bra{\background,\background}
    \pare{\Atw\middle|\eta}\eta,
\end{equation}
where
\begin{itemize}
    \item $a^{\tilde{\pi}_0,\eta}_1
    \in \cL^2\pare{\bR^2;\Gamma^1_{\geq0}}$ is real-valued and is explicitly defined as
    \begin{equation}\label{eq:apitilde0_eta1}
    a^{\tilde{\pi}_0,\eta}_1
    \bra{\background,\background}
    \pare{\Atw\middle|\eta;x,\xi}
    \defeq
    \lambda_1\pare{\eta;x,\xi}
    \left(
        q^{\tilde{\pi}_0}_0
        \bra{\background,\background}
        \pare{\Atw\middle|\eta;x,\xi}
        -
        \pare{
            b_{\textnormal{bg}}\bra{\background}
            \pare{\Atw\middle|\eta;x}
        }^2
    \right),
\end{equation}
    where
    \begin{itemize}
        \item $q^{\tilde{\pi}_0}_0
        \in \cL^2\pare{\bR^2;\Gamma^0_{\geq0}}$ is real-valued and is explicitly defined as
        \begin{equation}\label{eq:qpitilde0_bg}
    q^{\tilde{\pi}_0}_0
    \bra{\background,\background}
    \pare{\Atw\middle|\eta;x,\xi}
    \defeq
    \pare{4-\Atw^2}
    \left(
        A_{0}^{\mathpzc{B}_0,\omega;\perp}\pare{\eta;x,\xi}
        \cdot
        \tilde{\Omega}\bra{\background}\pare{\eta;x}
    \right)^2;
\end{equation}

        \item $b_{\textnormal{bg}}\bra{\background}
        \in \cL\pare{\bR^2;\cF^\bR}$ is real-valued and is explicitly defined as
        \begin{equation}\label{eq:bpitilde0_bg}
    b_{\textnormal{bg}}\bra{\background}
    \pare{\Atw\middle|\eta;x}
    \defeq
    \Atw\,
     v_3^{\star}
    \bra{\tilde{\Omega}\bra{\background}\pare{\eta}}
    \pare{\eta;x}
    -
    \frac{
        \nabla^\perp\eta\cdot
        \tilde{\Omega}\bra{\background}\pare{\eta;x}
    }{
        1+\av{\nabla\eta}^2
    };
\end{equation}
    \end{itemize}

    \item $a_{\bra{0,\geq1}}$ is a generic symbol as per
    \cref{notation:paralinearization};

    \item $r^\bR_{\bra{\varrho,\geq1}}$ is a smoothing remainder as per
    \cref{notation:paralinearization}.
    \end{itemize}
Consequently, the full zero-th order canonical term satisfies
\begin{equation}\label{eq:paralinearization_pican0}
\begin{aligned}
    \pi_{\mathrm{can},0}
    \pare{\rho^\pm,\background\middle|\eta}
    &=
    -\mathfrak{g}\Atw\,\eta
    +
    \OpBW{
        a^{\pi_0,\eta}_2
        \pare{\rho^\pm\middle|\eta;x,\xi}
        +
        a^{\tilde{\pi}_0,\eta}_1
        \bra{\background,\background}
        \pare{\Atw\middle|\eta;x,\xi}
        +
        a_{\bra{0,\geq1}}
        \bra{\background,\background}
        \pare{\sigma,\rho^\pm,\Atw\middle|\eta;x,\xi}
    }\eta
    \\
    &\quad
    +
    r^\bR_{\bra{\varrho,\geq1}}
    \bra{\background,\background}
    \pare{\sigma,\rho^\pm,\Atw\middle|\eta}\eta,
\end{aligned}
\end{equation}
where $a^{\pi_0,\eta}_2$ is defined in \cref{eq:symbol_a2_pi0star}.
\end{lemma}

\begin{proof}
Set
\[
    \Omega_{\textnormal{bg}}
    \defeq
    \tilde{\Omega}\bra{\background}\pare{\eta},
    \qquad
    \delta\Omega_{\textnormal{bg}}
    \defeq
    \Omega_{\textnormal{bg}}-\background.
\]
By \cref{eq:paralinearization_S_background,eq:ASeta1_background},
$\delta\Omega_{\textnormal{bg}}$ is the Sobolev correction whose principal
$\eta$-symbol is
$A^{S,\eta}_{1,\textnormal{bg}}\bra{\background}$. We first use the
bilinearity of $\pi_2^\star$ to split
\begin{equation}\label{eq:pitilde0_background_bilinear_split}
\begin{aligned}
\pi_2^\star\pare{\Atw\middle|\eta}
\bra{\Omega_{\textnormal{bg}},\Omega_{\textnormal{bg}}}
&=
\pi_2^\star\pare{\Atw\middle|\eta}\bra{\background,\background}
\\
&\quad
+2\pi_2^\star\pare{\Atw\middle|\eta}
\bra{\background,\delta\Omega_{\textnormal{bg}}}
\\
&\quad
+\pi_2^\star\pare{\Atw\middle|\eta}
\bra{\delta\Omega_{\textnormal{bg}},\delta\Omega_{\textnormal{bg}}}.
\end{aligned}
\end{equation}
The constant--constant block equals $c_{\background}$ at $\eta=0$ by
\cref{eq:c_background}; its difference from $c_{\background}$ vanishes at
the base point and is paralinearized in $\eta$ with $\background$ as an
external parameter. In the two remaining blocks we use the bilinear formula
of \cref{cor:paralinearization_pistar2_bilinear}, applying its Sobolev
estimates only to $\delta\Omega_{\textnormal{bg}}$ and using
\cref{rem:action_paradiff_constants} for every direct paradifferential action
on $\background$. Composing the Sobolev factors with the canonical
paralinearization just cited and resumming the three blocks gives
\begin{equation}\label{eq:pitilde0_raw_symbol}
\begin{aligned}
\tilde{\pi}_0\bra{\background,\background}\pare{\Atw\middle|\eta}
&=
c_{\background}
+
\OpBW{
    a^{\pi^\star_2,\eta}_1
    \bra{\Omega_{\textnormal{bg}},\Omega_{\textnormal{bg}}}
    \pare{\Atw\middle|\eta;x,\xi}
}
\eta
\\
&\quad
+
\OpBW{
    A^{\pi^\star_2,\Omega}_0
    \bra{\Omega_{\textnormal{bg}}}
    \pare{\Atw\middle|\eta;x,\xi}
    \cdot
    A^{S,\eta}_{1,\textnormal{bg}}
    \bra{\background}
    \pare{\Atw\middle|\eta;x,\xi}
}
\eta
\\
&\quad
+
\OpBW{
    a_{\bra{0,\geq1}}
    \bra{\background,\background}
    \pare{\Atw\middle|\eta;x,\xi}
}
\eta
+
r^\bR_{\bra{\varrho,\geq1}}
\bra{\background,\background}
\pare{\Atw\middle|\eta}\eta.
\end{aligned}
\end{equation}

We now simplify the order-one symbol in \cref{eq:pitilde0_raw_symbol}. From
\cref{eq:pistar2_Omega_symbol} and \cref{eq:ASeta1_background}, one has
\begin{equation}\label{eq:pitilde0_composition_symbol}
\begin{aligned}
    &A^{\pi^\star_2,\Omega}_0
    \bra{\Omega_{\textnormal{bg}}}
    \cdot
    A^{S,\eta}_{1,\textnormal{bg}}
    \bra{\background}
\\
=
&-\Atw
\left(
    A_{0}^{\mathpzc{B}_0,\omega;\perp}\cdot\Omega_{\textnormal{bg}}
    -
    \ii v_3^{\star}\bra{\Omega_{\textnormal{bg}}}
\right)
\left[
     V^{\star}\bra{\Omega_{\textnormal{bg}}}\cdot\xi
    +
    \frac{\Atw}{
        4\pare{1+\av{\nabla\eta}^2}
    }
    \bm J
    \pare{
        \Omega_{\textnormal{bg}}
        +
        \pare{\nabla\eta\cdot\Omega_{\textnormal{bg}}}\nabla\eta
    }
    \cdot\xi
\right.
\\
   &\left. +
    \ii\lambda_1
    \left(
        \Atw v_3^{\star}\bra{\Omega_{\textnormal{bg}}}
        -
        \frac{
            \nabla^\perp\eta\cdot\Omega_{\textnormal{bg}}
        }{
            1+\av{\nabla\eta}^2
        }
    \right)
\right],
\end{aligned}
\end{equation}
where all the symbols in the right-hand side are evaluated at
$\pare{\Atw\middle|\eta;x,\xi}$ and all functions at $x$.

The two pointwise identities
\begin{equation}\label{eq:pitilde0_identity_1}
    \frac{
        \bm J
        \pare{
            \Omega_{\textnormal{bg}}
            +
            \pare{\nabla\eta\cdot\Omega_{\textnormal{bg}}}\nabla\eta
        }
        \cdot\xi
    }{
        4\pare{1+\av{\nabla\eta}^2}
    }
    =
    \lambda_1\pare{\eta;x,\xi}
    \left(
        A_{0}^{\mathpzc{B}_0,\omega;\perp}\pare{\eta;x,\xi}
        \cdot
        \Omega_{\textnormal{bg}}
    \right)
\end{equation}
and
\begin{equation}\label{eq:pitilde0_identity_2}
\begin{aligned}
&\frac{1}{4\pare{1+\av{\nabla\eta}^2}}
\left(
    \pare{\Omega_{\textnormal{bg}}\cdot\nabla\eta}\Omega_{\textnormal{bg}}
    -
    \frac{
        \av{\Omega_{\textnormal{bg}}}^2
        +
        \av{\Omega_{\textnormal{bg}}\cdot\nabla\eta}^2
    }{
        1+\av{\nabla\eta}^2
    }
    \nabla\eta
\right)\cdot\xi
=
\lambda_1\pare{\eta;x,\xi}
\left(
    A_{0}^{\mathpzc{B}_0,\omega;\perp}\pare{\eta;x,\xi}
    \cdot
    \Omega_{\textnormal{bg}}
\right)
\frac{
    \nabla^\perp\eta\cdot\Omega_{\textnormal{bg}}
}{
    1+\av{\nabla\eta}^2
}
\end{aligned}
\end{equation}
follow directly from the explicit formula for $A_{0}^{\mathpzc{B}_0,\omega}$ and \cref{eq:lambda1}.

Combining \cref{eq:pistar2_eta_symbol} and
\cref{eq:pitilde0_composition_symbol}, the imaginary part of the resulting
order-one symbol is
\begin{equation}
\begin{aligned}
&-\Atw
\left[
\frac{1}{4\pare{1+\av{\nabla\eta}^2}}
\left(
    \pare{\Omega_{\textnormal{bg}}\cdot\nabla\eta}\Omega_{\textnormal{bg}}
    -
    \frac{
        \av{\Omega_{\textnormal{bg}}}^2
        +
        \av{\Omega_{\textnormal{bg}}\cdot\nabla\eta}^2
    }{
        1+\av{\nabla\eta}^2
    }
    \nabla\eta
\right)\cdot\xi
\right.
\left.
-
\lambda_1
\left(
    A_{0}^{\mathpzc{B}_0,\omega;\perp}\cdot\Omega_{\textnormal{bg}}
\right)
\frac{
    \nabla^\perp\eta\cdot\Omega_{\textnormal{bg}}
}{
    1+\av{\nabla\eta}^2
}
\right]
\\
&-\Atw^2
 v_3^{\star}\bra{\Omega_{\textnormal{bg}}}
\left[
    \lambda_1
    \left(
        A_{0}^{\mathpzc{B}_0,\omega;\perp}\cdot\Omega_{\textnormal{bg}}
    \right)
    -
    \frac{
        \bm J
        \pare{
            \Omega_{\textnormal{bg}}
            +
            \pare{\nabla\eta\cdot\Omega_{\textnormal{bg}}}\nabla\eta
        }
        \cdot\xi
    }{
        4\pare{1+\av{\nabla\eta}^2}
    }
\right],
\end{aligned}
\end{equation}
which vanishes by \cref{eq:pitilde0_identity_1} and
\cref{eq:pitilde0_identity_2}. Hence the principal symbol is real.

Taking the real part and using again \cref{eq:pitilde0_identity_1}, we obtain
\begin{align}
&\Re
\left[
    a^{\pi^\star_2,\eta}_1
    \bra{\Omega_{\textnormal{bg}},\Omega_{\textnormal{bg}}}
    +
    A^{\pi^\star_2,\Omega}_0
    \bra{\Omega_{\textnormal{bg}}}
    \cdot
    A^{S,\eta}_{1,\textnormal{bg}}
    \bra{\background}
\right]
=
\lambda_1
\left[
    \pare{4-\Atw^2}
    \left(
        A_{0}^{\mathpzc{B}_0,\omega;\perp}\cdot\Omega_{\textnormal{bg}}
    \right)^2
    -
    \left(
        \Atw v_3^{\star}\bra{\Omega_{\textnormal{bg}}}
        -
        \frac{
            \nabla^\perp\eta\cdot\Omega_{\textnormal{bg}}
        }{
            1+\av{\nabla\eta}^2
        }
    \right)^2
\right].
\end{align}
This is exactly \cref{eq:apitilde0_eta1}. The generic order-zero symbol and the
smoothing remainder in \cref{eq:paralinearization_pitilde0} collect the lower
symbolic terms from \cref{eq:paralinearization_S_background}, from
\cref{lem:paralinearization_pistar2}, and from the Bony--Weyl composition
formula. Finally, \cref{eq:paralinearization_pican0} follows by adding
\cref{eq:paralinearization_pi0}: the $c_{\background}$ in
\cref{eq:paralinearization_pitilde0} cancels the explicit
$-c_{\background}$ in the definition of $\pi_{\mathrm{can},0}$ from
\cref{eq:pi0tilde}.
\end{proof}

\paragraph[Paralinearization of pi-tilde 1]{Paralinearization of $\tilde{\pi}_1$}

We now paralinearize the term which is linear in the canonical perturbation
$\tilde{\psi}$. Recalling the decomposition \cref{eq:pi0tilde}, we set
\begin{equation}\label{eq:pitilde1_recall}
    \tilde{\pi}_1
    \pare{\Atw,\background\middle|\eta}
    \bra{\tilde{\psi}}
    \defeq
    2\,
    \pi^\star_2\pare{\Atw\middle|\eta}
    \bra{
        \tilde{\Omega}\bra{\background}\pare{\eta},
        S\pare{\Atw\middle|\eta}\tilde{\psi}
    } .
\end{equation}

\begin{lemma}[Paralinearization of $\tilde{\pi}_1$]\label{lem:paralinearization_pitilde1}
There exists $\varrho_0 > 0$ such that for any $\varrho \geq \varrho_0$, there exist
$s_0>0$ and $r_0>0$ such that for any
$\eta,\tilde{\psi}\in \BallR{s_0}{r_0}$ and any $\background\in\bR^2$, the term
$\tilde{\pi}_1\pare{\Atw,\background\middle|\eta}\bra{\tilde{\psi}}$
defined in \cref{eq:pitilde1_recall} admits the following paralinearization:
\begin{equation}\label{eq:paralinearization_pitilde1}
\begin{aligned}
    \tilde{\pi}_1
    \pare{\Atw,\background\middle|\eta}
    \bra{\tilde{\psi}}
    &=
    \OpBW{
        A^{\tilde{\pi}_1}_1
        \pare{\Atw,\background\middle|\eta,\tilde{\psi};x,\xi}
        +
        A_{\bra{0,\geq1}}
        \pare{\Atw,\background\middle|\eta,\tilde{\psi};x,\xi}
    }\cdot\vect{\eta}{\tilde{\psi}}
    \\
    &\quad
    +
    R^\bR_{\bra{\varrho,\geq1}}
    \pare{\Atw,\background\middle|\eta,\tilde{\psi}}
    \cdot\vect{\eta}{\tilde{\psi}},
\end{aligned}
\end{equation}
where
\begin{itemize}
    \item $A^{\tilde{\pi}_1}_1\in\pare{\Gamma^1_{\geq0}}^2$ and is explicitly defined as
    \begin{equation}\label{eq:Apitilde1_1}
        A^{\tilde{\pi}_1}_1
        \pare{\Atw,\background\middle|\eta,\tilde{\psi};x,\xi}
        \defeq
        \vect{
            a^{\tilde{\pi}_1,\eta}_1
            \bra{\background,\tilde{\psi}}
            \pare{\Atw\middle|\eta;x,\xi}
        }{
            a^{\tilde{\pi}_1,\tilde{\psi}}_1
            \bra{\background}
            \pare{\Atw\middle|\eta;x,\xi}
        },
    \end{equation}
    where
    \begin{itemize}
        \item $a^{\tilde{\pi}_1,\eta}_1
        \in \cL\pare{\bR^2;\cL\pare{H^{s_0};\Gamma^1_{\geq0}}}$ is real-valued and is explicitly defined as
    \begin{equation}\label{eq:apitilde1_eta1}
\begin{aligned}
a^{\tilde{\pi}_1,\eta}_1
\bra{\background,\tilde{\psi}}
\pare{\Atw\middle|\eta;x,\xi}
&\defeq
2\lambda_1\pare{\eta;x,\xi}
\Bigg[
    \pare{4-\Atw^2}
    \pare{
        A_{0}^{\mathpzc{B}_0,\omega;\perp}\pare{\eta;x,\xi}
        \cdot
        \tilde{\Omega}\bra{\background}\pare{\eta;x}
    }
\\
&\qquad\qquad\qquad\qquad
    \times
    \pare{
        A_{0}^{\mathpzc{B}_0,\omega;\perp}\pare{\eta;x,\xi}
        \cdot
        S\pare{\Atw\middle|\eta}\tilde{\psi}\pare{x}
    }
\\
&\qquad
    -
    b_{\textnormal{bg}}\bra{\background}
    \pare{\Atw\middle|\eta;x}
    b_{\tilde{\psi}}\bra{\tilde{\psi}}
    \pare{\Atw\middle|\eta;x}
\Bigg].
\end{aligned}
\end{equation}
        In particular,
        $a^{\tilde{\pi}_1,\eta}_1
        \bra{\background,\tilde{\psi}}
        \pare{\Atw\middle|\eta;x,\xi}$
        is real-valued and even in $\xi$;

        \item $a^{\tilde{\pi}_1,\tilde{\psi}}_1
        \in \cL\pare{\bR^2;\Gamma^1_{\geq0}}$ is explicitly defined as
        \begin{equation}\label{eq:apitilde1_psitilde1}
    a^{\tilde{\pi}_1,\tilde{\psi}}_1
    \bra{\background}
    \pare{\Atw\middle|\eta;x,\xi}
    \defeq
    -\lambda_1\pare{\eta;x,\xi}
    b_{\textnormal{bg}}\bra{\background}
    \pare{\Atw\middle|\eta;x}
    +
    \ii\,
    \mathpzc V_{\textnormal{bg}}
    \bra{\background}
    \pare{\Atw\middle|\eta;x}
    \cdot\xi .
\end{equation}
        In particular, the first term in
        \cref{eq:apitilde1_psitilde1} is real-valued and even in $\xi$,
        while the second one is the transport symbol generated by the
        canonical background;

        \item $b_{\textnormal{bg}}\bra{\background}
        \in \cL\pare{\bR^2;\cF^\bR}$ is the real-valued coefficient defined
        in \cref{eq:bpitilde0_bg}. In particular,
        \[
            b_{\textnormal{bg}}\bra{\background}
            \pare{\Atw\middle|0;x}
            =
            0;
        \]

        \item $b_{\tilde{\psi}}\bra{\tilde{\psi}}
        \in \cL\pare{H^{s_0};\cF^\bR}$ is real-valued and is explicitly defined as
        \begin{equation}\label{eq:bpitilde1_psitilde}
            b_{\tilde{\psi}}\bra{\tilde{\psi}}
            \pare{\Atw\middle|\eta;x}
            \defeq
            \Atw\,
            v_3\bra{\tilde{\psi}}
            \pare{\Atw\middle|\eta;x}
            -
            \frac{
                \nabla^\perp\eta\cdot
                S\pare{\Atw\middle|\eta}\tilde{\psi}\pare{x}
            }{
                1+\av{\nabla\eta}^2
            }\cdot
        \end{equation}
        In particular,
        $b_{\tilde{\psi}}\bra{\tilde{\psi}}
        \pare{\Atw\middle|\eta;x}$ is independent of $\xi$;

        \item $\mathpzc V_{\textnormal{bg}}\bra{\background}\pare{\Atw\middle|\eta}
        \in \cL\pare{\bR^2;\pare{\cF^\bR_{\geq0}}^2}$ is real-valued and is explicitly defined as
        \begin{equation}\label{eq:Vpitilde1_bg}
\begin{aligned}
    \mathpzc V_{\textnormal{bg}}
    \bra{\background}
    \pare{\Atw\middle|\eta;x}
    &\defeq
     V^{\star}
    \bra{\tilde{\Omega}\bra{\background}\pare{\eta}}
    \pare{\eta;x}
    \\
    &\quad
    +
    \frac{\Atw}{4\pare{1+\av{\nabla\eta}^2}}
    \bm J
    \pare{
        \tilde{\Omega}\bra{\background}\pare{\eta;x}
        +
        \pare{
            \nabla\eta\cdot
            \tilde{\Omega}\bra{\background}\pare{\eta;x}
        }\nabla\eta
    } .
\end{aligned}
        \end{equation}
        Its flat part and nonlinear correction are
        \begin{equation}\label{eq:Vpitilde1_bg_split}
            \mathpzc V_{\textnormal{bg}}\bra{\background}
            \pare{\Atw\middle|\eta}
            =
            \frac{\Atw}{4}\background^\perp
            +
            \mathpzc V_{\textnormal{bg}}^{\mathrm{nl}}
            \bra{\background}\pare{\Atw\middle|\eta},
            \qquad
            \mathpzc V_{\textnormal{bg}}^{\mathrm{nl}}
            \in\cL\pare{\bR^2;\pare{\cF_{\geq1}^\bR}^2}.
        \end{equation}
        In particular,
        $\mathpzc V_{\textnormal{bg}}
        \bra{\background}\pare{\Atw\middle|\eta;x}$ is independent of $\xi$ and equals
        $\frac{\Atw}{4}\background^\perp$ at $\eta=0$.
    \end{itemize}

    \item $A_{\bra{0,\geq1}}$ is a generic vector symbol as per
    \cref{notation:paralinearization};

    \item $R^\bR_{\bra{\varrho,\geq1}}$ is a smoothing remainder as per
    \cref{notation:paralinearization}.
\end{itemize}
\end{lemma}

\begin{proof}
For the duration of the proof set
\[
    \Omega_{\textnormal{bg}}
    \defeq
    \tilde{\Omega}\bra{\background}\pare{\eta},
    \qquad
    \delta\Omega_{\textnormal{bg}}
    \defeq
    \Omega_{\textnormal{bg}}-\background,
    \qquad
    \Omega_{\tilde{\psi}}
    \defeq
    S\pare{\Atw\middle|\eta}\tilde{\psi}.
\]
Then \cref{eq:pitilde1_recall} reads
\[
    \tilde{\pi}_1
    \pare{\Atw,\background\middle|\eta}
    \bra{\tilde{\psi}}
    =
    2\,
    \pi^\star_2\pare{\Atw\middle|\eta}
    \bra{\Omega_{\textnormal{bg}},\Omega_{\tilde{\psi}}}.
\]

\begin{proofpart}[Application of the bilinear paralinearization]
For the target gain $\varrho$, use
\cref{cor:paralinearization_pistar2_bilinear,prop:paralinearization_S} at
gain $\varrho+1$ and choose the Sobolev threshold one derivative above the
one in the corollary. This is legitimate because
$\Omega_{\tilde\psi}=S\pare{\Atw\middle|\eta}\tilde\psi$ loses one
derivative, and the subsequent order-one composition turns gain
$\varrho+1$ into gain $\varrho$. We also split
\begin{equation}\label{eq:pitilde1_background_bilinear_split}
\pi_2^\star\pare{\Atw\middle|\eta}
\bra{\Omega_{\textnormal{bg}},\Omega_{\tilde\psi}}
=
\pi_2^\star\pare{\Atw\middle|\eta}
\bra{\background,\Omega_{\tilde\psi}}
+
\pi_2^\star\pare{\Atw\middle|\eta}
\bra{\delta\Omega_{\textnormal{bg}},\Omega_{\tilde\psi}}.
\end{equation}
Using the bilinear formula with $\background$ as an external parameter,
applying its Sobolev estimates to
$\delta\Omega_{\textnormal{bg}}$ and $\Omega_{\tilde\psi}$, and using
\cref{rem:action_paradiff_constants} for the direct actions on
$\background$, then resumming, gives
\begin{equation}\label{eq:pitilde1_bilinear_raw}
\begin{aligned}
\tilde{\pi}_1
\pare{\Atw,\background\middle|\eta}
\bra{\tilde{\psi}}
&=
2\OpBW{
    a^{\pi^\star_2,\eta}_1
    \bra{\Omega_{\textnormal{bg}},\Omega_{\tilde{\psi}}}
    \pare{\Atw\middle|\eta;x,\xi}
}\eta
\\
&\quad
+
\OpBW{
    A^{\pi^\star_2,\Omega}_0
    \bra{\Omega_{\tilde{\psi}}}
    \pare{\Atw\middle|\eta;x,\xi}
}\cdot\Omega_{\textnormal{bg}}
\\
&\quad
+
\OpBW{
    A^{\pi^\star_2,\Omega}_0
    \bra{\Omega_{\textnormal{bg}}}
    \pare{\Atw\middle|\eta;x,\xi}
}\cdot\Omega_{\tilde{\psi}}
\\
&\quad
+
\OpBW{
    A_{\bra{0,\geq1}}
    \pare{\Atw,\background\middle|\eta,\tilde{\psi};x,\xi}
}\cdot\vect{\eta}{\tilde{\psi}}
+
R^\bR_{\bra{\varrho,\geq1}}
\pare{\Atw,\background\middle|\eta,\tilde{\psi}}
\cdot\vect{\eta}{\tilde{\psi}} .
\end{aligned}
\end{equation}
For the background factor we use directly
\cref{eq:paralinearization_S_background,eq:ASeta1_background}; in
particular, the principal $\eta$-symbol of
$\delta\Omega_{\textnormal{bg}}$ is
$A^{S,\eta}_{1,\textnormal{bg}}\bra{\background}$.
Likewise, by \cref{prop:paralinearization_S},
\begin{equation}\label{eq:pitilde1_Omegapsi_paralinearization}
    \Omega_{\tilde{\psi}}
    =
    \OpBW{
        \vectT{
            A^{S,\eta}_1\bra{\tilde{\psi}}
            \pare{\Atw\middle|\eta;x,\xi}
        }{
            \ii\xi^\perp
        }
        +
        A_{\bra{0,\geq1}}
        \pare{\Atw\middle|\eta,\tilde{\psi};x,\xi}
    }\vect{\eta}{\tilde{\psi}}
    +
    R^\bR_{\bra{\varrho,\geq1}}
    \pare{\Atw\middle|\eta,\tilde{\psi}}
    \cdot\vect{\eta}{\tilde{\psi}},
\end{equation}
where
\begin{equation}\label{eq:pitilde1_ASeta}
    A^{S,\eta}_1\bra{\tilde{\psi}}
    \pare{\Atw\middle|\eta;x,\xi}
    =
    -\Atw
    \left(
        A_{0}^{\mathpzc{B}_0,\omega;\perp}\pare{\eta;x,\xi}\cdot\Omega_{\tilde{\psi}}\pare{x}
        -
        \ii\,
        v_3\bra{\tilde{\psi}}\pare{\Atw\middle|\eta;x}
    \right)
    \xi^\perp .
\end{equation}
The constant vector $\background$ in
\cref{eq:paralinearization_S_background} is annihilated by Bony--Weyl
quantization. Therefore, composing
\cref{eq:pitilde1_bilinear_raw} with
\cref{eq:paralinearization_S_background} and
\cref{eq:pitilde1_Omegapsi_paralinearization}, and using the symbolic
composition formula, gives
\begin{equation}\label{eq:pitilde1_raw_vector_symbol}
\begin{aligned}
\tilde{\pi}_1
\pare{\Atw,\background\middle|\eta}
\bra{\tilde{\psi}}
&=
\OpBW{
    \vect{
        2a^{\pi^\star_2,\eta}_1
        \bra{\Omega_{\textnormal{bg}},\Omega_{\tilde{\psi}}}
        +
        A^{\pi^\star_2,\Omega}_0
        \bra{\Omega_{\tilde{\psi}}}
        \cdot
        A^{S,\eta}_{1,\textnormal{bg}}
        \bra{\background}
        +
        A^{\pi^\star_2,\Omega}_0
        \bra{\Omega_{\textnormal{bg}}}
        \cdot
        A^{S,\eta}_1\bra{\tilde{\psi}}
    }{
        A^{\pi^\star_2,\Omega}_0
        \bra{\Omega_{\textnormal{bg}}}
        \cdot \ii\xi^\perp
    }
}
\cdot\vect{\eta}{\tilde{\psi}}
\\
&\quad
+
\OpBW{
    A_{\bra{0,\geq1}}
    \pare{\Atw,\background\middle|\eta,\tilde{\psi};x,\xi}
}\cdot\vect{\eta}{\tilde{\psi}}
+
R^\bR_{\bra{\varrho,\geq1}}
\pare{\Atw,\background\middle|\eta,\tilde{\psi}}
\cdot\vect{\eta}{\tilde{\psi}} .
\end{aligned}
\end{equation}
All symbols in \cref{eq:pitilde1_raw_vector_symbol} are evaluated at
$\pare{\Atw\middle|\eta;x,\xi}$.
\end{proofpart}

\begin{proofpart}[$\tilde{\psi}$-component]
Using \cref{eq:pistar2_Omega_symbol_bilinear} with
$\Omega_j=\Omega_{\textnormal{bg}}$, we compute
\begin{equation}\label{eq:pitilde1_psicomponent_raw}
\begin{aligned}
A^{\pi^\star_2,\Omega}_0
\bra{\Omega_{\textnormal{bg}}}
\pare{\Atw\middle|\eta;x,\xi}
\cdot \ii\xi^\perp
&=
\ii\,
\mathpzc V_{\textnormal{bg}}
\bra{\background}
\pare{\Atw\middle|\eta;x}
\cdot\xi
+
\pare{1+\av{\nabla\eta}^2}
b_{\textnormal{bg}}\bra{\background}
\pare{\Atw\middle|\eta;x}
A_{0}^{\mathpzc{B}_0,\omega;\perp}\pare{\eta;x,\xi}\cdot\xi^\perp .
\end{aligned}
\end{equation}
Since
\begin{equation}\label{eq:pitilde1_H0_identity}
    \pare{1+\av{\nabla\eta}^2}
    A_{0}^{\mathpzc{B}_0,\omega;\perp}\pare{\eta;x,\xi}\cdot\xi^\perp
    =
    -\lambda_1\pare{\eta;x,\xi},
\end{equation}
we obtain
\begin{equation}\label{eq:pitilde1_psicomponent_final}
A^{\pi^\star_2,\Omega}_0
\bra{\Omega_{\textnormal{bg}}}
\pare{\Atw\middle|\eta;x,\xi}
\cdot \ii\xi^\perp
=
-\lambda_1\pare{\eta;x,\xi}
b_{\textnormal{bg}}\bra{\background}
\pare{\Atw\middle|\eta;x}
+
\ii\,
\mathpzc V_{\textnormal{bg}}
\bra{\background}
\pare{\Atw\middle|\eta;x}
\cdot\xi .
\end{equation}
This is exactly \cref{eq:apitilde1_psitilde1}.
\end{proofpart}

\begin{proofpart}[$\eta$-component]
We simplify the first component of
\cref{eq:pitilde1_raw_vector_symbol}. The two elementary identities
\begin{equation}\label{eq:pitilde1_identity_1_bg}
    \frac{
        \bm J
        \pare{
            \Omega_{\textnormal{bg}}
            +
            \pare{\nabla\eta\cdot\Omega_{\textnormal{bg}}}\nabla\eta
        }\cdot\xi
    }{
        4\pare{1+\av{\nabla\eta}^2}
    }
    =
    \lambda_1\pare{\eta;x,\xi}
    \pare{
        A_{0}^{\mathpzc{B}_0,\omega;\perp}\pare{\eta;x,\xi}\cdot\Omega_{\textnormal{bg}}
    },
\end{equation}
and
\begin{equation}\label{eq:pitilde1_identity_1_psi}
    \frac{
        \bm J
        \pare{
            \Omega_{\tilde{\psi}}
            +
            \pare{\nabla\eta\cdot\Omega_{\tilde{\psi}}}\nabla\eta
        }\cdot\xi
    }{
        4\pare{1+\av{\nabla\eta}^2}
    }
    =
    \lambda_1\pare{\eta;x,\xi}
    \pare{
        A_{0}^{\mathpzc{B}_0,\omega;\perp}\pare{\eta;x,\xi}\cdot\Omega_{\tilde{\psi}}
    }
\end{equation}
follow from the explicit formula for $A_{0}^{\mathpzc{B}_0,\omega}$ and from \cref{eq:lambda1}. We shall also use
\begin{equation}\label{eq:pitilde1_identity_2_cross}
\begin{aligned}
&\frac{1}{4\pare{1+\av{\nabla\eta}^2}}
\left(
    \frac{1}{2}
    \pare{\Omega_{\textnormal{bg}}\cdot\nabla\eta}\Omega_{\tilde{\psi}}
    +
    \frac{1}{2}
    \pare{\Omega_{\tilde{\psi}}\cdot\nabla\eta}\Omega_{\textnormal{bg}}
    -
    \frac{
        \Omega_{\textnormal{bg}}\cdot\Omega_{\tilde{\psi}}
        +
        \pare{\Omega_{\textnormal{bg}}\cdot\nabla\eta}
        \pare{\Omega_{\tilde{\psi}}\cdot\nabla\eta}
    }{
        1+\av{\nabla\eta}^2
    }
    \nabla\eta
\right)\cdot\xi
\\
&\qquad =
\frac{\lambda_1\pare{\eta;x,\xi}}{2}
\left[
    \pare{
        A_{0}^{\mathpzc{B}_0,\omega;\perp}\pare{\eta;x,\xi}\cdot\Omega_{\textnormal{bg}}
    }
    \frac{
        \nabla^\perp\eta\cdot\Omega_{\tilde{\psi}}
    }{
        1+\av{\nabla\eta}^2
    }
    +
    \pare{
        A_{0}^{\mathpzc{B}_0,\omega;\perp}\pare{\eta;x,\xi}\cdot\Omega_{\tilde{\psi}}
    }
    \frac{
        \nabla^\perp\eta\cdot\Omega_{\textnormal{bg}}
    }{
        1+\av{\nabla\eta}^2
    }
\right].
\end{aligned}
\end{equation}
This is the polarization of \cref{eq:pitilde0_identity_2}.

Combining \cref{eq:pistar2_eta_symbol_bilinear},
\cref{eq:ASeta1_background}, and \cref{eq:pitilde1_ASeta}, the
imaginary part of the first component in
\cref{eq:pitilde1_raw_vector_symbol} is a sum of the two brackets
\[
\lambda_1
\pare{
    A_{0}^{\mathpzc{B}_0,\omega;\perp}\cdot\Omega_{\textnormal{bg}}
}
\frac{\nabla^\perp\eta\cdot\Omega_{\tilde{\psi}}}
{1+\av{\nabla\eta}^2}
+
\lambda_1
\pare{
    A_{0}^{\mathpzc{B}_0,\omega;\perp}\cdot\Omega_{\tilde{\psi}}
}
\frac{\nabla^\perp\eta\cdot\Omega_{\textnormal{bg}}}
{1+\av{\nabla\eta}^2}
\]
minus the corresponding expression in
\cref{eq:pitilde1_identity_2_cross}, together with the two transport
brackets cancelled by
\cref{eq:pitilde1_identity_1_bg} and
\cref{eq:pitilde1_identity_1_psi}. Hence the imaginary part vanishes.

The real part is
\begin{equation}\label{eq:pitilde1_eta_component_final}
\begin{aligned}
&2a^{\pi^\star_2,\eta}_1
\bra{\Omega_{\textnormal{bg}},\Omega_{\tilde{\psi}}}
+
A^{\pi^\star_2,\Omega}_0
\bra{\Omega_{\tilde{\psi}}}
\cdot
A^{S,\eta}_{1,\textnormal{bg}}
\bra{\background}
+
A^{\pi^\star_2,\Omega}_0
\bra{\Omega_{\textnormal{bg}}}
\cdot
A^{S,\eta}_1\bra{\tilde{\psi}}
\\
=&
2\lambda_1\pare{\eta;x,\xi}
\Bigg[
    \pare{4-\Atw^2}
    \pare{
        A_{0}^{\mathpzc{B}_0,\omega;\perp}\pare{\eta;x,\xi}\cdot\Omega_{\textnormal{bg}}
    }
    \pare{
        A_{0}^{\mathpzc{B}_0,\omega;\perp}\pare{\eta;x,\xi}\cdot\Omega_{\tilde{\psi}}
    }
\\
-&
    \left(
        \Atw v_3^{\star}\bra{\Omega_{\textnormal{bg}}}\pare{\eta;x}
        -
        \frac{
            \nabla^\perp\eta\cdot\Omega_{\textnormal{bg}}
        }{
            1+\av{\nabla\eta}^2
        }
    \right)
    \left(
        \Atw v_3\bra{\tilde{\psi}}\pare{\Atw\middle|\eta;x}
        -
        \frac{
            \nabla^\perp\eta\cdot\Omega_{\tilde{\psi}}
        }{
            1+\av{\nabla\eta}^2
        }
    \right)
\Bigg].
\end{aligned}
\end{equation}
By the definitions
\cref{eq:bpitilde0_bg} and \cref{eq:bpitilde1_psitilde}, this is exactly
\cref{eq:apitilde1_eta1}.
\end{proofpart}

\begin{proofpart}[Synthesis]
The identities
\cref{eq:pitilde1_psicomponent_final} and
\cref{eq:pitilde1_eta_component_final}, inserted into
\cref{eq:pitilde1_raw_vector_symbol}, give
\cref{eq:paralinearization_pitilde1}. The lower order symbols and smoothing
remainders collect the order-zero part of the paralinearizations of
$S\pare{\Atw\middle|\eta}$ and of
$\tilde{\Omega}\bra{\background}\pare{\eta}$, the order $-1$ column of
\cref{cor:paralinearization_pistar2_bilinear}, and the symbolic remainders in
the Bony--Weyl composition formula.
\end{proofpart}
\end{proof}

\paragraph[Paralinearization of pi-tilde 2]{Paralinearization of $\tilde{\pi}_2$}
Recall that $\tilde{\pi}_2$ is defined in \cref{eq:pi2tilde}, we recall the
definition here for
convenience
\begin{equation}
     \tilde{\pi}_{2}\pare{\Atw\middle|\eta}\bra{\tilde{\psi},\tilde{\psi}}
    =
    \pi_2\pare{\Atw\middle|\eta}\bra{\tilde{\psi},\tilde{\psi}}
\end{equation}
so the paralinearization of $\tilde{\pi}_2$ is provided by \cref{eq:pi2_raw_vector_symbol}.

\begin{proof} [Proof of \cref{prop:paralinearization_pitilde}]

We combine the three homogeneous components in \cref{eq:tildepi_decomposition},
using the shorthand
\begin{equation}
    \Omega_{\textnormal{bg}}
    \defeq
    \tilde{\Omega}\bra{\background}\pare{\eta} ,
    \qquad
    \Omega_{\tilde{\psi}}
    \defeq
    S\pare{\Atw\middle|\eta}\tilde{\psi} ,
\end{equation}
so that, by \cref{eq:Omegacan_pitilde},
$\Omega\pare{\background\middle|\eta,\tilde{\psi}}
= \Omega_{\textnormal{bg}} + \Omega_{\tilde{\psi}}$.
Since the formulas \cref{eq:bpitilde_can} and \cref{eq:Vpitilde_can}
are linear in $\Omega$, this decomposition gives
\begin{equation}
    b_{\textnormal{can}}\pare{\Atw,\background\middle|\eta,\tilde{\psi}}
    =
    b_{\textnormal{bg}}\bra{\background}\pare{\Atw\middle|\eta}
    +
    b_{\tilde{\psi}}\bra{\tilde{\psi}}\pare{\Atw\middle|\eta} ,
    \qquad
    \mathpzc V_{\textnormal{can}}\pare{\Atw,\background\middle|\eta,\tilde{\psi}}
    =
    \mathpzc V_{\textnormal{bg}}\bra{\background}\pare{\Atw\middle|\eta}
    +
    \mathpzc V\bra{\tilde{\psi}}\pare{\Atw\middle|\eta} .
\end{equation}

\begin{proofpart}[Order-two symbol]
The only order-two contribution comes from $\tilde\pi_0$,
which is paralinearized in \cref{eq:paralinearization_pican0};
the order-two symbol is therefore exactly \cref{eq:Apitilde_2}.
\end{proofpart}
\begin{proofpart}[Order-one $\eta$-component]
By \cref{lem:paralinearization_pitilde0,lem:paralinearization_pitilde1} and
\cref{prop:paralinearization_pi}, the three contributions to the
order-one $\eta$-symbol are
$a^{\tilde{\pi}_0,\eta}_1$,
$a^{\tilde{\pi}_1,\eta}_1$, and
$a^{\pi,\eta}_1$ (with $\psi$ replaced by $\tilde{\psi}$).
Summing and using $\Omega = \Omega_{\textnormal{bg}}+\Omega_{\tilde{\psi}}$,
\[
\begin{aligned}
    &
    a^{\tilde{\pi}_0,\eta}_1
    +
    a^{\tilde{\pi}_1,\eta}_1
    +
    a^{\pi,\eta}_1
    =
    \lambda_1\pare{\eta;x,\xi}
    \bra{
        \pare{4-\Atw^2}
        \pare{
            A_{0}^{\mathpzc{B}_0,\omega;\perp}\pare{\eta;x,\xi}
            \cdot
            \Omega\pare{\background\middle|\eta,\tilde{\psi};x}
        }^2
        -
        \pare{
            b_{\textnormal{can}}
            \pare{\Atw,\background\middle|\eta,\tilde{\psi};x}
        }^2
    } ,
\end{aligned}
\]
which is precisely
$\lambda_1\bra{q^{\tilde{\pi}}_0 - b_{\textnormal{can}}^2}$,
i.e., \cref{eq:apitilde_eta1} by \cref{eq:qpitilde0_definition}.
\end{proofpart}
\begin{proofpart}[Order-one $\tilde{\psi}$-component]
The two contributions are $a^{\tilde{\pi}_1,\tilde{\psi}}_1$ from
\cref{eq:apitilde1_psitilde1} and $a^{\pi,\psi}_1$ from \cref{eq:Api_psi1}
with $\psi=\tilde{\psi}$.
Summing and using the additivity of $b_{\textnormal{can}}$
and $\mathpzc V_{\textnormal{can}}$,
\[
    a^{\tilde{\pi}_1,\tilde{\psi}}_1
    +
    a^{\pi,\psi}_1
    =
    -\lambda_1\pare{\eta;x,\xi}\,
    b_{\textnormal{can}}
    \pare{\Atw,\background\middle|\eta,\tilde{\psi};x}
    +
    \ii\,
    \mathpzc V_{\textnormal{can}}
    \pare{\Atw,\background\middle|\eta,\tilde{\psi};x}
    \cdot\xi ,
\]
which is exactly \cref{eq:apitilde_psitilde1}.
\end{proofpart}
\begin{proofpart}[Lower-order terms]
All remaining contributions have order at most zero:
they come from the generic order-zero symbols in
\cref{eq:paralinearization_pican0} and \cref{eq:paralinearization_pitilde1},
from the lower-order parts of the paralinearizations of
$\tilde{\Omega}\bra{\background}\pare{\eta}$ and
$S\pare{\Atw\middle|\eta}\tilde{\psi}$,
and from the Bony--Weyl composition remainders.
They are absorbed in $A_{\bra{0,\geq1}}$ and $R^\bR_{\bra{\varrho,\geq1}}$,
giving \cref{eq:paralinearization_pitilde} and concluding the proof of
\cref{prop:paralinearization_pitilde}.
\end{proofpart}
\end{proof}

\subsection{Synthesis: proof of \cref{thm:paraCDEetapsi}}
\normalfont

We now assemble the paralinearizations obtained in the previous subsections in
the canonical variables of \cref{proprecast1}. The point is that the equation
for $\eta$ contains both the canonical background contribution
$\mathpzc g^\star\pare{\eta}
\bra{\tilde{\Omega}\bra{\background}\pare{\eta}}$
and the Dirichlet--Neumann contribution
$\mathpzc g\pare{\Atw\middle|\eta}\tilde{\psi}$, whereas the equation for
$\tilde{\psi}$ is the shifted Bernoulli equation paralinearized in
\cref{prop:paralinearization_pitilde}. We keep the expansion of the
Dirichlet--Neumann column in $\tilde{\psi}$ up to order zero, namely
$\lambda_1+\Atw\lambda_0$, and absorb only the terms of order at most $-1$
in the corresponding generic symbol.

\begin{proof}
Set
\begin{equation}\label{eq:Omega_bg_and_tildepsi_synthesis}
    \Omega_{\textnormal{bg}}
    \defeq
    \tilde{\Omega}\bra{\background}\pare{\eta},
    \qquad
    \Omega_{\tilde{\psi}}
    \defeq
    S\pare{\Atw\middle|\eta}\tilde{\psi}.
\end{equation}
Thus, by \cref{eq:Omegacan_paraCDE},
\begin{equation}\label{eq:Omegacan_synthesis}
    \Omega
    \pare{\background\middle|\eta,\tilde{\psi}}
    =
    \Omega_{\textnormal{bg}}
    +
    \Omega_{\tilde{\psi}}.
\end{equation}

\begin{proofpart}[Paralinearization of the $\eta$-equation]
By \cref{eq:CDEeta}, the first equation is
\begin{equation}\label{eq:eta_equation_synthesis_start}
    \eta_t
    =
    \mathpzc g^\star\pare{\eta}
    \bra{\Omega_{\textnormal{bg}}}
    +
    \mathpzc g\pare{\Atw\middle|\eta}
    \bra{\tilde{\psi}}.
\end{equation}
We first paralinearize the background term. Write
\[
    \Omega_{\textnormal{bg}}
    =
    \background+\delta\Omega_{\textnormal{bg}},
    \qquad
    \delta\Omega_{\textnormal{bg}}
    \defeq
    \Omega_{\textnormal{bg}}-\background.
\]
The canonical paralinearization and its principal symbol are exactly
\cref{eq:paralinearization_S_background,eq:ASeta1_background}; we do not
redefine them here. Since $\mathpzc g^\star\pare\eta$ is linear in its
vorticity argument, split
\[
\mathpzc g^\star\pare\eta\bra{\Omega_{\textnormal{bg}}}
=
\mathpzc g^\star\pare\eta\bra{\background}
+
\mathpzc g^\star\pare\eta\bra{\delta\Omega_{\textnormal{bg}}}.
\]
For the first term we use \cref{rem:action_paradiff_constants}: every direct
Bony--Weyl action on $\background$ vanishes, while the dependence on
$\background$ remains in the symbol acting on $\eta$. We apply
\cref{prop:Gstar_paralinearization} to the Sobolev correction in the second
term and compose it with
\cref{eq:paralinearization_S_background,eq:ASeta1_background}. Resumming the
two pieces gives
\begin{equation}\label{eq:gstar_background_raw_synthesis}
\begin{aligned}
\mathpzc g^\star\pare{\eta}\bra{\Omega_{\textnormal{bg}}}
&=
\OpBW{
    a^{\mathpzc g^\star,\eta}_1
    \bra{\Omega_{\textnormal{bg}}}
    \pare{\eta;x,\xi}
    +
    \ii\,
    A^{\mathpzc g^\star,\Omega}_{0}
    \pare{\eta;x,\xi}
    \cdot
    A^{S,\eta}_{1,\textnormal{bg}}
    \bra{\background}
    \pare{\Atw\middle|\eta;x,\xi}
}
\eta
\\
&\quad
+
\OpBW{
    a_{\bra{0,\geq1}}\bra{\background}
    \pare{\Atw\middle|\eta;x,\xi}
}
\eta
+
r^{\lin,\bR}_{\bra{\varrho,\geq1}}\bra{\background}
\pare{\Atw\middle|\eta}\eta .
\end{aligned}
\end{equation}
Using \cref{eq:AGstarW0}, \cref{eq:aGstar_1_corrected}, and
\[
    \pare{1+\av{\nabla\eta}^2}
    A_{0}^{\mathpzc{B}_0,\omega;\perp}\pare{\eta;x,\xi}\cdot\xi^\perp
    =
    -\lambda_1\pare{\eta;x,\xi},
\]
the principal symbol in \cref{eq:gstar_background_raw_synthesis} becomes
\begin{equation}\label{eq:gstar_background_symbol_synthesis}
\begin{aligned}
&a^{\mathpzc g^\star,\eta}_1
    \bra{\Omega_{\textnormal{bg}}}
    \pare{\eta;x,\xi}
    +
    \ii\,
    A^{\mathpzc g^\star,\Omega}_{0}
    \pare{\eta;x,\xi}
    \cdot
    A^{S,\eta}_{1,\textnormal{bg}}
    \bra{\background}
    \pare{\Atw\middle|\eta;x,\xi}\\
=&
\lambda_1\pare{\eta;x,\xi}
\left(
    \Atw\,
     v_3^{\star}
    \bra{\Omega_{\textnormal{bg}}}
    \pare{\eta;x}
    -
    \frac{
        \nabla^\perp\eta\cdot\Omega_{\textnormal{bg}}
    }{
        1+\av{\nabla\eta}^2
    }
\right)
\\
+&
\ii
\left(
     V^{\star}
    \bra{\Omega_{\textnormal{bg}}}
    \pare{\eta;x}
    +
    \frac{\Atw}{4\pare{1+\av{\nabla\eta}^2}}
    \bm J
    \pare{
        \Omega_{\textnormal{bg}}
        +
        \pare{\nabla\eta\cdot\Omega_{\textnormal{bg}}}\nabla\eta
    }
\right)
\cdot\xi .
\end{aligned}
\end{equation}
Equivalently, using \cref{eq:bpitilde0_bg} and \cref{eq:Vpitilde1_bg}, we have
\begin{equation}\label{eq:gstar_background_final_synthesis}
    \mathpzc g^\star\pare{\eta}\bra{\Omega_{\textnormal{bg}}}
    =
    \OpBW{
        \lambda_1
        b_{\textnormal{bg}}
        +
        \ii\,
        \mathpzc V_{\textnormal{bg}}\bra{\background}\pare{\Atw\middle|\eta;x}\cdot\xi
        +
        A_{\bra{0,\geq1}}^\eta
    }\eta
    +
    r^{\lin,\bR}_{\bra{\varrho,\geq1}}\eta .
\end{equation}
By \cref{eq:Vpitilde1_bg_split}, the transport coefficient in this formula
is the sum of the flat velocity
$\frac{\Atw}{4}\background^\perp$ and a nonlinear correction in
$\pare{\cF_{\geq1}^\bR}^2$.

On the other hand, \cref{prop:G_paralinearization} applied to
$\tilde{\psi}$ gives, with the order-zero DN contribution kept explicit,
\begin{equation}\label{eq:g_tildepsi_synthesis}
\begin{aligned}
\mathpzc g\pare{\Atw\middle|\eta}\tilde{\psi}
&=
\OpBW{
    \vect{
        \lambda_1
        b_{\tilde{\psi}}
        +
        \ii\,
        \mathpzc V\bra{\tilde{\psi}}\pare{\Atw\middle|\eta;x}\cdot\xi
    }{
        \lambda_1+\Atw\lambda_0
    }
    +
    \vect{
        A_{\bra{0,\geq1}}^\eta
        \pare{\Atw\middle|\eta,\tilde{\psi};x,\xi}
    }{
        A_{\bra{-1,\geq1}}^\psi
        \pare{\Atw\middle|\eta,\tilde{\psi};x,\xi}
    }
}
\cdot
\vect{\eta}{\tilde{\psi}}\\&
+
R^{\lin,\bR}_{\bra{\varrho,\geq1}}
\pare{\Atw\middle|\eta,\tilde{\psi}}
\cdot
\vect{\eta}{\tilde{\psi}}.
\end{aligned}
\end{equation}
Here the symbol $\lambda_0$ is the real order-zero symbol defined in
\cref{eq:lambda0}; the generic symbol in the $\tilde{\psi}$-column is now of
order $-1$ because the order-zero part of that column is precisely the term
$\Atw\lambda_0$ displayed in \cref{eq:g_tildepsi_synthesis}.

Adding \cref{eq:gstar_background_final_synthesis} and
\cref{eq:g_tildepsi_synthesis}, and using the linearity of
${\bm v}^\star$ in the vorticity argument, gives
\begin{equation}\label{eq:eta_equation_final_symbol_synthesis}
\begin{aligned}
\eta_t
&=
\OpBW{
    \vect{
        \lambda_1
        b_{\textnormal{can}}
        \pare{\Atw,\background\middle|\eta,\tilde{\psi};x}
        +
        \ii\,
        \mathpzc V_{\textnormal{can}}
        \pare{\Atw,\background\middle|\eta,\tilde{\psi};x}
        \cdot\xi
    }{
        \lambda_1\pare{\eta;x,\xi}
        +
        \Atw\lambda_0\pare{\eta;x,\xi}
    }
    +
    \vect{
        A_{\bra{0,\geq1}}^\eta
        \pare{\Atw,\background\middle|\eta,\tilde{\psi};x,\xi}
    }{
        A_{\bra{-1,\geq1}}^\psi
        \pare{\Atw,\background\middle|\eta,\tilde{\psi};x,\xi}
    }
}
\cdot
\vect{\eta}{\tilde{\psi}}
\\&+
R^{\lin,\bR}_{\bra{\varrho,\geq1}}
\pare{\Atw,\background\middle|\eta,\tilde{\psi}}
\cdot
\vect{\eta}{\tilde{\psi}} .
\end{aligned}
\end{equation}
\end{proofpart}

\begin{proofpart}[Paralinearization of the $\tilde{\psi}$-equation]
By \cref{eq:evolution_psi}, the second equation is
\[
    \tilde{\psi}_t
    =
    \tilde{\pi}
    \pare{\rho^\pm,\background\middle|\eta,\tilde{\psi}}.
\]
Therefore \cref{prop:paralinearization_pitilde} gives
\begin{equation}\label{eq:psitilde_equation_final_symbol_synthesis}
\begin{aligned}
\tilde{\psi}_t
&=
-\mathfrak g\Atw\,\eta
+
\OpBW{
    \vect{
        a^{\pi_0,\eta}_2
        \pare{\rho^\pm\middle|\eta;x,\xi}
        +
        \lambda_1\pare{\eta;x,\xi}
        \left(
            q^{\tilde{\pi}}_0
            \pare{\Atw,\background\middle|\eta,\tilde{\psi};x,\xi}
            -
            b_{\textnormal{can}}^2
            \pare{\Atw,\background\middle|\eta,\tilde{\psi};x}
        \right)
    }{
        -\lambda_1\pare{\eta;x,\xi}
        b_{\textnormal{can}}
        \pare{\Atw,\background\middle|\eta,\tilde{\psi};x}
        +
        \ii\,
        \mathpzc V_{\textnormal{can}}
        \pare{\Atw,\background\middle|\eta,\tilde{\psi};x}
        \cdot\xi
    }
}
\cdot
\vect{\eta}{\tilde{\psi}}
\\
&\quad
+
\OpBW{
    \vect{
        A_{\bra{0,\geq1}}^\eta
        \pare{\Atw,\background\middle|\eta,\tilde{\psi};x,\xi}
    }{
        A_{\bra{-1,\geq1}}^\psi
        \pare{\Atw,\background\middle|\eta,\tilde{\psi};x,\xi}
    }
}
\cdot
\vect{\eta}{\tilde{\psi}}
\\
&\quad
+
R^{\lin,\bR}_{\bra{\varrho,\geq1}}
\pare{\Atw,\background\middle|\eta,\tilde{\psi}}
\cdot
\vect{\eta}{\tilde{\psi}}.
\end{aligned}
\end{equation}
All symbols in \cref{eq:psitilde_equation_final_symbol_synthesis} are
evaluated at the displayed arguments whenever the dependence is not written
explicitly.
\end{proofpart}

\begin{proofpart}[Matrix synthesis]
Combining \cref{eq:eta_equation_final_symbol_synthesis} and
\cref{eq:psitilde_equation_final_symbol_synthesis}, the terms that are kept
explicit in the quasilinear matrix are the following.

First, the symbols independent of $b_{\textnormal{can}}$ and of the transport
assemble into
\[
    {\bm Q}
    =
    \begin{bmatrix}
        0
        &
        \lambda_1+\Atw\lambda_0
        \\
        -\mathfrak g\Atw
        +
        a^{\pi_0,\eta}_2
        +
        \lambda_1 q^{\tilde{\pi}}_0
        &
        0
    \end{bmatrix}.
\]
Second, the real order-one terms involving $b_{\textnormal{can}}$ assemble
into the Alinhac block
\[
    {\bm B}
    =
    -
    \begin{bmatrix}
        -b_{\textnormal{can}}
        &
        0
        \\
        b_{\textnormal{can}}^2
        &
        b_{\textnormal{can}}
    \end{bmatrix}
    \lambda_1 .
\]
Indeed, the first row of ${\bm B}$ gives
$\lambda_1 b_{\textnormal{can}}$ in the $\eta$-column, while the second row
gives
$-\lambda_1 b_{\textnormal{can}}^2$ in the $\eta$-column and
$-\lambda_1 b_{\textnormal{can}}$ in the $\tilde{\psi}$-column.

Finally, the imaginary first-order transport symbols in the two equations are
the same, namely
\[
    \ii\,
    \mathpzc V_{\textnormal{can}}
    \pare{\Atw,\background\middle|\eta,\tilde{\psi};x}
    \cdot\xi,
\]
and therefore form the diagonal transport matrix
\[
    \ii
    \pare{
        \mathpzc V_{\textnormal{can}}
        \pare{\Atw,\background\middle|\eta,\tilde{\psi};x}
        \cdot\xi
    }
    \Id_{\bC^2}.
\]
The remaining terms in the $\eta$-column are collected in
$A_{\bra{0,\geq1}}^\eta$, whereas the remaining terms in the $\tilde{\psi}$-column
are collected in $A_{\bra{-1,\geq1}}^\psi$, consistently with the fact that the
order-zero DN contribution $\Atw\lambda_0$ has been displayed explicitly in
${\bm Q}$. The smoothing remainders are collected in
${\bm R}^{\lin,\bR}_{\bra{\varrho,\geq1}}$. This gives exactly
\cref{eq:paralinearazationetapsi}.
\end{proofpart}
\end{proof}

\section{Diagonal reduction and proof of \cref{thm:main}}
\label{sec:diagonal_reduction_main}

The purpose of this section is to prove \cref{thm:main}.
The argument proceeds by a cascade of paradifferential conjugations that
progressively simplify the canonical CDE \cref{eq:CDE}: the real
paradifferential system \cref{eq:paralinearazationetapsi} is brought to a
fully diagonal complex form with purely imaginary positive-order symbols, and
a modified energy for the resulting scalar system satisfies a closed
Riccati-type differential inequality.
A prescribed-background regularization followed by a Kato--Bona--Smith
argument then yields existence, uniqueness, and continuous dependence.
The section is organized as follows.

\begin{itemize}

    \item \cref{sec:real_alinhac_good_unknown_3d} introduces the real Alinhac
    good unknown $Y_{\rm Al} \defeq \vect{\eta_{\rm Al}}{\tilde\psi_{\rm Al}}$,
    defined by the nilpotent triangular conjugation ${\bm G}$ of
    \cref{def:real_alinhac_good_unknown_3d}.
    The conjugation is exact by \cref{eq:real_alinhac_nilpotency_3d}, leaves
    the elevation component unchanged, and removes the off-diagonal Alinhac
    block ${\bm B}$ from \cref{eq:paralinearazationetapsi}; the reduced system
    is stated in \cref{prop:real_alinhac_good_unknown_3d}.

    \item \cref{sec:complex_formulation_3d_vs} reformulates the
    Alinhac-reduced system in complex scalar variables.
    Under the strict Atwood--Weber stability condition
    \cref{hyp:strict_atwood_weber_stability}, the principal symbol of the
    resulting system is block-diagonal with purely imaginary entries of order
    $3/2$; the complex paradifferential system is recorded in
    \cref{prop:complex_paradifferential_system_3d}.

    \item \cref{sec:positive_order_diagonalization_3d} diagonalizes all
    positive-order symbols of the complex system via two successive
    paradifferential conjugations.
    The quasilinear order-$3/2$ block is treated in
    \cref{sec:principal_order_diagonalization_3d}; the explicit order-$1/2$
    off-diagonal block is treated in
    \cref{sec:subprincipal_order_diagonalization_3d}.
    The outcome is the fully diagonal positive-order system of
    \cref{prop:positive_order_diagonalization_3d}, with diagonalized unknown
    $W$.

    \item \cref{sec:modified_energy_estimates_3d} constructs a modified energy
    $E_s$ for the diagonalized system.
    The multiplier ${\bm \Lambda}_s$ is diagonal, with each entry a scalar
    function of the corresponding diagonal symbol $d_j$ of ${\bm D}^{\bC}_{>0}$.
    The resulting Poisson-bracket identity eliminates the potentially
    derivative-losing subprincipal commutator and yields the closed Riccati
    estimate of \cref{prop:modified_energy_estimate_diagonalized_3d}; the
    associated uniform a priori bound is
    \cref{cor:uniform_apriori_diagonalized_energy_3d}.

    \item \cref{sec:proof_main_theorem_3d} applies the same
    reductions as linear operators on an increment over a prescribed
    background. The resulting flow estimate and  the tame
    dependence of \(\mathfrak F_p\pare Y\)  yield the local
    Cauchy theory. This local theory is then combined with
    the inverse-size a priori estimate of
    \cref{prop:uniform_apriori_sobolev_energy_3d}.  The Hamiltonian part
    remains of size \(O\pare{\varepsilon}\), whereas the missing
    undifferentiated low-frequency component of \(\tilde\psi\) is controlled
    by \cref{lem:rough_low_frequency_potential_propagation_3d}; the
    continuation criterion then gives
    \(T_\varepsilon\gtrsim\varepsilon^{-1}\).

\end{itemize}

\subsection{A real Alinhac good unknown}
\label{sec:real_alinhac_good_unknown_3d}

The purpose of this subsection is to remove the Alinhac block
\({\bm B}\) displayed in \cref{eq:mathpzcB}.  We work directly on the real paradifferential system
\cref{eq:paralinearazationetapsi}, written for
\[
    Y\defeq \vect{\eta}{\tilde\psi}.
\]
The transformation is triangular and nilpotent.  It leaves the elevation
component unchanged and replaces the potential component by the usual
Alinhac good unknown.

\begin{definition}[Real Alinhac good unknown]
\label{def:real_alinhac_good_unknown_3d}
We define the matrix symbol
\begin{equation}\label{eq:real_alinhac_matrix_symbol_3d}
    {\bm N}
    \pare{\Atw,\background\middle|\eta,\tilde\psi;x}
    \defeq
    \begin{bmatrix}
        0 & 0
        \\
        b_{\textnormal{can}}
        \pare{\Atw,\background\middle|\eta,\tilde\psi;x}
        & 0
    \end{bmatrix}.
\end{equation}
In particular,
\begin{equation}\label{eq:real_alinhac_nilpotency_3d}
    {\bm N}^2=0 .
\end{equation}
The associated paradifferential operator is
\begin{equation}\label{eq:real_alinhac_operator_3d}
    {\bm G}
    \pare{\Atw,\background\middle|\eta,\tilde\psi}
    \defeq
    \Id_{\bC^2}
    +
    \OpBW{
        {\bm N}
        \pare{\Atw,\background\middle|\eta,\tilde\psi;x}
    } .
\end{equation}
Equivalently, the Alinhac variable is
\begin{equation}\label{eq:real_alinhac_unknown_3d}
    Y_{\rm Al}
    \defeq
    \vect{\eta_{\rm Al}}{\tilde\psi_{\rm Al}}
    \defeq
    {\bm G}
    \pare{\Atw,\background\middle|\eta,\tilde\psi}
    \vect{\eta}{\tilde\psi},
\end{equation}
that is,
\begin{equation}\label{eq:real_alinhac_components_3d}
    \eta_{\rm Al}
    =
    \eta,
    \qquad
    \tilde\psi_{\rm Al}
    =
    \tilde\psi
    +
    \OpBW{
        b_{\textnormal{can}}
        \pare{\Atw,\background\middle|\eta,\tilde\psi;x}
    }\eta .
\end{equation}
By \cref{eq:real_alinhac_nilpotency_3d}, the inverse is exact:
\begin{equation}\label{eq:real_alinhac_operator_inverse_3d}
    {\bm G}
    \pare{\Atw,\background\middle|\eta,\tilde\psi}^{-1}
    =
    \Id_{\bC^2}
    -
    \OpBW{
        {\bm N}
        \pare{\Atw,\background\middle|\eta,\tilde\psi;x}
    } .
\end{equation}
\end{definition}

\begin{prop}[Real Alinhac reduction]
\label{prop:real_alinhac_good_unknown_3d}
Let the assumptions of \cref{thm:paraCDEetapsi} hold.  If
\(\pare{\eta,\tilde\psi}\) solves \cref{eq:CDE}, and if \(Y_{\rm Al}\) is
defined by \cref{eq:real_alinhac_unknown_3d}, then
\begin{equation}\label{eq:real_system_after_alinhac_good_unknown_3d}
\begin{aligned}
    \partial_t Y_{\rm Al}
    =
    &\OpBW{
        {\bm Q}
        \pare{\Atw,\background\middle|\eta,\tilde\psi;x,\xi}
        +
        \ii
        \pare{
            \mathpzc V_{\textnormal{can}}
            \pare{\Atw,\background\middle|\eta,\tilde\psi;x}
            \cdot\xi
        }
        \Id_{\bC^2}
        +
        {\bm A}^{\bR,{\rm Al}}_{\rm lo}
        \pare{\Atw,\background\middle|\eta,\tilde\psi;x,\xi}
    }
    Y_{\rm Al}
    \\
    +&
    {\bm R}^\bR_{\bra{\varrho,\geq1}}
    \pare{\Atw,\background\middle|\eta,\tilde\psi}
    Y_{\rm Al}.
\end{aligned}
\end{equation}
The objects entering \cref{eq:real_system_after_alinhac_good_unknown_3d} are
defined as follows.
\begin{itemize}
       \item \(Y_{\rm Al}\) is the real Alinhac unknown defined in \cref{eq:real_alinhac_unknown_3d}; equivalently, its
    components are given in \cref{eq:real_alinhac_components_3d}.  In
    particular, the elevation variable is unchanged, while the potential component is
    corrected by the Alinhac coefficient \(b_{\textnormal{can}}\), defined
    in \cref{eq:bpitilde_main}.

    \item \({\bm Q}\) is the quasilinear off-diagonal symbol defined in
    \cref{eq:mathpzcQ}.  

    \item \(\mathpzc V_{\textnormal{can}}\) is the canonical transport
    velocity defined in \cref{eq:Vpitilde_main}. 

    \item
    \({\bm A}^{\bR,{\rm Al}}_{\rm lo}\) is a generic lower-order matrix
    symbol  satisfying
    \begin{equation}\label{eq:real_alinhac_lower_order_structure_3d}
        {\bm A}^{\bR,{\rm Al}}_{\rm lo}
        \pare{\Atw,\background\middle|\eta,\tilde\psi;x,\xi}
        \in
        \begin{bmatrix}
            \Gamma_{\geq1}^0 & \Gamma_{\geq1}^{-1}
            \vspace{0.1cm}\\
            \Gamma_{\geq1}^0 & \Gamma_{\geq1}^0
        \end{bmatrix}.
    \end{equation}

    \item
    \({\bm R}^\bR_{\bra{\varrho,\geq1}} \in
        \pare{\cR_{\geq1}^\varrho}^{2\times2}\) is a generic smoothing
    remainder.
    \end{itemize}
\end{prop}

\begin{proof}
By \cref{eq:paralinearazationetapsi,eq:mathpzcQ,eq:mathpzcB}, the real
paradifferential system can be written as
\begin{equation}\label{eq:real_system_before_alinhac_proof_3d}
    \partial_t Y
    =
    \OpBW{
        {\bm Q}
        +
        {\bm B}
        +
        {\bm T}
        +
        {\bm A}^{\bR}_{\rm lo}
    }Y
    +
    {\bm R}^\bR_{\bra{\varrho,\geq1}}Y,
\end{equation}
where, in this proof only,
\begin{equation}\label{eq:real_transport_symbol_for_alinhac_proof_3d}
    {\bm T}
    \defeq
    \ii
    \pare{
        \mathpzc V_{\textnormal{can}}
        \pare{\Atw,\background\middle|\eta,\tilde\psi;x}
        \cdot\xi
    }
    \Id_{\bC^2},
\end{equation}
and
\[
    {\bm A}^{\bR}_{\rm lo}
    \in
    \begin{bmatrix}
            \Gamma_{\geq1}^0 & \Gamma_{\geq1}^{-1}
            \vspace{0.1cm}\\
            \Gamma_{\geq1}^0 & \Gamma_{\geq1}^0
        \end{bmatrix}.
\]
The goal is to prove that the triangular good unknown
\[
    {\bm G}
    =
    \Id_{\bC^2}
    +
    \OpBW{{\bm N}}
\]
removes the Alinhac block \({\bm B}\) from the positive-order part. We shall use repeatedly the following consequence of
\cref{def:asympt_symbol,prop:composition_paradiff,lem:composition_rules}:
after choosing the truncation order $\varrho$ in the Bony--Weyl product sufficiently
large,
\begin{equation}\label{eq:real_alinhac_composition_convention_3d}
    \OpBW{{\bm A}}\OpBW{{\bm C}}
    =
    \OpBW{{\bm A}\#_\varrho{\bm C}}
    +
    {\bm R}^\bR_{\bra{\varrho,\geq1}},
\end{equation}
for matrix symbols \({\bm A},{\bm C}\) of the orders occurring below.
Moreover, if \(f\in\cF\subset\Gamma_{\geq1}^0\) is independent of \(\xi\) and
\(a\in\Gamma_{\geq0}^m\), then
\begin{equation}\label{eq:real_alinhac_basic_composition_orders_3d}
    a\#_\varrho f=af+a_{\bra{m-1,\geq1}},
    \qquad
    f\#_\varrho a=fa+a_{\bra{m-1,\geq1}},\qquad a_{\bra{m-1,\geq1}}\in\Gamma_{\geq1}^{m-1}.
\end{equation}
All smoothing errors generated by \cref{eq:real_alinhac_composition_convention_3d}
will be absorbed into the final \({\bm R}^\bR_{\bra{\varrho,\geq1}}\).

\begin{proofpart}[Conjugation of the equation]
Since
\[
    Y_{\rm Al}={\bm G}Y,
    \qquad
    Y={\bm G}^{-1}Y_{\rm Al},
\]
we have
\begin{equation}\label{eq:real_alinhac_time_derivative_conjugation_3d}
    \partial_tY_{\rm Al}
    =
    \pare{\partial_t{\bm G}}{\bm G}^{-1}Y_{\rm Al}
    +
    {\bm G}\partial_tY .
\end{equation}
Inserting \cref{eq:real_alinhac_time_derivative_conjugation_3d} into \cref{eq:real_system_before_alinhac_proof_3d} and using
\cref{eq:real_alinhac_operator_inverse_3d}, we obtain
\begin{equation}\label{eq:real_alinhac_full_operator_conjugation_3d}
\begin{aligned}
    \partial_tY_{\rm Al}
    =
    \Big[
        \pare{\partial_t{\bm G}}{\bm G}^{-1}
        +
        {\bm G}
        \OpBW{
            {\bm Q}+{\bm B}+{\bm T}+{\bm A}^{\bR}_{\rm lo}
        }
        {\bm G}^{-1}
    \Big]Y_{\rm Al}
    +
    {\bm G}{\bm R}^\bR_{\bra{\varrho,\geq1}}{\bm G}^{-1}Y_{\rm Al}.
\end{aligned}
\end{equation}
The first term in the square bracket will be treated separately.  The
second one is the spatial conjugation to which the paradifferential formula
below applies.
\end{proofpart}

\begin{proofpart}[Algebraic cancellation]
Let us first highlight the Alinhac cancellation of the ${\bm B}$ term. We write
\[
    b
    \defeq
    b_{\textnormal{can}},
    \qquad
    \lambda_1
    \defeq
    \lambda_1\pare{\eta;x,\xi},
    \qquad
    q_{12}
    \defeq
    \lambda_1\pare{\eta;x,\xi}
    +
    \Atw\lambda_0\pare{\eta;x,\xi}
\]
and
\begin{equation}\label{eq:real_alinhac_q21_shorthand_3d}
    q_{21}
    \defeq
    -\mathfrak g\Atw
    +
    a^{\pi_0,\eta}_2
    \pare{\rho^\pm\middle|\eta;x,\xi}
    +
    \lambda_1\pare{\eta;x,\xi}
    q^{\tilde{\pi}}_0
    \pare{\Atw,\background\middle|\eta,\tilde\psi;x,\xi}.
\end{equation}
Then \cref{eq:mathpzcQ,eq:mathpzcB} read
\begin{equation}\label{eq:real_QB_shorthand_alinhac_3d}
    {\bm Q}
    =
    \begin{bmatrix}
        0 & q_{12}
        \\
        q_{21} & 0
    \end{bmatrix},
    \qquad
    {\bm B}
    =
    \lambda_1
    \begin{bmatrix}
        b & 0
        \\
        -b^2 & -b
    \end{bmatrix}.
\end{equation}
Moreover, by \cref{eq:real_alinhac_matrix_symbol_3d},
\[
    \Id_{\bC^2}+{\bm N}
    =
    \begin{bmatrix}
        1 & 0
        \\
        b & 1
    \end{bmatrix},
    \qquad
    \Id_{\bC^2}-{\bm N}
    =
    \begin{bmatrix}
        1 & 0
        \\
        -b & 1
    \end{bmatrix}.
\]
At the pointwise symbolic level, isolating the principal upper-right entry
\(\lambda_1\), one has the exact cancellation
\begin{equation}\label{eq:real_alinhac_exact_lambda1_cancellation_3d}
    \begin{bmatrix}
        1 & 0
        \\
        b & 1
    \end{bmatrix}
    \pare{
    \begin{bmatrix}
        0 & \lambda_1
        \\
        q_{21} & 0
    \end{bmatrix}
    +
    \lambda_1
    \begin{bmatrix}
        b & 0
        \\
        -b^2 & -b
    \end{bmatrix}
    }
    \begin{bmatrix}
        1 & 0
        \\
        -b & 1
    \end{bmatrix}
    =
    \begin{bmatrix}
        0 & \lambda_1
        \\
        q_{21} & 0
    \end{bmatrix}.
\end{equation}
For the full upper-right entry
\(q_{12}=\lambda_1+\Atw\lambda_0\), the same pointwise multiplication gives
\begin{equation}\label{eq:real_alinhac_lambda0_residual_3d}
\begin{aligned}
    &
    \begin{bmatrix}
        1 & 0
        \\
        b & 1
    \end{bmatrix}
    \pare{{\bm Q}+{\bm B}}
    \begin{bmatrix}
        1 & 0
        \\
        -b & 1
    \end{bmatrix}
    =
    {\bm Q}
    +
    \Atw\lambda_0 
    \begin{bmatrix}
        -b & 0
        \\
        -b^2 & b
    \end{bmatrix}.
\end{aligned}
\end{equation}
Since \(\lambda_0\in\Gamma_{\geq1}^0\) and \(b\in\cF\subset\Gamma_{\geq1}^0\), the residual
term in \cref{eq:real_alinhac_lambda0_residual_3d} belongs to
\[
    \begin{bmatrix}
            \Gamma_{\geq1}^0 & 0
            \vspace{0.1cm}\\
            \Gamma_{\geq1}^0 & \Gamma_{\geq1}^0
        \end{bmatrix},
\]
and will be absorbed into the lower-order symbol.
\end{proofpart}

\begin{proofpart}[Paradifferential formula for the conjugation]
Let \({\bm A}\) be any matrix symbol of the orders occurring below.  Since
\(\OpBW{{\bm N}}\) has only one nonzero lower-left entry, one has
\[
    \pare{\OpBW{{\bm N}}}^{2}=0,
    \qquad
    {\bm G}^{-1}
    =
    \Id_{\bC^2}-\OpBW{{\bm N}},
\]
as recorded in \cref{eq:real_alinhac_operator_inverse_3d}.  Hence, at the
operator level,
\begin{equation}\label{eq:real_alinhac_exact_operator_conjugation_3d}
    {\bm G}\OpBW{{\bm A}}{\bm G}^{-1}
    =
    \OpBW{{\bm A}}
    +
    \OpBW{{\bm N}}\OpBW{{\bm A}}
    -
    \OpBW{{\bm A}}\OpBW{{\bm N}}
    -
    \OpBW{{\bm N}}\OpBW{{\bm A}}\OpBW{{\bm N}} .
\end{equation}
Applying \cref{eq:real_alinhac_composition_convention_3d} to each product,
and using \cref{prop:composition_paradiff,lem:composition_rules}, gives
\begin{align}
    \OpBW{{\bm N}}\OpBW{{\bm A}}
    &=
    \OpBW{{\bm N}\#_\varrho{\bm A}}
    +
    {\bm R}^\bR_{\bra{\varrho,\geq1}},
    \label{eq:real_alinhac_NA_composition_3d}
    \\
    \OpBW{{\bm A}}\OpBW{{\bm N}}
    &=
    \OpBW{{\bm A}\#_\varrho{\bm N}}
    +
    {\bm R}^\bR_{\bra{\varrho,\geq1}},
    \label{eq:real_alinhac_AN_composition_3d}
    \\
    \OpBW{{\bm N}}\OpBW{{\bm A}}\OpBW{{\bm N}}
    &=
    \OpBW{\pare{{\bm N}\#_\varrho{\bm A}}\#_\varrho{\bm N}}
    +
    {\bm R}^\bR_{\bra{\varrho,\geq1}}.
    \label{eq:real_alinhac_NAN_composition_3d}
\end{align}
In the last line, the smoothing remainder produced by the first composition
is composed once more with the order-zero operator \(\OpBW{{\bm N}}\), and
is still absorbed by \({\bm R}^\bR_{\bra{\varrho,\geq1}}\) after decreasing and
relabelling \(\varrho\).  Here \(\#_\varrho\) denotes the matrix product
whose entries are composed with the scalar Bony--Weyl product, and we denote
\[
    {\bm N}\#_\varrho{\bm A}\#_\varrho{\bm N}
    \defeq
    \pare{{\bm N}\#_\varrho{\bm A}}\#_\varrho{\bm N}.
\]
Therefore
\begin{equation}\label{eq:real_alinhac_operator_conjugated_symbol_3d}
    {\bm G}\OpBW{{\bm A}}{\bm G}^{-1}
    =
    \OpBW{{\bm A}_{\rm conj}}
    +
    {\bm R}^\bR_{\bra{\varrho,\geq1}},
\end{equation}
where
\begin{equation}\label{eq:real_alinhac_symbolic_conjugation_formula_3d}
    {\bm A}_{\rm conj}
    \defeq
    {\bm A}
    +
    {\bm N}\#_\varrho{\bm A}
    -
    {\bm A}\#_\varrho{\bm N}
    -
    {\bm N}\#_\varrho{\bm A}\#_\varrho{\bm N}.
\end{equation}
All smoothing remainders generated in
\cref{eq:real_alinhac_NA_composition_3d,eq:real_alinhac_AN_composition_3d,eq:real_alinhac_NAN_composition_3d}
are absorbed in the generic matrix remainder
\({\bm R}^\bR_{\bra{\varrho,\geq1}}\), after decreasing and relabelling \(\varrho\).
\end{proofpart}

\begin{proofpart}[The quasilinear and Alinhac blocks]
Set
\[
    {\bm P}\defeq{\bm Q}+{\bm B}.
\]
Applying
\cref{eq:real_alinhac_operator_conjugated_symbol_3d,eq:real_alinhac_symbolic_conjugation_formula_3d}
with \({\bm A}={\bm P}\), we have
\[
    {\bm G}\OpBW{{\bm P}}{\bm G}^{-1}
    =
    \OpBW{{\bm P}_{\rm conj}}
    +
    {\bm R}^\bR_{\bra{\varrho,\geq1}},
\]
where
\begin{equation}\label{eq:real_alinhac_P_conjugation_formula_3d}
    {\bm P}_{\rm conj}
    \defeq
    {\bm P}
    +
    {\bm N}\#_\varrho{\bm P}
    -
    {\bm P}\#_\varrho{\bm N}
    -
    {\bm N}\#_\varrho{\bm P}\#_\varrho{\bm N}.
\end{equation}
Using \cref{eq:real_QB_shorthand_alinhac_3d}, we have
\begin{equation}\label{eq:real_P_entries_alinhac_3d}
    {\bm P}
    =
    \begin{bmatrix}
        \lambda_1 b
        &
        q_{12}
        \\
        q_{21}-\lambda_1 b^2
        &
        -\lambda_1 b
    \end{bmatrix}.
\end{equation}
Since
\[
    {\bm N}
    =
    \begin{bmatrix}
        0&0\\
        b&0
    \end{bmatrix},
\]
the matrix structure gives
\begin{align}
    {\bm N}\#_\varrho{\bm P}
    &=
    \begin{bmatrix}
        0 & 0
        \\
        b\#_\varrho\pare{\lambda_1 b}
        &
        b\#_\varrho q_{12}
    \end{bmatrix},
    \label{eq:FsharpP_real_alinhac_3d}
    \\
    {\bm P}\#_\varrho{\bm N}
    &=
    \begin{bmatrix}
        q_{12}\#_\varrho b
        &
        0
        \\
        \pare{-\lambda_1 b}\#_\varrho b
        &
        0
    \end{bmatrix},
    \label{eq:PsharpF_real_alinhac_3d}
    \\
    {\bm N}\#_\varrho{\bm P}\#_\varrho{\bm N}
    &=
    \begin{bmatrix}
        0 & 0
        \\
        \pare{b\#_\varrho q_{12}}\#_\varrho b
        &
        0
    \end{bmatrix}.
    \label{eq:FsharpPsharpF_real_alinhac_3d}
\end{align}
The important point is that \(q_{21}\), although of order \(2\), is never
differentiated by the Bony--Weyl corrections generated by the triangular
conjugation.  Only \(q_{12}\), \(\lambda_1 b\), and \(b\) enter the composed
symbols. Combining
\cref{eq:real_alinhac_P_conjugation_formula_3d,eq:real_P_entries_alinhac_3d,eq:FsharpP_real_alinhac_3d,eq:PsharpF_real_alinhac_3d,eq:FsharpPsharpF_real_alinhac_3d}, we get
\begin{equation}\label{eq:conjugated_P_entries_real_alinhac_3d}
\begin{aligned}
    {\bm P}_{\rm conj}
    ={}&
    \begin{bmatrix}
        \lambda_1 b-q_{12}\#_\varrho b
        &
        q_{12}
        \\
        q_{21}-\lambda_1 b^2
        +b\#_\varrho\pare{\lambda_1 b}
        -\pare{-\lambda_1 b}\#_\varrho b
        -\pare{b\#_\varrho q_{12}}\#_\varrho b
        &
        -\lambda_1 b+b\#_\varrho q_{12}
    \end{bmatrix}.
\end{aligned}
\end{equation}
We now expand each product by
\cref{eq:real_alinhac_basic_composition_orders_3d}.  Since \(b\in\cF\subset\Gamma_{\geq1}^0\),
\(\lambda_1\in\Gamma_{\geq 0}^1\), and
\(q_{12}=\lambda_1+\Atw\lambda_0\in\Gamma_{\geq 0}^1\), one obtains
\begin{align}
    q_{12}\#_\varrho b
    &=
    q_{12} b
    +
    a_{\bra{0,\geq1}},
    \label{eq:q12sharpb_real_alinhac_order_3d}
    \\
    b\#_\varrho q_{12}
    &=
q_{12}b
    +
    a_{\bra{0,\geq1}},
    \label{eq:bsharpq12_real_alinhac_order_3d}
    \\
    b\#_\varrho\pare{\lambda_1 b}
    &=
    \lambda_1 b^2
    +
    a_{\bra{0,\geq1}},
    \label{eq:bsharplambdab_real_alinhac_order_3d}
    \\
    \pare{-\lambda_1 b}\#_\varrho b
    &=
    -\lambda_1 b^2
    +
    a_{\bra{0,\geq1}},
    \label{eq:lambdabsharpb_real_alinhac_order_3d}
    \\
    \pare{b\#_\varrho q_{12}}\#_\varrho b
    &=
    q_{12}b^2
    +
    a_{\bra{0,\geq1}}.
    \label{eq:bsharpq12sharpb_real_alinhac_order_3d}
\end{align}
Here \(a_{\bra{0,\geq1}}\in\Gamma_{\geq1}^0\) is generic: the first neglected Weyl
correction contains one \(\xi\)-derivative of an order-one symbol, while
all further corrections have lower order. 
Substituting
\cref{eq:q12sharpb_real_alinhac_order_3d,eq:bsharpq12_real_alinhac_order_3d,eq:bsharplambdab_real_alinhac_order_3d,eq:lambdabsharpb_real_alinhac_order_3d,eq:bsharpq12sharpb_real_alinhac_order_3d}
into \cref{eq:conjugated_P_entries_real_alinhac_3d} yields
\begin{equation}\label{eq:real_alinhac_P_conjugated_final_symbol_3d}
    {\bm P}_{\rm conj}
    =
    {\bm Q}
    +
    \Atw\lambda_0
    \begin{bmatrix}
        -b & 0
        \\
        -b^2 & b
    \end{bmatrix}
    +
    \begin{bmatrix}
        \Gamma_{\geq1}^0 & 0
        \vspace{0.1cm}\\
        \Gamma_{\geq1}^0 & \Gamma_{\geq1}^0
    \end{bmatrix}={\bm Q}
    +
    \begin{bmatrix}
        \Gamma_{\geq1}^0 & 0
        \vspace{0.1cm}\\
        \Gamma_{\geq1}^0 & \Gamma_{\geq1}^0
    \end{bmatrix}.
\end{equation}
% Since \(\lambda_0\in\Gamma^0\), the second term on the right-hand side of
% \cref{eq:real_alinhac_P_conjugated_final_symbol_3d} is also of order \(0\).
% Therefore
% \[
%     {\bm P}_{\rm conj}
%     =
%     {\bm Q}
%     +
%     {\bm A}^{\bR,{\rm Al},P}_{\rm lo},
%     \qquad
%     {\bm A}^{\bR,{\rm Al},P}_{\rm lo}
%     \in
%     \begin{bmatrix}
%         \Gamma^0 & \Gamma^{-1}
%         \\
%         \Gamma^0 & \Gamma^0
%     \end{bmatrix},
% \]
% where the zero upper-right entry is harmlessly included in \(\Gamma^{-1}\).
\end{proofpart}

\begin{proofpart}[The transport block]
Let
\[
    \tau
    \defeq
    \ii
    \mathpzc V_{\textnormal{can}}
    \pare{\Atw,\background\middle|\eta,\tilde\psi;x}
    \cdot\xi,
    \qquad
    {\bm T}\defeq \tau\Id_{\bC^2}.
\]
Applying
\cref{eq:real_alinhac_operator_conjugated_symbol_3d,eq:real_alinhac_symbolic_conjugation_formula_3d}
with \({\bm A}={\bm T}\), we have
\[
    {\bm G}\OpBW{{\bm T}}{\bm G}^{-1}
    =
    \OpBW{{\bm T}_{\rm conj}}
    +
    {\bm R}^\bR_{\bra{\varrho,\geq1}},
\]
where
\begin{equation}\label{eq:real_alinhac_T_conjugation_formula_3d}
    {\bm T}_{\rm conj}
    \defeq
    {\bm T}
    +
    {\bm N}\#_\varrho{\bm T}
    -
    {\bm T}\#_\varrho{\bm N}
    -
    {\bm N}\#_\varrho{\bm T}\#_\varrho{\bm N}.
\end{equation}
Since \({\bm T}\) is scalar in the matrix indices, it commutes with
\({\bm N}\) at the pointwise principal-symbol level:
\[
    {\bm N}{\bm T}-{\bm T}{\bm N}=0,
    \qquad
    {\bm N}{\bm T}{\bm N}=\tau{\bm N}^2=0.
\]
At the Bony--Weyl symbolic level,
\[
    {\bm N}\#_\varrho{\bm T}
    -
    {\bm T}\#_\varrho{\bm N}
    =
    \begin{bmatrix}
        0&0\\
        b\#_\varrho\tau-\tau\#_\varrho b&0
    \end{bmatrix},
\]
and \({\bm N}\#_\varrho{\bm T}\#_\varrho{\bm N}=0\) by the same lower-left
matrix structure. Since we have that \(b\) is independent of \(\xi\) and $ \partial_\xi\tau
    =
    \ii\mathpzc V_{\textnormal{can}}
    \pare{\Atw,\background\middle|\eta,\tilde\psi;x}$, then using the first Weyl correction in
\cref{def:asympt_symbol} and then the composition rules in
\cref{prop:composition_paradiff}, we get
\[
    b\#_\varrho\tau-\tau\#_\varrho b
    =
    -
    \mathpzc V_{\textnormal{can}}
    \pare{\Atw,\background\middle|\eta,\tilde\psi;x}
    \cdot\nabla b
    +
    a_{\bra{-1,\geq1}}.
\]
Consequently, by \cref{eq:real_alinhac_T_conjugation_formula_3d},
\begin{equation}\label{eq:transport_conjugation_real_alinhac_3d}
    {\bm T}_{\rm conj}
    =
    {\bm T}
    +
    \begin{bmatrix}
        0&0\\
        -\mathpzc V_{\textnormal{can}}
        \pare{\Atw,\background\middle|\eta,\tilde\psi;x}\cdot\nabla b&0
    \end{bmatrix}
    +
    \begin{bmatrix}
        0&0\\
        \Gamma_{\geq1}^{-1}&0
    \end{bmatrix}={\bm T}
    +\begin{bmatrix}
        0&0\\
        \Gamma_{\geq1}^{0}&0
    \end{bmatrix}.
\end{equation}
Thus the transport symbol is unchanged at order \(1\), and all new
transport contributions are lower order.
\end{proofpart}

\begin{proofpart}[The original lower-order symbol]
Applying
\cref{eq:real_alinhac_operator_conjugated_symbol_3d,eq:real_alinhac_symbolic_conjugation_formula_3d}
with \({\bm A}={\bm A}^{\bR}_{\rm lo}\), we have
\[
    {\bm G}\OpBW{{\bm A}^{\bR}_{\rm lo}}{\bm G}^{-1}
    =
    \OpBW{{\bm A}^{\bR,{\rm Al},{\rm lo}}_{\rm conj}}
    +
    {\bm R}^\bR_{\bra{\varrho,\geq1}},
\]
where
\begin{equation}\label{eq:real_alinhac_Alo_conjugation_formula_3d}
    {\bm A}^{\bR,{\rm Al},{\rm lo}}_{\rm conj}
    \defeq
    {\bm A}^{\bR}_{\rm lo}
    +
    {\bm N}\#_\varrho{\bm A}^{\bR}_{\rm lo}
    -
    {\bm A}^{\bR}_{\rm lo}\#_\varrho{\bm N}
    -
    {\bm N}\#_\varrho{\bm A}^{\bR}_{\rm lo}\#_\varrho{\bm N}.
\end{equation}
Write
\[
    {\bm A}^{\bR}_{\rm lo}
    =
    \begin{bmatrix}
        a_{11} & a_{12}
        \\
        a_{21} & a_{22}
    \end{bmatrix},
    \qquad
    a_{11},a_{21},a_{22}\in\Gamma_{\geq1}^0,
    \qquad
    a_{12}\in\Gamma_{\geq1}^{-1}.
\]
The pointwise products generated by the triangular matrix are
\[
    {\bm N}{\bm A}^{\bR}_{\rm lo}
    =
    \begin{bmatrix}
        0&0\\
        ba_{11}&ba_{12}
    \end{bmatrix},
    \qquad
    {\bm A}^{\bR}_{\rm lo}{\bm N}
    =
    \begin{bmatrix}
        a_{12}b&0\\
        a_{22}b&0
    \end{bmatrix},
    \qquad
     {\bm N}{\bm A}^{\bR}_{\rm lo}{\bm N}
    =
    \begin{bmatrix}
        0&0\\
        ba_{12}b&0
    \end{bmatrix}.
\]
Since \(b\in\cF\subset\Gamma_{\geq1}^0\), these products have the same matrix
order structure as \({\bm A}^{\bR}_{\rm lo}\).  The Bony--Weyl corrections
only differentiate symbols in \(x\) and \(\xi\), and every
\(\xi\)-derivative lowers the symbolic order.  Therefore, by
\cref{prop:composition_paradiff,lem:composition_rules},
\begin{equation}\label{eq:lower_order_conjugation_real_alinhac_3d}
    {\bm A}^{\bR,{\rm Al},{\rm lo}}_{\rm conj}
    \in
    \begin{bmatrix}
        \Gamma_{\geq1}^0 & \Gamma_{\geq1}^{-1}
        \vspace{0.1cm}\\
        \Gamma_{\geq1}^0 & \Gamma_{\geq1}^0
    \end{bmatrix}.
\end{equation}
This contribution is included in
\({\bm A}^{\bR,{\rm Al}}_{\rm lo}\).
\end{proofpart}

\begin{proofpart}[The time derivative of the good unknown]
It remains to treat the first term in
\cref{eq:real_alinhac_full_operator_conjugation_3d}.  Note that
\[
    {\bm N}
    =
    \begin{bmatrix}
        0&0\\
        b&0
    \end{bmatrix},
    \qquad
    \partial_t{\bm N}
    =
    \begin{bmatrix}
        0&0\\
        \partial_t b&0
    \end{bmatrix}\in\begin{bmatrix}
        0 & 0\\
        \Gamma_{\geq1}^0 & 0
    \end{bmatrix}.
\]
In particular, the product of the two corresponding lower-left matrix operators is exactly
zero:
\[
    \OpBW{\partial_t{\bm N}}\OpBW{{\bm N}}=0.
\]
Thus, using the exact inverse
\({\bm G}^{-1}=\Id_{\bC^2}-\OpBW{{\bm N}}\),
\begin{equation}\label{eq:time_derivative_G_real_alinhac_3d}
    \pare{\partial_t{\bm G}}{\bm G}^{-1}
    =
    \OpBW{\partial_t{\bm N}}.
\end{equation}
This term is included in
\(\OpBW{{\bm A}^{\bR,{\rm Al}}_{\rm lo}}\).\\

\end{proofpart}

\begin{proofpart}[Smoothing remainders and conclusion]
The map
\[
    {\bm A}
    \mapsto
    {\bm A}
    +
    {\bm N}\#_\varrho{\bm A}
    -
    {\bm A}\#_\varrho{\bm N}
    -
    {\bm N}\#_\varrho{\bm A}\#_\varrho{\bm N}
\]
appearing in \cref{eq:real_alinhac_symbolic_conjugation_formula_3d} is
linear in \({\bm A}\).  Therefore the spatial conjugation in
\cref{eq:real_alinhac_full_operator_conjugation_3d} is obtained by summing
the three contributions computed above with
\({\bm A}={\bm P},{\bm T},{\bm A}^{\bR}_{\rm lo}\).

It remains to control the smoothing remainder in
\cref{eq:real_system_before_alinhac_proof_3d}.  The operators
\({\bm G}^{\pm1}\) are the identity plus an order-zero paradifferential
operator with coefficient controlled by the same tame norms as
\(b_{\textnormal{can}}\).  Therefore, by \cref{lem:composition_rules},
\[
    {\bm G}{\bm R}^\bR_{\bra{\varrho,\geq1}}{\bm G}^{-1}
    \in
    \pare{\cR_{\geq1}^\varrho}^{2\times2},
\]
after decreasing and relabelling \(\varrho\).  The same convention absorbs
all smoothing remainders produced by the symbolic compositions in
\cref{eq:real_alinhac_composition_convention_3d}.

Collecting
\cref{eq:real_alinhac_P_conjugated_final_symbol_3d,eq:transport_conjugation_real_alinhac_3d,eq:lower_order_conjugation_real_alinhac_3d,eq:time_derivative_G_real_alinhac_3d}
inside the conjugation identity
\cref{eq:real_alinhac_full_operator_conjugation_3d}, we obtain
\begin{equation}\label{Alinhac:final}
    \partial_tY_{\rm Al}
    =
    \OpBW{
        {\bm Q}
        +
        {\bm T}
        +
        {\bm A}^{\bR,{\rm Al}}_{\rm lo}
    }
    Y_{\rm Al}
    +
    {\bm R}^\bR_{\bra{\varrho,\geq1}}Y_{\rm Al},
\end{equation}
with
\[
    {\bm A}^{\bR,{\rm Al}}_{\rm lo}
    \in
    \begin{bmatrix}
        \Gamma_{\geq1}^0 & \Gamma_{\geq1}^{-1}
        \vspace{0.1cm}\\
        \Gamma_{\geq1}^0 & \Gamma_{\geq1}^0
    \end{bmatrix},
    \qquad
    {\bm R}^\bR_{\bra{\varrho,\geq1}}\in
    \pare{\cR_{\geq1}^\varrho}^{2\times2}.
\]
Recalling the definition of \({\bm T}\) in
\cref{eq:real_transport_symbol_for_alinhac_proof_3d}, this is exactly
\cref{eq:real_system_after_alinhac_good_unknown_3d}.
\end{proofpart}
\end{proof}

\subsection{Complex formulation}
\label{sec:complex_formulation_3d_vs}
Throughout this subsection we assume
\cref{hyp:strict_atwood_weber_stability}.  Hence, by
\cref{prop:purely_imaginary_linear_frequencies},
\(
    \mathfrak d\pare{\Atw,\background\middle|\xi}<0
\)
for every \(\xi\in\bR^2\setminus\set{0}\), and the two branches of the flat
linearized CDE are purely imaginary and simple.  We use the positive frequency
introduced below to complexify the real Alinhac unknown and diagonalize the
positive-order flat dynamics.

\begin{definition}[Complex Fourier multipliers and non-homogeneous coordinate]
\label{def:complex_linear_multipliers_3d}
\label{def:complex_diagonalizing_matrix_3d}
Assume \cref{hyp:strict_atwood_weber_stability}, fix the even low-frequency
cut-off \(\chi\) of \cref{def-chi}, and let
\(\xi\in\bR^2\setminus\set{0}\).  We introduce the following objects.
\begin{enumerate}
    \item \emph{Homogeneous dispersive multipliers.}
    The order-\(3/2\) dispersive frequency is
    \begin{equation}\label{eq:complex_frequency_order_three_halves}
        \omega_{\frac32}
        \pare{\Atw,\background\middle|\xi}
        \defeq
        \sqrt{-\mathfrak d_{\Atw,\background}\pare{\xi}}
        \in\tilde\Gamma^{\frac32}.
    \end{equation}
    In particular, by \cref{eq:def_linearized_discriminant},
    \begin{equation}\label{eq:complex_frequency_square_identity}
        \omega_{\frac32}^{2}
        \pare{\Atw,\background\middle|\xi}
        =
        -\frac{\av{\xi}}{2}
        \mathfrak q\pare{\Atw,\background\middle|\xi}.
    \end{equation}
    Thus the two branches in \cref{eq:linearized_eigenvalues} are obtained by
    replacing \(\sqrt{\mathfrak d}\) with \(\ii\omega_{\frac32}\).

    The order-\((-1/4)\) symmetrizing multiplier is
    \begin{equation}\label{eq:complex_multiplier_p}
        \mathfrak m_{-\frac14}
        \pare{\Atw,\background\middle|\xi}
        \defeq
        \sqrt{
            \frac{\av{\xi}}{
                2\omega_{\frac32}
                \pare{\Atw,\background\middle|\xi}
            }
        }
        \in\tilde\Gamma^{-\frac14}.
    \end{equation}
    Equivalently,
    \begin{equation}\label{eq:complex_multiplier_p_equivalent}
        \mathfrak m_{-\frac14}^{2}
        =
        \frac{\frac12\av{\xi}}{\omega_{\frac32}}
        =
        -\frac{\omega_{\frac32}}{\mathfrak q},
        \qquad
        \mathfrak m_{-\frac14}^{-2}
        =
        \frac{2\omega_{\frac32}}{\av{\xi}},
    \end{equation}
    where all quantities in \cref{eq:complex_multiplier_p_equivalent} are
    evaluated at \(\pare{\Atw,\background\middle|\xi}\).

    The homogeneous leading capillary frequency and symmetrizer are
    \begin{equation}\label{eq:complex_leading_capillary_frequency_3d}
        \dot\omega_{\frac32}
        \pare{\Atw,\background\middle|\xi}
        \defeq
        \sqrt{\frac{\sigma}{2\aangle{\rho^\pm}}}
        \av{\xi}^{\frac32}
        \in\tilde\Gamma^{\frac32}
    \end{equation}
    and
    \begin{equation}\label{eq:complex_leading_multiplier_p_3d}
        \dot{\mathfrak m}_{-\frac14}
        \pare{\Atw,\background\middle|\xi}
        \defeq
        \sqrt{
            \frac{\av{\xi}}{
                2\dot\omega_{\frac32}
                \pare{\Atw,\background\middle|\xi}
            }
        }
        \in\tilde\Gamma^{-\frac14}.
    \end{equation}
    In particular,
    \begin{equation}\label{eq:complex_leading_multiplier_expansions_3d}
        \omega_{\frac32}-\dot\omega_{\frac32}
        \in\tilde\Gamma^{\frac12},
        \qquad
        \mathfrak m_{-\frac14}^{2}
        -\dot{\mathfrak m}_{-\frac14}^{2}
        \in\tilde\Gamma^{-\frac32},
        \qquad
        \mathfrak m_{-\frac14}^{-2}
        -\dot{\mathfrak m}_{-\frac14}^{-2}
        \in\tilde\Gamma^{-\frac12}.
    \end{equation}

    \item \emph{Non-homogeneous complex coordinate and flat defect.}
    Define, initially for \(\xi\neq0\),
    \begin{equation}\label{eq:complex_nonhomogeneous_multiplier_p_3d}
        \mathfrak m^{\rm nh}_{-\frac14}
        \pare{\Atw,\background\middle|\xi}
        \defeq
        \chi\pare{\xi}
        +
        \pare{1-\chi\pare{\xi}}
        \mathfrak m_{-\frac14}
        \pare{\Atw,\background\middle|\xi},
    \end{equation}
    and extend it smoothly through \(\xi=0\), where it equals \(1\).  In the
    following two matrix displays, every occurrence of
    \(\mathfrak m^{\rm nh}_{-\frac14}\) is evaluated at
    \(\pare{\Atw,\background\middle|\xi}\).  The non-homogeneous complex
    diagonalizing matrix is
    \begin{equation}\label{eq:complex_matrix_M}
        {\bm M}
        \pare{\Atw,\background\middle|\xi}
        \defeq
        \frac{1}{\sqrt2}
        \begin{bmatrix}
            \pare{\mathfrak m^{\rm nh}_{-\frac14}}^{-1}
            &
            \ii\mathfrak m^{\rm nh}_{-\frac14}
            \\
            \pare{\mathfrak m^{\rm nh}_{-\frac14}}^{-1}
            &
            -\ii\mathfrak m^{\rm nh}_{-\frac14}
        \end{bmatrix}.
    \end{equation}
    Its inverse is
    \begin{equation}\label{eq:complex_matrix_M_inverse}
        {\bm M}
        \pare{\Atw,\background\middle|\xi}^{-1}
        =
        \frac{1}{\sqrt2}
        \begin{bmatrix}
            \mathfrak m^{\rm nh}_{-\frac14}
            &
            \mathfrak m^{\rm nh}_{-\frac14}
            \\
            -\ii\pare{\mathfrak m^{\rm nh}_{-\frac14}}^{-1}
            &
            \ii\pare{\mathfrak m^{\rm nh}_{-\frac14}}^{-1}
        \end{bmatrix}.
    \end{equation}
    The complex coordinate is
    \begin{equation}\label{eq:complex_variable_definition_3d}
        U
        \defeq
        \vect{u}{\bar u}
        \defeq
        {\bm M}\pare{\Atw,\background\middle|D}Y_{\rm Al},
        \qquad
        Y_{\rm Al}
        \defeq
        \vect{\eta_{\rm Al}}{\tilde\psi_{\rm Al}},
    \end{equation}
    where \(Y_{\rm Al}\) is defined in
    \cref{eq:real_alinhac_unknown_3d}.

    Finally, define the diagonal sign matrix
    \begin{equation}\label{eq:complex_diagonal_sign_matrix}
        {\bm J}_{\rm d}
        \defeq
        \begin{bmatrix}
            -1&0
            \\
            0&1
        \end{bmatrix}
    \end{equation}
    and the flat low-frequency defect
    \begin{equation}\label{eq:linear_low_frequency_defect_complex_3d}
    \begin{aligned}
        {\bm A}^{\bC,{\rm lin}}_{\rm lo}
        \pare{\Atw,\background\middle|\xi}
        \defeq{}&
        {\bm M}\pare{\Atw,\background\middle|\xi}
        \begin{bmatrix}
            0&\dfrac{\av{\xi}}{2}
            \\
            \mathfrak q\pare{\Atw,\background\middle|\xi}&0
        \end{bmatrix}
        {\bm M}\pare{\Atw,\background\middle|\xi}^{-1}
        \\
        &-
        \pare{1-\chi\pare{\xi}}
        \ii\omega_{\frac32}
        \pare{\Atw,\background\middle|\xi}
        {\bm J}_{\rm d}.
    \end{aligned}
    \end{equation}
\end{enumerate}
\end{definition}

For later use, we recover the real Alinhac variable from \(U\) by
\begin{equation}\label{eq:alinhac_variables_as_functions_of_U}
    \vect{\eta_{\rm Al}\pare{U}}{\tilde\psi_{\rm Al}\pare{U}}
    \defeq
    {\bm M}\pare{\Atw,\background\middle|D}^{-1}U.
\end{equation}
The coefficients of the real Alinhac-reduced system are evaluated on the
physical variables
\begin{equation}\label{eq:eta_psitilde_as_functions_of_U}
    \vect{\eta\pare{U}}{\tilde\psi\pare{U}}
    \defeq
    {\bm G}
    \pare{
        \Atw,\background
        \middle|
        \eta\pare{U},\tilde\psi\pare{U}
    }^{-1}
    \vect{\eta_{\rm Al}\pare{U}}{\tilde\psi_{\rm Al}\pare{U}},
\end{equation}
where \({\bm G}\) is defined in \cref{eq:real_alinhac_operator_3d}.  More
explicitly, by
\cref{eq:real_alinhac_components_3d,eq:real_alinhac_operator_inverse_3d},
\begin{equation}\label{eq:physical_inverse_alinhac_components_complex_3d}
    \eta\pare{U}
    =
    \eta_{\rm Al}\pare{U},
    \qquad
    \tilde\psi\pare{U}
    =
    \tilde\psi_{\rm Al}\pare{U}
    -
    \OpBW{
        b_{\textnormal{can}}
        \pare{
            \Atw,\background
            \middle|
            \eta\pare{U},\tilde\psi\pare{U};x
        }
    }
    \eta_{\rm Al}\pare{U}.
\end{equation}

\begin{prop}[Complex form of the Alinhac-reduced paradifferential system]
\label{prop:complex_paradifferential_system_3d}
Let the assumptions of \cref{thm:paraCDEetapsi} and
\cref{hyp:strict_atwood_weber_stability} hold.  If
\(\pare{\eta,\tilde\psi}\) solves the canonical CDE \cref{eq:CDE}, if
\(Y_{\rm Al}\) is the real Alinhac unknown defined in
\cref{eq:real_alinhac_unknown_3d}, and if \(U\) is defined by
\cref{eq:complex_variable_definition_3d}, then
\begin{equation}\label{eq:complex_paradifferential_system_3d}
\begin{aligned}
    \partial_t U
    ={}&
    \OpBW{
     {\bm A}^{\bC,{\rm lin}}_{\rm lo}
        \pare{\Atw,\background\middle|\xi}
        +
        \pare{1-\chi\pare{\xi}}\,
        \ii\omega_{\frac32}\pare{\Atw,\background\middle|\xi}
        {\bm J}_{\rm d}
        \pare{
            \Id_{\bC^2}
            +
            \dot{\bm E}^{\bC}_{{\rm ql},0}
            \pare{U;x,\hat{\xi}}
        }
    }
    U
    \\
    &+
    \OpBW{
        \ii
        \pare{
            \mathpzc V_{\textnormal{can}}
            \pare{
                \Atw,\background
                \middle|
                \eta\pare{U},\tilde\psi\pare{U};
                x
            }
            \cdot\xi
        }
        \Id_{\bC^2}
    }
    U
    \\
    &+
    \OpBW{
        \pare{1-\chi\pare{\xi}}\,
        \dot{\bm A}^{\bC}_{0}
        \pare{\Atw,\background\middle|U;x,\hat{\xi}}
        \av{\xi}^{\frac12}
    }
    U
    \\
    &+
    \OpBW{
        {\bm A}_{\bra{0,\geq1}}
        \pare{\Atw,\background\middle|U;x,\xi}
    }
    U
    +
    {\bm R}_{\bra{\varrho,\geq1}}
    \pare{\Atw,\background\middle|U}U.
\end{aligned}
\end{equation}
where:
\begin{itemize}
    \item \(U=\vect{u}{\bar u}\) is the complex Alinhac-reduced unknown
    defined in \cref{eq:complex_variable_definition_3d}; the physical
    variables \(\eta\pare{U}\) and \(\tilde\psi\pare{U}\) are recovered
    through the inverse real Alinhac map in
    \cref{eq:eta_psitilde_as_functions_of_U,eq:physical_inverse_alinhac_components_complex_3d};
    \item
    \({\bm A}^{\bC,{\rm lin}}_{\rm lo}
    \pare{\Atw,\background\middle|\xi}\in\pare{\tilde\Gamma^0}^{2\times2}\)
    is the Fourier-localized symbol introduced in \cref{eq:linear_low_frequency_defect_complex_3d}.
    \item \(\omega_{\frac32}\pare{\Atw,\background\middle|\xi}
    \in\tilde\Gamma^{3/2}\) is the scalar Fourier multiplier defined in
    \cref{eq:complex_frequency_order_three_halves};

\item \({\bm J}_{\rm d}\) is the diagonal sign matrix defined in
    \cref{eq:complex_diagonal_sign_matrix};

   \item \(\dot{\bm E}^{\bC}_{{\rm ql},0}
    \pare{U;x,\hat{\xi}}
    \in\pare{\Gamma_{\geq1}^0}^{2\times2}\)
    is the homogeneous degree-zero quasilinear matrix coefficient
    \begin{equation}\label{eq:complex_ql_matrix_coefficient_3d}
        \dot{\bm E}^{\bC}_{{\rm ql},0}
        \pare{U;x,\hat{\xi}}
        \defeq
        \begin{bmatrix}
            \dot e_0 & \dot f_0
            \\
            \dot f_0 & \dot e_0
        \end{bmatrix}
        \pare{U;x,\hat{\xi}} ,
    \end{equation}
    where the scalar homogeneous coefficients are 
    \begin{itemize}
                \item $\dot e_0\in \Gamma_{\geq1}^0$ and is explicitly defined as
        \begin{equation}
             \dot e_0
        \pare{U;x,\hat{\xi}}
        \defeq
    \frac12
    \pare{
        \sqrt{
            1+\av{\nabla^\perp\eta\pare{U}\cdot\hat{\xi}}^2
        }
        +
        \frac{
            1+\av{\nabla^\perp\eta\pare{U}\cdot\hat{\xi}}^2
        }{
            \pare{1+\av{\nabla\eta\pare{U}}^2}^{\frac32}
        }
    }-1
    \label{eq:dote0_closed_principal_3d};
        \end{equation}

        \item  $\dot f_0\in \Gamma_{\geq1}^0$ and is explicitly defined as
        \begin{equation}
              \dot f_0
        \pare{U;x,\hat{\xi}}
        \defeq\frac12
    \pare{
        \frac{
            1+\av{\nabla^\perp\eta\pare{U}\cdot\hat{\xi}}^2
        }{
            \pare{1+\av{\nabla\eta\pare{U}}^2}^{\frac32}
        }
        -
        \sqrt{
            1+\av{\nabla^\perp\eta\pare{U}\cdot\hat{\xi}}^2
        }
    };\label{eq:dotf0_closed_principal_3d}
        \end{equation}
    \end{itemize}
    \item the transport is kept in the explicit scalar form
    \[
        \ii
        \pare{
            \mathpzc V_{\textnormal{can}}
            \pare{
                \Atw,\background
                \middle|
                \eta\pare{U},\tilde\psi\pare{U};
                x
            }
            \cdot\xi
        }
        \Id_{\bC^2},
    \]
    where \(\mathpzc V_{\textnormal{can}}\) is defined in
    \cref{eq:Vpitilde_main};

    \item
    \(\dot{\bm A}^{\bC}_{0}
    \pare{\Atw,\background\middle|U;x,\hat{\xi}}
    \in\pare{\Gamma_{\geq1}^0}^{2\times2}\)
    is the homogeneous degree-zero coefficient of the explicit order-\(1/2\)
    block. It is defined by
    \begin{equation}\label{eq:complex_half_order_symbol_definition_3d}
        \dot{\bm A}^{\bC}_{0}
        \pare{\Atw,\background\middle|U;x,\hat{\xi}}
        \defeq
        \ii
        \begin{bmatrix}
            \dot\alpha^{\bC}_{0}
            &
            \dot\beta^{\bC}_{0}
            \\
            -\dot\beta^{\bC}_{0}
            &
            -\dot\alpha^{\bC}_{0}
        \end{bmatrix}
        \pare{\Atw,\background\middle|U;x,\hat{\xi}} ,
    \end{equation}
    where
    \begin{small}
    \begin{subequations}
    \begin{align}
        \dot\alpha^{\bC}_{0}
        \pare{\Atw,\background\middle|U;x,\hat{\xi}}
        &\defeq
        \frac{
            \dot{\mathfrak b}^{\bC}_{0}
            -
            \dot{\mathfrak a}^{\bC}_{0}
        }{2}
        \pare{\Atw,\background\middle|U;x,\hat{\xi}},
        \label{eq:complex_half_alpha_coefficient_3d}
        \\
        \dot\beta^{\bC}_{0}
        \pare{\Atw,\background\middle|U;x,\hat{\xi}}
        &\defeq
        \frac{
            \dot{\mathfrak a}^{\bC}_{0}
            +
            \dot{\mathfrak b}^{\bC}_{0}
        }{2}
        \pare{\Atw,\background\middle|U;x,\hat{\xi}},
        \label{eq:complex_half_beta_coefficient_3d}
        \\
        \dot{\mathfrak a}^{\bC}_{0}
        \pare{\Atw,\background\middle|U;x,\hat{\xi}}
        &\defeq
        \Atw
        \sqrt{
            \frac{2\sigma}{\aangle{\rho^\pm}}
        }\,
        \lambda_0
        \pare{\eta\pare{U};x,\xi},
        \label{eq:complex_a_half_scalar_3d}
        \\
        \dot{\mathfrak b}^{\bC}_{0}
        \pare{\Atw,\background\middle|U;x,\hat{\xi}}
        &\defeq
        \frac{4-\Atw^2}{16}
        \sqrt{
            \frac{2\aangle{\rho^\pm}}{\sigma}
        }\,
        \sqrt{
            1+
            \av{
                \nabla^\perp\eta\pare{U}
                \cdot\hat{\xi}
            }^2
        }
        \nonumber
        \\
        &\quad\times
        \bra{
            4
            \pare{
                A_{0}^{\mathpzc{B}_0,\omega;\perp}
                \pare{\eta\pare{U};x,\xi}
                \cdot
                \Omega
                \pare{
                    \background
                    \middle|
                    \eta\pare{U},\tilde\psi\pare{U};
                    x
                }
            }^2
            -
            \frac{
                \pare{\background^\perp\cdot\hat{\xi}}^2
                \sqrt{
                    1+
                    \av{
                        \nabla^\perp\eta\pare{U}
                        \cdot\hat{\xi}
                    }^2
                }
            }{
                \pare{
                    1+
                    \av{\nabla\eta\pare{U}}^2
                }^{\frac32}
            }
        }.
        \label{eq:complex_b_half_scalar_3d}
    \end{align}
    \end{subequations}
    \end{small}
    Here \(\lambda_0\) is defined in \cref{eq:lambda0}, \(A_{0}^{\mathpzc{B}_0,\omega}\) is defined
    in \cref{eq:H0_definition}, and \(\Omega\) is defined in
    \cref{eq:Omegacan_paraCDE}.

    \item \({\bm A}_{\bra{0,\geq1}}
    \in\pare{\Gamma_{\geq1}^0}^{2\times2}\) is a generic matrix symbol as per
    \cref{notation:paralinearization};

    \item \({\bm R}_{\bra{\varrho,\geq1}}
    \in\pare{\cR_{\geq1}^\varrho}^{2\times2}\) is a generic smoothing remainder as
    per \cref{notation:paralinearization}.
\end{itemize}
\end{prop}

\begin{lemma}[High-frequency diagonalization of the constant-coefficient linear part]
\label{prop:linear_diagonalization_complex_3d}
Assume \cref{hyp:strict_atwood_weber_stability}.  Then
\({\bm A}^{\bC,{\rm lin}}_{\rm lo}\), defined in
\cref{eq:linear_low_frequency_defect_complex_3d}, belongs to
\(\pare{\tilde\Gamma^0}^{2\times2}\), and
\begin{equation}\label{eq:linear_diagonalization_identity_complex_3d}
    {\bm M}\pare{\Atw,\background\middle|\xi}
    {\bm L}_0\pare{\xi}
    {\bm M}\pare{\Atw,\background\middle|\xi}^{-1}
    =
    \ii\frac{\Atw}{4}
    \pare{\background^{\perp}\cdot\xi}
    \Id_{\bC^2}
    +
    \pare{1-\chi\pare{\xi}}
    \ii\omega_{\frac32}
    \pare{\Atw,\background\middle|\xi}
    {\bm J}_{\rm d}
    +
    {\bm A}^{\bC,{\rm lin}}_{\rm lo}
    \pare{\Atw,\background\middle|\xi}.
\end{equation}
On the low-frequency core where \(\chi\pare{\xi}=1\), one has
\begin{equation}\label{eq:linear_low_frequency_defect_core_complex_3d}
    {\bm A}^{\bC,{\rm lin}}_{\rm lo}
    \pare{\Atw,\background\middle|\xi}
    =
    \frac{\ii}{2}
    \begin{bmatrix}
        \mathfrak q\pare{\Atw,\background\middle|\xi}
        -\dfrac{\av{\xi}}{2}
        &
        \dfrac{\av{\xi}}{2}
        +\mathfrak q\pare{\Atw,\background\middle|\xi}
        \vspace{0.1cm}
        \\
        -\pare{
            \dfrac{\av{\xi}}{2}
            +\mathfrak q\pare{\Atw,\background\middle|\xi}
        }
        &
        \dfrac{\av{\xi}}{2}
        -\mathfrak q\pare{\Atw,\background\middle|\xi}
    \end{bmatrix}.
\end{equation}
Consequently, under the change of variables in
\cref{eq:complex_variable_definition_3d}, the linearized equation
\cref{eq:linearized_CDE_trivial_state} is equivalent to
\begin{equation}\label{eq:complex_linear_evolution_3d}
\begin{aligned}
    \partial_t U_{\rm lin}
    ={}&
    \ii\frac{\Atw}{4}
    \pare{\background^{\perp}\cdot D}U_{\rm lin}
    +
    \OpBW{
        \pare{1-\chi\pare{\xi}}
        \ii\omega_{\frac32}
        \pare{\Atw,\background\middle|\xi}
        {\bm J}_{\rm d}
    }U_{\rm lin}
    \\
    &+
    \OpBW{
        {\bm A}^{\bC,{\rm lin}}_{\rm lo}
        \pare{\Atw,\background\middle|\xi}
    }U_{\rm lin}.
\end{aligned}
\end{equation}
\end{lemma}

\begin{proof}
By
\cref{eq:def_linearized_discriminant,eq:complex_frequency_order_three_halves},
\[
    \mathfrak d\pare{\Atw,\background\middle|\xi}
    =
    \frac{\av{\xi}}{2}
    \mathfrak q\pare{\Atw,\background\middle|\xi},
    \qquad
    \omega_{\frac32}^2\pare{\Atw,\background\middle|\xi}
    =
    -\frac{\av{\xi}}{2}
    \mathfrak q\pare{\Atw,\background\middle|\xi}.
\]
Moreover, \cref{eq:complex_multiplier_p_equivalent} gives
\[
    \mathfrak q\mathfrak m_{-\frac14}^2
    =
    -\omega_{\frac32},
    \qquad
    \frac{\av{\xi}}{2}\mathfrak m_{-\frac14}^{-2}
    =
    \omega_{\frac32},
\]
where all multipliers are evaluated at
\(\pare{\Atw,\background\middle|\xi}\).  On the high-frequency region
where \(\chi\pare{\xi}=0\), one has
\(\mathfrak m^{\rm nh}_{-\frac14}=\mathfrak m_{-\frac14}\); hence ordinary
matrix multiplication yields
\begin{equation}\label{eq:linear_high_frequency_exact_diagonalization_3d}
    {\bm M}
    \begin{bmatrix}
        0&\frac{\av{\xi}}{2}
        \\
        \mathfrak q&0
    \end{bmatrix}
    {\bm M}^{-1}
    =
    \ii\omega_{\frac32}{\bm J}_{\rm d},
    \qquad
    \textnormal{on }\supp^\complement\pare{\chi}.
\end{equation}
The transport part
\(
    \ii\frac{\Atw}{4}
    \pare{\background^\perp\cdot\xi}\Id_{\bC^2}
\)
is scalar and therefore commutes with \({\bm M}\).  Combining
\cref{eq:linearized_CDE_matrix_symbol,eq:linear_low_frequency_defect_complex_3d,eq:linear_high_frequency_exact_diagonalization_3d}
proves \cref{eq:linear_diagonalization_identity_complex_3d}.  The defect
vanishes where \(\chi=0\), and is therefore supported in a fixed
low-frequency ball.  Since the symbolic calculus is restricted to
\(\av{\xi}\geq1/2\), it belongs to
\(\pare{\tilde\Gamma^0}^{2\times2}\).

On the core where \(\chi=1\), one has
\(\mathfrak m^{\rm nh}_{-\frac14}=1\).  Direct multiplication of the two
constant matrices in
\cref{eq:complex_matrix_M,eq:complex_matrix_M_inverse} gives
\cref{eq:linear_low_frequency_defect_core_complex_3d}.  Finally,
\cref{eq:complex_linear_evolution_3d} follows from
\cref{eq:complex_variable_definition_3d,eq:linearized_CDE_trivial_state,eq:linear_diagonalization_identity_complex_3d}.
\end{proof}

\begin{proof}[Proof of \cref{prop:complex_paradifferential_system_3d}]
Throughout the proof the physical variables
\(\eta\pare{U}\), \(\tilde\psi\pare{U}\) are those recovered from \(U\)
through
\cref{eq:eta_psitilde_as_functions_of_U,eq:physical_inverse_alinhac_components_complex_3d}.
We use the notation of \cref{notation:paralinearization}: symbols
\({\bm A}_{\bra{m,\geq1}}\) and smoothing remainders
\({\bm R}_{\bra{\varrho,\geq1}}\) are generic and may change from line to line.

\begin{proofpart}[Real order splitting]
\label{def:explicit_real_order_splitting_complex_3d}
We first split the nonlinear increments of the real off-diagonal entries
according to their symbolic order.

The homogeneous order-one Kelvin--Helmholtz shear symbol is
\begin{equation}\label{eq:complex_linear_shear_symbol}
    \dot q^{\rm sh}_{1}
    \pare{\Atw,\background\middle|\xi}
    \defeq
    \frac{4-\Atw^2}{8}
    \frac{\pare{\background^\perp\cdot\xi}^2}{\av{\xi}}
    \in\tilde\Gamma^1.
\end{equation}
The real order-one and order-zero increments in the upper-right entry are
\begin{align}
    \delta\dot q_{12,1}
    \pare{\Atw,\background\middle|U;x,\xi}
    &\defeq
    \lambda_1\pare{\eta\pare{U};x,\xi}
    -\frac{\av{\xi}}{2}
    \in\Gamma_{\geq1}^1,
    \label{eq:dq12_order_one_complex}
    \\
    \delta\dot q_{12,0}
    \pare{\Atw,\background\middle|U;x,\xi}
    &\defeq
    \Atw\lambda_0\pare{\eta\pare{U};x,\xi}
    \in\Gamma_{\geq1}^0,
    \label{eq:dq12_order_zero_complex}
\end{align}
where \(\lambda_1\) and \(\lambda_0\) are defined in
\cref{eq:lambda1,eq:lambda0}.  The real order-two and order-one increments in
the lower-left entry are
\begin{align}
    \delta\dot q_{21,2}
    \pare{\Atw,\background\middle|U;x,\xi}
    &\defeq
    a^{\pi_0,\eta}_2
    \pare{\rho^\pm\middle|\eta\pare{U};x,\xi}
    +
    \frac{\sigma}{\aangle{\rho^\pm}}\av{\xi}^2
    \in\Gamma_{\geq1}^2,
    \label{eq:dq21_order_two_complex}
    \\
    \delta\dot q_{21,1}
    \pare{\Atw,\background\middle|U;x,\xi}
    &\defeq
    \lambda_1\pare{\eta\pare{U};x,\xi}
    q^{\tilde\pi}_0
    \pare{
        \Atw,\background
        \middle|
        \eta\pare{U},\tilde\psi\pare{U};x,\xi
    }
    -
    \dot q^{\rm sh}_{1}
    \pare{\Atw,\background\middle|\xi}
    \in\Gamma_{\geq1}^1.
    \label{eq:dq21_order_one_complex}
\end{align}
Here \(a^{\pi_0,\eta}_2\) and \(q^{\tilde\pi}_0\) are defined in
\cref{eq:Api_star_2,eq:qpitilde_main}.

With this notation, the real off-diagonal quasilinear block in
\cref{eq:mathpzcQ} decomposes as
\begin{equation}\label{eq:dq_recomposition_complex}
    {\bm Q}\pare{\Atw,\background\middle|U;x,\xi}
    =
    \begin{bmatrix}
        0
        &
        \dfrac{\av{\xi}}{2}
        +\delta\dot q_{12,1}
        +\delta\dot q_{12,0}
        \\
        \mathfrak q
        +\delta\dot q_{21,2}
        +\delta\dot q_{21,1}
        &
        0
    \end{bmatrix}
    \pare{\Atw,\background\middle|U;x,\xi}.
\end{equation}
The flat entries \(\av{\xi}/2\) and
\(\mathfrak q\pare{\Atw,\background\middle|\xi}\) are defined in
\cref{eq:def_linearized_q_background}.
\end{proofpart}

\begin{proofpart}[Conjugation of the equation]
By
\cref{eq:complex_variable_definition_3d,eq:alinhac_variables_as_functions_of_U},
\[
    U
    =
    {\bm M}\pare{\Atw,\background\middle|D}Y_{\rm Al},
    \qquad
    Y_{\rm Al}
    =
    {\bm M}\pare{\Atw,\background\middle|D}^{-1}U .
\]
Since \({\bm M}\pare{\Atw,\background\middle|D}\) is a Fourier multiplier
independent of time,
\begin{equation}\label{eq:proof_complex_time_derivative}
    \partial_t U
    =
    {\bm M}\pare{\Atw,\background\middle|D}
    \partial_tY_{\rm Al}.
\end{equation}
Using
\cref{Alinhac:final},
the Alinhac-reduced real equation can be written as
\begin{equation}\label{eq:proof_real_alinhac_system_compact_complex}
    \partial_tY_{\rm Al}
    =
    \OpBW{
        {\bm Q}
        +
        {\bm T}
        +
        {\bm A}^{\bR,{\rm Al}}_{\rm lo}
    }Y_{\rm Al}
    +
    {\bm R}_{\bra{\varrho,\geq1}}
    Y_{\rm Al}.
\end{equation}
Substituting
\cref{eq:proof_real_alinhac_system_compact_complex} into
\cref{eq:proof_complex_time_derivative} gives
\begin{equation}\label{eq:proof_complex_conjugation_before_symbolic_product}
\begin{aligned}
    \partial_t U
    ={}&
    {\bm M}
    \OpBW{
        {\bm Q}
        +
        {\bm T}
        +
        {\bm A}^{\bR,{\rm Al}}_{\rm lo}
    }
    {\bm M}^{-1}U+
    {\bm M}
    {\bm R}_{\bra{\varrho,\geq1}}
    {\bm M}^{-1}U .
\end{aligned}
\end{equation}
We study each term separately.
\end{proofpart}
\begin{proofpart}[Paradifferential formula for the complex conjugation]
Let \({\bm A}\pare{U;x,\xi}\) be any matrix symbol whose structure will be specified below. We define
\begin{equation}\label{eq:complex_symbolic_conjugation_formula_3d}
    {\bm C}_{\bm M}\bra{{\bm A}}
    \pare{\Atw,\background\middle|U;x,\xi}
    \defeq
    {\bm M}\pare{\Atw,\background\middle|\xi}
    \#_\varrho
    {\bm A}\pare{\Atw,\background\middle|U;x,\xi}
    \#_\varrho
    {\bm M}\pare{\Atw,\background\middle|\xi}^{-1}.
\end{equation}
Note that if ${\bm A}={\bm A}\pare{\Atw,\background\middle|\xi}$ is a Fourier multiplier independent of $U$ and $x$, then the relation \cref{eq:complex_symbolic_conjugation_formula_3d} simplifies into the classical conjugation
\begin{equation}\label{simplif:FMconj}
    {\bm C}_{\bm M}\bra{{\bm A}}\pare{\Atw,\background\middle|\xi}={\bm M}\pare{\Atw,\background\middle|\xi}{\bm A}\pare{\Atw,\background\middle|\xi}{\bm M}\pare{\Atw,\background\middle|\xi}^{-1}
\end{equation}
and
\begin{equation}\label{Fourier:comp}
    {\bm M}
    \OpBW{{\bm A}}
    {\bm M}^{-1}=\OpBW{{\bm M}{\bm A}{\bm M}^{-1}}.
\end{equation}
In
\cref{Fourier:comp},
all matrices \({\bm M}^{\pm1}\) are evaluated at \(\pare{\Atw,\background\middle|D}\) when considered as the conjugating operators and at
\(\pare{\Atw,\background\middle|\xi}\) when considered as symbols. Since there is no possible confusion, in what follows, we alleviate the notations removing this dependence.
Otherwise, if the matrix symbol ${\bm A}$ belongs to $\pare{\Gamma_{\geq 1}^{m}}^{2\times 2}$, then, by the Bony--Weyl composition formula
\cref{def:asympt_symbol,prop:composition_paradiff} and by the composition
rules in \cref{lem:composition_rules},
\begin{equation}\label{eq:complex_operator_conjugation_formula_3d}
    {\bm M}
    \OpBW{{\bm A}}
    {\bm M}^{-1}
    =
    \OpBW{
        {\bm C}_{\bm M}\bra{{\bm A}}
    }
    +
    {\bm R}_{\bra{\varrho,\geq1}}.
\end{equation}
Indeed,
\begin{align}
    {\bm M}
    \OpBW{{\bm A}}
    &=
    \OpBW{
        {\bm M}\#_\varrho{\bm A}
    }
    +
    {\bm R}_{\bra{\varrho,\geq1}},
    \label{eq:complex_MA_composition_3d}
    \\
    \OpBW{
        {\bm M}\#_\varrho{\bm A}
    }
    {\bm M}^{-1}
    &=
    \OpBW{
        \pare{{\bm M}\#_\varrho{\bm A}}
        \#_\varrho
        {\bm M}^{-1}
    }
    +
    {\bm R}_{\bra{\varrho,\geq1}}.
    \label{eq:complex_MAMinv_composition_3d}
\end{align}
The smoothing remainder produced in
the first composition remains smoothing after composition with the Fourier
multiplier \({\bm M}^{-1}\), after decreasing
\(\varrho\) by a fixed amount.
\end{proofpart}

\begin{proofpart}[Conjugation of the linear Fourier multiplier and transport component]
By \cref{eq:dq_recomposition_complex},
\[
    {\bm Q}
    =
    \begin{bmatrix}
        0&\frac{\av{\xi}}{2}
        \\
        \mathfrak q&0
    \end{bmatrix}
    +
    \begin{bmatrix}
        0&\delta\dot q_{12,1}
        \\
        \delta\dot q_{21,2}&0
    \end{bmatrix}
    +
    \begin{bmatrix}
        0&\delta\dot q_{12,0}
        \\
        \delta\dot q_{21,1}&0
    \end{bmatrix}.
\]
Moreover,
\begin{equation}\label{eq:complex_real_transport_symbol_3d}
    {\bm T}
    =
    \ii
    \pare{
        \mathpzc V_{\textnormal{can}}
        \pare{
            \Atw,\background
            \middle|
            \eta\pare{U},\tilde\psi\pare{U};x
        }
        \cdot\xi
    }
    \Id_{\bC^2}.
\end{equation}
By
\cref{eq:Vpitilde_main,tot-can-vort_at_zero,tildeV and dv3 at zero}, define
\begin{equation}\label{eq:nonlinear_transport_remainder_complex_3d}
    \mathpzc V_{\geq1}
    \pare{\Atw,\background\middle|U;x}
    \defeq
    \mathpzc V_{\textnormal{can}}
    \pare{
        \Atw,\background
        \middle|
        \eta\pare{U},\tilde\psi\pare{U};x
    }
    -
    \frac{\Atw}{4}\background^\perp
    \in\pare{\Gamma_{\geq1}^0}^{2}.
\end{equation}
Consequently,
\[
    {\bm T}
    =
    \ii\frac{\Atw}{4}
    \pare{\background^\perp\cdot\xi}\Id_{\bC^2}
    +
    \ii\pare{\mathpzc V_{\geq1}\cdot\xi}\Id_{\bC^2},
\]
and, by \cref{eq:linearized_CDE_matrix_symbol},
\[
    {\bm Q}+{\bm T}
    =
    {\bm L}_0
    +
    \begin{bmatrix}
        0&\delta\dot q_{12,1}
        \\
        \delta\dot q_{21,2}&0
    \end{bmatrix}
    +
    \begin{bmatrix}
        0&\delta\dot q_{12,0}
        \\
        \delta\dot q_{21,1}&0
    \end{bmatrix}
    +
    \ii\pare{\mathpzc V_{\geq1}\cdot\xi}\Id_{\bC^2}.
\]
Applying
\cref{prop:linear_diagonalization_complex_3d,Fourier:comp} gives
\[
    {\bm M}\OpBW{{\bm L}_{0}\pare{\xi}}{\bm M}^{-1}
    =
    \OpBW{
        \ii\frac{\Atw}{4}
        \pare{\background^\perp\cdot\xi}\Id_{\bC^2}
        +
        \pare{1-\chi\pare{\xi}}
        \ii\omega_{\frac32}{\bm J}_{\rm d}
        +
        {\bm A}^{\bC,{\rm lin}}_{\rm lo}
    }.
\]
For the remaining transport component, apply
\cref{eq:complex_operator_conjugation_formula_3d,eq:complex_symbolic_conjugation_formula_3d}
with
\(
    {\bm A}
    =
    \ii\pare{\mathpzc V_{\geq1}\cdot\xi}\Id_{\bC^2}
\).
Its pointwise conjugation is unchanged because the symbol is scalar in the
matrix indices.  Every non-pointwise Bony--Weyl correction is of order zero
and at least one-homogeneous in \(U\).  Therefore
\begin{equation}\label{eq:transport_complex_conjugation_expanded}
    {\bm C}_{\bm M}
    \bra{
        \ii\pare{\mathpzc V_{\geq1}\cdot\xi}\Id_{\bC^2}
    }
    =
    \ii\pare{\mathpzc V_{\geq1}\cdot\xi}\Id_{\bC^2}
    +
    {\bm A}_{\bra{0,\geq1}}.
\end{equation}
Recombining the flat and nonlinear transport contributions yields
\begin{equation}\label{semi_complex_conj}
\begin{aligned}
    {\bm M}\OpBW{{\bm Q}+{\bm T}}{\bm M}^{-1}
    ={}&
    \OpBW{
        \pare{1-\chi\pare{\xi}}
        \ii\omega_{\frac32}{\bm J}_{\rm d}
        +
        {\bm A}^{\bC,{\rm lin}}_{\rm lo}
        +
        {\bm T}
        +
        {\bm A}_{\bra{0,\geq1}}
    }
    +
    {\bm R}_{\bra{\varrho,\geq1}}
    \\
    &+
    \OpBW{
        {\bm C}_{\bm M}
        \bra{
            \begin{bmatrix}
                0&\delta\dot q_{12,1}
                \\
                \delta\dot q_{21,2}&0
            \end{bmatrix}
            +
            \begin{bmatrix}
                0&\delta\dot q_{12,0}
                \\
                \delta\dot q_{21,1}&0
            \end{bmatrix}
        }
    }.
\end{aligned}
\end{equation}
It remains to analyze the last term in \cref{semi_complex_conj}.
\end{proofpart}

\begin{proofpart}[The high-frequency quasilinear block and the explicit order-\(\frac12\) block]
All positive-order terms are extracted on the high-frequency region
\(\supp^\complement\pare{\chi}\). On this region
\[
    \mathfrak m^{\rm nh}_{-\frac14}
    \pare{\Atw,\background\middle|\xi}
    =
    \mathfrak m_{-\frac14}
    \pare{\Atw,\background\middle|\xi}.
\]
The complementary low-frequency contribution is compactly supported in
\(\xi\), hence it corresponds to bounded operators. Moreover, since $\delta\dot{q}_{12,1}$, $\delta\dot{q}_{12,0}$, $\delta\dot{q}_{21,2}$ and $\delta\dot{q}_{21,1}$ are at least $1$-homogeneous in $U,$ such contributions are included in the generic bounded symbol
\({\bm A}_{\bra{0,\geq1}}\). In what follows, we shall only work within this high-frequency region. Given two
scalar symbols \(a,b\), we define
\[
    {\mathscr B}\pare{a,b}
    \defeq
    \begin{bmatrix}
        b-a&a+b
        \\
        -\pare{a+b}&-\pare{b-a}
    \end{bmatrix}.
\]
With this notation in hand, one readily has on $\supp^{\complement}\pare{\chi},$
\begin{equation}
    {\bm C}_{{\bm M}}\bra{\begin{bmatrix}
        0&a\\ b&0
    \end{bmatrix}}=\frac{\ii}{2}\mathscr{B}\pare{\mathfrak m_{-\frac14}^{-1}\#_{\varrho}a\#_{\varrho}\mathfrak m_{-\frac14}^{-1}\,,\, \mathfrak m_{-\frac14}\#_{\varrho}b\#_{\varrho}\mathfrak m_{-\frac14}}.
\end{equation}
Now, observe the following Weyl cancellation: If
\(\mathfrak m=\mathfrak{m}\pare{\xi}\in\tilde{\Gamma}^{\mu}\) is a scalar Fourier multiplier, then for any $c\in\Gamma_{\geq1}^m$, we have 
\begin{equation}\label{eq:same_multiplier_weyl_cancellation_complex}
    \mathfrak m\#_\varrho c\#_\varrho\mathfrak m
    =
    \mathfrak m^2c
    +
    \Gamma_{\geq1}^{2\mu+m-2}.
\end{equation}
Indeed, the first-order Weyl corrections cancel because the same multiplier
appears on both sides of \(c\). Applying \cref{eq:same_multiplier_weyl_cancellation_complex} with either the couple $\pare{\mathfrak m,c}=\pare{\mathfrak{m}_{-\frac14}^{-1},a}\in\tilde{\Gamma}^{\frac14}\times\Gamma_{\geq1}^{1}$ or a couple $\pare{\mathfrak m,c}=\pare{\mathfrak{m}_{-\frac14},b}\in\tilde{\Gamma}^{-\frac14}\times\Gamma_{\geq1}^{2}$, we find on 
\(\supp^\complement\pare{\chi}\),
\begin{equation}\label{eq:off_diagonal_high_frequency_conjugation_formula_complex_intermediaire}
    {\bm C}_{{\bm M}}\bra{\begin{bmatrix}
        0&a\\ b&0
    \end{bmatrix}}=
    \frac{\ii}{2}
    {\mathscr B}\pare{
        \mathfrak m_{-\frac14}^{-2}a,
        \mathfrak m_{-\frac14}^{2}b
    } +{\bm A}_{\bra{0,\geq1}}.
\end{equation}
According to the explicit formulae \cref{eq:complex_multiplier_p,eq:complex_multiplier_p_equivalent,eq:complex_leading_capillary_frequency_3d,eq:complex_leading_multiplier_p_3d}, one can check that
\begin{equation}\label{eq:frakm_nh_to_h}
    \mathfrak m_{-\frac14}^{-2}
    =
    \frac{
        \omega_{\frac32}
    }{
        \dot\omega_{\frac32}
    }
    \dot{\mathfrak m}_{-\frac14}^{-2}
\end{equation}
and
\begin{equation}\label{eq:p_square_half_order_defect_3d}
    \mathfrak m_{-\frac14}^{2}
    -
    \frac{
        \omega_{\frac32}
    }{
        \dot\omega_{\frac32}
    }
    \dot{\mathfrak m}_{-\frac14}^{2}
    =
    \dot{\mathfrak c}^{\rm sh}_{-\frac32}
    +
    \tilde\Gamma^{-\frac52},
\end{equation}
with \(\dot{\mathfrak c}^{\rm sh}_{-\frac32}\) defined by
\begin{equation}\label{eq:complex_shear_half_defect_3d}
        \dot{\mathfrak c}^{\rm sh}_{-\frac32}
        \pare{\Atw,\background\middle|\xi}
        \defeq
        \frac{4-\Atw^2}{32}
        \pare{
            \frac{2\aangle{\rho^\pm}}{\sigma}
        }^{\frac32}
        \frac{
            \pare{\background^\perp\cdot\xi}^2
        }{
            \av{\xi}^{\frac72}
        }
        \in\tilde\Gamma^{-\frac32}.
\end{equation}
Inserting \cref{eq:frakm_nh_to_h,eq:p_square_half_order_defect_3d} into \cref{eq:off_diagonal_high_frequency_conjugation_formula_complex_intermediaire}, we infer for any $a\in\Gamma_{\geq1}^1$ and $b\in\Gamma_{\geq1}^2$, on $\supp^{\complement}\pare{\chi},$
\begin{equation}\label{eq:off_diagonal_high_frequency_conjugation_formula_complex}
    {\bm C}_{{\bm M}}\bra{\begin{bmatrix}
        0&a\\ b&0
    \end{bmatrix}}=
    \frac{\ii}{2}\frac{
        \omega_{\frac32}
    }{
        \dot\omega_{\frac32}
    }
    {\mathscr B}\pare{
        \dot{\mathfrak m}_{-\frac14}^{-2}a,
        \dot{\mathfrak m}_{-\frac14}^{2}b
    }+\frac{\ii}{2}\mathscr{B}\pare{0,\dot{\mathfrak c}^{\rm sh}_{-\frac32}b} +{\bm A}_{\bra{0,\geq1}}.
\end{equation}
We now extract the quasilinear order-\(3/2\) correction. Applying
\cref{eq:off_diagonal_high_frequency_conjugation_formula_complex} with the pair
\[
    \pare{a,b}=\pare{
        \delta\dot q_{12,1},
        \delta\dot q_{21,2}
    }\in\Gamma_{\geq1}^{1}\times\Gamma_{\geq1}^2,
\]
we obtain by definition \cref{eq:complex_ql_matrix_coefficient_3d} that on $\supp^{\complement}\pare{\chi},$
\begin{equation}\label{eq:principal_high_frequency_ql_extraction_complex_3d}
\begin{aligned}
    {\bm C}_{\bm M}
    \bra{
        \begin{bmatrix}
            0&\delta\dot q_{12,1}\\
            \delta\dot q_{21,2}&0
        \end{bmatrix}
    }=
    \ii\omega_{\frac32}
    {\bm J}_{\rm d}
    \dot{\bm E}^{\bC}_{{\rm ql},0}+
    \frac{\ii}{2}
    {\mathscr B}\pare{
        0,
        \dot{\mathfrak c}^{\rm sh}_{-\frac32}
        \delta\dot q_{21,2}
    }
    +
    {\bm A}_{\bra{0,\geq1}},
\end{aligned}
\end{equation}
where
$$\dot{\bm E}^{\bC}_{{\rm ql},0}\defeq\begin{bmatrix}
            \dot e_0 & \dot f_0
            \\
            \dot f_0 & \dot e_0
        \end{bmatrix},$$
    with
    \begin{align}
    \dot e_0
        \defeq
        -
        \frac{
            \dot{\mathfrak m}_{-\frac14}^{2}
            \delta\dot q_{21,2}
            -
            \dot{\mathfrak m}_{-\frac14}^{-2}
            \delta\dot q_{12,1}
        }{
            2\dot\omega_{\frac32}
        },
        \label{eq:dote0}\\
        \dot f_0
        \defeq
        -
        \frac{
            \dot{\mathfrak m}_{-\frac14}^{-2}
            \delta\dot q_{12,1}
            +
            \dot{\mathfrak m}_{-\frac14}^{2}
            \delta\dot q_{21,2}
        }{
            2\dot\omega_{\frac32}
        }\cdot\label{eq:dotf0}
    \end{align}
By
\cref{eq:lambda1,eq:symbol_a2_pi0star,eq:dq12_order_one_complex,eq:dq21_order_two_complex,eq:complex_leading_capillary_frequency_3d,eq:complex_leading_multiplier_p_3d}, we compute explicitly
\begin{align}
    \dot{\mathfrak m}_{-\frac14}^{-2}
    \delta\dot q_{12,1}
    &=
    \dot\omega_{\frac32}
    \pare{
        \sqrt{
            1+\av{\nabla^\perp\eta\cdot\hat{\xi}}^2
        }
        -1
    },
    \label{eq:dotm_dq12_closed_principal_3d}
    \\
    \dot{\mathfrak m}_{-\frac14}^{2}
    \delta\dot q_{21,2}
    &=
    \dot\omega_{\frac32}
    \pare{
        1
        -
        \frac{
            1+\av{\nabla^\perp\eta\cdot\hat{\xi}}^2
        }{
            \pare{1+\av{\nabla\eta}^2}^{\frac32}
        }
    }.
    \label{eq:dotm_dq21_closed_principal_3d}
\end{align}
Inserting \cref{eq:dotm_dq12_closed_principal_3d,eq:dotm_dq21_closed_principal_3d} into \cref{eq:dote0,eq:dotf0}, we get \cref{eq:dote0_closed_principal_3d,eq:dotf0_closed_principal_3d}. Similarly, one can apply \cref{eq:off_diagonal_high_frequency_conjugation_formula_complex}to the pair
\[
    \pare{a,b}=\pare{
        \delta\dot q_{12,0},
        \delta\dot q_{21,1}}\in\Gamma_{\geq1}^0\times\Gamma_{\geq1}^1
\]
and get on $\supp^{\complement}\pare{\chi},$
\begin{equation}\label{eq:complex_lower_pair_conjugation_3d}
\begin{aligned}
    {\bm C}_{\bm M}
    \bra{
        \begin{bmatrix}
            0&\delta\dot q_{12,0}\\
            \delta\dot q_{21,1}&0
        \end{bmatrix}
    }=\frac{\ii}{2}
    {\mathscr B}\pare{
        \dot{\mathfrak m}_{-\frac14}^{-2}
        \delta\dot q_{12,0},
        \dot{\mathfrak m}_{-\frac14}^{2}
        \delta\dot q_{21,1}
    }
    +
    {\bm A}_{\bra{0,\geq1}} .
\end{aligned}
\end{equation}
Indeed, we obtain from \cref{eq:complex_leading_multiplier_expansions_3d} that $\frac{
        \omega_{\frac32}
    }{
        \dot\omega_{\frac32}
    }-1\in\tilde{\Gamma}^{-1}$, so the remaining terms create only bounded symbols. Now observe that the cut-off function being with compact support, for any $m\in\mathbb{R}$ and any ${\bm A},{\bm B}\in\pare{\Gamma_{\geq 1}^{m}}^{2\times2}$,
    \begin{equation}\label{support_to_cutoff}
        \pare{{\bm A}={\bm B}\quad\textnormal{on}\quad\supp^{\complement}\pare{\chi}}\qquad\Rightarrow\qquad{\bm A}=\pare{1-\chi\pare{\xi}}{\bm B}+{\bm A}_{\bra{0,\geq1}}.
    \end{equation}
	    Putting together \cref{semi_complex_conj,eq:principal_high_frequency_ql_extraction_complex_3d,eq:complex_lower_pair_conjugation_3d,support_to_cutoff}, we have obtained
	    \begin{equation}\label{complex_conj_QplusT}
	    \begin{aligned}
	        {\bm M}\OpBW{{\bm Q}+{\bm T}}{\bm M}^{-1}&={}\OpBW{\pare{1-\chi\pare{\xi}}
    \ii
    \omega_{\frac32}
    {\bm J}_{\rm d}\pare{\Id+\dot{\bm E}^{\bC}_{{\rm ql},0}}+{\bm A}^{\bC,{\rm lin}}_{\rm lo}+{\bm T}+\pare{1-\chi\pare{\xi}}{\bm A}_{\frac{1}{2}}^{\bC}+{\bm A}_{\bra{0,\geq1}}}+{\bm R}_{\bra{\varrho,\geq1}},
	    \end{aligned}
	    \end{equation}
where
\begin{equation}\label{def:bmA-half}
    {\bm A}_{\frac{1}{2}}^{\bC}\defeq\frac{\ii}{2}\pare{{\mathscr B}\pare{
        0,
        \dot{\mathfrak c}^{\rm sh}_{-\frac32}
        \delta\dot q_{21,2}
    }+{\mathscr B}\pare{
        \dot{\mathfrak m}_{-\frac14}^{-2}
        \delta\dot q_{12,0},
        \dot{\mathfrak m}_{-\frac14}^{2}
        \delta\dot q_{21,1}
    }}.
\end{equation}
We are left with the explicit computation of the symbol ${\bm A}_{\frac{1}{2}}^{\bC}$ in \cref{def:bmA-half}. Note that we can rewrite
$${\bm A}_{\frac{1}{2}}^{\bC}=\dot{\bm A}^{\bC}_{0}\av{\xi}^{\frac{1}{2}},\qquad \dot{\bm A}^{\bC}_{0}\triangleq\frac{\ii}{2}\mathscr{B}\pare{\dot{\mathfrak a}^{\bC}_{0},
        \dot{\mathfrak b}^{\bC}_{0}},$$
where
\begin{align}
    \dot{\mathfrak a}^{\bC}_{0}&\defeq\av{\xi}^{-\frac12}
    \dot{\mathfrak m}_{-\frac14}^{-2}
    \delta\dot q_{12,0},\\
    \dot{\mathfrak b}^{\bC}_{0}&\defeq\av{\xi}^{-\frac12}
    \bra{
        \dot{\mathfrak m}_{-\frac14}^{2}
        \delta\dot q_{21,1}
        +
        \dot{\mathfrak c}^{\rm sh}_{-\frac32}\delta\dot q_{21,2}}.
\end{align}
Combining \cref{eq:complex_leading_capillary_frequency_3d,eq:complex_leading_multiplier_p_3d,eq:dq12_order_zero_complex}, we immediately obtain the explicit expression \cref{eq:complex_a_half_scalar_3d} for $\dot{\mathfrak a}^{\bC}_{0}$. The explicit computation of $\dot{\mathfrak b}^{\bC}_{0}$ requires more explanations. We compute the two contributions separately. First,
\cref{eq:dq21_order_one_complex,eq:lambda1,eq:qpitilde_main,eq:complex_linear_shear_symbol}
give
\begin{align}
    \delta\dot q_{21,1}
    \pare{\Atw,\background\middle|U;x,\xi}
    &=
    \frac{4-\Atw^2}{2}
    \av{\xi}{
    \sqrt{
        1+
        \av{
            \nabla^\perp\eta\pare{U}
            \cdot\hat{\xi}
        }^2
    }
    \pare{
        A_{0}^{\mathpzc{B}_0,\omega;\perp}
        \pare{\eta\pare{U};x,\xi}
        \cdot
        \Omega
        \pare{
            \background
            \middle|
            \eta\pare{U},\tilde\psi\pare{U};
            x
        }
    }^2}
    \nonumber
    \\
    &\quad
    -
    \frac{4-\Atw^2}{8}
    \av{\xi}
    \pare{\background^\perp\cdot\hat{\xi}}^2 .
    \label{eq:dq21_order_one_expanded_half_proof_3d}
\end{align}
Moreover,
\[
    \av{\xi}^{-\frac12}
    \dot{\mathfrak m}_{-\frac14}^{2}
    \pare{\Atw,\background\middle|\xi}
    =
    \sqrt{
        \frac{\aangle{\rho^\pm}}{2\sigma}
    }\,
    \av{\xi}^{-1}.
\]
Therefore
\begin{align}
    &\av{\xi}^{-\frac12}
    \dot{\mathfrak m}_{-\frac14}^{2}
    \pare{\Atw,\background\middle|\xi}
    \delta\dot q_{21,1}
    \pare{\Atw,\background\middle|U;x,\xi}
    \nonumber
    \\
    &\quad =
    \frac{4-\Atw^2}{16}
    \sqrt{
        \frac{2\aangle{\rho^\pm}}{\sigma}
    }
    \bra{
        4
        \sqrt{
            1+
            \av{
                \nabla^\perp\eta\pare{U}
                \cdot\hat{\xi}
            }^2
        }
        \pare{
            A_{0}^{\mathpzc{B}_0,\omega;\perp}
            \pare{\eta\pare{U};x,\xi}
            \cdot
            \Omega
            \pare{
                \background
                \middle|
                \eta\pare{U},\tilde\psi\pare{U};
                x
            }
        }^2
        -
        \pare{\background^\perp\cdot\hat{\xi}}^2
    } .
    \label{eq:bhalf_first_contribution_3d}
\end{align}
Second, by \cref{eq:dq21_order_two_complex,eq:Api_star_2},
\[
    \delta\dot q_{21,2}
    \pare{\Atw,\background\middle|U;x,\xi}
    =
    \frac{\sigma}{\aangle{\rho^\pm}}
    \av{\xi}^2
    \bra{
        1
        -
        \frac{
            1+
            \av{
                \nabla^\perp\eta\pare{U}
                \cdot\hat{\xi}
            }^2
        }{
            \pare{
                1+\av{\nabla\eta\pare{U}}^2
            }^{\frac32}
        }
    } .
\]
Using \cref{eq:complex_shear_half_defect_3d}, this gives
\begin{equation}
\label{eq:bhalf_second_contribution_3d}
\av{\xi}^{-\frac12}
    \dot{\mathfrak c}^{\rm sh}_{-\frac32}
    \pare{\Atw,\background\middle|\xi}
    \delta\dot q_{21,2}
    \pare{\Atw,\background\middle|U;x,\xi}
  =
    \frac{4-\Atw^2}{16}
    \sqrt{
        \frac{2\aangle{\rho^\pm}}{\sigma}
    }\,
    \pare{\background^\perp\cdot\hat{\xi}}^2
    \bra{
        1
        -
        \frac{
            1+
            \av{
                \nabla^\perp\eta\pare{U}
                \cdot\hat{\xi}
            }^2
        }{
            \pare{
                1+\av{\nabla\eta\pare{U}}^2
            }^{\frac32}
        }
    } .
\end{equation}
Adding
\cref{eq:bhalf_first_contribution_3d} and \cref{eq:bhalf_second_contribution_3d}
cancels the pure flat shear contribution and yields precisely
\cref{eq:complex_b_half_scalar_3d}.  The square on
\(A_{0}^{\mathpzc{B}_0,\omega;\perp}\cdot\Omega\) comes from \(q^{\tilde\pi}_0\) in
\cref{eq:qpitilde_main}.
\end{proofpart}

\begin{proofpart}[The bounded real lower-order block]
We apply
\cref{eq:complex_operator_conjugation_formula_3d,eq:complex_symbolic_conjugation_formula_3d}
with ${\bm A}={\bm A}^{\bR,{\rm Al}}_{\rm lo}.$
Then
\[
    {\bm M}
    \OpBW{{\bm A}^{\bR,{\rm Al}}_{\rm lo}}
    {\bm M}^{-1}
    =
    \OpBW{{\bm A}^{\bC,{\rm lo}}_{\rm conj}}
    +
    {\bm R}_{\bra{\varrho,\geq1}},
    \qquad
    {\bm A}^{\bC,{\rm lo}}_{\rm conj}
    \defeq
   {\bm C}_{\bm M}
    \bra{{\bm A}^{\bR,{\rm Al}}_{\rm lo}}.
\]
By \cref{eq:real_alinhac_lower_order_structure_3d},
\[
    {\bm A}^{\bR,{\rm Al}}_{\rm lo}
    \in
    \begin{bmatrix}
        \Gamma_{\geq1}^0&\Gamma_{\geq1}^{-1}\vspace{0.1cm}\\
        \Gamma_{\geq1}^0&\Gamma_{\geq1}^0
    \end{bmatrix}.
\]
Moreover, by \cref{eq:complex_matrix_M,eq:complex_matrix_M_inverse}, the
entries of \({\bm M}\) and \({\bm M}^{-1}\) have orders \(-1/4\) or
\(1/4\) at high frequency, and are order zero at low frequency.  Therefore
the largest order produced by the pointwise products is \(0\).  All
non-pointwise Bony--Weyl corrections contain at least one \(\xi\)-derivative
and therefore have no larger order. Consequently
\begin{equation}\label{eq:complex_lower_order_conjugation_3d}
    {\bm A}^{\bC,{\rm lo}}_{\rm conj}
    \in
    \pare{\Gamma_{\geq1}^0}^{2\times2}.
\end{equation}
This contribution is included in the generic bounded symbol
\({\bm A}_{\bra{0,\geq1}}\).
\end{proofpart}

\begin{proofpart}[Smoothing remainders and conclusion]
The smoothing remainder in
\cref{eq:proof_complex_conjugation_before_symbolic_product} is controlled by
\cref{lem:composition_rules}:
\begin{equation}\label{complex_conj_smoothing}
    {\bm M}
    {\bm R}_{\bra{\varrho,\geq1}}
    {\bm M}^{-1}
    \in
    \pare{\cR_{\geq1}^\varrho}^{2\times2},
\end{equation}
after decreasing and relabelling \(\varrho\) and can be renamed ${\bm R}_{\bra{\varrho,\geq1}}$. The same convention absorbs
all remainders generated by the symbolic compositions in
\cref{eq:complex_MA_composition_3d,eq:complex_MAMinv_composition_3d}. Combining
\cref{eq:proof_complex_conjugation_before_symbolic_product,complex_conj_QplusT,eq:complex_lower_order_conjugation_3d,complex_conj_smoothing}, we have obtained \cref{eq:complex_paradifferential_system_3d}.
\end{proofpart}
\end{proof}

\begin{remark}[Non-homogeneous Sobolev scale and the explicit order-\(1/2\) block]
\label{rem:complex_half_order_block_after_GU_3d}
By \cref{eq:complex_variable_definition_3d},
\begin{equation}\label{eq:complex_variable_scalar_definition_3d}
    u
    =
    \frac{1}{\sqrt2}
    \pare{
        \pare{
            \mathfrak m^{\rm nh}_{-\frac14}
            \pare{\Atw,\background\middle|D}
        }^{-1}
        \eta_{\rm Al}
        +
        \ii\mathfrak m^{\rm nh}_{-\frac14}
        \pare{\Atw,\background\middle|D}
        \tilde\psi_{\rm Al}
    }.
\end{equation}
Conversely,
\begin{equation}\label{eq:inverse_complex_variable_definition_3d}
    \eta_{\rm Al}
    =
    \frac{1}{\sqrt2}
    \mathfrak m^{\rm nh}_{-\frac14}
    \pare{\Atw,\background\middle|D}
    \pare{u+\bar u},
    \qquad
    \tilde\psi_{\rm Al}
    =
    -\frac{\ii}{\sqrt2}
    \pare{
        \mathfrak m^{\rm nh}_{-\frac14}
        \pare{\Atw,\background\middle|D}
    }^{-1}
    \pare{u-\bar u}.
\end{equation}

By \cref{eq:complex_nonhomogeneous_multiplier_p_3d}, the multiplier
\(\mathfrak m^{\rm nh}_{-\frac14}\) is equal to \(1\) at low frequency and
coincides with \(\mathfrak m_{-\frac14}\) at high frequency.  Hence the
complex variable in \cref{eq:complex_variable_definition_3d} satisfies the
non-homogeneous equivalence
\[
    \norm{U}_{H^s}
    \sim
    \norm{\eta_{\rm Al}}_{H^{s+\frac14}}
    +
    \norm{\tilde\psi_{\rm Al}}_{H^{s-\frac14}} .
\]
In the capillary high-frequency regime,
\[
    \omega_{\frac32}
    \sim
    \dot\omega_{\frac32}
    =
    \sqrt{
        \frac{\sigma}{2\aangle{\rho^\pm}}
    }
    \av{\xi}^{\frac32},
    \qquad
    \mathfrak m_{-\frac14}
    \sim
    \dot{\mathfrak m}_{-\frac14}
    \sim
    C_{\rho^\pm,\sigma}\av{\xi}^{-\frac14}.
\]
Therefore the explicit order-\(\frac12\) block
\[
    \pare{1-\chi\pare{\xi}}\,
    \dot{\bm A}^{\bC}_{0}
    \pare{\Atw,\background\middle|U;x,\hat\xi}
    \av{\xi}^{\frac12}
\]
is displayed only on the high-frequency region.  The complementary
low-frequency part is a bounded contribution and is included in
\({\bm A}_{\bra{0,\geq1}}\).  The Alinhac order-\(\frac12\) terms are absent because
the Alinhac block has already been removed in
\cref{prop:real_alinhac_good_unknown_3d}.  \fine
\end{remark}

\subsection{Diagonalization of the positive-order symbols}
\label{sec:positive_order_diagonalization_3d}

In this subsection we diagonalize the positive-order part of the complex
Alinhac-reduced system obtained in
\cref{prop:complex_paradifferential_system_3d}.  The diagonalization is
performed in two steps.  First we diagonalize the quasilinear order-\(3/2\)
matrix
\[
    {\bm J}_{\rm d}
    \pare{
        \Id_{\bC^2}
        +
        \dot{\bm E}^{\bC}_{{\rm ql},0}
    },
\]
by a matrix paradifferential operator of order \(0\).  Then we remove the
off-diagonal part of the explicit order-\(1/2\) block by a correction of
order \(-1\).  The scalar transport part remains diagonal at all positive
orders.

Throughout this section we use the conventions of
\cref{notation:vector_convention}: scalar symbols are denoted by lower-case
letters, two-component variables by upper-case letters, and \(2\times2\)
matrix symbols and matrix operators by upper-case bold letters.  Homogeneous
symbols in \(\xi\) are denoted with a dot, and the subscript records the
degree of homogeneity.

\begin{prop}[Diagonalization of the positive-order symbols]
\label{prop:positive_order_diagonalization_3d}
Let the assumptions of
\cref{prop:complex_paradifferential_system_3d} hold, and let \(U\) solve
\cref{eq:complex_paradifferential_system_3d}.  There exists a
bounded and invertible operator
\[
    {\bm B}
    \pare{\Atw,\background\middle|U}
    \in
    \cL\pare{ H^s}
    \qquad
    \forall\,s\geq s_0,
\]
such that, defining
\begin{equation}\label{eq:positive_order_diagonalized_unknown_3d}
    W
    \defeq
    {\bm B}
    \pare{\Atw,\background\middle|U}^{-1}U ,
\end{equation}
the unknown \(W\) solves
\begin{equation}\label{eq:positive_order_diagonalized_system_3d}
    \partial_t W
    =
    \OpBW{
        {\bm A}^{\bC,{\rm lin}}_{\rm lo}
        \pare{\Atw,\background\middle|\xi}
    }W
    +
    \ii
    \OpBW{
        {\bm D}^{\bC}_{>0}
        \pare{\Atw,\background\middle|U;x,\xi}
        +
        {\bm A}_{\bra{0,\geq1}}
        \pare{\Atw,\background\middle|U;x,\xi}
    }W
    +
    {\bm R}_{\bra{\varrho,\geq1}}
    \pare{\Atw,\background\middle|U}W .
\end{equation}
The objects entering
\cref{eq:positive_order_diagonalized_unknown_3d,eq:positive_order_diagonalized_system_3d}
are defined as follows.
\begin{itemize}

    \item
    \({\bm A}^{\bC,{\rm lin}}_{\rm lo}
    \pare{\Atw,\background\middle|\xi}\)
    is the same flat low-frequency linear multiplier introduced in
    \cref{eq:linear_low_frequency_defect_complex_3d};

    \item
    \({\bm D}^{\bC}_{>0}
    \pare{\Atw,\background\middle|U;x,\xi}
    \in\pare{\Gamma_{\geq0}^{\frac32}}^{2\times2}\)
    is the diagonal positive-order matrix of real symbols
    \begin{equation}\label{eq:positive_order_diagonal_symbol_decomposition_3d}
        {\bm D}^{\bC}_{>0}
        \defeq
        {\bm D}^{\bC}_{\frac32}
        +
        \dot{\bm D}^{\bC}_{1}
        +
        \dot{\bm D}^{\bC}_{\frac12}.
    \end{equation}
    The symbols
    \({\bm D}^{\bC}_{\frac32}\),
    \(\dot{\bm D}^{\bC}_{1}\), and
    \(\dot{\bm D}^{\bC}_{\frac12}\) are real-valued diagonal matrix
    symbols for real \(U\) and are explicitly defined as follows:
    \begin{itemize}
        \item the high-frequency diagonalized quasilinear order-\(3/2\)
        symbol is
        \begin{equation}\label{eq:positive_order_diagonal_ql_symbol_3d}
            {\bm D}^{\bC}_{\frac32}
            \pare{\Atw,\background\middle|U;x,\xi}
            \defeq
            \pare{1-\chi\pare{\xi}}\,
            \omega_{\frac32}
            \pare{\Atw,\background\middle|\xi}
            \dot\theta_0
            \pare{U;x,\hat\xi}
            {\bm J}_{\rm d}
            \in
            \pare{\Gamma_{\geq0}^{\frac32}}^{2\times2};
        \end{equation}
        here \({\bm J}_{\rm d}\) is defined in
        \cref{eq:complex_diagonal_sign_matrix},
        \(\omega_{\frac32}\) is the order-\(3/2\) linear frequency from
        \cref{eq:complex_frequency_order_three_halves}, and
        \(\dot\theta_0
        \pare{U;x,\hat\xi}\)
        is the real-valued homogeneous degree-zero scalar symbol given by
        \begin{equation}\label{eq:principal_scalar_l_theta_3d}
            \dot\ell_0
            \pare{U;x,\hat\xi}
            \defeq
            1+
            \dot e_0
            \pare{U;x,\hat\xi},
            \qquad
            \dot\theta_0
            \pare{U;x,\hat\xi}
            \defeq
            \sqrt{
                \dot\ell_0^2
                -
                \dot f_0^2
            }
            \pare{U;x,\hat\xi}.
        \end{equation}
        The coefficients
        \(\dot e_0,\dot f_0\in\Gamma_{\geq1}^0\) are the real-valued homogeneous
        quasilinear coefficients introduced in
        \cref{eq:dote0,eq:dotf0};

        \item the scalar transport symbol is
        \begin{equation}\label{eq:positive_order_diagonal_transport_symbol_3d}
            \dot{\bm D}^{\bC}_{1}
            \pare{\Atw,\background\middle|U;x,\xi}
            \defeq
            \pare{
                \mathpzc V_{\textnormal{can}}
                \pare{
                    \Atw,\background
                    \middle|
                    \eta\pare{U},\tilde\psi\pare{U};
                    x
                }
                \cdot\xi
            }
            \Id_{\bC^2}
            \in
            \pare{\Gamma_{\geq0}^1}^{2\times2};
        \end{equation}
        here \(\mathpzc V_{\textnormal{can}}\) is the canonical transport
        velocity defined in \cref{eq:Vpitilde_main}.  This block is scalar in
        the complex components and is therefore unaffected by the
        diagonalizing matrices at positive order;

       \item the high-frequency diagonal order-\(1/2\) symbol is \begin{equation}\label{eq:positive_order_diagonal_half_symbol_3d} \dot{\bm D}^{\bC}_{\frac12} \pare{\Atw,\background\middle|U;x,\xi} \defeq - \pare{1-\chi\pare{\xi}}\, \dot\alpha^{\bC,{\rm pr}}_{0} \pare{\Atw,\background\middle|U;x,\hat\xi} \av{\xi}^{\frac12} {\bm J}_{\rm d} \in \pare{\Gamma_{\geq1}^{\frac12}}^{2\times2}. \end{equation} The coefficient \(\dot\alpha^{\bC,{\rm pr}}_0\in\Gamma_{\geq1}^0\) is the diagonal component of the full order-\(1/2\) coefficient after the principal diagonalization; it is defined in \cref{eq:principal_half_output_decomposition_3d,eq:principal_alpha_principal_half_3d}. 
    \end{itemize}

    \item
    \({\bm A}_{\bra{0,\geq1}}
    \pare{\Atw,\background\middle|U;x,\xi}
    \in\pare{\Gamma_{\geq1}^0}^{2\times2}\)
    is a generic bounded matrix symbol in the sense of
    \cref{notation:paralinearization}.

    \item
    \({\bm R}_{\bra{\varrho,\geq1}}
    \pare{\Atw,\background\middle|U}
    \in\pare{\cR_{\geq1}^\varrho}^{2\times2}\)
    is a generic smoothing remainder in the sense of
    \cref{notation:paralinearization}.
\end{itemize}
\end{prop}

    The operator ${\bm B}$ is constructed as a finite composition of paradifferential
parametrices, each eliminating a block of off-diagonal symbols at a prescribed order.
The diagonalization is carried out in two successive steps, corresponding to the
principal and subprincipal orders, in
\cref{prop:principal_order_diagonalization_3d,prop:subprincipal_order_diagonalization_3d}
respectively; \cref{prop:positive_order_diagonalization_3d} is then a direct
consequence, as detailed in \cref{sec:proof_positive_order_diagonalization_3d}.

\subsubsection{Diagonalization at principal order}
\label{sec:principal_order_diagonalization_3d}
Here, we perform the diagonalization of the order-$\frac{3}{2}$ in the complex paradifferential system \cref{eq:complex_paradifferential_system_3d}. The result reads as follows.

\begin{prop}[Principal-order diagonalization]
\label{prop:principal_order_diagonalization_3d}
Let the assumptions of
\cref{prop:complex_paradifferential_system_3d} hold, and let \(U\) solve
\cref{eq:complex_paradifferential_system_3d}.  There exists a
bounded and invertible operator
\[
    {\bm B}_{0}
    \pare{\Atw,\background\middle|U}
    \in
    \cL\pare{H^s}
    \qquad
    \forall\,s\geq s_0,
\]
such that, defining
\begin{equation}\label{eq:principal_diagonalized_unknown_3d}
    U_{0}
    \defeq
    {\bm B}_{0}
    \pare{\Atw,\background\middle|U}^{-1}U ,
\end{equation}
the unknown \(U_0\) solves
\begin{equation}\label{eq:principal_diagonalized_system_3d}
\begin{aligned}
    \partial_t U_{0}
    ={}&
    \OpBW{
        {\bm A}^{\bC,{\rm lin}}_{\rm lo}
        \pare{\Atw,\background\middle|\xi}
    }U_0
    \\
    &+
    \ii
    \OpBW{
        {\bm D}^{\bC}_{\frac32}
        \pare{\Atw,\background\middle|U;x,\xi}}U_{0}\\
        &+\ii\OpBW{
        \dot{\bm D}^{\bC}_{1}
        \pare{\Atw,\background\middle|U;x,\xi}
    }U_{0}
    \\
    &+
    \OpBW{
        \pare{1-\chi\pare{\xi}}\,
        \dot{\bm A}^{\bC,{\rm pr}}_{0}
        \pare{\Atw,\background\middle|U;x,\hat\xi}
        \av{\xi}^{\frac12}
    }U_{0}
    \\
    &+
    \OpBW{
        {\bm A}^{\bC,{\rm pr}}_{\bra{0,\geq1}}
        \pare{\Atw,\background\middle|U;x,\xi}
    }U_{0}
    +
    {\bm R}^{\bC,{\rm pr}}_{\bra{\varrho,\geq1}}
    \pare{\Atw,\background\middle|U}U_{0}.
\end{aligned}
\end{equation}
The objects entering
\cref{eq:principal_diagonalized_unknown_3d,eq:principal_diagonalized_system_3d}
are defined as follows.
\begin{itemize}
    \item
    \({\bm B}_{0}\pare{U}\) is the order-zero
    paradifferential operator
    \begin{equation}\label{eq:principal_diagonalizer_operator_3d}
        {\bm B}_{0}
        \pare{U}
        \defeq
        \OpBW{
            {\bm P}_{0}
            \pare{U;x,\xi}
        } ,
    \end{equation}
    where \({\bm P}_{0}\) is the non-homogeneous principal diagonalizer
    \begin{equation}\label{eq:principal_nonhomogeneous_diagonalizer_symbol_3d}
        {\bm P}_{0}
        \pare{U;x,\xi}
        \defeq
        \chi\pare{\xi}\Id_{\bC^2}
        +
        \pare{1-\chi\pare{\xi}}\,
        \dot{\bm P}_{0}
        \pare{U;x,\hat\xi}.
    \end{equation}
    The high-frequency homogeneous diagonalizer is
    \begin{equation}\label{eq:principal_diagonalizer_symbol_3d}
        \dot{\bm P}_{0}
        \pare{U;x,\hat\xi}
        \defeq
        \frac{1}{
            \sqrt{
                1-\dot\chi_0^2
            }
        }
        \begin{bmatrix}
            1 & \dot\chi_0
            \\
            \dot\chi_0 & 1
        \end{bmatrix}
        \pare{U;x,\hat\xi},\qquad
        \dot{\bm P}_{0}-\Id_{\bC^2}\in
        \pare{\Gamma_{\geq1}^0}^{2\times2}.
    \end{equation}
    Here \(\dot\chi_0\in\Gamma_{\geq1}^0\) is the real homogeneous scalar defined by \begin{equation}\label{eq:dotchi0}
           \begin{aligned}
            \dot\chi_0
            \pare{U;x,\hat\xi}
            &\defeq
            -
            \frac{
                \dot f_0
                \pare{U;x,\hat\xi}
            }{
                \dot\ell_0
                \pare{U;x,\hat\xi}
                +
                \dot\theta_0
                \pare{U;x,\hat\xi}
            }
            \\
            &=
            \frac{
                \pare{1+\av{\nabla\eta\pare U}^2}^{\frac34}
                -
                \pare{
                    1+\av{\nabla^\perp\eta\pare U\cdot\hat{\xi}}^2
                }^{\frac14}
            }{
                \pare{1+\av{\nabla\eta\pare U}^2}^{\frac34}
                +
                \pare{
                    1+\av{\nabla^\perp\eta\pare U\cdot\hat{\xi}}^2
                }^{\frac14}
            }
            \in\Gamma_{\geq1}^0;
        \end{aligned}
       \end{equation}

    \item
    \({\bm A}^{\bC,{\rm lin}}_{\rm lo}\) is the flat low-frequency linear
    multiplier explicitly defined in  \cref{eq:linear_low_frequency_defect_complex_3d};
    \item
    \({\bm D}^{\bC}_{\frac32}\in
    \pare{\Gamma_{\geq0}^{\frac32}}^{2\times2}\) and
    \(\dot{\bm D}^{\bC}_{1}\in
    \pare{\Gamma_{\geq0}^1}^{2\times2}\) are the diagonal quasilinear and transport
    symbols already defined in
    \cref{eq:positive_order_diagonal_ql_symbol_3d,eq:positive_order_diagonal_transport_symbol_3d}.

    \item
    \(\dot{\bm A}^{\bC,{\rm pr}}_0
    \pare{\Atw,\background\middle|U;x,\hat\xi}
    \in\pare{\Gamma_{\geq1}^0}^{2\times2}\)
    is defined as
    \begin{equation}\label{eq:principal_half_output_decomposition_3d}
        \dot{\bm A}^{\bC,{\rm pr}}_{0}
        \pare{\Atw,\background\middle|U;x,\hat\xi}
        \defeq
        \Big(
            \dot{\bm P}_{0}^{-1}
            \dot{\bm A}^{\bC}_{0}
            \dot{\bm P}_{0}
        \Big)
        \pare{\Atw,\background\middle|U;x,\hat\xi} =
        \ii
        \begin{bmatrix}
            \dot\alpha^{\bC,{\rm pr}}_{0}
            &
            \dot\beta^{\bC,{\rm pr}}_{0}
            \\
            -\dot\beta^{\bC,{\rm pr}}_{0}
            &
            -\dot\alpha^{\bC,{\rm pr}}_{0}
                    \end{bmatrix},
    \end{equation}
    where
    \begin{align}
        \dot\alpha^{\bC,{\rm pr}}_{0}
        &\defeq
        \frac{
            \pare{1+\dot\chi_0^2}\dot\alpha^{\bC}_{0}
            +
            2\dot\chi_0\dot\beta^{\bC}_{0}
        }{
            1-\dot\chi_0^2
        },
        \label{eq:principal_alpha_principal_half_3d}
        \\
        \dot\beta^{\bC,{\rm pr}}_{0}
        &\defeq
        \frac{
            2\dot\chi_0\dot\alpha^{\bC}_{0}
            +
            \pare{1+\dot\chi_0^2}\dot\beta^{\bC}_{0}
        }{
            1-\dot\chi_0^2
        }.
        \label{eq:principal_beta_principal_half_3d}
    \end{align}
    Here
    \(\dot\alpha^{\bC}_{0},\dot\beta^{\bC}_{0}\in \Gamma_{\geq1}^0\)
    are defined in
    \cref{eq:complex_half_alpha_coefficient_3d,eq:complex_half_beta_coefficient_3d}.

    \item
    \({\bm A}^{\bC,{\rm pr}}_{\bra{0,\geq1}}
    \pare{\Atw,\background\middle|U;x,\xi}
    \in\pare{\Gamma_{\geq1}^0}^{2\times2}\)
    is a generic bounded matrix symbol in the sense of
    \cref{notation:paralinearization}. 
    \item
    \({\bm R}^{\bC,{\rm pr}}_{\bra{\varrho,\geq1}}
    \pare{\Atw,\background\middle|U}
    \in\pare{\cR_{\geq1}^\varrho}^{2\times2}\)
    is a generic smoothing remainder in the sense of
    \cref{notation:paralinearization}.
\end{itemize}
\end{prop}

\begin{proof}
Throughout the proof all symbols are evaluated at
\[
    \pare{\Atw,\background\middle|U;x,\xi},
\]
unless explicitly stated otherwise.  
We use the convention of
\cref{notation:paralinearization}: the symbols
\({\bm A}_{\bra{m,\geq1}}\) and the remainders
\({\bm R}_{\bra{\varrho,\geq1}}\) are generic and may change from line to line.  
We
also use componentwise the symbolic calculus of
\cref{def:asympt_symbol,prop:composition_paradiff} and the remainder rules
\cref{lem:composition_rules}-\ref{item:composition_rules_paradiff_remainder},
\cref{lem:composition_rules}-\ref{item:composition_rules_remainder_remainder}.

\begin{proofpart}[The diagonalizer and its symbolic inverse]
We first verify that the high-frequency symbol defining the principal
diagonalizer is admissible.  By
\cref{eq:dote0_closed_principal_3d,eq:dotf0_closed_principal_3d,eq:principal_scalar_l_theta_3d},
 we have the explicit expressions
\begin{align}
    \dot\ell_0
    &=
    1+\dot e_0
    =
    \frac12
    \pare{
        \sqrt{
            1+\av{\nabla^\perp\eta\cdot\hat{\xi}}^2
        }
        +
        \frac{
            1+\av{\nabla^\perp\eta\cdot\hat{\xi}}^2
        }{
            \pare{1+\av{\nabla\eta}^2}^{\frac32}
        }
    },
    \label{eq:dotell0_closed_principal_3d}
    \\
    \dot\theta_0
    &=
    \pare{\frac{ 1+\av{\nabla^\perp\eta\cdot\hat{\xi}}^2
    }{1+\av{\nabla\eta}^2
    }}^{\frac34}.\label{eq:dottheta0_closed_principal_3d}
\end{align}
In particular,

\[
    \dot\ell_0-1\in\Gamma_{\geq1}^0,
    \qquad
    \dot\theta_0-1\in\Gamma_{\geq1}^0.
\]
After increasing \(s_0\) and decreasing \(r_0\), we get
\[
    \dot\ell_0\geq \frac34,
    \qquad
    \dot\theta_0\geq \frac34,
    \qquad
    \dot\ell_0+\dot\theta_0\geq \frac32\cdot
\]
This allows to define
\begin{equation}\label{,eq:dotell0_closed_principal_3d,eq:principal_chi0_symbol_class_refined_3d}
    \dot\chi_0\defeq
    -\frac{\dot f_0}{\dot\ell_0+\dot\theta_0}
    \in\Gamma_{\geq1}^0,
    \qquad
    \norm{\dot\chi_0}_{L^\infty}
    \leq C r_0 .
\end{equation}
In particular, decreasing \(r_0\) once more if necessary,
\begin{equation}\label{eq:principal_chi0_smallness_3d}
    1-\dot\chi_0^2\geq \frac12\cdot
\end{equation}
Moreover, after straightforward simplifications, we obtain from \cref{eq:dotf0_closed_principal_3d,eq:dottheta0_closed_principal_3d}, the explicit expression \cref{eq:dotchi0}. Now, using
\cref{eq:principal_diagonalizer_symbol_3d,eq:principal_chi0_smallness_3d}
and the substitution rule
\cref{lem:composition_rules}-\ref{item:composition_rules_substitution_symbol},
we get
\begin{equation}\label{eq:principal_P0_symbol_class_refined_3d}
    \dot{\bm P}_0-\Id_{\bC^2}\in\pare{\Gamma_{\geq1}^0}^{2\times2},
    \qquad
    \dot{\bm P}_0\in\pare{\Gamma^0_{\ge0}}^{2\times2}.
\end{equation}
Moreover its high-frequency pointwise inverse is
\begin{equation}\label{eq:principal_diagonalizer_inverse_symbol_3d}
    \dot{\bm P}_{0}^{-1}
    =
    \frac{1}{\sqrt{1-\dot\chi_0^2}}
    \begin{bmatrix}
        1 & -\dot\chi_0
        \\
        -\dot\chi_0 & 1
    \end{bmatrix}.
\end{equation}
It satisfies
\begin{equation}\label{eq:principal_P0_inverse_symbol_class_refined_3d}
    \dot{\bm P}_0^{-1}-\Id_{\bC^2}
    \in\pare{\Gamma_{\geq1}^0}^{2\times2},
    \qquad
    \dot{\bm P}_0^{-1}
    \in\pare{\Gamma^0_{\ge0}}^{2\times2}.
\end{equation}

We next record the differential algebra of the homogeneous diagonalizer,
which is the key point both for the symbolic inverse and for the cancellation
of the order-\(1/2\) Weyl correction below.  Set
\begin{equation}\label{eq:subprincipal_offdiagonal_matrix_3d}
    {\bm X}_{\rm od}
    \defeq
    \begin{bmatrix}
        0 & 1
        \\
        1 & 0
    \end{bmatrix}.
\end{equation}
Since
\begin{equation}\label{eq:principal_diagonalizer_exponential_form_3d}
    \dot{\bm P}_0
    =
    \exp\pare{\dot\varphi_0{\bm X}_{\rm od}},
    \qquad
    \dot\varphi_0\defeq\operatorname{arctanh}\dot\chi_0,
\end{equation}
we have, for every
\(y\in\set{x_1,x_2,\xi_1,\xi_2}\),
\begin{equation}\label{eq:principal_diagonalizer_derivative_identity_3d}
    \dot{\bm\Lambda}_y
    \defeq
    \dot{\bm P}_0^{-1}\partial_y\dot{\bm P}_0
    =
    \pare{\partial_y\dot\varphi_0}{\bm X}_{\rm od}.
\end{equation}
Thus all the matrices \(\dot{\bm\Lambda}_y\) commute.  Moreover
\(\partial_y\dot{\bm P}_0^{-1}
=-\dot{\bm\Lambda}_y\dot{\bm P}_0^{-1}\), and therefore
\begin{equation}\label{eq:P0_inverse_defect_principal_diag_3d}
\begin{aligned}
    \pbra{\dot{\bm P}_0^{-1}}{\dot{\bm P}_0}
    &=
    \sum_{j=1}^{2}
    \Big[
        \pare{\partial_{\xi_j}\dot{\bm P}_0^{-1}}
        \pare{\partial_{x_j}\dot{\bm P}_0}
        -
        \pare{\partial_{x_j}\dot{\bm P}_0^{-1}}
        \pare{\partial_{\xi_j}\dot{\bm P}_0}
    \Big]
    \\
    &=
    \sum_{j=1}^{2}
    \Big[
        -\dot{\bm\Lambda}_{\xi_j}\dot{\bm\Lambda}_{x_j}
        +
        \dot{\bm\Lambda}_{x_j}\dot{\bm\Lambda}_{\xi_j}
    \Big]
    =0.
\end{aligned}
\end{equation}
Hence the homogeneous symbolic inverse has no order-\(\pare{-1}\) defect
\begin{equation}\label{eq:principal_B0_symbolic_inverse_3d}
    \dot{\bm P}_0^{-1}
    \#_\varrho
    \dot{\bm P}_0
    =
    \Id_{\bC^2}
    +
    {\bm A}_{\bra{-2,\geq1}},
    \qquad
    {\bm A}_{\bra{-2,\geq1}}\in\pare{\Gamma_{\geq1}^{-2}}^{2\times2},
\end{equation}
and the same identity holds with the two factors exchanged.

We now pass from \(\dot{\bm P}_0\) to the final non-homogeneous symbol
\({\bm P}_0\) of
\cref{eq:principal_nonhomogeneous_diagonalizer_symbol_3d}.  Since
\(\dot{\bm P}_0-\Id_{\bC^2}=O(r_0)\) in the order-zero symbolic seminorms,
\[
    {\bm P}_0-\Id_{\bC^2}
    =
    \pare{1-\chi\pare{\xi}}\pare{\dot{\bm P}_0-\Id_{\bC^2}}
    \in\pare{\Gamma_{\geq1}^0}^{2\times2},
    \qquad
    \norm{{\bm P}_0-\Id_{\bC^2}}_{\Gamma^0}
    \leq C r_0.
\]
Thus \({\bm B}_0=\OpBW{{\bm P}_0}\) is bounded on \(H^s\), and after
decreasing \(r_0\) it is invertible on \(H^s\) by the Neumann series.

Let \({\bm P}_0^{[-1]}\) be the pointwise inverse of \({\bm P}_0\).  Then
\[
    {\bm P}_0^{[-1]}-\Id_{\bC^2}
    \in\pare{\Gamma_{\geq1}^0}^{2\times2}
    \qquad\textnormal{and}\qquad
    {\bm P}_0^{[-1]}=\dot{\bm P}_0^{-1}
    \quad\textnormal{on }\supp^{\complement}\pare{\chi}.
\]
On the region \(\chi=0\),
\cref{eq:P0_inverse_defect_principal_diag_3d} applies; on the region
\(\chi=1\), both \({\bm P}_0\) and \({\bm P}_0^{[-1]}\) are the identity.
Every discrepancy from the homogeneous calculation is supported in the
fixed transition annulus; since the symbols are smooth and compactly
supported in \(\xi\) there, it belongs to \(\Gamma^{-N}\) for every \(N\).
Consequently, after absorbing the transition contribution into the
order-\(\pare{-2}\) symbol,
\begin{equation}\label{eq:principal_B0_parametrix_inverse_3d}
    \OpBW{{\bm P}_0^{[-1]}}{\bm B}_0
    =
    \Id_{\bC^2}
    +
    \OpBW{{\bm A}_{\bra{-2,\geq1}}}
    +
    {\bm R}_{\bra{\varrho,\geq1}},
    \qquad
    {\bm A}_{\bra{-2,\geq1}}\in\pare{\Gamma_{\geq1}^{-2}}^{2\times2},
\end{equation}
and analogously on the other side.  Since \({\bm B}_0\) is already known
to be invertible, its true inverse has the form
\begin{equation}\label{eq:principal_B0_true_inverse_expansion_3d}
    {\bm B}_0^{-1}
    =
    \OpBW{
        {\bm P}_0^{[-1]}
        +
        {\bm A}^{\sharp}_{\bra{-2,\geq1}}
    }
    +
    {\bm R}_{\bra{\varrho,\geq1}},
    \qquad
    {\bm A}^{\sharp}_{\bra{-2,\geq1}}
    \in\pare{\Gamma_{\geq1}^{-2}}^{2\times2},
\end{equation}
after decreasing and relabelling \(\varrho\).

Let now \({\bm A}\in(\Gamma^m_{\ge0})^{2\times2}\), \(m\leq3/2\), be one
of the symbols treated below.  Combining
\cref{eq:principal_B0_true_inverse_expansion_3d} with
\cref{def:asympt_symbol,prop:composition_paradiff} and
\cref{lem:composition_rules}-\ref{item:composition_rules_paradiff_remainder},
we obtain
\begin{equation}\label{eq:principal_B0_conjugation_formula_3d}
\begin{aligned}
    {\bm B}_{0}^{-1}
    \OpBW{{\bm A}}
    {\bm B}_{0}
    =
    \OpBW{
        {\bm C}_{{\bm P}_0}\bra{{\bm A}}
        +
        {\bm A}_{\bra{m-2,\geq1}}
    }
    +
    {\bm R}_{\bra{\varrho,\geq1}},
\end{aligned}
\end{equation}
where
\begin{equation}\label{eq:principal_symbolic_conjugation_map_3d}
    {\bm C}_{{\bm P}_0}\bra{{\bm A}}
    \defeq
    {\bm P}_0^{[-1]}\#_\varrho {\bm A}\#_\varrho{\bm P}_0 .
\end{equation}
The additional term \({\bm A}_{\bra{m-2,\geq1}}\) is solely due to replacing the
true inverse by the symbolic parametrix.  Since \(m\leq3/2\), it has order
at most \(-1/2\), hence it is bounded and will be absorbed into
\({\bm A}_{\bra{0,\geq1}}\) below.
\end{proofpart}

\begin{proofpart}[Conjugation of the equation]
We write
\[
    U={\bm B}_0U_0 .
\]
Differentiating gives
\begin{equation}\label{eq:principal_diag_time_derivative_3d}
    \partial_tU_0
    =
    {\bm B}_0^{-1}\partial_tU
    -
    {\bm B}_0^{-1}\pare{\partial_t{\bm B}_0}U_0 .
\end{equation}
Set
\begin{equation}\label{eq:principal_qll_matrix_before_diag_3d}
    \dot{\bm C}_{{\rm ql},0}
    \defeq
    {\bm J}_{\rm d}
    \pare{
        \Id_{\bC^2}
        +
        \dot{\bm E}^{\bC}_{{\rm ql},0}
    }.
\end{equation}
By \cref{eq:complex_ql_matrix_coefficient_3d,eq:dote0,eq:dotf0},
\begin{equation}\label{eq:principal_diag_C0_entries_3d}
    \dot{\bm C}_{{\rm ql},0}
    =
    \begin{bmatrix}
        -\dot\ell_{0} & -\dot f_{0}
        \\
        \dot f_{0} & \dot\ell_{0}
    \end{bmatrix},
\end{equation}
where \(\dot\ell_0,\dot f_0\) are defined in
\cref{eq:principal_scalar_l_theta_3d}.  We rewrite
\cref{eq:complex_paradifferential_system_3d} as
\begin{equation}\label{eq:principal_diag_compact_system_3d}
    \partial_tU
    =
    \OpBW{
        {\bm L}_{\rm lo}
        +
        {\bm L}_{\frac32}
        +
        {\bm L}_{1}
        +
        {\bm L}_{\frac12}
        +
        {\bm L}_{0}
    }U
    +
    {\bm R}^{\bC}_{\bra{\varrho,\geq1}}U,
\end{equation}
with
\begin{align}
    {\bm L}_{\rm lo}
    &\defeq
    {\bm A}^{\bC,{\rm lin}}_{\rm lo}
    \pare{\Atw,\background\middle|\xi},
    \label{eq:principal_Llo_definition_3d}
    \\
    {\bm L}_{\frac32}
    &\defeq
    \pare{1-\chi\pare{\xi}}\,
    \ii\omega_{\frac32}\dot{\bm C}_{{\rm ql},0},
    \label{eq:principal_L32_definition_3d}
    \\
    {\bm L}_{1}
    &\defeq
    \ii\dot{\bm D}^{\bC}_{1},
    \label{eq:principal_L1_definition_3d}
    \\
    {\bm L}_{\frac12}
    &\defeq
    \pare{1-\chi\pare{\xi}}\,
    \dot{\bm A}^{\bC}_{0}\av{\xi}^{\frac12},
    \label{eq:principal_L12_definition_3d}
    \\
    {\bm L}_{0}
    &\defeq
    {\bm A}_{\bra{0,\geq1}}.
    \label{eq:principal_L0_definition_3d}
\end{align}
The classes are
\[
    {\bm L}_{\rm lo}\in\pare{\Gamma^0_{\ge0}}^{2\times2},
    \qquad
    {\bm L}_{\frac32}\in\pare{\Gamma^{\frac32}_{\ge0}}^{2\times2},
    \qquad
    {\bm L}_{1}\in\pare{\Gamma^1_{\ge0}}^{2\times2}, \qquad {\bm L}_{\frac12}\in\pare{\Gamma_{\geq1}^{\frac12}}^{2\times2},
    \qquad
    {\bm L}_{0}\in\pare{\Gamma_{\geq1}^0}^{2\times2}.
\]
%\[
    %{\bm L}_{\frac12}\in\pare{\Gamma_{\geq1}^{\frac12}}^{2\times2},
    %\qquad
    %{\bm L}_{0}\in\pare{\Gamma_{\geq1}^0}^{2\times2}.
%\]
Thus all five blocks are admissible inputs for
\cref{eq:principal_B0_conjugation_formula_3d}.  The factors
\(1-\chi\pare{\xi}\) encode that the positive capillary diagonalization is
performed only at high frequency.
\end{proofpart}

\begin{proofpart}[The order-\(3/2\) block]
We first work on the homogeneous high-frequency region \(\chi=0\), where
\({\bm P}_0=\dot{\bm P}_0\) and
\({\bm L}_{\frac32}=\ii\omega_{\frac32}\dot{\bm C}_{{\rm ql},0}\).
A direct multiplication gives
\begin{equation}\label{Conj:PC}
    \dot{\bm P}_{0}^{-1}
    \dot{\bm C}_{{\rm ql},0}
    \dot{\bm P}_{0}
    =
    \frac{1}{1-\dot\chi_0^2}
    \begin{bmatrix}
        -\dot\ell_0\pare{1+\dot\chi_0^2}
        -2\dot f_0\dot\chi_0
        &
        -2\dot\ell_0\dot\chi_0
        -\dot f_0\pare{1+\dot\chi_0^2}
        \\
        2\dot\ell_0\dot\chi_0
        +\dot f_0\pare{1+\dot\chi_0^2}
        &
        \dot\ell_0\pare{1+\dot\chi_0^2}
        +2\dot f_0\dot\chi_0
    \end{bmatrix}.
\end{equation}
The choice of \(\dot\chi_0\) in \cref{eq:dotchi0} gives
\begin{align}
    2\dot\ell_0\dot\chi_0
    +
    \dot f_0\pare{1+\dot\chi_0^2}
    &=0,
    \label{eq:principal_diag_off_identity_3d}
    \\
    \dot\ell_0\pare{1+\dot\chi_0^2}
    +
    2\dot f_0\dot\chi_0
    &=
    \dot\theta_0\pare{1-\dot\chi_0^2}.
    \label{eq:principal_diag_diag_identity_3d}
\end{align}
Inserting these identities into \cref{Conj:PC} yields
\begin{equation}\label{eq:pointwise_principal_diagonalization_3d}
    \dot{\bm P}_{0}^{-1}
    \dot{\bm C}_{{\rm ql},0}
    \dot{\bm P}_{0}
    =
    \dot\theta_0{\bm J}_{\rm d}.
\end{equation}
Set
\begin{equation}\label{eq:principal_Dsharp_definition_3d}
    {\bm D}_{\frac32}^{\sharp}
    \defeq
    \ii\omega_{\frac32}\dot\theta_0{\bm J}_{\rm d}.
\end{equation}
Equivalently,
\begin{equation}\label{eq:principal_homogeneous_block_factorization_3d}
    \ii\omega_{\frac32}\dot{\bm C}_{{\rm ql},0}
    =
    \dot{\bm P}_0
    {\bm D}_{\frac32}^{\sharp}
    \dot{\bm P}_0^{-1},
\end{equation}
and therefore
\begin{equation}\label{eq:principal_pointwise_L32_conjugation_3d}
    \dot{\bm P}_{0}^{-1}
    \pare{\ii\omega_{\frac32}\dot{\bm C}_{{\rm ql},0}}
    \dot{\bm P}_{0}
    =
    {\bm D}_{\frac32}^{\sharp}.
\end{equation}

We now verify the order-\(1/2\) cancellation in the principal
conjugation. On \(\chi=0\), the Bony--Weyl expansion gives
\begin{equation}\label{eq:principal_C0_L32_expansion_3d}
\begin{aligned}
    {\bm C}_{\dot{\bm P}_0}
    \bra{\ii\omega_{\frac32}\dot{\bm C}_{{\rm ql},0}}
    ={}&
    \dot{\bm P}_0^{-1}
    \pare{\ii\omega_{\frac32}\dot{\bm C}_{{\rm ql},0}}
    \dot{\bm P}_0
    \\
    &+
    \frac{1}{2\ii}
    \left(
        \pbra{\dot{\bm P}_0^{-1}}
              {\ii\omega_{\frac32}\dot{\bm C}_{{\rm ql},0}}
        \dot{\bm P}_0
        +
        \pbra{
            \dot{\bm P}_0^{-1}
            \pare{\ii\omega_{\frac32}\dot{\bm C}_{{\rm ql},0}}
        }{\dot{\bm P}_0}
    \right)
    +
    {\bm A}_{\bra{-\frac12,\geq1}}.
\end{aligned}
\end{equation}
Differentiating
\cref{eq:principal_homogeneous_block_factorization_3d} gives, for every
\(y\in\set{x_1,x_2,\xi_1,\xi_2}\),
\begin{equation}\label{eq:principal_differentiated_block_identity_3d}
    \dot{\bm P}_0^{-1}
    \partial_y\pare{\ii\omega_{\frac32}\dot{\bm C}_{{\rm ql},0}}
    \dot{\bm P}_0
    =
    \partial_y{\bm D}_{\frac32}^{\sharp}
    +
    \comm{ \dot{\bm\Lambda}_y}{{\bm D}_{\frac32}^{\sharp}}
      .
\end{equation}
By
\cref{eq:principal_diagonalizer_derivative_identity_3d}, every
\(\dot{\bm\Lambda}_y\) is a scalar multiple of
\({\bm X}_{\rm od}\), whereas
\({\bm D}_{\frac32}^{\sharp}\) and all its derivatives are scalar
multiples of \({\bm J}_{\rm d}\).  Since
\[
    {\bm X}_{\rm od}{\bm J}_{\rm d}
    +
    {\bm J}_{\rm d}{\bm X}_{\rm od}
    =0,
\]
we have
\begin{equation}\label{eq:principal_antidiagonal_anticommutation_3d}
    \dot{\bm\Lambda}_y
    \partial_z{\bm D}_{\frac32}^{\sharp}
    +
    \partial_z{\bm D}_{\frac32}^{\sharp}
    \dot{\bm\Lambda}_y
    =0,
    \qquad
    \forall\,y,z\in\set{x_1,x_2,\xi_1,\xi_2}.
\end{equation}
Furthermore,
\(\partial_y\dot{\bm P}_0^{-1}
=-\dot{\bm\Lambda}_y\dot{\bm P}_0^{-1}\), so that
\[
    \partial_y\!
    \left(
        \dot{\bm P}_0^{-1}
        \pare{\ii\omega_{\frac32}\dot{\bm C}_{{\rm ql},0}}
    \right)
    \partial_z\dot{\bm P}_0
    =
    \pare{
        \partial_y{\bm D}_{\frac32}^{\sharp}
        -
        {\bm D}_{\frac32}^{\sharp}\dot{\bm\Lambda}_y
    }
    \dot{\bm\Lambda}_z.
\]
Consequently, the formally order-\(1/2\) correction in
\cref{eq:principal_C0_L32_expansion_3d} is
\begin{equation}\label{eq:principal_weyl_correction_expanded_3d}
\begin{aligned}
    &\pbra{\dot{\bm P}_0^{-1}}
        {\ii\omega_{\frac32}\dot{\bm C}_{{\rm ql},0}}
        \dot{\bm P}_0
    +
    \pbra{
        \dot{\bm P}_0^{-1}
        \pare{\ii\omega_{\frac32}\dot{\bm C}_{{\rm ql},0}}
    }{\dot{\bm P}_0}
    \\
    &\quad=
    \sum_{j=1}^{2}
    \left[
        -\dot{\bm\Lambda}_{\xi_j}
        \pare{
            \partial_{x_j}{\bm D}_{\frac32}^{\sharp}
            +
            \left[
                \dot{\bm\Lambda}_{x_j},
                {\bm D}_{\frac32}^{\sharp}
            \right]
        }
        +
        \dot{\bm\Lambda}_{x_j}
        \pare{
            \partial_{\xi_j}{\bm D}_{\frac32}^{\sharp}
            +
            \left[
                \dot{\bm\Lambda}_{\xi_j},
                {\bm D}_{\frac32}^{\sharp}
            \right]
        } \right. 
        \\
        &\qquad\quad
        +\left. 
        \pare{
            \partial_{\xi_j}{\bm D}_{\frac32}^{\sharp}
            -
            {\bm D}_{\frac32}^{\sharp}
            \dot{\bm\Lambda}_{\xi_j}
        }
        \dot{\bm\Lambda}_{x_j}
        -
        \pare{
            \partial_{x_j}{\bm D}_{\frac32}^{\sharp}
            -
            {\bm D}_{\frac32}^{\sharp}
            \dot{\bm\Lambda}_{x_j}
        }
        \dot{\bm\Lambda}_{\xi_j}
    \right].
\end{aligned}
\end{equation}
The terms linear in \(\dot{\bm\Lambda}_y\) vanish by
\cref{eq:principal_antidiagonal_anticommutation_3d}.  The quadratic terms
vanish because all \(\dot{\bm\Lambda}_y\) commute with one another and
anti-commute with \({\bm D}_{\frac32}^{\sharp}\).  Hence
\begin{equation}\label{eq:principal_weyl_correction_vanishes_3d}
    \label{eq:principal_L32_half_weyl_correction_3d}
    \pbra{\dot{\bm P}_0^{-1}}
        {\ii\omega_{\frac32}\dot{\bm C}_{{\rm ql},0}}
        \dot{\bm P}_0
    +
    \pbra{
        \dot{\bm P}_0^{-1}
        \pare{\ii\omega_{\frac32}\dot{\bm C}_{{\rm ql},0}}
    }{\dot{\bm P}_0}
    =0.
\end{equation}

It remains to pass to the non-homogeneous symbols.  On \(\chi=0\) the
preceding calculation applies verbatim.  On \(\chi=1\), one has
\({\bm P}_0=\Id_{\bC^2}\) and \({\bm L}_{\frac32}=0\).  Every discrepancy
in the transition annulus, including every term in which a
\(\xi\)-derivative hits the cut-off, is smooth and compactly supported in
\(\xi\), hence belongs to \(\Gamma^{-N}\) for every \(N\).  Combining
\cref{eq:principal_C0_L32_expansion_3d,eq:principal_pointwise_L32_conjugation_3d,eq:principal_weyl_correction_vanishes_3d,eq:principal_B0_conjugation_formula_3d},
we obtain
\begin{equation}\label{eq:principal_block_conjugated_result_3d}
    {\bm B}_{0}^{-1}
    \OpBW{{\bm L}_{\frac32}}
    {\bm B}_{0}
    =
    \OpBW{
        \ii{\bm D}^{\bC}_{\frac32}
        +
        {\bm A}_{\bra{0,\geq1}}
    }
    +
    {\bm R}_{\bra{\varrho,\geq1}}.
\end{equation}
In particular, the order-\(3/2\) block produces no homogeneous
order-\(1/2\) symbol.
\end{proofpart}

\begin{proofpart}[The flat low-frequency linear block]
We treat the flat multiplier
\[
    {\bm L}_{\rm lo}
    =
    {\bm A}^{\bC,{\rm lin}}_{\rm lo}
    \pare{\Atw,\background\middle|\xi}.
\]
By \cref{eq:principal_nonhomogeneous_diagonalizer_symbol_3d},
\[
    {\bm P}_{0}-\Id_{\bC^2}
    =
    \pare{1-\chi\pare{\xi}}
    \pare{
        \dot{\bm P}_{0}-\Id_{\bC^2}
    }.
\]
Since \(\dot{\bm P}_{0}-\Id_{\bC^2}\in\Gamma_{\geq1}^0\) vanishes at least linearly
in \(U\), and since \({\bm L}_{\rm lo}\) is a constant-coefficient order-zero
symbol, the pointwise conjugation gives
\[
    {\bm P}_{0}^{[-1]}
    {\bm L}_{\rm lo}
    {\bm P}_{0}
    =
    {\bm L}_{\rm lo}
    +
    {\bm A}_{\bra{0,\geq1}}.
\]
The Bony--Weyl correction terms contain either an \(x\)-derivative of
\({\bm P}_{0}\) or an \(x\)-derivative of \({\bm P}_{0}^{[-1]}\), hence again
at least one nonlinear factor.  Since all symbols involved have order \(0\),
these corrections are bounded.  Therefore
\begin{equation}\label{eq:principal_low_linear_defect_conjugation_3d}
    {\bm B}_{0}^{-1}
    \OpBW{{\bm L}_{\rm lo}}
    {\bm B}_{0}
    =
    \OpBW{
        {\bm A}^{\bC,{\rm lin}}_{\rm lo}
        +
        {\bm A}_{\bra{0,\geq1}}
    }
    +
    {\bm R}_{\bra{\varrho,\geq1}}.
\end{equation}
\end{proofpart}

\begin{proofpart}[The scalar transport block]
We apply \cref{eq:principal_B0_conjugation_formula_3d} to
\({\bm A}={\bm L}_1\).  By
\cref{eq:positive_order_diagonal_transport_symbol_3d},
\[
    {\bm L}_1
    =
    \ii\dot{\bm D}^{\bC}_{1}
    =
    \ii
    \pare{
        \mathpzc V_{\textnormal{can}}
        \pare{
            \Atw,\background
            \middle|
            \eta\pare{U},\tilde\psi\pare{U};x
        }
        \cdot\xi
    }
    \Id_{\bC^2}.
\]
It is scalar in the matrix indices, so its pointwise conjugation is exactly
unchanged on the whole frequency space:
\begin{equation}\label{eq:principal_transport_pointwise_conjugation_3d}
    {\bm P}_0^{[-1]}{\bm L}_1{\bm P}_0
    =
    {\bm L}_1.
\end{equation}
Every non-pointwise term in
\({\bm C}_{{\bm P}_0}\bra{{\bm L}_1}\) has order at most \(0\): the
first Bony--Weyl correction contains one \(\xi\)-derivative of an order-one
symbol, and all higher corrections are lower order.  Terms in which a
\(\xi\)-derivative hits the cut-off are supported in the fixed transition
annulus and are smoothing.  Therefore
\begin{equation}\label{eq:principal_transport_conjugated_result_3d}
    {\bm B}_0^{-1}
    \OpBW{{\bm L}_1}
    {\bm B}_0
    =
    \OpBW{
        \ii\dot{\bm D}^{\bC}_{1}
        +
        {\bm A}_{\bra{0,\geq1}}
    }
    +
    {\bm R}_{\bra{\varrho,\geq1}}.
\end{equation}
Thus the scalar transport block is unchanged at every positive order.
\end{proofpart}

\begin{proofpart}[The explicit order-\(1/2\) block]
By \cref{eq:principal_weyl_correction_vanishes_3d}, the order-\(3/2\)
block produces no homogeneous order-\(1/2\) correction.  Hence the only
surviving source at order \(1/2\) is
\[
    {\bm L}_{\frac12}
    =
    \pare{1-\chi\pare{\xi}}\,
    \dot{\bm A}^{\bC}_{0}\av{\xi}^{\frac12}.
\]
On the high-frequency region \(\chi=0\), one has
\({\bm P}_0=\dot{\bm P}_0\).  Since
\(\dot{\bm P}_0^{\pm1}\in(\Gamma^0_{\ge0})^{2\times2}\), the first
non-pointwise Bony--Weyl correction in the conjugation of
\({\bm L}_{\frac12}\) has order \(-1/2\).  Every derivative of the cut-off
is supported in the fixed transition annulus and is smoothing.  Therefore
\begin{equation}\label{eq:principal_half_conjugation_reduction_3d}
    {\bm C}_{{\bm P}_0}\bra{{\bm L}_{\frac12}}
    =
    \pare{1-\chi\pare{\xi}}\,
    \dot{\bm P}_0^{-1}
    \dot{\bm A}^{\bC}_{0}
    \dot{\bm P}_0
    \av{\xi}^{\frac12}
    +
    {\bm A}_{\bra{0,\geq1}}.
\end{equation}
By \cref{eq:complex_half_order_symbol_definition_3d},
\[
    \dot{\bm A}^{\bC}_{0}
    =
    \ii
    \begin{bmatrix}
        \dot\alpha^{\bC}_{0}
        &
        \dot\beta^{\bC}_{0}
        \\
        -\dot\beta^{\bC}_{0}
        &
        -\dot\alpha^{\bC}_{0}
    \end{bmatrix}.
\]
Using
\cref{eq:principal_diagonalizer_symbol_3d,eq:principal_diagonalizer_inverse_symbol_3d},
we obtain
\begin{equation}\label{eq:principal_half_pointwise_conjugation_computation_3d}
\begin{aligned}
    \dot{\bm P}_0^{-1}
    \dot{\bm A}^{\bC}_{0}
    \dot{\bm P}_0
    & =
    \frac{\ii}{1-\dot\chi_0^2}
    \begin{bmatrix}
        1&-\dot\chi_0
        \\
        -\dot\chi_0&1
    \end{bmatrix}
    \begin{bmatrix}
        \dot\alpha^{\bC}_{0}
        &
        \dot\beta^{\bC}_{0}
        \\
        -\dot\beta^{\bC}_{0}
        &
        -\dot\alpha^{\bC}_{0}
    \end{bmatrix}
    \begin{bmatrix}
        1&\dot\chi_0
        \\
        \dot\chi_0&1
    \end{bmatrix}
    \\
    & =
    \ii
    \begin{bmatrix}
        \dot\alpha^{\bC,{\rm pr}}_{0}
        &
        \dot\beta^{\bC,{\rm pr}}_{0}
        \\
        -\dot\beta^{\bC,{\rm pr}}_{0}
        &
        -\dot\alpha^{\bC,{\rm pr}}_{0}
    \end{bmatrix}.
\end{aligned}
\end{equation}
where $\dot\alpha^{\bC,{\rm pr}}_{0}$ and $ \dot\beta^{\bC,{\rm pr}}_{0}$ are defined in 
\cref{eq:principal_alpha_principal_half_3d,eq:principal_beta_principal_half_3d}.
Thus
\begin{equation}\label{eq:principal_half_conjugated_result_3d}
    {\bm B}_0^{-1}
    \OpBW{{\bm L}_{\frac12}}
    {\bm B}_0
    =
    \OpBW{
        \pare{1-\chi\pare{\xi}}\,
        \dot{\bm A}^{\bC,{\rm pr}}_0
        \av{\xi}^{\frac12}
        +
        {\bm A}_{\bra{0,\geq1}}
    }
    +
    {\bm R}_{\bra{\varrho,\geq1}}.
\end{equation}
This is the unique homogeneous order-\(1/2\) block in the
principal-diagonalized system.
\end{proofpart}

\begin{proofpart}[Bounded symbols, time derivative, and smoothing remainders]
For \({\bm L}_0={\bm A}^{\bC}_{\bra{0,\geq1}}\),
\cref{eq:principal_B0_conjugation_formula_3d} gives
\begin{equation}\label{eq:principal_bounded_conjugation_result_3d}
    {\bm B}_0^{-1}
    \OpBW{{\bm L}_0}
    {\bm B}_0
    =
    \OpBW{{\bm A}_{\bra{0,\geq1}}}
    +
    {\bm R}_{\bra{\varrho,\geq1}}.
\end{equation}
Indeed the three factors are order-zero operators and all symbolic corrections
are of order at most \(0\).

We now treat the time derivative in
\cref{eq:principal_diag_time_derivative_3d}.  Since
\({\bm B}_0=\OpBW{{\bm P}_0}\), the symbolic composition formula gives
\begin{equation}\label{eq:principal_time_derivative_symbol_3d}
    {\bm B}_{0}^{-1}
    \pare{\partial_t{\bm B}_{0}}
    =
    \OpBW{
        {\bm P}_0^{[-1]}\partial_t{\bm P}_0
        +
        {\bm A}_{\bra{-1,\geq1}}
    }
    +
    {\bm R}_{\bra{\varrho,\geq1}}.
\end{equation}
The cut-off is time-independent and
\(\partial_t{\bm P}_0=(1-\chi)\partial_t\dot{\bm P}_0\).  Hence
\(\partial_t{\bm P}_0\) is of order zero in \(\xi\): the time derivative
acts only on the coefficients depending on \(U\) and does not alter their
high-frequency homogeneity.  Its tame bounds follow, after increasing
\(s_0\), from the equation for \(U\), from
\cref{eq:dote0,eq:dotf0,eq:principal_scalar_l_theta_3d,eq:dotchi0}, and from
\cref{lem:composition_rules}-\ref{item:composition_rules_substitution_symbol}.
Thus
\begin{equation}\label{eq:principal_time_derivative_bounded_result_3d}
    {\bm B}_{0}^{-1}
    \pare{\partial_t{\bm B}_{0}}
    =
    \OpBW{{\bm A}_{\bra{0,\geq1}}}
    +
    {\bm R}_{\bra{\varrho,\geq1}},
\end{equation}
and this contribution is absorbed into
\({\bm A}^{\bC,{\rm pr}}_{\bra{0,\geq1}}\).

Finally, by
\cref{lem:composition_rules}-\ref{item:composition_rules_paradiff_remainder},
\begin{equation}\label{eq:principal_smoothing_conjugation_result_3d}
    {\bm B}_{0}^{-1}
    {\bm R}^{\bC}_{\bra{\varrho,\geq1}}
    {\bm B}_{0}
    \in
    \pare{\cR^\varrho}^{2\times2},
\end{equation}
after decreasing and relabelling \(\varrho\).
\end{proofpart}

Collecting
\cref{eq:principal_low_linear_defect_conjugation_3d,eq:principal_block_conjugated_result_3d,eq:principal_transport_conjugated_result_3d,eq:principal_half_conjugated_result_3d,eq:principal_bounded_conjugation_result_3d,eq:principal_time_derivative_bounded_result_3d,eq:principal_smoothing_conjugation_result_3d}
inside \cref{eq:principal_diag_time_derivative_3d}, and using
\cref{eq:principal_weyl_correction_vanishes_3d}, the only homogeneous
order-\(1/2\) contribution is the pointwise-conjugated block defined in
\cref{eq:principal_half_output_decomposition_3d}.  Absorbing all bounded
symbols into \({\bm A}^{\bC,{\rm pr}}_{\bra{0,\geq1}}\) and all smoothing
contributions into \({\bm R}^{\bC,{\rm pr}}_{\bra{\varrho,\geq1}}\), we obtain
\[
\begin{aligned}
    \partial_t U_0
    ={}&
    \OpBW{{\bm A}^{\bC,{\rm lin}}_{\rm lo}}U_0
    +
    \OpBW{\ii{\bm D}^{\bC}_{\frac32}}U_0
    +
    \OpBW{\ii\dot{\bm D}^{\bC}_{1}}U_0
    \\
    &+
    \OpBW{
        \pare{1-\chi\pare{\xi}}\,
        \dot{\bm A}^{\bC,{\rm pr}}_0
        \av{\xi}^{\frac12}
    }U_0
    +
    \OpBW{{\bm A}^{\bC,{\rm pr}}_{\bra{0,\geq1}}}U_0
    +
    {\bm R}^{\bC,{\rm pr}}_{\bra{\varrho,\geq1}}U_0.
\end{aligned}
\]
Restoring the common symbol arguments
\[
    \pare{\Atw,\background\middle|U;x,\xi},
\]
and the multiplier arguments
\(\pare{\Atw,\background\middle|\xi}\), gives precisely
\cref{eq:principal_diagonalized_system_3d}.
\end{proof}

\subsubsection{Diagonalization at subprincipal order}
\label{sec:subprincipal_order_diagonalization_3d}

Here we perform the diagonalization of the lower positive order terms. The result reads as follows. 
\begin{prop}[Subprincipal diagonalization]
\label{prop:subprincipal_order_diagonalization_3d}
Let the assumptions of
\cref{prop:positive_order_diagonalization_3d} hold, and let \(U_0\) solve
\cref{eq:principal_diagonalized_system_3d}.  There exists a bounded and
invertible operator
\[
    {\bm B}_{-1}
    \pare{\Atw,\background\middle|U}
    \in
    \cL\pare{H^s}
    \qquad
    \forall\,s\geq s_0,
\]
such that, defining
\begin{equation}\label{eq:subprincipal_diagonalized_unknown_3d}
    W
    \defeq
    {\bm B}_{-1}
    \pare{\Atw,\background\middle|U}^{-1}U_0 ,
\end{equation}
the unknown \(W\) solves
\begin{equation}\label{eq:subprincipal_diagonalized_system_3d}\begin{aligned}
    \partial_t W
    =
    &\OpBW{
        {\bm A}^{\bC,{\rm lin}}_{\rm lo}
        \pare{\Atw,\background\middle|\xi}
    }W
    \\
    +
    &\ii
    \OpBW{
        {\bm D}^{\bC}_{\frac32}
        \pare{\Atw,\background\middle|U;x,\xi}
        +
        \dot{\bm D}^{\bC}_{1}
        \pare{\Atw,\background\middle|U;x,\xi}
        +
        \dot{\bm D}^{\bC}_{\frac12}
        \pare{\Atw,\background\middle|U;x,\xi}
        +
        {\bm A}^{\bC,{\rm d}}_{\bra{0,\geq1}}
        \pare{\Atw,\background\middle|U;x,\xi}
    }W
    \\
    +
    &{\bm R}^{\bC,{\rm d}}_{\bra{\varrho,\geq1}}
    \pare{\Atw,\background\middle|U}W .
\end{aligned}
\end{equation}
The objects entering
\cref{eq:subprincipal_diagonalized_unknown_3d,eq:subprincipal_diagonalized_system_3d}
are defined as follows.
\begin{itemize}
    \item
    \({\bm B}_{-1}\pare{\Atw,\background\middle|U}\) is the order-zero
    perturbation of the identity
    \begin{equation}\label{eq:subprincipal_diagonalizer_operator_3d}
        {\bm B}_{-1}
        \pare{\Atw,\background\middle|U}
        \defeq
        \Id_{\bC^2}
        +
        \OpBW{
            {\bm S}_{-1}
            \pare{\Atw,\background\middle|U;x,\xi}
        } .
    \end{equation}
    We use the fixed off-diagonal matrix
    \({\bm X}_{\rm od}\) introduced in
    \cref{eq:subprincipal_offdiagonal_matrix_3d}.  The final
    non-homogeneous generator is
    \begin{equation}\label{eq:subprincipal_nonhomogeneous_generator_symbol_3d}
        {\bm S}_{-1}
        \pare{\Atw,\background\middle|U;x,\xi}
        \defeq
        \pare{1-\chi\pare{\xi}}\,
        \dot{\bm S}_{-1}
        \pare{\Atw,\background\middle|U;x,\xi}
        \in
        \pare{\Gamma_{\geq1}^{-1}}^{2\times2},
    \end{equation}
    where the high-frequency homogeneous generator is
    \begin{equation}\label{eq:subprincipal_generator_symbol_3d}
        \dot{\bm S}_{-1}
        \pare{\Atw,\background\middle|U;x,\xi}
        \defeq
        \frac{
            \dot\beta^{\bC,{\rm pr}}_{0}
            \pare{\Atw,\background\middle|U;x,\hat\xi}
            \av{\xi}^{\frac12}
        }{
            2
            \dot\omega_{\frac32}
            \pare{\Atw,\background\middle|\xi}
            \dot\theta_0
            \pare{\Atw,\background\middle|U;x,\hat\xi}
        }
        {\bm X}_{\rm od}.
    \end{equation}
    Here
    \({\bm X}_{\rm od}\) is the off-diagonal matrix defined in
    \cref{eq:subprincipal_offdiagonal_matrix_3d},
    \(\dot\omega_{\frac32}\) is defined in
    \cref{eq:complex_leading_capillary_frequency_3d},
    \(\dot\theta_0\) is defined in
    \cref{eq:principal_scalar_l_theta_3d}, and
    \(\dot\beta^{\bC,{\rm pr}}_0\) is defined in
    \cref{eq:principal_beta_principal_half_3d};

    \item
    \({\bm A}^{\bC,{\rm lin}}_{\rm lo}\) is the same flat low-frequency
    linear multiplier from \cref{eq:linear_low_frequency_defect_complex_3d};

    \item
    the diagonal symbols
    \[
        {\bm D}^{\bC}_{\frac32},
        \qquad
        \dot{\bm D}^{\bC}_{1},
        \qquad
        \dot{\bm D}^{\bC}_{\frac12}
    \]
    are those defined in
    \cref{eq:positive_order_diagonal_ql_symbol_3d,eq:positive_order_diagonal_transport_symbol_3d,eq:positive_order_diagonal_half_symbol_3d}.
    They belong respectively to
    \[
        \pare{\Gamma_{\geq0}^{\frac32}}^{2\times2},
        \qquad
        \pare{\Gamma_{\geq0}^1}^{2\times2},
        \qquad
        \pare{\Gamma_{\geq1}^{\frac12}}^{2\times2},
    \]
    and are real-valued diagonal matrix symbols for real \(U\);

    \item
    \({\bm A}^{\bC,{\rm d}}_{\bra{0,\geq1}}
    \pare{\Atw,\background\middle|U;x,\xi}
    \in\pare{\Gamma_{\geq1}^0}^{2\times2}\)
    is a generic bounded matrix symbol in the sense of
    \cref{notation:paralinearization};

    \item
    \({\bm R}^{\bC,{\rm d}}_{\bra{\varrho,\geq1}}
    \pare{\Atw,\background\middle|U}
    \in\pare{\cR_{\geq1}^\varrho}^{2\times2}\)
    is a generic smoothing remainder in the sense of
    \cref{notation:paralinearization}.
\end{itemize}
\end{prop}

\begin{proof}
Throughout the proof all symbols are evaluated at
\[
    \pare{\Atw,\background\middle|U;x,\xi},
\]
unless explicitly stated otherwise; Fourier multipliers are evaluated at
\(\pare{\Atw,\background\middle|\xi}\).  
%We use the convention of
%\cref{notation:paralinearization}: the symbols
%\({\bm A}_{\bra{m,\geq1}}\) and the remainders
%\({\bm R}_{\bra{\varrho,\geq1}}\) are generic and may change from line to line.  
The
symbolic calculus is used componentwise, through
\cref{def:asympt_symbol,prop:composition_paradiff}, and the smoothing
remainder calculus is used through
\cref{lem:composition_rules}-\ref{item:composition_rules_paradiff_remainder}
and
\cref{lem:composition_rules}-\ref{item:composition_rules_remainder_remainder}.

\begin{proofpart}[The subprincipal diagonalizer and its inverse]
We first check that the symbol defining the change of variables is of order
\(-1\) on the high-frequency region.  By
\cref{eq:principal_beta_principal_half_3d},
\[
    \dot\beta^{\bC,{\rm pr}}_0\in\Gamma_{\geq1}^0 .
\]
Moreover, by \cref{eq:principal_scalar_l_theta_3d}, after increasing
\(s_0\) and shrinking \(r_0\),
\[
    \dot\theta_0-1\in\Gamma_{\geq1}^0,
    \qquad
    \dot\theta_0\geq \frac34,
    \qquad
    \dot\theta_0^{-1}\in\Gamma^0_{\ge0}.
\]
Since
\[
    \dot\omega_{\frac32}^{-1}\av{\xi}^{\frac12}
    \in\tilde\Gamma^{-1}
\]
by \cref{eq:complex_leading_capillary_frequency_3d}, the product rule for
\(\Gamma_{\ge0}\)-symbols and
\cref{lem:composition_rules}-\ref{item:composition_rules_substitution_symbol}
give
\begin{equation}\label{eq:subprincipal_scalar_sminusone_class_3d}
    \dot s_{-1}
    \defeq
    \frac{
        \dot\beta^{\bC,{\rm pr}}_{0}\av{\xi}^{\frac12}
    }{
        2\dot\omega_{\frac32}\dot\theta_0
    }
    \in\Gamma_{\geq1}^{-1}.
\end{equation}
The coefficient \(\dot s_{-1}\) is real-valued for real \(U\).  Therefore
\[
    \dot{\bm S}_{-1}
    =
    \dot s_{-1}{\bm X}_{\rm od}
    \in
    \pare{\Gamma_{\geq1}^{-1}}^{2\times2}.
\]
Multiplying by the fixed high-frequency cut-off gives
\begin{equation}\label{eq:subprincipal_Sminusone_class_refined_3d}
    {\bm S}_{-1}
    =
    \pare{1-\chi\pare{\xi}}\dot{\bm S}_{-1}
    \in
    \pare{\Gamma_{\geq1}^{-1}}^{2\times2}.
\end{equation}
In particular, \({\bm S}_{-1}\) is a non-homogeneous order-\(\pare{-1}\) symbol:
it vanishes at low frequency and coincides with the homogeneous generator
\(\dot{\bm S}_{-1}\) where \(1-\chi\pare{\xi}=1\).

By the action estimate for paradifferential operators,
\(\OpBW{{\bm S}_{-1}}\) is bounded on \(H^s\), for \(s\geq s_0\).  Since
\({\bm S}_{-1}\) vanishes at least linearly in \(U\), its low Sobolev
seminorm is \(O(r_0)\).  By \cref{thm:contS}, after decreasing \(r_0\), we
get
\[
    \norm{\OpBW{{\bm S}_{-1}}}_{\cL(H^s)}<1,
    \qquad s\ge s_0.
\]
Hence
\begin{equation}\label{eq:subprincipal_Bminusone_invertible_3d}
    {\bm B}_{-1}
    =
    \Id_{\bC^2}+\OpBW{{\bm S}_{-1}}
\end{equation}
is bounded and invertible on \(H^s\) by the Neumann series.

We next record the symbolic expansion of the inverse.  Since we have that
\({\bm S}_{-1}\in(\Gamma_{\geq1}^{-1})^{2\times2}\), then by virtue of
\cref{def:asympt_symbol,prop:composition_paradiff}, we infer
\begin{equation}\label{eq:subprincipal_parametrix_product_3d}
\begin{aligned}
    \pare{\Id_{\bC^2}-\OpBW{{\bm S}_{-1}}}
    \pare{\Id_{\bC^2}+\OpBW{{\bm S}_{-1}}}
    &=
    \Id_{\bC^2}
    -
    \OpBW{{\bm S}_{-1}\#_\varrho{\bm S}_{-1}}
    +
    {\bm R}_{\bra{\varrho,\geq1}}
    \\
    &=
    \Id_{\bC^2}
    +
    \OpBW{{\bm A}_{\bra{-2,\geq1}}}
    +
    {\bm R}_{\bra{\varrho,\geq1}}.
\end{aligned}
\end{equation}
The same identity holds with the two factors exchanged.  Combining the
invertibility statement \cref{eq:subprincipal_Bminusone_invertible_3d} with
the parametrix identity \cref{eq:subprincipal_parametrix_product_3d}, and
using
\cref{lem:composition_rules}-\ref{item:composition_rules_paradiff_remainder},
we may write the true inverse as
\begin{equation}\label{eq:subprincipal_Bm1_true_inverse_expansion_refined_3d}
    {\bm B}_{-1}^{-1}
    =
    \Id_{\bC^2}
    -
    \OpBW{{\bm S}_{-1}}
    +
    \OpBW{{\bm A}_{\bra{-2,\geq1}}}
    +
    {\bm R}_{\bra{\varrho,\geq1}},
\end{equation}
after decreasing and relabelling \(\varrho\).
\end{proofpart}

\begin{proofpart}[Symbolic conjugation formula]
Let \({\bm A}\in(\Gamma^m_{\ge0})^{2\times2}\), \(m\leq3/2\), be any of
the symbols to be conjugated below.  By
\cref{eq:subprincipal_Bm1_true_inverse_expansion_refined_3d},
\cref{def:asympt_symbol,prop:composition_paradiff}, and
\cref{lem:composition_rules}-\ref{item:composition_rules_paradiff_remainder},
\begin{equation}\label{eq:subprincipal_operator_conjugation_formula_refined_3d}
    {\bm B}_{-1}^{-1}\OpBW{{\bm A}}{\bm B}_{-1}
    =
    \OpBW{
        {\bm A}
        +
        {\bm A}\#_\varrho{\bm S}_{-1}
        -
        {\bm S}_{-1}\#_\varrho{\bm A}
    }
    +
    \OpBW{{\bm A}_{\bra{m-2,\geq1}}}
    +
    {\bm R}_{\bra{\varrho,\geq1}}.
\end{equation}
The pointwise products
\({\bm A}{\bm S}_{-1}\) and
\({\bm S}_{-1}{\bm A}\) have order \(m-1\).  The first non-pointwise
Bony--Weyl correction contains one \(\xi\)-derivative falling either on
\({\bm A}\) or on \({\bm S}_{-1}\), and therefore has order \(m-2\).
Consequently, at positive order,
\begin{equation}\label{eq:subprincipal_symbolic_conjugation_positive_part_refined_3d}
    {\bm B}_{-1}^{-1}\OpBW{{\bm A}}{\bm B}_{-1}
    =
    \OpBW{
        {\bm A}
        +
        {\bm A}{\bm S}_{-1}
        -
        {\bm S}_{-1}{\bm A}
        +
        {\bm A}_{\bra{m-2,\geq1}}
    }
    +
    {\bm R}_{\bra{\varrho,\geq1}}.
\end{equation}
For \(m=3/2\) the discarded terms have order at most \(-1/2\), and for
\(m\le1\) they have order at most \(-1\).  All such contributions are
bounded in the non-homogeneous scale and are absorbed into a generic
\({\bm A}_{\bra{0,\geq1}}\).  If a \(\xi\)-derivative hits the cut-off
\(1-\chi\pare{\xi}\), the resulting contribution is compactly supported in
frequency and is again absorbed into \({\bm A}_{\bra{0,\geq1}}\).
\end{proofpart}

\begin{proofpart}[Conjugation of the equation]
Set
\[
    U_0={\bm B}_{-1}W .
\]
Differentiating this identity gives
\begin{equation}\label{eq:subprincipal_time_derivative_refined_3d}
    \partial_t W
    =
    {\bm B}_{-1}^{-1}\partial_tU_0
    -
    {\bm B}_{-1}^{-1}
    \pare{\partial_t{\bm B}_{-1}}W .
\end{equation}
We write \cref{eq:principal_diagonalized_system_3d} as
\begin{equation}\label{eq:subprincipal_input_system_compact_refined_3d}
    \partial_tU_0
    =
    \OpBW{
{\bm L}_{\rm lo}
+
{\bm L}_{\frac32}
+
{\bm L}_{1}
+
{\bm L}_{\frac12}
+
{\bm L}_{0}
    }U_0
    +
    {\bm R}^{\bC,{\rm pr}}_{\bra{\varrho,\geq1}}U_0,
\end{equation}
where
\begin{align}
    {\bm L}_{\rm lo}
    &\defeq
    {\bm A}^{\bC,{\rm lin}}_{\rm lo}
    \pare{\Atw,\background\middle|\xi},
    \label{eq:subprincipal_Llo_definition_refined_3d}
\\
    {\bm L}_{\frac32}
    &\defeq
    \ii{\bm D}^{\bC}_{\frac32},
    \label{eq:subprincipal_L32_definition_refined_3d}
    \\
    {\bm L}_{1}
    &\defeq
    \ii\dot{\bm D}^{\bC}_{1},
    \label{eq:subprincipal_L1_definition_refined_3d}
    \\
    {\bm L}_{\frac12}
    &\defeq
    \pare{1-\chi\pare{\xi}}\,
    \dot{\bm A}^{\bC,{\rm pr}}_{0}
    \av{\xi}^{\frac12},
    \label{eq:subprincipal_L12_definition_refined_3d}
    \\
    {\bm L}_{0}
    &\defeq
    {\bm A}^{\bC,{\rm pr}}_{\bra{0,\geq1}}.
    \label{eq:subprincipal_L0_definition_refined_3d}
\end{align}
Moreover $ {\bm L}_{\rm lo}
    \in
    \pare{\Gamma^0_{\ge0}}^{2\times2}$ is constant-coefficient and does not vanish at \(U=0\).  It will be kept
explicitly through the conjugation.  The bounded nonlinear defects generated
by its commutators with the diagonalizer are absorbed into
\({\bm A}_{\bra{0,\geq1}}\).
The transport symbol \(\dot{\bm D}^{\bC}_{1}\) is the scalar diagonal
symbol defined in
\cref{eq:positive_order_diagonal_transport_symbol_3d}.  The factors
\(1-\chi\pare{\xi}\) in
\({\bm D}^{\bC}_{\frac32}\), \({\bm L}_{\frac12}\), and
\(\dot{\bm D}^{\bC}_{\frac12}\) encode that the subprincipal
diagonalization is performed only at high frequency.
\end{proofpart}

\begin{proofpart}[The flat low-frequency linear block]
The subprincipal diagonalizer has the form
\[
    {\bm B}_{-1}
    =
    \Id_{\bC^2}
    +
    \OpBW{{\bm S}_{-1}},
    \qquad
    {\bm S}_{-1}
    =
    \pare{1-\chi\pare{\xi}}\dot{\bm S}_{-1}.
\]
In particular \({\bm S}_{-1}\) vanishes at least linearly in \(U\).  Since
\({\bm L}_{\rm lo}={\bm A}^{\bC,{\rm lin}}_{\rm lo}\) is a
constant-coefficient order-zero multiplier, the conjugation formula gives
\[
    {\bm L}_{\rm lo}
    +
    {\bm L}_{\rm lo}{\bm S}_{-1}
    -
    {\bm S}_{-1}{\bm L}_{\rm lo}
    +
    {\bm A}_{\bra{-1,\geq1}}.
\]
The commutator term has order \(-1\) and vanishes at least linearly in \(U\).
It is therefore included in the bounded nonlinear symbol
\({\bm A}_{\bra{0,\geq1}}\).  Thus
\begin{equation}\label{eq:subprincipal_low_linear_defect_conjugation_3d}
    {\bm B}_{-1}^{-1}
    \OpBW{{\bm L}_{\rm lo}}
    {\bm B}_{-1}
    =
    \OpBW{
        {\bm A}^{\bC,{\rm lin}}_{\rm lo}
        +
        {\bm A}_{\bra{0,\geq1}}
    }
    +
    {\bm R}_{\bra{\varrho,\geq1}}.
\end{equation}
\end{proofpart}

\begin{proofpart}[The order-\(3/2\) block and the subprincipal correction]
We compute the order-\(1/2\) term produced by conjugating the already
diagonal order-\(3/2\) block with the subprincipal diagonalizer.  This is
the term designed to cancel the off-diagonal part of the explicit
order-\(1/2\) block.  First,
\begin{equation}\label{eq:Jd_Xod_commutator_refined_3d}
    {\bm J}_{\rm d}{\bm X}_{\rm od}
    -
    {\bm X}_{\rm od}{\bm J}_{\rm d}
    =
    2{\bm J}.
\end{equation}
By
\cref{eq:subprincipal_L32_definition_refined_3d,eq:subprincipal_nonhomogeneous_generator_symbol_3d,eq:subprincipal_generator_symbol_3d},
\[
    {\bm L}_{\frac32}
    =
    \ii\pare{1-\chi\pare{\xi}}\,
    \omega_{\frac32}\dot\theta_0{\bm J}_{\rm d},
    \qquad
    {\bm S}_{-1}
    =
    \pare{1-\chi\pare{\xi}}\,
    \dot s_{-1}{\bm X}_{\rm od},
\]
where
\[
    \dot s_{-1}
    =
    \frac{
        \dot\beta^{\bC,{\rm pr}}_0\av{\xi}^{\frac12}
    }{
        2\dot\omega_{\frac32}\dot\theta_0
    }\cdot
\]
The pointwise commutator in
\cref{eq:subprincipal_symbolic_conjugation_positive_part_refined_3d} is
therefore
\begin{equation}\label{eq:subprincipal_principal_pointwise_commutator_refined_3d}
\begin{aligned}
    {\bm L}_{\frac32}{\bm S}_{-1}
    -
    {\bm S}_{-1}{\bm L}_{\frac32}
    & =
    \ii\pare{1-\chi\pare{\xi}}^2
    \omega_{\frac32}\dot\theta_0\dot s_{-1}
    \pare{
        {\bm J}_{\rm d}{\bm X}_{\rm od}
        -
        {\bm X}_{\rm od}{\bm J}_{\rm d}
    }
    \\
    & =
    \ii\pare{1-\chi\pare{\xi}}^2
    \frac{\omega_{\frac32}}{\dot\omega_{\frac32}}
    \dot\beta^{\bC,{\rm pr}}_0
    \av{\xi}^{\frac12}{\bm J}.
\end{aligned}
\end{equation}
By \cref{eq:complex_leading_multiplier_expansions_3d},
\(\omega_{\frac32}/\dot\omega_{\frac32}-1\in\tilde\Gamma^{-1}\).
Hence its contribution after multiplication by
\(\dot\beta^{\bC,{\rm pr}}_0\av{\xi}^{1/2}\) has order \(-1/2\).
Moreover,
\(\pare{1-\chi}^2-\pare{1-\chi}\) is supported in the fixed transition
annulus and is smoothing.  Therefore
\begin{equation}\label{eq:subprincipal_principal_commutator_result_refined_3d}
    {\bm L}_{\frac32}{\bm S}_{-1}
    -
    {\bm S}_{-1}{\bm L}_{\frac32}
    =
    \ii\pare{1-\chi\pare{\xi}}\,
    \dot\beta^{\bC,{\rm pr}}_0
    \av{\xi}^{\frac12}{\bm J}
    +
    {\bm A}_{\bra{-\frac12,\geq1}}.
\end{equation}
The first non-pointwise Bony--Weyl term in
\({\bm L}_{\frac32}\#_\varrho{\bm S}_{-1}
-{\bm S}_{-1}\#_\varrho{\bm L}_{\frac32}\) also has order
\(3/2-1-1=-1/2\), and all subsequent terms are lower order.  Combining
\cref{eq:subprincipal_symbolic_conjugation_positive_part_refined_3d,eq:subprincipal_principal_commutator_result_refined_3d},
we obtain
\begin{equation}\label{eq:subprincipal_L32_conjugated_result_refined_3d}
    {\bm B}_{-1}^{-1}
    \OpBW{{\bm L}_{\frac32}}
    {\bm B}_{-1}
    =
    \OpBW{
        \ii{\bm D}^{\bC}_{\frac32}
        +
        \ii\pare{1-\chi\pare{\xi}}\,
        \dot\beta^{\bC,{\rm pr}}_0
        \av{\xi}^{\frac12}{\bm J}
        +
        {\bm A}_{\bra{0,\geq1}}
    }
    +
    {\bm R}_{\bra{\varrho,\geq1}}.
\end{equation}
Here \({\bm A}_{\bra{0,\geq1}}\) absorbs all order \(-1/2\) and lower terms.
\end{proofpart}

\begin{proofpart}[The explicit order-\(1/2\) block]
We now treat the explicit order-\(1/2\) block in the principal-diagonalized
system.  By
\cref{eq:principal_half_output_decomposition_3d,eq:positive_order_diagonal_half_symbol_3d},
\begin{equation}\label{eq:subprincipal_half_block_decomposition_refined_3d}
    {\bm L}_{\frac12}
    =
    \pare{1-\chi\pare{\xi}}\,
    \dot{\bm A}^{\bC,{\rm pr}}_0
    \av{\xi}^{\frac12}
    =
    \ii\dot{\bm D}^{\bC}_{\frac12}
    -
    \ii
    \pare{1-\chi\pare{\xi}}\,
    \dot\beta^{\bC,{\rm pr}}_0
    \av{\xi}^{\frac12}
    {\bm J}.
\end{equation}
Here the coefficient \(\dot\beta^{\bC,{\rm pr}}_0\) comes from the
pointwise conjugation of the explicit order-\(1/2\) block.  The
cancellation in \cref{eq:principal_weyl_correction_vanishes_3d} shows that
the order-\(3/2\) block adds no homogeneous order-\(1/2\) contribution.

Since
\({\bm L}_{\frac12}\in(\Gamma_{\geq1}^{1/2})^{2\times2}\) and
\({\bm S}_{-1}\in(\Gamma_{\geq1}^{-1})^{2\times2}\), the pointwise commutator
\[
    {\bm L}_{\frac12}{\bm S}_{-1}
    -
    {\bm S}_{-1}{\bm L}_{\frac12}
\]
has order \(-1/2\), and the first non-pointwise Bony--Weyl corrections have
even lower order.  Therefore the subprincipal diagonalizer does not modify
the order-\(1/2\) part of \({\bm L}_{\frac12}\).  Using
\cref{eq:subprincipal_symbolic_conjugation_positive_part_refined_3d}, we get
\begin{equation}\label{eq:subprincipal_L12_conjugated_result_refined_3d}
    {\bm B}_{-1}^{-1}
    \OpBW{{\bm L}_{\frac12}}
    {\bm B}_{-1}
    =
    \OpBW{
        \ii\dot{\bm D}^{\bC}_{\frac12}
        -
        \ii
        \pare{1-\chi\pare{\xi}}\,
        \dot\beta^{\bC,{\rm pr}}_0
        \av{\xi}^{\frac12}
        {\bm J}
        +
        {\bm A}_{\bra{0,\geq1}}
    }
    +
    {\bm R}_{\bra{\varrho,\geq1}}.
\end{equation}

The off-diagonal term in
\cref{eq:subprincipal_L12_conjugated_result_refined_3d} cancels exactly with
the off-diagonal order-\(1/2\) term generated by the conjugation of
\({\bm L}_{\frac32}\), computed in
\cref{eq:subprincipal_L32_conjugated_result_refined_3d}.  Therefore
\begin{equation}\label{eq:subprincipal_positive_half_diagonal_result_refined_3d}
\begin{aligned}
    &{\bm B}_{-1}^{-1}
    \OpBW{
        {\bm L}_{\frac32}
        +
        {\bm L}_{\frac12}
    }
    {\bm B}_{-1}
    =
    \OpBW{
        \ii{\bm D}^{\bC}_{\frac32}
        +
        \ii\dot{\bm D}^{\bC}_{\frac12}
        +
        {\bm A}_{\bra{0,\geq1}}
    }
    +
    {\bm R}_{\bra{\varrho,\geq1}}.
\end{aligned}
\end{equation}
\end{proofpart}

\begin{proofpart}[The scalar transport block]
We now conjugate
\({\bm L}_{1}=\ii\dot{\bm D}^{\bC}_{1}\). By
\cref{eq:positive_order_diagonal_transport_symbol_3d},
\[
    \dot{\bm D}^{\bC}_{1}
    =
    \pare{
        \mathpzc V_{\textnormal{can}}
        \pare{\Atw,\background\middle|\eta\pare U,\tilde\psi\pare U;x}
        \cdot\xi
    }\Id_{\bC^2}.
\]
Thus the pointwise commutator with \({\bm S}_{-1}\) vanishes. The first
non-pointwise Bony--Weyl commutator has order
$1+(-1)-1=-1$, and the defects produced by the
\({\bm A}_{\bra{-2,\geq1}}\)-part of the inverse have order \(1-2=-1\). Hence
\begin{equation}\label{eq:subprincipal_transport_conjugated_result_refined_3d}
    {\bm B}_{-1}^{-1}
    \OpBW{{\bm L}_{1}}
    {\bm B}_{-1}
    =
    \OpBW{
        \ii\dot{\bm D}^{\bC}_{1}
        +
        {\bm A}_{\bra{0,\geq1}}
    }
    +
    {\bm R}_{\bra{\varrho,\geq1}}.
\end{equation}
Thus the subprincipal diagonalizer does not change the scalar transport
block at positive order.
\end{proofpart}

\begin{proofpart}[Bounded terms, time derivative, and smoothing remainders]
Since
\({\bm L}_{0}={\bm A}^{\bC,{\rm pr}}_{\bra{0,\geq1}}\in(\Gamma_{\geq1}^0)^{2\times2}\),
\cref{eq:subprincipal_operator_conjugation_formula_refined_3d} gives
\begin{equation}\label{eq:subprincipal_bounded_conjugated_result_refined_3d}
    {\bm B}_{-1}^{-1}
    \OpBW{{\bm L}_{0}}
    {\bm B}_{-1}
    =
    \OpBW{{\bm A}_{\bra{0,\geq1}}}
    +
    {\bm R}_{\bra{\varrho,\geq1}}.
\end{equation}
Indeed, all additional products with \({\bm S}_{-1}\) have order at most
\(-1\), hence are bounded.

We next treat the time derivative in
\cref{eq:subprincipal_time_derivative_refined_3d}. Since the Bony--Weyl
cut-off is time-independent,
\[
    \partial_t{\bm B}_{-1}
    =
    \OpBW{\partial_t{\bm S}_{-1}}.
\]
The factor \(1-\chi\pare{\xi}\) is independent of time, and differentiating
the high-frequency homogeneous coefficient
\(\dot{\bm S}_{-1}\) along the solution does not change the \(\xi\)-order.
After replacing \(\partial_tU\) by the vector field in
\cref{eq:complex_paradifferential_system_3d}, and increasing
\(s_0\) if necessary to absorb the fixed loss of derivatives in the tame
bounds, we may regard
\begin{equation}\label{eq:subprincipal_time_derivative_symbol_class_refined_3d}
    \partial_t{\bm S}_{-1}
    \in
    \pare{\Gamma_{\geq1}^{-1}}^{2\times2}
    \subset
    \pare{\Gamma_{\geq1}^0}^{2\times2}
\end{equation}
as a bounded matrix symbol depending on \(U\) and on the fixed parameters.
Therefore
\begin{equation}\label{eq:subprincipal_time_derivative_symbol_refined_3d}
    {\bm B}_{-1}^{-1}
    \pare{\partial_t{\bm B}_{-1}}
    =
    \OpBW{{\bm A}_{\bra{0,\geq1}}}
    +
    {\bm R}_{\bra{\varrho,\geq1}}.
\end{equation}

Finally, by
\cref{lem:composition_rules}-\ref{item:composition_rules_paradiff_remainder}
and
\cref{lem:composition_rules}-\ref{item:composition_rules_remainder_remainder},
\begin{equation}\label{eq:subprincipal_smoothing_conjugated_result_refined_3d}
    {\bm B}_{-1}^{-1}
    {\bm R}^{\bC,{\rm pr}}_{\bra{\varrho,\geq1}}
    {\bm B}_{-1}
    \in
    \pare{\cR_{\geq1}^\varrho}^{2\times2},
\end{equation}
after decreasing and relabeling \(\varrho\). Indeed, the identity parts of
\({\bm B}_{-1}^{\pm1}\) preserve the smoothing order, whereas every factor
\(\OpBW{{\bm S}_{-1}}\) improves it, because \({\bm S}_{-1}\) has order
\(-1\).
\end{proofpart}

\begin{proofpart}[Conclusion]
Substituting
\cref{eq:subprincipal_low_linear_defect_conjugation_3d,eq:subprincipal_positive_half_diagonal_result_refined_3d,eq:subprincipal_transport_conjugated_result_refined_3d,eq:subprincipal_bounded_conjugated_result_refined_3d,eq:subprincipal_time_derivative_symbol_refined_3d,eq:subprincipal_smoothing_conjugated_result_refined_3d}
into \cref{eq:subprincipal_time_derivative_refined_3d}, we obtain
\[
    \partial_t W
    =
    \OpBW{
        {\bm A}^{\bC,{\rm lin}}_{\rm lo}
    }W
    +
    \OpBW{
        \ii{\bm D}^{\bC}_{\frac32}
        +
        \ii\dot{\bm D}^{\bC}_{1}
        +
        \ii\dot{\bm D}^{\bC}_{\frac12}
        +
        {\bm A}_{\bra{0,\geq1}}
    }W
    +
    {\bm R}^{\bC,{\rm d}}_{\bra{\varrho,\geq1}}W.
\]
Since \({\bm A}_{\bra{0,\geq1}}\in(\Gamma_{\geq1}^0)^{2\times2}\), the symbol
\(-\ii{\bm A}_{\bra{0,\geq1}}\) is again a generic bounded matrix symbol. Renaming
it as
\({\bm A}^{\bC,{\rm d}}_{\bra{0,\geq1}}\), and renaming the smoothing remainder as
\({\bm R}^{\bC,{\rm d}}_{\bra{\varrho,\geq1}}\), gives exactly
\cref{eq:subprincipal_diagonalized_system_3d}.
\end{proofpart}
\end{proof}

\subsubsection{Proof of \cref{prop:positive_order_diagonalization_3d}}
\label{sec:proof_positive_order_diagonalization_3d}

\begin{proofpart}[Construction of the complete diagonalizer]
Let \({\bm B}_0\) and \({\bm B}_{-1}\) be the principal and
subprincipal diagonalizing operators constructed in
\cref{prop:principal_order_diagonalization_3d,prop:subprincipal_order_diagonalization_3d},
respectively, and define
\begin{equation}\label{eq:positive_order_complete_diagonalizer_3d}
    {\bm B}
    \pare{\Atw,\background\middle|U}
    \defeq
    {\bm B}_0
    \pare{\Atw,\background\middle|U}
    {\bm B}_{-1}
    \pare{\Atw,\background\middle|U}.
\end{equation}
Both factors in \cref{eq:positive_order_complete_diagonalizer_3d} are
bounded and invertible on \(H^s\), for every \(s\geq s_0\). Hence
\({\bm B}\) has the same properties. Moreover, by
\cref{eq:principal_diagonalized_unknown_3d,eq:subprincipal_diagonalized_unknown_3d},
\begin{equation}\label{eq:positive_order_unknown_composition_3d}
    W
    =
    {\bm B}_{-1}^{-1}U_0
    =
    {\bm B}_{-1}^{-1}{\bm B}_0^{-1}U
    =
    {\bm B}^{-1}U,
\end{equation}
which proves \cref{eq:positive_order_diagonalized_unknown_3d}.
\end{proofpart}

\begin{proofpart}[Identification of the positive-order symbol]
By \cref{eq:subprincipal_diagonalized_system_3d}, the variable \(W\)
satisfies
\begin{equation}\label{eq:positive_order_assembled_system_3d}
\begin{aligned}
    \partial_t W
    =&
    \OpBW{
        {\bm A}^{\bC,{\rm lin}}_{\rm lo}
        \pare{\Atw,\background\middle|\xi}
    }W
    \\
    +
    &\ii\OpBW{
        {\bm D}^{\bC}_{\frac32}
        \pare{\Atw,\background\middle|U;x,\xi}
        +
        \dot{\bm D}^{\bC}_{1}
        \pare{\Atw,\background\middle|U;x,\xi}
        +
        \dot{\bm D}^{\bC}_{\frac12}
        \pare{\Atw,\background\middle|U;x,\xi}
        +
        {\bm A}^{\bC,{\rm d}}_{\bra{0,\geq1}}
        \pare{\Atw,\background\middle|U;x,\xi}
    }W
    \\
    +&
    {\bm R}^{\bC,{\rm d}}_{\bra{\varrho,\geq1}}
    \pare{\Atw,\background\middle|U}W.
\end{aligned}
\end{equation}
The three positive-order symbols in
\cref{eq:positive_order_assembled_system_3d} are
\begin{align}
    {\bm D}^{\bC}_{\frac32}
    \pare{\Atw,\background\middle|U;x,\xi}
    \defeq{}&
    \pare{1-\chi\pare{\xi}}\,
    \omega_{\frac32}
    \pare{\Atw,\background\middle|\xi}
    \dot\theta_0
    \pare{U;x,\hat\xi}
    {\bm J}_{\rm d},
    \label{eq:positive_order_proof_D32_3d}
    \\
    \dot{\bm D}^{\bC}_{1}
    \pare{\Atw,\background\middle|U;x,\xi}
    \defeq{}&
    \pare{
        \mathpzc V_{\textnormal{can}}
        \pare{
            \Atw,\background
            \middle|
            \eta\pare{U},\tilde\psi\pare{U};
            x
        }
        \cdot\xi
    }
    \Id_{\bC^2},
    \label{eq:positive_order_proof_D1_3d}
    \\
    \dot{\bm D}^{\bC}_{\frac12}
    \pare{\Atw,\background\middle|U;x,\xi}
    \defeq{}&
    -
    \pare{1-\chi\pare{\xi}}\,
    \dot\alpha^{\bC,{\rm pr}}_0
    \pare{\Atw,\background\middle|U;x,\hat\xi}
    \av{\xi}^{\frac12}
    {\bm J}_{\rm d}.
    \label{eq:positive_order_proof_D12_3d}
\end{align}
These formulas are precisely
\cref{eq:positive_order_diagonal_ql_symbol_3d,eq:positive_order_diagonal_transport_symbol_3d,eq:positive_order_diagonal_half_symbol_3d}.
In particular, the symbols are real-valued and diagonal for real \(U\),
and have orders \(3/2\), \(1\), and \(1/2\), respectively. Therefore,
by \cref{eq:positive_order_diagonal_symbol_decomposition_3d},
\begin{equation}\label{eq:positive_order_proof_Dgt0_3d}
    {\bm D}^{\bC}_{>0}
    =
    {\bm D}^{\bC}_{\frac32}
    +
    \dot{\bm D}^{\bC}_{1}
    +
    \dot{\bm D}^{\bC}_{\frac12}
    \in
    \pare{\Gamma_{\geq0}^{\frac32}}^{2\times2}.
\end{equation}
\end{proofpart}

\begin{proofpart}[Bounded symbols and smoothing remainders]
Set
\begin{equation}\label{eq:positive_order_final_remainders_3d}
    {\bm A}_{\bra{0,\geq1}}
    \pare{\Atw,\background\middle|U;x,\xi}
    \defeq
    {\bm A}^{\bC,{\rm d}}_{\bra{0,\geq1}}
    \pare{\Atw,\background\middle|U;x,\xi},
    \qquad
    {\bm R}_{\bra{\varrho,\geq1}}
    \pare{\Atw,\background\middle|U}
    \defeq
    {\bm R}^{\bC,{\rm d}}_{\bra{\varrho,\geq1}}
    \pare{\Atw,\background\middle|U}.
\end{equation}
By \cref{prop:subprincipal_order_diagonalization_3d},
\begin{equation}\label{eq:positive_order_final_remainder_classes_3d}
    {\bm A}_{\bra{0,\geq1}}
    \in
    \pare{\Gamma_{\geq1}^{0}}^{2\times2},
    \qquad
    {\bm R}_{\bra{\varrho,\geq1}}
    \in
    \pare{\cR_{\geq1}^{\varrho}}^{2\times2}.
\end{equation}
The flat multiplier
\({\bm A}^{\bC,{\rm lin}}_{\rm lo}\) remains explicit in
\cref{eq:positive_order_assembled_system_3d} and is not included in
\({\bm A}_{\bra{0,\geq1}}\). Substituting
\cref{eq:positive_order_proof_Dgt0_3d,eq:positive_order_final_remainders_3d}
into \cref{eq:positive_order_assembled_system_3d} gives
\cref{eq:positive_order_diagonalized_system_3d}.
\end{proofpart}
\qed

\subsection{Uniform energy estimates for the diagonalized system}
\label{sec:modified_energy_estimates_3d}

The purpose of this subsection is to prove a uniform a priori estimate for
classical solutions in the standard Sobolev phase space
\(\mathcal X^s\) of \cref{eq:main_phase_space}.  The proof is carried out in
three steps.  First we prove a closed estimate for a modified energy adapted
to the diagonalized variables.  Then we compare this modified energy with the
Sobolev norm, keeping track of the purely low-frequency component of
\(\tilde\psi\).  Finally we combine the two estimates.

\begin{prop}[Uniform inverse-size a priori estimate in the Sobolev scale]
\label{prop:uniform_apriori_sobolev_energy_3d}
Let \(K\Subset\mathcal S_{\rm st}\).  There exist \(s_0,r_0>0\) such that,
for every \(s\geq s_0\), the following holds uniformly for
\(\pare{\Atw,\background}\in K\).  Let
\[
    Y\pare t
    \defeq
    \vect{\eta\pare t}{\tilde\psi\pare t}
\]
be a real classical solution on an interval \(I\ni0\).  Assume that the
reductions of
\cref{prop:real_alinhac_good_unknown_3d,prop:complex_paradifferential_system_3d,prop:positive_order_diagonalization_3d}
are defined on \(I\), and that
\[
    \sup_{t\in I}
    \norm{Y\pare t}_{\mathcal X^{s_0}}
    \leq r_0 .
\]
Then there exist constants
\[
    \varepsilon_{s,K}>0,
    \qquad
    c_{s,K}>0,
    \qquad
    C_{s,K}>0,
\]
such that, if
\[
    \norm{Y\pare0}_{\mathcal X^s}
    \leq
    \varepsilon
    \leq
    \varepsilon_{s,K},
\]
then, for every \(t\in I\) such that
\[
    \av t\leq c_{s,K}\varepsilon^{-1},
\]
one has
\begin{equation}\label{eq:uniform_sobolev_apriori_lifespan_form_3d}
    \norm{Y\pare t}_{\mathcal X^s}
    \leq
    C_{s,K}\pare{\varepsilon+c_{s,K}} .
\end{equation}
Moreover, the Hamiltonian component satisfies the sharper bound
\begin{equation}\label{eq:uniform_hamiltonian_apriori_lifespan_form_3d}
    \norm{Y\pare t}_{\mathcal X^s_{\rm Ham}
    \pare{\Atw,\background}}
    \leq
    C_{s,K}\varepsilon .
\end{equation}
In particular, after choosing \(c_{s,K}\) sufficiently small, the solution
remains in a fixed small Sobolev neighborhood of the origin on a time
interval of size comparable to \(\varepsilon^{-1}\).  All constants are
uniform for \(\pare{\Atw,\background}\in K\).
\end{prop}

\paragraph{Roadmap of the proof.}
The proof of \cref{prop:uniform_apriori_sobolev_energy_3d} is organized in
four steps.
\begin{enumerate}
    \item The first step makes the dependence on the physical parameters quantitative.
For every \(K\Subset\mathcal S_{\rm st}\), compactness and the strict
Atwood--Weber inequality provide a uniform positive distance from the
stability threshold. In
\cref{proofpart:uniform_stable_parameter_control_3d} this is converted into
the uniform ellipticity estimate
\cref{eq:q_uniform_ellipticity_nonhomogeneous_3d} for the flat restoring
symbol. Consequently, the flat diagonalizing multipliers, their inverses,
and the corresponding symbolic estimates are uniform for
\(\pare{\Atw,\background}\in K\). This is the source of the uniformity in
the parameter-dependent constants used throughout the argument.
\item The second step works in the diagonalized complex variable constructed in
\cref{prop:real_alinhac_good_unknown_3d,%
prop:complex_paradifferential_system_3d,%
prop:positive_order_diagonalization_3d}.
At high frequency the positive-order part is diagonal at principal level.
The modified weight is chosen so that the potentially derivative-losing
commutator is cancelled by the Poisson-bracket identity
\cref{eq:modified_energy_poisson_cancellation_3d}. At low frequency the
energy is completed with the Hamiltonian low-frequency symmetrizer of
\cref{lem:flat_low_frequency_hamiltonian_symmetrizer_3d}. The crucial point
is the exact flat cancellation
\cref{eq:full_flat_linear_symmetrizer_identity_3d}: the isolated
low-frequency defect is not skew-adjoint, but the complete flat linear
symbol is. These two cancellations give the closed modified energy estimate
of \cref{prop:modified_energy_estimate_diagonalized_3d} and the inverse-size
modified energy bound of
\cref{cor:uniform_apriori_diagonalized_energy_3d}.
\item The third step compares the modified energy with the physical Sobolev
variables. This is done in
\cref{sec:modified_sobolev_energy_equivalence_3d}. The modified energy is
equivalent to the Hamiltonian Sobolev norm
\(\mathcal X^s_{\rm Ham}\pare{\Atw,\background}\), which controls all
high-frequency components and, at low frequency, the Hamiltonian variables
\[
    \eta,
    \qquad
    \av{D}^{\frac12}\tilde\psi .
\]
The only missing component of the full inhomogeneous Sobolev norm is the
undifferentiated low-frequency part of \(\tilde\psi\), isolated in
\cref{eq:low_frequency_potential_control_norm_3d}.
\item The fourth step propagates this missing component by the rough estimate of
\cref{lem:rough_low_frequency_potential_propagation_3d}. Since the
Hamiltonian energy remains \(O\pare{\varepsilon}\), this component grows at
most like \(O\pare{\varepsilon\av t}\). Therefore, on the time interval
\(\av t\leq c_{s,K}\varepsilon^{-1}\), it is bounded by
\(O\pare{\varepsilon+c_{s,K}}\). Choosing \(c_{s,K}\) sufficiently small
keeps the full Sobolev norm inside the perturbative regime, while the
Hamiltonian component remains \(O\pare{\varepsilon}\). This proves the
inverse-size a priori estimate stated in
\cref{prop:uniform_apriori_sobolev_energy_3d}.

\end{enumerate}
\subsubsection{Uniform modified energy estimates}
\label{sec:uniform_modified_energy_estimates_3d}

\begin{proofpart}[Uniform control on compact stable parameter sets]
\label{proofpart:uniform_stable_parameter_control_3d}
Fix \(K\Subset\mathcal S_{\rm st}\).  Since
\(\mathsf W_{\Atw,\background}\) is continuous in the physical parameters
and \(\mathsf W_{\Atw,\background}<1\) on \(\mathcal S_{\rm st}\), compactness
gives
\begin{equation}\label{eq:uniform_atwood_weber_gap_3d}
\delta_K
\defeq
1-
\sup_{\pare{\Atw,\background}\in K}
\mathsf W_{\Atw,\background}
>0.
\end{equation}
Moreover,
$$
    \inf_{\pare{\Atw,\background}\in K}
    \mathfrak g\Atw
    >0.
$$
Let \(\xi\neq0\).  From
\cref{eq:def_linearized_q_background,eq:strict_atwood_weber_number_intro},
\begin{equation}\label{eq:uniform_q_shear_control_3d}
-\mathfrak q
\pare{\Atw,\background\middle|\xi}
\geq
\mathfrak g\Atw
+
\frac{\sigma}{\aangle{\rho^\pm}}\av{\xi}^2
\
-
2\mathsf W_{\Atw,\background}
\sqrt{
\frac{\sigma\mathfrak g\Atw}
{\aangle{\rho^\pm}}
}
\av{\xi}.
\end{equation}
Using
$$
    2
    \sqrt{
        \frac{\sigma\mathfrak g\Atw}
        {\aangle{\rho^\pm}}
    }
    \av{\xi}
    \leq
    \mathfrak g\Atw
    +
    \frac{\sigma}{\aangle{\rho^\pm}}\av{\xi}^2
$$
and \cref{eq:uniform_atwood_weber_gap_3d}, we obtain
$$
    -\mathfrak q
    \pare{\Atw,\background\middle|\xi}
    \geq
    \delta_K
    \pare{
        \mathfrak g\Atw
        +
        \frac{\sigma}{\aangle{\rho^\pm}}\av{\xi}^2
    }.
$$
Since the coefficients in the last expression vary continuously on \(K\),
there exist \(0<c_K\leq C_K\) such that, uniformly for
\(\pare{\Atw,\background}\in K\),
\begin{equation}\label{eq:q_uniform_ellipticity_nonhomogeneous_3d}
c_K\angles{\xi}^2
\leq
-\mathfrak q
\pare{\Atw,\background\middle|\xi}
\leq
C_K\angles{\xi}^2,
\qquad
\xi\in\bR^2.
\end{equation}
By
\cref{eq:complex_frequency_square_identity,eq:q_uniform_ellipticity_nonhomogeneous_3d},
there exist \(0<c_K'\leq C_K'\) such that, uniformly for
\(\pare{\Atw,\background}\in K\),
\[
    c_K'
    \av{\xi}^{\frac12}\angles{\xi}
    \leq
    \omega_{\frac32}
    \pare{\Atw,\background\middle|\xi}
    \leq
    C_K'
    \av{\xi}^{\frac12}\angles{\xi},
    \qquad
    \xi\neq0.
\]
Hence, by \cref{eq:complex_multiplier_p},
\begin{equation}\label{eq:uniform_flat_complex_multiplier_bounds_3d}
    \mathfrak m_{-\frac14}
    \pare{\Atw,\background\middle|\xi}
    \sim_K
    \av{\xi}^{\frac14}\angles{\xi}^{-\frac12},
    \qquad
    \mathfrak m_{-\frac14}^{-1}
    \pare{\Atw,\background\middle|\xi}
    \sim_K
    \av{\xi}^{-\frac14}\angles{\xi}^{\frac12}.
\end{equation}
In particular, on \(\supp\pare{1-\chi}\), where \(\av{\xi}\) is bounded
away from zero,
\[
    \mathfrak m_{-\frac14}
    \sim_K
    \angles{\xi}^{-\frac14},
    \qquad
    \mathfrak m_{-\frac14}^{-1}
    \sim_K
    \angles{\xi}^{\frac14}.
\]
The corresponding symbol seminorms are uniform on \(K\): indeed,
\(\mathfrak q\) is given explicitly by
\cref{eq:def_linearized_q_background}, while
\(\omega_{\frac32}\) and \(\mathfrak m_{-\frac14}\) are obtained from it
through the smooth operations in
\cref{eq:complex_frequency_square_identity,eq:complex_multiplier_p};
on \(\supp\pare{1-\chi}\) their denominators are uniformly bounded away
from zero by the estimates above.

It remains only to check the low-frequency regularization. By
\cref{eq:complex_nonhomogeneous_multiplier_p_3d},
\(\mathfrak m^{\rm nh}_{-\frac14}=1\) on the core where \(\chi=1\), while
it coincides with \(\mathfrak m_{-\frac14}\) where \(\chi=0\).
On the fixed transition annulus \(0<\chi<1\),
\cref{eq:uniform_flat_complex_multiplier_bounds_3d} and compactness of \(K\)
give uniform upper and positive lower bounds for
\(\mathfrak m_{-\frac14}\). Therefore
\(\mathfrak m^{\rm nh}_{-\frac14}\) and its inverse have uniform symbol
seminorms for \(\pare{\Atw,\background}\in K\). By
\cref{eq:complex_matrix_M,eq:complex_matrix_M_inverse}, the same is true
for the flat complex change of variables and its inverse.

Thus the flat multipliers which enter the complex reduction are uniformly
controlled on \(K\). Together with the convention that the symbol and
smoothing classes have uniform seminorm bounds when the external parameters
range in a fixed compact set, this allows all constants arising in the
subsequent diagonalization, symmetrization, and energy estimates to be
chosen uniformly for \(\pare{\Atw,\background}\in K\).
\end{proofpart}

Recall the cutoff function \(\chi\) fixed in \cref{def-chi}.  We set
\begin{equation}\label{eq:modified_energy_cutoffs_3d}
    \chi_{\rm lo}\pare{\xi}
    \defeq
    \chi\pare{\frac{\xi}{2}},
    \qquad
    \chi_{\rm hi}\pare{\xi}
    \defeq
    1-\chi_{\rm lo}\pare{\xi}.
\end{equation}
In particular, recalling \cref{eq:linear_low_frequency_defect_complex_3d,eq:linear_high_frequency_exact_diagonalization_3d},
\begin{equation}\label{eq:modified_energy_cutoff_separation_3d}
    \chi_{\rm lo}=1
    \quad\textnormal{on }\supp\pare{\chi},
    \qquad
    \chi_{\rm hi}\,
    {\bm A}^{\bC,{\rm lin}}_{\rm lo}
    =
    0.
\end{equation}
Recalling \cref{eq:positive_order_diagonal_symbol_decomposition_3d}, we can write
\begin{equation}\label{eq:diagonal_entries_positive_symbol_3d}
    {\bm D}^{\bC}_{>0}
    \pare{\Atw,\background\middle|U;x,\xi}
    =
    \operatorname{diag}
    \pare{
        d_1\pare{\Atw,\background\middle|U;x,\xi},
        d_2\pare{\Atw,\background\middle|U;x,\xi}
    },
\end{equation}
where
\begin{align}
    d_1
    &\defeq
    \mathpzc V_{\textnormal{can}}
    \pare{
        \Atw,\background
        \middle|
        \eta\pare{U},\tilde\psi\pare{U};
        x
    }
    \cdot\xi
    -
    \pare{1-\chi\pare{\xi}}\,
    \omega_{\frac32}
    \pare{\Atw,\background\middle|\xi}
    \dot\theta_0
    \pare{\Atw,\background\middle|U;x,\hat\xi}
    \nonumber
    \\
    &\qquad+
    \pare{1-\chi\pare{\xi}}\,
    \dot\alpha^{\bC,{\rm pr}}_0
    \pare{\Atw,\background\middle|U;x,\hat\xi}
    \av{\xi}^{\frac12},
    \label{eq:first_diagonal_entry_positive_symbol_3d}
    \\
    d_2
    &\defeq
    \mathpzc V_{\textnormal{can}}
    \pare{
        \Atw,\background
        \middle|
        \eta\pare{U},\tilde\psi\pare{U};
        x
    }
    \cdot\xi
    +
    \pare{1-\chi\pare{\xi}}\,
    \omega_{\frac32}
    \pare{\Atw,\background\middle|\xi}
    \dot\theta_0
    \pare{\Atw,\background\middle|U;x,\hat\xi}
    \nonumber
    \\
    &\qquad-
    \pare{1-\chi\pare{\xi}}\,
    \dot\alpha^{\bC,{\rm pr}}_0
    \pare{\Atw,\background\middle|U;x,\hat\xi}
    \av{\xi}^{\frac12}.
    \label{eq:second_diagonal_entry_positive_symbol_3d}
\end{align}
Here \(\hat\xi\defeq\xi/\av{\xi}\), and every symbol depending on
\(\hat\xi\) is multiplied by \(1-\chi\pare{\xi}\). We split
\begin{equation}\label{eq:diagonalized_positive_symbol_linear_nonlinear_split_3d}
    {\bm D}^{\bC}_{>0}
    =
    {\bm D}^{\bC}_{>0,{\rm lin}}
    +
    {\bm D}^{\bC}_{>0,{\rm nl}},
\end{equation}
where
\begin{equation}\label{eq:flat_diagonal_symbol_3d}
    {\bm D}^{\bC}_{>0,{\rm lin}}
    \pare{\Atw,\background\middle|\xi}
    \defeq
    {\bm D}^{\bC}_{>0}
    \pare{\Atw,\background\middle|0;x,\xi},
\end{equation}
and
\begin{equation}\label{eq:nonlinear_diagonal_symbol_3d}
    {\bm D}^{\bC}_{>0,{\rm nl}}
    \pare{\Atw,\background\middle|U;x,\xi}
    \defeq
    {\bm D}^{\bC}_{>0}
    \pare{\Atw,\background\middle|U;x,\xi}
    -
    {\bm D}^{\bC}_{>0,{\rm lin}}
    \pare{\Atw,\background\middle|\xi}.
\end{equation}
Thus the diagonalized equation for \(W\) can be written as
\begin{equation}\label{eq:diagonalized_system_linear_nonlinear_split_3d}
\begin{aligned}
    \partial_tW
    &=
    \OpBW{
        {\bm A}^{\bC,{\rm lin}}_{\rm lo}
        \pare{\Atw,\background\middle|\xi}
    +\ii 
        {\bm D}^{\bC}_{>0,{\rm lin}}
        \pare{\Atw,\background\middle|\xi}
    }W
    \\
    &\quad+
    \ii
    \OpBW{
        {\bm D}^{\bC}_{>0,{\rm nl}}
        \pare{\Atw,\background\middle|U;x,\xi}
        +
        {\bm A}_{\bra{0,\geq1}}
        \pare{\Atw,\background\middle|U;x,\xi}
    }W
    +
    {\bm R}_{\bra{\varrho,\geq1}}
    \pare{\Atw,\background\middle|U}W.
\end{aligned}
\end{equation}
We next introduce the flat low-frequency symmetrizer.  Let
\begin{equation}\label{eq:real_low_frequency_symmetrizer_3d}
    {\bm S}^{\bR}_{\rm lo}
    \pare{\Atw,\background\middle|\xi}
    \defeq
    \begin{bmatrix}
        -\mathfrak q
        \pare{\Atw,\background\middle|\xi}
        &
        0
        \\
        0
        &
        \dfrac{\av{\xi}}{2}
    \end{bmatrix},
\end{equation}
and define
\begin{equation}\label{eq:complex_low_frequency_symmetrizer_3d}
    {\bm S}^{\bC}_{\rm lo}
    \pare{\Atw,\background\middle|\xi}
    \defeq
    \pare{
        {\bm M}
        \pare{\Atw,\background\middle|\xi}^{-1}
    }^*
    {\bm S}^{\bR}_{\rm lo}
    \pare{\Atw,\background\middle|\xi}
    {\bm M}
    \pare{\Atw,\background\middle|\xi}^{-1}.
\end{equation}
For \(s\geq0\), let
\begin{equation}\label{eq:modified_energy_weight_symbols_3d}
    \lambda_s\pare{z}
    \defeq
    \pare{1+z^2}^{\frac{s}{3}},
    \qquad
    {\bm \Lambda}_{s,{\rm hi}}
    \pare{\Atw,\background\middle|U;x,\xi}
    \defeq
    \chi_{\rm hi}\pare{\xi}
    \operatorname{diag}
    \pare{
        \lambda_s\pare{d_1},
        \lambda_s\pare{d_2}
    }.
\end{equation}
Since \(d_j\) is of order \(3/2\), the weight
\(\lambda_s\pare{d_j}\) is of order \(s\) on the high-frequency region.
For the unknowns \(U\) and \(W\), define
\begin{equation}\label{eq:modified_energy_unknowns_3d}
    W_{s,{\rm hi}}
    \defeq
    \OpBW{
        {\bm \Lambda}_{s,{\rm hi}}
        \pare{\Atw,\background\middle|U}
    }W,
    \qquad
    U_{s,{\rm lo}}
    \defeq
    \chi_{\rm lo}\pare{D}
    \angles{D}^{s-\frac34}U.
\end{equation}
The modified energy is
\begin{equation}\label{eq:modified_energy_definition_3d}
    E_s\pare{t}
    \defeq
    \frac12
    \norm{W_{s,{\rm hi}}\pare{t}}_{L^2}^2
    +
    \frac12
    \psc{
        {\bm S}^{\bC}_{\rm lo}
        \pare{\Atw,\background\middle|D}
        U_{s,{\rm lo}}\pare{t}
    }{
        U_{s,{\rm lo}}\pare{t}
    }.
\end{equation}
Our goal is to prove the following modified energy estimate of the diagonalized system.
\begin{prop}[Modified energy estimate for the diagonalized system]
\label{prop:modified_energy_estimate_diagonalized_3d}
Let the assumptions of
\cref{prop:positive_order_diagonalization_3d} hold.  There exist
\(s_0,r_0>0\) such that, for every \(s\geq s_0\), every classical real
solution satisfying
\[
    \sup_{t\in I}E_{s_0}^{\frac12}\pare{t}
    \leq r_0
\]
obeys
\begin{equation}\label{eq:modified_energy_closed_estimate_3d}
    \av{
        \frac{\dd}{\dd t}E_s\pare{t}
    }
    \leq
    C_s
    E_{s_0}^{\frac12}\pare{t}
    E_s\pare{t}.
\end{equation}
The constant \(C_s\) is uniform when
\(\pare{\Atw,\background}\) ranges in a compact subset of
\(\mathcal S_{\rm st}\).
\end{prop}
To prove \cref{prop:modified_energy_estimate_diagonalized_3d}, we start by proving the following two lemmas about flat low-frequency Hamiltonian symmetrizer and high-low frequency modified energy bounds.
\begin{lemma}[Flat low-frequency Hamiltonian symmetrizer]
\label{lem:flat_low_frequency_hamiltonian_symmetrizer_3d}
Assume \cref{hyp:strict_atwood_weber_stability}.  Set
\[
    {\bm Q}_0
    \pare{\Atw,\background\middle|\xi}
    \defeq
    \begin{bmatrix}
        0&\dfrac{\av{\xi}}2\\
        \mathfrak q\pare{\Atw,\background\middle|\xi}&0
    \end{bmatrix},
    \qquad
    \nu_1
    \pare{\Atw,\background\middle|\xi}
    \defeq
    \frac{\Atw}{4}
    \pare{\background^\perp\cdot\xi}.
\]
Then the following identities hold pointwise in \(\xi\).

\begin{enumerate}
    \item The complex low-frequency symmetrizer has the explicit form
    \begin{equation}\label{eq:complex_low_frequency_symmetrizer_explicit_3d}
    \begin{aligned}
        {\bm S}^{\bC}_{\rm lo}
        &=
        \frac12
        \begin{bmatrix}
            -\mathfrak q
            \pare{\mathfrak m^{\rm nh}_{-\frac14}}^2
            +
            \dfrac{\av{\xi}}2
            \pare{\mathfrak m^{\rm nh}_{-\frac14}}^{-2}
            &
            -\mathfrak q
            \pare{\mathfrak m^{\rm nh}_{-\frac14}}^2
            -
            \dfrac{\av{\xi}}2
            \pare{\mathfrak m^{\rm nh}_{-\frac14}}^{-2}
            \\
            -\mathfrak q
            \pare{\mathfrak m^{\rm nh}_{-\frac14}}^2
            -
            \dfrac{\av{\xi}}2
            \pare{\mathfrak m^{\rm nh}_{-\frac14}}^{-2}
            &
            -\mathfrak q
            \pare{\mathfrak m^{\rm nh}_{-\frac14}}^2
            +
            \dfrac{\av{\xi}}2
            \pare{\mathfrak m^{\rm nh}_{-\frac14}}^{-2}
        \end{bmatrix}
        \pare{\Atw,\background\middle|\xi}.
    \end{aligned}
    \end{equation}
    In particular \({\bm S}^{\bC}_{\rm lo}\) is Hermitian and
    nonnegative.  It is positive on the Hamiltonian low-frequency variables
    \[
        \eta,
        \qquad
        \av{D}^{\frac12}\tilde\psi,
    \]
    and it is degenerate only on the purely undifferentiated zero-frequency
    component of \(\tilde\psi\).

    \item The real Hamiltonian symmetrizer cancels the flat off-diagonal
    real symbol:
    \begin{equation}\label{eq:real_flat_hamiltonian_symmetrizer_identity_3d}
        {\bm S}^{\bR}_{\rm lo}{\bm Q}_0
        +
        {\bm Q}_0^*
        {\bm S}^{\bR}_{\rm lo}
        =
        0 .
    \end{equation}
    Equivalently, in the complex coordinate,
    \begin{equation}\label{eq:complex_flat_off_diagonal_symmetrizer_identity_3d}
        {\bm S}^{\bC}_{\rm lo}
        \pare{
            {\bm M}{\bm Q}_0{\bm M}^{-1}
        }
        +
        \pare{
            {\bm M}{\bm Q}_0{\bm M}^{-1}
        }^*
        {\bm S}^{\bC}_{\rm lo}
        =
        0 .
    \end{equation}

    \item The complete flat complex linear symbol satisfies
    \begin{equation}\label{eq:full_flat_linear_symmetrizer_identity_3d}
    \begin{aligned}
        &{\bm S}^{\bC}_{\rm lo}
        \pare{
            {\bm A}^{\bC,{\rm lin}}_{\rm lo}
            +
            \ii{\bm D}^{\bC}_{>0,{\rm lin}}
        }
       +
        \pare{
            {\bm A}^{\bC,{\rm lin}}_{\rm lo}
            +
            \ii{\bm D}^{\bC}_{>0,{\rm lin}}
        }^*
        {\bm S}^{\bC}_{\rm lo}
        =
        0 .
    \end{aligned}
    \end{equation}

    \item The low-frequency block has the diagonalization defect
    \begin{equation}\label{eq:low_frequency_linear_symmetrizer_defect_3d}
    \begin{aligned}
        &{\bm S}^{\bC}_{\rm lo}
        {\bm A}^{\bC,{\rm lin}}_{\rm lo}
        +
        \pare{
            {\bm A}^{\bC,{\rm lin}}_{\rm lo}
        }^*
        {\bm S}^{\bC}_{\rm lo}
       =
        -\ii
        \pare{1-\chi\pare{\xi}}\,
        \omega_{\frac32}
        \pare{\Atw,\background\middle|\xi}
        \pare{
            {\bm S}^{\bC}_{\rm lo}{\bm J}_{\rm d}
            -
            {\bm J}_{\rm d}{\bm S}^{\bC}_{\rm lo}
        } .
    \end{aligned}
    \end{equation}
\end{enumerate}
\end{lemma}

\begin{proof}
The explicit formula
\cref{eq:complex_low_frequency_symmetrizer_explicit_3d} follows from
\cref{eq:complex_low_frequency_symmetrizer_3d}.  Indeed,
\[
    {\bm M}^{-1}
    =
    \frac1{\sqrt2}
    \begin{bmatrix}
        \mathfrak m^{\rm nh}_{-\frac14}
        &
        \mathfrak m^{\rm nh}_{-\frac14}
        \\
        -\ii
        \pare{\mathfrak m^{\rm nh}_{-\frac14}}^{-1}
        &
        \ii
        \pare{\mathfrak m^{\rm nh}_{-\frac14}}^{-1}
    \end{bmatrix},
    \qquad
    {\bm S}^{\bR}_{\rm lo}
    =
    \begin{bmatrix}
        -\mathfrak q&0\\
        0&\dfrac{\av{\xi}}2
    \end{bmatrix},
\]
and
\[
    {\bm S}^{\bC}_{\rm lo}
    =
    \pare{{\bm M}^{-1}}^*
    {\bm S}^{\bR}_{\rm lo}
    {\bm M}^{-1}.
\]
A direct multiplication gives
\cref{eq:complex_low_frequency_symmetrizer_explicit_3d}.  Since
\cref{hyp:strict_atwood_weber_stability} gives
$-\mathfrak q
    \pare{\Atw,\background\middle|\xi}
    >0,$
the real quadratic form associated with \({\bm S}^{\bR}_{\rm lo}\) is
\[
    -\mathfrak q
    \av{\eta}^2
    +
    \frac{\av{\xi}}2
    \av{\tilde\psi}^2.
\]
It is nonnegative and degenerates only at \(\xi=0\) on the undifferentiated
\(\tilde\psi\)-component.  Since \({\bm S}^{\bC}_{\rm lo}\) is the pullback
of \({\bm S}^{\bR}_{\rm lo}\) by the invertible matrix \({\bm M}\), the same
positivity and degeneracy properties hold in the complex coordinate.

We now prove the symmetrizer identities.  Direct multiplication gives
\[
    {\bm S}^{\bR}_{\rm lo}{\bm Q}_0
    =
    \begin{bmatrix}
        0&
        -\dfrac{\av{\xi}}2\mathfrak q
        \\
        \dfrac{\av{\xi}}2\mathfrak q&
        0
    \end{bmatrix},
    \qquad
    {\bm Q}_0^*
    {\bm S}^{\bR}_{\rm lo}
    =
    \begin{bmatrix}
        0&
        \dfrac{\av{\xi}}2\mathfrak q
        \\
        -\dfrac{\av{\xi}}2\mathfrak q&
        0
    \end{bmatrix}.
\]
This proves
\cref{eq:real_flat_hamiltonian_symmetrizer_identity_3d}.  Conjugating this
identity by \({\bm M}^{-1}\), and using the definition of
\({\bm S}^{\bC}_{\rm lo}\), gives
\cref{eq:complex_flat_off_diagonal_symmetrizer_identity_3d}.

At the flat state the nonlinear principal corrections vanish:
\[
    \dot\theta_0=1,
    \qquad
    \mathpzc V_{\textnormal{can}}
    \pare{\Atw,\background\middle|0,0}
    =
    \frac{\Atw}{4}\background^\perp,
    \qquad
    \dot\alpha^{\bC,{\rm pr}}_0=0.
\]
Hence the flat diagonal positive-order part is
\[
    {\bm D}^{\bC}_{>0,{\rm lin}}
    =
    \nu_1\Id_{\bC^2}
    +
    \pare{1-\chi}
    \omega_{\frac32}{\bm J}_{\rm d}.
\]
On the other hand, by
\cref{eq:linear_low_frequency_defect_complex_3d},
\[
    {\bm A}^{\bC,{\rm lin}}_{\rm lo}
    =
    {\bm M}{\bm Q}_0{\bm M}^{-1}
    -
    \ii
    \pare{1-\chi}
    \omega_{\frac32}{\bm J}_{\rm d}.
\]
Therefore the two transition terms cancel inside
\[
    {\bm A}^{\bC,{\rm lin}}_{\rm lo}
    +
    \ii{\bm D}^{\bC}_{>0,{\rm lin}}.
\]
Consequently, the left-hand side of
\cref{eq:full_flat_linear_symmetrizer_identity_3d} reduces to
\[
\begin{aligned}
    &{\bm S}^{\bC}_{\rm lo}
    \pare{
        {\bm M}{\bm Q}_0{\bm M}^{-1}
    }
    +
    \pare{
        {\bm M}{\bm Q}_0{\bm M}^{-1}
    }^*
    {\bm S}^{\bC}_{\rm lo}
    \\
    &\qquad
    +
    {\bm S}^{\bC}_{\rm lo}
    \ii\nu_1\Id_{\bC^2}
    +
    \pare{\ii\nu_1\Id_{\bC^2}}^*
    {\bm S}^{\bC}_{\rm lo}.
\end{aligned}
\]
The first line vanishes by
\cref{eq:complex_flat_off_diagonal_symmetrizer_identity_3d}.  The second line
also vanishes, since \(\nu_1\) is real and the term
\(\ii\nu_1\Id_{\bC^2}\) is scalar and skew-adjoint:
\[
    {\bm S}^{\bC}_{\rm lo}
    \ii\nu_1\Id_{\bC^2}
    +
    \pare{\ii\nu_1\Id_{\bC^2}}^*
    {\bm S}^{\bC}_{\rm lo}
    =
    \ii\nu_1{\bm S}^{\bC}_{\rm lo}
    -
    \ii\nu_1{\bm S}^{\bC}_{\rm lo}
    =
    0.
\]
This proves
\cref{eq:full_flat_linear_symmetrizer_identity_3d}.

It remains to record what happens if the diagonal positive-order flat part is
removed.  Using again
\[
    {\bm A}^{\bC,{\rm lin}}_{\rm lo}
    =
    {\bm M}{\bm Q}_0{\bm M}^{-1}
    -
    \ii
    \pare{1-\chi}
    \omega_{\frac32}{\bm J}_{\rm d},
\]
and cancelling the off-diagonal term by
\cref{eq:complex_flat_off_diagonal_symmetrizer_identity_3d}, one obtains
\[
\begin{aligned}
    &{\bm S}^{\bC}_{\rm lo}
    {\bm A}^{\bC,{\rm lin}}_{\rm lo}
    +
    \pare{
        {\bm A}^{\bC,{\rm lin}}_{\rm lo}
    }^*
    {\bm S}^{\bC}_{\rm lo}=
    -\ii
    \pare{1-\chi}
    \omega_{\frac32}
    {\bm S}^{\bC}_{\rm lo}{\bm J}_{\rm d}
    +
    \ii
    \pare{1-\chi}
    \omega_{\frac32}
    {\bm J}_{\rm d}{\bm S}^{\bC}_{\rm lo}.
\end{aligned}
\]
This is exactly
\cref{eq:low_frequency_linear_symmetrizer_defect_3d}.
\end{proof}

\begin{lemma}[High- and low-frequency modified energy bounds]
\label{lem:adapted_diagonal_weights_3d}
Under the assumptions of
\cref{prop:modified_energy_estimate_diagonalized_3d}, one has
\begin{align}
    \av{
        \frac12
        \frac{\dd}{\dd t}
        \norm{W_{s,{\rm hi}}}_{L^2}^2
    }
    &\leq
    C_sE_{s_0}^{\frac12}E_s,
    \label{eq:high_frequency_energy_bound_3d}
    \\
    \av{
        \frac12
        \frac{\dd}{\dd t}
        \psc{
            {\bm S}^{\bC}_{\rm lo}
            U_{s,{\rm lo}}
        }{
            U_{s,{\rm lo}}
        }
    }
    &\leq
    C_sE_{s_0}^{\frac12}E_s.
    \label{eq:low_frequency_energy_bound_3d}
\end{align}
\end{lemma}

\begin{proof}
We omit the parameters
\(\Atw,\background\) and the time variable whenever no confusion is possible.

\begin{proofpart}[High-frequency contribution]
Differentiating the first term in
\cref{eq:modified_energy_definition_3d} and using
\cref{eq:positive_order_diagonalized_system_3d}, we obtain
\begin{align}
    \frac12
    \frac{\dd}{\dd t}
    \norm{W_{s,{\rm hi}}}_{L^2}^2
    &={}
    \Re\,
    \psc{
        \OpBW{
            \partial_t{\bm \Lambda}_{s,{\rm hi}}
        }W
    }{
        W_{s,{\rm hi}}
    }
    \nonumber
    \\
    &+
    \Re\,
    \psc{
        \OpBW{{\bm \Lambda}_{s,{\rm hi}}}
        \OpBW{
            {\bm A}^{\bC,{\rm lin}}_{\rm lo}
        }W
    }{
        W_{s,{\rm hi}}
    }
    \nonumber
    \\
    &+
    \Re\,
    \psc{
        \ii
        \OpBW{{\bm \Lambda}_{s,{\rm hi}}}
        \OpBW{{\bm D}^{\bC}_{>0}}W
    }{
        W_{s,{\rm hi}}
    }
    \nonumber
    \\
    &+
    \Re\,
    \psc{
        \ii
        \OpBW{{\bm \Lambda}_{s,{\rm hi}}}
        \OpBW{{\bm A}_{\bra{0,\geq1}}}W
    }{
        W_{s,{\rm hi}}
    }
    \nonumber
    \\
    &+
    \Re\,
    \psc{
        \OpBW{{\bm \Lambda}_{s,{\rm hi}}}
        {\bm R}_{\bra{\varrho,\geq1}}W
    }{
        W_{s,{\rm hi}}
    }.
    \label{eq:high_frequency_energy_decomposition_3d}
\end{align}

At \(U=0\), the weight
\({\bm \Lambda}_{s,{\rm hi}}\pare0\) is a Fourier multiplier supported away
from \(\supp {\bm A}^{\bC,{\rm lin}}_{\rm lo}\).  Hence its composition with
the flat low-frequency multiplier vanishes.  The remaining part of
\({\bm \Lambda}_{s,{\rm hi}}\) vanishes at least linearly in \(U\).  The
spectral localization of the Bony--Weyl quantization and
\cref{eq:modified_energy_cutoff_separation_3d} therefore give
\begin{equation}\label{eq:high_low_linear_interaction_bound_3d}
    \norm{
        \OpBW{{\bm \Lambda}_{s,{\rm hi}}}
        \OpBW{{\bm A}^{\bC,{\rm lin}}_{\rm lo}}W
    }_{L^2}
    \leq
    C_sE_{s_0}^{\frac12}E_s^{\frac12}.
\end{equation}

We now consider the positive-order diagonal contribution.  Since
\({\bm D}^{\bC}_{>0}\) is real and diagonal,
\[
    \Re\,
    \psc{
        \ii\OpBW{{\bm D}^{\bC}_{>0}}
        W_{s,{\rm hi}}
    }{
        W_{s,{\rm hi}}
    }
    =
    0.
\]
Thus only the commutator remains.  For \(j=1,2\), its potentially
order-losing symbol is
\begin{equation}\label{eq:modified_energy_principal_commutator_symbol_3d}
    \frac{1}{\ii}
    \pbra{
        \chi_{\rm hi}\pare{\xi}
        \lambda_s\pare{d_j}
    }{
        d_j
    }.
\end{equation}
Where \(\chi_{\rm hi}\) is constant, one has the crucial scalar cancellation
\begin{equation}\label{eq:modified_energy_poisson_cancellation_3d}
    \pbra{
        \lambda_s\pare{d_j}
    }{
        d_j
    }
    =
    \lambda_s'\pare{d_j}
    \pbra{d_j}{d_j}
    =
    0.
\end{equation}
When a \(\xi\)-derivative hits \(\chi_{\rm hi}\), the resulting symbol is
supported in a fixed annulus and contains \(\nabla_xd_j\).  It therefore
vanishes at \(U=0\).  After the cancellation
\cref{eq:modified_energy_poisson_cancellation_3d}, the first remaining
high-frequency Weyl term has order
\[
    s+\frac32-3
    =
    s-\frac32.
\]
Consequently,
\begin{equation}\label{eq:adapted_weight_commutator_bound_3d}
    \norm{
        \comm{
            \OpBW{{\bm \Lambda}_{s,{\rm hi}}}
        }{
            \OpBW{{\bm D}^{\bC}_{>0}}
        }W
    }_{L^2}
    \leq
    C_sE_{s_0}^{\frac12}E_s^{\frac12}.
\end{equation}

Similarly,
\(\partial_t{\bm \Lambda}_{s,{\rm hi}}\pare0=0\), and differentiating the
nonlinear coefficients produces one lower norm of \(\partial_tU\).  Using
the equation, the tame estimates, and
\cref{eq:small_factor_control_by_low_energy_3d}, we get
\begin{equation}\label{eq:time_derivative_symmetrizer_bound_3d}
    \norm{
        \OpBW{
            \partial_t{\bm \Lambda}_{s,{\rm hi}}
        }W
    }_{L^2}
    \leq
    C_sE_{s_0}^{\frac12}E_s^{\frac12}.
\end{equation}
Finally, the bounded symbol and the smoothing remainder vanish at \(U=0\).
Their tame estimates give the same right-hand side.  Substitution in
\cref{eq:high_frequency_energy_decomposition_3d} proves
\cref{eq:high_frequency_energy_bound_3d}.
\end{proofpart}

\begin{proofpart}[Low-frequency contribution]
We differentiate the low-frequency part of
\cref{eq:modified_energy_definition_3d} using the complex equation
\cref{eq:complex_paradifferential_system_3d}.  The flat contribution is
\[
    {\bm A}^{\bC,{\rm lin}}_{\rm lo}
    +
    \ii{\bm D}^{\bC}_{>0,{\rm lin}}.
\]
Since all low-frequency weights are constant-coefficient Fourier multipliers,
\cref{eq:full_flat_linear_symmetrizer_identity_3d} gives
\begin{equation}\label{eq:low_frequency_flat_linear_cancellation_3d}
    \Re\,
    \psc{
        {\bm S}^{\bC}_{\rm lo}
        \pare{
            {\bm A}^{\bC,{\rm lin}}_{\rm lo}
            +
            \ii{\bm D}^{\bC}_{>0,{\rm lin}}
        }
        U_{s,{\rm lo}}
    }{
        U_{s,{\rm lo}}
    }
    =
    0.
\end{equation}
This cancellation acts on the whole multiplier
\({\bm A}^{\bC,{\rm lin}}_{\rm lo}\), including its transition-frequency
part.

It remains to estimate the nonlinear terms.  Conjugating the low-frequency
quadratic form by \({\bm M}^{-1}\), these estimates reduce to the real
Alinhac system \cref{eq:real_system_after_alinhac_good_unknown_3d}.  The
off-diagonal quasilinear block has the same Hamiltonian order structure as
its flat part, while the lower-order block satisfies
\[
    {\bm A}^{\bR,{\rm Al}}_{\rm lo}
    \in
    \begin{bmatrix}
        \Gamma_{\geq1}^0&\Gamma_{\geq1}^{-1}\vspace{0.1cm}
        \\
        \Gamma_{\geq1}^0&\Gamma_{\geq1}^0
    \end{bmatrix}
\]
by \cref{eq:real_alinhac_lower_order_structure_3d}.  Together with the
spectral localization of the Bony--Weyl quantization, this structure is
bounded on the low-frequency Hamiltonian energy space
\[
    L^2
    \times
    H^{\frac12}.
\]
The nonlinear coefficients vanish at the flat state.  Since the CDE
coefficients depend on \(\eta\) and on differentiated components of
\(\tilde\psi\), their lower seminorms are controlled by the Hamiltonian
smallness
\[
    \norm{Y}_{\mathcal X^{s_0}_{\rm Ham}
    \pare{\Atw,\background}}
    \lesssim
    E_{s_0}^{\frac12}.
\]
The undifferentiated low-frequency component
\(\norm{Y}_{\mathcal Z^{s_0}_{\rm lo}}\) is therefore not needed in this
closed symmetrized estimate.  The smoothing remainder satisfies the same
bound after fixing \(\varrho\) sufficiently large.  Hence
\[
    \av{
        \Re\,
        \psc{
            {\bm S}^{\bC}_{\rm lo}
            \chi_{\rm lo}\pare{D}
            \angles{D}^{s-\frac34}
            \pare{
                \partial_tU
                -
                \pare{
                    {\bm A}^{\bC,{\rm lin}}_{\rm lo}
                    +
                    \ii{\bm D}^{\bC}_{>0,{\rm lin}}
                }U
            }
        }{
            U_{s,{\rm lo}}
        }
    }
    \leq
    C_sE_{s_0}^{\frac12}E_s.
\]
Combining this bound with
\cref{eq:low_frequency_flat_linear_cancellation_3d} proves
\cref{eq:low_frequency_energy_bound_3d}.
\end{proofpart}
\end{proof}
We are now ready to prove \cref{prop:modified_energy_estimate_diagonalized_3d}.
\begin{proof}[Proof of \cref{prop:modified_energy_estimate_diagonalized_3d}]
By \cref{eq:modified_energy_definition_3d},
\[
    E_s
    =
    \frac12\norm{W_{s,{\rm hi}}}_{L^2}^2
    +
    \frac12
    \psc{
        {\bm S}^{\bC}_{\rm lo}U_{s,{\rm lo}}
    }{
        U_{s,{\rm lo}}
    }.
\]
Adding
\cref{eq:high_frequency_energy_bound_3d,eq:low_frequency_energy_bound_3d}
gives
\[
    \av{
        \frac{\dd}{\dd t}E_s
    }
    \leq
    C_sE_{s_0}^{\frac12}E_s,
\]
which is \cref{eq:modified_energy_closed_estimate_3d}.
\end{proof}

\begin{cor}[Uniform a priori estimate in the modified energy]
\label{cor:uniform_apriori_diagonalized_energy_3d}
Under the assumptions of
\cref{prop:modified_energy_estimate_diagonalized_3d}, set
\[
    e_0\pare{t}
    \defeq
    E_{s_0}^{\frac12}\pare{t}.
\]
Then
\begin{equation}\label{eq:low_energy_inverse_size_bound_3d}
    e_0\pare{t}
    \leq
    \frac{
        e_0\pare{0}
    }{
        1-C_{s_0}e_0\pare{0}\av{t}
    }
\end{equation}
whenever the denominator is positive.  Moreover, for every \(s\geq s_0\),
\begin{equation}\label{eq:high_energy_inverse_size_bound_3d}
    E_s\pare{t}
    \leq
    E_s\pare{0}
    \pare{
        1-C_{s_0}e_0\pare{0}\av{t}
    }^{-\frac{C_s}{C_{s_0}}}.
\end{equation}
\end{cor}

\begin{proof}
Taking \(s=s_0\) in
\cref{eq:modified_energy_closed_estimate_3d} gives
\[
    \av{
        \frac{\dd}{\dd t}e_0
    }
    \leq
    C_{s_0}e_0^2.
\]
Integrating forward and backward from \(t=0\) proves
\cref{eq:low_energy_inverse_size_bound_3d}.  For \(s\geq s_0\),
\[
    \av{
        \frac{\dd}{\dd t}E_s
    }
    \leq
    C_se_0E_s.
\]
Integrating again and using
\cref{eq:low_energy_inverse_size_bound_3d} proves
\cref{eq:high_energy_inverse_size_bound_3d}.
\end{proof}

\subsubsection{Equivalence between the modified and Sobolev energies}
\label{sec:modified_sobolev_energy_equivalence_3d}

For \(s\in\bR\) and \(Y=\vect{\eta}{\tilde\psi}\), define the Hamiltonian
low-high energy norm
\begin{equation}\label{eq:hamiltonian_phase_space_norm_3d}
\begin{aligned}
    \norm{Y}_{\mathcal X^s_{\rm Ham}
    \pare{\Atw,\background}}^2
    \defeq{}&
    \norm{
        \chi_{\rm hi}\pare{D}\eta
    }_{H^{s+\frac14}}^2
    +
    \norm{
        \chi_{\rm hi}\pare{D}\tilde\psi
    }_{H^{s-\frac14}}^2
    \\
    &+
    \norm{
        \chi_{\rm lo}\pare{D}
        \angles{D}^{s-\frac34}
        \pare{
            -\mathfrak q
            \pare{\Atw,\background\middle|D}
        }^{\frac12}
        \eta
    }_{L^2}^2
    \\
    &+
    \frac12
    \norm{
        \chi_{\rm lo}\pare{D}
        \angles{D}^{s-\frac34}
        \av{D}^{\frac12}
        \tilde\psi
    }_{L^2}^2 .
\end{aligned}
\end{equation}
We also isolate the purely inhomogeneous low-frequency component of
\(\tilde\psi\):
\begin{equation}\label{eq:low_frequency_potential_control_norm_3d}
    \norm{Y}_{\mathcal Z^s_{\rm lo}}^2
    \defeq
    \norm{
        \chi_{\rm lo}\pare{D}
        \angles{D}^{s-\frac14}
        \tilde\psi
    }_{L^2}^2 .
\end{equation}

\begin{remark}[Hamiltonian scale versus the inhomogeneous Sobolev scale]
The norm in \cref{eq:hamiltonian_phase_space_norm_3d} is the scale naturally
produced by the flat Hamiltonian symmetrizer.  It is fully inhomogeneous on
\(\eta\) and on the high-frequency part of \(\tilde\psi\), but at low
frequencies it measures \(\av{D}^{1/2}\tilde\psi\), not \(\tilde\psi\)
itself.  This is already present in the real flat Hamiltonian energy through
the block
\[
    \begin{bmatrix}
        -\mathfrak q\pare{\Atw,\background\middle|\xi} & 0
        \\
        0 & \av{\xi}/2
    \end{bmatrix}.
\]
Thus the missing component in the standard Sobolev norm \(\mathcal X^s\) is
precisely \(\norm{Y}_{\mathcal Z^s_{\rm lo}}\).  The symmetrized energy gives
a closed estimate in the Hamiltonian scale; the missing low-frequency
Sobolev component is propagated by a rough estimate. \fine
\end{remark}

\begin{remark}[Why the enlarged low cutoff is used]
The choice \(\chi_{\rm lo}\pare{\xi}=\chi\pare{\xi/2}\) guarantees that
\(\chi_{\rm lo}=1\) on the whole support of the low-frequency linear defect
\({\bm A}^{\bC,{\rm lin}}_{\rm lo}\).  Therefore the low-frequency
Hamiltonian symmetrizer sees the entire flat low block.  The transition
annulus between \(\chi\) and \(\chi_{\rm lo}\) is compact in frequency; all
symbols supported there are treated as bounded symbols in the energy
estimate. \fine
\end{remark}

\begin{lemma}[Coercivity of the modified energy]
\label{lem:lower_norms_controlled_by_modified_energy_3d}
Let \(K\Subset\mathcal S_{\rm st}\). There exist
\(r_0, s_0 > 0\), such that for every \(s\geq s_0\) and every real classical solution satisfying
\[
    E_{s_0}^{\frac12}\pare{t}\leq r_0,
\]
one has
\begin{align}
    E_s\pare{t}
    &\sim_{s,K}
    \norm{
        \pare{
            \eta\pare{t},
            \tilde\psi\pare{t}
        }
    }_{\mathcal X^s_{\rm Ham}\pare{\Atw,\background}}^2,
    \label{eq:modified_energy_equivalence_3d}
    \\
    \norm{
        \pare{
            \eta\pare{t},
            \tilde\psi\pare{t}
        }
    }_{\mathcal X^s}^2
    &\sim_{s,K}
    E_s\pare{t}
    +
    \norm{
        \pare{
            \eta\pare{t},
            \tilde\psi\pare{t}
        }
    }_{\mathcal Z^s_{\rm lo}}^2.
    \label{eq:modified_energy_full_sobolev_equivalence_3d}
\end{align}
In particular,
\begin{equation}\label{eq:small_factor_control_by_low_energy_3d}
    \norm{
        \pare{\eta\pare{t},\tilde\psi\pare{t}}
    }_{\mathcal X^{s_0}_{\rm Ham}\pare{\Atw,\background}}
    \leq
    C_KE_{s_0}^{\frac12}\pare{t}.
\end{equation}
\end{lemma}

\begin{proof}
On \(\supp\chi_{\rm hi}\), one has \(\chi=0\), hence
\(\mathfrak m^{\rm nh}_{-\frac14}=\mathfrak m_{-\frac14}\).  The ellipticity
of the diagonal entries gives
\[
    \lambda_s\pare{d_j}
    \sim_{s,K}
    \angles{\xi}^{s},
    \qquad
    j=1,2.
\]
Using the high-frequency form of \({\bm M}\), the tame invertibility of
\({\bm B}\), and the real Alinhac transformation, we obtain
\begin{equation}\label{eq:modified_energy_high_equivalence_3d}
    \norm{W_{s,{\rm hi}}}_{L^2}^2
    \sim_{s,K}
    \norm{
        \chi_{\rm hi}\pare{D}\eta
    }_{H^{s+\frac14}}^2
    +
    \norm{
        \chi_{\rm hi}\pare{D}\tilde\psi
    }_{H^{s-\frac14}}^2.
\end{equation}
For the low-frequency part, using \(U={\bm M}Y_{\rm Al}\) and
\cref{eq:complex_low_frequency_symmetrizer_3d}, one has
\begin{equation}\label{eq:modified_energy_low_pullback_3d}
\begin{aligned}
    \psc{
        {\bm S}^{\bC}_{\rm lo}U_{s,{\rm lo}}
    }{
        U_{s,{\rm lo}}
    }
    ={}&
    \norm{
        \chi_{\rm lo}\pare{D}
        \angles{D}^{s-\frac34}
        \pare{
            -\mathfrak q
            \pare{\Atw,\background\middle|D}
        }^{\frac12}
        \eta_{\rm Al}
    }_{L^2}^2
    +
    \frac12
    \norm{
        \chi_{\rm lo}\pare{D}
        \angles{D}^{s-\frac34}
        \av{D}^{\frac12}
        \tilde\psi_{\rm Al}
    }_{L^2}^2.
\end{aligned}
\end{equation}
By \cref{eq:q_uniform_ellipticity_nonhomogeneous_3d},
\cref{eq:modified_energy_low_pullback_3d} is uniformly coercive for
\(\pare{\Atw,\background}\in K\).
The Alinhac transformation is triangular and tame on the Hamiltonian
low-frequency scale.  More explicitly, the correction to \(\tilde\psi\) is
an order-zero coefficient times \(\eta\), and therefore
\[
    \chi_{\rm lo}\pare{D}
    \angles{D}^{s-\frac34}
    \av{D}^{\frac12}
    \pare{
        \tilde\psi_{\rm Al}-\tilde\psi
    }
\]
is controlled by the corresponding low-frequency \(\eta\)-term.  Hence
\cref{eq:modified_energy_low_pullback_3d} is equivalent to the low-frequency
part of \(\norm{Y}_{\mathcal X^s_{\rm Ham}\pare{\Atw,\background}}\).
Together with
\cref{eq:modified_energy_high_equivalence_3d}, this proves
\cref{eq:modified_energy_equivalence_3d}.

Finally, by the definition of the standard inhomogeneous space
\(\mathcal X^s\) in \cref{eq:main_phase_space}, the only missing term in the
Hamiltonian norm is the low-frequency inhomogeneous component of
\(\tilde\psi\).  Therefore
\[
    \norm{Y}_{\mathcal X^s}^2
    \sim_{s,K}
    \norm{Y}_{\mathcal X^s_{\rm Ham}\pare{\Atw,\background}}^2
    +
    \norm{Y}_{\mathcal Z^s_{\rm lo}}^2.
\]
Using \cref{eq:modified_energy_equivalence_3d} gives
\cref{eq:modified_energy_full_sobolev_equivalence_3d}, and
\cref{eq:small_factor_control_by_low_energy_3d} follows immediately.
\end{proof}

\begin{lemma}[Rough propagation of the missing low-frequency component]
\label{lem:rough_low_frequency_potential_propagation_3d}
Let \(K\Subset\mathcal S_{\rm st}\).  There exist \(s_0,r_0>0\) such that,
for every \(s\geq s_0\), every real classical solution satisfying
\[
    \sup_{t\in I}
    E_{s_0}^{\frac12}\pare{t}
    \leq r_0
\]
obeys
\begin{equation}\label{eq:rough_low_frequency_potential_differential_3d}
    \frac{\dd}{\dd t}
    \norm{Y\pare{t}}_{\mathcal Z^s_{\rm lo}}
    \leq
    C_s
    E_s^{\frac12}\pare{t}
    +
    C_s
    E_{s_0}^{\frac12}\pare{t}
    \norm{Y\pare{t}}_{\mathcal Z^s_{\rm lo}} .
\end{equation}
Consequently, if
\[
    \sup_{\av{t}\leq T}
    E_{s_0}^{\frac12}\pare{t}
    \leq
    M_0\varepsilon,
    \qquad
    \sup_{\av{t}\leq T}
    E_s^{\frac12}\pare{t}
    \leq
    M_s\varepsilon,
\]
then, for \(\av{t}\leq T\),
\begin{equation}\label{eq:rough_low_frequency_potential_integrated_3d}
    \norm{Y\pare{t}}_{\mathcal Z^s_{\rm lo}}
    \leq
    \exp\pare{C_sM_0\varepsilon\av{t}}
    \pare{
        \norm{Y_0}_{\mathcal Z^s_{\rm lo}}
        +
        C_sM_s\varepsilon\av{t}
    }.
\end{equation}
\end{lemma}

\begin{proof}
The CDE vector field depends on \(\tilde\psi\) only through differentiated
quantities.  At the level of the second equation, its flat low-frequency part
has the form
\[
    \partial_t\tilde\psi
    =
    \mathfrak q\pare{\Atw,\background\middle|D}\eta
    +
    \textnormal{terms at least linear in }Y
    \textnormal{ and involving differentiated }\tilde\psi .
\]
Applying \(\chi_{\rm lo}\pare{D}\angles{D}^{s-\frac14}\) gives
\[
    \norm{
        \chi_{\rm lo}\pare{D}
        \angles{D}^{s-\frac14}
        \mathfrak q\pare{\Atw,\background\middle|D}\eta
    }_{L^2}
    \leq
    C_{s,K}
    \norm{Y}_{\mathcal X^s_{\rm Ham}
    \pare{\Atw,\background}},
\]
because \(\mathfrak q\) is a smooth bounded multiplier on the support of
\(\chi_{\rm lo}\), while the Hamiltonian low-frequency norm controls
\[
    \chi_{\rm lo}\pare{D}
    \angles{D}^{s-\frac34}
    \pare{-\mathfrak q\pare{\Atw,\background\middle|D}}^{\frac12}
    \eta .
\]
For the nonlinear terms, the tame paralinearization estimates of
\cref{thm:paraCDEetapsi} and the real Alinhac structure give
\[
    \norm{
        \chi_{\rm lo}\pare{D}
        \angles{D}^{s-\frac14}
        \mathfrak N\pare{Y}
    }_{L^2}
    \leq
    C_s
    \norm{Y}_{\mathcal X^{s_0}_{\rm Ham}
    \pare{\Atw,\background}}
    \pare{
        \norm{Y}_{\mathcal X^s_{\rm Ham}
        \pare{\Atw,\background}}
        +
        \norm{Y}_{\mathcal Z^s_{\rm lo}}
    }.
\]
The appearance of \(\mathcal Z^s_{\rm lo}\) here accounts only for bounded
zero-order low-frequency contributions.  Using
\cref{eq:modified_energy_equivalence_3d}, we obtain
\cref{eq:rough_low_frequency_potential_differential_3d}.  The integrated
bound \cref{eq:rough_low_frequency_potential_integrated_3d} follows by
Gronwall's lemma.
\end{proof}

\subsubsection{Proof of the uniform Sobolev estimate}
\label{sec:proof_uniform_sobolev_energy_estimate_3d}

\begin{proof}[Proof of \cref{prop:uniform_apriori_sobolev_energy_3d}]
By the bootstrap smallness assumption in \(\mathcal X^{s_0}\) and the
coercivity estimate \cref{eq:modified_energy_full_sobolev_equivalence_3d},
after decreasing \(r_0\) we may apply
\cref{prop:modified_energy_estimate_diagonalized_3d}.  Set
\[
    e_s\pare t
    \defeq
    E_s^{\frac12}\pare{t},
    \qquad
    z_s\pare t
    \defeq
    \norm{Y\pare t}_{\mathcal Z^s_{\rm lo}}.
\]
The initial bound in \(\mathcal X^s\), together with
\cref{eq:modified_energy_full_sobolev_equivalence_3d}, gives
\begin{equation}\label{eq:sobolev_proof_initial_energy_comparison_3d}
    e_s\pare0+z_s\pare0
    \leq
    C_{s,K}\varepsilon,
    \qquad
    e_{s_0}\pare0
    \leq
    C_{s_0,K}\varepsilon .
\end{equation}
Choose \(c_{s,K}>0\) and \(\varepsilon_{s,K}>0\) so small that
\[
    C_{s_0,K}e_{s_0}\pare0\av t
    \leq
    \frac12
    \qquad
    \textnormal{whenever }\av t\leq c_{s,K}\varepsilon^{-1}.
\]
Then \cref{cor:uniform_apriori_diagonalized_energy_3d} gives, for every
\(t\in I\) with \(\av t\leq c_{s,K}\varepsilon^{-1}\),
\begin{equation}\label{eq:sobolev_proof_modified_energy_bound_3d}
    e_s\pare t
    \leq
    C_{s,K}\varepsilon,
    \qquad
    e_{s_0}\pare t
    \leq
    C_{s_0,K}\varepsilon .
\end{equation}
In particular, by \cref{eq:modified_energy_equivalence_3d},
\cref{eq:sobolev_proof_modified_energy_bound_3d} gives the Hamiltonian bound
\cref{eq:uniform_hamiltonian_apriori_lifespan_form_3d}.

It remains to control the undifferentiated low-frequency component of
\(\tilde\psi\).  Applying
\cref{eq:rough_low_frequency_potential_integrated_3d} with
\[
    T=c_{s,K}\varepsilon^{-1},
    \qquad
    M_0=C_{s_0,K},
    \qquad
    M_s=C_{s,K},
\]
and using \cref{eq:sobolev_proof_initial_energy_comparison_3d}, we obtain
\begin{equation}\label{eq:sobolev_proof_low_frequency_bound_3d}
    z_s\pare t
    \leq
    C_{s,K}
    \pare{
        \varepsilon
        +
        \varepsilon\av t
    }
    \leq
    C_{s,K}\pare{\varepsilon+c_{s,K}}
\end{equation}
for every \(t\in I\) with \(\av t\leq c_{s,K}\varepsilon^{-1}\), after
possibly decreasing \(c_{s,K}\) once more to keep the exponential factor
uniformly bounded.  Combining
\cref{eq:modified_energy_full_sobolev_equivalence_3d,eq:sobolev_proof_modified_energy_bound_3d,eq:sobolev_proof_low_frequency_bound_3d}
proves \cref{eq:uniform_sobolev_apriori_lifespan_form_3d}.  The final
small-neighborhood statement follows by choosing \(c_{s,K}\) and
\(\varepsilon_{s,K}\) sufficiently small.
\end{proof}

\subsection{Proof of \cref{thm:main}}
\label{sec:proof_main_theorem_3d}

\begin{proof}
Fix  a compact subset \(K\) of the stable parameter set \(\mathcal S_{\rm st}\) introduced in \cref{def:main_stable_parameter_set}. Write
$$p\defeq\pare{\Atw,\background}\in K.$$
We construct the proof of \Cref{thm:main} following a standard approximation procedure for quasilinear-hyperbolic systems. 

\begin{proofpart}[Analyticity]
Recall that we use the notation $Y\defeq\begin{bmatrix}
    \eta\\
    \tilde{\psi}
\end{bmatrix}.$ Denote by \(\mathfrak F_p\pare Y\) the right-hand side of \cref{eq:CDE}.
Its formulas combine pointwise analytic functions, Fourier multipliers,
the kernels in \cref{eq:B0,eq:Q}, analytic in the normalized differences
of \(\eta\), and the analytic inverse
\(\pare{\Id+\Atw\tau_0\pare\eta}^{-1}\) of
\cref{lem:uniform_inverse_tau0}.
Thus, for \(s_0\) sufficiently large, \(r_\ast>0\) sufficiently small,
and \(q\geq s_0\),
\begin{equation}\label{eq:main_proof_analytic_exact_vector_field_3d}
    \mathfrak F_p
    \in
    \cC^\omega\pare{
        B_{\mathcal X^q}\pare{0,r_\ast};
        \mathcal X^{q-\frac32}
    },
\end{equation}
with uniform tame estimates for \(p\in K\).
The smallness of
\(\norm{\eta}_{H^{s_0+\frac14}}\) keeps the normalized kernel factors
and the inverse controlled when differentiating these formulas.

Define
\begin{equation}\label{eq:main_proof_nonlinear_diagonalizing_map_3d}
\begin{aligned}
    U_p\pare Y
    &\defeq
    {\bm M}\pare{p\middle|D}
    {\bm G}\pare{p\middle|Y}Y,
    \\
    \Phi_p\pare Y
    &\defeq
    {\bm B}\pare{p\middle|U_p\pare Y}^{-1}
    U_p\pare Y .
\end{aligned}
\end{equation}
By
\cref{eq:real_alinhac_unknown_3d,eq:complex_variable_definition_3d,eq:positive_order_diagonalized_unknown_3d},
the maps entering \cref{eq:main_proof_nonlinear_diagonalizing_map_3d}
are real-analytic and tame for
\(\norm{Y}_{\mathcal X^{s_0}}<r_\ast\), uniformly for \(p\in K\).
Hence, for \(q\geq s_0\),
\[
    \Phi_p
    \in
    \cC^\omega\pare{
        B_{\mathcal X^q}\pare{0,r_\ast};
        H^q
    }
\]
on the real subspace, with uniform tame bounds. Moreover,
\[
    \Phi_p\pare0=0,
    \qquad
    \dd\Phi_p\pare0
    =
    {\bm M}\pare{p\middle|D},
\]
and \({\bm M}\pare{p\middle|D}\) is a uniform isomorphism from
\(\mathcal X^q\) onto the real subspace of
\(H^q\pare{\bR^2;\bC^2}\).
After decreasing \(r_\ast\), the analytic inverse function theorem
therefore gives a real-analytic tame inverse
\[
    \Psi_p
    \defeq
    \Phi_p^{-1},
\]
uniformly for \(p\in K\).

We use \(\Phi_p\) as a nonlinear change of coordinates and set
\[
    W\defeq\Phi_p\pare Y,
    \qquad
    Y=\Psi_p\pare W.
\]
Accordingly, the vector field of the canonical CDE
\cref{eq:CDE} in the \(W\)-coordinates is
\begin{equation}\label{eq:main_proof_exact_transformed_vector_field_3d}
    \mathfrak S_p\pare W
    \defeq
    \dd\Phi_p\pare{\Psi_p\pare W}
    \bra{
        \mathfrak F_p\pare{\Psi_p\pare W}
    } .
\end{equation}
Since \(\mathfrak F_p\), \(\Phi_p\), and \(\Psi_p\) are real-analytic
and tame, for every \(q\geq s_0\),
\begin{equation}\label{eq:main_proof_analytic_transformed_vector_field_3d}
    \mathfrak S_p
    \in
    \cC^\omega\pare{
        B_{H^q}\pare{0,r_\ast};
        H^{q-\frac32}
    },
\end{equation}
on the real subspace, with uniform tame bounds.

If \(Y\) solves the canonical CDE
\[
    \partial_tY=\mathfrak F_p\pare Y,
\]
then, using \(W=\Phi_p\pare Y\) and the chain rule,
\begin{equation}\label{eq:evolution-W_LWP}
    \partial_tW
    =
    \dd\Phi_p\pare Y
    \bra{\partial_tY}
    =
    \dd\Phi_p\pare Y
    \bra{\mathfrak F_p\pare Y}
    =
    \mathfrak S_p\pare{\Phi_p\pare Y}
    =
    \mathfrak S_p\pare W .
\end{equation}
Conversely, since \(Y=\Psi_p\pare W\), the same argument with the inverse
change of coordinates shows that
\[
    \partial_tY=\mathfrak F_p\pare Y
    \quad\Longleftrightarrow\quad
    \partial_tW=\mathfrak S_p\pare W .
\]

We now write explicitly the quasilinear operator associated with
\(\mathfrak S_p\). Set
\[
    \mathscr U_p\pare W
    \defeq
    U_p\pare{\Psi_p\pare W},
\]
and, for a sufficiently regular \(Z\), define
\begin{equation}\label{eq:main_proof_frozen_operator_3d}
\begin{aligned}
    \mathscr L_p\pare W Z
    \defeq{}&
    \OpBW{
        {\bm A}^{\bC,{\rm lin}}_{\rm lo}
        \pare{p\middle|\xi}
    }Z
    \\
    &+
    \ii
    \OpBW{
        {\bm D}^{\bC}_{>0}
        \pare{p\middle|\mathscr U_p\pare W;x,\xi}
        +
        {\bm A}_{\bra{0,\geq1}}
        \pare{p\middle|\mathscr U_p\pare W;x,\xi}
    }Z
    \\
    &+
    {\bm R}_{\bra{\varrho,\geq1}}
    \pare{p\middle|\mathscr U_p\pare W}Z .
\end{aligned}
\end{equation}
The paradifferential representation of
\cref{eq:evolution-W_LWP} is precisely
\cref{eq:positive_order_diagonalized_system_3d}. By \cref{eq:main_proof_analytic_transformed_vector_field_3d},
\(\mathscr L_p\pare W\) depends real-analytically and tamely on \(W\),
and
\begin{equation}\label{eq:main_proof_vector_field_factorization_3d}
    \mathfrak S_p\pare W
    =
    \mathscr L_p\pare W\bra W .
\end{equation}
\end{proofpart}

\begin{proofpart}[The regularized problem]
We now construct the approximate solutions.  Let
\(\tilde{\chi}\in\cC_c^\infty\pare{\bR;\bra{0,1}}\) be real and even, with
\(\tilde{\chi}=1\) on \(\bra{-1,1}\), and set
\begin{equation}\label{eq:main_proof_fourier_truncation_3d}
    P_N
    \defeq
    \tilde{\chi}\pare{\frac{\angles D}{2^N}} .
\end{equation}
The operators \(P_N\) are real, self-adjoint smoothing Fourier
multipliers and commute with Fourier multipliers.
Let
\[
    {\bm D}_p\pare W
    \defeq
    \operatorname{diag}\pare{
        d_1\pare{W;x,\xi},
        d_2\pare{W;x,\xi}
    }
\]
denote the real diagonal positive-order part of
\(\mathscr L_p\pare W\) in the sense of \cref{eq:first_diagonal_entry_positive_symbol_3d,eq:second_diagonal_entry_positive_symbol_3d}, and write
\begin{equation}\label{eq:main_proof_linearized_bounded_part_3d}
    \mathscr B_p\pare W\bra Z
    \defeq
    \mathscr L_p\pare W\bra Z
    -
    \ii\OpBW{{\bm D}_p\pare W}Z .
\end{equation}
Thus \(\mathscr B_p\pare W\) contains the bounded part of the
linearized transformed vector field.
For \(N\geq0\), define the regularized diagonal symbols
\begin{equation}\label{eq:main_proof_regularized_diagonal_symbol_3d}
\begin{aligned}
    f_{N+1}\pare z
    &\defeq
    z\,
    \tilde{\chi}^2\pare{
        \frac{z}{2^{\frac32(N+1)}}
    },
    \\
    \widetilde d_{j,N+1}\pare{W;x,\xi}
    &\defeq
    f_{N+1}\pare{
        d_j\pare{W;x,\xi}
    },
    \qquad j=1,2,
    \\
    \widetilde{\bm D}_{p,N+1}\pare W
    &\defeq
    \operatorname{diag}\pare{
        \widetilde d_{1,N+1}\pare W,
        \widetilde d_{2,N+1}\pare W
    } .
\end{aligned}
\end{equation}
Since the cutoff is real and even, the reality relation between the two
diagonal branches is preserved.
Set
\begin{equation}\label{eq:main_proof_regularized_iteration_3d_0}
    W_{\rm in}\defeq\Phi_p\pare{Y_0},
    \qquad
    W_0\pare t\equiv0.
\end{equation}
Having constructed \(W_N\), define \(W_{N+1}\) as the solution of
\begin{equation}\label{eq:main_proof_regularized_iteration_3d}
\begin{cases}
    \partial_tW_{N+1}
    =
    \ii\OpBW{
        \widetilde{\bm D}_{p,N+1}\pare{W_N}
    }W_{N+1}
    +
    P_{N+1}
    \mathscr B_p\pare{W_N}
    \bra{
       W_{N+1}
    },
    \\[2mm]
    W_{N+1}\pare0
    =
    P_{N+1}W_{\rm in}.
\end{cases}
\end{equation}
For fixed \(N\), \cref{eq:main_proof_regularized_iteration_3d} is
linear in \(W_{N+1}\).  The symbol
\(\widetilde{\bm D}_{p,N+1}\pare{W_N}\) has compact frequency support,
so its Bony--Weyl quantization is bounded on every fixed Sobolev space,
with a bound which may depend on \(N\).  The second term in
\cref{eq:main_proof_regularized_iteration_3d} is bounded on the same
space by the Friedrichs regularization.

Since \(W_N\) is known and continuous in time, the right-hand side of
\cref{eq:main_proof_regularized_iteration_3d} is a continuous family of
bounded linear operators on \(H^s\).  Hence the linear
Cauchy--Lipschitz theorem gives a unique solution \(W_{N+1}\) on a
possibly \(N\)-dependent time interval.  The real subspace is invariant,
and compatibility of the same equation at higher Sobolev indices gives
smooth approximate solutions for the smoothed initial data.
\end{proofpart}

\begin{proofpart}[The uniform bounds]

We prove that the solutions constructed in Part~2 exist on a common
time interval, independent of \(N\).

For \(N\geq0\), set
\[
    U_N
    \defeq
    \mathscr U_p\pare{W_N}.
\]
For a real function \(Z\), let
\(\mathcal E_{q,N}\bra Z\) denote the quadratic energy obtained from
\cref{eq:modified_energy_definition_3d} by freezing the
state-dependent coefficients at \(U_N\) and replacing \(W\) by \(Z\).
In particular, its high-frequency weight is
\begin{equation}\label{eq:main_proof_frozen_modified_weight_3d}
    {\bm\Lambda}^{N}_{q,{\rm hi}}
    \defeq
    \chi_{\rm hi}\pare{\xi}
    \operatorname{diag}\pare{
        \lambda_q\pare{d_1\pare{U_N;x,\xi}},
        \lambda_q\pare{d_2\pare{U_N;x,\xi}}
    }.
\end{equation}
The low-frequency part is obtained by freezing the change
\(W\mapsto U\) at \(U_N\) and using the flat symmetrizer
\({\bm S}^{\bC}_{\rm lo}\) of
\cref{eq:complex_low_frequency_symmetrizer_3d}.
The coercivity argument of
\cref{lem:lower_norms_controlled_by_modified_energy_3d}
is uniform with respect to \(N\).

We verify that the regularization of Part~2 preserves the estimates used
in the proof of \cref{lem:adapted_diagonal_weights_3d}.

The crucial point is that the principal high-frequency cancellation is
preserved \emph{exactly}.  We use in
\cref{eq:main_proof_frozen_modified_weight_3d} the unregularized weight
\(\lambda_q\pare{d_j}\).  For \(j=1,2\), the potentially order-losing
symbol in the commutator with the regularized diagonal operator is
\[
    \frac{1}{\ii}
    \pbra{
        \chi_{\rm hi}\pare{\xi}\lambda_q\pare{d_j}
    }{
        \widetilde d_{j,N+1}
    }.
\]
Where \(\chi_{\rm hi}\) is constant, one has
\begin{equation}\label{eq:main_proof_regularized_poisson_cancellation_3d}
    \pbra{
        \lambda_q\pare{d_j}
    }{
        \widetilde d_{j,N+1}
    }
    =
    \pbra{
        \lambda_q\pare{d_j}
    }{
        f_{N+1}\pare{d_j}
    }
    =
    \lambda_q'\pare{d_j}
    f_{N+1}'\pare{d_j}
    \pbra{d_j}{d_j}
    =
    0.
\end{equation}
Thus the regularization preserves exactly the cancellation in
\cref{eq:modified_energy_poisson_cancellation_3d}.  When a
\(\xi\)-derivative hits \(\chi_{\rm hi}\), the resulting symbol is
supported in its fixed transition annulus and contains
\(\nabla_x d_j\).  In particular, it vanishes at \(U=0\) and is estimated
exactly as in the proof of \cref{lem:adapted_diagonal_weights_3d}.

Finally, the bounded part of
\cref{eq:main_proof_regularized_iteration_3d} is harmless:
\(P_{N+1}\) is real, self-adjoint and uniformly bounded on every Sobolev
space.  We choose once and for all the initial regularization scale in
Part~2 so that, for every \(N\geq0\), both \(P_{N+1}\) and
\(d_j\mapsto\widetilde d_{j,N+1}\) are the identity on the fixed
low-frequency region entering
\cref{eq:complex_low_frequency_symmetrizer_3d}.  Hence the flat
low-frequency cancellation
\cref{eq:low_frequency_flat_linear_cancellation_3d} is unchanged.

We may therefore repeat the proof of
\cref{lem:adapted_diagonal_weights_3d} with the coefficients frozen at
\(U_N\).  The only difference is that the coefficients depend on
\(W_N\), whereas the equation is solved by \(W_{N+1}\).
When the time derivative falls on a frozen coefficient it produces
\(\partial_tW_N\); using the equation in  Part~2 and the same
tame estimates as in
\cref{eq:time_derivative_symmetrizer_bound_3d} gives, for every
\(q\geq s_0\),
\begin{equation}\label{eq:main_proof_frozen_energy_estimate_3d}
\begin{aligned}
    \frac{\dd}{\dd t}
    \mathcal E_{q,N}^{\frac12}\bra{W_{N+1}}
    \leq{}&
    C_{q,K}
    \mathcal E_{s_0,N}^{\frac12}\bra{W_N}
    \mathcal E_{q,N}^{\frac12}\bra{W_{N+1}}+
    C_{q,K}
    \mathcal E_{q,N}^{\frac12}\bra{W_N}
    \mathcal E_{s_0,N}^{\frac12}\bra{W_{N+1}}.
\end{aligned}
\end{equation}
In particular, for \(q=s_0\),
\begin{equation}\label{eq:main_proof_frozen_low_energy_estimate_3d}
    \frac{\dd}{\dd t}
    \mathcal E_{s_0,N}^{\frac12}\bra{W_{N+1}}
    \leq
    C_{s_0,K}
    \mathcal E_{s_0,N}^{\frac12}\bra{W_N}
    \mathcal E_{s_0,N}^{\frac12}\bra{W_{N+1}}.
\end{equation}

The initial data satisfy, uniformly in \(N\),
\[
    \mathcal E_{s,N}^{\frac12}
    \bra{W_{N+1}\pare0}
    \leq
    C_{s,K}\varepsilon,
\]
by the tame bounds for \(\Phi_p\) and the uniform boundedness of
\(P_{N+1}\).  We set
\begin{equation}\label{eq:main_proof_target_interval_3d}
    T_\varepsilon
    \defeq
    c_{s,K}\varepsilon^{-1},
    \qquad
    I_\varepsilon
    \defeq
    \bra{-T_\varepsilon,T_\varepsilon}.
\end{equation}
We first apply
\cref{eq:main_proof_frozen_low_energy_estimate_3d}.
Assuming inductively that
\[
    \sup_{t\in I_\varepsilon}
    \mathcal E_{s_0,N}^{\frac12}
    \bra{W_N\pare t}
    \leq
    M_0\varepsilon,
\]
Gronwall's inequality gives
\[
    \sup_{t\in I_\varepsilon}
    \mathcal E_{s_0,N}^{\frac12}
    \bra{W_{N+1}\pare t}
    \leq
    C_K\varepsilon
    \exp\pare{C_KM_0\varepsilon T_\varepsilon}.
\]
Choosing first \(M_0\) and then \(c_{s,K}\) sufficiently small closes
the induction.  With this lower bound available,
\cref{eq:main_proof_frozen_energy_estimate_3d} gives in the same way,
at the index \(s\),
\begin{equation}\label{eq:main_proof_uniform_modified_energy_bounds_3d}
    \sup_{N\geq0}
    \sup_{t\in I_\varepsilon}
    \mathcal E_{s,N}^{\frac12}
    \bra{W_{N+1}\pare t}
    \leq
    C_{s,K}\varepsilon.
\end{equation}

It remains only to propagate the undifferentiated low-frequency
potential, which is not contained in the Hamiltonian energy.
The proof of
\cref{lem:rough_low_frequency_potential_propagation_3d}
applies to the frozen equation in the same way.  If
\(z_{q,N+1}\) denotes the corresponding frozen version of
\cref{eq:low_frequency_potential_control_norm_3d}, its tame
polarization gives
\begin{equation}\label{eq:main_proof_frozen_low_frequency_bound_3d}
\begin{aligned}
    \frac{\dd}{\dd t}z_{q,N+1}
    \leq{}&
    C_{q,K}
    \mathcal E_{q,N}^{\frac12}\bra{W_{N+1}}
    +
    C_{q,K}
    \mathcal E_{s_0,N}^{\frac12}\bra{W_N}
    z_{q,N+1}
    \\
    &+
    C_{q,K}
    \mathcal E_{q,N}^{\frac12}\bra{W_N}
    z_{s_0,N+1}.
\end{aligned}
\end{equation}
Using first \(q=s_0\), and then \(q=s\), together with
\cref{eq:main_proof_uniform_modified_energy_bounds_3d}, yields
\begin{equation}\label{eq:main_proof_uniform_low_frequency_bounds_3d}
    \sup_{N\geq0}
    z_{s,N+1}\pare t
    \leq
    C_{s,K}
    \pare{
        \varepsilon+\varepsilon\av t
    },
    \qquad
    t\in I_\varepsilon.
\end{equation}
The frozen analogue of
\cref{eq:modified_energy_full_sobolev_equivalence_3d} therefore gives
\begin{equation}\label{eq:main_proof_uniform_smooth_bounds_3d}
    \sup_{N\geq0}
    \sup_{t\in I_\varepsilon}
    \norm{W_N\pare t}_{H^s}
    \leq
    C_{s,K}
    \pare{
        \varepsilon+c_{s,K}
    }.
\end{equation}
After decreasing \(c_{s,K}\) and \(\varepsilon_0\), the right-hand side
of \cref{eq:main_proof_uniform_smooth_bounds_3d} remains inside the
coordinate neighborhood on which all the transformations and
coefficients are defined.  The continuation criterion for the linear
problems of Part~2 therefore extends every \(W_N\) to the whole interval
\(I_\varepsilon\).   
\end{proofpart}

\begin{proofpart}[Convergence of approximate solutions]
We prove that the sequence \(\pare{W_N}_N\) is Cauchy in
\(\cC\pare{I_\varepsilon;L^2}\).  For convenience, denote by
\begin{equation}\label{eq:main_proof_regularized_linear_operator_3d}
    \mathscr L_{p,N}\pare W\bra Z
    \defeq
    \ii\OpBW{
        \widetilde{\bm D}_{p,N}\pare W
    }Z
    +
    P_N
    \mathscr B_p\pare W
    \bra{P_NZ}
\end{equation}
the regularized linear operator introduced in Part~2.  Thus
\[
    \partial_tW_{N+1}
    =
    \mathscr L_{p,N+1}\pare{W_N}\bra{W_{N+1}}.
\]

Set
\[
    Z_N\defeq W_{N+1}-W_N .
\]
Subtracting the equations at the two consecutive levels gives
\begin{equation}\label{eq:main_proof_L2_difference_equation_3d}
\begin{aligned}
    \partial_tZ_N
    ={}&
    \mathscr L_{p,N+1}\pare{W_N}\bra{Z_N}
    \\
    &+
    \pare{
        \mathscr L_{p,N+1}\pare{W_N}
        -
        \mathscr L_{p,N+1}\pare{W_{N-1}}
    }\bra{W_N}
    \\
    &+
    \pare{
        \mathscr L_{p,N+1}\pare{W_{N-1}}
        -
        \mathscr L_{p,N}\pare{W_{N-1}}
    }\bra{W_N}.
\end{aligned}
\end{equation}
The first line is the homogeneous linearized equation.  Since the
diagonal symbols \(\widetilde d_{j,N+1}\) are real, their Bony--Weyl
quantizations are self-adjoint and hence
\[
    \Re\angles{
        \ii\OpBW{
            \widetilde{\bm D}_{p,N+1}\pare{W_N}
        }Z_N,
        Z_N
    }_{L^2}
    =0.
\]
The bounded part contributes at most
\(C\norm{Z_N}_{L^2}^2\), uniformly in \(N\).

By the tame Lipschitz estimates for the coefficients and the uniform
bounds of Part~3, for \(W,V\) in the structural ball and \(X\in H^s\),
\begin{equation}\label{eq:main_proof_L2_state_difference_bound_3d}
    \norm{
        \pare{
            \mathscr L_{p,N+1}\pare{W}
            -
            \mathscr L_{p,N+1}\pare V
        }\bra X
    }_{L^2}
    \leq
    C_{s,K}
    \norm{W-V}_{L^2}
    \norm X_{H^s}.
\end{equation}

For the regularization error, set
\[
    a_{j,N}\pare W
    \defeq
    f_{N+1}\pare{d_j\pare W}
    -
    f_N\pare{d_j\pare W},
    \qquad j=1,2.
\]
Since \(f_{N+1}-f_N\) vanishes for
\(\av z\lesssim 2^{\frac32N}\), the high-frequency ellipticity of
\(d_j\) on the $\BallR{s_0}{r_0}$ implies
\[
    \supp a_{j,N}\pare W
    \subset
    \set{
        \angles\xi\gtrsim_K 2^N
    }.
\]
Together with the uniform order-\(3/2\) symbol bounds proved in Part~3,
this gives
\begin{equation}\label{eq:main_proof_diagonal_regularization_difference_bound_3d}
    \norm{
        \OpBW{
            a_{j,N}\pare W
        }X
    }_{L^2}
    \leq
    C_{s,K}
    2^{-\pare{s-\frac32}N}
    \norm X_{H^s},
    \qquad j=1,2.
\end{equation}
For the bounded part, one has
\[
\begin{aligned}
    &P_{N+1}\mathscr B_p\pare W P_{N+1}
    -
    P_N\mathscr B_p\pare W P_N
    \\
    &\qquad=
    \pare{P_{N+1}-P_N}
    \mathscr B_p\pare W P_{N+1}
    +
    P_N\mathscr B_p\pare W
    \pare{P_{N+1}-P_N}.
\end{aligned}
\]
Since
\[
    \norm{
        \pare{P_{N+1}-P_N}X
    }_{L^2}
    \leq
    C_s2^{-sN}\norm X_{H^s},
\]
and \(\mathscr B_p\pare W\) is uniformly bounded on the Sobolev scale,
\[
    \norm{
        \pare{
            P_{N+1}\mathscr B_p\pare W P_{N+1}
            -
            P_N\mathscr B_p\pare W P_N
        }X
    }_{L^2}
    \leq
    C_{s,K}2^{-sN}\norm X_{H^s}.
\]
Consequently,
\begin{equation}\label{eq:main_proof_L2_regularization_difference_bound_3d}
    \norm{
        \pare{
            \mathscr L_{p,N+1}\pare W
            -
            \mathscr L_{p,N}\pare W
        }\bra X
    }_{L^2}
    \leq
    C_{s,K}
    2^{-\pare{s-\frac32}N}
    \norm X_{H^s}.
\end{equation}
Using Part~3 in
\cref{eq:main_proof_L2_difference_equation_3d} therefore gives for some $\sigma>0,$
\begin{equation}\label{eq:main_proof_L2_difference_estimate_3d}
    \frac{\dd}{\dd t}\norm{Z_N}_{L^2}
    \leq
    C_{s,K}\norm{Z_N}_{L^2}
    +
    C_{s,K}\norm{Z_{N-1}}_{L^2}
    +
    C_{s,K}2^{-\sigma N}.
\end{equation}
Moreover,
\[
    Z_N\pare0
    =
    \pare{P_{N+1}-P_N}W_{\rm in},
\]
and hence, after decreasing \(\sigma>0\) if necessary,
\begin{equation}\label{eq:main_proof_L2_initial_difference_3d}
    \norm{Z_N\pare0}_{L^2}
    \leq
    C_{s,K}2^{-\sigma N}.
\end{equation}
Integrating \cref{eq:main_proof_L2_difference_estimate_3d} and using
\cref{eq:main_proof_L2_initial_difference_3d}, for \(0\leq t\leq
T_\varepsilon\), we obtain
\[
    A_N\pare t
    \leq
    C_{s,K}2^{-\sigma N} (1+t)
    +
    C_{s,K}\int_0^t A_N\pare\tau\dd\tau
    +
    C_{s,K}\int_0^t A_{N-1}\pare\tau\dd\tau ,
\]
where
\[
    A_N\pare t
    \defeq
    \sup_{0\leq \tau\leq t}
    \norm{Z_N\pare\tau}_{L^2}.
\]
Summing in \(N\geq1\) and applying Gronwall's inequality yields
\[
    \sum_{N\geq1}
    \sup_{0\leq t\leq T_\varepsilon}
    \norm{Z_N\pare t}_{L^2}
    <+\infty .
\]
The same argument on \([-T_\varepsilon,0]\) gives, recalling the \cref{notations:intro}-3,
\begin{equation}\label{eq:main_proof_L2_summability_3d}
    \sum_{N\geq1}
    \norm{W_{N+1}-W_N}_{0,I_\varepsilon}
    <+\infty .
\end{equation}
Consequently there exists a real-valued
\[
    W\in\cC\pare{I_\varepsilon;L^2}
\]
such that
\begin{equation}\label{eq:main_proof_L2_convergence_3d}
    W_N\longrightarrow W
    \qquad\textnormal{in}\qquad
    \cC\pare{I_\varepsilon;L^2}.
\end{equation}

By Part~3 the sequence is uniformly bounded in
\(L^\infty\pare{I_\varepsilon;H^s}\).  Therefore, for every
\(0<q<s\), Sobolev interpolation gives
\[
    \norm{W_N-W_M}_{q,I_\varepsilon}
    \leq
    C
    \norm{W_N-W_M}_{0,I_\varepsilon}^{1-\frac qs}
    \pare{
        \norm{W_N}_{L^\infty H^s}
        +
        \norm{W_M}_{L^\infty H^s}
    }^{\frac qs}.
\]
Hence
\begin{equation}\label{eq:main_proof_intermediate_convergence_3d}
    W_N\longrightarrow W
    \qquad\textnormal{in}\qquad
    \cC\pare{I_\varepsilon;H^q},
    \qquad
    0<q<s.
\end{equation}
In particular, choosing \(s_0\leq q<s\) and using the local
Lipschitz continuity of \(\Psi_p:H^q\to\mathcal X^q\), we obtain, with
\[
    Y_N\defeq\Psi_p\pare{W_N},
    \qquad
    Y\defeq\Psi_p\pare W,
\]
that
\begin{equation}\label{eq:main_proof_intermediate_physical_convergence_3d}
    Y_N\longrightarrow Y
    \qquad\textnormal{in}\qquad
    \cC\pare{I_\varepsilon;\mathcal X^q},
    \qquad
    s_0\leq q<s.
\end{equation}
\end{proofpart}

\begin{proofpart}[The limit is a solution of the original equation]
Fix \(q\) with $s_0\leq q<s$.
By Part~4,
\begin{equation}\label{eq:main_proof_strong_lower_convergence_3d}
    W_N\longrightarrow W
    \qquad\textnormal{in}\qquad
    \cC\pare{I_{\varepsilon};H^q}.
\end{equation}
We pass to the limit in the integral form of
\cref{eq:main_proof_regularized_iteration_3d},
\begin{equation}\label{eq:main_proof_approximate_integral_equation_3d}
    W_{N+1}\pare t
    =
    P_{N+1}W_{\rm in}
    +
    \int_0^t
    \mathscr L_{p,N+1}\pare{W_N\pare\tau}
    \bra{W_{N+1}\pare\tau}
    \dd\tau .
\end{equation}

We first remove the regularization.  The uniform \(H^s\)-bounds of
Part~3, the Fourier localization of \(P_N-\Id\), and the definition
\(\widetilde d_{j,N}=f_N\pare{d_j}\) give
\begin{equation}\label{eq:main_proof_regularization_vanishes_3d}
\begin{aligned}
    &
    \norm{
        \mathscr L_{p,N+1}\pare{W_N}
        \bra{W_{N+1}}
        -
        \mathscr L_p\pare{W_N}
        \bra{W_{N+1}}
    }_{q-\frac32,I_{\varepsilon}}
   \leq
    C_{q,s,K}\,
    2^{-\pare{s-q}N}
    \norm{W_{N+1}}_{s,I_{\varepsilon}}
    \xrightarrow{N\to \infty} 0 .
\end{aligned}
\end{equation}
For the diagonal part this follows from the fact that
\(f_N\pare{d_j}-d_j\) is supported where
\(\angles\xi\gtrsim2^N\), since
\(\av{d_j}\sim\angles\xi^{3/2}\); the bounded part is treated directly
with \(P_N-\Id\).

 By
\cref{prop:positive_order_diagonalization_3d}, the concrete operator
\(\mathscr L_p\pare V\) consists of a symbol of order \(3/2\), bounded
symbols, the fixed low-frequency multiplier, and a smoothing remainder.
The symbol estimates in
\cref{def:sfr,simb-pro-geq-q-geq-q}, the Bony--Weyl boundedness
\cref{thm:contS,cont00}, and the smoothing estimate
\cref{def:smoothing,eq:smoothing-tame-geq-q} therefore give, uniformly
for \(V\) in $B_{s_0,\mathbb{R}}\pare{I_{\varepsilon},r_0}$,
\begin{equation}\label{eq:main_proof_unregularized_operator_uniform_bound_3d}
    \norm{
        \mathscr L_p\pare V\bra Z
    }_{H^{q-\frac32}}
    \leq
    C_{q,K}\norm Z_{H^q}.
\end{equation}
Here the lower-order paradifferential terms take values in spaces of
higher regularity and are then embedded into \(H^{q-\frac32}\).

We also use the dependence of the concrete coefficients on the state.
The same analytic dependence of the CDE and of the changes of variables
used in Part~1 to obtain
\cref{eq:main_proof_analytic_exact_vector_field_3d,%
eq:main_proof_analytic_transformed_vector_field_3d},
together with the polarized versions of
\cref{simb-pro-geq-q-geq-q,eq:smoothing-tame-geq-q}, gives
\begin{equation}\label{eq:main_proof_unregularized_operator_differential_bound_3d}
    \norm{
        \dd\mathscr L_p\pare V\bra H\bra Z
    }_{H^{q-\frac32}}
    \leq
    C_{q,s,K}
    \norm H_{H^q}
    \norm Z_{H^s},
\end{equation}
for \(V\) in the same ball.  In particular, since
\[
    \mathscr L_p\pare{V_1}
    -
    \mathscr L_p\pare{V_2}
    =
    \int_0^1
    \dd\mathscr L_p
    \pare{
        V_2+\theta\pare{V_1-V_2}
    }
    \bra{V_1-V_2}
    \dd\theta,
\]
we have, thanks to the analyticity of the nonlinear vector field
\begin{equation}\label{eq:main_proof_unregularized_operator_difference_bound_3d}
    \norm{
        \pare{
            \mathscr L_p\pare{V_1}
            -
            \mathscr L_p\pare{V_2}
        }\bra Z
    }_{H^{q-\frac32}}
    \leq
    C_{q,s,K}
    \norm{V_1-V_2}_{H^q}
    \norm Z_{H^s}.
\end{equation}

We now apply these estimates to
\[
\begin{aligned}
    &\mathscr L_p\pare{W_N}\bra{W_{N+1}}
    -
    \mathscr L_p\pare W\bra W
    \\
    &\qquad=
    \mathscr L_p\pare{W_N}\bra{W_{N+1}-W_N}
   +
    \pare{
        \mathscr L_p\pare{W_N}
        -
        \mathscr L_p\pare W
    }\bra{W_N}
    +
    \mathscr L_p\pare W\bra{W_N-W}.
\end{aligned}
\]
By \cref{eq:main_proof_unregularized_operator_uniform_bound_3d},
\[
\begin{aligned}
    \norm{
        \mathscr L_p\pare{W_N}
        \bra{W_{N+1}-W_N}
    }_{q-\frac32,I_{\varepsilon}}
    &\leq
    C_{q,K}
    \norm{
        W_{N+1}-W_N
    }_{q,I_{\varepsilon}}
    \longrightarrow0,
    \\
    \norm{
        \mathscr L_p\pare W
        \bra{W_N-W}
    }_{q-\frac32,I_{\varepsilon}}
    &\leq
    C_{q,K}
    \norm{
        W_N-W
    }_{q,I_{\varepsilon}}
    \longrightarrow0.
\end{aligned}
\]
For the middle term,
\cref{eq:main_proof_unregularized_operator_difference_bound_3d} and the
uniform \(H^s\)-bound of Part~3 give
\[
\begin{aligned}
    &\norm{
        \pare{
            \mathscr L_p\pare{W_N}
            -
            \mathscr L_p\pare W
        }\bra{W_N}
    }_{q-\frac32,I_{\varepsilon}}
    \leq
    C_{q,s,K}
    \norm{W_N-W}_{q,I_{\varepsilon}}
    \norm{W_N}_{s,I_{\varepsilon}}
    \longrightarrow0
\end{aligned}
\]
by \cref{eq:main_proof_strong_lower_convergence_3d}.  Hence
\begin{equation}\label{eq:main_proof_unregularized_vector_field_convergence_3d}
    \norm{
        \mathscr L_p\pare{W_N}\bra{W_{N+1}}
        -
        \mathscr L_p\pare W\bra W
    }_{q-\frac32,I_{\varepsilon}}
    \longrightarrow0 .
\end{equation}
Since
\[
    \mathscr L_p\pare W\bra W
    =
    \mathfrak S_p\pare W,
\]
\cref{eq:main_proof_regularization_vanishes_3d,eq:main_proof_unregularized_vector_field_convergence_3d}
give
\begin{equation}\label{eq:main_proof_approximate_vector_fields_converge_3d}
    \mathscr L_{p,N+1}\pare{W_N}\bra{W_{N+1}}
    \longrightarrow
    \mathfrak S_p\pare W
    \quad\textnormal{in}\quad
    \cC\pare{I_{\varepsilon};H^{q-\frac32}} .
\end{equation}

Moreover,
\[
    P_{N+1}W_{\rm in}\longrightarrow W_{\rm in}
    \qquad\textnormal{in}\qquad H^q.
\]
Passing to the limit in
\cref{eq:main_proof_approximate_integral_equation_3d} therefore gives
\begin{equation}\label{eq:main_proof_exact_W_integral_equation_3d}
    W\pare t
    =
    W_{\rm in}
    +
    \int_0^t
    \mathfrak S_p\pare{W\pare\tau}
    \dd\tau .
\end{equation}
Thus
\[
    W\in
    \cC\pare{I_{\varepsilon};H^q}
    \cap
    \cC^1\pare{I_{\varepsilon};H^{q-\frac32}}
\]
and solves the exact transformed equation
\[
    \partial_tW=\mathfrak S_p\pare W .
\]

Finally set
\[
    Y\defeq\Psi_p\pare W .
\]
Since \(W_{\rm in}=\Phi_p\pare{Y_0}\), one has \(Y\pare0=Y_0\).
By the exact equivalence established in Part~1,
\[
    \begin{cases}
        \partial_tW=\mathfrak S_p\pare W\\
        W\pare{0}=W_{\rm in}
    \end{cases}
    \quad\Longleftrightarrow\quad
    \begin{cases}
        \partial_tY=\mathfrak F_p\pare Y\\
        Y\pare{0}=Y_0.
    \end{cases}
\]
Hence
\begin{equation}\label{eq:main_proof_exact_Y_integral_equation_3d}
    Y\pare t
    =
    Y_0
    +
    \int_0^t
    \mathfrak F_p\pare{Y\pare\tau}
    \dd\tau ,
\end{equation}
and the limit is a solution of the original canonical CDE.
\end{proofpart}

\begin{proofpart}[Uniqueness of the solution]
Let \(Y_1,Y_2\) be two solutions in the class constructed above with
initial data \(Y_{1,0},Y_{2,0}\), and set
\[
    W_j\defeq\Phi_p\pare{Y_j},
    \qquad
    Z\defeq W_1-W_2 .
\]
Since
\[
    \partial_tW_j
    =
    \mathscr L_p\pare{W_j}\bra{W_j},
\]
one has
\[
    \partial_tZ
    =
    \mathscr L_p\pare{W_1}\bra Z
    +
    \pare{
        \mathscr L_p\pare{W_1}
        -
        \mathscr L_p\pare{W_2}
    }\bra{W_2}.
\]
The \(L^2\)-comparison estimate of Part~4 applies directly, now with no
regularization error, and gives
\begin{equation}\label{eq:main_proof_uniqueness_L2_3d}
    \norm{Z\pare t}_{L^2}
    \leq
    C_{s,K}
    \norm{Z\pare0}_{L^2}
    +
    C_{s,K}
    \int_0^{\av t}
    \norm{Z\pare\tau}_{L^2}\dd\tau ,
    \qquad t\in I_{\varepsilon} .
\end{equation}
Hence, by Gronwall's inequality,
\[
    \norm{W_1-W_2}_{0,I_{\varepsilon}}
    \leq
    C_{s,K,I_{\varepsilon}}
    \norm{W_1\pare0-W_2\pare0}_{L^2}.
\]
In particular, if \(Y_{1,0}=Y_{2,0}\), then \(W_1=W_2\), and therefore
\(Y_1=Y_2\) since \(\Phi_p\) is invertible.  This proves uniqueness.
\end{proofpart}

\begin{proofpart}[Regularity of the limit]
By Part~3 and
\cref{eq:main_proof_intermediate_physical_convergence_3d},
for every \(t\in I_\varepsilon\) weak lower semicontinuity gives
\begin{equation}\label{eq:main_proof_limit_high_bound_3d}
    \norm{Y\pare t}_{\mathcal X^s}
    \leq
    \liminf_{N\to\infty}
    \norm{Y_N\pare t}_{\mathcal X^s}
    \leq
    C_{s,K}\pare{\varepsilon+\varepsilon\av t}.
\end{equation}
Hence
\[
    Y\in L^\infty\pare{I_\varepsilon;\mathcal X^s}.
\]
Moreover, by Part~4,
\[
    Y\in
    \cC\pare{I_\varepsilon;\mathcal X^q},
    \qquad
    s_0\leq q<s.
\]

It remains to prove strong continuity with values in
\(\mathcal X^s\).
Fix a lower Sobolev index \(q\) such that
\[
    s_0\leq q<s-\frac32,
    \qquad
    \delta\defeq s-q>\frac32.
\]
For \(N\geq1\), recall from \cref{eq:main_proof_regularized_iteration_3d_0,eq:main_proof_regularized_iteration_3d} that
\begin{equation}\label{eq:main_proof_smoothed_transformed_data_3d}
    W_{N}\pare{0}=
    P_NW_{\rm in},
    \qquad
    W_{\rm in}
    =
    \Phi_p\pare{Y_0}.
\end{equation}
Since \(W_{N}\pare{0}\) is smooth, the solution construction of Parts~2--4
may be performed at the index \(s+\frac32\).  Its associated lifespan is still
\(I_\varepsilon\) up to slightly modifying the value of the constant $c_{s,K}$, since
\[
    \norm{W_{N}\pare{0}}_{H^s}
    \leq
    C\norm{W_{\rm in}}_{H^s}
    \leq
    C_{s,K}\norm{Y_0}_{\mathcal X^s},
\]
by the uniform boundedness of \(P_N\) and the tame bounds for
\(\Phi_p\) following
\cref{eq:main_proof_nonlinear_diagonalizing_map_3d}.
Moreover, the Fourier localization gives
\[
    \norm{W_{N}\pare 0}_{H^{s+\frac32}}
    \leq
    C2^{\frac32N}
    \norm{W_{\rm in}}_{H^s}
    \leq
    C_{s,K}2^{\frac32N}
    \norm{Y_0}_{\mathcal X^s}.
\]
Hence the higher-order estimate of Part~3 gives
\[
    \sup_{t\in I_\varepsilon}
    \norm{W_{N}\pare t}_{H^{s+\frac32}}
    \leq
    C_{s,K,I_\varepsilon}
    2^{\frac32N}
    \norm{Y_0}_{\mathcal X^s}.
\]

Take $M,N\in\mathbb{N}^{\star}$  with $M>N$ and set
\[
    \mathcal{Z}_{M,N}\defeq W_M-W_N,
    \qquad
    \tau_{M,N}
    \defeq
    \norm{\pare{P_M-P_N}W_{{\rm in}}}_{H^s}.
\]
Then
\[
    \tau_{M,N}\longrightarrow0,\qquad\textnormal{as }N\to\infty.
\]
The frequency
localization of \(P_M-P_N\) gives
\begin{equation}\label{eq:main_proof_smoothed_initial_difference_3d}
\begin{aligned}
    \norm{\mathcal{Z}_{M,N}\pare0}_{H^s}
    &=\tau_{M,N},
    \\
    \norm{\mathcal{Z}_{M,N}\pare0}_{H^q}
    &\leq
    C_{s,K}
    2^{-\delta N}\tau_{M,N}.
\end{aligned}
\end{equation}

We first compare the two solutions at the lower index \(q\).
The difference \(\mathcal{Z}_{M,N}\) satisfies the same equation used in
Part~5.  Since
\[
    q+\frac32<s,
\]
the order-\(\frac32\) coefficient-difference term is controlled by
the uniform \(H^s\)-bounds of Part~3, and the corresponding energy
estimate gives
\begin{equation}\label{eq:main_proof_lower_bona_smith_comparison_3d}
    \norm{\mathcal{Z}_{M,N}}_{q,I_{\varepsilon}}
    \leq
    C_{s,K,I_\varepsilon}
    \norm{\mathcal{Z}_{M,N}\pare0}_{H^q}
    \leq
    C_{s,K,I_\varepsilon}
    2^{-\delta N}\tau_{M,N}.
\end{equation}

We now repeat the same comparison at the index \(s\), using the
symmetrized energy estimate of
\cref{lem:adapted_diagonal_weights_3d}.  The homogeneous part is
controlled exactly as in
\cref{eq:high_frequency_energy_bound_3d,%
eq:low_frequency_energy_bound_3d}.
For the coefficient-difference term, recall from
\cref{prop:positive_order_diagonalization_3d} that the only
positive-order contribution has order \(3/2\).  The tame symbol estimate
\cref{simb-pro-geq-q-geq-q}, the Bony--Weyl boundedness
\cref{thm:contS,cont00}, and the smoothing estimate
\cref{eq:smoothing-tame-geq-q}, applied in their polarized form as in
\cref{eq:main_proof_unregularized_operator_difference_bound_3d},
give
\begin{equation}\label{eq:main_proof_bona_smith_tame_operator_bound_3d}
\begin{aligned}
    &\norm{
        \pare{
            \mathscr L_{p,M}\pare{W_M}
            -
            \mathscr L_{p,N}\pare{W_N}
        }
        \bra{W_N}
    }_{H^s}\leq
    C_{s,K}
    \norm{\mathcal{Z}_{M,N}}_{H^s}
    +
    C_{s,K}
    \norm{\mathcal{Z}_{M,N}}_{H^{s-2}}
    \pare{\norm{W_M}_{H^{s+\frac32}}+\norm{W_N}_{H^{s+\frac32}}} .
\end{aligned}
\end{equation}
Here the first term is the high--low contribution and is absorbed in the
Gronwall estimate, while in the second one the high derivatives are
placed on the smoother solution.  This is the same tame splitting used
in \cref{eq:main_proof_frozen_energy_estimate_3d}. Recalling that
\[
    \sup_{t\in I_\varepsilon}
    \norm{W_N\pare t}_{H^{s+\frac32}}
    \leq
    C_{s,K,I_\varepsilon}
    2^{\frac32N}
    \norm{Y_0}_{\mathcal X^s},
\]
the coefficient-difference term is controlled by the lower-order
comparison multiplied by \(2^{\frac32N}\).
\[
\begin{aligned}
    \norm{\mathcal{Z}_{M,N}}_{s,I_\varepsilon
    }
    \leq
    C_{s,K,I_\varepsilon}
    \pare{
        \norm{\mathcal{Z}_{M,N}\pare0}_{H^s}
        +
        \norm{\mathcal{Z}_{M,N}}_{
            q,I_\varepsilon
        }
        \pare{\norm{W_M}_{s+\frac32,I_\varepsilon
        }+\norm{W_N}_{s+\frac32,I_\varepsilon
        }}
    }.
\end{aligned}
\]
Using
\cref{eq:main_proof_smoothed_initial_difference_3d,eq:main_proof_lower_bona_smith_comparison_3d},
we obtain
\begin{equation}\label{eq:main_proof_bona_smith_comparison_3d}
    \norm{\mathcal{Z}_{M,N}}_{
        s,I_\varepsilon
    }
    \leq
    C_{s,K,I_\varepsilon}
    \pare{
        1+
        2^{-\pare{\delta-\frac32}N}
    }
    \tau_{M,N}.
\end{equation}
Since \(\delta>\frac32\) and \(\tau_N\to0\),
\cref{eq:main_proof_bona_smith_comparison_3d} implies that the sequence $\pare{W_N}_N$ is Cauchy in $L^{\infty}\pare{I_\varepsilon;H^s}.$ By uniqueness of the limit, Part 4 gives
\[
    W_N
    \longrightarrow
    W
    \qquad\textnormal{in}\qquad
    L^{\infty}\pare{I_\varepsilon;H^s}.
\]
The local Lipschitz
continuity of \(\Psi_p:H^s\to\mathcal X^s\) yields
\[
    Y_N
    \longrightarrow
    Y
    \qquad\textnormal{in}\qquad
    L^{\infty}\pare{I_\varepsilon;\mathcal X^s}.
\]
Each \(Y_N\) belongs to
\(\cC\pare{I_\varepsilon;\mathcal X^s}\); hence, as a uniform limit of continuous functions,
\begin{equation}\label{eq:main_proof_strong_top_continuity_3d}
    Y\in
    \cC\pare{I_\varepsilon;\mathcal X^s}.
\end{equation}

Finally, by
\cref{eq:main_proof_analytic_exact_vector_field_3d},
\[
    \mathfrak F_p:
    \mathcal X^s
    \longrightarrow
    \mathcal X^{s-\frac32}
\]
is continuous.  Therefore
\[
    \mathfrak F_p\pare Y
    \in
    \cC\pare{
        I_\varepsilon;
        \mathcal X^{s-\frac32}
    },
\]
and the integral equation
\cref{eq:main_proof_exact_Y_integral_equation_3d} gives
\[
    Y\in
    \cC\pare{I_\varepsilon;\mathcal X^s}
    \cap
    \cC^1\pare{
        I_\varepsilon;
        \mathcal X^{s-\frac32}
    }.
\]
Together with Part~3 and Part~6, this completes the proof.\end{proofpart}

\end{proof}

\appendix

\section{Technical miscellanea}

\label{sec:technical lemma}

\begin{lemma}[Radial distribution identities]
\label{lem:abstract-radial-distribution-identities}
The following identities hold true in
\(\mathcal S'\pare{\bR}\). One has
\begin{equation}\label{eq:abstract-radial-distribution-identities}
\begin{aligned}
    \int_0^\infty
    \frac{\sin\pare{rt}}r\dd r
    &=
    \frac\pi2\sgn\pare t,
    &
    \int_0^\infty
    e^{-\ii rt}\dd r
    &=
    \pi\delta\pare t-\ii\pv\frac1t,
    &
    \int_0^\infty
    \frac{1-\cos\pare{rt}}{r^2}\dd r
    &=
    \frac\pi2\av t.
\end{aligned}
\end{equation}

\end{lemma}

\begin{proof}
We prove the three identities separately.
\begin{proofpart}
For \(\epsilon>0\), define
\[
    I_\epsilon\pare t
    \defeq
    \int_0^\infty
    e^{-\epsilon r}
    \frac{\sin\pare{rt}}r\dd r.
\]
The integral is absolutely convergent and may be differentiated with
respect to \(t\). Thus
\[
    \partial_t I_\epsilon\pare t
    =
    \int_0^\infty
    e^{-\epsilon r}\cos\pare{rt}\dd r
    =
    \frac{\epsilon}{\epsilon^2+t^2}\cdot
\]
Since \(I_\epsilon\pare0=0\), integration in \(t\) gives
\[
    I_\epsilon\pare t
    =
    \arctan\pare{t/\epsilon}.
\]
Hence
\[
    I_\epsilon\pare t
    \longrightarrow
    \frac\pi2\sgn\pare t
    \qquad
    \textnormal{for }t\neq0.
\]
Moreover,
\(
    \av{I_\epsilon\pare t}\leq\pi/2
\)
uniformly in \(\epsilon>0\) and \(t\in\bR\). Therefore, for every
\(\varphi\in\mathcal S\pare{\bR}\), dominated convergence yields
\[
    \lim_{\epsilon\downarrow0}
    \int_{\bR}
    I_\epsilon\pare t\varphi\pare t\dd t
    =
    \frac\pi2
    \int_{\bR}
    \sgn\pare t\varphi\pare t\dd t.
\]
This proves the first identity in
\cref{eq:abstract-radial-distribution-identities} in
\(\mathcal S'\pare{\bR}\).
\end{proofpart}
\begin{proofpart}
We next consider
\[
    J_\epsilon\pare t
    \defeq
    \int_0^\infty
    e^{-\epsilon r-\ii rt}\dd r
    =
    \frac1{\epsilon+\ii t}
    =
    \frac{\epsilon}{\epsilon^2+t^2}
    -
    \ii\frac{t}{\epsilon^2+t^2}\cdot
\]
We treat the real and imaginary parts separately. Let
\(\varphi\in\mathcal S\pare{\bR}\). For the real part, the change of
variables \(t=\epsilon s\) gives
\[
\begin{aligned}
    \int_{\bR}
    \frac{\epsilon}{\epsilon^2+t^2}
    \varphi\pare t\dd t
    &=
    \int_{\bR}
    \frac{\varphi\pare{\epsilon s}}{1+s^2}\dd s
    \\
    &\longrightarrow
    \varphi\pare0
    \int_{\bR}
    \frac1{1+s^2}\dd s
    =
    \pi\varphi\pare0.
\end{aligned}
\]
Thus
\[
    \frac{\epsilon}{\epsilon^2+t^2}
    \longrightarrow
    \pi\delta\pare t
    \qquad
    \textnormal{in }\mathcal S'\pare{\bR}.
\]

For the imaginary part, fix an even function
\(\chi\in\mathcal C^\infty_c\pare{\bR}\) such that
\(\chi\equiv1\) in a neighborhood of the origin. Since
\(t\pare{\epsilon^2+t^2}^{-1}\chi\pare t\) is odd,
\[
\begin{aligned}
    \int_{\bR}
    \frac{t}{\epsilon^2+t^2}
    \varphi\pare t\dd t
    &=
    \int_{\bR}
    \frac{t}{\epsilon^2+t^2}
    \bra{
        \varphi\pare t-\varphi\pare0\chi\pare t
    }\dd t
    \\
    &=
    \int_{\bR}
    \frac{t^2}{\epsilon^2+t^2}
    \frac{
        \varphi\pare t-\varphi\pare0\chi\pare t
    }t\dd t.
\end{aligned}
\]
The function
\[
    t
    \longmapsto
    \frac{
        \varphi\pare t-\varphi\pare0\chi\pare t
    }t
\]
extends smoothly through \(t=0\) and belongs to \(L^1\pare{\bR}\).
Since $0\leq
    \frac{t^2}{\epsilon^2+t^2}
    \leq1$, dominated convergence gives
\[
\begin{aligned}
    \lim_{\epsilon\downarrow0}
    \int_{\bR}
    \frac{t}{\epsilon^2+t^2}
    \varphi\pare t\dd t
    &=
    \int_{\bR}
    \frac{
        \varphi\pare t-\varphi\pare0\chi\pare t
    }t\dd t
    \\
    &=
    \pv\int_{\bR}\frac{\varphi\pare t}t\dd t,
\end{aligned}
\]
where the last equality uses the evenness of \(\chi\). Consequently,
\[
    \frac{t}{\epsilon^2+t^2}
    \longrightarrow
    \pv\frac1t
    \qquad
    \textnormal{in }\mathcal S'\pare{\bR}.
\]
Combining the real and imaginary parts proves the second identity in
\cref{eq:abstract-radial-distribution-identities}.
\end{proofpart}
\begin{proofpart}
Finally, for every \(t\in\bR\), the integral
\[
    K\pare t
    \defeq
    \int_0^\infty
    \frac{1-\cos\pare{rt}}{r^2}\dd r
\]
is absolutely convergent. Indeed,
\[
    1-\cos\pare{rt}
    =\begin{cases}
    O\pare{r^2t^2}
    \qquad\textnormal{as }r\to0,\\
    O\pare{1}
    \quad\,\,\,\,\,\qquad\textnormal{as }r\to\infty.\\
    \end{cases}
\]
For \(t\neq0\), the change of variables \(s=\av t r\) gives
\[
    K\pare t
    =
    \av t
    \int_0^\infty
    \frac{1-\cos s}{s^2}\dd s.
\]
It remains to evaluate the constant on the right-hand side. For
\(0<\delta<R\), integration by parts gives
\[
\begin{aligned}
    \int_\delta^R
    \frac{1-\cos r}{r^2}\dd r
    &=
    \bra{
        -\frac{1-\cos r}{r}
    }_\delta^R
    +
    \int_\delta^R
    \frac{\sin r}{r}\dd r.
\end{aligned}
\]
The boundary term tends to zero as
\(\delta\downarrow0\) and \(R\to\infty\). The remaining improper
integral converges by the Dirichlet criterion. Denoting
\[
    F\pare R
    \defeq
    \int_0^R\frac{\sin r}{r}\dd r,
\]
let
\(
    L\defeq\lim_{R\to\infty}F\pare R
\).
Since \(F\) is bounded, integration by parts gives, for every
\(\epsilon>0\),
\[
\begin{aligned}
    \int_0^\infty
    e^{-\epsilon r}
    \frac{\sin r}{r}\dd r
    &=
    \epsilon
    \int_0^\infty
    e^{-\epsilon r}F\pare r\dd r
    \\
    &=
    \int_0^\infty
    e^{-s}F\pare{s/\epsilon}\dd s.
\end{aligned}
\]
Dominated convergence therefore implies that the last expression tends
to \(L\) as \(\epsilon\downarrow0\). On the other hand, by the
calculation for \(I_\epsilon\),
\[
    \int_0^\infty
    e^{-\epsilon r}
    \frac{\sin r}{r}\dd r
    =
    \arctan\pare{1/\epsilon}
    \longrightarrow
    \frac\pi2\cdot
\]
Hence \(L=\pi/2\), and therefore
\[
    \int_0^\infty
    \frac{1-\cos r}{r^2}\dd r
    =
    \frac\pi2\cdot
\]
It follows that
\[
    K\pare t
    =
    \frac\pi2\av t,
\]
including \(t=0\), where both sides vanish. This proves the third
identity in \cref{eq:abstract-radial-distribution-identities} and concludes the proof.
\end{proofpart}
\end{proof}

\section{Abstract paralinearization for subcritical homogeneous graph kernels}
\label{app:abstract-locally-integrable-paralinearization-Rd}

Let
\(
    0<n<\mathtt d
\)
and let \(\pare{h,k}\) be an admissible homogeneous pair of order \(-n\)
and vanishing order \(\ell\) in the sense of
\cref{def:abstract-admissible-angular-pair-Rd}.  Throughout this appendix we
assume
\begin{equation}\label{eq:abstract-subcritical-orders-Rd}
    \beta\defeq\mathtt d-n>0,
    \qquad
    k\pare0=0,
    \qquad
    \ell>\beta.
\end{equation}
The condition \(k\pare0=0\) removes the flat homogeneous convolution, while
\(\ell>\beta\) makes the nonlinear
far-field contribution integrable. 

For every real \(\eta\) satisfying
\(
    \sup_{z\neq0}\norm{\Delta_z\eta}_{L^\infty}<r_k
\), define
\begin{equation}\label{eq:abstract-locally-integrable-graph-operator-Rd}
\begin{aligned}
    \mathpzc t_{h,k}^{\pare n}\pare\eta\bra\omega\pare x
    &\defeq
    \int_{\bR^{\mathtt d}}
    \frac{h\pare{\hat z}}{\av z^n}
    k\pare{\Delta_z\eta\pare x}
    \omega\pare{x-z}\dd z,
    \\
    \Delta_z\eta\pare x
    &\defeq
    \frac{\eta\pare x-\eta\pare{x-z}}{\av z},
    \qquad
    \hat z\defeq\frac z{\av z}\cdot
\end{aligned}
\end{equation}
The integral in \cref{eq:abstract-locally-integrable-graph-operator-Rd} is
ordinary both at the origin and at infinity: local integrability follows
from \(n<\mathtt d\), while
\cref{eq:abstract-subcritical-orders-Rd} and the vanishing order of \(k\)
give the required far-field decay. Our main goal is to prove the following theorem.

\begin{theorem}[Abstract subcritical paralinearization]
\label{thm:abstract-locally-integrable-paralinearization-Rd}
Let \(0<n<\mathtt d\), and let \(\pare{h,k}\) be an admissible homogeneous pair of order \(-n\) in the sense of
\cref{def:abstract-admissible-angular-pair-Rd} and satisfy
\cref{eq:abstract-subcritical-orders-Rd}.  There exists
\(\varrho_0>0\) such that, for every \(\varrho\geq\varrho_0\), there exist
\(s_0>\mathtt d/2+1\), \(\mu\geq0\), and \(r_0>0\) such that, for every
\(s\geq s_0\), every real
\[
    \eta\in H^{s+1+\mu}\pare{\bR^{\mathtt d};\bR},
    \qquad
    \omega\in H^{s+\mu}\pare{\bR^{\mathtt d};\bR},
    \qquad
    \norm\eta_{H^{s_0+1+\mu}}<r_0,
\]
the operator in
\cref{eq:abstract-locally-integrable-graph-operator-Rd} can be recast as
\begin{equation}\label{eq:abstract-locally-integrable-paralinearization-Rd}
    \mathpzc t_{h,k}^{\pare n}\pare\eta\bra\omega
    =
    \OpBW{
        A^{h,k}_{1,-\beta}\bra\omega
        \pare{\eta;x,\xi}
        +
        \vect{
            a_{\bra{0,\geq\ell}}\bra\omega\pare{\eta;x,\xi}
        }{
            a_{\bra{-\beta-2,\geq\ell}}\pare{\eta;x,\xi}
        }
    }
    \cdot
    \vect{\eta}{\omega}
    +
    R^{\lin,\bR}_{\bra{\varrho,\geq\ell}}\pare{\eta,\omega}
    \cdot
    \vect{\eta}{\omega}.
\end{equation}
 The objects in
\cref{eq:abstract-locally-integrable-paralinearization-Rd} are as follows.
\begin{itemize}

    \item The vector symbol
    \begin{equation}\label{eq:abstract-locally-integrable-vector-symbol-Rd}
        A^{h,k}_{1,-\beta}\bra\omega
        \pare{\eta;x,\xi}
        \defeq
        \vect{
            \omega\pare x
            a_1^{h,k;\eta}
            \pare{\nabla\eta\pare x,\xi}
        }{
            a_{-\beta}^{h,k;\omega}
            \pare{\nabla\eta\pare x,\xi}
        }
    \end{equation}
    has the following components.
    \begin{itemize}

        \item
        \(a_1^{h,k;\eta}\pare{\nabla\eta,\xi}
        \in\Gamma^1_{\geq\ell-1}\) is the symbol
        \begin{equation}\label{eq:abstract-locally-integrable-transported-symbol-Rd}
            a_1^{h,k;\eta}\pare{p,\xi}
            \defeq
            \int_{\bR^{\mathtt d}}
            \chi\pare z
            \frac{h\pare{\hat z}}{\av z^{n+1}}
            k'\pare{p\cdot\hat z}
            \pare{1-e^{-\ii z\cdot\xi}}\dd z;
        \end{equation}

        \item
        \(a_{-\beta}^{h,k;\omega}\pare{\nabla\eta,\xi}
        \in\Gamma^{-\beta}_{\geq\ell}\) is defined by
        \begin{equation}\label{eq:abstract-locally-integrable-quasilinear-symbol-Rd}
            a_{-\beta}^{h,k;\omega}\pare{p,\xi}
            \defeq
            \int_{\bR^{\mathtt d}}
            \chi\pare z
            \frac{h\pare{\hat z}}{\av z^n}
            k\pare{p\cdot\hat z}
            e^{-\ii z\cdot\xi}\dd z.
        \end{equation}

    \end{itemize}
    Both components satisfy
    \begin{equation}\label{eq:abstract-locally-integrable-principal-reality-Rd}
    \begin{aligned}
        a_1^{h,k;\eta}\pare{p,-\xi}
        &=
        \overline{a_1^{h,k;\eta}\pare{p,\xi}},
        &
        a_{-\beta}^{h,k;\omega}\pare{p,-\xi}
        &=
        \overline{a_{-\beta}^{h,k;\omega}\pare{p,\xi}}.
    \end{aligned}
    \end{equation}

    \item The generic vector symbol
    \[
        \vect{
            a_{\bra{0,\geq\ell}}\bra\omega\pare{\eta;x,\xi}
        }{
            a_{\bra{-\beta-2,\geq\ell}}\pare{\eta;x,\xi}
        }
    \]
    has the following components.
    \begin{itemize}

        \item
        \(
            a_{\bra{0,\geq\ell}}\bra\cdot
            \in
            \cL\pare{H^{s_0};\Gamma^0}
        \), and, for every real \(\omega\),
        \[
            a_{\bra{0,\geq\ell}}\bra\omega\pare{\eta;x,-\xi}
            =
            \overline{
                a_{\bra{0,\geq\ell}}\bra\omega\pare{\eta;x,\xi}
            }.
        \]
        If \(\ell\geq2\), the range of
        \(a_{\bra{0,\geq\ell}}\bra\cdot\) improves to
        \(\Gamma^0_{\geq\ell-1}\);

        \item
        \(a_{\bra{-\beta-2,\geq\ell}}
        \in\Gamma^{-\beta-2}_{\geq\ell}\), and
        \[
            a_{\bra{-\beta-2,\geq\ell}}\pare{\eta;x,-\xi}
            =
            \overline{
                a_{\bra{-\beta-2,\geq\ell}}\pare{\eta;x,\xi}
            }.
        \]

    \end{itemize}

    \item
    \(R^{\lin,\bR}_{\bra{\varrho,\geq\ell}}\pare{\eta,\omega}\) is a real-to-real
    linear smoothing remainder of the form
    \begin{equation}\label{eq:abstract-locally-integrable-linear-smoothing-remainder-Rd}
        R^{\lin,\bR}_{\bra{\varrho,\geq\ell}}\pare{\eta,\omega}
        \defeq
        \vect{
            r^\bR_{\bra{\varrho,\geq\ell}}\bra\omega\pare\eta
        }{
            r^\bR_{\bra{\varrho,\geq\ell}}\pare\eta
        },
    \end{equation}
    where
    \[
        r^\bR_{\bra{\varrho,\geq\ell}}\bra\cdot
        \in
        \cL\pare{H^{s_0};\cR^\varrho},
        \qquad
        r^\bR_{\bra{\varrho,\geq\ell}}
        \in
        \cR^\varrho_{\geq\ell}.
    \]
    Moreover, if \(\ell\geq2\), then
    \[
        r^\bR_{\bra{\varrho,\geq\ell}}\bra\cdot
        \in
        \cL\pare{
            H^{s_0};
            \cR^\varrho_{\geq\ell-1}
        }.
    \]

\end{itemize}
\end{theorem}

\subsection{Local--tail and Bony reduction}

For fixed \(z\neq0\), write
\(
    \omega_z\pare x\defeq\omega\pare{x-z}
\).
Since \(k\pare0=0\), one has \(\widetilde k=k\) in the notation of
\cref{eq:abstract-admissible-K-vanishing}.

\begin{lemma}[Local-tail and paraproduct reduction]
\label{lem:abstract-locally-integrable-Bony-reduction-Rd}
For every prescribed \(\varrho\geq0\), one has
\begin{equation}\label{eq:abstract-locally-integrable-Bony-decomposition-Rd}
\begin{aligned}
    \mathpzc t_{h,k}^{\pare n}\pare\eta\bra\omega
    ={}&
    \mathpzc t_{h,k}^{\QL}\pare\eta\bra\omega
    +
    \mathpzc t_{h,k}^{\mathrm{tr}}\pare\eta\bra\omega
    +
    \mathpzc t_{h,k}^{R}\pare\eta\bra\omega
    +
    \OpBW{
        f_{\geq\ell}^{\tail,\bR}
        \bra\omega\pare{\eta;x}
    }\eta
    +
    R^{\lin,\bR}_{\bra{\varrho,\geq\ell}}\pare{\eta,\omega}
    \cdot\vect{\eta}{\omega}.
\end{aligned}
\end{equation}
The objects in \cref{eq:abstract-locally-integrable-Bony-decomposition-Rd}
are as follows.
\begin{itemize}

    \item The quasilinear branch is
    \begin{equation}\label{eq:abstract-locally-integrable-quasilinear-branch-Rd}
        \mathpzc t_{h,k}^{\QL}\pare\eta\bra\omega
        \defeq
        \int_{\bR^{\mathtt d}}
        \chi\pare z
        \frac{h\pare{\hat z}}{\av z^n}
        \OpBW{k\pare{\Delta_z\eta}}\omega_z\dd z.
    \end{equation}

    \item The transported branch is
    \begin{equation}\label{eq:abstract-locally-integrable-transported-branch-Rd}
        \mathpzc t_{h,k}^{\mathrm{tr}}\pare\eta\bra\omega
        \defeq
        \int_{\bR^{\mathtt d}}
        \chi\pare z
        \frac{h\pare{\hat z}}{\av z^n}
        \OpBW{\omega_z}\bra{k\pare{\Delta_z\eta}}\dd z.
    \end{equation}

    \item The localized Bony remainder is
    \begin{equation}\label{eq:abstract-locally-integrable-localized-Bony-remainder-Rd}
        \mathpzc t_{h,k}^{R}\pare\eta\bra\omega
        \defeq
        \int_{\bR^{\mathtt d}}
        \chi\pare z
        \frac{h\pare{\hat z}}{\av z^n}
        \RBW{k\pare{\Delta_z\eta}}{\omega_z}\dd z.
    \end{equation}

    \item The tail function is
    \begin{equation}\label{eq:abstract-locally-integrable-tail-function-Rd}
        f_{\geq\ell}^{\tail,\bR}
        \bra\omega\pare{\eta;x}
        \defeq
        \int_{\bR^{\mathtt d}}
        \pare{1-\chi\pare z}
        \frac{h\pare{\hat z}}{\av z^{n+1}}
        k'\pare{\Delta_z\eta\pare x}
        \omega_z\pare x\dd z.
    \end{equation}

    \item The smoothing vector
    \(R^{\lin,\bR}_{\bra{\varrho,\geq\ell}}\pare{\eta,\omega}\)
    is as in
    \cref{eq:abstract-locally-integrable-linear-smoothing-remainder-Rd}.

\end{itemize}
\end{lemma}

\begin{proof}
The localized contribution is decomposed as in the proof of
\cref{lem:abstract-common-local-tail-Bony-reduction-Rd} by applying
\cref{lem:paraprod_Bony} to
\(k\pare{\Delta_z\eta}\omega_z\).

For the tail contribution, by
\cref{eq:abstract-subcritical-orders-Rd},
\[
    n+\ell
    =
    \mathtt d-\beta+\ell
    >
    \mathtt d.
\]
Moreover,
\[
    b_h\pare z
    \defeq
    \frac{h\pare{\hat z}}{\av z^n}
\]
satisfies the far-field bounds required in the definition of
\(\KF_{\textnormal{tail},\geq\ell}\). Hence
\[
    \pare{1-\chi\pare z}
    \frac{h\pare{\hat z}}{\av z^n}
    k\pare{\Delta_z\eta}
    \in
    \KF_{\textnormal{tail},\geq\ell}.
\]
Applying \cref{prop:tail-kernel-smoothing} gives
\[
\begin{aligned}
    &\int_{\bR^{\mathtt d}}
    \pare{1-\chi\pare z}
    \frac{h\pare{\hat z}}{\av z^n}
    k\pare{\Delta_z\eta}
    \omega_z\dd z
    \\
    &\qquad=
    \OpBW{
        f_{\geq\ell}^{\tail,\bR}
        \bra\omega\pare{\eta;x}
    }\eta
    +
    R^{\lin,\bR}_{\bra{\varrho,\geq\ell}}\pare{\eta,\omega}
    \cdot\vect{\eta}{\omega},
\end{aligned}
\]
after absorbing the two smoothing remainders in the generic smoothing
vector. Combining the local and tail contributions proves
\cref{eq:abstract-locally-integrable-Bony-decomposition-Rd}.
\end{proof}

\begin{lemma}[Localized Bony remainder]
\label{lem:abstract-locally-integrable-localized-Bony-remainder-Rd}
For every prescribed \(\varrho\geq0\), the localized Bony remainder
satisfies
\begin{equation}\label{eq:abstract-locally-integrable-Bony-remainder-final-Rd}
    \mathpzc t_{h,k}^{R}\pare\eta\bra\omega
    =
    r^\bR_{\bra{\varrho,\geq\ell}}\pare\eta\omega,
    \qquad
    r^\bR_{\bra{\varrho,\geq\ell}}
    \in
    \cR^\varrho_{\geq\ell}.
\end{equation}
In particular, the remainder is real-to-real.
\end{lemma}

\begin{proof}
Since \(n<\mathtt d\), the localized homogeneous factor
\(
    \chi\pare z h\pare{\hat z}\av z^{-n}
\)
has local order \(-n>-\mathtt d\).  By
\cref{rem:smoothBONYremainder}, after increasing the baseline index
according to the prescribed \(\varrho\), the family
\[
    \RBW{k\pare{\Delta_z\eta}}{\omega_z}
\]
is a \(z\)-dependent smoothing operator with background vanishing order
\(\ell\).  Multiplication by the localized homogeneous factor preserves a
locally integrable kernel-smoothing class.  Hence
\cref{lem:integration_KR} gives
\cref{eq:abstract-locally-integrable-Bony-remainder-final-Rd}.  The
real-to-real property follows from that of the Bony remainder and from the
reality of the kernel.
\end{proof}

\subsection{Quasilinear contribution}

\begin{lemma}[Quasilinear contribution]
\label{lem:abstract-locally-integrable-quasilinear-contribution-Rd}
For every prescribed \(\varrho\geq0\), the quasilinear branch in
\cref{eq:abstract-locally-integrable-quasilinear-branch-Rd} satisfies
\begin{equation}\label{eq:abstract-locally-integrable-quasilinear-final-Rd}
    \mathpzc t_{h,k}^{\QL}\pare\eta\bra\omega
    =
    \OpBW{
        a_{-\beta}^{h,k;\omega}\pare{\nabla\eta,\xi}
        +
        a_{\bra{-\beta-2,\geq\ell}}\pare{\eta;x,\xi}
    }\omega
    +
    r^\bR_{\bra{\varrho,\geq\ell}}\pare\eta\omega.
\end{equation}
The objects in this equation have the following properties.
\begin{itemize}

    \item
    \(a_{-\beta}^{h,k;\omega}\) is defined in
    \cref{eq:abstract-locally-integrable-quasilinear-symbol-Rd} and belongs
    to \(\Gamma^{-\beta}_{\geq\ell}\);

    \item
    \(a_{\bra{-\beta-2,\geq\ell}}\in\Gamma^{-\beta-2}_{\geq\ell}\) and
    \(r^\bR_{\bra{\varrho,\geq\ell}}\in\cR^\varrho_{\geq\ell}\) are generic terms in
    the sense of \cref{notation:paralinearization}.  Both are
    real-to-real.

\end{itemize}
\end{lemma}

\begin{proof}
\begin{proofpart}[Proof of \cref{eq:abstract-locally-integrable-local-Taylor-splitting-Rd}]
We first decompose, analogously as it was done in the critical case in \cref{sec:abstract-structure-preserving-paralinearization-Rd}
\begin{equation}\label{eq:abstract-locally-integrable-local-Taylor-splitting-Rd}
    \mathpzc t_{h,k}^{\QL}
    =
    \mathpzc t_{h,k}^{\QL,-\beta}
    +
    \mathpzc t_{h,k}^{\QL,-\beta-1}
    +
    \mathpzc t_{h,k}^{\QL,\leq-\beta-2},
\end{equation}
where
\begin{subequations}\label{eq:abstract-locally-integrable-local-Taylor-models-Rd}
\begin{align}
    \mathpzc t_{h,k}^{\QL,-\beta}\pare\eta\bra\omega
    &\defeq
    \int_{\bR^{\mathtt d}}
    \chi\pare z
    \frac{h\pare{\hat z}}{\av z^n}
    \OpBW{k\pare{\mathsf{v}\bra{\eta}}}\omega_z\dd z,
    \label{eq:abstract-locally-integrable-leading-model-Rd}
    \\
    \mathpzc t_{h,k}^{\QL,-\beta-1}\pare\eta\bra\omega
    &\defeq
    -\frac12
    \int_{\bR^{\mathtt d}}
    \chi\pare z
    \frac{h\pare{\hat z}}{\av z^{n-1}}
    \OpBW{k'\pare{\mathsf{v}\bra{\eta}}
        \mathsf{h}\bra{\eta}}\omega_z\dd z,
    \label{eq:abstract-locally-integrable-subprincipal-model-Rd}
    \\
    \mathpzc t_{h,k}^{\QL,\leq-\beta-2}\pare\eta\bra\omega
    &\defeq
    \int_{\bR^{\mathtt d}}
    \chi\pare z
    \frac{h\pare{\hat z}}{\av z^n}
    \OpBW{k^{\loc,\bR}_{\bra{2,\geq\ell}}}\omega_z\dd z.
    \label{eq:abstract-locally-integrable-lower-model-Rd}
\end{align}
\end{subequations}
\end{proofpart}

\begin{proofpart}[Analysis of \cref{eq:abstract-locally-integrable-leading-model-Rd}]
We prove
\begin{equation}\label{eq:abstract-locally-integrable-leading-branch-Rd}
\begin{aligned}
    \mathpzc t_{h,k}^{\QL,-\beta}\pare\eta\bra\omega
    ={}&
    \OpBW{
        a_{-\beta}^{h,k;\omega}\pare{\nabla\eta,\xi}
        +
        \frac\ii2
        \pare{\nabla_x\cdot\nabla_\xi}
        a_{-\beta}^{h,k;\omega}\pare{\nabla\eta,\xi}
        +
        a_{\bra{-\beta-2,\geq\ell}}
    }\omega
    +
    r^\bR_{\bra{\varrho,\geq\ell}}\pare\eta\omega,
\end{aligned}
\end{equation}
where
\[
    a_{\bra{-\beta-2,\geq\ell}}
    \in\Gamma^{-\beta-2}_{\geq\ell},
    \qquad
    r^\bR_{\bra{\varrho,\geq\ell}}
    \in\cR^\varrho_{\geq\ell}.
\]

The point at which the subcritical assumption differs from the critical
case is that \(n<\mathtt d\). Hence the kernel in
\cref{eq:abstract-locally-integrable-leading-model-Rd} is locally
integrable at \(z=0\), and \cref{prop:int-opbw} applies directly to the
leading branch. In particular, by
\cref{eq:int-opbw-centered-symbol},
\begin{equation}\label{eq:abstract-locally-integrable-leading-centered-symbol-Rd}
    \mathpzc t_{h,k}^{\QL,-\beta}\pare\eta\bra\omega
    =
    \OpBW{
        b_{-\beta}^{h,k}\pare{\eta;x,\xi}
    }\omega,
\end{equation}
where
\begin{equation}\label{eq:abstract-locally-integrable-leading-centered-symbol-definition-Rd}
\begin{aligned}
    b_{-\beta}^{h,k}\pare{\eta;x,\xi}
    \defeq
    \int_{\bR^{\mathtt d}}
    \chi\pare z
    \frac{h\pare{\hat z}}{\av z^n}
    k\pare{
        \mathsf{v}\bra{\eta}
        \pare{x+\frac z2,z}
    }
    e^{-\ii z\cdot\xi}\dd z.
\end{aligned}
\end{equation}
For fixed \(z\neq0\), set
\[
    F_z\pare x
    \defeq
    k\pare{\mathsf{v}\bra{\eta}\pare{x,z}}.
\]
Taylor's formula at \(x\) gives
\begin{equation}\label{eq:abstract-locally-integrable-centered-Taylor-Rd}
\begin{aligned}
    F_z\pare{x+\frac z2}
    ={}&
    F_z\pare x
    +
    \frac12 z\cdot\nabla_xF_z\pare x
    +
    \frac14
    \int_0^1
    \pare{1-t}
    \pare{z\cdot\nabla_x}^2
    F_z\pare{x+\frac t2z}
    \dd t.
\end{aligned}
\end{equation}
Substituting the zeroth Taylor term into
\cref{eq:abstract-locally-integrable-leading-centered-symbol-definition-Rd}
gives exactly
\[
    a_{-\beta}^{h,k;\omega}
    \pare{\nabla\eta\pare x,\xi}
\]
by
\cref{eq:abstract-locally-integrable-quasilinear-symbol-Rd}.
For the first Taylor term, differentiating under the integral sign gives
\begin{align}
    &\frac12
    \int_{\bR^{\mathtt d}}
    \chi\pare z
    \frac{h\pare{\hat z}}{\av z^n}
    z\cdot\nabla_x
    \bra{
        k\pare{\mathsf{v}\bra{\eta}\pare{x,z}}
    }
    e^{-\ii z\cdot\xi}\dd z
    =
    \frac\ii2
    \pare{\nabla_x\cdot\nabla_\xi}
    a_{-\beta}^{h,k;\omega}
    \pare{\nabla\eta\pare x,\xi}.
\end{align}

It remains to consider the last term in
\cref{eq:abstract-locally-integrable-centered-Taylor-Rd}.
The two explicit factors of \(z\) improve the local order of the kernel
from \(-n\) to \(-n+2\), while differentiation in \(x\) preserves its
background vanishing order \(\ell\). Hence the corresponding kernel
belongs to $ \KF^{-n+2}_{\geq\ell}$. 
Since \(-n+2>-\mathtt d\), another application of
\cref{prop:int-opbw} gives a symbol in $\Gamma^{-\pare{\pare{-n+2}+\mathtt d}}_{\geq\ell}
    =
    \Gamma^{-\beta-2}_{\geq\ell}$. 
This yields the term
\(a_{\bra{-\beta-2,\geq\ell}}\) in
\cref{eq:abstract-locally-integrable-leading-branch-Rd} and proves the
claim.
\end{proofpart}
\begin{proofpart}[Analysis  of \cref{eq:abstract-locally-integrable-subprincipal-model-Rd}]
We prove
\begin{equation}\label{eq:abstract-locally-integrable-geometric-branch-Rd}
\begin{aligned}
    \mathpzc t_{h,k}^{\QL,-\beta-1}\pare\eta\bra\omega
    ={}&
    \OpBW{
        -\frac\ii2
        \pare{\nabla_x\cdot\nabla_\xi}
        a_{-\beta}^{h,k;\omega}\pare{\nabla\eta,\xi}
        +
        a_{\bra{-\beta-2,\geq\ell}}
    }\omega
    +
    r^\bR_{\bra{\varrho,\geq\ell}}\pare\eta\omega,
\end{aligned}
\end{equation}
with the classes in
\cref{eq:abstract-locally-integrable-leading-branch-Rd}.
Differentiating
\cref{eq:abstract-locally-integrable-quasilinear-symbol-Rd} under the
integral sign gives
\begin{equation}\label{eq:abstract-locally-integrable-principal-derivative-Rd}
\begin{aligned}
    &\pare{\nabla_x\cdot\nabla_\xi}
    a_{-\beta}^{h,k;\omega}
    \pare{\nabla\eta\pare x,\xi}
    =
    -\ii
    \int_{\bR^{\mathtt d}}
    \chi\pare z
    \frac{h\pare{\hat z}}{\av z^{n-1}}
    k'\pare{\mathsf{v}\bra{\eta}\pare{x,z}}
    \mathsf{h}\bra{\eta}\pare{x,z}
    e^{-\ii z\cdot\xi}\dd z.
\end{aligned}
\end{equation}
Thus the uncentered Fourier integral of the kernel in
\cref{eq:abstract-locally-integrable-subprincipal-model-Rd} is the
explicit term in
\cref{eq:abstract-locally-integrable-geometric-branch-Rd}.
The kernel has local order \(-n+1\), so every nonconstant term in its
expansion at \(x+z/2\) has local order at least \(-n+2\).
Applying \cref{prop:int-opbw} proves
\cref{eq:abstract-locally-integrable-geometric-branch-Rd}.
Consequently, the two order-\(\pare{-\beta-1}\) terms in
\cref{eq:abstract-locally-integrable-leading-branch-Rd,eq:abstract-locally-integrable-geometric-branch-Rd}
cancel:
\begin{equation}\label{eq:abstract-locally-integrable-subprincipal-cancellation-Rd}
    \frac\ii2
    \pare{\nabla_x\cdot\nabla_\xi}
    a_{-\beta}^{h,k;\omega}
    -
    \frac\ii2
    \pare{\nabla_x\cdot\nabla_\xi}
    a_{-\beta}^{h,k;\omega}
    =0.
\end{equation}
\end{proofpart}

\begin{proofpart}[Proof of \cref{eq:abstract-locally-integrable-lower-remainder-symbol-Rd}]
We prove
\begin{equation}\label{eq:abstract-locally-integrable-lower-remainder-symbol-Rd}
    \mathpzc t_{h,k}^{\QL,\leq-\beta-2}\pare\eta\bra\omega
    =
    \OpBW{
        a_{\bra{-\beta-2,\geq\ell}}\pare{\eta;x,\xi}
    }\omega
    +
    r^\bR_{\bra{\varrho,\geq\ell}}\pare\eta\omega,
\end{equation}
where
\[
    a_{\bra{-\beta-2,\geq\ell}}
    \in
    \Gamma^{-\beta-2}_{\geq\ell},
    \qquad
    r^\bR_{\bra{\varrho,\geq\ell}}
    \in
    \cR^\varrho_{\geq\ell}.
\]
Since \(-n+2>-\mathtt d\), this follows from
\cref{prop:int-opbw}.
The real-to-real conclusions follow from
\cref{eq:int-opbw-reality} and the reality of the kernel.
\end{proofpart}

\begin{proofpart}[Proof of \cref{eq:abstract-locally-integrable-quasilinear-final-Rd}]
Combining
\cref{eq:abstract-locally-integrable-local-Taylor-splitting-Rd,%
eq:abstract-locally-integrable-leading-branch-Rd,%
eq:abstract-locally-integrable-geometric-branch-Rd,%
eq:abstract-locally-integrable-lower-remainder-symbol-Rd}
and using
\cref{eq:abstract-locally-integrable-subprincipal-cancellation-Rd}
proves \cref{eq:abstract-locally-integrable-quasilinear-final-Rd}.
The symbol estimates in \cref{prop:int-opbw}, applied to the Fourier
integral in
\cref{eq:abstract-locally-integrable-quasilinear-symbol-Rd}, give $$  a_{-\beta}^{h,k;\omega}\pare{\nabla\eta,\xi}
    \in\Gamma^{-\beta}_{\geq\ell}.$$ All kernels are real-valued; hence
\cref{eq:int-opbw-reality} gives the real-to-real symbols and remainders
in \cref{eq:abstract-locally-integrable-quasilinear-final-Rd}.
\end{proofpart}
\end{proof}

\subsection{Transported contribution}

\begin{lemma}[Transported contribution]
\label{lem:abstract-locally-integrable-transported-contribution-Rd}
For every prescribed \(\varrho\geq0\), the transported branch in
\cref{eq:abstract-locally-integrable-transported-branch-Rd} satisfies
\begin{equation}\label{eq:abstract-locally-integrable-transported-final-Rd}
    \mathpzc t_{h,k}^{\mathrm{tr}}\pare\eta\bra\omega
    =
    \OpBW{
        \omega\pare x
        a_1^{h,k;\eta}\pare{\nabla\eta\pare x,\xi}
        +
        a_{\bra{0,\geq\ell}}\bra\omega\pare{\eta;x,\xi}
    }\eta
    +
    r^\bR_{\bra{\varrho,\geq\ell}}\bra\omega\pare\eta\eta.
\end{equation}
The objects in this equation have the following properties.
\begin{itemize}

    \item
    \(a_1^{h,k;\eta}\) is defined in
    \cref{eq:abstract-locally-integrable-transported-symbol-Rd} and belongs
    to \(\Gamma^1_{\geq\ell-1}\);

    \item
    \(
        a_{\bra{0,\geq\ell}}\bra\cdot
        \in
        \cL\pare{H^{s_0};\Gamma^0}
    \)
and
    \(
        r^\bR_{\bra{\varrho,\geq\ell}}\bra\cdot
        \in
        \cL\pare{H^{s_0};\cR^\varrho}
    \).
    If \(\ell\geq2\), both ranges improve to vanishing order
    \(\ell-1\).  The symbol and the remainder are real-to-real.

\end{itemize}
\end{lemma}

\begin{proof}
\begin{proofpart}[Proof of \cref{eq:abstract-locally-integrable-transported-Bony-splitting-Rd}]
Applying \cref{lem paralin} to \(k\pare{\Delta_z\eta}\) gives
\begin{equation}\label{eq:abstract-locally-integrable-transported-Bony-splitting-Rd}
    k\pare{\Delta_z\eta}
    =
    \OpBW{k'\pare{\Delta_z\eta}}\Delta_z\eta
    +
    \mathpzc q_k\pare{\Delta_z\eta}\bra{\Delta_z\eta},
\end{equation}
where \(\mathpzc q_k\defeq R_k\). Substitution into
\cref{eq:abstract-locally-integrable-transported-branch-Rd} gives
\[
    \mathpzc t_{h,k}^{\mathrm{tr}}
    =
    \mathpzc t_{h,k}^{\mathrm{tr},\mathrm I}
    +
    \mathpzc t_{h,k}^{\mathrm{tr},\mathrm {II}},
\]
with
\begin{subequations}\label{eq:abstract-locally-integrable-transported-two-branches-Rd}
\begin{align}
    \mathpzc t_{h,k}^{\mathrm{tr},\mathrm I}\pare\eta\bra\omega
    &\defeq
    \int_{\bR^{\mathtt d}}
    \chi\pare z
    \frac{h\pare{\hat z}}{\av z^n}
    \OpBW{\omega_z}
    \OpBW{k'\pare{\Delta_z\eta}}
    \Delta_z\eta\dd z,
    \\
    \mathpzc t_{h,k}^{\mathrm{tr},\mathrm {II}}\pare\eta\bra\omega
    &\defeq
    \int_{\bR^{\mathtt d}}
    \chi\pare z
    \frac{h\pare{\hat z}}{\av z^n}
    \OpBW{\omega_z}
    \mathpzc q_k\pare{\Delta_z\eta}\bra{\Delta_z\eta}\dd z.
\end{align}
\end{subequations}
\end{proofpart}

\begin{proofpart}[Proof of \cref{eq:abstract-locally-integrable-transported-principal-almost-Rd}]
We prove
\begin{equation}\label{eq:abstract-locally-integrable-transported-principal-almost-Rd}
    \mathpzc t_{h,k}^{\mathrm{tr},\mathrm I}\pare\eta\bra\omega
    =
    \mathpzc t_{h,k}^{\mathrm{tr},1,\loc}\pare\eta\bra\omega
    +
    \OpBW{a_{\bra{0,\geq\ell}}\bra\omega}\eta
    +
    r^\bR_{\bra{\varrho,\geq\ell}}\bra\omega\pare\eta\eta,
\end{equation}
where the first term is
\begin{equation}\label{eq:abstract-locally-integrable-transported-local-principal-model-Rd}
    \mathpzc t_{h,k}^{\mathrm{tr},1,\loc}\pare\eta\bra\omega
    \defeq
    \int_{\bR^{\mathtt d}}
    \chi\pare z
    \frac{h\pare{\hat z}}{\av z^{n+1}}
    \OpBW{\omega\pare x k'\pare{\mathsf{v}\bra{\eta}}}
    \delta_z\eta\dd z,
\end{equation}
and the generic terms have the classes stated in
\cref{lem:abstract-locally-integrable-transported-contribution-Rd}.
Recall \(\mathsf{v}\bra{\eta}\) from
\cref{eq:abstract-quasilinear-local-building-blocks}.
The identities
\[
    \omega_z-\omega\pare x
    =
    -\av z
    \int_0^1
    \nabla\omega\pare{x-tz}\cdot\hat z\dd t,\quad\,
    \,
    k'\pare{\Delta_z\eta}-k'\pare{\mathsf{v}\bra{\eta}}
    =
    \pare{\Delta_z\eta-\mathsf{v}\bra{\eta}}
    \int_0^1
    k''\pare{
        \mathsf{v}\bra{\eta}
        +t\pare{\Delta_z\eta-\mathsf{v}\bra{\eta}}
    }\dd t,
\]
together with
\cref{eq:abstract-quasilinear-difference-Taylor-expansion}, give
\begin{equation}\label{eq:abstract-locally-integrable-transported-freezing-Rd}
    \omega_z k'\pare{\Delta_z\eta}
    =
    \omega\pare x k'\pare{\mathsf{v}\bra{\eta}}
    +
    e_{\bra1}^{\mathrm{tr}}\bra\omega,
    \qquad
    e_{\bra1}^{\mathrm{tr}}\bra\cdot
    \in
    \cL\pare{H^{s_0};\KF^1_{\loc}}.
\end{equation}
If \(\ell\geq2\), the range improves to
\(\KF^1_{\loc,\geq\ell-1}\).
Using the composition defect
\cref{eq:abstract-transported-composition-defect} and
\(\Delta_z\eta=\delta_z\eta/\av z\), we obtain the exact splitting
\begin{equation}\label{eq:abstract-locally-integrable-transported-principal-splitting-Rd}
\begin{aligned}
    \mathpzc t_{h,k}^{\mathrm{tr},\mathrm I}\pare\eta\bra\omega
    ={}&
    \mathpzc t_{h,k}^{\mathrm{tr},1,\loc}\pare\eta\bra\omega
    \\
    &+
    \int_{\bR^{\mathtt d}}
    \chi\pare z
    \frac{h\pare{\hat z}}{\av z^{n+1}}
    \OpBW{e_{\bra1}^{\mathrm{tr}}\bra\omega}
    \delta_z\eta\dd z
    \\
    &+
    \int_{\bR^{\mathtt d}}
    \chi\pare z
    \frac{h\pare{\hat z}}{\av z^n}
    \mathpzc c_z\bra\omega\pare\eta\bra{\Delta_z\eta}\dd z.
\end{aligned}
\end{equation}
We analyze the three terms in
\cref{eq:abstract-locally-integrable-transported-principal-splitting-Rd}
in their order of appearance.

\begin{proofsubpart}[Proof of \cref{eq:abstract-locally-integrable-transported-principal-symbol-reduction-Rd}]
We first prove
\begin{equation}\label{eq:abstract-locally-integrable-transported-principal-symbol-reduction-Rd}
\begin{aligned}
    \mathpzc t_{h,k}^{\mathrm{tr},1,\loc}\pare\eta\bra\omega
    ={}&
    \OpBW{
        \omega\pare x
        a_1^{h,k;\eta}\pare{\nabla\eta\pare x,\xi}
    }\eta
    \\
    &+
    \OpBW{a_{\bra{0,\geq\ell}}\bra\omega}\eta
    +
    r^\bR_{\bra{\varrho,\geq\ell}}\bra\omega\pare\eta\eta.
\end{aligned}
\end{equation}
The first difference in
\cref{eq:abstract-locally-integrable-transported-local-principal-model-Rd}
cannot in general be split as \(\delta_z\eta=\eta-\eta_z\):
the two terms would be separately locally integrable only for
\(\beta>1\). We therefore retain the first difference.
In the Bony--Weyl Fourier representation
\cref{BonyWeyl,eq:int-opbw-centered-symbol-hat}, let \(\varsigma\) be
the Fourier variable of the coefficient and \(\xi\) the midpoint
frequency. We use
\begin{equation}\label{eq:abstract-locally-integrable-transported-compensated-splitting-Rd}
\begin{aligned}
    1-e^{-\ii\pare{\xi-\varsigma/2}\cdot z}
    ={}&
    1-e^{-\ii\xi\cdot z}
    \\
    &+
    e^{-\ii\xi\cdot z}
    \pare{1-e^{\ii\varsigma\cdot z/2}}.
\end{aligned}
\end{equation}
By \cref{eq:abstract-locally-integrable-transported-symbol-Rd}, the
first term in
\cref{eq:abstract-locally-integrable-transported-compensated-splitting-Rd}
produces exactly the principal symbol in
\cref{eq:abstract-locally-integrable-transported-principal-symbol-reduction-Rd}.
For the second term,
\[
    1-e^{\ii\varsigma\cdot z/2}
    =
    -\frac\ii2\pare{\varsigma\cdot z}
    \int_0^1 e^{\ii t\varsigma\cdot z/2}\dd t,
    \qquad
    \av{1-e^{\ii\varsigma\cdot z/2}}
    \lesssim\av\varsigma\,\av z.
\]
The factor \(\varsigma\) differentiates the coefficient, while
\(\av z\) lowers the local singularity from \(-\pare{n+1}\) to
\(-n>-\mathtt d\). Thus \cref{prop:int-opbw} gives an order-zero
symbol and a smoothing remainder, as in
\cref{eq:abstract-locally-integrable-transported-principal-symbol-reduction-Rd}.

To record the order of \(a_1^{h,k;\eta}\), decompose the localized kernel
in \cref{eq:abstract-locally-integrable-transported-symbol-Rd}
dyadically. On an annulus \(\av z\sim2^j\), its \(L^1_z\)-size is
\(O\pare{2^{j\pare{\beta-1}}}\), while
\(1-e^{-\ii z\cdot\xi}\) gives
\(\min\set{2^j\av\xi,1}\). The resulting bounds, and their
\(p\)- and \(\xi\)-derivatives, imply
\begin{equation}\label{eq:abstract-locally-integrable-transported-symbol-bound-Rd}
    \av{
        \partial_p^\alpha\partial_\xi^\gamma
        a_1^{h,k;\eta}\pare{p,\xi}
    }
    \lesssim_{\alpha,\gamma}
    \angles\xi^{1-\av\gamma},
    \qquad
    \av p\leq r_*.
\end{equation}
Together with
\cref{eq:abstract-admissible-K-vanishing-order}, this proves
\(a_1^{h,k;\eta}\pare{\nabla\eta,\xi}\in\Gamma^1_{\geq\ell-1}\).
\end{proofsubpart}

\begin{proofsubpart}[Proof of \cref{eq:abstract-locally-integrable-transported-freezing-error-Rd}]
For the second term in
\cref{eq:abstract-locally-integrable-transported-principal-splitting-Rd},
we prove
\begin{equation}\label{eq:abstract-locally-integrable-transported-freezing-error-Rd}
\begin{aligned}
    \int_{\bR^{\mathtt d}}
    \chi\pare z
    \frac{h\pare{\hat z}}{\av z^{n+1}}
    \OpBW{e_{\bra1}^{\mathrm{tr}}\bra\omega}
    \delta_z\eta\dd z
    &=
    \OpBW{
        a_{\bra{0,\geq\ell}}\bra\omega
        \pare{\eta;x,\xi}
    }\eta\\
    &+
    r^\bR_{\bra{\varrho,\geq\ell}}\bra\omega\pare\eta\eta.
\end{aligned}
\end{equation}
By \cref{eq:abstract-locally-integrable-transported-freezing-Rd}, the
kernel after multiplication by \(\chi h\av z^{-n-1}\) has local order
\(-n>-\mathtt d\). We may therefore split
\(\delta_z\eta=\eta-\eta_z\).
Then \cref{lem:integration_KF} applies to the unshifted term and
\cref{prop:int-opbw} to the translated one, with \(\omega\) as a
separately estimated linear parameter and \(\eta\) as the argument.
The latter symbol has order \(-\beta<0\) and is included in the generic
order-zero term. This proves
\cref{eq:abstract-locally-integrable-transported-freezing-error-Rd},
with the classes stated in
\cref{lem:abstract-locally-integrable-transported-contribution-Rd}.
\end{proofsubpart}

\begin{proofsubpart}[Proof of \cref{eq:abstract-locally-integrable-transported-principal-reduction-Rd}]
For the last term in
\cref{eq:abstract-locally-integrable-transported-principal-splitting-Rd},
we prove the equivalent reduction
\begin{equation}\label{eq:abstract-locally-integrable-transported-principal-reduction-Rd}
\begin{aligned}
    \mathpzc t_{h,k}^{\mathrm{tr},\mathrm I}\pare\eta\bra\omega
    ={}&
    \int_{\bR^{\mathtt d}}
    \chi\pare z
    \frac{h\pare{\hat z}}{\av z^{n+1}}
    \OpBW{\omega_z k'\pare{\Delta_z\eta}}
    \delta_z\eta\dd z
    \\
    &+
    r^\bR_{\bra{\varrho,\geq\ell}}\bra\omega\pare\eta\eta.
\end{aligned}
\end{equation}
The two symbols in
\cref{eq:abstract-transported-composition-defect} are independent of
\(\xi\). Hence every non-principal term in their Bony--Weyl product
vanishes, and the composition defect is arbitrarily smoothing by
\cref{prop:composition_paradiff,lem:composition_rules_z}.
Since its homogeneous factor has local order \(-n>-\mathtt d\),
\cref{lem:integration_KR} applies directly, without a frozen
subtraction. This proves
\cref{eq:abstract-locally-integrable-transported-principal-reduction-Rd};
the composition defect contributes only to the smoothing remainder in
\cref{eq:abstract-locally-integrable-transported-principal-almost-Rd}.
\end{proofsubpart}

Combining
\cref{eq:abstract-locally-integrable-transported-principal-reduction-Rd,%
eq:abstract-locally-integrable-transported-freezing-Rd,%
eq:abstract-locally-integrable-transported-freezing-error-Rd}
proves \cref{eq:abstract-locally-integrable-transported-principal-almost-Rd}.
\end{proofpart}

\begin{proofpart}[Proof of \cref{eq:abstract-locally-integrable-transported-Bony-remainder-final-Rd}]
We prove
\begin{equation}\label{eq:abstract-locally-integrable-transported-Bony-remainder-final-Rd}
    \mathpzc t_{h,k}^{\mathrm{tr},\mathrm {II}}\pare\eta\bra\omega
    =
    r^\bR_{\bra{\varrho,\geq\ell}}\bra\omega\pare\eta\eta.
\end{equation}
By \cref{lem paralin}, choosing a smoothing gain larger than the
prescribed \(\varrho\), we have
\[
    \mathpzc q_k\in\cR^\varrho_{\geq1}.
\]
Moreover, \(\mathpzc q_k\pare u\bra u\) vanishes to order at least
\(\ell\), while its differential vanishes to order at least
\(\ell-1\).
Since the homogeneous factor has local order \(-n>-\mathtt d\),
\cref{rem:smoothBONYremainder,lem:integration_KR} apply directly to the
second line of
\cref{eq:abstract-locally-integrable-transported-two-branches-Rd} and
give
\cref{eq:abstract-locally-integrable-transported-Bony-remainder-final-Rd}.
If \(\ell\geq2\), the range vanishes to order at least \(\ell-1\).
\end{proofpart}

\begin{proofpart}[Proof of \cref{eq:abstract-locally-integrable-transported-final-Rd}]
Combining
\cref{eq:abstract-locally-integrable-transported-two-branches-Rd,%
eq:abstract-locally-integrable-transported-principal-almost-Rd,%
eq:abstract-locally-integrable-transported-principal-symbol-reduction-Rd,%
eq:abstract-locally-integrable-transported-Bony-remainder-final-Rd}
proves \cref{eq:abstract-locally-integrable-transported-final-Rd}.
The symbol classes and their vanishing orders follow from
\cref{eq:abstract-locally-integrable-transported-freezing-Rd,%
eq:abstract-locally-integrable-transported-symbol-bound-Rd}
and the reductions above.
Conjugating
\cref{eq:abstract-locally-integrable-transported-symbol-Rd}
and replacing \(\xi\) with \(-\xi\) gives the principal
real-to-real symmetry. For the other terms, the same property follows
from \cref{eq:int-opbw-reality,eq:compositionsymbol1,lem:integration_KR}.
\end{proofpart}
\end{proof}

\subsection{Proof of \cref{thm:abstract-locally-integrable-paralinearization-Rd}}
\begin{proof}
For smooth inputs, combining
\cref{eq:abstract-locally-integrable-Bony-decomposition-Rd,%
eq:abstract-locally-integrable-quasilinear-final-Rd,%
eq:abstract-locally-integrable-transported-final-Rd,%
eq:abstract-locally-integrable-Bony-remainder-final-Rd}
gives
{
\begin{equation}\label{eq:abstract-locally-integrable-synthesis-before-relabeling-Rd}
\begin{aligned}
    \mathpzc t_{h,k}^{\pare n}\pare\eta\bra\omega
    ={}&
    \OpBW{
        \omega\pare x
        a_1^{h,k;\eta}\pare{\nabla\eta,\xi}
        +
        a_{\bra{0,\geq\ell}}\bra\omega
        +f_{\geq\ell}^{\tail,\bR}\bra\omega\pare\eta
    }\eta
    \\
    &+
    \OpBW{
        a_{-\beta}^{h,k;\omega}\pare{\nabla\eta,\xi}
        +
        a_{\bra{-\beta-2,\geq\ell}}
    }\omega
    \\
    &+
    \pare{
        r^\bR_{\bra{\varrho,\geq\ell}}\bra\omega\pare\eta
        +r^{\tail,\bR}_{\bra{\varrho,\geq\ell}}\bra\omega\pare\eta
    }\eta
    +
    r^\bR_{\bra{\varrho,\geq\ell}}\pare\eta\omega.
\end{aligned}
\end{equation}
}
The first principal term is identified in
\cref{eq:abstract-locally-integrable-transported-principal-symbol-reduction-Rd}.
The second follows from
\cref{eq:abstract-locally-integrable-leading-branch-Rd}, while
\cref{eq:abstract-locally-integrable-subprincipal-cancellation-Rd}
removes the intermediate order \(-\beta-1\).
The symbol and remainder classes, their vanishing orders, and their
real-to-real properties follow from
\cref{lem:abstract-locally-integrable-quasilinear-contribution-Rd,%
lem:abstract-locally-integrable-transported-contribution-Rd,%
lem:abstract-locally-integrable-localized-Bony-remainder-Rd,%
lem:abstract-locally-integrable-Bony-reduction-Rd}.

{
Absorbing the tail function in the generic order-zero coefficient and the
tail remainders in the corresponding smoothing components, and using
\cref{rem:linear-parameter-vanishing-order}, the bracketed terms have total
vanishing order at least $\ell$. Collecting the terms according to
\cref{eq:abstract-locally-integrable-vector-symbol-Rd,%
eq:abstract-locally-integrable-linear-smoothing-remainder-Rd}
then gives \cref{eq:abstract-locally-integrable-paralinearization-Rd}.
}
Choose a common baseline \(s_0\), finite loss \(\mu\), and radius
\(r_0\) for the preceding reductions. Smooth approximation and their
tame estimates extend the identity to the Sobolev inputs of the theorem.
\end{proof}

	\begin{footnotesize}
		\bibliography{references}
		\bibliographystyle{plain}
	\end{footnotesize}

\end{document}